\documentclass{memo-l}

\usepackage{exscale,amssymb,color,graphicx,accents,mathrsfs,mathtools,hyperref,comment}
\usepackage{enumitem}
\input xy
\xyoption{all}

\setlist[itemize]{leftmargin=1cm}

\newtheorem{theorem}{Theorem}[chapter]
\newtheorem{lemma}[theorem]{Lemma}
\newtheorem{proposition}[theorem]{Proposition}
\newtheorem{corollary}[theorem]{Corollary}

\theoremstyle{definition}
\newtheorem{definition}[theorem]{Definition}
\newtheorem{construction}[theorem]{Construction}
\newtheorem{example}[theorem]{Example}
\newtheorem{problem}[theorem]{Problem}

\theoremstyle{remark}
\newtheorem{remark}[theorem]{Remark}

\newtheorem{conjecture}[theorem]{Conjecture}

\numberwithin{section}{chapter}
\numberwithin{equation}{chapter}

\numberwithin{figure}{chapter}

\newcommand{\abs}[1]{\left\lvert #1 \right\rvert}

\makeatletter
 \DeclareRobustCommand{\checkarg}{\@ifnextchar[{\@witharg}{}}
 \DeclareRobustCommand{\@witharg}[1][]{\ensuremath{\left(#1\right)}}
 
 \DeclareRobustCommand{\scaleGen}[1]{\@ifnextchar[{\@scalewithargs{#1}}{\odot^{}_{#1}}}
 \def\@scalewithargs#1[#2][#3]{#2 \odot^{}_{#1} #3}
\makeatother

\newcommand{\IPspace}{\mathcal{I}}
\newcommand{\dI}{d_{\IPspace}}

\newcommand{\Exc}{\mathcal{E}}

\newcommand{\mBxc}{\nu_{\texttt{BESQ}}}	

\def\len{\textnormal{len}}			
\def\life{\zeta}					
\def\dis{\textnormal{dis}}			
\def\IPmag#1{\left\|\vphantom{I}#1\right\|}		

\def\skewer{\textsc{skewer}}		
\def\skewerP{\widebar{\skewer}}		

\def\Dirac#1{\delta\left( #1 \right)}

\def\ShiftRestrict#1#2{#1\big|^{\from}_{#2}} 
\def\shiftrestrict#1#2{#1|^{\from}_{#2}}

\def\Restrict#1#2{#1\big|_{#2}}
\def\restrict#1#2{#1|_{#2}}

\def\Concat{ \mathop{ \raisebox{-2pt}{\Huge$\star$} } }

\def\concat{\star}

\def\BR{\mathbb{R}}				
\def\BN{\mathbb{N}}				
\def\BQ{\mathbb{Q}}				

\def\Leb{\textnormal{Leb}}		
\def\downto{\downarrow}

\def\from{\leftarrow}

\def\cf{\mathbf{1}}				
\def\Pr{\mathbf{P}}				
\def\BQ{\mathbb{Q}}				
\def\EV{\mathbf{E}}				

\newcommand{\distribfont}[1]{\ensuremath{\mathtt{#1}}}

\def\ExpDist{\distribfont{Exponential}\checkarg}
\def\GammaDist{\distribfont{Gamma}\checkarg}

\def\BetaDist{\distribfont{Beta}\checkarg}

\def\PoiDir{\distribfont{PD}\checkarg}

\def\PRM{\distribfont{PRM}\checkarg}

\def\Stable{\distribfont{Stable}\checkarg}
\def\BESQ{\distribfont{BESQ}\checkarg}
\def\PDIP{\distribfont{PDIP}\checkarg}

\DeclareRobustCommand{\CRP}{\texttt{CRP}\checkarg}
\DeclareRobustCommand{\oCRP}{\texttt{oCRP}\checkarg}
\def\oCRPAT{\oCRP\ensuremath{(\alpha,\theta)}}
\def\CRPAT{\CRP\ensuremath{(\alpha,\theta)}}

\def\cadlag{c\`adl\`ag}

\newcommand{\cev}[1]{\accentset{\leftharpoonup}{#1}}

\makeatletter
\let\save@mathaccent\mathaccent
\newcommand*\if@single[3]{%
  \setbox0\hbox{${\mathaccent"0362{#1}}^H$}%
  \setbox2\hbox{${\mathaccent"0362{\kern0pt#1}}^H$}%
  \ifdim\ht0=\ht2 #3\else #2\fi
  }
\newcommand*\rel@kern[1]{\kern#1\dimexpr\macc@kerna}
\newcommand{\widebar}{}
\DeclareRobustCommand*\widebar[1]{\@ifnextchar^{\wide@bar{#1}{0}}{\wide@bar{#1}{1}}}
\newcommand*\wide@bar[2]{\if@single{#1}{\wide@bar@{#1}{#2}{1}}{\wide@bar@{#1}{#2}{2}}}
\newcommand*\wide@bar@[3]{%
  \begingroup
  \def\mathaccent##1##2{%
    \let\mathaccent\save@mathaccent
    \if#32 \let\macc@nucleus\first@char \fi
    \setbox\z@\hbox{$\macc@style{\macc@nucleus}_{}$}%
    \setbox\tw@\hbox{$\macc@style{\macc@nucleus}{}_{}$}%
    \dimen@\wd\tw@
    \advance\dimen@-\wd\z@
    \divide\dimen@ 3
    \@tempdima\wd\tw@
    \advance\@tempdima-\scriptspace
    \divide\@tempdima 10
    \advance\dimen@-\@tempdima
    \ifdim\dimen@>\z@ \dimen@0pt\fi
    \rel@kern{0.6}\kern-\dimen@
    \if#31
      \overline{\rel@kern{-0.6}\kern\dimen@\macc@nucleus\rel@kern{0.4}\kern\dimen@}%
      \advance\dimen@0.4\dimexpr\macc@kerna
      \let\final@kern#2%
      \ifdim\dimen@<\z@ \let\final@kern1\fi
      \if\final@kern1 \kern-\dimen@\fi
    \else
      \overline{\rel@kern{-0.6}\kern\dimen@#1}%
    \fi
  }%
  \macc@depth\@ne
  \let\math@bgroup\@empty \let\math@egroup\macc@set@skewchar
  \mathsurround\z@ \frozen@everymath{\mathgroup\macc@group\relax}%
  \macc@set@skewchar\relax
  \let\mathaccentV\macc@nested@a
  \if#31
    \macc@nested@a\relax111{#1}%
  \else
    \def\gobble@till@marker##1\endmarker{}%
    \futurelet\first@char\gobble@till@marker#1\endmarker
    \ifcat\noexpand\first@char A\else
      \def\first@char{}%
    \fi
    \macc@nested@a\relax111{\first@char}%
  \fi
  \endgroup
}
\makeatother

\setlist[enumerate,1]{leftmargin=1cm}

\newcommand{\wt}[1]{\widetilde{#1}}
\newcommand{\wh}[1]{\widehat{#1}}

\newcommand{\bD}{\mathbb{D}}
\newcommand{\bE}{\mathbb{E}}

\newcommand{\bN}{\mathbb{N}}
\newcommand{\bP}{\mathbb{P}}
\newcommand{\bR}{\mathbb{R}}
\newcommand{\bT}{\mathbb{T}}

\newcommand{\cB}{\mathcal{B}}

\newcommand{\cD}{\mathcal{D}}
\newcommand{\cE}{\mathcal{E}}
\newcommand{\cF}{\mathcal{F}}
\newcommand{\cG}{\mathcal{G}}

\newcommand{\cI}{\mathcal{I}}
\newcommand{\cJ}{\mathcal{J}}

\newcommand{\cL}{\mathcal{L}}
\newcommand{\cM}{\mathcal{M}}
\newcommand{\cN}{\mathcal{N}}

\newcommand{\cR}{\mathcal{R}}
\newcommand{\cS}{\mathcal{S}}
\newcommand{\cT}{\mathcal{T}}
\newcommand{\cU}{\mathcal{U}}
\newcommand{\cV}{\mathcal{V}}
\newcommand{\cW}{\mathcal{W}}

\newcommand{\sD}{\mathscr{D}}

\newcommand{\ff}{\mathbf{f}}
\newcommand{\fn}{\mathbf{n}}
\newcommand{\fm}{\mathbf{m}}

\newcommand{\ft}{\mathbf{t}}

\newcommand{\fN}{\mathbf{N}}
\newcommand{\fP}{\mathbf{P}}
\newcommand{\fT}{\mathbf{T}}

\newcommand{\fX}{\mathbf{X}}

\newcommand{\parent}[1]{\accentset{\leftarrow}{#1}}
\newcommand{\longparent}[1]{\accentset{\longleftarrow}{#1}}
\DeclareRobustCommand{\Ast}[1]{\accentset{*}{#1}}

\newcommand{\wi}{{\tilde{\textit{\i}}}}
\newcommand{\wj}{{\tilde{\textit{\j}}}}

\newcommand{\TInt}{\bT^{\textnormal{int}}}
\newcommand{\TShape}{\bT^{\textnormal{shape}}}
\newcommand{\bTInt}{\widebar{\bT}^{\textnormal{int}}}
\newcommand{\tdTInt}{\widetilde{\bT}^{\textnormal{int}}}

\newcommand{\block}{\textsc{block}}

\newcommand{\wbD}{\widebar{D}}
\newcommand{\wbT}{\widebar{\cT}}

\newcommand{\bTMarkk}{\Ast{\widebar{\bT}}^{\textnormal{int}}_k}
\newcommand{\tdTMarkk}{\Ast{\widetilde{\bT}}^{\textnormal{int}}_k}
\newcommand{\TMarkk}{\Ast\bT_{k}}

\newcommand{\besq}{{\tt BESQ}}

\newcommand{\skewerbar}{\ensuremath{\overline{\normalfont\textsc{skewer}}}}
\newcommand{\clade}{\ensuremath{\normalfont\textsc{clade}}}

\newcommand{\StableA}{\ensuremath{\distribfont{Stable}\big(\frac32\big)}}
\newcommand{\besqA}{\ensuremath{\distribfont{BESQ}(-1)}}

\newcommand{\PDIPAT}{\ensuremath{\distribfont{PDIP}\big(\frac12,\frac12\big)}}

\makeindex

\begin{document}

\frontmatter

\title{The Aldous diffusion: \\ a stationary evolution\\ of the Brownian CRT}


\author{Noah Forman}
\address{Department of Mathematics and Statistics, McMaster University, 1280 Main St W, Hamilton, ON L8S 4L8, Canada}
\curraddr{}
\email{noah.forman@gmail.com}
\thanks{}
\author{Soumik Pal}
\address{Department of Mathematics, University of Washington, Seattle, WA 98195, USA}
\curraddr{}
\email{soumikpal@gmail.com}
\thanks{}
\author{Douglas Rizzolo}
\address{Department of Mathematical Sciences, University of Delaware, Newark, DE 19716, USA}
\curraddr{}
\email{drizzolo@udel.edu}
\thanks{}
\author{Matthias Winkel}
\address{Department of Statistics, University of Oxford, 24--29 St Giles', Oxford OX1 3LB, UK}
\curraddr{}
\email{winkel@stats.ox.ac.uk}
\thanks{}

\date{\today}

\subjclass[2020]{Primary 60J80, 60J25; Secondary 60J60, 60G18, 60C05}

\keywords{Brownian CRT, weighted $\mathbb{R}$-tree, Aldous diffusion, interval partition, Chinese restaurant process, 
Wright--Fisher diffusion, stable process, Poisson--Dirichlet distribution, scaling limit, intertwining}


\begin{abstract}
Motivated by a down-up Markov chain on cladograms, David Aldous conjectured in 1999 that there exists a ``diffusion on continuum trees'' whose mass partitions at any finite number of branch points evolve as certain Wright--Fisher diffusions with some negative mutation rates, until some branch point disappears. Building on previous work on interval-partition-valued processes, we construct this conjectured process via a consistent system of stationary evolutions of binary trees with $k$ labeled leaves and edges decorated with interval partitions. The interval partitions are scaled Poisson--Dirichlet interval partitions whose interval lengths record subtree masses. They also 
possess a diversity property that captures certain distances in the continuum tree. Continuously evolving diversities give access to continuously evolving continuum tree distances. 

The pathwise construction allows us to study this ``Aldous diffusion'' in the Gromov--Hausdorff--Prokhorov space of rooted, weighted $\mathbb{R}$-trees. We establish the (simple) Markov property and path-continuity. The Aldous diffusion is stationary with the distribution of the Brownian continuum random tree. While the Brownian continuum random tree is binary almost surely, we show that there is a dense null set of exceptional times when the Aldous diffusion has a ternary branch point, and this set includes stopping times at which the strong Markov property fails. 

Our construction relates to the two-parameter Chinese restaurant process, branching processes, and stable L\'evy processes, among other connections. Wright--Fisher diffusions and the aforementioned processes of Poisson--Dirichlet interval partitions arise as interesting projections of the Aldous diffusion.

Finally, the Aldous diffusion with its consistent system of evolving binary trees embedded, allows us to embed Aldous's stationary down-up Markov chain on 
cladograms in the Aldous diffusion and hence address a related conjecture by David Aldous by establishing a scaling limit theorem.
\end{abstract}

\maketitle

\tableofcontents


\mainmatter

\chapter{Introduction}\label{ch:intro}
Tree-valued dynamics arise in applications in computer science \cite{PropWils98}, machine learning \cite{AdamGhahJord10,BleiGrifJord10,BrodGram11,GramacyLee08,Neal03}, and phylogenetics \cite{BarkerThesis,DinhEtAl,DrumRamb07,HuelRonq01,LargetSimon,WhidMats15}, often in the context of Markov chain Monte Carlo inference. The immense size of phylogenetic trees has motivated a growing literature on asymptotic properties of such Markov chains, notably mixing times \cite{Aldous00,Caraceni20,EppsFris22,McShineTetali,Schweinsberg02} and continuum analogs and scaling limits; see the lecture notes of Evans \cite{EvansStFlour} and Zambotti \cite{Zambotti17}.

The purpose of this memoir is to construct and study a path-continuous Markov process on a space of continuum trees whose existence was conjectured by David Aldous \cite{ADProb2, AldousDiffusionProblem} in 1999, and to establish a scaling limit theorem also conjectured by David Aldous \cite{ADProb2}, which we believe is an instance of an invariance principle that would frame this new process as a universal limit object. We call this continuum-tree-valued limiting process ``the Aldous diffusion.''

Aldous's original motivation for the conjecture and a reason for its continuing significance is that while there are a number of different Markov chains in common use, the one suggested by Aldous was among those that seemed most accessible to a fuller asymptotic analysis, as Aldous had illustrated by mixing time calculations \cite{Aldous00} and some further observations that led to his conjecture \cite{ADProb2,AldousDiffusionProblem}. In the meantime, several continuum-tree-valued processes relating to different Markov chains have been studied by a variety of authors \cite{EtheLabb15,EPW,EW,Zambotti01} and a variant of Aldous's conjecture has been resolved by L\"ohr, Mytnik and Winter \cite{LohrMytnWint20}.

This work develops and explores relationships between a range of classical stochastic processes, 
including the two-parameter Chinese restaurant process, branching processes, squared Bessel processes, Wright--Fisher diffusions, and stable L\'evy processes and subordinators. 
These connections have already borne additional fruit: en route to resolving Aldous's conjecture, we have resolved a 2009 conjecture of Feng and Sun \cite{FengBook,FengSun10} on (measure-valued) Fleming--Viot processes with Poisson--Dirichlet stationary distributions \cite{Paper1-3,PaperFV,ShiWinkel-2}.

We use this introduction not only to state the main results that solve the conjectures, but also to point out connections to related stochastic processes and complex random structures, and we make observations that generalize to other settings or shed light on the difficulties encountered in other approaches. 

This document builds upon the authors' previous work on interval-partition- and (combinatorial) tree-valued processes \cite{Paper1-1,Paper1-0,Paper2,Paper1-2}, supersedes two unpublished preprints \cite{Paper4,Paper3}, and develops new material.

\section{The Aldous diffusion conjecture}\label{sec:adconj}
\begin{figure}\centering\vspace{-0.4cm}
 \scalebox{.88}{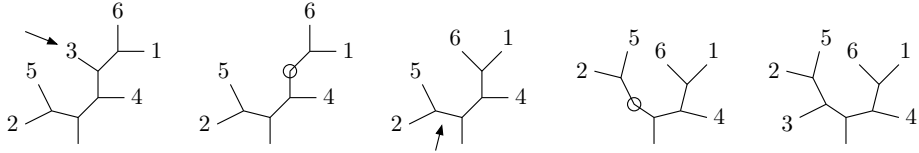}
 \caption{From left to right, one Aldous down-up move.\label{fig:AC_move}}
\end{figure}

The Aldous chain is a Markov chain on the space of rooted binary trees with $n$ labeled leaves. Each transition of the {Aldous chain}, called a down-up move, has two steps. In the down-move a uniform random leaf is deleted and its parent branch point is contracted away; in the up-move a uniform random edge is selected, a branch point is inserted into the middle of the edge, and the leaf is reattached at that point. See Figure \ref{fig:AC_move}. It is not difficult to prove that this chain is stationary with the uniform distribution on such trees, indeed reversible. Aldous \cite{Aldous00} studied the analog of this chain on unrooted trees.

Suppose that for each $n\ge1$: $T_n$ is a uniform random rooted binary tree with $n$ labeled leaves, as in the first panel in Figure \ref{fig:AC_move}; $\rho_n$ is its root; $d_n$ is the graph distance metric on $T_n$; and $\mu_n$ is the uniform probability distribution on the leaves of $T_n$. Then the sequence $(T_n,d_n / \sqrt{n},\rho_n,\mu_n)$, $n\ge1$, converges in distribution \cite{CuriHaas13,HMPW} under the rooted Gromov--Hausdorff--Prokhorov metric \cite{M09}, which we will discuss in Section \ref{sec:Rtrees}. The limit $(\cT,d,\rho,\mu)$ is the \emph{Brownian continuum random tree} (\emph{BCRT}) \cite{AldousCRT1}, which we discuss in Section \ref{sec:BCRT}.

Aldous \cite{ADProb2, AldousDiffusionProblem} suggested taking a diffusive limit of the Markov chain as $n\rightarrow\infty$, yielding a diffusion on some space of continuum trees as the limiting object. He also did some calculations showing how some tree statistics of the limiting process should evolve.  Specifically, he considered partitioning the initial tree around some subset of its branch points and studying the fluctuating leaf counts in each connected component, as the tree evolves according to the Aldous chain. This induces a Markov chain until one of the partitioning branch points is contracted away. Aldous observed that there are two types of components: ``internal'' ones adjacent to two branch points and ``external'' ones adjacent to one branch point. External components with $n_i\ge 1$ leaves have $2n_i-1$ edges, while internal components with $n_i\ge 0$ leaves have $2n_i+1$ edges. In our setting of rooted trees, the component containing the root, which is not considered a leaf (and hence cannot be deleted in a down-move), behaves like an internal component. Aldous further conjectured that the induced Markov chain has a space-time scaling limit as $n\to\infty$: a Wright--Fisher-like multi-dimensional diffusion process on the simplex, run until some component vanishes. Interestingly, external components have negative Wright--Fisher ``mutation rates.''

Wright--Fisher diffusions are Markov processes in the $d$-dimensional simplex $\Delta_d=\{(w_1,\ldots,w_d)\in[0,1]^d\colon\sum_{1\le i\le d}w_i=1\}$. Specifically, let $\theta_1,\ldots,\theta_d$ be
non-negative and/or negative real parameters. We consider the infinitesimal generator
\begin{equation}\label{eq:WFgen}
  \cG=2\sum_{1\le i\le d}w_i\frac{\partial^2}{\partial w_i^2}-2\sum_{1\le i,j\le d}w_iw_j\frac{\partial^2}{\partial w_i\partial w_j}-2\sum_{1\le i\le d}(\theta_+w_i-\theta_i)\frac{\partial}{\partial w_i},
\end{equation}
where $\theta_+=\sum_{1\le i\le d}\theta_i$. In population genetics, diffusions with such generators (up to a constant factor) arise with mutation parameters $\theta_1,\ldots,\theta_d\ge 0$, when the boundaries where a coordinate vanishes are reflecting (or absorbing when the corresponding parameter is zero). See e.g.\ Ethier and Kurtz \cite{EthKurtzBook}. Pal \cite{Pal13} extended this to negative parameters 
and constructed diffusions ${\tt WF}(\theta_1,\ldots,\theta_d)$ that are stopped when a coordinate with a negative parameter vanishes (or the diffusion may continue on a lower-dimensional simplex). Pal was motivated by Aldous's observation about the induced Markov chain on masses $\big(X^{(n)}_m(1),\ldots,X^{(n)}_m(2k-1)\big)$, 
$0\le m\le D^{(n)}$, recording the proportions of leaves in the components around a finite number $k-1$ of branch points: as $n\to\infty$, we have
\begin{equation}\label{eq:3massconv}
  \Big(\Big(X^{(n)}_{\lfloor n^2 t\rfloor\wedge D^{(n)}}(1),\ldots,X^{(n)}_{\lfloor n^2t\rfloor\wedge D^{(n)}}(2k\!-\!1)\Big),\,t\ge0\Big)\stackrel{d}{\longrightarrow}{\tt WF}(\theta_1,\ldots,\theta_{2k-1}),\!\!
\end{equation}
where $(\theta_1,\ldots,\theta_{2k-1})$ is a vector of $\frac12$ for each internal component between two branch points (or the root) and $-\frac12$ for each external component adjacent to a single branch point.

\begin{conjecture}[Aldous, 1999 \cite{ADProb2, AldousDiffusionProblem}]
\label{conj:Aldous}
  There exists a ``diffusion on continuum trees'' that is stationary with the law of the BCRT $(\cT,d,\rho,\mu)$ and for which the evolution of $\mu$-masses around some finite subsets of its branch points are Wright--Fisher diffusions distributed like the limits in \eqref{eq:3massconv}.  Moreover, this process is the diffusive limit of the Aldous chain.
\end{conjecture}

We call the conjectured process the \emph{Aldous diffusion}.  
The present work resolves Aldous's conjecture by constructing such a process and showing that it is the scaling limit, in the sense of finite-dimensional distributions, of the Aldous chain run according to a Poisson clock. The Aldous diffusion has continuous paths and the simple Markov property but, surprisingly, not the strong Markov property.

\subsection*{Motivation for our approach}

\begin{figure}
\begin{center}
\begin{minipage}{.5\linewidth}
\hspace{0.8cm}\scalebox{.20}{ \includegraphics{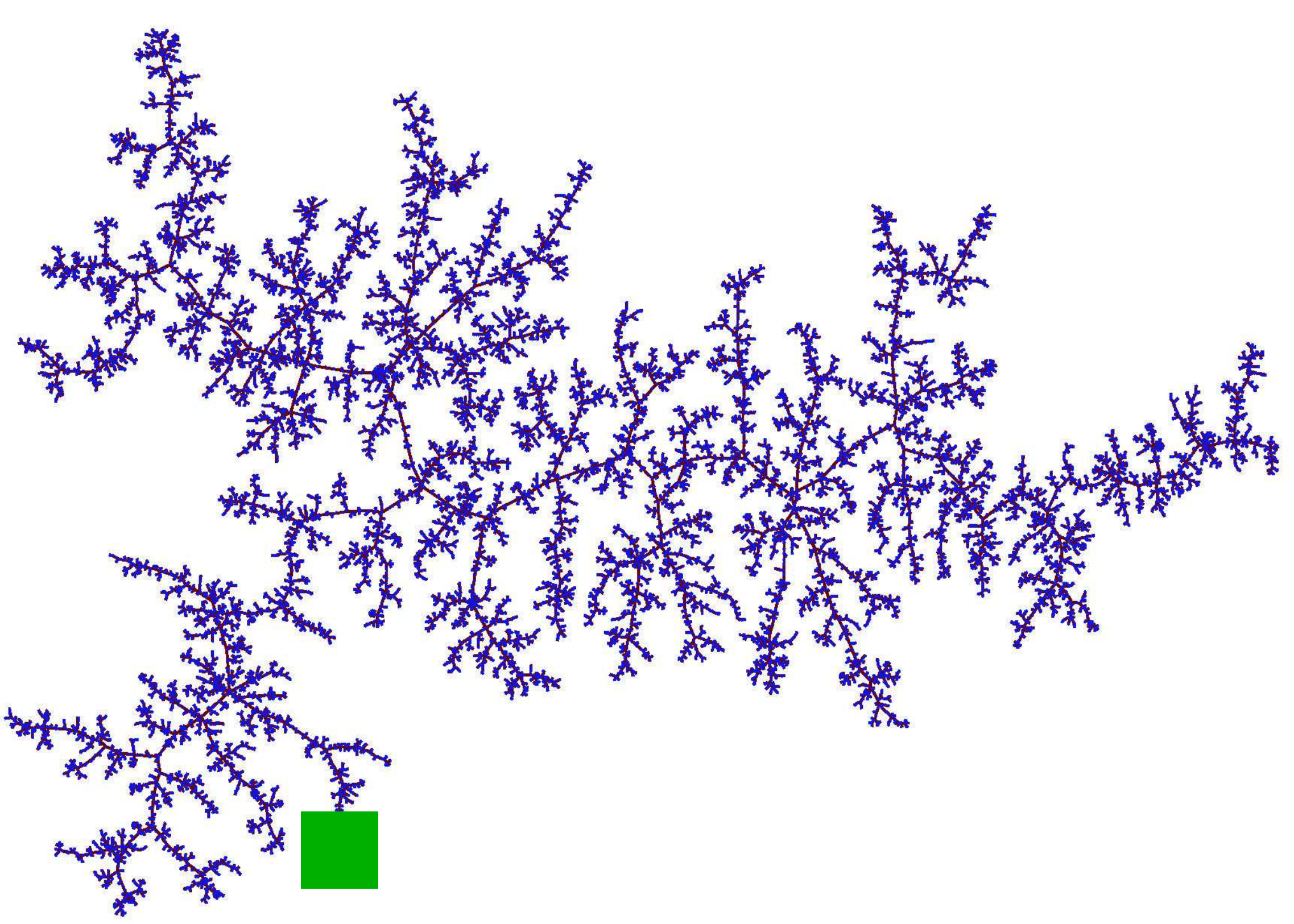}}
\end{minipage}%
\begin{minipage}{.5\linewidth}
\scalebox{1}{ \includegraphics{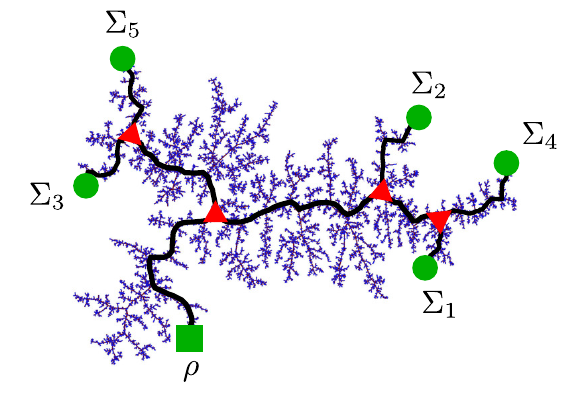}}
\end{minipage}
\end{center}
 \caption{Left: Simulation of a BCRT $(\cT,d,\rho,\mu)$ courtesy of Igor Kortchemski, with a root vertex (green square). Right: $\cR_5^+$ is the tree spanned by the root $\rho$ and leaves
  $\Sigma_1,\ldots,\Sigma_5$, here depicted as black branches with (red) triangles as branch points.\label{fig:5treeEdges}}
\end{figure}

 Our construction of the process on continuum trees is based on a dynamic variant of Aldous's original construction of the BCRT.  In \cite{AldousCRT3}, Aldous constructed the BCRT $(\cT,d,\rho,\mu)$ as the limit, in an appropriate metric space, of a consistent sequence $(\cR_k^+,\,k\geq 1)$  of rooted, leaf-labeled, binary trees with edge lengths such that the tree shape of $\cR_k^+$ is uniformly distributed on rooted binary trees with $k$ labeled leaves and, conditionally given the shape of  $\cR_k^+$, the lengths of its $2k-1$ edges have the joint density on $(0,\infty)^{2k-1}$
\[ f(x_1,\dots, x_{2k-1}) =\left(\prod_{i=1}^{k-1}(2i-1)\right) s\exp(-s^2/2) \quad \textrm{where } s=\sum_{i=1}^{2k-1} x_i.\]
The consistency of $(\cR_k^+,\,k\geq 1)$ means that $\cR_k^+$ can be obtained from $\cR_{k+1}^+$ by removing the leaf labeled $k+1$ along with the branch connecting it to the rest of the tree.  The intuitive idea is that if one takes an i.i.d sequence of leaves in $\cT$, then $\cR_k^+$ is the subtree of $\cT$ spanned by the first $k$ leaves and the root; see Figure \ref{fig:5treeEdges}.  This becomes formally true once $\cT$ has been constructed. 
A natural aim would be to develop a dynamic version of this by constructing a consistent family $(\cR_k^+(t),\, t\geq 0)$, ${k\geq 1}$, of evolving trees such that $\cR_k^+(t)$ evolves as a sampled subtree should evolve in the conjectured process and then to construct the process as the limit of this sequence. Similar projective consistency has been used previously to construct the limit of the root growth with re-grafting process \cite{EPW}.  However, in the present case, several novel challenges arise.

The first challenge arises because $\cR_k^+$ is spanned by randomly sampled leaves and the root, but 
it is precisely the leaves that are being moved 
in each step of the Aldous chain.  From our discussion of the discrete mass split around branch points above, we see that the appropriate scaling is to have $n^2$ steps per unit time.  But there are only $n$ leaves, so all leaves have been moved approximately after $n\log(n)$ steps.  Thus, in the limiting process, the leaves are moving too quickly for a dynamic version of $\cR_k^+$ to behave nicely (or even be well-defined).

Our discussion of the evolving mass split around branch points suggests that instead of a dynamic version of $\cR_k^+$, we consider a dynamic version of the subtree $\cR_k$ spanned by the branch points and root of $\cR_k^+$. 
The leaves move rapidly, but the components of the tree around branch points have identities that are stable over time until the branch point is contracted away. We adopt this approach.


The second challenge is that of describing the evolving lengths of edges between branch points.  As discussed in \cite{LohrMytnWint20}, distances should not have finite quadratic variation, making them difficult to describe using classical techniques like martingale problems or stochastic differential equations.  Indeed, a consequence of our construction is that these distances evolve like the spatial variation of the local time process of a stable L\'evy process.  

One natural approach to recording distances would be to use the strings of beads developed in \cite{PitmWink09} and further studied in \cite{RW}.  We turn each edge of $\cR_k$ into a string of beads, which captures in an atomic measure a point mass at each location where a subtree of the limiting tree will branch off from the edge, recording in the atom size the mass of the subtree. The length of the edge can be recovered from the support of the measure.  We find, however, that as the locations of atoms evolve, it can (and will) happen that two atoms reflect off each other, instantaneously occupying the same location before bouncing off.  This behavior can be found at stopping times, so the evolution of the string of beads is not a strong Markov process: one cannot tell from the string of beads at that time that one of the atoms was two atoms an instant before and will again be two atoms an instant later.  This is related to the failure of the strong Markov property for the Aldous diffusion as well as the existence of ternary branch points.  

To avoid this problem, rather than working with strings of beads, we work with interval partitions, which are collections of disjoint open intervals whose lengths record subtree masses in the same manner as the point masses in the string of beads, ordered left-to-right by decreasing distance from the subtree to the root.  Instead of a dynamic version of $\cR_k$, we introduce interval partition trees $R_k$; see Section \ref{sec:BCRT} for a formal definition.  The combinatorial tree shape of $R_k$ is the shape of $\cR_k$.  To each edge we associate the interval partition of masses of subtrees that attach to that edge, as well as the masses of the subtrees above branch points that contain only a single sampled leaf; see Figure \ref{fig:B_k_tree_proj} in Section \ref{sec:BCRT} below.  With this choice of $R_k$, we are able to construct a consistent family $(R_k(t),\, t\geq 0)$, $k\geq 1$, of evolving interval partition trees by leveraging our recent progress on understanding evolving interval partitions \cite{Paper1-1,Paper1-2} and adapting, from the discrete to the continuum setting, our strategy in \cite{Paper2} for selecting new branch points when one disappears.  The Aldous diffusion is then defined as an appropriate $k\to\infty$ limit.  The process that we construct is a path-continuous simple Markov process, reversible with respect to the distribution of the BCRT, but it is not strongly Markovian, and thus it is not a diffusion in the strict sense.

%
%
%
%

%

\medskip

The remainder of this introduction is structured as follows. We give an introduction to the BCRT and its reduced $k$-trees in Section \ref{sec:BCRT}, before we state our main results in Section \ref{sec:mainresults}. We provide a literature overview in Section \ref{sec:related_work}. In Section \ref{sec:discrete_to_cts}, we explain in the more elementary setting of the Aldous chain the approach we develop in this memoir for reduced subtrees of continuum trees. We conclude the introduction by giving a chapter overview in Section \ref{sec:intro:overview}.

\section{The Brownian continuum random tree and its reduced $k$-trees}
\label{sec:BCRT}

Since its introduction in the 1990s, the BCRT has become a central object in probability theory, with 
a variety of representations \cite{AldousCRT1,EPW,EW,LeGall93,LohrWint21}, connections and ramifications. 
Connections include branching processes \cite{LeGall91}, coalescents and fragmentation processes \cite{AldoPitm98,Bertoin06}, Dirichlet and Poisson--Dirichlet distributions \cite{Aldous94,PitmWink09}. BCRTs serve as building blocks for scaling limits of Erd\H{o}s--R\'enyi random graphs in the critical window \cite{ABBG,Aldous97}, the Brownian map \cite{LeGallMier12,MarcMokk} and Liouville quantum gravity \cite{MatingTrees}. There is a large (universality) class of random trees that converge to the BCRT \cite{AldousCRT3,ArchShal23arxiv,BertMier13,CuriHaasKort15,HaasMier12,HaasStep21,MarcMier11}. 

We recall two classical constructions of the BCRT: the line-breaking construction and the construction from a Brownian excursion. We also discuss the representation of reduced
$k$-trees that is crucial for our construction of the Aldous diffusion.

\begin{definition} \label{def:real_tree}
 An \emph{$\BR$-tree} (\emph{real tree}) is a complete, separable metric space $(T,d)$ with the property that: (i) for each $x,y\in T$, there is a unique non-self-intersecting path in $T$ from $x$ to $y$, called a \emph{segment} $[\![x,y]\!]_{T}$, and (ii) each segment is isometric to a closed real interval.
 
 A \emph{rooted, weighted $\BR$-tree} is a quadruple $(T,d,\rho,\mu)$, where $(T,d)$ is an $\BR$-tree, $\rho\in T$ is a distinguished vertex called the \emph{root}, and $\mu$ is a probability distribution on the $\sigma$-algebra of Borel sets of $(T,d)$. In cases where the weight measure $\mu$ is supported on the leaves of $T$, it is also called \emph{leaf mass}. Here, a \em leaf \em is any $x\in T$ such that $T\setminus\{x\}$ is connected.
 
 For our purposes, a \emph{continuum random tree} (\emph{CRT}) is a random rooted, weighted $\BR$-tree with the a.s.\ properties that the weight measure is diffuse and supported on the leaves of the tree, and every neighborhood of every leaf has positive weight. See Section \ref{sec:Rtrees}.
\end{definition}

CRTs were introduced by Aldous in \cite{AldousCRT1, AldousCRT2, AldousCRT3}. The BCRT is the most famous example. In general, the properties stipulated for a CRT imply that 
the set of leaves has to be uncountable (in order to carry a diffuse weight measure). 


\subsection*{The line-breaking construction of the BCRT \cite{AldousCRT1}}%
One way to construct $\BR$-trees is by embedding them in the vector space of summable sequences of real numbers equipped with  the $\ell_1$-norm $\|(x_j,j\!\ge\!1)\| = \sum_{j\ge1}|x_j|$ and the associated $\ell_1$-distance $d_{\ell_1}(x,y)=\|x-y\|$. 
For $j\ge 1$, let $\mathbf{e}_j$ denote the sequence with 1 as its $j^{\text{th}}$ entry and all other entries 0. These sequences are of course independent of each other as vectors, though they do not form a basis for the space. 

Given a sequence of non-negative branch lengths $D_1,D_2,\ldots$ with suitable properties, we sequentially construct a random $\mathbb{R}$-tree as follows. Begin with $\cT_0 = \{\mathbf{0}\}$, the zero sequence in $\ell_1$. Then, iteratively for $j\ge1$, let $\mathbf{x}^{(j)}$ be a random point in $\cT_{j-1}$ sampled from the normalized length measure on the tree, and let
$$\cT_j = \cT_{j-1} \cup \big\{\mathbf{x}^{(j)} + tD_j\mathbf{e}_j\colon t\in [0,1]\big\}.$$
In other words, at each step we add a new branch of length $D_j$, extending in the $j^{\text{th}}$ coordinate direction from $\mathbf{x}^{(j)}$. Finally, we define $\cT$ to be the topological closure of the increasing union $\bigcup_{j\ge0}\cT_j$.

With suitably chosen random branch lengths, this construction gives rise to the random $\mathbb{R}$-tree $(\cT,d_{\ell_1},\mathbf{0})$ that we will further equip with a weight
measure to obtain the Brownian CRT $(\cT,d_{\ell_1},\mathbf{0},\mu)$. 
Specifically, consider a Poisson process $(N(t),\,t\ge0)$ with variable intensity $tdt$ on $[0,\infty)$; i.e.\ $\EV[N(v) - N(u)] = \int_{u}^vtdt$ for $0\le u<v$. 
Let $S_0 := 0$ and for $j\ge 1$ denote the $j^{\text{th}}$ arrival time by $S_j := \inf\{t\ge0\colon N(t)\ge j\}$. Then breaking the line $[0,\infty)$ at $S_j$, $j\ge 1$, produces branch
lengths $D_j = S_j - S_{j-1}$, $j\ge 1$. These branch lengths, plugged into the previous construction, produce Aldous's Brownian CRT $(\cT,d_{\ell_1},\mathbf{0},\mu)$, if we
further consider approximating uniform probability measures $\mu_j$ on $\cT_j$ that converge weakly almost surely as random measures on $\ell_1$ to a limiting measure $\mu$ that is supported by the leaves of $\cT$. See \cite[Theorem 3]{AldousCRT1}.

\subsection*{The construction of the BCRT from a Brownian excursion} 
\begin{figure}
 \centering
 \includegraphics[scale=1.5]{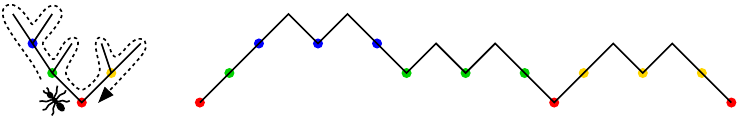}\\
 \includegraphics[scale=1.5]{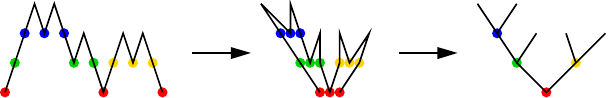}
 \caption{A rooted plane tree and corresponding Dyck path, with branch points of the tree color-coded to match corresponding points along the path. 
 Top: obtaining the path as the contour process of the tree. Bottom: recovering the tree by squeezing the path together laterally.
 \label{fig:Dyck_tree}}
\end{figure}

\cite{AldousCRT3,LeGall91,LeGall93} building on \cite{LeGall89,NeveuPitman89}. 
An unlabeled \emph{rooted plane tree} is a rooted tree in which the children of each non-leaf vertex are assigned a left-to-right order. For $n\ge1$, there is a classical bijection from the set of unlabeled rooted plane trees with $n+1$ vertices to the set of \emph{Dyck paths} of length $2n$: sequences $(a_i,\,i\in [0,2n])$ of non-negative integers with: (i) $a_0 = a_{2n} = 0$, and (ii) $|a_{i+1} - a_i| = 1$ for $0\le i\le 2n-1$. In the literature, the path associated with a tree is also known as its associated \emph{contour process} or \emph{Harris path} \cite{Dwass1975,Harris52,StanleyVol2}. Both directions of this bijection are illustrated in the top and bottom panels of Figure \ref{fig:Dyck_tree}.
%
%
%

\begin{figure}
 \centering
 \includegraphics[scale=1]{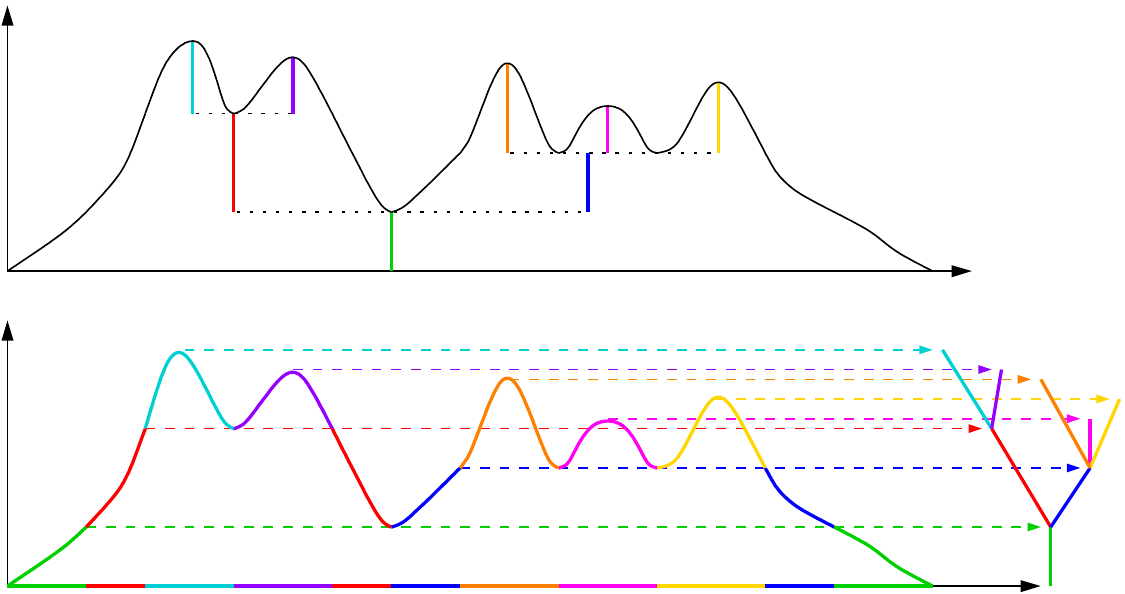}
 \caption{A continuous excursion function $H$ and the associated $\BR$-tree, $([0,1] / \sim_H , d_H)$. Top: The $\BR$-tree shown ``inscribed'' under the path. Horizontal dotted lines represent branch points. Bottom: The path and unit interval beneath color-coded to indicate which branch of the tree each segment corresponds to, with the correspondingly colored tree shown on the right.\label{fig:R_tree_from_exc}}
\end{figure}

This map can be extended from discrete Dyck paths to apply to continuous functions $H\colon [0,1]\to [0,\infty)$ with $H(0) = H(1) = 0$ and $H(t)\ge 0$ for $t\in (0,1)$. Given such a function, we define an associated pseudometric $d_H$ on $[0,1]$ by
\begin{equation}\label{eq:tree_from_path}
 d_H(a,b) = H(a) + H(b) - 2\min_{t\in [a,b]}H(t).
\end{equation}
This fails to be a proper metric at points $a<b$ where $H(a) = H(b) = \min_{t\in [a,b]}H(t)$, and we get $d_H(a,b) = 0$. Define an equivalence relation $\sim_H$ on $[0,1]$ by $a\sim_H b$ if and only if $d_H(a,b) = 0$. Then the quotient space $([0,1] /\!\!\sim_H , d_H , [0]_{\sim_H}, \Leb)$, where ${\rm Leb}$ denotes (the $\sim_H$-image of) Lebesgue measure on $[0,1]$, is a weighted $\BR$-tree, rooted at the equivalence class $[0]_{\sim_H}$ of $0$. See Figure \ref{fig:R_tree_from_exc}.

It is not hard to show that each branch point in the resulting tree corresponds to one or more local minima of $H$, depending on the degree of the branch point. Similarly, local maxima of $H$ all correspond to leaves, however, in general, not all leaves correspond to local maxima. In particular, assuming that $H$ is not locally constant anywhere, a point $r\in [0,1]$ corresponds to a leaf if and only if $H$ is neither locally non-increasing to the left of $r$ nor locally non-decreasing to the right:
\begin{equation}\label{eq:leaves_from_exc}
 \forall \delta > 0,\,\exists u\in (r-\delta,r),\,v\in(r,r+\delta)\text{\ such that\ }H(u)\!<\!H(r)\mbox{ and }H(v)\!<\!H(r).
\end{equation}

A standard Brownian excursion can be thought of informally as standard Brownian motion conditioned to: (i) escape upwards from 0 at time zero (instead of immediately having an accumulation of visits to 0) and (ii) then make its first subsequent return to 0 at time one. Formally, a standard Brownian excursion $(B^{\text{ex}}_r,\,r\in[0,1])$ can be constructed from a standard Brownian motion $(B_t,\,t\ge 0)$, as follows. For $s\ge 0$, let $g_s=\sup\{u\in[0,s]\colon B_u=0\}$ and $d_s=\inf\{u\ge s\colon B_u=0\}$ be the nearest zeroes before and after time $s$. As a consequence of the scaling property of Brownian motion, the excursion $(B_{g_s+t},\,t\in[0,d_s-g_s])$ straddling time $s$ is such that the
distribution of the excursion
\[
\left(\frac{1}{\sqrt{d_s-g_s}}|B_{g_s+(d_s-g_s)r}|,\ r\in[0,1]\right)
\]
scaled to unit time does not depend on $s$. This common distribution is the distribution of a standard Brownian excursion. See e.g.\ \cite{ItoMcKean1965,CSP,RevuzYor}.

The BCRT is the rooted, weighted $\mathbb{R}$-tree associated with $\big(2B^{\rm ex}_r,\,r\in [0,1]\big)$ via the map described around \eqref{eq:tree_from_path} and illustrated in Figure \ref{fig:R_tree_from_exc}.


\subsection*{Some properties of the BCRT}

We refer to Aldous \cite[Corollary 22]{AldousCRT3} for the property that these constructions are distributionally equivalent up to root- and weight-preserving isometry and to Pitman \cite[Chapter 7]{CSP} for a detailed discussion of this equivalence and its consequences. Relevant for us at this stage is that these two constructions grant easy access to some interesting properties of the BCRT. 

The line-breaking construction can be viewed in continuous time $t\ge 0$ as unit-rate continuous growth of unit length per unit time \cite{EPW}, with the $j^{\rm th}$ 
branch growing between times $S_{j-1}$ and $S_j$. Then the Poisson process with intensity $tdt$ and the uniformly random points $\mathbf{x}^{(j)}$, 
$j\ge 2$, make branch points appear on any existing branch at unit rate per unit length, and branch lengths $D_j$ tend to zero a.s.\ as $j\rightarrow\infty$. In 
particular, branch points and leaves are both dense in the BCRT, in the topological sense. 

From the excursion construction, we see that the BCRT is compact and its weight measure is indeed concentrated on the leaves, as the property \eqref{eq:leaves_from_exc} holds at ${\rm Leb}$-almost every $t\in [0,1]$ almost surely (for Brownian motion and hence) for the Brownian excursion. This construction also makes the BCRT inherit various instances of self-similarity from the Brownian excursion \cite{AlbG,Aldous94,AldoMierPitm04,AldoPitm98,Bertoin02,Bismut85,HM04}. 
Specifically, all excursions above a fixed level, or the two excursions adjacent to the mimimum between two independent uniform times, or all excursions above the past minimum process after (and above the future minimum before) a uniform random time are scaled independent Brownian excursions each encoding a scaled BCRT. 

Identifying $\mathbb{R}$-trees that are equal up to root- and weight-measure-preserving isometry classes is the continuum analog of identifying rooted combinatorial trees 
(or more general graphs) up to graph isomorphisms that preserve the root vertex, i.e.\ graphs that only differ in their vertex names. Indeed, the coordinates in the line-breaking construction are important in the construction but obscure the self-similarities of the resulting objects. On the other hand, the equivalence classes on $[0,1]$ encode additional 
planarity structure. It is instructive to explore the Aldous diffusion conjecture in conjunction with these and other representations of the BCRT, and we will do so in Section \ref{sec:openprob}.

\subsection*{Interval partitions and the Brownian reduced $k$-tree}
\label{sec:intro:B_k_tree}


To introduce our notion of a $k$-tree, we first require interval partitions, in the sense of \cite{Aldous85,GnedPitm05,PitmWink09}.

\begin{definition}\label{def:intro:IP}
 An \emph{interval partition} (IP) is a set $\beta$ of disjoint, open subintervals of some interval $[0,M]$, $M\ge 0$, that cover $[0,M]$ up to a Lebesgue null set. We refer to $M =: \|\beta\|$ as the \emph{mass} of $\beta$. The subintervals comprising the interval partition are called its \emph{blocks}. We refer to their lengths as \em block masses \em or \em sizes\em. 
\end{definition}

\begin{figure}
 \centering
 \includegraphics[height=1.5in, width=3.5in]{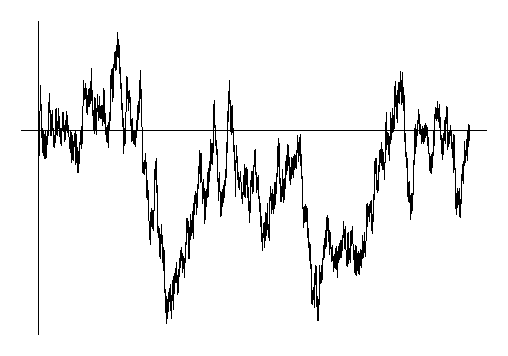}
 \caption{A standard Brownian bridge. In between any two excursions away from zero are infinitely many small excursions. The time intervals of these excursions form a $\PDIP\big(\frac12,\frac12\big)$.\label{fig:BBr}}
\end{figure}

Simple examples include finite partitions, such as $\{(0,2),(2,3)\}$, or infinite sequential partitions, such as $\{(0,\frac12),(\frac12,\frac34),(\frac34,\frac78),\ldots\}$. However, the blocks of a partition can also be ordered in a more complicated manner, such as the open intervals, called \emph{excursion intervals}, that comprise the complement of the zero set of a standard Brownian bridge; see Figure \ref{fig:BBr}. The left-to-right ordering of these intervals is isomorphic, as an ordered set, to $\BQ$. This interval partition is called a \emph{Poisson--Dirichlet$\big(\frac12,\frac12\big)$ interval partition}, or $\PDIP\big(\frac12,\frac12\big)$. In particular, the sequence of interval lengths, written in decreasing order, is Poisson--Dirichlet$\big(\frac12,\frac12\big)$-distributed \cite[Corollary 4.9]{CSP}; we will discuss this family of interval partitions at greater length in Section \ref{sec:IP}.

\begin{figure}[t!]\centering
  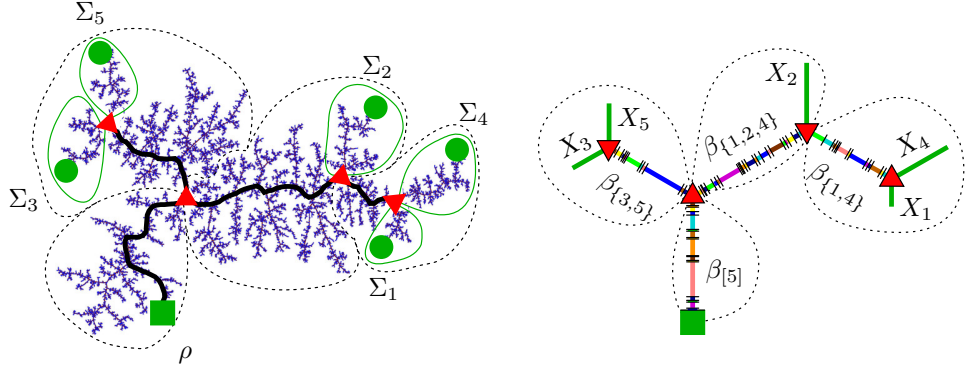
  \caption{Left: Simulation of a BCRT $(\cT,d,\rho,\mu)$, courtesy of Igor Kortchemski, 
  with root $\rho$, and $k=5$ leaves $\Sigma_1,\ldots,\Sigma_5$. The bold lines and triangles are the branches and vertices of $\cR_5$.  Right: The associated Brownian reduced $k$-tree.\label{fig:B_k_tree_proj}}
\end{figure}

The concept of a $k$-tree is best understood by generating a random $k$-tree from the BCRT, as in Figure \ref{fig:B_k_tree_proj}. Consider a BCRT $(\cT,d,\rho,\mu)$. 
Let $\Sigma_n$, $n\!\ge\!1$, denote a sequence of leaves sampled conditionally i.i.d.\ with law $\mu$. Denote by $\cR_k^+$ the subtree of $\cT$ spanned by $\rho,\Sigma_1,\ldots,\Sigma_k$, and by $\cR_k$ the subtree of $\cR_k^+$ spanned by the set ${\rm Br}(\cR_k^+)$ of branch points and the root $\rho$ of $\cR_k^+$. A.s., $\cR_k^+$ is a binary tree. 

Let $[k] := \{1,\ldots,k\}$ and suppose that $k\geq 1$. The \emph{Brownian reduced $k$-tree}, denoted by $R_k=\big(\ft_k,(X_j^{(k)},j\in [k]),(\beta_E^{(k)}, E \in{\rm IntEdge}(\ft_k))\big)$, is defined as follows.
\begin{itemize}
  \item 
  Let $\ft_k$ denote the rooted, binary tree with vertices ${\rm vert}(\ft_k)=\{\Sigma_1,\ldots,\Sigma_k\}\cup{\rm Br}(\cR_k^+)\cup \{\rho\}$ and edges connecting pairs of vertices if and only if the path between those points in $\cT$ does not pass through any other vertices of $\ft_k$.
  \item For $j\in[k]$, the \emph{top mass} $X_j^{(k)}$ is the $\mu$-mass of the component of $\cT\setminus\cR_k$ containing $\Sigma_j$.
  \item For each internal edge $E\in \text{IntEdge}(\ft_k)$, i.e.\ each edge between non-leaf vertices $v_1$ and $v_2$, let $\cB_E$ denote the unique non-self-intersecting path from $v_1$ to $v_2$ in $\cT$. 
  We assign an interval partition $\beta_E^{(k)}$ to this edge as follows. 
  Consider the connected components of $\cT\setminus\cR_k$ that attach to $\cR_k$ along the interior of $\cB_E$. These are ordered by decreasing distance between their attachment points on $\cB_E$ and the root. This order is not sequential since branch points of $\cT$ on $\cB_E$ are dense. We define $\beta_E^{(k)}$ to be the interval partition whose block sizes equal the $\mu$-masses of these components, in this order.
\end{itemize}
See Figure \ref{fig:B_k_tree_proj} for an illustration of such a $k$-tree.

In Section \ref{sec:killed_def} we formally define sets of $k$-trees, for each $k\ge1$, that support the laws of the Brownian reduced $k$-trees. The probability distribution of the Brownian reduced $k$-tree can be described in terms of some familiar objects.

\begin{proposition}[Section 3.3 of \cite{PW13}]\label{prop:B_ktree}
 Fix $k\ge1$. 
 \begin{itemize}[topsep=3pt]
  \item Let $\ft_k$ denote a uniform random rooted binary tree with $k$ labeled leaves.  
  \item Independently, let $(M_i,\ i\in [2k-1])\sim \distribfont{Dirichlet}\big(\frac12,\ldots,\frac12\big)$.
  \item Independently, let $\beta_i$, $k+1\le i\le 2k-1$, be independent $\PDIP\big(\frac12,\frac12\big)$.
 \end{itemize}
 Denote by $\phi_{\ft_k}\colon \textnormal{IntEdge}(\ft_k)\to\{k+1,\ldots,2k-1\}$ a bijection, e.g.\ via depth-first search of edges. Then $R_k=\big(\ft_k,(M_i,i\in [k]),(M_{\phi_{\ft_k}(E)}\beta_{\phi_{\ft_k}(E)},E\!\in\!\textnormal{IntEdge}(\ft_k))\big)$ is a Brownian reduced $k$-tree. In particular, its distribution is invariant under the permutation of labels.
\end{proposition}

\section{Main results: $k$-tree evolutions and the Aldous diffusion}
\label{sec:mainresults}

Let $(R_k,\,k\ge1)$ denote the Brownian reduced $k$-trees described above. There is a natural projection map $\pi_{k}$ from $(k+1)$-trees to $k$-trees for every $k\ge 1$ such that $\pi_k(R_{k+1})=R_k$. While $\pi_k$ affects the combinatorial tree shape $\ft_{k+1}$ by just removing vertex $\Sigma_{k+1}$ and the adjacent branch point, all top mass(es) and interval partition(s) adjacent to this branch point in $R_{k+1}$ are also suitably combined to form the resulting $k$-tree $R_k$. There is also a natural map $S_k$ taking $k$-trees to rooted, weighted $\BR$-trees such that, almost surely, $S_k(R_k)$ is the rooted, weighted $\BR$-tree that results from projecting the leaf mass measure of $\cT$ onto $\cR_k$.  To see this, observe that since the combinatorial tree shape of $\cR_k$ is $\ft_k$ and atoms of the projected measure are given by the top masses and intervals of the interval partitions, one only needs to be able to recover the metric structure of $\cR_k$ from $R_k$.  This is accomplished using the \emph{diversity} of the interval partitions $\beta$:
\begin{equation}\label{eq:diversity:intro}
 \sD(\beta) := \sqrt{\pi}\lim_{h\downarrow 0}\sqrt{h}\#\{(a,b)\in\beta\colon |b-a| > h\}.
\end{equation}
We will say more about this formula in Section \ref{sec:IP}. For now, we note that for the $\PDIP\big(\frac12,\frac12\big)$ formed by the excursion intervals of a Brownian bridge, as in Figure \ref{fig:BBr}, this limit exists almost surely. In fact, up to a scaling constant, the diversity of the interval partition equals the local time of the bridge at level 0. 
It was observed in \cite[Lemma 9 and equation (12)]{AldoMierPitm04} and \cite[Section 3.3]{PW13} that the diversities of interval partitions obtained by projecting the BCRT, as in Proposition \ref{prop:B_ktree}, equal distances in the tree almost surely. Hence, we use diversity to recover the metric structure along the branches $\cB_E$, $E\in\ft_k$, and define $S_k(R_k)$.

Given any system of $k$-trees $\mathrm{R}=(R_k,\, k\geq 1)$, we define $S(\mathrm{R}) = \lim_{k\to\infty} S_k(R_k)$, where the limit is taken with respect to the Gromov--Hausdorff--Prokhorov metric \cite{M09}, provided this limit exists, and $S(\mathrm{R})$ is the trivial one-point tree if the limit does not exist.  From the definitions of $\cR_k$ and $R_k$, $k\ge 1$, and from the relationship between diversities and distances, in the BCRT construction above, $S((R_k,k\geq 1)) = \cT$ almost surely. We formalize this 
construction in Section \ref{sec:Rtrees}.

The following three theorems summarize the main contributions of this memoir, resolving Conjecture \ref{conj:Aldous}.

\begin{theorem}\label{thm:intro:k_tree}
 There is a projective system $(\widebar\cT^s_{\!\bullet},\, s\geq 0)=\big((\widebar\cT^s_{\!k},k\geq 1), s\geq 0\big)$ such that the following hold.
 \begin{enumerate}[label=(\roman*), ref=(\roman*), itemsep=.1cm]
  \item \label{main markov} For each $k$, $(\widebar{\cT}_{\!k}^s,\,s\geq 0)$ is a $k$-tree-valued Markov process. 
  \item \label{main cons} The processes are consistent in the sense that $(\widebar\cT^s_{\!k},\,s\!\geq\! 0) = (\pi_k(\widebar\cT^s_{\!k+1}),\,s\!\geq\! 0)$. 
  \item \label{main stat} The law of the consistent family of Brownian reduced $k$-trees $(R_{k},\, k\geq 1)$ is a stationary law for the process $(\widebar\cT^s_{\!\!\bullet}, s\geq 0)$.
  \item \label{main Wright--Fisher} With notation $\widebar\cT^s_{\!k} = \big(\ft^s_k,(X^s_j,\,j\in[k]),(\beta^s_E,\,E\in\ft^s_k)\big)$, let $\|\beta^s_E\|$ be the mass of the interval partition $\beta^s_E$, and $\tau$ the first time a top mass vanishes.  Then $\ft^s_k=\ft^0_k$ for $s<\tau$ and $\big((( X^{(s\wedge)\tau/4-}_j,\,j\in[k]),(\|\beta^{(s\wedge\tau)/4-}_E\|,\,E\in\ft^0_k)), s\geq 0\big)$ is a Wright--Fisher diffusion ${\tt WF}(\theta_1,\ldots,\theta_{2k-1})$ with generator \eqref{eq:WFgen}, with the parameters $\theta_i=-\frac12$ for coordinates corresponding to top masses and $\theta_i=\frac12$ for coordinates corresponding to masses of interval partitions, stopped when a
top mass coordinate vanishes.
 \end{enumerate}
\end{theorem}


\begin{theorem}\label{thm:intro:AD}
 Let $(\widebar\cT^s_{\!\bullet},\, s\geq 0)$ denote the projective system of Theorem \ref{thm:intro:k_tree}, running in stationarity. Then the process $\big(S(\widebar\cT^s_{\!\bullet}),\, s\geq 0\big)$ of rooted, weighted $\BR$-tree-valued projective limits of this evolving system
 \begin{enumerate}[label=(\roman*),ref=(\roman*), itemsep=.1cm]
  \item  \label{main:CRT:stat} is stationary with the law of the BCRT,
  \item  \label{main:CRT:cont} has a modification that is path-continuous under the Gromov--Hausdorff--Prokhorov metric, and
  \item  \label{main:CRT:Markov} is a simple Markov process.
 \end{enumerate}
\end{theorem}

\begin{definition}\label{def:Aldousdiff}
 The \textit{Aldous diffusion} is a path-continuous modification of the process $\big(S(\widebar\cT_{\!\bullet}^s), s\geq 0\big)$ described in Theorem \ref{thm:intro:AD}.
\end{definition}

Together with \eqref{eq:3massconv}, this relates back to the Aldous chain via scaling limits of leaf proportions in finitely many components. In the next theorem, we establish the Aldous diffusion as a diffusive limit of a \em continuous-time Aldous chain \em with jump rate $1/2n(n-2)$ and jumps according to the non-trivial Aldous chain transitions (see Section \ref{sec:emb}).

\begin{theorem}\label{thm:scalinglim}
 The continuous-time Aldous chain, running in stationarity, represented as a process of $\mathbb{R}$-trees with edge lengths $1/\sqrt{n}$ and uniform weight measure on the leaves, converges to the Aldous diffusion as $n\to\infty$, in the sense of finite-dimensional distributions in the Gromov--Hausdorff--Prokhorov sense.
\end{theorem}

We prove this by embedding the continuous chain in the Aldous diffusion. 
Showing tightness to obtain a functional scaling limit theorem appears to be a difficult problem that relates to an only partially resolved problem about local times of stable processes. Such a result would immediately yield a scaling limit theorem for the discrete-time Aldous chain. 
We expect that the Aldous diffusion is also the scaling limit of many other down-up Markov chains with different stationary distributions in the large domain of attraction of the BCRT, but proving such results is beyond the scope of this memoir. 
We provide a fuller discussion in Section \ref{sec:openprob}. 


Our construction of the projective system of Theorem \ref{thm:intro:k_tree} yields a consistent family of stationary $k$-tree-valued Markov processes that turn out not to be reversible. It is not clear how the
reversibility of the Aldous diffusion can be proved directly from this construction. However, this is a direct consequence of the reversibility of the (continuous-time) Aldous chain and Theorem \ref{thm:scalinglim}.

\begin{corollary}\label{cor:rev:intro} The Aldous diffusion is reversible.
\end{corollary}

We will also explain in Section \ref{sec:discrete_to_cts} how our construction in the continuum relates to the discrete process. Before doing this, we discuss related literature, including another approach to Aldous's conjecture by L\"ohr, Mytnik and Winter \cite{LohrMytnWint20}, which does achieve process-level scaling limits of a version of the Aldous chain in a rather different state space of trees without the full continuum tree structure.

\section{Related literature}
\label{sec:related_work}

\subsection*{The Aldous diffusion project}

The present memoir is the culmination of ideas that we have developed in our previous joint work. Although it is not necessary to read all the previous papers to follow the mathematics here, in the interest of the ``big picture,'' Table \ref{tbl:project_outline} outlines their dependence structure. 
\begin{table}
 \centering

 \begin{tabular}{cccccc}
  Metrics on sets of interval&&&&& Uniform control of local times \\ 
  partitions with diversity \cite{Paper1-0}&&&&& of stable processes \cite{Paper0} 
\end{tabular}
 $\underbrace{\hphantom{Interval partition evolutions with emigration [20] of the Aldous chain [20] fill}}$\\
 $\downarrow$\\
 Construction of interval-partition-valued diffusions \cite{Paper1-1}\\
\begin{tabular}{ccccc}
  $\downarrow$&&&&\\
  Stationary Poisson--Dirichlet &&&& \\
  interval partition diffusions \cite{Paper1-2} &&&& \\
  $\downarrow$
  		&&&&	Consistent projections of \\
  Chapters \ref{ch:type-2}--\ref{dePoiss}: Two-tree evolutions
  		&&&&	the Aldous chain \cite{Paper2}
 \end{tabular}
 
 $\underbrace{\hphantom{Interval partition evolutions with emigration [20] of the Aldous chain [20] fill}}$\\
 $\downarrow$\\
 Chapters \ref{ch:constr}--\ref{ch:consistency}: Construction of the projective system of $k$-tree evolutions\\
 $\downarrow$\\
 Chapters \ref{ch:properties}--\ref{chap8}: Properties of the projective limit -- the Aldous diffusion
 
 \vphantom{text}
 
 \caption{Outline of the present authors' construction of the Aldous diffusion.\label{tbl:project_outline}} \vspace{-.7cm}
\end{table}%

The $k$-tree evolutions of Theorem \ref{thm:intro:k_tree} require Markovian evolutions on spaces of interval partitions with diversities. The state spaces were introduced in \cite{Paper1-0}, while the existence and properties of some of the evolutions have been worked out in \cite{Paper1-1,Paper1-2}. As mentioned above, the leaf masses of the continuum-tree-valued process are captured by the interval lengths of the evolving interval partitions. What is less explicit and not obvious is that the evolution of the metric structure of the continuum-tree-valued process is also captured by the diversity of the evolving interval partitions, as defined in \eqref{eq:diversity:intro}. This has been dealt with in \cite{Paper0} establishing and exploiting links between diversities and the local times of stable processes.

The consistency of the $k$-tree evolutions is subtle and requires a non-trivial labeling of the $k$-tree shapes and a \textit{resampling} mechanism when the masses of certain components reach 0. In \cite{Paper2}, a similar labeling and resampling scheme has been worked out for the Aldous chain where it has been proved to lead to consistent Markovian projections. Although the proofs of \cite{Paper2} cannot be generalized directly, this provides the basis for our approach to the consistency claimed in Theorem \ref{thm:intro:k_tree} and ultimately allows us to construct in the present work the continuum-tree-valued process of Theorem \ref{thm:intro:AD} and indeed of Conjecture \ref{conj:Aldous}, the Aldous diffusion.

Our published work on interval partitions is summarized in Chapter \ref{ch:prelim}; the remainder of this document presents previously unpublished work. Chapters 3-4 and Chapters 5-6 supersede two unpublished preprints, \cite{Paper3} and \cite{Paper4} respectively.
\subsection*{Related work by L\"ohr, Mytnik, and Winter}%
%
%

%
Recently L\"ohr, Mytnik, and Winter \cite{LohrMytnWint20} used a martingale problem to find the diffusive limit of the unrooted Aldous chain on a new space of trees, which they call binary algebraic measure trees.  They named their process the \textit{Aldous diffusion on binary algebraic non-atomic measure trees}, which they abbreviated to \textit{Aldous diffusion}, but for clarity, we will abbreviate as the \textit{algebraic Aldous diffusion}.  Algebraic measure trees \cite{LohrWint21} can be thought of as the structures that remain when one forgets the metric on a weighted $\BR$-tree but retains knowledge of the branch points; cf.\ mass-structural equivalence in \cite{IPTrees}. 
More formally, an algebraic measure tree is a triple $(t,c,\mu)$, where $t$ is a vertex set, $\mu$ is a measure on $t$, and $c\colon t^3\to t$ is a map that identifies the branch point separating any three vertices. 
Equivalence classes of algebraic measure trees form the state space for the algebraic Aldous diffusion.  The topology on the set of binary algebraic measure trees is most easily thought of as being the one generated by sample subtree convergence: for a binary algebraic measure tree $t$, let $t^{[k]}$ be the subtree spanned by $k$ vertices drawn independently from the measure on $t$.  Convergence of binary algebraic measure trees is essentially defined as $t_n\to t$ if and only if $t^{[k]}_n \stackrel{d}{\to} t^{[k]}$ for every $k$, on the space of graph-theoretic binary trees with $k$ leaves.  The actual definition of the topology is more subtle, but reduces to this for the trees that appear when studying the Aldous chain on algebraic trees; see \cite[Proposition 2.8]{LohrMytnWint20}.

There are advantages and disadvantages to using this state space and topology compared to our choice of the Gromov--Hausdorff--Prokhorov setting with distances given by the rescaled graph metric. The most significant advantage is that constructing the algebraic Aldous diffusion and proving convergence to it can be done using classical martingale problem techniques.  While the calculations are involved, they are surprisingly simple relative to what is required in the Gromov--Hausdorff--Prokhorov setting of the present work.  Additionally, they are able to show that the algebraic Aldous diffusion is an ergodic Feller process whose unique stationary distribution is the algebraic Brownian CRT.  One disadvantage is that, because of the topology, the only statistics that can be computed are averages of quantities over uniformly sampled subtrees. Because of this, the Wright--Fisher diffusions that Aldous described appear only in an annealed sense. 
Another disadvantage is that the aforementioned statistics do not capture lengths in the tree. Indeed, the algebraic tree setting was chosen by \cite{LohrMytnWint20} 
in order to sidestep the difficulties posed by distances in the tree: 
 the quadratic variation of the averaged distance process scales like $N^{3/2}$ instead of $N^2$, suggesting that distances in the Aldous chain might fluctuate too wildly for tightness to hold; see \cite[p. 2567]{LohrMytnWint20}.

An advantage of our approach is that we can show that the stationary tree-valued process converges, in the sense of finite-dimensional distributions under the Gromov--Hausdorff--Prokhorov metric, to a continuous limiting process, in which the evolution of distances is described in terms of local times of stable L\'evy processes. 
See Section \ref{sec:Holder} for the connection to local times.  The ability to understand the evolution of distances is a primary advantage of our approach but also a great source of complexity, leading us to decorate our trees with interval partitions. 
Another advantage of our setting is that the Wright--Fisher processes described by Aldous appear directly as mass evolutions around selected branch points.  A disadvantage of our approach is that we only construct the diffusion in stationarity. 
 Additionally, we note that in contrast to \cite{LohrMytnWint20} where the limiting process is Feller, in our setting the limiting process has the simple, but not the strong, Markov property; see Section \ref{sec:smp}. This is a disadvantage of our choice of state space in the sense of doing calculations with the limiting process. However, this is the state space most naturally implied by Aldous's conjecture, and our method reveals the unexpected failure of the strong Markov property. See Section \ref{sec:openprob} for further discussion and open problems related to these matters.

There is a natural conjecture relating the processes we call the Aldous diffusion and the algebraic Aldous diffusion.  In particular, in Theorem \ref{thm:intro:k_tree} we construct a consistent system of $k$-tree evolutions $(\widebar{\cT}^s_{\!k},\,s\geq 0)$, $k\geq 2$ that capture the Wright--Fisher diffusions proposed by Aldous.  In Definition \ref{def:Aldousdiff} we define the Aldous diffusion by mapping $\widebar{\cT}^s_{\!k}$ to a weighted $\BR$-tree using the diversities and block sizes in the interval partitions to determine branch lengths and masses of atoms, then taking a projective limit as $k\to \infty$ in the Gromov--Hausdorff--Prokhorov metric.  If instead of mapping $\widebar{\cT}^s_{\!k}$ to a weighted $\BR$-tree we map it to an algebraic measure tree and take the limit in the topology of \cite{LohrMytnWint20} -- this limit is easily seen to exist -- we conjecture that the resulting process is the diffusive limit of the rooted Aldous chain in a space of rooted algebraic trees.%
\subsection*{Continuum-tree-valued Markov processes}
The Aldous diffusion is not the first continuum-tree-valued Markov process. Notable previous examples include \cite{EvansStFlour,EPW,EW,Zambotti01}. Of these, only \cite{EW} has as its state space the Gromov--Hausdorff--Prokhorov (GHP) space of (unrooted) weighted $\mathbb{R}$-trees. See \cite[Section 6.3]{M09} for a proof that the metric used in \cite{EW} is indeed equivalent to the GHP metric. \cite{EPW} uses a Gromov--Hausdorff state space of unweighted $\mathbb{R}$-trees, while \cite{Zambotti01} uses a
space of unit length excursions, which relates to weighted $\mathbb{R}$-trees as explained around \eqref{eq:tree_from_path}. All have the BCRT (or the associated random 
$\mathbb{R}$-tree without weight measure, or a Brownian excursion) as their stationary distribution. 

Specifically, Evans, Pitman, and Winter \cite{EPW} study \em root growth with re-grafting\em. Under these dynamics, an $\mathbb{R}$-tree grows just above the root pushing up 
the tree at unit speed, while a Poisson point process at unit rate per unit length places cut-points onto the growing tree. Cut subtrees are re-grafted to the root. 
Indeed, starting from a one-point tree, the tree just before the $j^{\rm th}$ re-grafting has the same distribution as the $j^{\rm th}$ tree in the line-breaking construction. The dynamics are consistent across all connected subsets containing the root and can be defined on the Gromov--Hausdorff space of rooted compact $\mathbb{R}$-trees. This gives rise
to a recurrent Feller process that converges to stationarity, the law of the BCRT. This discontinuous evolution of the BCRT is the scaling limit of a continuous-time Markov chain on discrete trees with $n$ vertices called the \em Aldous--Broder chain\em, see \cite[Section 5]{Broder89} and Aldous \cite[end of Section 2]{Aldous90}. This chain selects each non-root vertex at rate 1, cuts the adjacent edge towards the root, adds an edge from the selected vertex to the root and re-roots at the selected vertex.

Evans and Winter \cite{EW} are motivated by a Markov chain on binary trees with $n$ leaves whose transitions are \em subtree prune and re-graft \em moves, in which a cut-edge and a re-graft edge are chosen uniformly at random. Specifically, in the component containing the re-graft edge, the branch point of the cut-edge is contracted away, while the other component is re-grafted by connecting the cut-edge to a new branch point in the middle of the re-graft edge. They use Dirichlet-form techniques to build a recurrent Hunt process on the Gromov--Hausdorff--Prokhorov space of weighted $\mathbb{R}$-trees that they expect to be the scaling limit of this Markov chain. Informally, a Poisson point process places cut-points at unit rate per unit length, while the re-graft point is then selected from the weight measure. Again, the component not containing the re-graft point is re-grafted. This discontinuous evolution is reversible with respect to the distribution of the BCRT. 

Zambotti \cite{Zambotti01} studies a class of stochastic partial differential equations driven by a Brownian sheet $W$ on $[0,\infty)\times[0,1]$, which includes the special case
\begin{equation}\label{eq:Zambotti}
\left\{\begin{array}{l}\displaystyle\frac{\partial u}{\partial s}=\frac12\frac{\partial^2u}{\partial^2r}+\frac{\partial^2W}{\partial s\,\partial r}+\frac{\eta(ds,dr)}{ds\,dr}\\[0.3cm]
                                  u(0,r)=x(r),\quad u(s,0)=u(s,1)=0,\quad
                                  u\ge 0,\quad u=0\quad\eta\mbox{-a.e.},
        \end{array}
\right.
\end{equation}
whose solutions $(u,\eta)$ we interpret as evolutions in time $s\ge 0$ of continuous functions $u(s,\,\cdot\,)\colon[0,1]\rightarrow[0,\infty)$ that vanish at 0 and 1, with initial
function $u(0,\cdot)=x$, while the auxiliary measure $\eta$ is to achieve reflection at the boundary. Zambotti shows that there is a solution $(u,\eta)$ such that $(u(s,\,\cdot\,),\,s\ge 0)$ is a path-continuous Markov process (adapted to the filtration of the Brownian sheet, generated by $W|_{[0,s]\times[0,1]}$, $s\ge 0$) for which the law of the standard Brownian excursion is an invariant measure. 
Zambotti \cite{Zambotti02} establishes a Dirichlet form for this process and shows that $\eta(\,\cdot\,,(0,1))$-a.e.\ $u(s,\cdot)$ has precisely one zero in $(0,1)$. The study of the reflection at the boundary is further refined in \cite{Zambotti04} establishing occupation densities, while \cite{Zambotti08} introduces related dynamics of Brownian excursions conditioned on their area. See also \cite{Zambotti17} for lecture notes on this material. 

Finally, \cite[Theorem 3]{EtheLabb15} showed that Zambotti's Markov process \cite{Zambotti01} is the scaling limit of a Markov chain on Dyck paths of length $2n$ that at each step chooses a flipping point on the path uniformly at random. If this point is a local minimum it is flipped into a local maximum by adding 2 at this point, and if it is a local maximum above height 2, it is flipped into a local minimum by subtracting 2. Zambotti \cite[Section 5.6.4]{Zambotti17} poses the problem of providing a description of the limiting process directly as an evolution of trees. There is no heuristic reason to believe that the projection of Zambotti's process into the space of rooted weighted $\mathbb{R}$-trees would be the same as the Aldous diffusion. However, due to the difficulty of computing statistics related to these processes, it is not straightforward to distinguish them. 
In Section \ref{sec:nonbin} we argue, but stop short of proving, that these processes are indeed distinct. In Section \ref{sec:openprob} we explore the possibility of a variant of our Aldous diffusion on a space of excursions.%

\subsection*{Dynamics for Poisson--Dirichlet distributions}

A one-parameter family of Poisson--Dirichlet distributions on the Kingman simplex was introduced by Kingman, in 1978, in his work on population genetics \cite{Kingman1975}. This was extended to two parameters in 1997 by Pitman and Yor \cite{PitmYor97}.

In 1981, Ethier and Kurtz introduced a measure-valued Fleming--Viot diffusion called the infinitely-many-neutral-alleles model \cite{EthiKurt81} to describe fluctuations in allele frequency in a large population under no fitness preferences (hence ``neutral''), with the possibility of mutation into infinitely many new types. Under the projection that maps a purely atomic measure to its ranked sequences of atom masses, this process maps to a diffusion on the Kingman simplex with $\PoiDir(0,\theta)$ stationary distribution. In 2009, Petrov generalized this latter diffusion to the two-parameter $\PoiDir(\alpha,\theta)$-setting \cite{Petrov09}. Around that same time, it was conjectured \cite{FengBook,FengSun10} that the aforementioned Fleming--Viot diffusion should be similarly generalizable: that there should exist a two-parameter family of Fleming--Viot diffusions that would project down to Petrov's diffusions on the Kingman simplex. Both Petrov's work and Feng and Sun's conjecture spurred a great deal of research; see e.g.\ \cite{CdBERS17,FengSunWangXu11,RuggWalk09,Zhou15}.

En route to the present work, the research program outlined in Table \ref{tbl:project_outline} also provided the tools to resolve Feng and Sun's conjecture. In particular, the stationary interval partition diffusions mentioned in the table have Poisson--Dirichlet stationary laws. We generalized and adapted these into measure-valued processes in \cite{PaperFV,ShiWinkel-2}. In \cite{Paper1-3}, we resolved the conjecture by showing that the constructed Fleming--Viot diffusions 
project down to Petrov's diffusion on the Kingman simplex.

Poisson--Dirichlet and similar mass splits ordered along a geodesic or represented as an interval partition appear in wider families of random metric spaces and related structures. See 
for example \cite{GoldHaas15,HPW,PitmWink09,RW2,RW} for classes of random $\mathbb{R}$-trees and \cite{ABACFG,ABBG,ConcGold23,GoldHaasSeni22} for critical random graphs with some cycles as they arise in the critical window of the Erd\H{o}s--R\'enyi random graph, stable graphs and random planar maps related to Voronoi tesselations and large uniform planar maps with a fixed finite number of faces on higher-genus surfaces.

\section{Consistent $k$-tree down-up chains and skewer representations}
\label{sec:discrete_to_cts}

In this section, we discuss several features of our approach to the Aldous diffusion in the discrete setting of the Aldous chain. While reading this section is not strictly necessary to
understand the mathematics that we develop in the remaining chapters, the notions introduced here provide useful motivation and insights into some of the obstacles that we encounter and some of the concepts that we develop as building blocks in our construction of the Aldous diffusion. In the discrete setting here, some technical complexity is absent, and this is helpful to more directly expose some of the challenges before we turn to the setting of interval partition evolutions and $k$-tree evolutions. In Section \ref{sec:intro:overview}, we discuss the structure of this memoir using some terminology for the building blocks that we introduce here, but that will be formally set up in the further chapters.

\subsection*{Label swapping and consistent Markovian projections}

Aldous \cite{Aldous00} and Schweinsberg \cite{Schweinsberg02} showed that the relaxation time for the unrooted Aldous chain is $\Theta(n^2)$, where $n$ is the number of leaves, and indeed, $n^2$ is the correct number of steps per unit time for mass fluctuations around branch points to have scaling limits in \eqref{eq:3massconv}. However, it only takes $n\log(n)$ steps for every single leaf to be displaced. 
Hence, leaf labels mix much faster than tree structure in the Aldous chain. 
The key challenge here is that we often describe continuum trees based on subtrees spanned by a finite sequence of leaves. 

In \cite{Paper2} we take the liberty of modifying the Aldous chain in order to slow down the motion of low leaf labels in such a way that it does not affect the dynamics of the underlying unlabeled tree. This allows us to identify persistent \emph{branch points} -- for example, the branch point separating leaf label 1 from label 2. In the continuum, these are the branch points that make up the $k$-tree structures of Theorem \ref{thm:intro:k_tree}.

For $n\ge 1$, we denote by $\bT^{\rm graph}_{n}$ the set of rooted binary trees with $n$ labeled leaves, with internal vertices treated as being unlabeled; the trees in the first and last panels of Figure \ref{fig:AC_move} are members of $\bT^{\rm graph}_6$. For $k\ge 1$, a \emph{decorated $k$-tree} is a binary tree $\overline{\ft}\in\bT^{\rm graph}_k$ with external edges (those incident to leaves) decorated with positive integer weights and internal edges (i.e.\ non-external ones) decorated with (potentially trivial) sequences of positive integer weights.

\begin{figure}
 \centering
 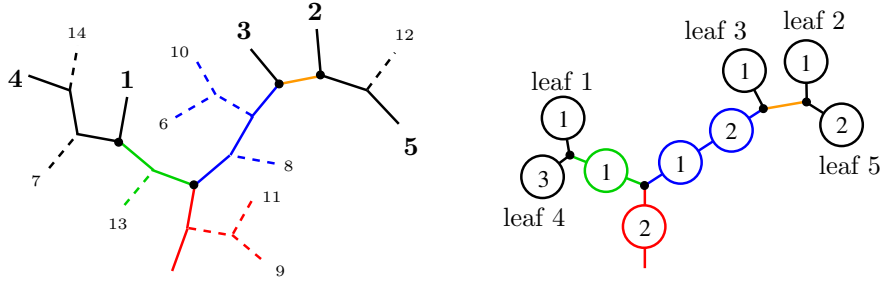
 \caption{The projection onto decorated 5-trees, applied to a 14-leaf tree $\overline{\mathbf{t}}$, with internal edges color-coded to indicate corresponding subtrees of $\overline{\ft}$. The (sequence of) weight(s) on each external (resp.\ internal) edge of the projected tree is inscribed in a (sequence of) circle(s) at the end (resp.\ in the middle) of the edge.}
\label{fig:decorated_k_tree}
\end{figure}

Figure \ref{fig:decorated_k_tree} illustrates the \emph{decorated $k$-tree projection} of a tree $\overline{\ft} \in \bT^{\rm graph}_{14}$, in the case $k=5$. This projection, which we denote by $\pi^\circ_k(\overline{\ft})$, is obtained by the following steps.
\begin{enumerate}
\item Consider the subtree spanned by leaves $1,\ldots,k$ and the root. 
\item Contract away all degree-2 branch points to obtain a binary tree $\overline{\ft}_k\in\bT^{\rm graph}_k$.
\item Decorate each external edge of $\overline{\ft}_k$ with weight equal to the number of leaves in the subtree of $\overline{\ft}$ corresponding to that edge.
\item Decorate each internal edge of $\overline{\ft}_k$ with a sequence of weights equal to the leaf counts in each of the subtrees grafted to the corresponding path in $\overline{\ft}$, in order of decreasing distance from the root.
\end{enumerate}

In the example in Figure \ref{fig:decorated_k_tree}, in the right panel, the internal edge colored blue is decorated with the sequence $(2,1)$, in order of decreasing distance from the root, corresponding to the two subtrees of the blue path in the left panel: the one containing leaves 6 and 10, and the one containing only leaf 8. The orange edge, to the right of the blue edge, is decorated by the null sequence, as no subtrees are grafted to the orange path in the panel on the left. The green and red edges are each decorated with singleton sequences corresponding to single subtrees.

\begin{remark}
 In \cite{Paper2}, the term ``decorated $k$-tree'' is first defined differently, with internal edges decorated by a single, non-negative mass, equal to the sum of the sequence of masses considered here. However, the ``decorated $k$-trees'' of this section are taken up in Section 5 of that paper.
\end{remark}

These projections are the discrete analogs to the $k$-trees of Theorem \ref{thm:intro:k_tree}. The benefit of modifying the leaf label dynamics of the Aldous chain is two-fold. We achieve a slowing-down of the movement of small labels, and the projections of the modified chain to decorated trees, for each $k\in[n]$, are Markov chains themselves.  Indeed, it is not hard to see that projections while using the original label dynamics of the Aldous chain are not Markovian, in general: e.g.\ for $k=2$, if leaf 1 is moved twice in a row, first to the root edge and then elsewhere, the initial state constrains the state after two steps while the intermediate state consists of a null sequence and weights $1$ and $n-1$ as decorations of the two external edges.   

 %
\begin{figure}\centering
 \scalebox{.75}{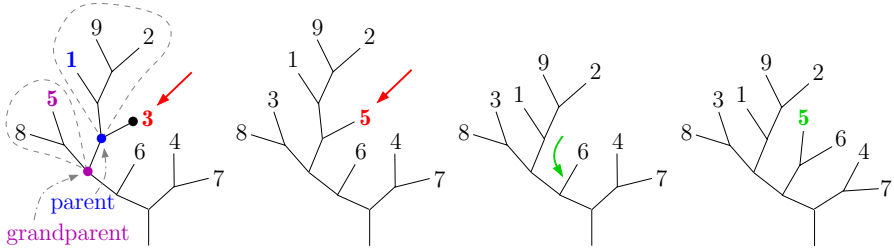}
 \caption{A down-up move with label swapping as in Definition \ref{def:modified_AC}. Here, $(i,a,b) = (3,1,5)$, so $\wi\!=\!5$. We swap labels $i\!=\!3$ and $\wi\!=\!5$ before deleting the chosen leaf. Label $5$ then regrows.\label{fig:downmove}\vspace{-0.1cm}}
\end{figure}

\begin{definition}[Modified Aldous chain, Definition 1 of \cite{Paper2}]\label{def:modified_AC}
 Fix $n\ge 3$. The \em modified Aldous chain \em is the down-up Markov chain on $\bT^{\rm graph}_{n}$, in which each transition has the following two steps. 
  \begin{enumerate}[label=(\roman*), ref=(\roman*)]    
  \item Down-move: Sample a uniform leaf $I$. Suppose $I=i$. Let $a$ (resp., $b$) denote the smallest leaf label from the other subtree of the parent branch point (resp., grandparent branch point, unless this is the root, in which case, by convention, $b=0$). Let $\wi=\max\{i,a,b\}$. Swap labels $i$ and $\wi$, unless $i=\wi$. Next, remove the leaf now labeled $\wi$ (which had been labeled $i$) and contract away its parent branch point. 
        
      \item Up-move: insert a new leaf $\wi$ at an edge chosen uniformly at random.
  \end{enumerate}
\end{definition}

Figure \ref{fig:downmove} illustrates a transition of the modified chain in which the leaf labeled 3 is selected for deletion, so label 3 swaps places with label 5, then the newly labeled leaf 5 is deleted and regrown elsewhere in the subsequent up-move. Crucially, the dynamics for the unlabeled tree are unchanged from those indicated around Figure \ref{fig:AC_move}; only leaf labels are affected.

These label dynamics may seem arbitrary; indeed, for brevity's sake, we have told this story backwards. In fact, this scheme emerged from first devising sensible Markovian dynamics for projected 2-trees, then 3-trees, etc., finally looking to the top of the resulting projective tower to find dynamics for individual leaf labels.

In \cite[Theorem 1]{Paper2} we observed that this modified chain is, like the original Aldous chain, stationary under the uniform distribution on $\bT^{\rm graph}_n$. Interestingly, despite the asymmetrical handling of labels in the transitions, leaf labels are exchangeable at fixed times, in stationarity. In \cite[Theorem 19]{Paper2}, we observed that the $\pi^\circ_{k}$-projections of a stationary modified Aldous chain are themselves stationary Markov processes, for each $k$.

A key feature of this label-swapping scheme is that for each $k\in[n]$, the decorated $k$-tree projection $\pi^\circ_{k}(T(j))$, $j\ge0$, of a modified Aldous chain $(T(j),\,j\ge 0)$ retains its shape in $\mathbb{T}^{\rm graph}_k$ much longer than the projection of the original chain. In particular, this shape only changes after an external weight has been reduced to 1 and the internal decoration immediately below to the null sequence, and then the single leaf of the external weight is selected for deletion in a down-move.

Note that for $k_1 > k_2$, the $\pi^\circ_{k_2}$-projection of a tree $\overline{\ft}\in\bT^{\rm graph}_{[n]}$ can be obtained in a natural manner as a further projection of the $\pi^\circ_{k_1}$-projection. Hence, the sequence of these Markov chains, for varying $k$, form a consistent projective system. Our aim is to mimic this system in the continuum, but rather than going from the ``top down,'' beginning with a continuum-tree-valued diffusion and then describing its projections, we must go in reverse: we will define $k$-tree-valued processes and use them to obtain a continuum-tree-valued process as their projective limit.

The label swapping of Definition \ref{def:modified_AC} solves the problem of identifying persistent branch points; in the remainder of this section, we describe a representation of the dynamics for our projected trees in a way that can be passed to the continuum.

\subsection*{Dynamics for decorated 2-trees and ordered Chinese restaurants}

Let us focus on the simplest case for our projections: $k=2$. Consider the branch point that separates leaves 1 and 2 from each other and the root. We decompose the tree into 
two \emph{top subtrees} above this branch point and a 
sequence of \emph{spinal subtrees} grafted to the path, called the 
\emph{spine}, from the branch point to the root. The decorated 2-tree comprises a pair of \emph{top masses} 
$(m_1,m_2)$, followed by a (possibly null) finite sequence of \emph{spinal masses} 
$(b_1,\ldots,b_{k})$, ordered by decreasing distance from the root.

\begin{figure}\centering
 \input{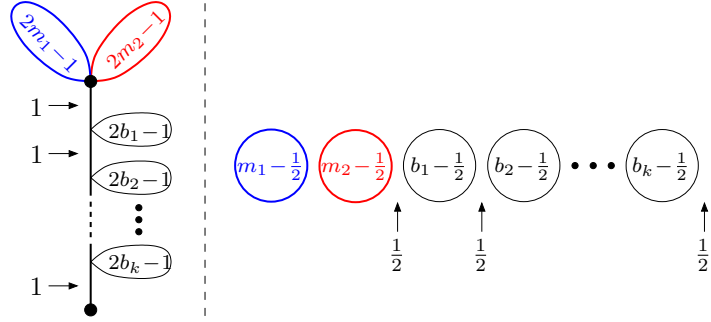}
 \caption{Left: Up-move weights for a 2-tree projection of the Aldous chain. The tree is decomposed into two top subtrees with $m_1$ and $m_2$ leaves, and $k$ spinal subtrees, with $b_1,\ldots,b_k$ leaves, respectively. Right: Seating weights for $\oCRP\big(\frac12,-\frac12\big)$, when the tables have $m_1,m_2,b_1,b_2,\ldots,b_k$ customers, in left-to-right order.\label{fig:2tree_CRP}}
\end{figure}

The down-up moves of Definition \ref{def:modified_AC} act on the decorated 2-tree as follows. In the 
down-move, we make a size-biased pick among the masses and reduce that mass
by one. 
For up-moves, we choose a mass $m$ with probability proportional to $2m - 1$,
or choose any edge along the spine with probability proportional to 1; see
Figure \ref{fig:2tree_CRP}. If a mass is chosen, it is incremented by 1; if
a spinal edge is chosen, a `1' is inserted into the sequence of spinal
masses at that point, representing the appearance of a new spinal subtree.

There are three cases for the down-move that require more explanation. 
\begin{enumerate}[label=(D\arabic*),ref=(D\arabic*)]
 \item If one of the spinal masses is reduced to 0 in the down-move, then it is deleted from the list; this corresponds to a single-leaf subtree and its parent branch point being contracted away, with the spine consequently shrinking in length.\smallskip%
 \label{item:2tdm:spine_shorter}%
 \item If one of the two top masses is reduced to zero at a time when there is at least one spinal mass, then the first spinal mass replaces it as a new top mass. This results from the label swapping of Definition \ref{def:modified_AC}, which will swap the label in the targeted top subtree down into the ``uncle'' subtree, one step down the spine. In this case, we say that the label is \emph{moving down the spine}.\smallskip%
 \label{item:2tdm:moving_down}%
 \item If one of the top masses is reduced to zero at a time when the spine is trivial, meaning that all but one unit of mass was already in the other top mass, then the chain jumps to a new state sampled from the stationary law for this chain, which is the $\pi^\circ_2$-image of the uniform law on $\bT^{\rm discr}_n$. We call this \emph{resampling}.
 \label{item:2tdm:resampling}
 
 This arises from the degenerate $b=0$ case of Definition \ref{def:modified_AC}. If the top mass being reduced to 0 bears label 2, then this will not swap at all, and will be deleted and regrow at a random new edge. If it bears label 1, then this will swap with label 2 before deletion, but the effect is the same.
\end{enumerate}

The up-move weights in the left panel of Figure \ref{fig:2tree_CRP} resemble the seating rule for an \emph{ordered Chinese restaurant process} (\oCRP) \cite{PitmWink09,RogeWink20}. The \oCRPAT\ begins with a single customer sitting alone at a table. New customers enter one-by-one. Upon entering, the $n+1^{\text{st}}$ customer chooses to join a table that already has $m$ customers with probability $(m\!-\!\alpha)/(n\!+\!\theta)$; sits alone at a new table inserted at the far left end of the restaurant with probability $\theta/(n\!+\!\theta)$; or sits alone at a new table, inserted to the right of any particular table already present, with weight $\alpha/(n\!+\!\theta)$, so that the total probability to sit alone is $(k\alpha\!+\!\theta)/(n\!+\!\theta)$, where $k$ is the number of tables already present. 
If we ignore the left-to-right order of these tables, then this is the well-known (unordered) \CRPAT\ due to Pitman \cite[\S 3.2]{CSP}, which generalizes the $\alpha=0$ case first studied by Blackwell and MacQueen \cite{BlacMacQ73}. 

If we take $(\alpha,\theta) = \big(\frac12,0\big)$, then this seating rule differs from the up-move probabilities in Figure \ref{fig:2tree_CRP} only in that, in the \oCRP, a new table can be introduced between the two leftmost tables, whereas in the 2-tree no new mass can be inserted in between the two leftmost masses, representing the two top subtrees, which are not separated by an edge but only by a branch point. We refer to the probabilities in Figure \ref{fig:2tree_CRP} as the seating rule for the \oCRP$\big(\frac12,-\frac12\big)$, as, under this rule, if there are a total of $k$ masses (2 top masses and $k-2$ spinal masses), then the probability for insertion of a new mass `1' is $(k-1)/(2n-1) = \big(k\frac12 - \frac12\big)/\big(n-\frac12\big)$.

This is outside of the usual parameter range considered for the \CRP. Indeed, if we start a \CRP$\big(\frac12,-\frac12\big)$ with a single customer, as described above, then all subsequent customers will be forced to join the first at a single table, as the probability to sit alone will be zero. However, if we start with two customers sitting separately, then the \oCRP$\big(\frac12,-\frac12\big)$ seating rule produces a non-trivial configuration with the same distribution as the decorated 2-tree projection of a uniform random rooted binary tree with labeled leaves. In this analogy, the down-up moves of the modified Aldous chain become \emph{re-seating}: a uniform random customer leaves their seat; their table is removed if empty; and they choose a new seat according to the seating rule, as if entering for the first time.


\subsection*{Discrete scaffolding, spindles, and skewer}

We simplify matters by \emph{Poissonizing} the re-seating \oCRP\ described above. In the Poissonized process, each customer exits the restaurant after an independent exponential time with rate 1, and each of the seating weights in the right panel of Figure \ref{fig:2tree_CRP} is taken as an exponential rate at which customers will either join a given table or sit alone at a new table in a given position. Thus, the total number of customers in the restaurant will fluctuate according to a birth-and-death Markov chain. In particular, we no longer think of this as \emph{re}-seating -- rather, old customers exit the restaurant independently of new customers entering. 
This technique was previously used 
in \cite{Pal13} to rigorously establish the connection proposed by Aldous \cite{AldousDiffusionProblem} between the Aldous chain and Wright--Fisher diffusions, which inspired Conjecture \ref{conj:Aldous}. As in that paper, after we pass to the continuum, we will apply a \emph{de-Poissonization} transformation that normalizes the total mass of the process and applies a corresponding time-change in order to preserve the Markov property.

In the Poissonized process, the table populations evolve independently of each other. Each one is a birth-and-death chain, having deaths with rate $m$ and births with rate $m-\frac12$ when the population is $m$, until absorption at population 0. Meanwhile, to the right of any table except for the leftmost (i.e.\ not between the two leftmost), a new table of population 1 appears with rate $\frac12$.

The Poissonized down-up \oCRPAT\ admits a Ray--Knight representation of a form that was introduced in \cite{Paper1-1,IPPAT} and used in \cite{Paper1-2,IPPAT} to construct interval partition diffusions that arise as scaling limits \cite{ShiWinkel-2} of the Poissonized down-up \oCRP$(\alpha,\theta)$ in the regime $\theta\ge 0$. 
An expanded discussion of this representation in the discrete regime can be read in \cite{RogeWink20}. See also \cite{RivRiz20} for scaling limit results on down-up \oCRPAT\ without Poissonization.

\begin{figure}
 \centering
 \scalebox{.75}{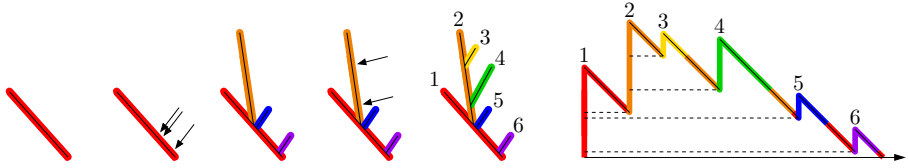}
 \caption{Iterative construction of splitting tree representation of tables in Poissonized down-up \oCRP[\frac12,0], started with one customer, and JCCP.\label{fig:splitting}}
\end{figure}

In Chapter \ref{ch:type-2} we will extend this representation from the $\theta\ge0$ regime to the required case $(\alpha,\theta) = \big(\frac12,-\frac12\big)$. For now, we will discuss the $\big(\frac12,0\big)$ case. 

We think of the tables that appear and vanish in the evolving \oCRP\ as members of a family: when a new table is born, the table immediately to its left at that time is its parent. The number of tables is then evolving over time as a homogeneous Crump--Mode--Jagers (CMJ) branching process \cite{Jagers69}. The genealogy among these tables, and their lifetimes, can be represented in a splitting tree \cite{GeigKers97}. For our purposes, this can be formalized as a rooted plane tree with edge lengths.

Figure \ref{fig:splitting} depicts the construction of a splitting tree representation of the Poissonized down-up $\oCRP\big(\frac12,0\big)$ started with a single customer.
\begin{enumerate}[label=(\arabic*), ref=(\arabic*)]
 \item Draw a line of random length, sampled from the probability distribution $\mu$ of the lifetime of a table started with population 1; this represents the first table. One end of this line will be the root of the tree, representing time 0.
 \item Now, mark that line with Poisson points along its length with rate $\frac12$, representing birth events.
 \item At each marked point, attach a new ``child'' line, branching off to the right from its parent, with length independently sampled from $\mu$. Each such line represents a table ``born,'' at some time, immediately to the right of the first table.
 \item Repeat steps (2), (3), and (4) on each of the newly drawn lines, if any.
\end{enumerate}
It can be shown that this procedure almost surely terminates for this choice of $\mu$.

This tree can be represented by a \emph{jumping chronological contour process} (JCCP) \cite{Geiger95,GeigKers97}, shown in Figure \ref{fig:splitting}. Imagine a flea traveling around the splitting tree. It begins to the left of the root and immediately jumps up to the top of the leftmost branch, representing the first table. It then slides down the right hand side of that branch at unit speed until its path is blocked by a branch sticking out to the right. When that happens, it jumps to the top of the new branch, and carries on in the same manner, until it finally reaches the root. The JCCP records the distance from the flea to the root, as a function of time.

The tables that arise in the evolving \oCRP\ are in bijective correspondence with the jumps of the JCCP, with the levels of the bottom and top of each jump equaling the birth and death times of the corresponding table. The genealogy among tables can be recovered from the JCCP by looking to the bottom of each jump (the birth-time of a child), and drawing a horizontal line to the left from that point, seeing where it crosses another jump (its parent). 

JCCP representations of splitting trees like ours are L\'evy processes of positive jumps and negative drift \cite{Lambert10}. Our particular JCCP has drift $-1$ and L\'evy measure $\frac12\mu$. \emph{Levels} in the JCCP correspond to \emph{times} in the evolving \oCRP. On the other hand, times in the JCCP have no simple meaning in the \oCRP, and serve mainly to record the left-to-right order of tables.

What is missing from this JCCP picture is the evolving table populations. Recall that each table population evolves as a birth-and-death chain with lifetime distribution $\mu$. This is also the law of jump heights in our JCCP. We incorporate both the genealogy among tables and the evolving table populations into a single formal object by marking each jump with such a birth-and-death chain, with lifetime equal to the height of the jump.

\begin{figure}[t]\centering
  \includegraphics[scale=.7]{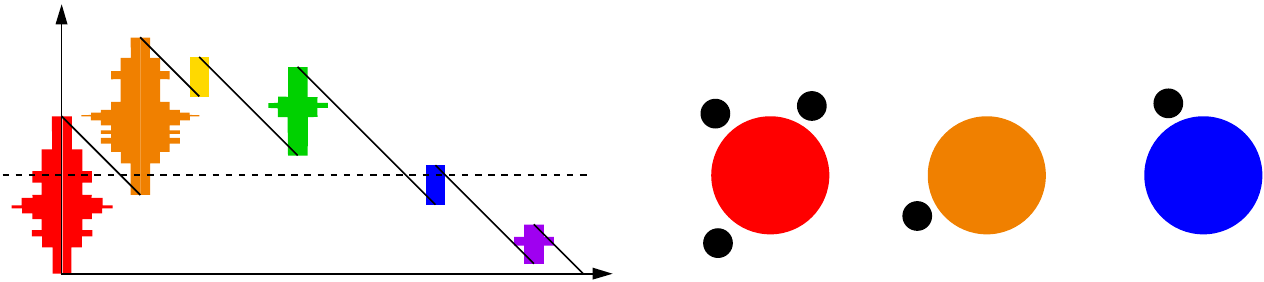}
  \caption{Left: Scaffolding with spindles. Right: a composition, visualized as an ordered Chinese restaurant state, corresponding to the skewer at a fixed level.}
  \label{fig:discrete_scaff_skewer}
\end{figure}

We depict this by representing each birth-and-death chain as a laterally symmetric ``spindle'' shape, beginning at the bottom of the jump and evolving towards its top, with width at each level describing the value of the chain at the corresponding time. In the context of this construction, we refer to the JCCP as \emph{scaffolding} and the markings as \emph{spindles}. See Figure \ref{fig:discrete_scaff_skewer}. In the CMJ framework \cite{Jagers69,JagersBranching}, our spindles are individual ``characteristics'' changing during lifetimes.

Then, to recover the Poissonized down-up \oCRP$\big(\frac12,0\big)$ from the scaffolding and spindles representation, we apply a \emph{skewer map}: for any $y\ge0$, we draw a horizontal line through the picture at that level, and look at the cross-sections of spindles pierced by the line. The widths of these cross-sections represent populations of tables, and their left-to-right order corresponds to that in the \oCRP. If we slide this horizontal line up continuously, then the cross-sections gradually change in width, with some dying out as the horizontal line passes the top of a jump, and new ones appearing as it reaches the bottom of a jump.

\subsection*{Putting the pieces together}

Recall the discussion of ``resampling'' as a special case \ref{item:2tdm:resampling} that may arise in a down-move acting on the decorated 2-tree. Analogous behavior arises in all decorated $k$-tree projections, when an external edge is reduced to mass 0 at a time when the internal edge below it also has mass 0. For fixed $k$, the down-up chain on decorated $k$-trees proceeds through $O(n^2)$ steps in between resampling events \cite{Pal13}, cf.\ \eqref{eq:3massconv}. In our continuum analog, we will construct $k$-tree-valued processes with certain continuous-time dynamics that are interrupted at discrete resampling times, when the process jumps away from a degenerate state in which an external component and the internal component below it simultaneously hit mass 0. The label-swapping of Definition \ref{def:modified_AC} governs the behavior of our continuum analog at resampling times, while a continuum analog to the scaffolding-and-spindles construction of Figure \ref{fig:discrete_scaff_skewer} governs its behavior in between these times.

\begin{figure}
 \centering
 \scalebox{0.9}{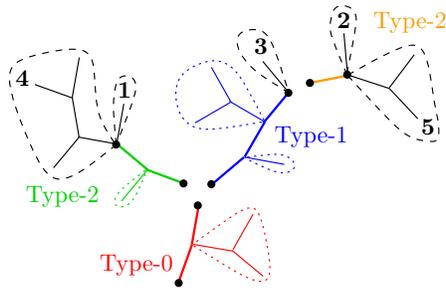}
 \caption{Decomposition of the tree of Figure \ref{fig:decorated_k_tree} into components so that, under the Poissonized Aldous chain, the masses of the spinal and top subtrees in each component evolve according to a Poissonized down-up \oCRP\ with parameters either $\big(\frac12,\frac12\big)$, called type-0; $\big(\frac12,0\big)$, called type-1; or $\big(\frac12,-\frac12\big)$, called type-2.\label{fig:discr_T5_decomp}}
\end{figure}

Note how, in Figure \ref{fig:B_k_tree_proj}, the continuum tree and its $k$-tree projection are partitioned with dashed black lines, so that each external subtree in the left panel (respectively, each leaf edge in the right panel), is grouped together with the internal subtree (resp.\ edge) below it. In that example, leaves 3 and 5 are grouped together with the internal component below them; similarly for leaves 1 and 4; leaf 2 has an internal component to itself; and the root is in an internal component that is not grouped with any external components. Figure \ref{fig:discr_T5_decomp} illustrates the analogous grouping of components for the $5$-tree projection of the combinatorial tree in Figure \ref{fig:decorated_k_tree}. This grouping has to do with case \ref{item:2tdm:moving_down} of the dynamics, when an external component hits mass 0 in a down-move and its label moves down the spine. 
Resampling times are discrete, but the times when a label moves down a spine have accumulation points.

Hence, rather than describing external components and internal components separately, we use a single stochastic process to describe each internal component together with any external components immediately above it. We call these type-0, type-1, and type-2 (interval partition) evolutions, with reference to the number of external components included in the process. Referring back to Figure \ref{fig:B_k_tree_proj}, the external masses containing leaves 3 and 5 and the internal interval partition below them would collectively evolve according to a type-2 evolution; likewise for leaves 1 and 4; leaf 2 would be included in a type-1 evolution; and the root component would evolve by a type-0 evolution.
Type-0 and type-1 evolutions are continuum analogs to Poissonized down-up \oCRP s with $(\alpha,\theta) = \big(\frac12,\frac12\big)$ or $\big(\frac12,0\big)$, respectively. These were constructed via a continuous scaffolding-and-spindles setup in \cite{Paper1-2}.

\section{Structure of this memoir}
\label{sec:intro:overview}

We review type-0 and type-1 evolutions in Chapter \ref{ch:prelim}, before, in Chapter \ref{ch:type-2}, we construct type-2 evolutions, which are the continuum analogs of the Poissonized down-up \oCRP$\big(\frac12,-\frac12\big)$. In Chapter \ref{dePoiss}, we de-Poissonize type-2 evolutions, study unit-mass 2-tree evolutions and prove the statements in Theorem \ref{thm:intro:k_tree} that only require $k=2$. In Chapter \ref{ch:constr} we then define and study $k$-trees and self-similar and unit-mass $k$-tree evolutions based on equipping each branch with a type-0, type-1, or type-2 evolution, resampling and de-Poissonization. In Chapter \ref{ch:consistency}, we study projective consistency and prove the remainder of Theorem \ref{thm:intro:k_tree}. 

At last, Chapter \ref{ch:properties} is devoted to proving Theorem \ref{thm:intro:AD} concerning the Aldous diffusion: the evolving continuum-tree-valued projective limit of projectively consistent $k$-tree evolutions. We derive further properties of the Aldous diffusion in Chapter \ref{chap8}: we study the existence of branch points of higher multiplicity at exceptional times and disprove the strong Markov property, we embed a continuous-time version of the stationary Aldous chain in the Aldous diffusion to prove Theorem \ref{thm:scalinglim}, and we collect open questions related to the Aldous diffusion.

\section{Acknowledgements}

The authors are grateful to Jim Pitman, Nick Bingham, and Quan Shi for their insightful feedback on early drafts of this manuscript. This research has been partially supported by NSF grants DMS-1204840, DMS-1308340, DMS-1612483, and DMS-1855568, UW-RRF grant A112251, EPSRC grant EP/K029797/1, and NSERC grant RGPIN-2020-06907.

\chapter{Preliminaries on type-0 and type-1 interval partition evolutions}
\label{ch:prelim}

In this chapter we recall the constructions and main properties of the type-0 and type-1 interval partition (IP) evolutions introduced in \cite{Paper1-1,Paper1-2}, and we introduce a variant of type-1 evolution that records the mass of a leftmost block separately.

\section{Interval partitions with diversity}
\label{sec:IP}

In this section we specify the metric space of interval partitions with diversity as introduced in \cite{Paper1-0}. Recall from Definition \ref{def:intro:IP} the notion of an interval partition that we use. Also recall the setting  illustrated in Figure \ref{fig:B_k_tree_proj} of a Brownian reduced $k$-tree $R_k = \big(\ft_k,(X_j^{(k)}\!,j\!\in\![k]),(\beta_E^{(k)}\!,E\!\in\!\text{edge}(\ft_k))\big)$, which features interval partitions  $\beta_E^{(k)}$ that are (scaled) $\PDIP\big(\frac12,\frac12\big)$ interval partitions. We introduced $\PDIP\big(\frac12,\frac12\big)$ as the law of the random interval partition $\bar\beta$ formed as the collection of disjoint open intervals in $\{t\in[0,1]\colon B^{\rm br}_t\neq 0\}$, i.e. the set of excursion intervals of a standard one-dimensional Brownian bridge $B^{\rm br}$ as in  Figure \ref{fig:BBr}.  

We begin by extending the notion of (total) diversity $\sD(\beta)$ of an interval partition $\beta$ introduced in \eqref{eq:diversity:intro} to a diversity function 
$(\sD_\beta(t),\,t\ge 0)$.

\begin{definition}\label{def:diversity}
 If it exists, the following limit is called the \emph{diversity} of an interval partition $\beta$ to the left of $t\ge 0$:
 \begin{equation}\label{eq:diversity}
  \sD_\beta(t) := \sqrt{\pi}\lim_{h\downarrow 0}\sqrt{h}\#\{(a,b)\in\beta\colon b\le t,\ |b-a|>h\}.
 \end{equation}
 If this limit exists for all $t\ge 0$, then $\beta$ is said to possess the \emph{diversity property}. We write $\cI$ to denote the set of all interval partitions with this property. As the diversity is constant across any given interval of $\beta$, we will write $\sD_\beta(U)$ for $U\in\beta$ to denote this constant value: $\sD_\beta(U) = \sD_\beta(t)$ for all $t\in U$. We refer to the constant value on $[\|\beta\|,\infty)$ as the \em total diversity of $\beta$ \em and write $\sD_\beta(\infty)=\sD_\beta(\|\beta\|)=\sD(\beta)$.
\end{definition}

In fact, this is the $\alpha=\frac12$ case of what is more generally known as $\alpha$-diversity \cite{CSP}. As this is the only case that we consider, we suppress $\alpha$ in our terminology. 

We define two operations on interval partitions: scaling and concatenation. For $c> 0$ let $c\beta := \{(cx,cy)\colon (x,y)\in\beta\}$ and we define $0\beta=\emptyset$. A collection of interval partitions $(\beta_a)_{a\in\mathcal{A}}$ is \emph{summable} if $\sum_{a\in\mathcal{A}}\IPmag{\beta_a} < \infty$. If $(\mathcal{A},\preceq)$ is a totally ordered set, then we can define 
 $S(a-) := \sum_{b\prec a}\IPmag{\beta_{a}}$ for $a\in\mathcal{A}$,
 and the \emph{concatenation}
 \begin{equation}\label{eq:IP:concat_def}
  \Concat_{a\in\mathcal{A}}\beta_a := \{(S(a-)+x,S(a-)+y)\colon\ a\in\mathcal{A},\ (x,y)\in \beta_a\}.
 \end{equation} 
We also write $\beta\concat\beta^\prime$ to concatenate two interval partitions and simplify notation to $(0,x)\concat\beta^\prime\!:=\!\{(0,x)\}\concat\beta^\prime$ when concatenating a single block into an interval partition. 

\begin{proposition}[\cite{CSP}, Proposition 2.2 of \cite{Paper1-2}]\label{prop:PDIP}
  \begin{enumerate}
    \item[(i)] The random interval \linebreak partition $\bar\beta\sim\PDIP\big(\frac12,\frac12\big)$ a.s.\ possesses the diversity property with positive total diversity, $\sD_{\bar\beta}(\infty) > 0$. 
       Moreover, the diversity process $\big(\sD_{\bar\beta}(t),\,t\in [0,1]\big)$ equals the level-0 local time process of the Brownian bridge, up to scaling.
    \item[(ii)] The ranked block masses of $\bar\beta\sim{\tt PDIP}\big(\frac12,\frac12\big)$ are Poisson--Dirichlet distributed, ${\tt PD}\big(\frac{1}{2},\frac{1}{2}\big)$. 
       Indeed, we can represent ${\tt PDIP}\big(\frac12,\frac12\big)$, as follows. Consider jointly independent $(P_i,\,i\ge 1)\sim{\tt PD}\big(\frac12,\frac12\big)$ and 
       $U_i\sim{\tt Unif}(0,1)$, $i\ge 1$. Then    
       $$\qquad\bar\beta=\Concat_{i\in\mathcal{A}}(0,P_i)\sim{\tt PDIP}\Big(\frac12,\frac12\Big),\quad\mbox{where }\mathcal{A}=\bN,\mbox{ and }i\preceq j\iff U_i\le U_j,$$
       i.e.\ $\bar{\beta}$ consists of intervals of lengths $P_i$, $i\ge 1$, in exchangeable random order. 
    \item[(iii)] Let $M$ be a ${\tt Stable}\big(\frac12\big)$ subordinator with Laplace exponent $\Phi(\theta)=\sqrt{\theta}$, and let $Z\sim{\tt Exponential}(\lambda)$ be independent of $M$. Then
       the partition formed by the jump sizes of $M$ prior to exceeding $Z$ at time $T:=\inf\{t\ge 0\colon M(t)>Z\}$,
       $$\beta:=\{(M(t-),M(t))\colon t\in[0,T),M(t-)<M(t)\},$$
       is a ${\tt PDIP}\big(\frac12,\frac12\big)$ scaled by an independent ${\tt Gamma}\big(\frac12,\lambda\big)$ variable $M(T-)$.
  \end{enumerate}
\end{proposition} 
In the Brownian reduced $k$-tree of Figure \ref{fig:B_k_tree_proj}, blocks $U_1,U_2\in\beta_E^{(k)}$ capture masses of connected components of $\cT\setminus\cR_k$, by construction, while diversities were shown in \cite{PitmWink09} to capture the heights of their attachment points on $\cB_E\subset\cT$ in the BCRT $\cT$, and hence the distance $|\sD_\beta(U_1)-\sD_\beta(U_2)|$ between them.    

We metrize the set $\cI$ of interval partitions with diversity as follows.

\begin{definition} \label{def:IP:metric}
 We adopt the notation $[n] := \{1,2,\ldots,n\}$. 
 For $\beta,\gamma\in \IPspace$, a \emph{correspondence} from $\beta$ to $\gamma$ is a finite sequence of ordered pairs of intervals $(U_1,V_1),\ldots,(U_n,V_n) \in \beta\times\gamma$, $n\geq 0$, where the sequences $(U_j)_{j\in [n]}$ and $(V_j)_{j\in [n]}$ are each strictly increasing in the left-to-right ordering of the interval partitions.

 The \emph{distortion} of a correspondence $(U_j,V_j)_{j\in [n]}$ from $\beta$ to $\gamma$, which we denote by $\dis(\beta,\gamma,(U_j,V_j)_{j\in [n]})$, is defined to be the maximum of the following: 
 \begin{enumerate}[label=(\roman*), ref=(\roman*)]
  \item $\sup_{j\in [n]}|\sD_{\beta}(U_j) - \sD_{\gamma}(V_j)|$,
  \item $|\sD_{\beta}(\infty) - \sD_{\gamma}(\infty)|$,
  \item $\sum_{j\in [n]}|\Leb(U_j)-\Leb(V_j)| + \IPmag{\beta} - \sum_{j\in [n]}\Leb(U_j)$, \label{item:IP_m:mass_1}
  \item $\sum_{j\in [n]}|\Leb(U_j)-\Leb(V_j)| + \IPmag{\gamma} - \sum_{j\in [n]}\Leb(V_j)$. \label{item:IP_m:mass_2}
 \end{enumerate}
 For $\beta,\gamma\in\IPspace$ we define
 \begin{equation}\label{eq:IP:metric_def}
  \dI(\beta,\gamma) := \inf \dis\big(\beta,\gamma,(U_j,V_j)_{j\in [n]}\big),
 \end{equation}
 where the infimum is over all correspondences from $\beta$ to $\gamma$.
\end{definition}

\begin{proposition}[Theorem 2.4 of \cite{Paper1-0}]\label{prop:IPspace:Lusin}
 The map $\dI$ is a metric, and $(\cI,\dI)$ is a Lusin space, i.e.\ homeomorphic to a Borel subset of a compact metric space.
\end{proposition}

\section{Definitions and properties of type-0 and type-1 evolutions}
\label{sec:prelim:type01_def}

We will specify two diffusions on $(\cI,\dI)$, called type-1 and type-0 evolutions, via their semigroups. These are the $\alpha=\frac12$ case of the diffusions introduced in \cite{Paper1-2}.

\begin{figure}[t]
 \centering
 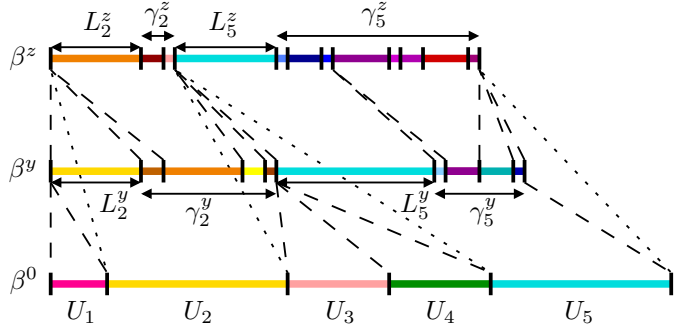
 \caption{Illustration of the transition kernel $\kappa^{1}_y$: $\beta^0$ has five blocks $U_1,\ldots,U_5$. Some blocks contribute $\emptyset$ for time $y$, here 
   $U_1$, $U_3$ and $U_4$; others yield non-trivial partitions, here $U_2$ and $U_5$, hence $\beta^y=(0,L_2^y)\concat\gamma_2^y\concat(0,L_5^y)\concat\gamma_5^y$.   
   The semigroup property requires consistency of the dotted transition from 0 to $z$ and the composition of the dashed transitions from $0$ to $y$ and from $y$ to $z$.\label{fig:semigroup}}
\end{figure}

The semigroup $(\kappa_y^{1},\,y\ge0)$ for type-1 evolutions satisfies a \emph{branching property}: given any interval partition $\beta$, the blocks $U\in\beta$ will give rise to independent random interval partitions $\gamma_U$ (possibly empty) at time $y>0$, and $\kappa_y^{1}(\beta,\,\cdot\,)$ will be the distribution of their concatenation. This is illustrated in Figure \ref{fig:semigroup}.

Specifically, fix $b,r>0$. Let $B_{r} \sim \GammaDist\big(\frac12,r\big)$, $\bar\gamma \sim \PDIP\big(\frac12,\frac12\big)$, and let $L_{b,r}$ be a $(0,\infty)$-valued random variable with Laplace transform
\begin{equation}
  \bE\left[e^{-\lambda L_{b,r}}\right] = \sqrt{\frac{r+\lambda}{r}}\frac{\exp(br^2/(r+\lambda))-1}{\exp(br)-1},
  \label{LMBintro}
\end{equation}
all assumed to be independent. Consider the following distribution on $\cI$:
\begin{equation}\label{eq:mu}
 \mu_{b,r}(\,\cdot\,) = e^{-br}\delta_\emptyset(\,\cdot\,) + (1-e^{-br})\bP\big\{\big(0,L_{b,r}\big) \concat B_{r}\bar\gamma\in\cdot\,\big\}.
\end{equation}
For $y>0$ and $\beta$ any interval partition, let $\kappa_y^{1}(\beta,\,\cdot\,)$ denote the law of
\begin{equation}
 \label{eq:intro:transn_1}
 \Concat_{U\in\beta}\gamma_U \quad \text{where} \quad \gamma_U\sim \mu_{\Leb(U),1/2y}\text{ independently for each }U\in\beta.
\end{equation}
Also set $\kappa_0^1(\beta,\,\cdot\,)=\delta_\beta(\,\cdot\,)$. Similarly, let $\kappa_0^0(\beta,\,\cdot\,)=\delta_\beta(\,\cdot\,)$ and $\kappa_y^{0}(\beta,\,\cdot\,)$ the law of 
\begin{equation}
 \label{eq:intro:transn_0}
 B\bar\gamma \concat \Concat_{U\in\beta}\gamma_U,
\end{equation}
where the $\gamma_U$ are as in \eqref{eq:intro:transn_1} jointly independent with $B \sim \GammaDist\big(\frac12,\frac{1}{2y}\big)$ and $\bar\gamma\sim\PDIP\big(\frac12,\frac12\big)$, for each $y>0$.

A Markov process $(\beta^y,y\ge0)$ is said to be \em self-similar\em, or \emph{1-self-similar} in the sense of Lamperti \cite{Lamperti72}, if it has the same semigroup as the scaled process $(c\beta^{y/c},y\ge 0)$ for all $c>0$.

\begin{proposition}[Theorems 1.2--1.3 of \cite{Paper1-2}]
\label{prop:type01:diffusion}
 The maps $\beta\mapsto \kappa_y^{1}(\beta,\cdot\,)$, $y\ge 0$, are weakly continuous and form the transition semigroup of a self-similar path-continuous Hunt process $(\beta^y,\,y\ge 0)$ on $(\IPspace,\dI)$, and likewise for $\kappa_y^{0}$.
\end{proposition}

It is not at all obvious why one would choose these transition kernels, nor that they satisfy the semigroup property. In fact, these kernels fall out of a Poissonian construction that we will describe in Sections \ref{sec:prelim:SSS}--\ref{sec:prelim:clade_type01}.

We refer to a diffusion with transition semigroup $(\kappa_y^{0},\,y\ge0)$ as a \emph{type-0 evolution}. We refer to a diffusion with semigroup $(\kappa_y^{1},\,y\ge0)$ as an \emph{$\cI$-valued type-1 evolution}. Previously \cite{Paper1-2}, we have simply called this a ``type-1 evolution,'' but for our present purpose of constructing continuum-tree-valued diffusions, we wish to set up type-1 evolutions on an equivalent state space.

It is not hard to show that the concatenations in \eqref{eq:intro:transn_1}--\eqref{eq:intro:transn_0} are almost surely finite, in the sense that all but finitely many of the components being concatenated will be null partitions \cite[Lemma 6.1]{Paper1-1}. Even so, as each \PDIPAT\ has infinitely many blocks, so too do type-0 and $\cI$-valued type-1 evolutions. More formally, given such a process $(\beta^y,y\ge 0)$ on the event $\{\beta^z\neq\emptyset\}$ for some $z>0$, the partition $\beta^z$ a.s.\ has infinitely many blocks. Moreover, there is a.s.\ no rightmost block but rather $\beta^z$ has infinitely many blocks to the right of $\IPmag{\beta^z}-\epsilon$, for every $\epsilon\!>\!0$. 
 However, as only finitely many of the $\gamma_U$ in \eqref{eq:intro:transn_1} are non-empty, in the type-1 case, $\beta^z$ comprises a finite alternating sequence of the leftmost blocks of those $\gamma_U$ and rescaled \PDIPAT. In particular, $\beta^z$ a.s.\ has a leftmost block when $\beta^z\neq\emptyset$. 

We define 
\begin{equation}\label{eq:type1space_def}
\begin{split}
 &\cJ^\bullet:=\{(m,\gamma)\in[0,\infty)\times\cI\colon m>0\mbox{ or }\gamma\mbox{ has no leftmost block}\},\\
 &d^\bullet((m_1,\gamma_1),(m_2,\gamma_2)) = |m_1-m_2|+d_\cI(\gamma_1,\gamma_2),
\end{split}
\end{equation}
and consider the continuous bijection $\varphi(m,\gamma)=(0,m)\concat\gamma$ from $\cJ^\bullet$ to $\cI$, which has a (discontinuous) measurable inverse. 

We define a (pair-valued) \emph{type-1 evolution} to be the image of an $\cI$-valued type-1 evolution under $\varphi^{-1}$, as a process on $\cJ^\bullet$. It follows immediately from Proposition \ref{prop:type01:diffusion} that this is a Markov process. We will expand on this result in Corollary \ref{corpairtype1}.

Squared Bessel processes are $[0,\infty)$-valued diffusions described by the SDE 
\begin{equation}\label{eq:BESQ}
 dZ(y) = r\,dy + 2\sqrt{Z(y)}dB(y), \quad 0\le y\le \zeta,
\end{equation}
where $r\in\BR$ is a parameter, $B$ is standard one-dimensional Brownian motion and either $\zeta = \infty$ if $r > 0$, or $\zeta = \inf\{y\ge0\colon Z(y)\le 0\}$ if $r\le 0$. In the latter case, we adopt the convention that the process is absorbed at 0, though this convention is not universal \cite{GoinYor03}. For $x\ge 0$ we write $\BESQ_x(r)$ to denote the law of such a process started from initial state $x$. 
When $r = 0$, this is a continuous-state branching process known as the Feller diffusion. Otherwise, this can be viewed as a branching process with immigration (when $r >0$; see \cite{KawaWata71}) or emigration (when $r <0$). See \cite[Chapter XI]{RevuzYor} or \cite{GoinYor03} for more discussion of these diffusions.

\begin{proposition}[Theorem 1.4 of \cite{Paper1-2}]\label{type1totalmass}
 The total mass process associated with a type-0 evolution $(\beta^y,\,y\ge0)$ is $(\|\beta^y\|,\,y\ge0) \sim \besq_{\|\beta^0\|}(1)$, while the total mass of a type-1 evolution $((m^y,\gamma^y),\,y\ge0)$ is $(m^y\!+\!\|\gamma^y\|,\,y\ge0)\sim\besq_{m^0+\|\gamma^0\|}(0)$. In particular, the type-1 evolution a.s.\ is absorbed at $(0,\emptyset)$ in finite time, whereas the type-0 evolution visits $\emptyset$ but is reflected rather than absorbed.
\end{proposition}

We refer to the absorption time of a type-1 evolution as its \emph{degeneration time}, and we say that a type-0 evolution never degenerates, or that it has degeneration time $D=\zeta=\infty$.


\begin{proposition}[Theorem 1.5 of \cite{Paper1-2}]\label{prop:01:pseudo}
 Let $A\sim \BetaDist\big(\frac12,\frac12\big)$, $\bar\beta\sim\PDIP\big(\frac12,\frac12\big)$, and let $Z_0\sim\besq(1)$ and $Z_1\sim \besq(0)$ with arbitrary initial distributions, all jointly independent. If $(\beta^y,\,y\ge0)$ is a type-0 evolution starting from $\beta^0$ distributed as the independent random multiple $Z_0(0)\bar\beta$ of $\bar\beta$, then 
$\beta^y$ is distributed as the independent random multiple $Z_0(y)\bar\beta$ of $\bar\beta$ for each $y>0$. Similarly, if $((m^y,\gamma^y),\,y\ge0)$ is a type-1 evolution with $(m^0,\gamma^0)\stackrel{d}{=}(Z_1(0)A,Z_1(0)(1-A)\bar\beta)$, then $(m^y,\gamma^y)\stackrel{d}{=}(Z_1(y)A,Z_1(y)(1-A)\bar\beta)$ for each $y>0$.
\end{proposition}

For this reason, the law of any independently randomly scaled $\PDIP\big(\frac12,\frac12\big)$ is called a \emph{pseudo-stationary law for the type-0 evolution}, and similarly, the law of $(A,(1-A)\bar\beta)$ times an independent scaling random variable is a \emph{pseudo-stationary law for the type-1 evolution}. This proposition has the following key special case.

\begin{proposition}[Proposition 4.1 of \cite{Paper1-2}]\label{prop:01:pseudo_g}
 Fix $\lambda>0$. Let $A$ and $\bar\beta$ be as in Proposition \ref{prop:01:pseudo} and, independently, consider $M_i\sim\GammaDist\big(\frac{1+i}{2},\lambda\big)$ for $i=0,1$. If $(\beta^y,\,y\ge0)$ is a type-0 evolution with $\beta^0\stackrel{d}{=}M_0\bar\beta$, then $\beta^y\stackrel{d}{=}(2y\lambda+1)M_0\bar\beta$. \linebreak If $((m^y,\gamma^y),\,y\ge0)$ is a type-1 evolution with $(m^0,\gamma^0)\stackrel{d}{=}M_1(A,(1-A)\bar\beta)$, then given that this process does not degenerate prior to time $y>0$, the conditional law of $(m^y,\gamma^y)$ equals the (unconditional) law of $(2y\lambda+1)M_1(A,(1-A)\bar\beta)$. Moreover, 
 \begin{equation}\label{eq:pseudo:type_1_survival}
  \Pr\{(m^y,\gamma^y) = (0,\emptyset)\} = \frac{1}{2y\lambda+1}.
 \end{equation}
\end{proposition}



\section{Scaffolding, spindles, and skewer}
\label{sec:prelim:SSS}

\begin{figure}
 \centering
  \includegraphics[width = 4.8in,height=2.88in]{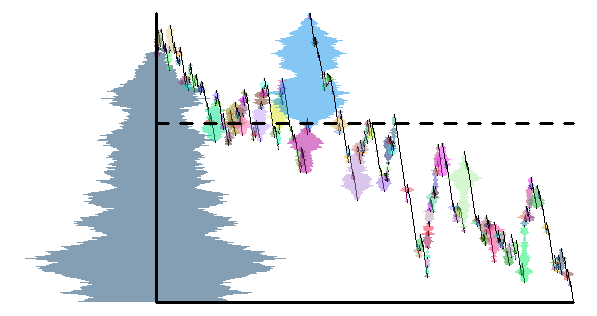}\\
  \hspace{.22in}\includegraphics[width = 2.2in,height=.15in]{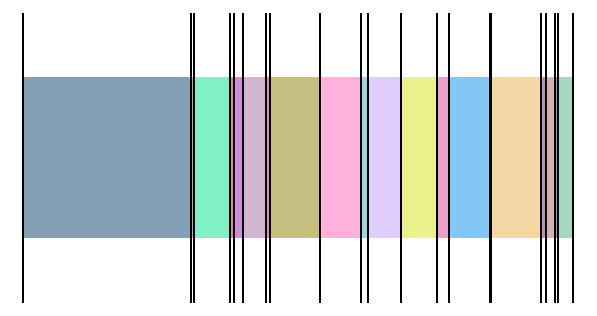}\hspace{.5in}
  \caption{A simulated approximation of the scaffolding and spindles construction of a type-1 interval partition evolution in the simplest case: starting from a one block partition $\{(0,a)\}$. In this case, the scaffolding and spindles comprise a \emph{clade}; see Section \ref{sec:prelim:clade_type01}. Compare this to the discrete rendering in Figure \ref{fig:discrete_scaff_skewer}.\label{fig:scaff_spind_skewer_1}}
\end{figure}

In this section and the next, we recall from \cite{Paper1-1,Paper1-2} the setup of scaffolding, spindles and the skewer map, as well as associated constructions of type-0 and type-1 interval partition evolutions. In this memoir, this setup serves three purposes. 
The first is to stress parallels to the discrete regime introduced in Section \ref{sec:discrete_to_cts}. Figure \ref{fig:scaff_spind_skewer_1} depicts a simulated approximation to the 
construction that we will undertake here, in the continuum. 
The second is to acknowledge the role this setup has played in proving the results stated in the previous section. Indeed, transition semigroups provide an efficient way to introduce type-0 and type-1 evolutions and to state the main results that also form the interface for the use of type-0 and type-1 evolutions in the construction and study of $k$-tree evolutions and the Aldous diffusion in Chapters \ref{ch:constr}--\ref{ch:properties}. However, we have not been able to develop a theory of type-0 and type-1 evolutions directly from the semigroups or even to show that they are semigroups without this setup. Last, but not least, we also need type-2 evolutions. Their constructions in Sections \ref{sec:sym} and \ref{sec:interweaving} will build explicitly on the constructions for types 0 and 1 that we recall and enhance here.    
To begin with, here is an informal introduction to the key terminology, to be made formal later.\medskip 

\emph{Scaffolding.} This is a random c\`adl\`ag process $X\colon [0,T]\to \BR$ 
formed by concatenating path segments of stopped spectrally positive \StableA\ L\'evy processes. In Figure \ref{fig:scaff_spind_skewer_1}, (an approximation of) the scaffolding is plotted in thin, steeply downwards-sloping black lines.

\emph{Spindles.} These are random $[0,\infty)$-valued excursions, typically continuous, that mark the jumps of the scaffolding in such a way that each jump of height $z$ is marked by an excursion with lifetime $z$. In particular, our construction uses \besqA\ excursions. In Figure \ref{fig:scaff_spind_skewer_1}, the spindles are depicted as laterally symmetric, shaded blobs inscribed into the jumps. The fluctuating width of the blob, as it progresses up from the bottom towards the top of the jump, depicts the fluctuating value of the excursion function, starting from and ending at width 0.

\emph{Skewer.} This setup gives rise to an interval-partition-valued process via the skewer map: in Figure \ref{fig:scaff_spind_skewer_1}, we see a horizontal dashed line cutting through the scaffolding-and-spindles plot. The interval partition shown below the plot is made up of a single block corresponding to each spindle bisected by the dashed line, with the color of the block matching that of the corresponding spindle and the size of the block equaling the width of the spindle at the point at which it is bisected. In order to obtain a continuous process of interval partitions, rather than a single partition as in the figure, we begin with a skewer at height 0, then move the dashed line continuously up the page. Color in the figure is just for illustration.

Our formal setup diverges from the heuristic description above in one key respect: rather than beginning with a scaffolding process and adding spindle markings, we find it more parsimonious to begin with Poisson random measures of spindles and then associate an intrinsic scaffolding that has the spindle lifetimes as jump heights. 

Before we turn to the formal setup, let us further motivate the terminology. The name ``spindle'' is in recognition of their visual appearance in illustrations such as Figure \ref{fig:scaff_spind_skewer_1} after we chose the laterally symmetric shapes for aesthetic reasons \cite{Paper1-1}. The vertical placement of each spindle of the Poisson random measure is implicit, and the ``associated scaffolding'' makes this explicit -- we can view the spindles as being placed onto the scaffolding. 
The ``skewer,'' at a given level, is pushed through the spindles from the left. Each spindle straddling the level has a certain width and will occupy a corresponding length of skewer. We leave no gaps on the skewer. Having a designated name ``scaffolding'' allows us to refer unambiguously to horizontal ``scaffolding time'' and vertical ``scaffolding level.'' The latter is ``skewer time'' and indeed interval partition evolution time and also relates to ``spindle time.'' 

\medskip

\emph{Spindles, formally.} We define the set $\cE$ of spindles as the subset 
\begin{equation}
 \cE := \left\{f\in\cD\ \middle| \begin{array}{c}
   \displaystyle \exists\ z\in(0,\infty)\textrm{ such that }\restrict{f}{(-\infty,0)\cup [z,\infty)} = 0,\\[0.2cm]
   \displaystyle \text{ and }f\text{ is positive and continuous on }(0,z)
  \end{array}\right\}\label{eq:cts_exc_space_def}
\end{equation}
of the space $\cD$ of c\`adl\`ag functions from $\BR$ to $[0,\infty)$. In words, spindles $f\in\cE$ are positive c\`adl\`ag excursions whose only jumps may be at birth and death.
We refer to spindles that have a jump at birth and/or death as \em broken spindles\em. For any spindle $f\in\Exc$, we define the \emph{lifetime} or \emph{absorption time} by 
\begin{equation}\label{eq:exc:lifetime}
 \zeta(f) := \inf\{s>0\colon f(s) = 0\}.
\end{equation}

\begin{lemma}[Equation (13) in \cite{GoinYor03}]\label{lem:BESQ:length}
  Let $Z=(Z(y),\,y\ge 0)\sim\BESQ_x(-1)$. Then the lifetime $\zeta(Z) = \inf\{y\ge0\colon Z(y) = 0\}$ has the same distribution as $x/2G$ where $G\sim\GammaDist(\frac32,1)$.
\end{lemma}

Pitman and Yor \cite{PitmYor82} gave a general construction of $\sigma$-finite excursion measures, which applies to $\BESQ(-1)$ even though there is no It\^o excursion measure as 0 is not an entrance boundary for $\BESQ(-1)$. In their terminology, the construction for the 0-diffusion $\BESQ(-1)$ uses the associated $\uparrow$-diffusion $\BESQ(5)$ up 
to a \emph{first passage time}, which we denote by $H^b\colon\cE \to [0,\infty]$ via $H^b(f)\!=\!\inf\{s\!\ge\! 0\colon f(s)\!=\!b\}$, $b>0$. 

\begin{lemma}[Description (3.1) in \cite{PitmYor82}]\label{lem:BESQ:existence}
  There is measure $\Lambda$ on $\cE$ such that: 
  \begin{enumerate}\item[(i)] $\Lambda\{f\!\in\!\cE\colon f(0)\!\neq\! 0\}=0$, and $\Lambda\{H^b\!<\!\infty\}=b^{-3/2}$, $b\!>\!0$,
  \item[(ii)] under $\Lambda(\,\cdot\,|\,H^b\!<\!\infty)$, the restricted canonical process 
  $f|_{[0,H^b]}$ is a $\BESQ_0(5)$ stopped at its first passage at $b$, independent of $f(H^b+\cdot\,)\sim\BESQ_b(-1)$.
  \end{enumerate} 
\end{lemma}  
Following \cite{Paper0,Paper1-2}, we scale the measure $\Lambda$ of Lemma \ref{lem:BESQ:existence} and set
$$ \mBxc = \frac{3}{2\sqrt{\pi}}\Lambda.$$



\emph{Scaffolding, formally.}
Let $\fN$ denote a Poisson random measure on $[0,\infty)\times\cE$ with intensity measure ${\rm Leb}\otimes\mBxc$, abbreviated $\fN\sim{\tt PRM}({\rm Leb}\otimes\nu_{\tt BESQ})$. This is a point process in which spindles occur at a dense set of times, but spindles with lifetime greater than $\epsilon$, for any $\epsilon>0$, occur at discrete times. More precisely, mapping spindles onto their lifetimes yields a Poisson 
random measure on $[0,\infty)\times(0,\infty)$ whose intensity measure ${\rm Leb}\otimes\nu_{\tt BESQ}(\zeta\in\,\cdot\,)$ we can use to build a L\'evy process that incorporates
all spindle lifetimes as jump heights in a compensated limit. Keep this in mind as a key example for the following definition.

\begin{definition}
\label{def:scaff}
Given a point measure $N$ on $[0,\infty)\times\cE$, we define
\begin{equation}\label{eq:scaff:def}
 \xi_{N}(t) := \lim_{z\downto 0} \left(\int_{[0,t]\times\{g\in\Exc\colon\zeta(g) > z\}}\zeta(f)dN(u,f)-\frac{3tz^{-1/2}}{\pi\sqrt{2}}\right)\quad \text{for }t\ge0.
\end{equation}
If this limit exists for all $t\ge0$, then we call this the \emph{scaffolding} associated with $N$ and we abbreviate $\xi(N) := (\xi_{N}(t),\,t\ge0)$ and often write $X:=\xi(N)$.

If $N$ is supported on a bounded time interval $[0,T]\times\cE$, then the \emph{length} $\len(N)$ is defined to be the least such $T$. In this case, we only require the limit to converge for $t\in [0,\len(N)]$ in order to call $X:=\xi(N) := (\xi_{N}(t),\,t\in [0,\len(N)])$ the scaffolding.
\end{definition}

The term outside of the integral in \eqref{eq:scaff:def} equals the expected value of the integral if we substitute $\fN$ for $N$.

\begin{lemma}[Proposition 2.12 of \cite{Paper1-1}, (2.11) of \cite{Paper1-2}]\label{lm:scaff}
 For $\fN$ as above, the associated scaffolding $\fX:=\xi(\fN)$ is a spectrally positive \StableA\ L\'evy process with L\'evy measure and Laplace exponent given by
 \begin{equation}\label{eq:scaff:Levy_Laplace}
  \Pi(dx) = \nu_{\tt BESQ}(\zeta\in dx) = \frac{3}{2\sqrt{2}\pi}x^{-5/2}dx\ \ \text{and}\ \ \psi(\lambda) = \sqrt{\frac{2}{\pi}}\lambda^{3/2}.
 \end{equation} 
\end{lemma}

\emph{Skewer, formally.} 
  Let $N=\sum_{i\in I}\delta(t_i,f_i)$ be a point measure on $[0,\infty)\times\cE$ with associated scaffolding $X=\xi(N)=(\xi_N(t),\,t\in[0,\len(N)])$. By construction,
  each $t_i$, $i\in I$, is a jump time of $X=\xi(N)$ of jump height $\zeta(f_i)$, and we associate spindle times $s\in[0,\zeta(f_i)]$ with the scaffolding levels 
  $\xi_N(t_i-)+s\in[\xi_N(t_i-),\xi_N(t_i)]$ crossed by the jump at scaffolding time $t_i$. This means that $f_i(y-\xi_N(t_i-))$ is associated with level $y\in\mathbb{R}$, which is 
  positive if $y$ is crossed at time $t_i$ and zero otherwise. If $f_i$ is continuous, this quantity equals $f_i((y-\xi_N(t_i-))-)$, but in order to achieve the desired effect in all cases 
  (including cutoff point measures needed for Lemma \ref{lem:type1Markov} where spindle lifetimes are cut short by c\`adl\`ag jumps down to 0), we consider 
  $\max\{f_i(y-\xi_N(t_i-)),f_i((y-\xi_N(t_i-))-)\}$.

\begin{definition}
 \label{def:skewer} Let $N=\sum_{i\in I}\delta(t_i,f_i)$ be a point measure on $[0,\infty)\times\cE$ with scaffolding $\xi(N)$ and $y\in\mathbb{R}$. Then the 
  \emph{aggregate mass} (sum of spindle widths) of $N$ at scaffolding level $y$ up to scaffolding time $t\in[0,\len(N)]$ is
 \begin{align}\label{eq:skewer_def}
  M^y_{N}(t) :=& \sum_{i\in I\colon t_i\le t} \max\Big\{ f_i\big(y-\xi_{N}(t_i-)\big),f_i\big((y-\xi_{N}(t_i-))-\big)\Big\}\\
                      =& \int_{[0,t]\times\Exc} \max\Big\{  f\big(y-\xi_{N}(u-)\big),f\big((y-\xi_{N}(u-))-\big)\Big\}dN(u,f),\nonumber
\end{align}
and the \emph{skewer} of $N$ at level $y$ is the interval partition formed by the range of $M_N^y$
\begin{equation}\label{eq:skewer_def2}
  \skewer(y,N) := \left\{\left(M^y_{N}(t-),M^y_{N}(t)\right)\!\colon t\in[0,\len(N)],\ M^y_{N}(t-) < M^y_{N}(t)\right\}\!.\!\!\!\!
 \end{equation}
 We abbreviate $\skewerbar(N) := \big(\skewer(y,N),\, y\geq 0\big)$.
\end{definition}

A simulation of this construction is depicted in Figure \ref{fig:scaff_spind_skewer_1}. Compare this to the analogous discrete construction in Section \ref{sec:discrete_to_cts} and the depiction in Figure \ref{fig:discrete_scaff_skewer}.

In \cite[Definitions 2.13 and 3.4]{Paper1-1} we defined a measurable space $\cN_{\rm fin}^{\rm sp,*}$ of point measures $N$ on $[0,\infty)\times\cE$, supported on bounded time intervals $[0,T]\times\cE$, for which $\skewerbar(N)$ is well-defined and $d_\cI$-continuous, and space-time local times of $\xi(N)$ equal diversities for all scaffolding times and levels: $\sD_{\skewer(y,N|_{[0,t]})}(\infty)=\ell_{\xi(N)}^y(t)$.  For brevity, in this memoir we simply denote this space by $\cN$.

Consider the restriction of $\fN$ to a bounded time interval $\restrict{\fN}{[0,T]\times\Exc}$ for some random time $T\in (0,\infty)$. From \cite[Proposition 3.8]{Paper1-1} applied with $\alpha=\frac12$, $q=1$, the process $\skewerbar\big(\restrict{\fN}{[0,T]\times\Exc}\big)$ is a.s.\ continuous on $(\cI,\dI)$, and is in fact a.s.\ H\"older-$\theta$ for any $\theta \in \big(0,\frac14\big)$. Moreover, $\restrict{\fN}{[0,T]\times\Exc}$ belongs to $\cN$ almost surely.

\section{Clades and the construction of type-0 and type-1 evolutions}
\label{sec:prelim:clade_type01}

We stated at the end of Section \ref{sec:discrete_to_cts} that the tree evolving under the Aldous chain can be decomposed, until a stopping time, into so-called type-0, type-1, and type-2 components that evolve, under spinal projections, as down-up ordered Chinese restaurant processes with parameters $\big(\frac12,\frac12\big)$, $\big(\frac12,0\big)$, and $\big(\frac12,-\frac12\big)$, respectively. The scaffolding-and-spindles construction initiated above is based on imagining a genealogy among the tables in the restaurant: whenever a new table appears (``is born'') in the restaurant, it is the ``child'' of the table immediately to its left at its time of birth. Then the scaffolding defined above is (the continuum version of) a contour process representation of this family tree of tables (or forest, in the case of having multiple tables at time 0 rooting multiple trees).

In evolutionary biology, a ``clade'' is the set of all descendants of a single individual. The ``clades'' defined below are continuum analogues to genealogical clades in our (imagined) genealogy among tables in the restaurant.

\begin{definition}\label{def:clade_from_PRM}
 Fix $x>0$ and consider $\fN\sim\PRM({\rm Leb}\otimes\mBxc)$ as above and an independent spindle $\ff\sim\besq_{x}(-1)$. A \emph{clade of initial mass $x$} is a random point measure $\fn\in\cN$, distributed as
\begin{equation}
\label{eq:clade_from_PRM}
\begin{split}
 \fn = \clade(\ff,\fN) &:= \delta(0,\ff)+\Restrict{\fN}{(0,T_{-\zeta(\ff)}(\fN)]\times\cE},\\
 \text{where} \quad T_{-y}(\fN) &:= \inf\{t\ge 0\colon\xi_{\fN}(t)=-y\}.
\end{split}
\end{equation}
\end{definition}

This construction gives the continuum analogue of a clade of tables in the down-up \oCRP: the spindle $\ff$ (``broken'' because it starts positive, rather than entering continuously from 0) represents an ancestor table started with positive population. The restricted point measure of spindles $\Restrict{\fN}{(0,T_{-\zeta(\ff)}(\fN)]\times\cE}$ describes the descendant tables. The scaffolding and spindles in Figure \ref{fig:scaff_spind_skewer_1} comprise a single clade, with the large grey broken spindle $\ff$ on the left, followed by a \StableA\ scaffolding $\xi(\fn)$ marked by descendant spindles, stopped when the scaffolding reaches level 0.

In the following clade construction and elsewhere, the notion of ``concatenation,'' denoted by $\concat$, is in the sense of excursion theory. Let $(N_a)_{a\in\mathcal{A}}$ denote a family of point processes of spindles indexed by a totally ordered set $(\mathcal{A},\preceq)$. For the purpose of the following, set
 \begin{equation}\label{eq:length_partial_sums}
  S(a) := \sum_{b\preceq a}\len(N_b) \quad \text{and} \quad S(a-) := \sum_{b\prec a}\len(N_b) \quad \text{for each }a\in \mathcal{A}.
 \end{equation}
 If $S(a-) < \infty$ for every $a\in \mathcal{A}$, 
 then we define the \emph{concatenation} of $(N_a)_{a\in\mathcal{A}}$ to be the point measure
 \begin{equation}\label{eqn:concat}
   \Concat_{a\in\mathcal{A}}N_a := \sum_{a\in\mathcal{A}}\int\Dirac{S(a-)+t,f}dN_a(t,f).
 \end{equation}

\begin{proposition}[Theorem 1.8 of \cite{Paper1-2}]\label{def:type01_meas} 
  Fix $\beta\in\cI$. Consider a family $N_U$ of independent clades with initial mass $\Leb(U)$, $U\in\beta$. Denote by $\fP_\beta^1$ the law of the  \em 
  type-1 point measure\em
 \begin{equation*}
  \fN_\beta := \Concat_{U\in\beta} N_U.
 \end{equation*}
  Then $\skewerbar(\fN_\beta)$ is an $\cI$-valued type-1 evolution with initial state $\beta$.
 %
\end{proposition}

\begin{corollary}\label{type01plustype1}
  Consider 
  $\beta_1,\beta_2\in\cI$ and two independent type-1 point measures $\fN_{\beta_1}$ and $\fN_{\beta_2}$ as in Proposition \ref{def:type01_meas}. Then $\fN_{\beta_1}\concat\fN_{\beta_2}$ is also a type-1 point measure. In particular,
  the concatenation of the associated skewer processes, $(\beta_1^y\concat\beta_2^y,y\ge 0)$, is an $\cI$-valued type-1 evolution starting from $\beta:=\beta_1\concat\beta_2$.
\end{corollary}

Now, consider the point measure $\cev\fN$ on $(-\infty,0)\!\times\!\cE$ formed by concatenating a sequence of independent copies of $\restrict{\fN}{(0,T_{-1}(\fN)]\times\cE}$, with each copy being concatenated to the left of the previous copies. We denote by $\cev\cN$ the space of all point measures of spindles constructed in this manner, by concatenating a sequence of point measures in $\cN$ whose associated scaffoldings, as in Definition \ref{def:scaff}, are first-passage descents to $-1$. In particular, $\cev\fN$ is a random element of this space. 

We adapt \eqref{eq:scaff:def}, \eqref{eq:skewer_def} and \eqref{eq:skewer_def2} to this setting. We define pre-0 scaffolding
\begin{equation}\xi_{\cev\fN}(t) := \lim_{z\downto 0}\!\left(\!-\int_{(t,0)\times\{g\in\Exc\colon\zeta(g) > z\}}\life(f)d\cev\fN(s,f) + \frac{3|t|z^{-1/2}}{\pi\sqrt{2}}\right)\!,\quad t\le0,
 \label{eq:JCCP_def_type0}
\end{equation}
and set $\xi\big(\cev\fN\big) := (\xi_{\cev\fN}(t),\,t\le 0)$. Informally, this is a spectrally positive \StableA\ first-passage descent from $\infty$ down to 0, arranged to arrive at 0 at time zero. For $y\ge 0$, we write as $T_y\big(\cev{\fN}\big)=\inf\{t\le 0\colon\xi_{\cev{\fN}}(t)=y\}$ the pre-0 downward first passage time at level $y$, and for $t\in[T_y(\cev\fN),0]$, the pre-0 aggregate mass process as
$$M_{\cev{\fN}}^y(t):=\int_{[T_y(\cev\fN),t]\times\cE} \max\Big\{f\big(y-\xi_{\cev\fN}(u-)\big), f\big((y-\xi_{\cev\fN}(u-))-\big)\Big\}d\cev\fN(u,f).$$
Finally, we define the associated skewer at level $y$
$$ \skewer(y,\cev\fN) := \left\{\left(M^y_{\cev\fN}(t-),M^y_{\cev\fN}(t)\right)\!\colon t\in[T_y(\cev\fN),0],\ M^y_{\cev\fN}(t-) < M^y_{\cev\fN}(t)\right\}\!.$$

\begin{construction}[Type 0]\label{type0:construction} Let $\beta\in\cI$. Consider $\cev\fN$ as above and, independently, $\fN_\beta\sim\fP_\beta^1$ as in Proposition \ref{def:type01_meas}. We denote by $\fP^0_\beta$ the distribution of $(\cev\fN,\fN_\beta)$. We define an $\cI$-valued evolution $(\beta^y,\,y\ge 0)$ as
  $$\beta^y:=\skewer(y,\cev\fN)\concat\skewer(y,\fN_\beta),\quad y\ge 0.$$ 
\end{construction}

  This was proposed as a construction of type-0 evolutions in \cite[Remark 3.9]{Paper1-2}. In \cite[Definition 3.8 and Proposition 3.10]{Paper1-2}, we constructed a type-0 evolution 
  with time interval $[0,j]$ starting from $\fN\sim{\tt PRM}({\rm Leb}\otimes\nu_{\tt BESQ})$, for each $j\ge 0$. 
  Specifically, we take the skewer process of $\restrict{\fN}{[0,T_{-j}(\fN)]}$ on the associated scaffolding shifted up by $j$ so as to yield a first passage descent from $j$ to $0$. 
  Since the point measure $\restrict{\cev\fN}{[T_j(\cev\fN),0]}$ in the setting of Construction \ref{type0:construction} also has a first passage descent from $j$ to $0$ as its 
  associated scaffolding, we obtain the following result.

\begin{proposition}\label{prop:type0:construction} Let $\beta\in\cI$. Then the $\cI$-valued process resulting from Construction \ref{type0:construction} is a type-0 evolution starting from $\beta$.
\end{proposition}

  Recall from Section \ref{sec:prelim:type01_def} the definition of a (pair-valued) type-1 evolution starting from $(x,\gamma)\in\cJ^\bullet$. In the setting of Proposition
  \ref{def:type01_meas}, for the corresponding $\cI$-valued initial state $\beta=(0,x)\concat\gamma$, we naturally split off the first clade and write $\fN_\beta$ in the  
  form $\fN_\beta=N_{(0,x)}\concat\fN_\gamma$, where $N_{(0,x)}$ is a clade of initial mass $x$ as in Definition \ref{def:clade_from_PRM}. The identification of first passage 
  descents from $\cev\fN$ and $\fN$ above Proposition \ref{prop:type0:construction} yields the following construction of type-1 evolutions that replaces 
  $\clade(\ff,\fN)$ by
  \begin{equation}
    \label{eq:clade_from_descent}
    \clade(\ff,\cev{\fN}):=\delta(0,\ff)+\ShiftRestrict{\cev{\fN}}{[T_{\zeta(\ff)}(\cev{\fN}),0)\times\cE}\stackrel{d}{=}\clade(\ff,\fN),
  \end{equation}
where the notation in the middle expression is in the following sense. We define the \emph{shifted restriction} of a point measure $N$, denoted by $\shiftrestrict{N}{[a,b]\times\Exc}$ to be the point measure obtained by first restricting its support to the indicated region, and then shifting the resulting point measure to be supported on $[0,b-a]\times\Exc$.

\begin{construction}[Type 1]\label{type1:construction} For $(x,\gamma)\in\cJ^\bullet$, consider
  $$(\ff,\cev\fN,\fN_\gamma)\sim{\tt BESQ}_x(-1)\otimes\fP^0_\gamma=:\fP_{x,\gamma}^1,$$
  where $\fP^0_\gamma$ is as in Construction \ref{type0:construction}. 
  Let $\fN_*:=\clade(\ff,\cev{\fN})\concat\fN_\gamma$. We define a $\cJ^\bullet$-valued evolution $((m^y,\gamma^y),y\ge 0)$ as
  \begin{equation}\label{eq:type1:construction}
     (0,m^y)\concat\gamma^y:=\skewer(y,\fN_*),\quad y\ge 0,
  \end{equation}
  where $m^y=0$ if and only if the skewer in the last expression has no leftmost block.  
\end{construction}

By combining \ref{def:type01_meas} and \eqref{eq:clade_from_descent}, we obtain the following result.

\begin{proposition}\label{prop:type1:construction} Let $(x,\gamma)\in\cJ^\bullet$. Then the $\cJ^\bullet$-valued process resulting from Construction \ref{type1:construction} is a type-1 evolution starting from $(x,\gamma)$.
\end{proposition}
\begin{remark}\label{rem:LMB}
We will use Construction \ref{type1:construction} in the construction and analysis of type-2 evolutions in Chapter \ref{ch:type-2}. Specifically, it is instructive to explore the 
behaviour of the leftmost block in this construction.

\begin{figure}
 \centering
 \includegraphics[height=2in,width=0.45in]{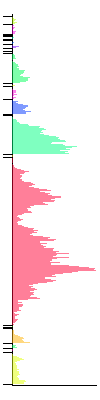}\qquad 
 \includegraphics[height=2in,width=4in]{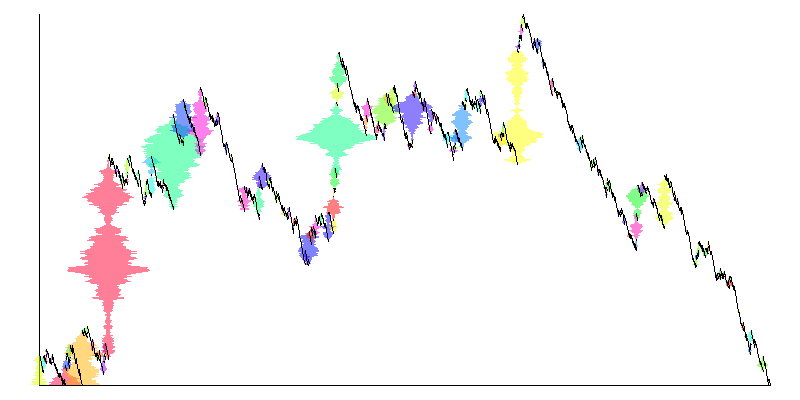}
 \caption{Left: The leftmost block process in a type-1 evolution, plotted with a vertical time-axis, with tick marks to the left of the axis indicating jump times in this process. Right: The scaffolding and spindles giving rise to this type-1 evolution.\label{fig:LMB_process}}
\end{figure}

When $x>0$, the evolution of the leftmost block is initially just $\ff\sim{\tt BESQ}_x(-1)$ independently of the evolution of the 
interval partition component, which by Construction \ref{type0:construction} is a type-0 evolution until time $\zeta(\ff)$. This time is an independent time for the type-0 evolution encoded by 
$(\cev\fN,\fN_\gamma)\sim\fP_\gamma^0$, so it has no leftmost block. 

Indeed, the scaffolding $\fX_*:=\xi(\fN_*)$ begins with (a first passage descent of) a spectrally positive 
${\tt Stable}\big(\frac32)$ L\'evy process starting from $\zeta(\ff)$ until reaching level zero. This L\'evy process has unbounded variation. By \cite[Corollary VII.5]{BertoinLevy}), it enters $(\zeta(\ff),\infty)$ 
immediately via an accumulation of small jumps. The skewer map extracts the leftmost block at these levels from the parts of these spindles that exceed the running supremum. 

As a consequence, the evolution of the leftmost block, depicted in Figure \ref{fig:LMB_process}, exhibits an accumulation of small jumps up from zero each continued by a ${\tt BESQ}(-1)$ evolution to take it
back to zero. Whenever the leftmost block jumps from $m^{y-}=0$ to $m^y>0$, there is a corresponding jump from $\gamma^{y-}=(0,m^y)\concat\gamma^y$ to $\gamma^y$. In other words, each jump of the leftmost block corresponds to the removal of an interval from the interval partition component. This delicate behaviour is efficiently encoded by scaffolding and spindles, and our understanding stems from the theory of L\'evy processes.

We care about the separate leftmost block component in a type-1 evolution, because type-2 evolutions will have two of them, and they will both interact in the same way with the 
same interval partition component. In the (discrete or) interval partition tree context, we will apply type-1 and type-2 evolutions to obtain evolutions of the decompositions around a binary branch point into one or two subtree masses and one edge partition, as in (Figure \ref{fig:discr_T5_decomp} or) Figure \ref{fig:B_k_tree_proj}.  
\end{remark}

We showed in \cite[Proposition 5.4]{Paper1-1} that $\beta\mapsto\fP_\beta^1$ is a stochastic kernel and in \cite[Proposition 6.11]{Paper1-1} that all of these distributions are measures on the space $\cN$ discussed at the end of Section \ref{sec:prelim:SSS}. It follows that $\beta\mapsto\fP_\beta^0$ and $(x,\gamma)\mapsto\fP_{x,\gamma}^1$ are likewise stochastic kernels, that $\fP_\beta^0$ is a measure on $\cev\cN\times\cN$ and $\fP_{x,\gamma}^1$ a measure on $\cE\times\cev\cN\times\cN$.
For any probability distribution $\nu$ on $\cI$ we denote by $\fP_\nu^0$ the $\nu$-mixture of the laws $\fP_\beta^0$. Similarly, for any probability 
distribution $\mu$ on $\cJ^\bullet$, we denote by $\fP_\mu^1$ the $\mu$-mixture of $\fP_{x,\gamma}^1$.




Previously in \cite{Paper1-2} we only studied $\cI$-valued type-1 evolutions, so here we collect some analogous properties of pair-valued type-1 evolutions.

\begin{proposition}\label{corpairtype1}
 Type-1 evolutions are $\cJ^\bullet$-valued Borel right Markov processes, but not Hunt processes. 
\end{proposition}

Recall Sharpe's definition \cite{Sharpe} (see also \cite[Definition A.18]{Li2011}) of Borel right Markov processes:
\begin{enumerate}[label=\arabic*., ref=\arabic*.]
  \item Lusin state space (homeomorphic to a Borel set in a compact metric space),
  \item right-continuous sample paths,
  \item Borel measurable semigroup and strong Markov property.
\end{enumerate}
Hunt processes are additionally required to be quasi-left-continuous, i.e.
\begin{enumerate}
  \item[4.] left-continuous along all increasing sequences of stopping times.
\end{enumerate}

\begin{proof}
 1. By Proposition \ref{prop:IPspace:Lusin}, the space $(\cJ^\bullet,d^\bullet)$ 
    is a Borel subset of a product of Lusin spaces and is therefore Lusin.

 2. Consider $\fN_*$ and $(0,m^y)\concat\gamma^y=\skewer(y,\fN_*)$ as in Construction \ref{type1:construction}. Then $y\mapsto m^y$ is c\`adl\`ag and  
    the only jumps are up from zero, one at the starting level of each excursion of $\xi(\fN_*)$ below the supremum. Recall from \cite[Theorem VII.4 or Lemma VIII.1]{BertoinLevy} 
    that ${\tt Stable}\big(\frac32\big)$ processes have ${\tt Stable}\big(\frac12\big)$ ladder height subordinators with zero drift coefficient. By \cite[Theorem VI.19]{BertoinLevy}, it is a.s.\ the case 
    that no two such excursions share an endpoint. It is not difficult to show that $(m^y,\gamma^y)$ is also c\`adl\`ag since for $m_n\rightarrow m_0$ and 
    $(m_n,\gamma_n)\in\cJ^\bullet$ for all $n\ge 0$, we have 
    $$d_\cI(\gamma_n,\gamma_0)\rightarrow 0\quad\mbox{if and only if}\quad d_\cI((0,m_n)\concat\gamma_n,(0,m_0)\concat\gamma_0)\rightarrow 0.$$ 
  
 3. Since $\varphi$ and $\varphi^{-1}$ are measurable bijections, the measurability of the semigroup and the strong Markov property follow 
    from Proposition \ref{prop:type01:diffusion}. 
  
 4. Consider two independent $\cI$-valued type-1 evolutions $(\beta^y,y\ge 0)$ and $(\gamma^y,\,y\ge 0)$. By Corollary \ref{type01plustype1}, the concatenation  $\beta^y\concat\gamma^y$ defines an $\cI$-valued type-1 evolution. Consider $\eta_n=\inf\{y\ge 0\colon \|\beta^y\|<\frac1n\}$. Then $\eta_n$ increases to $\eta=\inf\{y\ge 0\colon\beta^y=\emptyset\}$. Then the leftmost block at level $\eta_n$ converges to 0, but the leftmost block of $\gamma^\eta$ is non-zero with positive probability.
\end{proof}
%
%
%

\begin{figure}
 \centering
 \scalebox{.9}{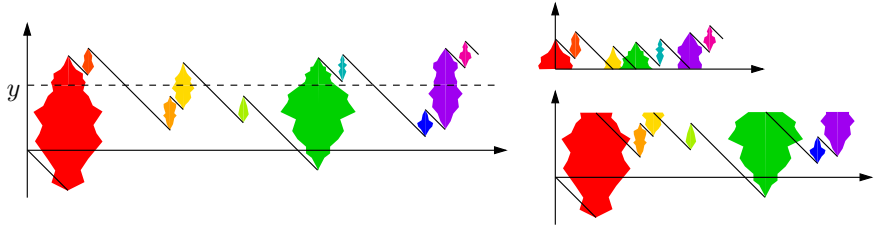}
 \caption{Illustration \cite[Figure 4]{Paper1-1} of spindles cut at level $y$. Left: $N$. Right: $N^{\ge y}$ and $N^{\le y}$.\label{fig:cutoff}}
\end{figure}

The Markov property of a ($\cJ^\bullet$-valued or $\cI$-valued) type-1 evolution corresponds to a Markov-like property of the scaffolding-and-spindles construction. 
In the setting of Proposition \ref{def:type01_meas}, this was developed in \cite{Paper1-1} and \cite[Appendix B]{Paper1-2}. 
Here, we decompose $\fN_\beta\sim\fP_\beta^1$ for each $y\ge 0$ into a point measure $\fN_\beta^{\ge y}$ of spindles (some broken) above level $y$ and a point measure 
$\fN_\beta^{\le y}$ of spindles (some broken) below level $y$, as illustrated in Figure \ref{fig:cutoff}.

More formally, consider any $N=\sum_{i\in I}\delta(t_i,f_i)\in\cN$ and any level $y\ge 0$. Let $I^y=\{i\in I\colon\xi_N(t_i-)<y<\xi_N(t_i)\}$. Then each spindle $f_i$ for $i\in I^y$ marks an upward passage of level $y$. We break $f_i$ into $f_i^{\le y}(z)=f_i(z)\cf\{z\in[0,y-\xi_N(t_i-)]\}$ and $f_i^{\ge y}(z)=f_i(y-\xi_N(t_i-)+z)\cf\{z\in[0,\infty)\}$. Since $\xi(N)$ is c\`adl\`ag with no negative jumps, each such $t_i$ is followed by a point measure $N_i^{\ge y}=\shiftrestrict{N}{(t_i,s_i^{\ge y}]}$ such that $\delta(0,f_i^{\ge y})+N_i^{\ge y}$ is a clade (or a concatenation of clades), i.e. 
$\zeta(f_i^{\ge y})+\xi(N_i^{\ge y})$ is nonnegative and ends at the next passage of $\zeta(f_i^{\ge y})+\xi(\shiftrestrict{N}{(t_i,\len(N)]})$ below 0 (or at $\len(N)-t_i$). We similarly define $N_i^{\le y}$. Then 
\begin{equation}\label{eq:upper}
  N^{\ge y}=\Concat_{i\in\mathcal{A}^y}\!\left(\delta(0,f_i^{\ge y})+N_i^{\ge y}\right)\!,\ \ \mbox{where }\mathcal{A}^y=I^y,\mbox{ and }i\preceq j\iff t_i\le t_j,
\end{equation}
and $N^{\le y}=\Concat_{i\in\mathcal{A}^y}\big(\delta(0,f_i^{\le y})+N_i^{\le y}\big)$ are point measures in $\cN$ that decompose $N$. 


\begin{lemma}[Proposition 6.6 of \cite{Paper1-1}]\label{lem:type1Markov}
 Let $y\ge 0$. In the setting of Proposition \ref{def:type01_meas} with $\fN_\beta\sim\fP_\beta^1$, conditionally given $\fN_\beta^{\le y}$, the point measure $\fN_\beta^{\ge y}$ has regular conditional distribution $\fP_{\beta^y}^1$, 
where $\beta^y=\skewer(y,\fN_\beta)$. This includes the degenerate case $\beta^y=\emptyset$.
\end{lemma}
%

Finally, we give a simpler construction of type-0 evolutions in the pseudo-stationary case that (unlike the construction based on Proposition \ref{def:type01_meas}) does not require concatenating infinitely many clades. Consider $\fN\sim\PRM[\Leb\otimes\nu_{\tt BESQ}]$ and its aggregate mass process $M_{\fN}^0$ defined as in \eqref{eq:skewer_def}. By \cite[Proposition 3.2]{Paper1-1}, this is a $\Stable(\frac12)$ subordinator (up to a time-change that does not affect its range). By Proposition \ref{prop:PDIP}(iii), the jump sizes of $M_{\fN}^0$ prior to exceeding an independent $\ExpDist[\lambda]$ threshold form a ${\tt Gamma}(\frac12,\lambda)$ multiple of a ${\tt PDIP}\big(\frac12,\frac12\big)$, which is pseudo-stationary for type-0 evolutions, by Proposition \ref{prop:01:pseudo_g}. In this context, \cite[Proposition 5.6]{Paper1-1} yields the following result.

\begin{lemma}\label{lem:type1:pseudo_constr}
 Fix $\lambda>0$.  Consider independent $\fN\sim\PRM[\Leb\otimes\nu_{\tt BESQ}]$ and $Z\sim\ExpDist[\lambda]$. Define $T := \inf\{t > 0\colon M^0_{\fN}(t) > Z\}$. Then $\skewer\big(0,\restrict{\fN}{[0,T)}\big)$ is a \PDIP[\frac12,\frac12] scaled by an independent \GammaDist[\frac12,\lambda], and, recalling the notation \eqref{eq:upper}, $\big(\restrict{\fN}{[0,T)}\big)^{\ge 0}$ is a type-1 point measure with initial state $\gamma$, in the sense of Proposition \ref{def:type01_meas}. Moreover, applying Construction \ref{type0:construction} to $\big(\cev\fN,\big(\restrict{\fN}{[0,T)}\big)^{\ge 0}\big)$ yields a pseudo-stationary type-0 evolution with \GammaDist[\frac12,\lambda] initial mass.
%
\end{lemma}


\chapter{Type-2 evolutions}\label{ch:type-2}\label{clocking}

\begin{figure}
 \centering
 \includegraphics[angle=90,scale=0.9]{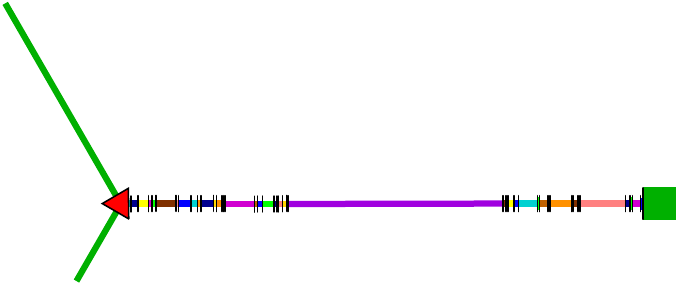}
 \caption{Tree representation of $(a,b,\gamma)\in\cJ^\circ$.}
\label{fig:2tree:int}
\end{figure}

A type-2 evolution 
is a process that takes values in a product space that combines two masses and an interval partition. More precisely, we will establish them as Markov processes in the state space 
\begin{equation}\label{type2spaces}
 \cJ^\circ :=\left\{ (a,b,\gamma)\in [0,\infty)^2 \times \cI,\;  a+b >0  \right\} \cup \{(0,0,\emptyset)\}.
\end{equation} 
Let $d^\circ$ denote the metric on $\cJ^\circ$ given by 
\[
d^\circ\left( (a_1, b_1, \gamma_1), (a_2, b_2, \gamma_2)  \right)= \abs{a_1-a_2} + \abs{b_1 - b_2} + d_{\cI}(\gamma_1, \gamma_2). 
\]
We think of $(a,b,\gamma)\in\cJ^\circ$ as a tree with a branch point separating two masses at the top of a spine down to the root of the tree, with the intervals representing an ordered collection of further masses on the spine. Equivalently, we can view $(a,b,\gamma)$ as an interval partition $\gamma$ with two additional blocks of sizes $a$ and $b$ that 
we consider both adjacent to the left end of $\gamma$. See Figure \ref{fig:2tree:int} for an illustration.

Recall that we denote by $\besq_a(-1)$ the distribution of a squared Bessel diffusion of dimension $-1$ starting from $a\ge 0$, killed upon hitting zero, and that $\zeta(\ff)$ denotes the \emph{lifetime} of the process $\ff\sim\besq_a(-1)$. Intuitively, the transition mechanism of type-2 evolutions is such that the interval partition evolves as a type-0 evolution independent of the two top masses that evolve as independent $\besq(-1)$, up until one of the top masses reaches 0. At that time, that top mass interacts with the interval partition component in the same manner in which the top mass and interval partition interact in a type-1 evolution. This is a delicate notion which we will revisit throughout the early stages of this chapter. Here is a formal definition.

\begin{definition}\label{def:type2:v1}
 Let $(a,b,\gamma)\in\cJ^\circ$. A type-2 evolution starting from $(a,b,\gamma)$ is a $\cJ^\circ$-valued process of the form 
$((m_1^y,m_2^y,\beta^y),y\ge 0)$, with $(m_1^0,m_2^0,\beta^0)=(a,b,\gamma)$. Its distribution is specified by the following iterative construction.
 
 Let $\left(\fm^{(0)}, \gamma^{(0)} \right)$ be a type-1 evolution starting with the initial condition $( b, \gamma)$ and independent of $\ff^{(0)}\sim \besq_a(-1)$, and let $Y_0=0$. Prior to time $Y_1:=\zeta(\ff^{(0)})$, we define the type-2 evolution as
 \[
 \big( m_1^y, m_2^y, \beta^y  \big):=\big(\ff^{(0)}(y), \fm^{(0)}(y), \gamma^{(0)}(y)  \big), \quad 0\le y \le Y_1,
 \]
 and proceed inductively. Suppose, for some $n\ge 1$, the process has been constructed until time $Y_n$ with $m_1^{Y_n}\!+m_2^{Y_n}\!>\!0$. Conditionally given this history, consider a type-1 evolution $(\fm^{(n)}, \gamma^{(n)})$ starting from $(0,\beta^{Y_n}) = (0,\gamma^{(n-1)}(Y_n\!-\!Y_{n-1}))$ that is independent of $\ff^{(n)}$, a $\besq(-1)$ diffusion with initial value $\fm^{(n-1)}(Y_n\!-\!Y_{n-1})$. The latter equals $m_2^{Y_n}$ if $n$ is odd or $m_1^{Y_n}$ if $n$ is even.
Set $Y_{n+1}=Y_n\!+\zeta(\ff^{(n)})$. For $y\in(0,Y_{n+1}\!-\!Y_n]$, define
 \begin{align*}
  &\left(m_1^{Y_n+y},m_2^{Y_n+y},\beta^{Y_n+y}\right):= \begin{dcases}
  (\fm^{(n)}(y),\ff^{(n)}(y),\gamma^{(n)}(y)), &\mbox{if $n$ is odd},\\
  (\ff^{(n)}(y),\fm^{(n)}(y),\gamma^{(n)}(y)), &\mbox{if $n$ is even.}
  \end{dcases}
 \end{align*}
 If, for some $n\!\ge\! 1$, $m_1^{Y_n}\!+\!m_2^{Y_n}\!=\!0$, set $(m_1^y,m_2^y,\beta^y)\!:=\!(0,0,\emptyset)$, $y\!>\!Y_n$, and $Y_{n+1} \!:=\! \infty$.
\end{definition}

We refer to the alternation between even and odd $n$ as \em regime switching\em. In the even (respectively odd) regime, only the second (respectively first) top mass can interact with the interval partition, and does so in the same delicate way as in a type-1 evolution extracting the masses from the left end of the interval partition one at a time and implicitly handling accumulations of small intervals as explained in Remark \ref{rem:LMB}. A regime change is triggered each time the other top mass vanishes. In this chapter we establish the following two theorems -- analogues of
Proposition \ref{prop:type01:diffusion} or Corollary \ref{corpairtype1}, and of Proposition \ref{type1totalmass}, respectively.

\begin{theorem}\label{thm:diffusion}
 Type-2 evolutions are Borel right Markov processes on $(\cJ^\circ\!,d^\circ)$. 
\end{theorem}

\begin{theorem}\label{thm:type2:total_mass}
 For a type-2 evolution $((m_1^y,m_2^y,\beta^y),\,y\ge0)$, the total mass process $(m_1^y+m_2^y+\|\beta^y\|,\,y\ge0)$ is a \BESQ[-1] process.
\end{theorem}

Before we can turn to proving any of the main claims including c\`adl\`ag sample paths, strong Markov property and $\besq(-1)$ total mass in Sections \ref{secbrm} and \ref{sectm}, we need to first verify that type-2 evolutions are well-defined. Specifically, we start by showing that 
\begin{enumerate}\item[(1)] the distribution of a type-2 evolution does not depend on the starting regime; 
  \item[(2)] regime change times do not accumulate to a finite limit, rather the number of regime changes is almost surely finite, thereby ensuring for instance that total mass approaches zero continuously.
\end{enumerate} To do this, in Section \ref{sec:sym}, we provide a scaffolding-and-spindles construction of type-2 evolutions. 
We also provide 
a further construction in Section \ref{sec:interweaving} that yields type-2 evolutions with special initial distributions that are relevant in Chapter \ref{dePoiss} to establish stationary unit-mass 2-tree evolutions and to study induced 3-mass processes.  

\section{Symmetry and non-accumulation of regime changes}\label{sec:sym}


Definition \ref{def:type2:v1} builds a type-2 evolution from sequences of type-1 evolutions and $\besq(-1)$ processes, ensuring after each regime change that there will again be two top masses. Specifically, one top mass is obtained by extracting the leftmost block (if the type-1 evolution has not degenerated yet) and letting it evolve as a $\besq(-1)$ evolution. The other top mass is obtained since type-1 evolutions even when starting without a leftmost block (i.e.\ from an interval partition that has an accumulation of small blocks at the left end) will again give rise to a further top mass (until they degenerate). See the discussion after Proposition \ref{prop:type01:diffusion}. More precisely, this ensures the persistence of two positive top masses at ({\rm Leb}-almost) all times up to the degeneration of a type-1 evolution, leading to a state that consists of just a single non-zero top mass and an empty spinal partition. 

\begin{figure}
 \centering
 \scalebox{.8}{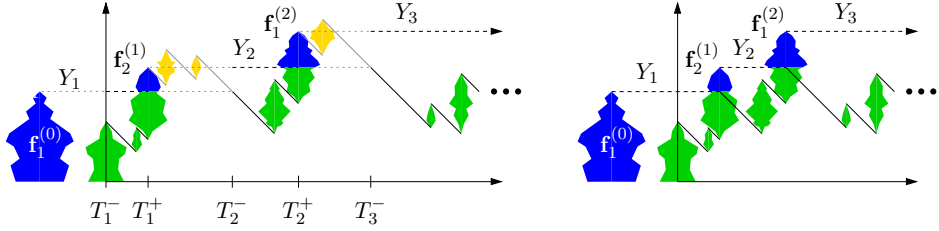}
 \caption{Construction \ref{deftype2} identifies a succession of part-spindles $(\ff^{(n)})$, colored blue, followed by spindles that are ignored, colored yellow. Other spindles are green.\label{fig:d_clk}}
\end{figure}

In fact, we can construct a type-2 evolution starting from $(a,b,\gamma)\in\cJ^\circ$ using independent $\ff_1\!\sim\!\besq_a(-1)$, $\ff_2\!\sim\!\besq_b(-1)$, and the point measures $(\cev{\fN},\fN_\gamma)\!\sim\!\fP_{\gamma}^0$ in Construction \ref{type0:construction} of an associated type-0 evolution. The type-1 evolution associated with
$(\ff_2,\cev{\fN},\fN_\gamma)\sim\besq_b(-1)\otimes\fP_\gamma^0=\fP_{b,\gamma}^1$ as in Construction \ref{type1:construction} is as required for $(\mathbf{m}^{(0)},\gamma^{(0)})$ in Definition \ref{def:type2:v1}, up to the time $Y_1=\zeta(\ff_1)$ of the first regime change. See Figure \ref{fig:d_clk} up to level $Y_1$. 

The following construction and proposition will show that, not only can $(\mathbf{m}^{(0)},$ $\gamma^{(0)})$ be derived as a function of $(\ff_2,\cev{\fN},\fN_\gamma)$, but in fact, all subsequent \BESQ\ and type-1 evolutions required in Definition \ref{def:type2:v1}, $\ff^{(n)}$, $(\fm^{(n)},\gamma^{(n)})$, $n\ge1$, can be extracted as functions of this same scaffolding and spindles.

Definition \ref{def:type2:v1} does not make use of $(\mathbf{m}^{(0)}(y),\gamma^{(0)}(y))$ for $y>Y_1$. 
By Corollary \ref{type01plustype1}, we can decompose this as the concatenation of two type-1 evolutions starting respectively from the single leftmost block, $(\mathbf{m}^{(0)}(Y_1),\emptyset)$, and the remaining partition, $(0,\gamma^{(0)}(Y_1))$. We define $(\mathbf{m}^{(1)},\gamma^{(1)})$ to be the type-1 evolution starting from $(0,\gamma^{(0)}(Y_1))$. As for the type-1 evolution starting from $(\mathbf{m}^{(0)}(Y_1),\emptyset)$, we define $\ff^{(1)}$ to be the mass evolution of the initial left-most block of this process; the remaining blocks in this process are not used in the construction. 

Applying this procedure inductively, re-framed in terms of scaffolding-and-spindles, yields the following.



\begin{construction}[Type 2]\label{deftype2}
 For $(a,b,\gamma)\in\cJ^\circ$, consider 
 $$(\ff_1,\ff_2,\cev\fN,\fN_{\gamma}) \sim \besq_a(-1)\otimes\besq_b(-1)\otimes\fP^0_{\gamma} =: \fP^2_{a,b,\gamma}.$$
 Let $\fN_* := \textsc{clade}(\ff_2,\cev\fN)\concat\fN_{\gamma}$ be the point measure of spindles of Construction \ref{type1:construction}, with $\ff_2$ taking the role of $\ff$, and let $\fX_* := \xi(\fN_*)$ be the associated scaffolding of $\fN_*$. 
 We define a $\cJ^\circ$-valued evolution 
 $((m_1^y,m_2^y,\beta^y),\,y\ge0)$ in three steps.
 
 \textbf{Step 1}. We define levels $(Y_n)$ and passage times $(T_n^{\pm})$ for $\fX_*$ inductively. See Figure \ref{fig:d_clk} for an illustration. Set $Y_0= 0$, $T_1^- = 0$, $Y_1 := \zeta(\ff_1)$, and for $n\ge 1$,
 \begin{equation}
 \begin{array}{l}\displaystyle
   T^+_n:=\inf\{t \ge T^-_n\colon \fX_*(t) > Y_n\},\qquad 
   Y_{n+1}:=\fX_*(T^+_n),\\[4pt]
   \displaystyle
   T^-_{n+1}:=\inf\{t > T^+_n\colon \fX_*(t) \leq Y_n\},
  \label{eq:d_clk:levels}
 \end{array}
 \end{equation}
 with the conventions $\inf\emptyset\!=\!\infty$ and $\fX_*(\infty)\!=\!\infty$. 
 
 \textbf{Step 2}. We define spindles $(\ff^{(n)})$ to provide one top mass on each of the intervals $[Y_n,Y_{n+1})$. Let $\ff^{(0)}=\ff_1$. For $n\ge 1$, let $\ff^{(n)}$ denote the cut-off top part $f^{\ge Y_{n}}=f(Y_n-\fX_*(T_n^+-)+\cdot\,)$ of the spindle $f$ that occurs at time $T^+_{n}$ in $\fN_*$. 
 
 \textbf{Step 3}. We finally define the evolution. For $n\ge 0$ even and $y\in [Y_n,Y_{n+1})$,
 \begin{equation}
  m_1^y := \ff^{(n)}(y-Y_{n}), \quad (0,m_{2}^y)\concat\beta^y := \skewer\big(y-Y_{n},\ShiftRestrict{\fN_*}{(T_{n+1}^-,\infty)\times\cE}\big),\label{eq:type_2_def}
 \end{equation}  
 where $m_2^y = 0$ if and only if the skewer in the last expression has no leftmost block. For $n\ge 1$ odd, the definition is the same, but with $m_1^y$ and $m_2^y$ swapped.
\end{construction}
\medskip

The effect of this construction is to skip over intervals of spindles from $\fN_*$, ensuring that they never contribute blocks to the skewer in \eqref{eq:type_2_def}. Specifically, for each $n\ge 1$, the process $\restrict{\fN_*}{(T_n^+,T_{n+1}^-]\times\cE}$ is redundant. 
This is illustrated in Figure \ref{fig:d_clk}. The (blue) part-spindles form a succession $(\ff^{(n)})$ with precisely one contribution for each $y$. The (yellow) deletions next to each $\ff^{(n)}$ are naturally interpreted as \em emigration. \em Each such family of spindles is associated with
a first passage descent of a \Stable[\frac32] scaffolding process $\xi(\fN_*)|_{(T_n^+,T_{n+1}^-]}$ as was used around \eqref{eq:JCCP_def_type0}.

\begin{proposition}
 Let $(a,b,\gamma)\in\cJ^\circ$. Then the $\cJ^\circ$-valued process resulting from Construction \ref{deftype2} is a type-2 evolution starting from $(a,b,\gamma)$.
\end{proposition}
\begin{proof} 
%
  Consider $(\ff_1,\ff_2,\cev{\fN},\fN_\gamma)\sim\fP_{a,b,\gamma}^2$ and $\fN_*:=\textsc{clade}(\ff_2,\cev\fN)\concat\fN_{\gamma}$, as well as the filtration 
  $(\cF^y,y\ge 0)$ generated by $(\ff_1(y),\fN_*^{\le y})$, using notation introduced above Lemma \ref{lem:type1Markov} for point measures that have 
  been cut off at scaffolding level $y$. We use the notation of Construction \ref{deftype2} to inductively
  set up all random variables as needed for Definition \ref{def:type2:v1}, and we will show that Definition
  \ref{def:type2:v1} and Construction \ref{deftype2}, in this setup, yield pathwise the same process $((m_1^y,m_2^y,\beta^y),y\ge 0)$. For the purpose 
  of this proof we will mark all random variables appearing in Definition \ref{def:type2:v1} by an underscore. We also define 
  $\fN_*^{(0)}:= \fN_*$ and 
 \begin{equation}\label{eq:d_clk:shift}
   \fN_*^{(n)}:=\left(\ShiftRestrict{\fN_*}{(T_{n+1}^-,\infty)\times\cE}\right)^{\ge 0},
 \end{equation}
 for $n\ge1$, where we recall our notation for shifted restrictions introduced for \eqref{eq:clade_from_descent}, and where the superscript $^{\ge 0}$ on the right-hand side is in the sense of the cutoff processes $N^{\ge y}$ defined in \eqref{eq:upper}, in which spindles below a given level are removed (or cut off if they straddle the level). Each $\fN_*^{(n)}$ is a type-1 point measure for one top mass and spinal masses during the interval $[Y_{n},Y_{n+1})$, which we understand as a time interval for $((m_1^y,m_2^y,\beta^y),y\ge 0)$ and a level interval for $\fN_*$ and its associated scaffolding $\fX_*=\xi(\fN_*)$.

  Now, $\underline{\ff}^{(0)}:=\ff_1$ and 
  $(0,\underline{\fm}^{(0)})\star\underline{\gamma}^{(0)}:=\skewerbar(\fN_*^{(0)})$ 
  have the appropriate joint distribution and achieve 
  $((\underline{m}_1^y,\underline{m}_2^y,\underline{\beta}^y),\,0\!\le\! y\!\le\! \underline{Y}_1)=((m_1^y,m_2^y,\beta^y),\,0\!\le\! y\!\le\! Y_1)$. \linebreak
  Suppose we have completed the construction up to $(\underline{\ff}^{(n-1)},\underline{\fm}^{(n-1)},$ 
  $\underline{\gamma}^{(n-1)})$ and identified 
    $((\underline{m}_1^y,\underline{m}_2^y,\underline{\beta}^y),\,0\le y\le \underline{Y}_{n})=((m_1^y,m_2^y,\beta^y),\,0\le y\le Y_{n})$ 
  for some $n\ge 1$. Then given $\cF^{Y_{n-1}}$, we apply Lemma \ref{lem:type1Markov}, which is the Markov-like property of the type-1 point measure 
  $\fN_*^{(n-1)}$ at the (conditionally independent) level $\zeta(\ff^{(n-1)})=Y_n-Y_{n-1}$, to find an above-$\zeta(\ff^{(n-1)})$ point measure $(\fN_*^{(n-1)})^{\ge\zeta(\ff^{(n-1)})}$. On the event that this is non-trivial, the first clade of this point measure has an initial spindle $\ff^{(n)}$ and the remaining clades form a type-1 point measure $\fN_*^{(n)}$. Noting that $\ff^{(n)}$ and $\fN_*^{(n)}$ are 
  conditionally independent given the below-$\zeta(\ff^{(n-1)})$ information of $(\ff^{(n-1)},\fN_*^{(n)})$, indeed given $\cF^{Y_n}$, we proceed as follows. Suppose $n$ is even. First, 
    $\underline{\ff}^{(n)}:=\ff^{(n)}\sim{\tt BESQ}_{m_1^{Y_{n}}}(-1)$, 
  is as needed for Definition \ref{def:type2:v1}, since $\underline{Y}_{n}=Y_{n}$. 
  Second, $\fN_*^{(n)}\sim\fP_{\beta^{Y_{n}}}^1=\fP_{\underline{\beta}^{\underline{Y}_n}}^1$ gives rise to a type-1 evolution 
  $(0,\underline{\fm}^{(n)})\star\underline{\gamma}^{(n)}:=\skewerbar(\fN_*^{(n)})$ started from $(0,\underline{\beta}^{\underline{Y}_n})$, 
  as required, as $\underline{m}_{2}^{\underline{Y}_n}=\underline{\ff}^{(n-1)}(\zeta(\underline{\ff}^{(n-1)}))=0$. This also implies that 
  for all $y\in[0,Y_{n+1}-Y_n)$
  \begin{align*}
   &(\underline{m}_{1}^{\underline{Y}_n+y},(0,\underline{m}_{2}^{\underline{Y}_n+y})\star\underline{\beta}^{\underline{Y}_n+y})
   		=(\underline{\ff}^{(n)}(y),(0,\underline{\fm}^{(n)}(y))\star\underline{\gamma}^{(n)}(y))\\
   &\qquad\qquad =(\ff^{(n)}(y),\skewer(y,\fN_*^{(n)}))
   =(m_{1}^{Y_n+y},(0,m_{2}^{Y_n+y})\star\beta^{Y_n+y}).
  \end{align*}
  The same argument applies for $n$ odd, with the roles of 1 and 2 interchanged.
\end{proof}


\begin{lemma}[Symmetry]\label{lem:type2_symm}
  If we modify Construction \ref{deftype2} by letting $Y_1=\zeta(\ff_2)$ and $\fN_*=\textsc{clade}(\ff_1,\cev\fN)\concat\fN_{\gamma}$ and accordingly swapping the parity in Step 3., we obtain a type-2 evolution that is pathwise the same as in Construction \ref{deftype2}, with identical sets $\{\zeta(\ff_1),\zeta(\ff_2)\}\cup\{Y_n,n\ge 0\}$. In particular, the point measure $\restrict{\cev{\fN}}{(-\infty,T_{\min\{\zeta(\ff_1),\zeta(\ff_2)\})}(\cev{\fN}))\times\cE}$ is redundant in this construction of the type-2 evolution.
\end{lemma}

\begin{proof} 
 For the purposes of this proof, we add underscores and write $\underline{Y}_{\,j}$, $\underline{\ff}^{(j)}$, $j\ge 0$, and  
  $((\underline{m}_{\,1}^y,\underline{m}_{\,2}^y,\underline{\beta}^y),y\ge 0)$ in the modification of Construction \ref{deftype2}. We remark that the underscores here are unrelated to those in the previous proof. 
  The main aim of this proof is to show the pathwise equality 
  $((\underline{m}_{\,1}^y,\underline{m}_{\,2}^y,\underline{\beta}^y),y\ge 0)=((m_1^y,m_2^y,\beta^y),y\ge 0)$. We only discuss the case where $a>0$ and $b>0$. The cases where $a=0$ or $b=0$ can
  then be checked similarly.
  
  On the event $\{\zeta(\ff_1)<\zeta(\ff_2)\}$, we have $\underline{Y}_{\,0}=0=Y_0<Y_1=\zeta(\ff_1)<Y_2=\zeta(\ff_2)=\underline{Y}_{\,1}$, and we see inductively that
  $\underline{Y}_{\,j}=Y_{j+1}$, $\underline{\ff}^{(j)}=\ff^{(j+1)}$ and $\ShiftRestrict{\underline{\fN}_*}{(\underline{T}_{j}^-,\infty)\times\cE}=\ShiftRestrict{\fN_*}{(T_{j+1}^-,\infty)\times\cE}$ for all $j\ge 1$. It is now easy to see that the pathwise equality holds on this event. Similarly, on $\{\zeta(\ff_1)>\zeta(\ff_2)\}$, we have $\underline{Y}_{\,1}=\zeta(\ff_2)$ and 
  $\underline{Y}_{\,j+1}=Y_{j}$ for all $j\ge 1$, and the same argument applies.

  In particular, the sets $\{Y_n,n\ge 0\}$ and $\{\underline{Y}_{\,n},n\ge 0\}$ differ precisely by the omission of either $\zeta(\ff_2)$ from the former or of $\zeta(\ff_1)$ 
  from the latter. The last statement of the lemma follows using the original definition on $\{\zeta(\ff_1)<\zeta(\ff_2)\}$ and the modified definition on 
  $\{\zeta(\ff_1)>\zeta(\ff_2)\}$. 
\end{proof}   

It is not a priori clear in Definition \ref{def:type2:v1}, nor from Construction \ref{deftype2}, that regime changes cannot accumulate at a 
finite $Y_\infty=\sup_{n\ge 0}Y_n<\infty$. This would leave the type-2 evolution undefined for times $y\ge Y_\infty$, so we address this point before turning to any further properties. 

\begin{lemma}\label{type2welldef}
 For all $(a,b,\gamma)\in\cJ^\circ$, the type-2 evolution of Definition \ref{def:type2:v1} or equivalently of Construction \ref{deftype2} is such that there is a.s.\ some finite $n\ge 0$ for which $Y_n<Y_{n+1}=\infty$ and as $y$ increases to $Y_n$, the evolution $(m_1^y,m_2^y,\beta^y)$ approaches $(0,0,\emptyset)$ continuously. 
\end{lemma}
\begin{proof}
%
 %
 First, we prove the claimed convergence to $(0,0,\emptyset)$. 
 The three events $\{T_1^+ = \infty\}$, $\{Y_2 = \infty\}$, and $\{\zeta(\skewerbar(\fN_*)) < Y_1\}$ are equal up to null sets. 
On these events, $m_1^y$ converges to 0 as $y$ increases to $Y_1$, and $(m_2^y,\beta^y)$ is already absorbed at $(0,\emptyset)$ beforehand. We proceed inductively. On $\{T_{n-1}^+ <\infty\}$, this time $T_{n-1}^+$ is a time when the type-1 scaffolding $\fX_*$ exceeds level $Y_{n-1}$. Since this scaffolding eventually reaches level 0, we get $T_{n}^- <\infty$ a.s. Now, on the event $\{T_{n}^+ = \infty\}$, we apply the same argument as before to $\shiftrestrict{\fN_*}{(T_{n}^-,\infty)\times\cE}$ in place of $\fN_*$, to conclude that $(m_1^y,m_2^y,\beta^y)$ approaches $(0,0,\emptyset)$ as $y$ increases to $Y_{n}$.
 
 It remains to show that $Y_n<Y_{n+1}=\infty$ for some $n\ge 1$. 
  We claim that it suffices to prove the following.
  \begin{center}
   $(*)$\ \ \ \parbox[c]{\textwidth - 2cm}{
  Consider any two spindles of heights $\zeta(f_1)=c_1$ and $\zeta(f_2)=c_2$ with $c_1<c_2$. Apply Construction \ref{deftype2} to $(\ff_1,\ff_2,\cev{\fN},\fN_\emptyset)=(f_1,f_2,\cev{\fN},0)$ for $(\cev{\fN},0)\sim\fP_\emptyset^0$. Then there is some $n\ge 1$ for which $Y_n < Y_{n+1}=\infty$.}\ \ \ \hphantom{(*)}
  \end{center}
 Indeed, once this is shown, $a+b>0$ in the general case implies $Y_2>0$, and only finitely many clades of $\fN_\gamma$ survive to level $Y_2$ \cite[Lemma 6.1]{Paper1-1}. 
 We apply $(*)$ to these clades one by one, with $c_1$ as the final regime switch level of the preceding clades and $c_2$ as the next level beyond $c_1$ at which the top mass of the next clade vanishes, to see that each clade contributes a finite number of regime changes.
 
 To prove $(*)$, we note that this can be read as a statement about the 
  \Stable[\frac32]\ L\'evy process $X=c_2+\xi\big(\shiftrestrict{\cev{\fN}}{(T_{c_2}(\cev{\fN}),0)\times\cE}\big)$, cf.\ Lemma \ref{lm:scaff} and the discussion leading up to
  \eqref{eq:clade_from_descent}.  
  Specifically, note that $\zeta\big(\ff^{(j)}\big)$ is the overshoot $Y_{j+1}-Y_j$ of $X$ when first crossing level $Y_j$ after the stopping time $T_j^-$.

  Now we extend $X$ to a \Stable[\frac32]\ process with infinite lifetime so that $T_j^+<\infty$ for all $j\ge 1$, and we show that $Y_j\rightarrow\infty$. To this end, let $\Delta_n=Y_{n+1}-Y_n$ and $R_n=\Delta_{n+1}/\Delta_n$
  for $n\ge 1$. By the strong Markov property of \Stable[\frac32], the conditional distribution of $\Delta_{n+1}$ given $\Delta_1,\ldots,\Delta_n$ equals the law of the overshoot of a $\Stable[\frac32]$ process
  when first crossing $\Delta_n$, which is the same as the overshoot of its $\Stable[\frac12]$ ladder height subordinator, see e.g.  \cite[Theorem VII.4 or Lemma VIII.1]{BertoinLevy}. By stable scaling, for each $n$, $R_n$ is independent of $\Delta_n$ and  
  distributed like the overshoot of a $\Stable[\frac12]$ (ladder height) subordinator across 1. So the sequence $(R_n,\,n\geq 1)$ is i.i.d.\ and\vspace{-0.2cm}
  $$\Delta_{n+1} = \Delta_1 \cdot \prod_{i=1}^{n} R_{i} \qquad \text{for }n\geq 1.$$
  Thus, $(\log(\Delta_n),\,n\geq 1)$ is a random walk. It suffices to show that the increments $\log(R_n)$, $n\ge 1$, of this walk have non-negative expected value.
 
  We can get at the law of $R_{n}$ by taking advantage of the $\Stable[\frac12]$ inverse local time subordinator associated with one-dimensional Brownian motion, $(B(t),\,t\geq 0)$, see e.g. \cite[Proposition III.(3.8) and Corollary VI.(2.3)]{RevuzYor}. In this 
  setting, $R_{n}$ is distributed like $T-1$, where $T$ is the time of the first return of $B$ to zero, after time 1. By a calculation based on the reflection principle, we find
  $\bP(T < t) = \frac{2}{\pi}\arctan(t-1)$.
 Thus,
  $\bE\left(  \log(R_n) \right) = \int_1^\infty \log(t-1)\bP(T\in dt) = \frac{2}{\pi}\int_1^{\infty} \log(t-1)\frac{1}{(t-1)^2 + 1}dt = 0.
  $ 
\end{proof}
 
We record the following consequence of the proof of Lemma \ref{type2welldef}.

\begin{corollary}\label{type2:degenlife} Let $((m_1^y,m_2^y,\beta^y),\,y\ge 0)$ be a type-2 evolution starting from $(a,b,\gamma)\!\in\!\cJ^\circ\!\setminus\!\{(0,0,\emptyset)\}$
  with regime changes $(Y_n)$. Then there is a.s.\ $n\!\ge\! 1$ such that the \em lifetime \em $\zeta\!=\!\inf\{y\!\ge\! 0\colon(m_1^y,m_2^y,\beta^y)\!=\!(0,0,\emptyset)\}$ and 
  the \em degeneration time \em $D\!=\!\inf\{y\!\ge\! 0\colon(m_1^y,\beta^y)\!=\!(0,\emptyset)\mbox{ or }(m_2^y,\beta^y)\!=\!(0,\emptyset)\}$ satisfy 
  $0\!\le\! Y_{n-1}\!\le\! D\!<\!\zeta\!=\!Y_n\!<\!\infty$, and furthermore $Y_{n-1}\!<\!D$ unless the initial state is already degenerate with $D\!=\!0$.
\end{corollary}

\section{Type-2 evolutions as Borel right Markov processes}\label{secbrm}

In this section we will prove Theorem \ref{thm:diffusion}, i.e.\ that type-2 evolutions are Borel right Markov processes, and we will further show that their semigroup is continuous. 
We listed the properties 1.--3.\ that this comprises after Corollary \ref{corpairtype1}.

\begin{proof}[Proof of Theorem \ref{thm:diffusion}]
 1. By Lemma \ref{type2welldef}, type-2 evolutions take values in the space $(\cJ^\circ,d^\circ)$ of (\ref{type2spaces}), 
    which is Lusin as Borel subset of a product of Lusin spaces (see Proposition \ref{prop:IPspace:Lusin}). 
 
 2. To confirm right-continuous paths, note that in the notation of Definition \ref{def:type2:v1}, for any $n\ge 0$, on the interval $[Y_n,Y_{n+1})$, the type-2 evolution inherits 
    c\`adl\`ag paths from the $\cJ^\bullet$-valued type-1 evolution of Corollary \ref{corpairtype1}, and from a path-continuous \BESQ[-1] process. Furthermore, we argue that 
    type-2 evolutions are continuous at $Y_n$, as follows. For even $n$, the top mass $m_1^y$ approaches $0$ continuously as $y\uparrow Y_n$, while $m_1^{Y_n}=0$ holds 
    since the transition kernels of type-1 evolutions in \eqref{eq:intro:transn_1} are concentrated on interval partitions without a second leftmost 
    block, for each fixed $y$, but hence also when mixed over the distribution of the conditionally independent random lifetime $Y_n-Y_{n-1}$ of the \BESQ[-1] process $\ff^{(n-1)}$. The argument for odd $n\ge 1$ is the
    same, with $m_2^y$ and $m_2^{Y_n}$ in the place of $m_1^y$ and $m_1^{Y_n}$.
 
 3. The method of construction undertaken in Definition \ref{def:type2:v1}, in which a right Markov process with finite lifetime is reborn at the end of the 
    lifetime according to a probability kernel, has been studied by Meyer \cite{Mey75}. 
    Type-1 evolutions and \BESQ[-1] processes are Borel right Markov processes (see Corollary \ref{corpairtype1}), and thus so too is the process $\big((m_1^y,m_2^y,\beta^y,1),\,0\le y\le Y_1\big)$ for any initial $(a,b,\gamma)\in\cJ^\circ$ 
with $a>0$.  By swapping the parity as in the statement of Lemma \ref{lem:type2_symm}, we can similarly define $\big((m_1^y,m_2^y,\beta^y,2),\,0\!\le\! y\!\le\! Y_1\big)$ for initial 
$(a,b,\gamma)\in\cJ^\circ$ 
with $b>0$, where the added fourth component $I(y)=1$ or $I(y)=2$ records which of the two top masses is set up to evolve according to \BESQ[-1] and which is forming a type-1 evolution with $\beta^y$. We define the deterministic kernel $\widebar{N}((0,x,\gamma,1);\,\cdot\,)=\delta_{(0,x,\gamma,2)}$, $\widebar{N}((x,0,\gamma,2);\,\cdot\,)=\delta_{(x,0,\gamma,1)}$. As noted in \cite[Definition 8.1]{Sharpe}, Borel right Markov processes are right Markov processes satisfying the hypoth\`eses droites, in Meyer's sense. Therefore, we can apply \cite[Th\'eor\`eme 1 and Remarque on p.474]{Mey75} to conclude that if we alternate killed processes with $I(y)=1$ and $I(y)=2$, using transitions according to $\widebar{N}$ to determine initial states from the previous killing state,
\begin{equation}\label{augment}\mbox{the process $\big((m_1^y,m_2^y,\beta^y,I(y)),\,y\ge 0\big)$ is a right Markov process,}
\end{equation}
satisfying the strong Markov property. It is not hard to show that the semigroup of this process is Borel, see e.g. the last point in the proof of \cite[Th\'eor\`eme (3.18)]{Bec07}. In Proposition \ref{contini} we strengthen this to continuity in the initial state. 
 
    Lemma \ref{lem:type2_symm} verifies Dynkin's criterion (see Appendix \ref{sec:DynkinIntertwining}) to show that the type-2 evolution is a right Markov process as well.
 %
\end{proof}


In order to establish continuity of the semigroup of type-2 evolution in the initial condition we require some intermediate results.

\begin{lemma}\label{lemma slightly stronger continuity}
Suppose that $((a_n,\beta_n),\,n\geq 1)$ is a sequence in $(\cJ^\bullet,d^\bullet)$ that converges to $(a,\beta)$ and that $(x_n,\,n\ge 1)$ is a 
sequence of times converging to $x>0$. Let $((m^y_n, \gamma^y_n), y\geq 0)$ and  $((m^y, \gamma^y), y\geq 0)$ be type-1 evolutions started from $(a_n,\beta_n)$ and $(a,\beta)$ respectively.  If $f\colon\cJ^\bullet\to \bR$ is bounded and continuous, then 
$$\bE\left[ f(m^{x_n}_n, \gamma^{x_n}_n)\right] \to  \bE\left[f(m^x, \gamma^x)\right].$$
\end{lemma}

\begin{proof}
If $g\colon\cI \to \bR$ is bounded and continuous, then the fact that
$$\bE\big[g\left((0,m^{x_n}_n)\concat \gamma^{x_n}_n\right)\big] \to  \bE\big[g\left((0,m^{x})\concat \gamma^x\right)\big].$$
is established in the proof of \cite[Proposition 6.15]{Paper1-1}.  The slightly stronger version that separates out convergence of the top mass 
follows from the coupling used in that proof. Specifically, that proof uses Proposition \ref{def:type01_meas} and reduces the argument to finitely many clades, each of which is composed of 
an initial spindle and an independent $\Stable[\frac32]$ L\'evy process. Furthermore, as noted in the proof of Lemma \ref{type2welldef}, the ladder height process of a $\Stable[\frac32]$ L\'evy 
process, in which the leftmost spindle at each level can be found, is a $\Stable[\frac12]$ subordinator. The probability that $x$ is in its range 
is zero, so that the evolution of the leftmost mass is continuous around scaffolding level $x$, i.e.\ time $x$ of the type-1 evolution, with probability one.    
\end{proof}

It will be convenient to augment the type-2 evolution $\Gamma^y:=(m_1^y,m_2^y,\beta^y)$, $y\ge 0$, by the counting process $J(y)=\inf\{j\ge 0\colon Y_{j+1}>y\}$ counting 
its regime changes. This process $((\Gamma^y,J(y)),y\ge0)$ can be constructed as a strong Markov process in the same way as in (\ref{augment}) and similarly relates to $(\Gamma^y,y\ge 0)$ by 
Dynkin's criterion. 
The state space for the evolution 
$((\Gamma^y,J(y)), y\geq 0)$ is the set
$$ \cJ^+ = \{ ((m_1,m_2,\beta), j) \in \cJ^\circ \times \mathbb{N}_0\colon \ m_{p(j+1)}>0\mbox{ or }(m_1,m_2,\beta)=(0,0,\emptyset)\}.$$
In the following lemma, we write $\bE_{\Gamma,j}$ to denote the expectation for the augmented process starting from $(\Gamma,j)\in\cJ^+$. We often write $(m_1,m_2,\beta,j)\in\cJ^+$ to mean $((m_1,m_2,\beta),j)\in\cJ^+$.

\begin{lemma}\label{lemma markovprops}
 Suppose that $(\Gamma^y,\,y\ge 0)$ is a type-2 evolution with regime changes at $(Y_n)$ and right-continuous natural filtration $(\cF^y,\,y\ge 0)$. Let $j\ge 0$. 
 Then
 \begin{enumerate}[label=(\roman*), ref=(\roman*)]
  \item for  all $f\colon\cJ^+ \to \bR$ bounded and continuous and $z\ge 0$ 
   $$\bE \left[ f(\Gamma^{Y_j+z}, J(Y_j+z)) \middle| \cF^{{Y_j}} \right] = \bE_{\Gamma^{Y_j}, j}\left[f(\Gamma^{z}, J(z))\right],\qquad\bP\mbox{-a.s.,}$$
  \item for all $h\colon\cJ^\circ \to \bR$ bounded and continuous, $y\ge 0$, and for $\bP$-a.e.\ $\omega$,
 \end{enumerate}
   \begin{align*}  &\bE\left[ h(\Gamma^y) \mathbf{1}\{Y_j\leq y <Y_{j+1}\} \middle| \cF^{Y_j}\right](\omega)\\ 
				&=\mathbf{1}\{Y_{j}(\omega)\leq y\} 
				  \bE_{\Gamma^{Y_j(\omega)}(\omega),j}\! \left[ h(\Gamma^{y\vee Y_j(\omega)- Y_j(\omega)}) \mathbf{1}\{ y\vee Y_j(\omega)-Y_j(\omega) <Y_{1}\}\right]\!.
   \end{align*}
\end{lemma}

\begin{proof}
The first claim is immediate from the definition of type-2 evolutions and the second follows from the proof of \cite[Theorem 2.3.3]{chung2006markov} applied to the augmented Markov process $((\Gamma^y,J(y)),y\ge 0)$. The book \cite{chung2006markov} (and indeed the earlier \cite{chung1982} on which the relevant parts of \cite{chung2006markov} are based) assumes that the Markov process takes values in a locally compact state space, but that is not needed in the proof of Theorem 2.3.3.  The right-continuous dependence of the semigroup on time needed in the proof follows from the right-continuity of sample paths.
\end{proof}

Next we establish weak continuity at regime changes.

\begin{lemma}\label{lemma clockcont}
Suppose that $(a_n,b_n,\gamma_n) \rightarrow (a,b,\gamma)$ in $(\cJ^\circ,d^\circ)$ with $a>0$.  Let $(\Gamma^y_n,\,y\geq 0)$ and $(\Gamma^y,\,y\geq 0)$
be type-2 evolutions started from $(a_n,b_n,\gamma_n)$ and $(a,b,\gamma)$ respectively with respective regime changes $(Y^n_k)$ and $(Y_k)$.  Then for all $j\ge 1$ and all bounded continuous functions $f\colon\cJ^\circ\times[0,\infty)\rightarrow\bR$,
 $$\bE\left[f\left(\Gamma_n^{Y^n_j}, Y^n_j\right)\right] \rightarrow \bE\left[f\left(\Gamma^{Y_j}, Y_j\right)\right].$$
\end{lemma}

\begin{proof} We first establish the claim for $j=1$. Let $((\Gamma^y_n, J_n(y)), y\geq 0)$ and $((\Gamma^y,J(y)),y\geq 0)$ be the augmented type-2 evolutions started from $(a_n,b_n,\gamma_n,0)$ and $(a,b,\gamma,0)$. Let $\ff^{(0)}$ be a $\besq(-1)$ started from $a$, let $(\mathbf{m}^{(0)}, \gamma^{(0)})$ be an independent type-1 evolution started from $(b, \gamma)$, and let $(\mathbf{m}^{(0)}_n,\gamma^{(0)}_n)$ be a type-1 evolution, independent from $\ff^{(0)}$, and started from $(b_n, \gamma_n)$.  From the definition of type-2 evolutions and the scaling property of $\besq(-1)$, we see that 
$$ (\Gamma^y,\,0\leq y\leq Y_1) \stackrel{d}{=}  \left( \left(\ff^{(0)}(y), \mathbf{m}^{(0)}(y), \gamma^{(0)}(y)\right),\,0\leq y\leq Y_1\right)$$
and for $n$ sufficiently large, $a_n>0$ and 
\begin{equation} (\Gamma^y_n,\,0\leq y\leq Y^n_1) \stackrel{d}{=} 
                  \Big( \Big(\frac{a_n}{a}\ff^{(0)}\Big(\frac{a}{a_n}y\Big), \mathbf{m}_n^{(0)}(y), \gamma^{(0)}_n(y)\Big),\,0\leq y\leq Y^n_1\Big).
  \label{scalebesq}
\end{equation}
Note that, from this construction, $Y^n_1 = (a_n/a)Y_1$.  Furthermore, from 
Lemma \ref{lem:BESQ:length} we see that $Y_1$ is distributed like $a/(2G)$ where 
$G\sim\GammaDist[\frac32,1]$.  In particular, $Y_1$ has a continuous density $q$ on $(0,\infty)$. Disintegrating
 based on the value of $Y^n_1$, i.e.\ conditioning on $Y^n_1$, we see that
\begin{equation*}
\begin{split}
&\bE\big[ f( \ff^{(0)}_{n}(Y^n_1), \mathbf{m}^{(0)}_n(Y^n_1),\gamma^{(0)}_n(Y^n_1), Y^n_1)\big]\\
 &\quad = \int_0^\infty\! \bE\!\left[ f\left( 0, \mathbf{m}^{(0)}_n\left(x\right), \gamma^{(0)}_n\left(x\right),x\right)\right]\!\frac{a}{a_n}q\!\left(\frac{ax}{a_n}\right)\! dx.
\end{split}
\end{equation*}
It follows from Lemma \ref{lemma slightly stronger continuity} and a version of the dominated convergence theorem (e.g.\ \cite[Theorem 1.21]{Kallenberg}) that
\begin{equation}\label{base1}\bE_{(a_n,b_n,\gamma_n),0}\left[f\left(\Gamma^{Y_1}, Y_1\right)\right] \rightarrow \bE_{(a,b,\gamma),0}\left[f\left(\Gamma^{Y_1}, Y_1\right)\right].
\end{equation}
This completes the proof for $j=1$, for all $a>0$, $b\ge 0$ and $\gamma\in\cI$. The same proof applied to augmented type-2 evolutions started from $(a_n,b_n,\gamma_n,1)$ and $(a,b,\gamma,1)$ shows 
\begin{equation}\label{base2}\bE_{(a_n,b_n,\gamma_n),1}\left[f\left(\Gamma^{Y_1},Y_1\right)\right] \rightarrow \bE_{(a,b,\gamma),1}\left[f\left(\Gamma^{Y_1}, Y_1\right)\right],
\end{equation}
for all $a\ge 0$, $b>0$ and $\gamma\in\cI$. The inductive step $j\rightarrow j+1$ follows from the strong Markov property of the augmented type-2 
evolutions at regime changes $Y^n_j$ and $Y_j$, applying (\ref{base2}) for odd $j$ and (\ref{base1}) for even $j$.
\end{proof}

In the following result, which is the continuity of type-2 evolutions in the initial condition, we write $\mathcal{P}(\cJ^\circ)$ for the space of Borel probability measures on $(\cJ^\circ,d^\circ)$, equipped with the topology of weak convergence.

\begin{proposition}\label{contini}
Fix  $y\!\geq\! 0$ and define $F_y\colon\cJ^\circ \to \mathcal{P}(\cJ^\circ)$, by letting $F_y(a,b,\gamma)$ be the law at time $y$ of a type-2 evolution starting from 
the initial state $(a,b,\gamma)\in \cJ^\circ$. Then $(a,b,\gamma)\mapsto F_y(a,b,\gamma)$ is weakly continuous. 
\end{proposition}

\begin{proof}
Suppose that $(a_n,b_n,\gamma_n)\rightarrow(a,b,\gamma)$ in $(\cJ^\circ,d^\circ)$, i.e. $a_n\to a$, 
$b_n\to b$ and $d_\cI(\gamma_n,\gamma)\to 0$. We may assume without loss of generality that $a>0$. Once the proof is complete for this subcase, we 
can apply Lemma \ref{lem:type2_symm} to deduce the subcase $a=0$, $b>0$; the subcase $a=b=0$, $\gamma=\emptyset$ is trivial. Let 
$(\Gamma^y_n,y\ge 0)$ and $(\Gamma^y,y\ge 0)$ be 
type-2 evolutions started from $(a_n,b_n,\gamma_n)$ and $(a,b,\gamma)$, 
respectively, with respective regime changes $(Y^n_j)_{j\ge0}$ and $(Y_j)_{j\ge0}$.  Observe that for all bounded continuous $f\colon\cJ^\circ\rightarrow\bR$
\begin{equation}\label{serieseq}
  \bE\left[ f\left(\Gamma^y_n\right)\right] = \sum_{j=0}^\infty \bE \left[ f\left(\Gamma^y_n\right) \mathbf{1}\{Y^n_j\leq y <Y^n_{j+1}\}\right].
\end{equation}
By Lemma \ref{lemma clockcont} and the Skorohod representation theorem, we may now assume $Y_j^n(\omega)\rightarrow Y_j(\omega)$, 
$d^\circ\big(\Gamma^{Y_j^n(\omega)}_n(\omega),\Gamma^{Y_j(\omega)}(\omega)\big)\rightarrow 0$, and since $\bP(Y_j=y)=0$, also 
$\mathbf{1}\{Y_j^n(\omega)\le y\}\rightarrow\mathbf{1}\{Y_j(\omega)\le y\}$ for $\bP$-a.e.\ $\omega$. 
Recall that $Y_1$ under $\bP_{\Gamma,j}$ is the lifetime of the top mass labeled 1 when $j$ is even and labeled $2$ when $j$ is odd. For $\Gamma = \Gamma^{Y_j(\omega)}(\omega)$ or $\Gamma=\Gamma_n^{Y^n_j(\omega)}(\omega)$, in either case, this is the non-zero top mass of $\Gamma$. 
Recall also from (\ref{scalebesq}) that $\besq(-1)$ 
processes with converging initial states can be coupled to converge uniformly together with their lifetimes. In particular, we 
can use their convergence in distribution together with Lemma \ref{lemma slightly stronger continuity} for the convergence of the other top mass 
and interval partitions at times $y\vee Y_j^n(\omega)-Y_j^n(\omega)\rightarrow y\vee Y_j(\omega)-Y_j(\omega)$ to obtain for $\bP$-a.e.\ $\omega$
\begin{align*}
  &\bE_{\Gamma^{Y^n_j(\omega)}_n(\omega),j} \left[ f\left(\Gamma^{y\vee Y^n_j(\omega)- Y^n_j(\omega)}\right) 
												 \mathbf{1}\left\{ y\vee Y^n_j(\omega)-Y^n_j(\omega) <Y_{1}\right\}\right]\\
  &\rightarrow \bE_{\Gamma^{Y_j(\omega)}(\omega),j} \left[ f\left(\Gamma^{y\vee Y_j(\omega)- Y_j(\omega)}\right) 
											     \mathbf{1}\left\{ y\vee Y_j(\omega)-Y_j(\omega) <Y_{1}\right\}\right].
\end{align*} 
By Lemma \ref{lemma markovprops}(ii) and applying the previous convergences and dominated convergence, we find 
\begin{equation}\label{conveq}
  \bE \left[ f(\Gamma^y_n) \mathbf{1}\{Y^n_j\leq y <Y^n_{j+1}\}\right]\rightarrow \bE \left[ f(\Gamma^y) \mathbf{1}\{Y_j\leq y <Y_{j+1}\}\right].
\end{equation}
A further application of the dominated convergence theorem yields $\bE[f(\Gamma^y_n)]$ 
$\rightarrow \bE[f(\Gamma^y)]$, completing the proof. 
\end{proof}


\section{The total mass process}\label{sectm}

In this section, we prove Theorem \ref{thm:type2:total_mass}, that the total mass process of any type-2 evolution is a $\besq(-1)$. We use the notation of Definition \ref{def:type2:v1} and work with the $\besq(-1)$ processes $\ff^{(n)}$, and with the type-1 evolutions $(\fm^{(n)},\gamma^{(n)})$, which have $\besq(0)$ total mass, by Proposition \ref{type1totalmass}, $n\ge 0$. Since the type-2
total mass process is built from the sum of these, the following additivity lemma will be useful. This extends the well-known additivity of
\besq\ processes with nonnegative parameters and has been taken up in higher generality in \cite{PW18}, where we refer for a proof.

\begin{lemma}[Proposition 1.1 of \cite{PW18}]\label{lmadd}
 Let $X\sim \besq_a(-1)$, $W \sim \besq_b(0)$ and $\underline{Z}\sim \besq_1(-1)$ be independent. 
 Let $\tau = \inf\left\{t\ge 0\colon X_t=0\right\} \wedge \inf\left\{t\ge 0\colon W_t=0  \right\}$.
 Define a process 
 \[
  V_t=\begin{dcases}
   X_t + W_t, \qquad t \le \tau,\\
   Z_{t - \tau}, \qquad\quad\ t> \tau,
  \end{dcases}
 \]
 where $Z_s=(X_\tau + W_{\tau})\underline{Z}_{s/ (X_\tau+ W_\tau)}$, $s\ge 0$. Then $V \sim \besq_{a+b}(-1)$. 
\end{lemma}
%
%

\begin{proof}[Proof of Theorem \ref{thm:type2:total_mass}]
 Consider a type-2 evolution $((m_1^y,m_2^y,\beta^y),y\ge 0)$ as defined in Definition \ref{def:type2:v1}, with initial state $(a,b,\gamma)\in\cJ^\circ$. If $(a,b,\gamma)$ equals 
$(a,0,\emptyset)$ or $(0,b,\emptyset)$, then the result is trivial from the definition, so assume not. Then by Lemma \ref{type2welldef} and Corollary \ref{type2:degenlife}, there 
is a.s.\ some finite $K\ge 0$ such that the degeneration time $D=\inf\{y\ge 0\colon(m_1^y,\beta^y)=(0,\emptyset)\mbox{ or }(m_2^y,\beta^y)=(0,\emptyset)\}$ and the lifetime 
$\zeta=\inf\{y\ge 0\colon(m_1^y,m_2^y,\beta^y)=(0,0,\emptyset)\}$ satisfy $Y_K<D<\zeta=Y_{K+1}$. 
 
 By the strong Markov property and Definition \ref{def:type2:v1}, after time $D$, the type-2 evolution comprises a single non-zero component $m_i^y$, with $i$ being either 1 or 2, evolving as a \BESQ[-1] until its absorption at zero. Let $\underline{Z}$ denote the $\besq_1(-1)$ process obtained by applying \BESQ\ scaling to normalize mass of this component at degeneration: $\underline{Z}_{\,y} := (m_i^{D})^{-1}m_i^{D+m_i^{D}y}$, $y\ge0$. By the strong Markov property, $\underline{Z}$ is independent of the type-2 evolution run up until time $D$.
 
 We define $D_n:=\min\{Y_n,D\}$, $n\ge 0$, so that $D_n=D$ for $n$ sufficiently large, and set
  $$V_y:=m_1^y+m_2^y+\|\beta^y\|,\quad
  V_y^{(n)}:=\left\{\begin{array}{ll}
  		V_y&\mbox{if } y\le D_n,\\
        Z^{(n)}_{y-D_n}&\mbox{if } y>D_n,
    \end{array}\right.$$
    where $Z_s^{(n)}=V_{D_n} \underline{Z}_{\,s/V_{D_n}}$ for $s\ge 0$.
 We will show inductively that all $V^{(n)}$, $n\ge 1$, and hence the a.s.\ limit $V=\lim_{n\rightarrow\infty}V^{(n)}$, are $\besq_{a+b+\|\gamma\|}(-1)$.
 
 For $n=1$, we have $V_y=X_y+W_y$, $0\le y\le D_1$, with $X=\ff^{(0)}\sim\besq_{a}(-1)$, $W=\fm^{(0)} + \|\gamma^{(0)}\|$ independent, and $D_1= \inf\left\{t\ge 0\colon X_t=0\right\} \wedge \inf\left\{t\ge 0\colon W_t=0  \right\}$ is $\tau$ of Lemma \ref{lmadd}. Since $W\sim\besq_{b+\|\gamma\|}(0)$ by Proposition \ref{type1totalmass}, Lemma \ref{lmadd} yields $V^{(1)}\sim\besq_{a+b+\|\gamma\|}(-1)$.

 Now, assume for induction that for some $n\ge1$, $\widehat{V}^{(n)}\sim\besq_{\widehat{a}+\widehat{b}+\|\widehat{\gamma}\|}(-1)$ for all type-2 evolutions  $((\widehat{m}_1^y,\widehat{m}_2^y,\widehat{\beta}^y),y\ge 0)$ starting from any $(\widehat{a},\widehat{b},\widehat{\gamma})\in\cJ^\circ$.  By the strong Markov property, we can apply the inductive hypothesis to 
 the process 
 $(\widehat{m}_1^y,\widehat{m}_2^y,\widehat{\beta}^y) := (m_1^{D_1+y},m_2^{D_1+y},\beta^{D_1+y})$, $y\!\ge\! 0$, on the event $\{Y_1\!=\!D_1\} = \{D\!>\!Y_1\}$.  Then $\widehat D_n = D_{n+1}-D_1$ and $\widehat{\underline{Z}}=\underline{Z}$. We see that
 \begin{align*}
  V_y^{(n+1)} &= \left\{\begin{array}{ll}
  		V_y&\mbox{if } y\le D_{n+1},\\
		Z^{(n+1)}_{y-D_{n+1}}&\mbox{if } y>D_{n+1},
	\end{array}\right.\\
  &= \left\{\begin{array}{ll}
   		V_y&\mbox{if } y\le D_1,\\
		\widehat{V}_{y-D_1}&\mbox{if }D_1<y\le D_1+\widehat{D}_n\\
    	\widehat{Z}^{(n)}_{y-D_1-\widehat{D}_n}&\mbox{if } y>D_1+\widehat{D}_n,
    \end{array}\right\}       
    = \left\{\begin{array}{ll}
      	V_y&\mbox{if } y\le D_1,\\
		\widehat{V}^{(n)}_{y-D_1}&\mbox{if } y>D_1.\end{array}\right.
 \end{align*}
 By the inductive hypothesis, $\widehat{V}^{(n)}\sim\besq_{m_1^{D_1}+m_2^{D_1}+\|\beta^{D_1}\|}(-1)$. By the strong Markov property and \BESQ\ scaling, $\big(\big(\widehat{V}^{(n)}_0\big)^{-1}\widehat{V}^{(n)}_{s\widehat{V}^{(n)}_0},s\ge 0\big)\sim\besq_1(-1)$ is unconditionally independent of $((m_1^y,m_2^y,\beta^y),\, 0\le y\le D_1)$ and in particular of $((X_y,W_y),\,0\le y\le D_1)$.
  Then, by the $n=1$ case already established, we conclude that $V^{(n+1)}\sim\besq_{a+b+\|\gamma\|}(-1)$, as required.
\end{proof}

\section{Type-2 evolutions via interweaving two type-1 point measures}\label{sec:interweaving}

In this section we present another construction of type-2 evolutions from initial states in which the interval partition component is an independent multiple of a $\PDIP[\frac12,\frac12]$ random variable. 
Such interval partitions appear as pseudo-stationary distributions of type-0 and type-1 evolutions, and indeed, we will use this construction in Chapter \ref{dePoiss} to study (pseudo-)stationarity properties of (type-2 evolutions and) unit-mass 2-tree evolutions, as well as projections to 3-mass processes that only retain the evolution of the two top masses and the total mass of the interval partition.

The construction builds a $\cJ^\circ$-valued process from two type-1 evolutions in such a way that the two top masses are taken from the respective type-1 evolution until one of the type-1 evolutions degenerates, while the interval partition is obtained by a procedure that alternates parts from the two interval partitions. We call the 
mechanism that generates this alternation and is based on scaffolding and spindles, \em interweaving\em. This construction is illustrated in Figure \ref{fig:interweaving}.

In the following, we will use notation $\mu^1_m$ for the type-1 pseudo-stationary distribution on $\cJ^\bullet$ with total mass $m\ge 0$ identified in Proposition \ref{prop:01:pseudo}, i.e.\ the distribution of $(mA,m(1-A)\widebar{\gamma})$ for independent $A\sim{\tt Beta}\big(\frac12,\frac12\big)$ and $\widebar{\gamma}\sim{\tt PDIP}\big(\frac12,\frac12\big)$. 
We will write $\mu_{a,c}^1$ for the distribution of $(a,c\widebar{\gamma})$ for $a,c\ge 0$. 
By abuse of notation, we will denote by $\mu_M^1$ and $\mu_{A,C}^1$ the associated mixture distributions, mixed according to the distributions of a random mass $M\ge 0$ or a random pair $(A,C)$ with $A,C\ge 0$, respectively. Finally, we denote by $\mathbf{P}^1_\mu$ the distribution on $\cE\times\cev\cN\times\cN$ of the triple $(\ff,\cev\fN,\fN_\gamma)$ where, conditionally given $(A,\gamma)\sim\mu$, the triple has distribution $\fP^1_{A,\gamma}$ as defined in Construction \ref{type1:construction}.   

\begin{construction}[Interweaving]\label{interweaving} Consider independent $A$ and $B$ for which $\bP(A+B>0)=1$. Also consider independent 
$C_1,C_2\sim\GammaDist(\frac12,\lambda)$ and $\widebar{\gamma}_1,\widebar{\gamma}_2\sim{\tt PDIP}\big(\frac12,\frac12\big)$ independent of $(A,B)$. Let 
\begin{equation}\label{interweav:pointmeas}
(\ff_1,\cev\fN_1,\fN_{\gamma_1})\sim\fP^1_{\mu_{A,C_1\widebar{\gamma}_1}^1}\quad\mbox{and}\quad (\ff_2,\cev\fN_2,\fN_{\gamma_2})\sim\fP^1_{\mu_{B,C_2\widebar{\gamma}_2}^1}
\end{equation}
be independent. Let $\fN_1 := \clade(\ff_1,\cev\fN_1)\concat\fN_{\gamma_1}$, and 
correspondingly define $\fN_2$. We will combine these to define a $\cJ^\circ$-valued process $((\wt m_1^y,\wt m_2^y,\wt\beta^y),\,y\ge 0)$.

\begin{figure}
 \centering
 \scalebox{.9}{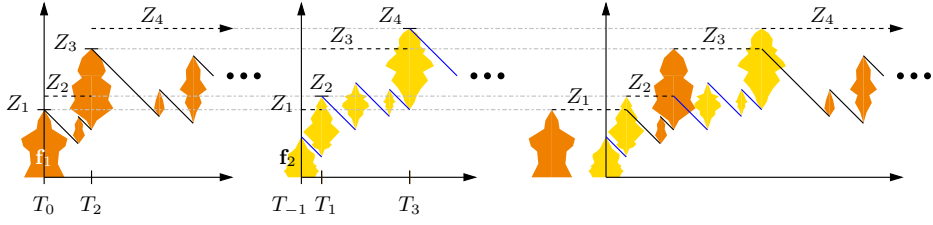}
 \caption{Interweaving is alternating restrictions $\shiftrestrict{\fN_i}{(T_{j-2},T_j]}$ from two type-1 scaffoldings with spindles, $(\fN_1,\fN_2)$. We begin with a single spindle, $\ff_1$, from $\fN_1$. Then, we include spindles from $\fN_2$ until the time $T_1$ at which its scaffolding exceeds the death level $Z_1 = \zeta(\ff_1)$, reaching some higher level $Z_2$. To this, we add spindles from $\fN_1$ until the time $T_2$ at which its scaffolding exceeds level $Z_2$, reaching some higher level $Z_3$, and so on.\label{fig:interweaving}}
\end{figure}

Let $\fX_1 := \xi(\fN_1)$ and $\fX_2 := \xi(\fN_2)$. We set $T_{-1} := T_0 := 0$, $Z_0 := 0$ and $Z_1 := \zeta(\ff_1)$. We define times for each scaffolding, and common levels: for $i\ge 1$, 
\begin{equation}\label{eq:interweave_vars}
\begin{split}
 T_{2i-1} := \inf\{t\ge T_{2i-3}\colon \fX_2(t) > Z_{2i-1}\}, &\qquad Z_{2i} := \fX_{2}(T_{2i-1}),\\
 T_{2i} := \inf\{t\ge T_{2i-2}\colon \fX_1(t) > Z_{2i}\}, &\qquad Z_{2i+1} := \fX_{1}(T_{2i}),
\end{split}
\end{equation}
with the conventions that $\inf(\emptyset) = \infty$ and $\fX_1(\infty)=\infty$ and $\fX_2(\infty)=\infty$. Also note that this includes setting $T_1=0$ if $\zeta(\ff_2)>\zeta(\ff_1)$. Let $p$ denote the \emph{parity} map, sending even numbers to 2 and odd numbers to 1. For $y\ge 0$ we define
\begin{gather}
 \wt I(y) := p\big(\inf\{j\ge 0\colon Z_{j+1} > y\}\big), \qquad \wt J(\infty) := \inf\{j\ge 1\colon T_{j} = \infty\},\notag\\
 \label{eq:interweaving_skewer}
 \left(0,\wt m_{3-\wt I(y)}^y\right) \concat \left(0,\wt m_{\wt I(y)}^y\right) \concat\wt\beta^y := (0,\ff_1(y))\concat (0,\ff_2(y))\concat\widetilde{\theta}(y),
\end{gather}
where
\begin{equation*}
\begin{split}
  \widetilde{\theta}(y) &
     := \skewer\!\left(y \!-\! \zeta(\ff_2), \Restrict{\fN_{2}}{(0,T_{1}]\times\cE}\right)\\[-5.5pt]
  	&\qquad 
       \concat\Concat_{2\le j\le\wt J(\infty)}\! \skewer\!\left(y \!-\! Z_{j-1}, \ShiftRestrict{\fN_{p(j+1)}}{(T_{j-2},T_{j}]\times\cE}\right)\!.
 \end{split}
\end{equation*}
By this we mean that, 
\begin{enumerate}\item[(i)] if the expression on the right of \eqref{eq:interweaving_skewer} has a leftmost block (which equals $(0,\ff_1(y))$ if and only if $y < \zeta(\ff_1)$), then we take $\wt m_{3-\wt I(y)}^y$ to denote the mass of this block, otherwise setting $\wt m_{3-\wt I(y)}^y :=0$; 
  \item[(ii)] if said expression has a second-to-leftmost block, then we denote its mass by $\wt m_{\wt I(y)}^y$, otherwise setting $\wt m_{\wt I(y)}^y := 0$.
\end{enumerate}
 Then $\wt\beta^y$ denotes what remains of $\widetilde{\theta}(y)$ after removing leftmost blocks as required to form $\widetilde{m}_1^y$ and $\widetilde{m}_2^y$, and, if necessary, shifting the remaining interval partition down to line up with 0 on its left end.
\end{construction}

\begin{proposition}\label{prop:interweaving}
 The process $((\wt m_1^y,\wt m_2^y,\wt\beta^y),y\ge 0)$ defined in \eqref{eq:interweaving_skewer} is a type-2 evolution with initial state  
 $(\wt m_1^0,\wt m_2^0,\wt\beta^0)=(A,B,C\widebar{\gamma})$, where $A$, $B$, $C$ and $\widebar{\gamma}$ are jointly independent, with 
 $C\sim\GammaDist(\frac12,\lambda)$ and $\widebar{\gamma}\sim\PDIP(\frac12,\frac12)$.
\end{proposition}

Before we prove this proposition, let us recall from Lemma \ref{lem:type1:pseudo_constr} the construction of a type-1 point measure $\fN_\gamma\sim\fP^1_{\mu_{0,C}^1}$ for $C\sim{\tt Gamma}\big(\frac12,\lambda\big)$. Specifically, $\fN_\gamma=\big(\restrict{\fN}{[0,T)}\big)^{\ge 0}$, where $\fN$ is a ${\tt PRM}({\rm Leb}\otimes\nu_{\tt BESQ})$ on $[0,\infty)\times\cE$, stopped just before its level-0 aggregate mass $M_\fN^0$ exceeds an independent threshold $S\sim{\tt Exponential}(\lambda)$. 

Now, let $\ff$ denote a \BESQ[-1] independent of the other objects, with any random initial mass, and define $\fN_* := \clade(\ff,\cev\fN)\concat\fN_{\gamma}$. In the special case that $\ff(0)\sim\GammaDist[\frac12,\lambda]$, the measure $\fN_*$ describes a pseudo-stationary type-1 evolution with \ExpDist[\lambda] initial mass, as in Proposition \ref{prop:01:pseudo_g}. For any distribution of $\ff(0)$, this construction has the following consequence, by way of the strong Markov property of $\fN$ and the memorylessness of $S$. To state the result, recall 
the notation $T_x(\cev{\fN})$ for the first passage time of $x>0$ by the scaffolding $\xi(\cev{\fN})$ of \eqref{eq:JCCP_def_type0}. 

\begin{lemma}[Memorylessness for some type-1 point measures]\label{lem:pseudostat:memoryless}
 Fix $\lambda>0$ and let $\fN_\gamma$, $\cev{\fN}$ and $\fN_*$ be as above. Let $R$ be a stopping time in the right-continuous time filtration $(\cF_t,\,t\ge0)$ generated by $\fN_*$, i.e.\ the least right-continuous filtration in which $\restrict{\fN_*}{[0,t]}$ is $\cF_t$-measurable for every $t\ge 0$. Given $\restrict{\fN_*}{[0,R]\times\cE}$ with $\xi_{\fN_*}(R) = x$, and further conditioning on 
 $\{\restrict{\fN_*}{(R,\infty)\times\cE}\neq 0\}$, the conditional distribution of $\shiftrestrict{\fN_*}{(R,\infty)\times\cE}$ equals 
 the (unconditioned) distribution of $\ShiftRestrict{\cev\fN}{[T_{x}(\cev\fN),0)\times\cE} \concat \fN_{\gamma}$.
\end{lemma}

\begin{proof}[Proof of Proposition \ref{prop:interweaving}]
 Let $(\ff_1,\ff_2,\cev\fN,\fN_{\gamma})\sim\mathbf{P}_\mu^2$, where $\mu$ is the distribution of $(A,B,C\widebar{\gamma})$ in the setting of the statement of this 
 proposition, and let $\big((m_1^y,m_2^y,\beta^y),\,y\ge0)$ be the type-2 evolution obtained from $(\ff_1,\ff_2,\cev\fN,\fN_\gamma)$ via Construction \ref{deftype2}. We follow the notation of Construction \ref{deftype2}, and further recall
from Lemma \ref{lemma markovprops} and the proof of Theorem \ref{thm:diffusion} notation $J(y)$ for the regime change counter at $y\ge 0$, and $I(y)$ to capture its parity, which we can write here as follows:
\begin{equation}\label{eq:IJ_index_def}
\begin{array}{l}
 J(y) = \inf\{j\ge0\colon Y_{j+1}>y\},\\[3pt]
 I(y) = \text{1 if $J(y)+1$ is odd,\qquad or } I(y)=\text{2 if $J(y)+1$ is even}.
\end{array}
\end{equation}
 Additionally, we define $J(\infty) := \inf\{j\ge 1\colon T_{j}^+ = \infty\}$.
 We prove our assertion by showing:
 \begin{equation}
  \left(0,m_{I(y)}^y\right)\concat \left(0,m_{3-I(y)}^y\right) \concat \beta^y = (0,\ff_1(y))\concat (0,\ff_2(y))\concat \theta(y),\quad y\ge0,\label{eq:alt_deletion_clocking}\\
 \end{equation}
 where
 \begin{equation*}
  \begin{split}
  \theta(y) &:= \skewer\left(y-\zeta(\ff_2), \Restrict{\fN_{*}}{(0,T^+_1]\times\cE} \right)\\ 
  		&\qquad \concat\Concat_{2\le i\le J(\infty)} \skewer\left(y-Y_{i-1}, \ShiftRestrict{\fN_{*}}{(T^-_i,T^+_i]\times\cE} \right)\!;
  \end{split}
 \end{equation*}
 and
 \begin{equation}
 \begin{split}
   &\left(\ff_1,\ff_2,\!\left( \ShiftRestrict{\fN_{p(i+1)}}{(T_{i-2},T_{i}]\times\cE}, Z_{i} \right)\!,1\!\le\! i\!\le\!\wt J(\infty) \right)\\
   &\qquad \stackrel{d}{=} 
  	\left( \ff_1,\ff_2,\!\left( \ShiftRestrict{\fN_{*}}{(T^-_i,T^+_i]\times\cE}, Y_{i} \right)\!,1\le i\!\le\!J(\infty) \right)\!.
  \end{split}\label{eq:inter_clocking_intervals}
 \end{equation}
 These formulas, together with \eqref{eq:interweaving_skewer}, complete the proof.
 
 First, we prove \eqref{eq:inter_clocking_intervals}. For $i\ge 1$, we note the equality of events
 \begin{equation}\label{eq:last_clock_adapted}
 \begin{split}
  \left\{\wt J(\infty)\! =\! i\right\} &= \left\{T_{i}\! =\! \infty;\,\wt J(\infty)\! \ge\! i\right\}\\
  &= \left\{\sup\nolimits_t\xi_{\shiftrestrict{\fN_{p(i+1)}}{(T_{i-2},T_{i}]\times\cE}}(t) < Z_i\!-\!Z_{i-1};\,\wt J(\infty)\! \ge\! i\right\}\!.
 \end{split}
 \end{equation}
 We conclude, by a recursive argument, that the indicator $\mathbf{1}\{\wt J(\infty) \le j\}$ is a function of the $1\le i\le j\wedge \wt J(\infty)$ terms on the left in \eqref{eq:inter_clocking_intervals}. By a corresponding argument, the indicator $\mathbf{1}\{J(\infty) \le j\}$ is a function of the $1\le i\le j\wedge J(\infty)$ terms on the right.
 
 We now establish the base case for an induction. By definition, $T_{-1} = T_1^- = 0$, $Z_1 = Y_1$, and $\fN_2\stackrel{d}{=}\fN_*$. Recall from \eqref{eq:interweave_vars} that $T_1$ is the time when $\fX_2$ first exceeds $Z_1$, while $T_1^+$ in \eqref{eq:d_clk:levels} is the time when $\fX_*$ exceeds $Y_1$. This proves equality in distribution for the $i=1$ terms of \eqref{eq:inter_clocking_intervals}.
 
 Assume for induction that, for some $j\ge 1$, \eqref{eq:inter_clocking_intervals} holds when we substitute $j\wedge\wt J(\infty)$ for the $\wt J(\infty)$ bound on the left and substitute $j\wedge J(\infty)$ for $J(\infty)$ on the right. By the argument following \eqref{eq:last_clock_adapted}, $\bP\{\wt J(\infty) \le j\} = \bP\{J(\infty) \le j\}$. 
 We now show that the conditional distribution of the $(j+1)^{\text{st}}$ term on the left in \eqref{eq:inter_clocking_intervals}, given the preceding terms and the event $\{j < \wt J(\infty)\}$, equals the conditional law of the corresponding term on the right given the preceding terms and the event $\{j < J(\infty)\}$. 
 %
 %
 
 Note that
 \begin{equation*}
 \begin{split}
  Z_{j+1} &= \fX_{p(j+1)}(T_{j}) = Z_{j-1} + \xi_{\shiftrestrict{\fN_{p(j+1)}}{(T_{j-2},T_{j}]\times\cE}}(T_{j}\!-\!T_{j-2})\\
  	&=: G\left(\left(\shiftrestrict{\fN_{p(i+1)}}{(T_{i-2},T_{i}]\times\cE},Z_i\right),\,i\le j\right)
 \end{split}
 \end{equation*}
 and $Y_{j+1} = G\left(\left( \ShiftRestrict{\fN_{*}}{(T^-_i,T^+_i]\times\cE}, Y_{i} \right),\,i\le j\right)$. 
 Next, observe that $\fX_{p(j+1)}(T_{j-2})$ 
 $= Z_{j-1}$ while, correspondingly, $\fX_*(T_{j}^-) = Y_{j-1}$. Since we have conditioned on $\{j < \wt J(\infty)\}$, which means $T_j < \infty$, we may apply Lemma \ref{lem:pseudostat:memoryless} to $\fN_{p(j+1)}$ at this time. In particular, by the independence of $\fN_1$ and $\fN_2$, and by this lemma, given $Z_{j}$, the restricted process $\shiftrestrict{\fN_{p(j+1)}}{(T_{j-2},\infty)\times\cE}$ is conditionally independent of all preceding terms on the left in \eqref{eq:inter_clocking_intervals}. Correspondingly, $\shiftrestrict{\fN_*}{(T_{j}^-,\infty)\times\cE}$ is conditionally independent of all preceding terms on the right in \eqref{eq:inter_clocking_intervals}, given $Y_j$, and these restricted point processes have the same conditional distribution. Finally, $T_{j}^+-T_j^-$ is the first time that $\xi(\shiftrestrict{\fN_{*}}{(T_{j}^-,\infty)\times\cE})$ exceeds $Y_{j}-Y_{j-1}$, and correspondingly for $\fN_{p(j+1)}$. This completes our induction and proves \eqref{eq:inter_clocking_intervals}.
 
 We now prove \eqref{eq:alt_deletion_clocking}. Recall Construction \ref{deftype2}. We distinguish four cases covering the four ways $y$ can be positioned with respect to $\zeta(\ff_1)$ and $\zeta(\ff_2)$.
 
 Case 1: $y < \min\{\zeta(\ff_1),\zeta(\ff_2)\}$. Then $J(y) = 0$, $I(y) = 2$, so the two leftmost blocks on the left hand side of \eqref{eq:alt_deletion_clocking} are $(0,m_1^y)\concat (0,m_2^y)$, which equal $(0,\ff_1(y))\concat (0,\ff_2(y))$, as claimed. By definition, $\fX_*$ is bounded below by $Y_{i}\ge Y_1 > y$ on each interval $(T_{i}^+,T_{i+1}^-]$. Therefore,
 \begin{equation}
  \beta^y = \skewer\left(y - \zeta(\ff_2),\restrict{\fN_*}{(0,\infty)\times\cE}\right) = \theta(y),
 \end{equation}
 as desired. Indeed, $\theta(y)$, as defined following \eqref{eq:alt_deletion_clocking}, simply skips over certain intervals of $\fN_*$ that cannot contribute to the skewer at levels below $\zeta(\ff_1)$.
 
 Case 2: $\zeta(\ff_2)\le y < \zeta(\ff_1)$. Then, again, $J(y) = 0$ and $I(y) = 2$. As before, $(0,m_1^y) = (0,\ff_1(y))$, in agreement with \eqref{eq:alt_deletion_clocking}. However, now $\ff_2(y) = 0$. Thus,
 \begin{equation*}
  (0,m_2^y)\concat\beta^y = \skewer(y,\fN_*) = \skewer\left(y - \zeta(\ff_2),\restrict{\fN_*}{(0,\infty)\times\cE}\right) = \theta(y),
 \end{equation*}
 since, as in Case 1, $\theta(y)$ skips over intervals that do not contribute.
 
 Case 3: $\zeta(\ff_1)\le y < \zeta(\ff_2)$. Then $J(y) = I(y) = 1$ and $T_1^+ = 0$. Then $(0,m_{3-I(y)}^y) = (0,m_2^y) = (0,\ff_2(y))$, while $\ff_1(y)=0$, in agreement with \eqref{eq:alt_deletion_clocking}. Moreover,
 \begin{equation*}
  (0,m_1^y)\concat\beta^y = \skewer(y-Y_1,\shiftrestrict{\fN_*}{(T_2^-,\infty)\times\cE}).
 \end{equation*}
 In this case, since $T_1^- = T_1^+=0$, the first term in the formula for $\theta(y)$ is empty. Then, the concatenation of subsequent terms in $\theta(y)$ equals the above expression, since $\fX_*$ is bounded below by $Y_i \ge Y_2 > y$ on each interval $(T_{i}^+,T_{i+1}^-]$ with $i\ge 2$.
 
 Case 4: $\max\{\zeta(\ff_1),\zeta(\ff_2)\} \le y$. Then $J(y)\ge 1$ and $T_{J(y)+1}^- > 0$. Moreover, $\ff_1(y)=\ff_2(y)=0$, so all that remains on the right in \eqref{eq:alt_deletion_clocking} is $\theta(y)$. Note that $\fX_*$ is bounded above by $Y_{J(y)} \le y$ on each interval $(T_i^-,T_i^+]$ with $i < J(y)$, as well as on $(T_{J(y)}^-,T_{J(y)}^+)$. Then $\fX_*$ jumps up across level $y$ at time $T_{J(y)}^+$, giving rise to the broken spindle $\ff^{(J(y))}$ relating to top mass label $3-I(y)$. Thus, the terms in $\theta(y)$ with $i < J(y)$ do not contribute, and the $i=J(y)$ term contributes only a single block:
 \begin{equation*}
 \begin{split}
  \left(0,m_{3-I(y)}^y\right) &= \left(0,\ff^{(J(y))}(y-Y_{J(y)})\right)\\
  	&= \skewer\left(y-Y_{J(y)-1},\shiftrestrict{\fN_*}{(T_{J(y)}^-,T_{J(y)}^+]\times\cE}\right).
 \end{split}
 \end{equation*}
 Then
 \begin{equation*}
  \left(0,m_{I(y)}^y\right) \concat \beta^y = \skewer\left(y-Y_{J(y)},\shiftrestrict{\fN_*}{(T_{J(y)+1}^-,\infty)\times\cE}\right),
 \end{equation*}
 which equals the concatenation of terms in $\theta(y)$ over $i>J(y)$, since, similarly to the previous cases, this expression skips over intervals where $\fX_*$ is bounded below by $Y_{J(y)+1}>y$.
\end{proof}

\begin{remark}\label{rem:type0fromtype2} After Construction \ref{deftype2} we interpreted the spindles that we ignored/removed during the construction (yellow in Figure \ref{fig:d_clk}) as emigration. Where is the emigration in the interweaving construction, Construction \ref{interweaving}?
The interweaving construction is based on two type-1 evolutions (without emigration). The one that degenerates earlier is completely incorporated into the type-2
evolution, while the other one will only be partially incorporated. Specifically, in the left-to-right order of its scaffolding, the spindles are incorporated up and including its first spindle that exceeds the highest level attained by the scaffolding in the construction of the former type-1 evolution. Following this spindle and starting at its top is a \Stable[\frac32] scaffolding process run until it reaches level 0. The corresponding spindles allow an analogous interpretation of emigration as the (yellow) spindles identified ignored in Construction \ref{deftype2}. 

  Indeed, this part of the marked scaffolding process encodes a type-0 evolution as in Construction \ref{type0:construction} up to its starting height, continued as a type-1 evolution
  as the cutoff point measure above this level, defined as in \ref{eq:upper}, is of the form of Proposition \ref{def:type01_meas}. Furthermore, by the memoryless property of the exponential distribution, the initial distribution of this type-0 evolution is a ${\tt Gamma}(\frac12,\frac12)$ multiple of ${\tt PDIP}(\frac12,\frac12)$, and the type-0 evolution is 
  conditionally independent of the type-2 evolution given the lifetime of the type-2 evolution. 
\end{remark}


\chapter{Unit-mass 2-tree evolutions and stationarity}\label{dePoiss}

In this chapter we establish a stationary variant of the type-2 evolution of Chapter \ref{ch:type-2}, which we will later identify as the $k=2$ case in the consistent system of $k$-tree evolutions of Theorem \ref{thm:intro:k_tree}. Beyond the Markov property and stationarity, we also prove in this chapter that the associated stopped three-mass process is a Wright--Fisher diffusion, hence establishing everything that Theorem \ref{thm:intro:k_tree} claims for $k=2$. In Chapters \ref{ch:constr}--\ref{ch:consistency}, we use and generalize the construction and results of this chapter to obtain more general $k$-tree evolutions and to establish projective consistency properties. Indeed, the special case of this chapter allows us to introduce several of the main techniques in a 
simpler setting so that we can build some familarity before combining them with further structure in the later chapters. 

Specifically, to define the stationary variant, recall from Theorem \ref{thm:type2:total_mass} that type-2 evolutions have $\besq(-1)$ total mass processes, which eventually get absorbed at zero. As a consequence, a type-2 evolution has the same finite \em lifetime \em and is not stationary. We modify the process in two ways: de-Poissonization and resampling. 

De-Poissonization means that we normalize so that the total mass remains constant at one, and then we apply a time-change. De-Poissonization was used in \cite{Paper1-2} to obtain stationary variants of type-0 and type-1 evolutions and has previously been applied in related settings in \cite{Pal11,Pal13,Shiga1990,WarrYor98} and also in further, more recent, related work in \cite{IPPAT,PaperFV,ShiWinkel-2}. 

Resampling is a new idea in this context. 
We saw in Corollary \ref{type2:degenlife} that strictly before reaching the absorbing state $(0,0,\emptyset)$ at the end of their life, type-2 evolutions have a \em degeneration time \em  when they enter a state $(a,b,\gamma)$ with a single block: either $a=\|\gamma\| = 0 < b$ or $b= \|\gamma\| = 0 < a$. In either case, the post-degeneration type-2 evolution will just be a one-dimensional $\besq(-1)$ in the single non-zero component. Resampling will have the process jump instead of degenerating, into an independent state sampled from the law of a Brownian reduced 2-tree; see Proposition \ref{prop:B_ktree}. The state space of these (resampling de-Poissonized) unit-mass 2-tree evolutions is the subspace
\begin{gather}\label{eq:IP_spaces}
 \cJ^*_1 := \{(a,b,\gamma)\in\cJ^\circ\colon a+b+\|\gamma\| = 1;\,a,b,\|\gamma\|<1\}
\end{gather}
of the space $(\cJ^\circ,d^\circ)$ introduced in \eqref{type2spaces}. We also consider the intermediate space $\cJ_1^\circ=\{(a,b,\gamma)\!\in\!\cJ^\circ\colon a\!+\!b\!+\!\|\gamma\|\!=\!1\}$.
Formally, $\cJ^\circ_1$-valued de-Poissonized type-2 evolutions and $\cJ^*_1$-valued (resampling) unit-mass 2-tree evolutions are defined, as follows.

%
Let $\boldsymbol{\mathcal{T}}=(\mathcal{T}^y,y\ge 0)=\big((m_1^y,m_2^y,\beta^y),y\ge 0\big)$ be a type-2 evolution as in Definition \ref{def:type2:v1}. We now consider the distribution $\bP^2_{a,b,\gamma}$ of $\boldsymbol{\mathcal{T}}$ on the space $\bD([0,\infty),\cJ^\circ)$ of c\`adl\`ag functions from $[0,\infty)$ to $\cJ^\circ$. For $T=(a,b,\gamma)\in\cJ^\circ$, we consider the total mass $\|T\|=a+b+\|\gamma\|$. For $\boldsymbol{\mathbf{T}}=(T^y,\,y\ge 0)\in\bD([0,\infty),\cJ^\circ)$, we define a time-change function $\rho_{\boldsymbol{\mathbf{T}}}\colon[0,\infty)\rightarrow[0,\infty]$ by
\begin{equation}\label{eq:dePoi_time_change}
 \rho_{\boldsymbol{\mathbf{T}}}(u)=\inf\left\{y\ge 0\colon\int_0^y\|T^x\|^{-1}dx>u\right\},\quad u\ge 0,
\end{equation}
which is continuous and strictly increasing until a potential absorption at $\infty$. Recall from Theorem \ref{thm:type2:total_mass} that for a type-2 evolution $\boldsymbol{\mathcal{T}}$ starting from $(a,b,\gamma)\in\cJ^\circ\setminus\{(0,0,\emptyset)\}$, we have $\|\boldsymbol{\mathcal{T}}\|:=(\|\mathcal{T}^y\|,\,y\ge 0)\sim\besq_{a+b+\|\gamma\|}(-1)$. By \cite[p.\,314-5]{GoinYor03}, $\rho_{\boldsymbol{\mathcal{T}}}$ is bijective from $[0,\infty)$ onto $[0,\zeta)$ a.s., where $\zeta=\inf\{y\ge 0\colon\|\mathcal{T}^y\|=0\}$.  

\begin{definition}\label{defdePoiss}
 Let $\widebar{\nu}$ be a distribution on $\cJ^\circ_1$. Given a type-2 evolution $\boldsymbol{\mathcal{T}}\sim\bP^2_{\widebar\nu}$ starting according to $\widebar{\nu}$, we associate the \emph{de-Poissonized type-2 evolution} 
$$\widebar{\mathcal{T}}^u=\mathcal{T}^{\rho_{\boldsymbol{\mathcal{T}}}(u)} / \|\mathcal{T}^{\rho_{\boldsymbol{\mathcal{T}}}(u)}\|,\qquad u\ge 0.$$ 
We denote its distribution on $\bD([0,\infty),\cJ^\circ_1)$ by $\widebar{\bP}^{2,-}_{\widebar{\nu}}$.
 \end{definition}

By the bijective property of $\rho_{\boldsymbol{\mathcal{T}}}\colon[0,\infty)\rightarrow[0,\zeta)$ noted above, the degeneration time $D$ of $\boldsymbol{\mathcal{T}}$ gives rise to an a.s.\ finite degeneration time 
$\widebar{D}$ of $\overline{\boldsymbol{\mathcal{T}}}$, which satisfies $\rho_{\boldsymbol{\mathcal{T}}}(\widebar{D})=D$ a.s. Denote by $\widebar\mu$ the distribution on $\cJ^*_1$ of a Brownian reduced 2-tree, i.e. the distribution of $(A_1,A_2,A_3\widebar{\gamma})$, where $(A_1,A_2,A_3)\sim{\tt Dirichlet}\big(\frac{1}{2},\frac{1}{2},\frac{1}{2}\big)$ is independent of the interval partition $\widebar{\gamma}\sim\PDIP\big(\frac12,\frac12\big)$. 

\begin{definition}\label{def:resampling}
 Let $(a,b,\gamma)\in\cJ^*_1$. Let $(\widebar{\mathcal{T}}^u_{\!\!(j)},\,0\le u<\widebar{D}_{(j)})$, $j\ge 0$, be a sequence of 
  independent de-Poissonized type-2 evolutions run until degeneration, with $\widebar{\mathcal{T}}^0_{\!\!(0)} = (a,b,\gamma)$ and $\widebar{\mathcal{T}}^0_{\!\!(j)} \sim \widebar\mu$ for $j\ge1$.  
 Set $V_0=0$ and define the \em resampling times \em
  $V_j=\widebar{D}_{(0)}+\cdots+\widebar{D}_{(j-1)}$, $j\ge 1$. Then the concatenation
  $$\widebar{\mathcal{T}}^{V_j+u}_{\!\!+}=\widebar{\mathcal{T}}^u_{\!\!(j)},\qquad 0\le u<\widebar{D}_{(j)},\ j\ge 0,$$
  is called a \emph{(resampling) unit-mass 2-tree evolution} starting from $(a,b,\gamma)$. We denote its distribution on 
  $\bD([0,\infty),\cJ^*_1)$ by $\widebar{\bP}^{2,+}_{a,b,\gamma}$. For clarity, we continue to use notation $(\widebar{\mathcal{T}}_{\!\!+}^u,u\ge 0)$ for the 
  canonical process on $\bD([0,\infty),\cJ^*_1)$ when working under $\widebar{\bP}^{2,+}_{a,b,\gamma}$.
\end{definition}

Let us state the main results of this chapter here.

\begin{theorem}\label{absorption}
 De-Poissonized type-2 evolutions as defined in Definition \ref{defdePoiss} are Borel right Markov processes absorbed in finite time in either $(1,0,\emptyset)$ or $(0,1,\emptyset)$. 
\end{theorem}

\begin{theorem}\label{thm:stationary}
The unit-mass 2-tree evolutions of Definition \ref{def:resampling} are Borel right Markov process on $(\cJ_1^*,d^\circ)$. 
  Consider $(A_1,A_2,A_3)\sim {\tt Dirichlet}\big(\frac12,\frac12,\frac12\big)$ and an independent interval partition $\widebar\gamma\sim\PDIP\left(\frac12,\frac12\right)$. The law of $(A_1,A_2,A_3\widebar\gamma)$ is the unique stationary distribution for the unit-mass 2-tree evolution.
\end{theorem}



Consider the map $\pi^\bullet_2$ on $\cJ_1^*$ given by $(a,b,\gamma)\mapsto (a,b,\|\gamma\|)$. The range of this map is $\Delta:=\{ (p_1, p_2, p_3) \in [0,1)^3, \sum_{i=1}^3 p_i=1   \}$. Let $\Lambda$ be the stochastic kernel from $\Delta$ to $\cJ_1^*$ that maps $(p_1,p_2,p_3)$ to the law of $(p_1,p_2,p_3\widebar\gamma)$, where $\widebar\gamma\sim\PDIP[\frac12,\frac12]$. Given $(p_1,p_2,p_3)\in\Delta$, run a unit-mass 2-tree evolution $(\widebar{\mathcal{T}}^u,\,u\ge0)$ with initial distribution $\Lambda(p_1,p_2,p_3)$. The \emph{induced $3$-mass process} is then $(X_1(u),X_2(u),X_3(u)) := \pi^\bullet_2(\widebar{\mathcal{T}}^u)$, $u\ge0$. In Appendix \ref{sec:DynkinIntertwining}, we review (Dynkin's criterion and) the Rogers--Pitman intertwining criterion for when a function of a Markov process is again a Markov process. The following result exhibits our first instance of intertwining.

\begin{theorem}\label{thm:wright-fisher}
 The induced $3$-mass process is a recurrent Markovian extension of the Wright--Fisher$(-\frac{1}{2}, -\frac{1}{2}, \frac{1}{2})$ diffusion with generator \eqref{eq:WFgen}, in the following sense. Let 
\[
U=\inf\{u\ge 0\colon X_1(u)=0\mbox{ or }X_2(u)=0\}
\]
be the first time when one of the two top masses vanishes. Then the process killed at $U$ is the killed Wright--Fisher diffusion. The $3$-mass process is intertwined with the unit-mass 2-tree evolution, and it converges to its unique stationary law ${\tt Dirichlet}\left( \frac{1}{2}, \frac{1}{2}, \frac{1}{2} \right)$.
\end{theorem}

Notice that the 3-mass process jumps back into the interior of the simplex immediately after either of the first two coordinates vanish. This extension of the generalized Wright--Fisher diffusion is natural from the perspective of the (modified) Aldous chain (Definition \ref{def:modified_AC}) as the continuum analogue of the construction in \cite{Paper2}. Indeed it suggests an extension of the scaling limit result \eqref{eq:3massconv}, here in the (rooted) case with $k=2$ top masses,  
that was observed by Aldous \cite{AldousDiffusionProblem} and Pal \cite{Pal13}, relating killed 3-mass processes in the discrete and continuum settings.

Before turning to de-Poissonized processes, in Section \ref{sec:type2:pseudostat} we prepare for stationarity arguments by establishing pseudo-stationary behaviour of the type-2 evolutions of Definition \ref{def:type2:v1}. Due to degeneration, this takes a slightly different form to the corresponding results for type-0 and type-1 evolutions in Proposition \ref{prop:01:pseudo}, and will be complemented by further pseudo-stationary behaviour at degeneration in Section \ref{sec:type2:pseudodeg}. We will also use both sets of results in Chapters \ref{ch:constr}--\ref{ch:consistency}. We then turn to de-Poissonized type-2 evolutions and prove Theorem \ref{absorption} in Section \ref{sec:type2:dePoiss}, and to (resampling) unit-mass 2-tree evolutions and the proofs of Theorems \ref{thm:stationary} and \ref{thm:wright-fisher} in Section \ref{sec:type2:resampling}. Finally, we establish some H\"older estimates for interval-partition-valued variants of type-2 evolutions in Section \ref{sec:Holder}, which we will use in Chapter \ref{ch:properties} to establish the path-continuity of the induced continuum-tree-valued process.


\section{Pseudo-stationary type-2 evolutions}
\label{sec:type2:pseudostat}

Recall from Corollary \ref{type2:degenlife} that type-2 evolutions degenerate to a single block of positive mass before reaching zero total mass, while type-1 evolutions degenerate when they reach zero total mass and type-0 evolutions do not degenerate (and are not absorbed) when they reach zero total mass. In Propositions \ref{prop:01:pseudo}--\ref{prop:01:pseudo_g}, we recalled the pseudo-stationarity behaviour for types 0 and 1, which we may read as conditional on non-degeneration. In this section, we establish pseudo-stationarity of type-2 evolutions, again conditionally given that degeneration has not yet happened. 

\begin{proposition}[Pseudo-stationarity of type-2 evolution]\label{prop:type2:pseudo}
 Let $\widebar\gamma\!\sim\!\PDIP(\frac12,\frac12)$, $(A_1,A_2,A_3)\!\sim\!{\tt Dirichlet}\big(\frac12,\frac12,\frac12\big)$ and $M(0)\!>\!0$ be independent and $(m_1^y,m_2^y,\beta^y)$, $y\ge 0$, a type-2 evolution started from $(M(0)A_1,M(0)A_2,M(0)A_3\widebar\gamma)$. Let $M(y)$, $y\ge0$, denote its total mass process. For fixed $y>0$, given $\{D>y\}$, the total mass $M(y)$ is conditionally independent of $(m_1^y/M(y),m_2^y/M(y),\beta^y/M(y))$. The latter is conditionally distributed according to the (unconditioned) law of $(A_1,A_2,A_3\widebar\gamma)$.
\end{proposition}

In light of this result, we refer to the law of $(M(0)A_1,M(0)A_2,M(0)A_3\widebar\gamma)$ above as the \emph{pseudo-stationary law for a type-2 evolution} with mass $M(0)$. 
Following the strategy of proof in \cite{Paper1-2} of Proposition \ref{prop:01:pseudo} above, we first prove this for $M(0)\sim\GammaDist[\frac{3}{2},\lambda]$, $\lambda>0$, and then generalize via Laplace inversion. 

\begin{proposition}\label{prop:type2:pseudo_g}
 Consider a type-2 evolution $((m_1^y,m_2^y,\beta^y),y\ge 0)$ with initial blocks $(m_1^0,m_2^0)$ independent of $\beta^0 = M\widebar\gamma$,
 where $M\sim \GammaDist[\frac{1}{2},\lambda]$ and $\widebar\gamma\sim\PDIP(\frac12,\frac12)$ are independent. 
 Then for $y>0$, given $\{D>y\}$, the interval partition $\beta^y$ is conditionally independent of $(m_1^y,m_2^y)$, conditionally distributed according to the (unconditional) law of $(2y\lambda+1)M\widebar\gamma$.
 
 If, additionally, $m_1^0$ and $m_2^0$ are i.i.d.\ \GammaDist[\frac12,\lambda], then given $\{D>y\}$, $m_1^y$ and $m_2^y$ are conditionally i.i.d.\ \GammaDist[\frac12,\lambda/(2y\lambda+1)].
\end{proposition}

\begin{proof}
 Let $(\ff_1,\ff_2,\cev\fN,\fN_{\gamma})$ be point measures and spindles for such an evolution, as in Construction \ref{deftype2}. From Lemma \ref{lem:type1:pseudo_constr}, we may assume $\fN_\gamma = \big(\restrict{\fN}{[0,T)}\big)^{\ge 0}$, where $\fN$ is a \PRM[\Leb\otimes\mBxc]\ on $[0,\infty)\times\cE$ and $T$ is the time at which the aggregate mass of spindles crossing level 0, as defined in \eqref{eq:skewer_def}, first exceeds an independent mass threshold $S\sim \ExpDist[\lambda]$. 
 
 We follow the notation of \eqref{eq:IJ_index_def}, in which $J(y)$ is the number of regime switches up to time $y$ and $I(y)$ denotes the index, 1 or 2, that records the alternating regime of the construction at scaffolding level $y$. So $m_{3-I(y)}^y$ is the top mass that is part of a type-1 evolution in this construction at that level, $\beta^y$ is the interval partition of remaining, ``spinal'' masses, while $m_{I(y)}^y$ is the further top mass. We set $\fN_* := \clade(\ff_2,\cev\fN)\concat\fN_{\gamma}$. Let $((m_*^z,\beta_*^z),z\ge0)$ denote the type-1 evolution $\skewerbar(\fN_*)$.
  It follows from Construction \ref{deftype2} that, on $\{D>y\}$, one top mass at level $y$ is the mass of a spindle found in $\restrict{\fN_*}{[0,T)}$ at the stopping time $$R = \inf\{t>T_{J(y)+1}^- \colon \xi_{\fN_*}(t)>y\}<T,$$ and the interval partition $\beta^y$ equals $\skewer(y - \xi_{\fN_*}(R),\restrict{\fN_*}{(R,T)})$. 
 
 Let $R' := \inf\{t\!>\!R\colon \xi_{\fN_*}(t)\!=\!y\}$ and $T_y(\fN_*) := \inf\{t\!\ge\!0\colon \xi_{\fN_*}(t)\! =\! y\}$. 
 By Lemma \ref{lem:pseudostat:memoryless}, the conditional law of $\shiftrestrict{\fN_*}{(R',\infty)}$ given $(\ff_1,\restrict{\fN_*}{[0,R']})$ and $\{D>y\}$ equals the conditional law of $\shiftrestrict{\fN_*}{(T_y(\fN_*),\infty)}$ given $\{m_*^y+\|\beta_*^y\|>0\}$. Passing to the skewers, the correspondingly conditioned laws of $\beta^y$ and $\beta^y_*$ are equal. 
 By Proposition \ref{prop:01:pseudo_g}, this is an independent \GammaDist[\frac12,\lambda/(2y\lambda+1)]\ multiple of a \PDIP[\frac12,\frac12]. This also implies that $\beta^y$ is conditionally independent of $(m_1^y,m_2^y)$ given $\{D>y\}$, proving the first assertion of the proposition.
 
 To prove the second assertion, we apply Proposition \ref{prop:interweaving}. In the representation there, with notation in and after \eqref{interweav:pointmeas}, we can express 
 $D = \min\{\zeta_1,\zeta_2\}$ in terms of the degeneration times $\zeta_i$ of the two type-1 evolutions $\skewerP(\fN_i)$, $i=1,2$. In particular, conditioning on $\{D>y\}$ is the same as conditioning on $\{\zeta_1>y,\zeta_2>y\}$. \linebreak By Proposition \ref{prop:01:pseudo_g} and the independence of the two pseudo-stationary type-1 evolutions in that construction, $m^y_1$ and $m^y_2$ are conditionally independent given $\{\zeta_1>y,\zeta_2>y\}$, with common distribution $\GammaDist[\frac{1}{2},\lambda/(2y\lambda+1)]$.
\end{proof}

\begin{proposition}\label{prop:pseudo_f}
 For $a,b,c>0$ and $\widebar\gamma\sim\PDIP[\frac12,\frac12]$, consider a type-2 evolution starting from $(a,b,c\widebar{\gamma})$. Let $\widebar\gamma'$ be an independent $\PDIP[\frac12,\frac12]$, and let $\wt\beta^y$ denote $\beta^y/\|\beta^y\|$ when $\beta^y\neq\emptyset$ (this holds a.s.\ given $y<D$), or $\widebar\gamma'$ otherwise. Then for $y>0$, $\wt\beta^y$ is independent of $(m_1^y,m_2^y,\|\beta^y\|)$ and has law $\PDIP[\frac12,\frac12]$.
\end{proposition}
\begin{proof}
 For $\lambda>0$, consider $B_{\lambda}\sim \GammaDist[\frac12,\lambda]$ independent of all other objects. By decomposing according to the events $\{D>y\}$ and $\{D\le y\}$, and applying the first assertion of Proposition \ref{prop:type2:pseudo_g} in the former case, we see that 
 for all continuous $f\colon \BR^3\to[0,\infty)$ and $g\colon\cI\to[0,\infty)$,
  \begin{align*}
    &\int_0^\infty \sqrt{\frac{\lambda}{\pi x}}e^{-\lambda x}\bE^2_{a,b,x\widebar\gamma}\big[ f(m_1^y,m_2^y,\|\beta^y\|)g(\wt\beta^y) \big]dx\\
    &\qquad\qquad = \bE^2_{a,b,B_{\lambda}\widebar\gamma}\big[ f(m_1^y,m_2^y,\|\beta^y\|)g(\wt\beta^y) \big]\\
    &\qquad\qquad = \bE\big[g(\widebar\gamma)\big]\int_0^\infty \sqrt{\frac{\lambda}{\pi x}}e^{-\lambda x} \bE^2_{a,b,x\widebar\gamma}\big[ f(m_1^y,m_2^y,\|\beta^y\|) \big]dx.
  \end{align*}
  We cancel factors of $\sqrt{\lambda}$. By the uniqueness of Laplace transforms,
  $$\bE^2_{a,b,x\widebar\gamma}\big[ f(m_1^y,m_2^y,\|\beta^y\|)g(\wt\beta^y) \big] = \bE^2_{a,b,x\widebar\gamma}\big[ f(m_1^y,m_2^y,\|\beta^y\|) \big]\bE\big[g(\widebar\gamma)\big]$$
  for a.e.\ $x>0$. By Proposition \ref{contini}, the right-hand side is continuous in $x$. 
  Note that $(a,b,\gamma)\mapsto f(a,b,\|\gamma\|)g(\gamma/\|\gamma\|)\cf\{\gamma\neq\emptyset\}$ is $\bP^2_{a,b,x\widebar\gamma}$-a.s.\ continuous at $(m_1^y,m_2^y,\beta^y)$. Thus, the left-hand side is continuous in $x$ as well; see e.g.\  \cite[Theorem 4.27]{Kallenberg}. We conclude that the above formula holds for every $x$.
\end{proof}

\begin{proof}[Proof of Proposition \ref{prop:type2:pseudo}]
 Let $((m_1^z,m_2^z,\beta^z),\,z\ge0)$ be as in the statement of the proposition, and fix $y>0$. The conditional law of $(m_1^y,m_2^y,\beta^y)$ given $D>y$ can be obtained as a mixture, over the law of the vector $(MA_1,MA_2,MA_3)$ of initial masses, of the conditional laws described in Proposition \ref{prop:pseudo_f}. In particular, conditionally given $\{\beta^y\neq\emptyset\}$, $\beta^y/\|\beta^y\|\sim\PDIP[\frac12,\frac12]$, conditionally independent of $(m_1^y,m_2^y,\|\beta^y\|)$. To prove that $(m_1^y,m_2^y,\|\beta^y\|)/M(y)$ then has conditional law ${\tt Dirichlet}\big(\frac12,\frac12,\frac12\big)$, we make an argument similar to that in the proof of Proposition \ref{prop:pseudo_f}.
 
 Recall the standard beta-gamma algebra that a ${\tt Dirichlet}(x_1,\ldots,x_n)$ vector, multiplied by an independent \GammaDist[x_1+\cdots+x_n,\lambda] scalar, gives rise to a vector of independent variables, with the $j^{\text{th}}$ having law \GammaDist[x_j,\lambda]. Let $(\widetilde m_1^y,\widetilde m_2^y,\widetilde m_3^y)$ denote $(m_1^y/M(y),m_2^y/M(y),\|\beta^y\|/M(y))$ when $y<D$ or $(A_1',A_2',A_3')$ otherwise, where the latter is an independent ${\tt Dirichlet}\big(\frac12,\frac12,\frac12\big)$. By the second assertion of Proposition \ref{prop:type2:pseudo_g}, for $\lambda>0$ and measurable $f\colon \BR^3\to[0,\infty)$ we have
  \begin{align*}
    \int_0^\infty 2\sqrt{\frac{x\lambda^3}{\pi}}e^{-\lambda x}\bE^2_{A_1x,A_2x,A_3x\widebar\gamma}\big[ f(\widetilde m_1^y,\widetilde m_2^y,\widetilde m_3^y) \big]dx
    = \bE\big[f(A_1,A_2,A_3)\big].
  \end{align*}
  Multiplying the right-hand side by $\int_0^\infty 2\sqrt{x\lambda^3/\pi}e^{-\lambda x}dx=1$, canceling factors of $\lambda^{3/2}$, and appealing to uniqueness of Laplace transforms and Proposition \ref{contini}, as in the previous proof, gives the desired result.
\end{proof}

For our next results, we require a scaling invariance property of type-2 evolutions. We recall the scaling invariance of type-1 evolutions from Proposition \ref{prop:type01:diffusion}, which states that for any type-1 evolution $((m^y,\beta^y),y\ge 0)$ and any $c>0$, the process $((cm^{y/c},c\beta^{y/c}),y\ge 0)$ is also a
type-1 evolution.

Together with the well-known scaling invariance of squared Bessel processes (see e.g.\ \cite[Appendix A.3]{GoinYor03}), the corresponding result for type-2 evolutions
follows straight from their definition in Definition \ref{def:type2:v1}.

\begin{lemma}\label{lem:type2:scaling}
For any type-2 evolution $((m_1^y,m_2^y,\beta^y),y\ge 0)$ and any $c>0$, the process $((cm_1^{y/c},cm_2^{y/c},c\beta^{y/c}),y\ge 0)$ is also a
type-2 evolution.
\end{lemma}

We denote by $\mu_m^2$ the pseudo-stationary distribution on $\cJ^\circ$ with total mass $m$, and by $\mu_{a,b,c}^2$ the distribution on $\cJ^\circ$ of $(a,b,c\widebar{\gamma})$, with $\widebar{\gamma}\sim\PDIP[\frac12,\frac12]$,
for all $(a,b,c)\in[0,\infty)$ with either $a+b>0$ or $a=b=c=0$. 

\begin{lemma}[Strong pseudo-stationarity]\label{pretotal}
 \begin{enumerate}[label=(\roman*), ref=(\roman*)]
  \item Let $\widebar{\gamma}\sim\PDIP[\frac12,\frac12]$, and let  
  $(A,B,C)$ be an independent vector for which, with probability 1, at least two components are positive. Consider a type-2 evolution 
  $((m_1^y,m_2^y,\beta^y),\,y\ge 0)$ with initial state $(A,B,C\widebar{\gamma})$. Denote by $(A(y),B(y),C(y))=(m_1^y,m_2^y,||\beta^y||)$ the associated 3-mass process and by 
  $(\cF_{\rm 3-mass}^y,\,y\ge 0)$ the right-continuous filtration it generates. Let $Y$ be a stopping time in this filtration. Then for all 
  $\cF^Y_{\rm 3-mass}$-measurable $\eta\colon\Omega\rightarrow[0,\infty)$ and all measurable 
  $H\colon\cJ^\circ\rightarrow[0,\infty)$,
  \begin{equation}\label{eq:3massstrongpseudo}
    \bE\left[\eta H(m_1^Y,m_2^Y,\beta^Y)\right]=\bE\big[\eta\mu^2_{A(Y),B(Y),C(Y)}[H]\big].
  \end{equation}
  I.e.\ conditionally given $\cF_{\rm 3-mass}^Y$, the state of the type-2 evolution at time $Y$ is distributed as  
  $(A(Y),B(Y),C(Y)\widebar{\gamma}^\prime)$ for independent $\widebar{\gamma}^\prime\sim\PDIP[\frac12,\frac12]$.  \label{item:pretotal:3}
  \item  Now consider instead a type-2 evolution $(\Gamma^y,y\ge 0)=((m_1^y,m_2^y,\beta^y),\,y\ge 0)$, whose initial state is an independent multiple of a random state with unit-mass pseudo-stationary distribution. 
  Denote by $M(y) = m_1^y+m_2^y+\|\beta^y\|$, $y\ge 0$, the associated total mass process and by $(\cF^y_{\rm mass},\,y\ge 0)$ the right-continuous filtration it generates. Let $Y$ be a stopping time in this filtration. Then for all $\cF^Y_{\rm mass}$-measurable $\eta\colon\Omega\!\rightarrow\![0,\infty)$ and measurable $H\colon\cJ^\circ\!\rightarrow\![0,\infty)$,
  $$\bE\left[\left.\eta H(M(Y)^{-1}\Gamma^Y)\,\right|D>Y\right]=\bE\left[\left.\eta\,\right|D>Y\right]\mu_1[H].\vspace{-0.6cm}$$\label{item:pretotal:1}%
 \end{enumerate}
\end{lemma}
\begin{proof} Results similar to (ii) for type-0 and type-1 evolutions were obtained in \cite[Lemma 4.7 and Theorem 4.8]{Paper1-2}. In the following, we adapt the proofs of those
   results to the present setting. 
 
   (i) If we further condition on $\{D\le Y\}$, then the statement follows trivially as $\|\beta^Y\|=0$. Now suppose first that $Y=y$ is non-random and consider the event $\{D>y\}$.
   We will simplify notation and write $ABC(y)=(A(y),B(y),C(y))$. By the pseudo-stationarity of the interval partition in Proposition \ref{prop:pseudo_f}, we have for all bounded measurable $f_0\colon\mathbb{R}^3\rightarrow[0,\infty)$ and
   $H\colon\mathcal{J}^\circ\rightarrow[0,\infty)$
  \begin{align}\label{eq:nequals1}
    &\mathbb{E}_{\mu_{A,B,C}^2}^2\Big[f_0\left(ABC(0)\right)\mathbf{1}\{D>y\}H(m_1^y,m_2^y,\beta^y)\Big]\\
    &=\int_{\mathcal{I}}\mathbb{E}_{\mu_{A,B,C}^2}^2\Big[f_0(ABC(0))\mathbf{1}\{D>y\}H(A(y),B(y),C(y)\gamma)\Big]\mathbb{P}(\widebar{\gamma}\in d\gamma)\nonumber\\
    &=\mathbb{E}_{\mu_{A,B,C}^2}^2\Big[f_0(ABC(0))\mathbf{1}\{D>y\}\mu^2_{ABC(y)}[H]\Big].\nonumber
  \end{align}
   Consider $0=y_0<y_1<\cdots<y_n<y_{n+1}$ and further bounded measurable functions $f_1,\ldots,f_{n+1}\colon\mathbb{R}^3\rightarrow[0,\infty)$. We will apply the 
   Markov property of the type-2 evolution at time $y_1$ and write $\widebar{y}_j=y_{j+1}-y_1$, $0\le j\le n$. Inductively, consider
  \begin{align*}
   &\mathbb{E}_{\mu_{A,B,C}^2}^2\Bigg[\prod_{j=0}^{n+1}f_j\left(ABC(y_j)\right)\mathbf{1}\{D>y_{n+1}\}H(m_1^{y_{n+1}},m_2^{y_{n+1}},\beta^{y_{n+1}})\Bigg]\\
   &=\mathbb{E}_{\mu_{A,B,C}^2}^2\Bigg[f_0\left(ABC(0)\right)\mathbf{1}\{D>y_1\}\\
   &\qquad\qquad\quad\ \ \mathbb{E}_{m_1^{y_1},m_2^{y_1},\beta^{y_1}}^2\Bigg[\prod_{j=0}^{n}f_{j+1}\left(ABC(\widebar{y}_j)\right)\mathbf{1}\{D\!>\!\widebar{y}_{n}\}H(m_1^{\widebar{y}_{n}},m_2^{\widebar{y}_{n}},\beta^{\widebar{y}_{n}})\Bigg]\Bigg].
  \end{align*}
  Applying \eqref{eq:nequals1}, this further equals
  \begin{align*}
   &\int_{\mathcal{I}}\mathbb{E}_{\mu_{A,B,C}^2}^2\Bigg[f_0\left(ABC(0)\right)\mathbf{1}\{D>y_1\}\\
   &\quad\ \ \mathbb{E}_{A(y_1),B(y_1),C(y_1)\gamma}^2\!\Bigg[\prod_{j=0}^{n}f_{j+1}(ABC(\widebar{y}_j))\mathbf{1}\{D\!>\!\widebar{y}_{n}\}H(m_1^{\widebar{y}_{n}},m_2^{\widebar{y}_{n}},\beta^{\widebar{y}_{n}})\Bigg]\Bigg]\mathbb{P}(\widebar{\gamma}\!\in\! d\gamma)\\
   &=\mathbb{E}_{\mu_{A,B,C}^2}^2\Bigg[f_0\left(ABC(0)\right)\mathbf{1}\{D>y_1\}\\
   &\qquad\qquad\quad\ \ \mathbb{E}_{\mu^2_{ABC(y_1)}}^2\Bigg[\prod_{j=0}^{n}f_{j+1}\left(ABC(\widebar{y}_j)\right)\mathbf{1}\{D>\widebar{y}_{n}\}H(m_1^{\widebar{y}_{n}},m_2^{\widebar{y}_{n}},\beta^{\widebar{y}_{n}})\Bigg]\Bigg].
  \end{align*}
  Applying the induction hypothesis and taking the other steps in reverse, we conclude that
  \begin{align*}
    &\mathbb{E}_{\mu_{A,B,C}^2}^2\left[\prod_{j=0}^{n+1}f_j\left(ABC(y_j)\right)\mathbf{1}\{D>y_{n+1}\}H(m_1^{y_{n+1}},m_2^{y_{n+1}},\beta^{y_{n+1}})\right]\\
   &=\mathbb{E}_{\mu_{A,B,C}^2}^2\Bigg[f_0\left(ABC(0)\right)\mathbf{1}\{D>y_1\}\\
   &\qquad\qquad\quad\ \ \mathbb{E}_{\mu^2_{ABC(y_1)}}^2\Bigg[\prod_{j=0}^{n}f_{j+1}\left(ABC(\widebar{y}_j)\right)\mathbf{1}\{D>\widebar{y}_{n}\}\mu^2_{ABC(\widebar{y}_{n})}[H]\Bigg]\Bigg]\\
   &=\mathbb{E}_{\mu_{A,B,C}^2}^2\left[\prod_{j=0}^{n+1}f_{j}\left(ABC(y_j)\right)\mathbf{1}\{D>y_{n+1}\}\mu_{ABC(y_{n+1})}^2[H]\right]. 
  \end{align*}
  This yields \eqref{eq:3massstrongpseudo} in the case $Y=y$ by a monotone class theorem. The generalization to stopping times is standard, first considering continuous $H$ 
  and approximating $Y$ by $Y_n=2^{-n}\lfloor 2^nY+1\rfloor\wedge 2^n$, and noting that firstly, $\eta_k:=\eta\mathbf{1}\{Y_n=k2^{-n}\}$ is 
  $\mathcal{F}^{k2^{-n}}_{\rm 3-mass}$-measurable, $k\in[2^{2n}-1]$; secondly, $\{D>Y_n\}$ increases to $\{D>Y\}$ up to a null set; thirdly, type-2 evolutions
  are (right-)continuous along $Y_n\downarrow Y$. 

  (ii) The proof of (i) is easily adapted, using the total mass process instead of the 3-mass process, applying Proposition \ref{prop:type2:pseudo} instead of Proposition \ref{prop:pseudo_f} and finally conditioning on $\{D>Y\}$.
%
\end{proof}

\begin{proposition}\label{prop:3mass:Poi}
 Consider a type-2 evolution $((m_1^y,m_2^y,\beta^y),y\ge 0)$ starting from $(x_1,x_2,x_3\widebar{\gamma})$ for 
  $\widebar{\gamma}\sim\PDIP(\frac12,\frac12)$. Then the associated 3-mass process $((m_1^y,m_2^y,\|\beta^y\|),y\ge 0)$ is a Markov process starting from $(x_1,x_2,x_3)$. 
\end{proposition}
\begin{proof} We check the Rogers--Pitman intertwining criterion (see Appendix \ref{sec:DynkinIntertwining}).
  Consider the map $\phi(x_1,x_2,\beta) = (x_1,x_2,\|\beta\|)$ and the stochastic kernel $\Lambda((x_1,x_2,x_3),A)=\bP((x_1,x_2,x_3\widebar{\gamma})\in A)$, where $\widebar{\gamma}\sim\PDIP(\frac12,\frac12)$. Then clearly $\Lambda((x_1,x_2,x_3),\phi^{-1}(\{(x_1,x_2,x_3)\}))=1$ and
  by Lemma \ref{pretotal}\ref{item:pretotal:3}, we also have 
  $$\bP_\mu^2\big((m_1^y,m_2^y,\beta^y)\in A\ \big|\ (m_1^y,m_2^y,\|\beta^y\|),(m_1^0,m_2^0,\|\beta^0\|)\big)=\Lambda((m_1^y,m_2^y,\|\beta^y\|),A)\quad a.s.,$$
  for all initial distributions $\mu$ of the form $\Lambda((x_1,x_2,x_3),\,\cdot\,)$, as required.
\end{proof}

The semigroup of the 3-mass process can be described as ``replace the third component by a scaled $\PDIP(\frac12,\frac12)$, make type-2 evolution transitions, and then project the interval partition onto its mass.''

\section{Degeneration in pseudo-stationarity}
\label{sec:type2:pseudodeg}

In the interweaving construction, Construction \ref{interweaving}, particularly in the pseudo-stationary case with ${\tt Gamma}(\frac32,\lambda)$ initial total mass, it is easy to describe the distribution of the degeneration time.

\begin{proposition}\label{prop:degeneration}
 Fix $\lambda>0$. Let $D$ be the degeneration time of a type-2 evolution starting from $(A,B,C\widebar{\gamma})$, where $A$, $B$, $C$ and 
 $\widebar{\gamma}$ are jointly independent, with $C\sim\GammaDist(\frac12,\lambda)$ and $\widebar{\gamma}\sim\PDIP(\frac12,\frac12)$. Then
 $\bP(D>y)=\bP(\zeta_1>y)\bP(\zeta_2>y)$ for all $y>0$, where $\zeta_1$ and $\zeta_2$ are the lifetimes (which are also the degeneration times) of two independent type-1 evolutions starting from
 $(A,C\widebar{\gamma})$ and $(B,C\widebar{\gamma})$, respectively. 

 If also $A,B\sim\GammaDist(\frac12,\lambda)$, then $\bP\{D>y\} = (2y\lambda+1)^{-2}$ for all $y>0$.
\end{proposition}


\begin{proof}
 Construction \ref{interweaving}, interweaving, is such that on $\{\wt J(\infty)\mbox{ even}\}$,
 $$\left(0,\widetilde{m}_1^y\right)
	\concat\Concat_{2\le j\le\wt J(\infty)\textrm{ even}} 
             \skewer\left(y - Z_{j-1}, \ShiftRestrict{\fN_{p(j+1)}}{(T_{j-2},T_{j}]\times\cE}\right)$$
 equals $\skewer\left(y,\fN_1\right)$ and
 \begin{equation*}
 \begin{split}
  \left(0,\widetilde{m}_2^y\right)
	&\concat\skewer\left(y - \zeta(\ff_2), \Restrict{\fN_{2}}{(0,T_{1}]\times\cE}\right)\\
	&\concat\Concat_{2\le j\le\wt J(\infty)\textrm{ odd}} 
			 \skewer\left(y - Z_{j-1}, \ShiftRestrict{\fN_{p(j+1)}}{(T_{j-2},T_{j}]\times\cE}\right)
 \end{split}
 \end{equation*}
 equals $\skewer\big(y,\restrict{\fN_2}{[0,T_{\wt J(\infty)-1}]}\big)$. On the other hand, on $\{\wt J(\infty)\mbox{ odd}\}$, the first displayed expression is $\skewer\big(y,\restrict{\fN_1}{[0,T_{\wt J(\infty)-1}]}\big)$, the second is $\skewer\left(y,\fN_2\right)$.
 
 On $\{\wt J(\infty)\mbox{ even}\}$, the definitions of $Z_{\wt J(\infty)}$ and $T_{\wt J(\infty)}$ imply that
 $\zeta_2\ge\zeta=Z_{\wt J(\infty)}>\zeta_1=D$, where $\zeta_i$ is the lifetime of the type-1 evolution $\skewerbar(\fN_i)$,
 $i=1,2$, and $\zeta$ is the lifetime of the type-2 evolution $((\wt m_1^y,\wt m_2^y,\wt\beta^y),y\ge 0)$. Together with corresponding observations on 
 $\{\wt J(\infty)\mbox{ odd}\}$, we see that $D$ is the minimum of the lifetimes $\zeta_1$ and $\zeta_2$ of two type-1 evolutions.

 If we apply the interweaving construction to independent $A,B\sim\GammaDist[\frac12,\lambda]$, we 
 obtain the type-1 pseudo-stationary initial distribution. As these are i.i.d.\ with $\ExpDist(\lambda)$ initial mass, from \eqref{eq:pseudo:type_1_survival} they each have lifetime at least $y$ with probability $(2y\lambda+1)^{-1}$. The minimum of two i.i.d.\ variables with 
 this law has probability $(2y\lambda+1)^{-2}$ of exceeding $y$, as claimed.
\end{proof}

\begin{proposition}\label{degdist}
 Consider a type-2 evolution $((m_1^y,m_2^y,\beta^y),\,y\ge 0)$ starting from the initial condition of Proposition \ref{prop:type2:pseudo} with $M\sim{\tt Gamma}\big(\frac{3}{2},\lambda\big)$, with degeneration time $D$. Let $A = \{I(D)=1\}$; this is the event that $(m_1^y,\,y\ge0)$ is the surviving top mass process at the time of degeneration. On this event, $(m_1^D,m_2^D,\beta^D)=(M^D,0,\emptyset)$; on the complementary event, $(m_1^D,m_2^D,\beta^D)=(0,M^D,\emptyset)$. Then $\bP(A) =\frac12$, the event $A$ is independent of $(D,M^D)$, and $\GammaDist[\frac{1}{2},\lambda/(2\lambda y+1)]$ is a regular conditional distribution for $M^D$ given $D=y$.
\end{proposition}
\begin{proof}
 We are interested in the joint distribution of $(D,m_1^D,m_2^D,\beta^D)$. Using Construction \ref{interweaving} from two independent type-1 evolutions with lifetimes $\zeta_1$ and $\zeta_2$ and top mass processes $m_1$ and $m_2$, we have $D=\min\{\zeta_1,\zeta_2\}$, $A = \{\zeta_1>\zeta_2\}$, and
  $$(D,M^D)=(\zeta_2,m_1^{\zeta_2})\cf_A+(\zeta_1,m_2^{\zeta_1})\cf_{A^c}.$$
  Under the stated initial conditions, these two type-1 evolutions are in fact i.i.d. From this, it is clear by symmetry that $\bP(A)=\frac12$ and $A$ is independent of $(D,M^D)$, as claimed.
  
  For all nonnegative measurable $f$ and $g$ on $\bR$,
  $$\bE\big[f(D)g\big(M^D\big)\big]=\bE\big[f(\zeta_2)g\big(m_1^{\zeta_2}\big)\cf_A\big]+\bE\big[f(\zeta_1)g\big(m_2^{\zeta_1}\big)\cf_{A^c}\big].$$ 
  We use Proposition \ref{prop:01:pseudo} to rewrite the first term on the right hand side as
  \begin{equation*}
  \begin{split}
   &\int_0^\infty f(y)\bE[g(m_1^y)\mathbf{1}{\{\zeta_1>y\}}]\bP(\zeta_2\in dy)\\
       &=\int_0^\infty f(y)\bP\{\zeta_1>y\}\int_0^\infty g(x)\sqrt{\frac{\lambda}{\pi x(2\lambda y+1)}}\exp\left(-\frac{\lambda x}{2\lambda y+1}\right)dx\bP(\zeta_2\in dy)\\
    &=\bE\!\left[f(\zeta_2)\mathbf{1}{\{\zeta_1>\zeta_2\}}\int_0^\infty g(x)\sqrt{\frac{\lambda}{\pi x(2\lambda \zeta_2+1)}}\exp\left(-\frac{\lambda x}{2\lambda\zeta_2+1}\right)dx\right]\!.
  \end{split}
  \end{equation*}
The second term can be written similarly, by symmetry, and together they give
  $$\bE\big[f(D)g\big(M^D\big)\big] = \bE\!\left[f(D)\int_0^\infty g(x)\sqrt{\frac{\lambda}{\pi x(2\lambda D+1)}}\exp\!\left(-\frac{\lambda}{2\lambda D+1}x\right)\!dx\right]\!.$$
 This proves the claimed regular conditional distribution for $M^D$.
\end{proof}

Note that this result (and proof) formalizes an extension of the second part of Proposition \ref{prop:type2:pseudo_g} to the random time $y=D$, the degeneration time, and yields the same conditional distribution for the single surviving top mass as for the two surviving top masses when conditioning on $y<D$.

\section{De-Poissonized type-2 evolutions}
\label{sec:type2:dePoiss}


Recall that Theorem \ref{absorption} claims that de-Poissonized type-2 evolutions (i.e. time-changed and normalized to unit mass, but without resampling) are $\cJ_1^\circ$-valued Borel right Markov processes and that they reach one of the two absorbing states, $(1,0,\emptyset)$ or $(0,1,\emptyset)$, in finite time.

\begin{proof}[Proof of Theorem \ref{absorption}]
 $\cJ^\circ_1$ is a Borel subset of a Lusin space, and is therefore Lusin. Both continuous time changes and normalization on $\cJ^\circ\setminus\{(0,0,\emptyset)\}$ preserve the property of sample paths being c\`adl\`ag.  The strong Markov property of Theorem \ref{thm:diffusion} and the continuity in the initial state of Proposition \ref{contini} transfer to the de-Poissonized processes as in \cite[Proposition 4.6, proof of Theorem 1.6]{Paper1-2}. 

  Degeneration occurs at the time $\widebar{D}$ that satisfies $\rho_{\boldsymbol{\mathcal{T}}}(\widebar{D}) = D$ since $D<\zeta$ a.s.; note that all states $(m,0,\emptyset)$, $m\in(0,\infty)$, are normalized to $(1,0,\emptyset)$, and similarly for $(0,1,\emptyset)$. 
Finally, the time-change is such that $\rho_{\boldsymbol{\mathcal{T}}}(u)<\zeta$ for all $u\in[0,\infty)$. 
\end{proof}

Pal \cite{Pal11,Pal13} studied Wright--Fisher diffusions with positive and negative real parameters $\theta_1,\ldots,\theta_n$ as de-Poissonized processes associated with vectors of independent $Z_i\sim\besq(2\theta_i)$, $1\le i\le n$. Combining the arguments of \cite[Proposition 11]{Pal11} and \cite[Theorem 4]{Pal13}, we may define generalized Wright--Fisher diffusions (running at 4 times the speed of \cite{Pal11,Pal13}) via generators \eqref{eq:WFgen}, or as weak solutions to certain systems of stochastic differential equations, or, as is relevant for us, as
\begin{equation}\label{eq:WF_constr}
 \widebar{Z}_i(u)=\frac{Z_i(\rho(u))}{Z_+(\rho(u))},\ \ 1\le i\le n,\ \ 0\le u\le\widebar{\tau},
\end{equation}
where $\widebar\tau := \inf\{s\ge 0\colon \exists i\text{ s.t.\ }\widebar{Z}_i(s)=0\}$, $Z_+(y) := \sum_{i=1}^nZ_i(y)$, and $\rho(u)$ is as in \eqref{eq:dePoi_time_change}, but with $Z_+(x)$ in place of $\|\mathcal{T}^x\|$ inside the integral. 
See also \cite{Paper1-3}.

\begin{proposition}\label{prop:3mass:dePoi}
 Let $\overline{\boldsymbol{\mathcal{T}}}=((X_1(u),X_2(u),\gamma(u)),u\ge 0)\sim\widebar{\bP}^{2,-}_{a,b,\gamma}$ be a de-Poissonized type-2 evolution starting from any initial state $(a,b,\gamma)\in\cJ^\circ_1$. Let 
 \linebreak
 $U=\inf\{u\ge 0\colon X_1(u)=0\mbox{ or }X_2(u)=0\}$. Then the 3-mass process 
$$((X_1(u),X_2(u),1-X_1(u)-X_2(u)),\,0\le u\le U)$$ 
is a generalized Wright--Fisher process with parameter vector $(-\frac{1}{2},-\frac{1}{2},\frac{1}{2})$.
  
  If furthermore the initial state is taken as $\gamma=(1-a-b)\widebar{\gamma}$ for $\widebar{\gamma}\sim\PDIP(\frac12,\frac12)$, then the 3-mass process $((X_1(u),X_2(u),1-X_1(u)-X_2(u)),\,u\ge 0)$ is a Markovian extension of the generalized Wright--Fisher process.
\end{proposition}

\begin{proof}
 For the first claim, we assume without loss of generality that $\overline{\boldsymbol{\mathcal{T}}}$ is constructed as in Definition \ref{defdePoiss} from a type-2 evolution $\boldsymbol{\mathcal{T}}=(\mathcal{T}^y,\,y\ge 0)$ arising from $(\ff_1,\ff_2,\cev{\fN},\fN_\gamma)\sim\fP_{a,b,\gamma}^2$ as in Construction \ref{deftype2}. By Proposition \ref{type1totalmass}, we have $\|\skewerbar(\cev{\fN}\concat\fN_\gamma)\|\sim\besq_{\|\gamma\|}(1)$. This process, together with the independent processes $\ff_1\sim\besq_a(-1)$ and $\ff_2\sim\besq_b(-1)$ forms a triple of \besq\ processes, as in the paragraph above the proposition. Thus, we can construct a generalized Wright--Fisher process from $Z_1=\ff_1$, $Z_2=\ff_2$ and $Z_3=\|\skewerbar(\cev{\fN}\concat\fN_\gamma)\|$, as in \eqref{eq:WF_constr}. 
 Since $Z_+(y)=\|\mathcal{T}^y\|$ for $0\le y\le\tau:=\inf\{y\ge 0:\exists i\text{ s.t.\ }Z_i(y)=0\}$, we have $\rho(u)=\rho_{\boldsymbol{\mathcal{T}}}(u)$, and hence $X_i(u)=\widebar{Z}_i(u)$, $i=1,2$, and $X_3(u)=1-X_1(u)-X_2(u)=1-\widebar{Z}_1(u)-\widebar{Z}_2(u)=\widebar{Z}_3(u)$, $0\le u\le U=\widebar{\tau}$. This completes the proof.   
  
  The second claim follows from Proposition \ref{prop:3mass:Poi} and the observation that the (Poissonized) 3-mass process of that proposition can be de-Poissonized by the same scaling/time-change operation as the type-2 evolution, as the scaling and time change only depend on the common total mass process.
\end{proof}

\section{Resampling and stationarity of unit-mass 2-tree evolutions}
\label{sec:type2:resampling}

As we have seen in Theorem \ref{absorption}, de-Poissonized type-2 evolutions degenerate at a finite random time $\widebar{D}<\infty$ in one of the two absorbing states $(1,0,\emptyset)$ and $(0,1,\emptyset)$. In this section we consider resampling evolutions that are restarted, as defined in Definition \ref{def:resampling}, instead of entering the absorbing states. Informally and with a stationary Aldous diffusion in mind, we take the opportunity to sample afresh from the reduced Brownian CRT at each degeneration time. Recall the state space $\cJ_1^*$ of \eqref{eq:IP_spaces}, and recall that Theorem \ref{thm:stationary} claims that unit-mass 2-tree evolutions (resampling de-Poissonized type-2 evolutions) are Borel right Markov with the Brownian reduced 2-tree as their
unique stationary distribution.


\begin{proof}[Proof of Theorem \ref{thm:stationary}]
 To confirm that the unit-mass 2-tree evolution is a Borel right Markov process, we only need to check the strong Markov property. Given the construction of Definition \ref{def:resampling}, with resampling times $(V_n)$, this follows as an application of general results about resurrecting Markov processes \cite{Mey75}.

  Now, we prove that $\widebar\mu$, defined before Definition \ref{def:resampling}, is the unique stationary distribution, and that the process converges to it. Applying Lemma \ref{pretotal}\ref{item:pretotal:1} to the $(\cF^y_{\rm mass},y\ge 0)$-stopping time $Y=\rho_{\boldsymbol{\mathcal{T}}}(u)$, we find
  \begin{equation}\label{statpre}
   \bE_{\widebar\mu}^2\left[g\big(\mathcal{T}^{\rho_{\boldsymbol{\mathcal{T}}}(u)}/\|\mathcal{T}^{\rho_{\boldsymbol{\mathcal{T}}}(u)}\|\big)\mathbf{1}_{\{D>\rho_{\boldsymbol{\mathcal{T}}}(u)\}}\right] =
   \bP_{\widebar\mu}^2\big\{D>\rho_{\boldsymbol{\mathcal{T}}}(u)\big\}{\widebar\mu}\big[g\big].
  \end{equation}
  Now consider a (resampling) unit-mass 2-tree evolution $(\widebar{\mathcal{T}}_{\!\!+}^u,u\ge 0)$ with initial distribution ${\widebar\mu}$. We use the notation of Definition \ref{def:resampling}. Let $U\sim\ExpDist(\lambda)$ independent of the unit-mass 2-tree evolution. Then \eqref{statpre} yields
  $$\bE\left[g\left(\widebar{\mathcal{T}}_{\!\!+}^U\right)\mathbf{1}_{\{U<V_1\}}\right] = \bP\{U<V_1\}{\widebar\mu}[g].$$ 
  For $m\ge 1$,
  \begin{align*}
  &\bE\left[g\left(\widebar{\mathcal{T}}_{\!\!+}^U\right)\mathbf{1}_{\{V_m\le U<V_{m+1}\}}\right]
	=\int_0^\infty\lambda e^{-\lambda u}\bE\left[g\left(\widebar{\mathcal{T}}^u_{\!\!+}\right)\mathbf{1}_{\{V_m\le u<V_{m+1}\}}\right]du\\
	&\qquad =\bE\left[e^{-\lambda V_m}\int_0^\infty\lambda e^{-\lambda s}g\left(\widebar{\mathcal{T}}_{\!\!+}^{V_m+s}\right)\mathbf{1}_{\{V_m+s<V_{m+1}\}}ds\right]\\
	&\qquad =\int_0^\infty\lambda e^{-\lambda s}\bE\left[e^{-\lambda V_m}g\left(\widebar{\mathcal{T}}_{\!\!+}^{V_m+s}\right)\mathbf{1}_{\{V_m+s<V_{m+1}\}}\right]ds\\
	&\qquad =\int_0^\infty\lambda e^{-\lambda s}\bE\left[e^{-\lambda V_m}\bE\left[g(\widebar{\mathcal{T}}_{\!\!+}^s)\mathbf{1}_{\{s<V_1\}}\right]\right]ds\\
	&\qquad =\bE\left[e^{-\lambda V_m}\bE\left[g\left(\widebar{\mathcal{T}}_{\!\!+}^U\right)\mathbf{1}_{\{U<V_1\}}\right]\right]\\
	&\qquad =\bE\left[e^{-\lambda V_m}\right]\bP\big\{U<V_1\big\}{\widebar\mu}\big[g\big] = \bP\big\{V_m\le U<V_{m+1}\big\}{\widebar\mu}\big[g\big].
  \end{align*}
  Summing over $m$ and inverting Laplace transforms in $\lambda$, we find that\linebreak
   $\widebar{\bE}^{2,+}_{\widebar\mu}[g(\widebar{\mathcal{T}}^u)]={\widebar\mu}[g]$ for all $u\ge 0$, i.e.\ ${\widebar\mu}$ is stationary for $\overline{\boldsymbol{\mathcal{T}}}_{\!\!+}$. Furthermore, since resampling is according to the stationary distribution ${\widebar\mu}$, we have for any other initial distribution ${\widebar\nu}$ that for all bounded measurable $g\colon\cJ^*_1\rightarrow[0,\infty)$
  $$\widebar{\bE}_{\widebar\nu}^{2,+}\big[g(\widebar{\mathcal{T}}^u)\big] = \widebar{\bE}_{\widebar\nu}^{2,+}\left[g(\widebar{\mathcal{T}}^u)\mathbf{1}_{\{u<V_1\}}\right]+\widebar{\bP}_{\widebar\nu}^{2,+}\big\{V_1\le u\big\}{\widebar{\mu}}\big[g\big]\rightarrow{\widebar{\mu}}\big[g\big],$$ 
  since $V_1=\widebar{D}$ is finite $\widebar{\bP}_{\widebar\nu}^{2,+}$-a.s. In particular, the stationary distribution is unique.
\end{proof}

As in Propositions \ref{prop:3mass:Poi} and \ref{prop:3mass:dePoi}, we can project a unit-mass 2-tree evolution to a \emph{resampling 3-mass process},  $((\widebar m_1^y,\widebar m_2^y,\widebar\beta^y),\,y\ge0) \mapsto \big((\widebar m_1^y,\widebar m_2^y,\|\widebar\beta^y\|),\,y\ge0\big)$. 
We can now prove Theorem \ref{thm:wright-fisher}, which identifies this projected process as a Markovian extension of a Wright--Fisher diffusion.

\begin{proof}[Proof of Theorem \ref{thm:wright-fisher}]
Proposition \ref{prop:3mass:dePoi} and Theorem \ref{thm:stationary} imply that 
the resampling 3-mass process is a Borel right Markov process that extends the generalized Wright--Fisher process to a recurrent process on the simplex 
$\{(a,b,c)\in [0,1)^3:$ $a+b+c=1\}$, which has 
${\tt Dirichlet}\big(\frac12,\frac12,\frac12\big)$ stationary distribution, and converges to stationarity. 

The intertwining relationship was noted in the proof of Proposition \ref{prop:3mass:Poi} for the type-2 evolution. Using Lemma \ref{pretotal}(i) and the fact that resampling is into stationarity, this extends to the present de-Poissonized setting
with resampling.
\end{proof}

Note that the Wright--Fisher diffusion with parameters $\left(\frac12,\frac12,\frac12\right)$ has this same invariant law. See e.g.\ Ethier and Kurtz \cite[Lemma 4.1]{EthiKurt81}. 

Recall the definition in \eqref{eq:diversity} of the diversity $\sD_\beta$ of an interval partition $\beta\in\cI$.

\begin{corollary}
 Under $\widebar{\bP}_{\widebar\nu}^{2,+}$, let $(\widebar\beta^u,\,u\ge 0)$ denote the evolution of the interval partition component. Then the total diversity process $(\sD_{\widebar\beta^u}(\infty),\,u\ge 0)$ is continuous except at the resampling times $V_m$, $m\ge 1$.
\end{corollary}
\begin{proof}
First consider a type-2 evolution starting according to $\widebar{\nu}$. Up until the first regime change, the interval partition component differs from an 
$\cI$-valued type-1 evolution by at most one interval. The continuity of total diversity follows from Proposition \ref{prop:type01:diffusion}. At regime changes, the interval partition component loses a single interval, without affecting diversity. An induction extends this up to degeneration. De-Poissonization maintains the continuity of total diversity. The same argument proves continuity between any two consecutive resampling times under  $\widebar{\bP}_{\widebar\nu}^{2,+}$.
\end{proof}

\section{Other state spaces for 2-trees and H\"older estimates}\label{sec:Holder}

The spaces $\cJ^\circ$ and $\cJ^\circ_1$ of interval partitions with two top masses are not the only state spaces in which 2-trees such as those obtained as reduced trees in a BCRT
$(\cT,d,\rho,\mu)$ can be represented. Indeed, recall from Section \ref{sec:intro:B_k_tree} that we defined the Brownian reduced $2$-tree in three steps. In the first step, we followed Aldous \cite{AldousCRT3} and considered the subtree $\cR_2^+$ of $\cT$ spanned by $\rho$ and two leaves $\Sigma_1$ and $\Sigma_2$ randomly sampled from $\mu$, a Y-shaped tree with a unique branch point $v$. In \cite{PitmWink09}, $\cR_2^+$ was equipped with the measure $\mu_2^+$ obtained by projecting $\mu$ onto $\cR_2^+$. Before doing such a projection, our second step was to further reduce to just the trunk $\cR_2=[\![\rho,v]\!]_\cT$. The third step is a projection of $\mu$ to $\cR_2$. But the image measure $\mu_2$ under this projection has an atom at $v$ that adds the masses $X_1^{(2)}$ and $X_2^{(2)}$ of the two connected components of 
$\cT\setminus\cR_2$ containing $\Sigma_1$ and $\Sigma_2$, respectively. We instead recorded these masses separately as top masses, and we represented the remainder of 
$\mu_2$ by the interval partition $\beta_{\{1,2\}}^{(2)}$ that contains intervals of lengths corresponding to the atoms sizes of $\mu_2$ on $]\!]v,\rho]\!]_\cT$, which in turn capture the masses of the other connected components of $\cT\setminus\cR_2$.

In this section we consider $(\cR_2,\mu_2)$ and related representations of $(a,b,\gamma)\in\cJ^\circ$ as 
\begin{equation}\label{eq:metric2tree}
M_2(a,b,\gamma):=\bigg([0,\sD_\gamma(\infty)]\ ,\ \ (a\!+\!b)\delta(\sD_\gamma(\infty))+\sum_{U\in\gamma}\Leb(U)\delta(W_\gamma(U))\bigg),
\end{equation}
where $W_\gamma(U)=\sD_\gamma(\infty)-\sD_\gamma(U)$, $U\in\gamma$.
As a consequence of \cite[Theorem 2.5(a)--(b)]{Paper1-0}, the map $M_2$ is Lipschitz continuous from $(\cJ^\circ,d^\circ)$ into the space
\begin{equation}\label{eq:onebranch}
\cM=\{(C,\nu)\colon C\subset [0,\infty)\text{ compact},\nu\text{ finite Borel measure on }[0,\infty)\},
\end{equation}
equipped with the sum $d_{\rm HP}$ of the Hausdorff metric on compact subsets of $[0,\infty)$ and the Prokhorov metric on finite Borel measures on $[0,\infty)$. We write 
$\cM^\circ:=M_2(\cJ^\circ)$ and $\cM_1^\circ:=M_2(\cJ_1^\circ)$.

It is evident that $M_2\colon\cJ^\circ\rightarrow\cM^\circ$ is not one-to-one. Moreover, the $\cM^\circ$-valued projection $(M_2(\Gamma^y),y\ge 0)$ of a
$\cJ^\circ$-valued type-2 evolution $(\Gamma^y,y\ge 0)$ cannot be expected to be Markovian. Indeed, recall from Remark \ref{rem:LMB} how a top mass in a type-1 evolution 
(and hence in a type-2 evolution) interacts with the interval partition. In $\cM^\circ$, both top masses contribute to the atom at 0. Informally, while both top masses are large, the 
(type-0) evolution of the interval partition does not contribute atoms at 0 almost surely at any fixed time, but when a top mass vanishes, the interval partition provides new top 
masses. When this last happened and how large the last atom was can be seen in the history of $(M_2(\Gamma^y),y\ge 0)$, is not recorded in the current state in $\cM^\circ$, but 
is relevant for the further evolution. 

\begin{remark}\label{nostrongMarkov} The reader may wonder if this loss of the simple Markov property could be avoided by studying evolutions of $(\cR_2,\mu_2)$ rather than 
  $(\cR_2^\circ,\mu_2^\circ)$. We think this is true. Indeed, the natural starting point for this would be an evolution of the string of beads $(\cR_1,\mu_1)$ of 
  \cite{PitmWink09} obtained by projecting $\mu$ onto $[\![\rho,\Sigma_1]\!]_\cT$. However, we have been unable to devise an evolution that is compatible with the Aldous 
  chain and handles the delicate behaviour in the neighborhood of a leaf. 

  In any case, the strong Markov property would still be lost. Even a type-0 evolution $(\gamma^y,y\ge 0)$ similarly represented as $(M_0(\gamma^y),y\ge 0)$, where 
  $M_0(\gamma):=\big([0,\sD_\gamma(\infty)]\,,\ \sum_{U\in\gamma}\Leb(U)\delta(\sD_\gamma(U))\big)$, would fail to be strong Markov. Specifically, suppose that
  $\gamma^0=\beta_0\concat(0,a_1)\concat\beta_1\concat(0,a_2)\concat\beta_2$. Consider independent type-0 evolutions $(\beta_i^y,y\ge 0)$ starting from $\beta_i$, $i=0,1,2$, 
  as in Construction \ref{type0:construction} and ${\tt BESQ}_{a_i}(-1)$ evolutions $\ff_i$, $i=1,2$. By Proposition \ref{prop:type1:construction}, $(0,\ff_i(y))\concat\beta_i^y$, 
  $0\le y\le\zeta(\ff_i)$, $i=1,2$, are independent type-1 evolutions stopped when the top mass vanishes. Now consider the random time 
  $\eta=\inf\{y\ge 0\colon\beta_1^y=\emptyset\}$. Then $\bP(\eta<\min\{\zeta(\ff_1),\zeta(\ff_2)\})>0$, so at time $\eta$, on this event, the interval partition 
  $\gamma^\eta=\beta_0^\eta\concat(0,\ff_1(\eta))\concat(0,\ff_2(\eta))\concat\beta_2^\eta$ has two blocks with the same diversity, hence their atoms in $M_0(\gamma^\eta)$ 
  add up to a single atom. After time $\eta$, the type-0 evolution $(\beta_1^{\eta+z},z\ge 0)$ separates these two atoms again in the same way as before. Only the sum of atoms
  is recorded in $M_0(\gamma^\eta)$, but the split, part of the history of the process, is relevant for the future.  
\end{remark}

But, again informally, the $\cM^\circ$-valued process $(M_2(\Gamma^y),y\ge 0)$ is $d_{\rm HP}$-continuous since atom sizes and atom locations (diversities) evolve continuously. We will make this precise in a framework that is (still) easier to handle, interval partitions.

\begin{definition}\label{def:type2:IPvalued} Let $((m_1^y,m_2^y,\beta^y),y\ge 0)$ be a type-2 evolution with regime change times $(Y_n)$ as defined in Definition \ref{def:type2:v1}, and let $(I(y),y\ge0)$ be the $\{1,2\}$-valued parity process that flips at each regime change as in \eqref{eq:IJ_index_def}. Then we refer to 
  $$\gamma^y=(0,m_{I(y)}^y)\concat(0,m_{3-I(y)}^y)\concat\beta^y,\qquad y\ge 0,$$
  as an \emph{$\cI^\circ$-valued type-2 evolution}, where 
  $$\cI^\circ:=\{(0,m)\concat\gamma\colon m\in(0,\infty),\gamma\in\cI\}\cup\{\emptyset\}\subset\cI.$$
\end{definition}

\begin{remark}\label{prop:type2:IPvalued} $\cI^\circ$-valued type-2 evolutions are, in fact, path-continuous Hunt processes in their own right. As we will not need this result, we do not provide a formal proof, but we point out that the (Borel right) Markov property follows from \eqref{augment} by Dynkin's argument, and we will establish path-continuity when starting from initial conditions that correspond to the pseudo-stationary distributions of $\cJ^\circ$-valued type-2 evolutions. 
\end{remark}

Let us denote by $\widetilde\mu$ the distribution of the 
$\cI^\circ$-valued interval partition 
$$(0,\widebar{A})\star(0,\widebar{B})\star \widebar{G}\,\widebar{\beta}$$
for independent $(\widebar{A},\widebar{B},\widebar{G})\sim{\tt Dirichlet}\big(\frac12,\frac12,\frac12\big)$ and 
$\widebar{\beta}\sim{\tt PDIP}\big(\frac12,\frac12\big)$. This distribution is not pseudo-stationary in the strong sense that the distribution of an 
$\cI^\circ$-valued type-2 evolution starting from $\widetilde\mu$ has as marginal distributions the distributions of random multiples of this interval 
partition -- intuitively, the leftmost block is stochastically larger than the second block. However, we will be able to appeal to the 
pseudo-stationarity of $\cJ^\circ$-valued type-2 evolutions starting from $(\widebar{A},\widebar{B},\widebar{G}\,\widebar{\beta})$ in 
situations that treat the two top masses symmetrically.


\begin{proposition}\label{type2holder}
 Let $(\widetilde{\gamma}^y,y\ge 0)$ be an $\cI^\circ$-valued type-2 evolution starting according to $\widetilde\mu$. Let $\theta\in(0,\frac14)$ and $y>0$. Then 
 there is a random H\"older constant $L=L_{\theta,y}$ with moments of all orders such that
 $$d_\cI(\widetilde{\gamma}^a,\widetilde{\gamma}^b)\le L|b-a|^\theta\qquad\mbox{for all }0\le a<b\le y.$$ 
\end{proposition} 

The remainder of this subsection is devoted to the proof of this proposition. We begin by some preliminary considerations. Let us first consider the type-2 evolution $\gamma^y=M\widetilde{\gamma}^{y/M}$, $y\ge0$, with \GammaDist[\frac32,\lambda] initial mass $M$ for some $\lambda>0$, cf.\ Proposition \ref{prop:type2:pseudo_g}. Recall that 
\begin{itemize}
  \item the evolution $(\gamma^y,y\!\ge\! 0)$ can be constructed by interweaving two
    independent pseudo-stationary 
    type-1 evolutions of initial mass 
    ${\tt Exponential}(\lambda)$, see Propositions \ref{prop:interweaving} and \ref{prop:type2:pseudo};
  \item such pseudo-stationary type-1 evolutions consist of a type-1 evolution starting from a single 
    interval $(0,A)$ with $A\sim\GammaDist\big(\frac12,\lambda\big)$ concatenated left-to-right with an independent type-1 evolution starting from a $\PDIP\big(\frac12,\frac12\big)$ scaled by an independent $\GammaDist\big(\frac12,\lambda\big)$ mass, see Proposition \ref{prop:01:pseudo_g}; 
  \item a type-1 evolution starting from a single $\GammaDist\big(\frac12,\lambda\big)$-distributed interval $(0,A)$ can be constructed
    from a $\besq_A(-1)$ process, with death level $\zeta$ and an independent $\Stable\big(\frac32\big)$ process $\mathbf{X}$  
    with $\besq(-1)$ excursions in its jumps and run until it first descends to $-\zeta$, see Proposition \ref{def:type01_meas}; adding $\zeta$ we obtain a \em descent from $\zeta$ to $0$\em;
  \item a type-1 evolution starting from $\PDIP\big(\frac12,\frac12\big)$ scaled by mass $\GammaDist\big(\frac12,\lambda\big)$ can be constructed from a $\Stable\big(\frac32\big)$ process $\widetilde{\fX}$ starting from 0, with $\besq(-1)$ excursions in its jumps stopped at a time $\widetilde{T}$, which is the left endpoint of the excursion away from 0 where the mass at level 0 exceeds an independent ${\tt Exponential}(\lambda)$ threshold, see Lemma \ref{lem:type1:pseudo_constr}.
\end{itemize}
Access to H\"older bounds is via local times of ${\tt Stable}\big(\frac32\big)$ processes. Recall that
\begin{itemize}
  \item the local times $(\widetilde{\ell}^y(t),0\le t\le\widetilde{T},y\ge 0)$ of 
    $(\widetilde{\fX}(t),0\le t\le\widetilde{T})$ have the property that for each $a\ge 0$ and $\theta\in(0,\frac14)$, the random variable 
    $$\widetilde{D}_\theta^a=\sup_{0\le t\le\widetilde{T},0\le x<y\le a}\frac{|\widetilde{\ell}^x(t)-\widetilde{\ell}^y(t)|}{|y-x|^\theta}$$
    has moments of all orders (\cite[Theorem 3]{Paper0});
  \item 
  for a type-1 evolution $(\beta^y,y\!\ge\!0)$ arising from ${\tt Stable}\big(\frac32\big)$ scaffolding $\fX$ marked by \BESQ[-1] spindles, 
  it is a.s.\ the case that for every $y$ and every block $U\in\beta^y$, the diversity $\mathscr{D}_{\beta^y}(U)$ equals the local time $\ell^y(t)$ in $\fX$, up to the time $t$ at which the spindle corresponding to block $U$ arises (\cite[Theorem 1]{Paper0}).
\end{itemize}
Consider the $\Stable\big(\frac32\big)$ process starting from $\zeta$ obtained by concatenating the descent $\fX+\zeta$ from $\zeta$ to $0$ before $\wt \fX$. Denote this process by $(\widehat{\fX}(t),0\le t\le\widehat{T})$ and its local times by $(\hat{\ell}^y(t),0\le t\le\widehat{T},y\ge 0)$. 
In this context, \cite[Theorem 3]{Paper0} has the following consequence. 

\begin{lemma}\label{lem:LT_Holder} 
 The following random variable has moments of all orders:
  $$\widehat{D}_\theta^a=\sup_{0\le t\le\widehat{T},0\le x<y\le a}\frac{|\hat{\ell}^x(t)-\hat{\ell}^y(t)|}{|y-x|^\theta}.$$
\end{lemma}
\begin{proof}
 Let $\widetilde{H}_\zeta=\inf\{t\ge 0\colon\widetilde{\fX}(t)=\zeta\}$. Then the event $\{\widetilde{H}_\zeta<\widetilde{T}\}$ has positive probability. 
 %
  As with the memorylessness property of Lemma \ref{lem:pseudostat:memoryless}, the conditional distribution given $\widetilde{H}_\zeta<\widetilde{T}$ of the process
  $(\widetilde{\fX}(\widetilde{H}_\zeta+s),0\le s\le\widetilde{T}-\widetilde{H}_\zeta)$ is the same as the unconditional distribution of $(\widehat{\fX}(t),0\le t\le \widehat{T})$. 
  The idea is the same as in Lemma \ref{lem:pseudostat:memoryless}; one just needs to show that the negative parts of $\widetilde{\fX}$ do not affect the argument. 
  Then the associated local times 
  $(\wt\ell^y(\widetilde{H}_\zeta+s)-\wt\ell^y(\widetilde{H}_\zeta),0\le s\le\widetilde{T}-\widetilde{H}_\zeta,y\ge 0)$ have as their conditional distribution the
  distribution of $(\hat{\ell}^y(t),0\le t\le\widehat{T},y\ge 0)$. By the triangle inequality,
  $$\bE\big[(\widehat{D}_\theta^a)^p\big] \le \bE\big[(2\widetilde{D}_\theta^a)^p\big|\widetilde{H}_\zeta<\widetilde{T}\big]<\infty.\vspace{-0.4cm}$$
\end{proof}
This allows us to bound terms (i) and (ii) of Definition \ref{def:IP:metric} of $d_\cI$, which deal with diversity.
\begin{lemma}\label{lem:LT_matching_bd}
  There is a random variable $L_\theta$ with moments of all orders such that uniformly over all
  correspondences $((U_j,U_j^\prime),1\le j\le m)$ from $\gamma^0$ and $\gamma^y$ that are taken from the same 
  $\besq(-1)$ excursion, we have
  \begin{equation*}
   \big|\sD_{\gamma^y}(\infty)-\sD_{\gamma^0}(\infty)\big| \le L_\theta y^\theta \quad \text{and} \quad \max_{1\le j\le m}\big|\sD_{\gamma^y}(U'_j)-\sD_{\gamma^0}(U_j)\big| \le L_\theta y^\theta.
  \end{equation*}
\end{lemma}
\begin{proof}
  Think of $(\gamma^y,\,y\ge0)$ as arising from an interweaving construction, as in Section \ref{sec:interweaving}, so for each $y\ge0$, $\gamma^y$ is formed as in \eqref{eq:interweaving_skewer}, by concatenating alternating intervals of the skewers of two i.i.d.\ copies $(\fX_{1},\fX_{2})$ of $\widehat\fX$ with jumps marked by \BESQ[-1] spindles. Now, consider a block $U\in\gamma^0$; this corresponds to one such spindle, marking a jump at some time $t$ in either $\fX_{1}$ or $\fX_{2}$. Suppose, for example, that this spindle appears in $\fX_{1}$ with $t\in [T_2,T_4)$, in the notation of \eqref{eq:interweaving_skewer}. Then by \cite[Theorem 3.3]{Paper1-2}, $\sD_{\gamma^0}(U) = \ell^0_{1}(t) + \ell^0_{2}(T_3)$, and if $U'\in\gamma^y$ corresponds to the same spindle, then $\sD_{\gamma^y}(U') = \ell^y_{1}(t) + \ell^y_{2}(T_3)$. Such comparisons can be made for spindles coming from any interval $[T_{j-2},T_j)$ in either $\fX_{1}$ or $\fX_{2}$. Thus, the claimed bounds follow from Lemma \ref{lem:LT_Holder} by the triangle inequality, with the $p^{\rm th}$ moment of $L_\theta$ being bounded by twice that of $\widehat{D}_\theta^y$.
\end{proof}
It remains to bound terms (iii) and (iv) in Definition \ref{def:IP:metric}, which deal with mass. Consider a sequence of $m$ distinct size-biased picks among the blocks of $\gamma^0$, and match these with the blocks arising from the same spindle at time $y$, $((U_j,U_j^\prime),1\le j\le m)$, allowing that ${\rm Leb}(U_j^\prime)$ may equal zero for some $j$ if the spindle does not survive. We can separately control
\begin{itemize}
  \item total discrepancy between matched blocks $\sum_{1\le j\le m}\big|{\rm Leb}(U_j)-{\rm Leb}(U_j^\prime)\big|$,
  \item unmatched level-0 mass $\|\gamma^0\|-\sum_{1\le j\le m}{\rm Leb}(U_j)\cf\{{\rm Leb}(U_j')> 0\}$,
  \item and unmatched level-$y$ mass $\|\gamma^y\|-\sum_{1\le j\le m}{\rm Leb}(U_j^\prime)$.
\end{itemize}

Denote by $\widetilde\mu_\lambda$ the distribution of $M\widetilde{\gamma}^0$ for independent $\widetilde{\gamma}^0\sim\widetilde\mu$ and $M\sim\GammaDist[\frac32,\lambda]$.

\begin{lemma}\label{gammaholder}
 Let $(\gamma^y,y\ge 0)$ be an $\cI^\circ$-valued type-2 evolution starting according to $\widetilde\mu_\lambda$. Let $\theta\in(0,\frac14)$ and $p>0$. Then there is a constant $C=C_{\lambda,\theta,p}$ such that
 \begin{equation}\label{eq:gammaholder}
  \bE\left[(d_\cI(\gamma^0,\gamma^y))^p\right]\le Cy^{\theta p}\qquad\mbox{for all }0\le y\le 1.
 \end{equation}
\end{lemma}
\begin{proof}
 Consider such a process $(\gamma^y,\,y\ge0)$. Its initial state is of the form $\beta^0 = (0,A)\concat (0,B)\concat G\bar\beta$, where $A,B,G$ are i.i.d.\ \GammaDist[\frac12,\lambda] random variables, independent of $\widebar\beta\sim\PDIP[\frac12,\frac12]$. Further let $\beta^*=(\beta^*_1,\beta^*_2,\dots) \in [0,1]^\infty$ denote a size-biased random ordering of the masses of $\bar\beta$. 
 
 To construct the correspondence, we take the blocks $U_1=A$, $U_2=B$, together with the blocks $U_i$ for $3\leq i\leq m=\lfloor 3y^{-1/4}\rfloor$, where $U_i$ is the block corresponding to $G\beta^*_{i-2}$ and match them with the blocks $U_1',U_2',U_3',\ldots$ that arise from the corresponding spindles at level $y$. Consequently, ${\rm Leb}(U_i') = \mathbf{g}_i(y)$ where $\mathbf{g}_i\sim\besq_{{\rm Leb}(U_i)}(-1)$ given ${\rm Leb}(U_i)$. 
%
 Note that this means that some of our blocks will be matched with empty blocks and should thus be omitted from the correspondence and accounted for in the remaining mass component of the metric.  We will handle this later.
 
 Let $p \geq 2$ and $y\in(0,1]$. Using \cite[Lemma 33]{Paper0} and the fact that $M:=A+B+G$ has finite moments of all orders, there are constants $C_1,C_2, C_3,C_4$ and $C_5$, depending only on $p$, such that
 \[ \begin{split}
  &\bE \Bigg[ \Bigg(\sum_{1\le j\le m}\!|{\rm Leb}(U_j)-{\rm Leb}(U_j^\prime)|\Bigg)^{\!\!p\,}\Bigg]
  	\leq m^{{p}-1} \bE \Bigg[ \sum_{1\le j\le m}\!\Big|{\rm Leb}(U_j)-{\rm Leb}(U_j^\prime)\Big|^{p}\Bigg]\\
	&\qquad\leq m^{{p}-1}y^{{p}/2} \sum_{1\le j\le m}\bE\Big[\Big(C_1 +\sqrt{C_2{\rm Leb}(U_j)+C_3}\Big)^{\!p\,}\Big] \\
	&\qquad\leq m^{{p}-1}y^{{p}/2} \sum_{1\le j\le m}\bE\Big[\Big(C_1 +\sqrt{C_2{\rm Leb}(M)+C_3}\Big)^{\!p\,}\Big] \\
	&\qquad\leq C_4 m^{{p}} y^{{p}/2} 
   \leq C_5 y^{p/4},
 \end{split} \]
 using in the last step that $m\le 3y^{-1/4}$ and absorbing $3^p$ into the constant.
 
 The unmatched mass at level $0$ is $G\sum_{j=m-1}^\infty \beta^*_i$.  Let $Y_n\sim{\tt Beta}\big(\frac12, (n+1)/2\big)$, $n\ge 1$, be a sequence of independent random variables, also independent of $G$.  Using the stick-breaking construction of the Poisson--Dirichlet distribution, we see that for some $C_6\ge C_5$ and all $y\in(0,1]$
 \[\begin{split}
  &\bE\Bigg[ \Bigg(G \!\sum_{j=m-1}^\infty\! \beta^*_i\Bigg)^{\!p\,}\Bigg]
 	= \bE\big[G^p\big] \bE\Bigg[ \Bigg(1 - \sum_{j=1}^{m-2} \beta^*_i\Bigg)^{\!p\,}\Bigg]
	= \bE\big[G^p\big]\bE\Bigg[ \Bigg(\prod_{j=1}^{m-2} (1\!-\!Y_j)\Bigg)^{\!p\,}\Bigg] \\
	&\qquad = \bE\big[G^p\big]\prod_{j=1}^{m-2}\bE\big[ \left(1-Y_j\right)^{p}\big] 
	=\bE\big[G^p\big]\Gamma\left(1+{p} \right)\frac{ \Gamma\left(\frac{m}{2}\right)}{\Gamma\left(\frac{m}{2} +{p}\right)} \\
	&\qquad\sim \bE\big[G^p\big]\Gamma\left(1+{p} \right)m^{-{p}} 
	= \bE\big[G^p\big]\Gamma(1+{p}) \lfloor 3y^{-1/4}\rfloor^{-{p}}
 	\le C_6 y^{p/4}.
 \end{split}\]

 It remains to estimate the unmatched mass at time $y$. By the triangle inequality,
 \[ \|\gamma^y\| - \sum_{j=1}^{m} {\rm Leb}(U'_j)  \leq \big| \|\gamma^y\|-M\big| + \Bigg|G\sum_{j=m-1}^\infty \beta^*_j\Bigg|+ \sum_{1\le j\le m}\big|{\rm Leb}(U_j)-{\rm Leb}(U_j^\prime)\big|.\]
 Furthermore, by Theorem \ref{thm:total_mass}, $(\|\gamma^y\|)_{y\geq 0}$ is a \besq$_{M}(-1)$ process to which \cite[Lemma 33]{Paper0} applies, as above. Consequently, we have for some $C_7>0$ that
 \[\begin{split}
  \bE \Bigg[\Bigg|\|\gamma^y\| - \sum_{j=1}^{m} {\rm Leb}(U'_j)\Bigg|^{p\,}\Bigg] 
   &\leq 3^{{p}-1}\!\left(\!\!
   		\begin{array}{l}
   			\displaystyle \bE\Big[ \big| \|\gamma^y\|-M\big|^{p}\Big] + \bE\Bigg[\Bigg|G\sum_{j=m-1}^{\infty} \beta^*_j\Bigg|^{p\,}\Bigg]\\[4pt]
   			\displaystyle\ \ +\ \bE\Bigg[\Bigg(\sum_{1\le j\le m}|{\rm Leb}(U_j)-{\rm Leb}(U_j^\prime)|\Bigg)^{{\!p\,}}\Bigg]
   		\end{array}\!\!\right)\\
   & \leq C_7(y^{{p}/2} +y^{p/4} + y^{p/4})
   \leq 3C_7 y^{{p} /4}.
 \end{split}\]
 To account for the fact that some ${\rm Leb}(U'_i)$ may be $0$, and thus the corresponding $U_i$ should count towards unmatched mass at time $0$, we bound the metric $d_{\mathcal{I}}$ using the correspondence defined above (and bounding the maximum in the definition of $d_{\mathcal{I}}$ by a sum) to see that $d_{\mathcal{I}}(\gamma^0,\gamma^y)$ is bounded by
 \[\begin{split}
   &\sum_{1\le j\le m}|{\rm Leb}(U_j)-{\rm Leb}(U_j^\prime)|\mathbf{1}\{{\rm Leb}(U'_j) >0\}\\[-.2cm]
   &\quad + \|\gamma^y\| - \sum_{j=1}^{m} {\rm Leb}(U'_j)\mathbf{1}\{{\rm Leb}(U'_j)\! >\!0\}
    + M - \sum_{j=1}^m {\rm Leb}(U_j)\mathbf{1}\{{\rm Leb}(U'_j)\! >\!0\} \\
   &\quad + |\ell^0(T)-\ell^y(T)|+  \max_{1\le j\le m}|\ell^0(U_j)-\ell^y(U_j^\prime)|\mathbf{1}\{{\rm Leb}(U'_j) >0\} \\
   &=   \sum_{1\le j\le m}|{\rm Leb}(U_j)-{\rm Leb}(U_j^\prime)|
    + \|\gamma^y\| - \sum_{j=1}^{m} {\rm Leb}(U'_j)
    + M -\sum_{j=1}^m {\rm Leb}(U_j)\\
   &\quad + |\ell^0(T)-\ell^y(T)|+  \max_{1\le j\le m}|\ell^0(U_j)-\ell^y(U_j^\prime)|\mathbf{1}\{{\rm Leb}(U'_j) >0\}.
 \end{split}\]
 Dropping the indicator on the last term and combining this with our calculations above and Lemma \ref{lem:LT_matching_bd} shows that for $0<\theta<\frac14$ there exists some constant $C_{\lambda,\theta,{p}}$ depending only on $\lambda$, $\theta$ and $p$ that satisfies \eqref{eq:gammaholder}.
\end{proof}

\begin{proof}[Proof of Proposition \ref{type2holder}]
Let $(\gamma^y,y\ge 0)$ be an $\cI^\circ$-valued type-2 evolution with initial distribution $\widetilde\mu_\lambda$ and degeneration time $D$. Denote the total mass evolution by 
$Z(y)=\|\gamma^y\|$, $y\ge 0$. Let $0\le a<b\le 1$. Then
\begin{align*}
  \bE_{\widetilde\mu_\lambda}\!\left[(d_\cI(\gamma^a,\gamma^b))^p\right]
  &=\bE_{\widetilde{\mu}_\lambda}\!\left[\cf\{D>a\}\bE_{\gamma^a}\!\left[(d_\cI(\gamma^0,\gamma^{b-a}))^p\right]\right]\\
    &\qquad +\bE_{\widetilde\mu_\lambda}\!\big[|Z(a)-Z(b)|^p\cf\{D<a\}\big].
\end{align*}
For the first term, we condition on $D>a$ and apply the pseudo-stationarity of Proposition \ref{prop:type2:pseudo}. While the distribution of $\gamma^a$ given $D>a$ may not be $\widetilde\mu_{\lambda/(2\lambda a+1)}$, it is $\widetilde\mu_{\lambda/(2\lambda a+1)}$ up to a potential swap of the two leftmost blocks, and the matching set up in the proof of Lemma \ref{gammaholder} is unaffected by such a swap so that scaling by $2\lambda a+1$ and applying the bound of Lemma \ref{gammaholder} yields the upper bound:
\begin{align*}
  \bE_{\widetilde\mu_{\lambda/(2\lambda a+1)}}\!\left[\big(d_\cI\big(\gamma^0,\gamma^{b-a}\big)\big)^{p\,}\right]
  &=(2\lambda a+1)^p\bE_{\widetilde\mu_\lambda}\!\left[\big(d_\cI\big(\gamma^0,\gamma^{(b-a)/(2\lambda a+1)}\big)\big)^p\right]\\
  &\le (2\lambda a+1)^pC(b-a)^{\theta p}/(2\lambda a+1)^{\theta p}\\
  &\le (2\lambda+1)^{p(1-\theta)}C|b-a|^{\theta p}.
\end{align*}
For the second term, we apply \cite[Lemma 33]{Paper0} to find the upper bound
$$|b-a|^{p/2}\bE_{\widetilde\mu_\lambda}\!\left[\left(1+2(p-1)+2\sqrt{p-1}\sqrt{Z(0)+2(p-1)}\right)^{p\,}\right],$$
which is easily seen to be a finite multiple of $|b-a|^{p/2}\le|b-a|^{\theta p}$.

By the Kolmogorov--Chentsov theorem \cite[Theorem I.(2.1)]{RevuzYor}, this shows that for all $0<\theta<\frac14$ and $p>0$,
$$\bE_{\widetilde\mu_\lambda}\!\left[\left(\sup_{0\le a<b\le 1}\frac{d_\cI(\gamma^a,\gamma^b)}{|b-a|^\theta}\right)^{\!p\,}\right]<\infty.$$
We can write the left-hand side by integrating out the random initial mass. Canceling $\lambda^{3/2}/\Gamma(\frac32)$ gives
$$\int_0^\infty e^{-\lambda x}\sqrt{x}\bE_x\!\left[\left(\sup_{0\le a<b\le 1}\frac{d_\cI(\gamma^a,\gamma^b)}{|b-a|^\theta}\right)^{\!p\,}\right]dx<\infty.$$
By Fubini's theorem, this yields for a.e.\ $x\in(0,\infty)$ that
$$\bE_x\!\left[\left(\sup_{0\le a<b\le 1}\frac{d_\cI(\gamma^a,\gamma^b)}{|b-a|^\theta}\right)^{\!p\,}\right] < \infty.$$
But for any $x,y\in(0,\infty)$, we can find $c<1/y$ so that this expectation is finite for initial mass $cx$. By scaling, 
\begin{align*}
 \infty&>\bE_{c x}\!\left[\sup_{0\le a<b\le 1}\left(\frac{d_\cI(\gamma^a,\gamma^b)}{|b-a|^\theta}\right)^{\!p\,}\right]
       = c^p\bE_{x}\!\left[\sup_{0\le a<b\le 1}\left(\frac{d_\cI(\gamma^{a/c},\gamma^{b/c})}{|b-a|^\theta}\right)^{\!p\,}\right]\\
       &=c^{p(1-\theta)}\bE_{x}\!\left[\sup_{0\le a^\prime<b^\prime\le 1/c}\left(\frac{d_\cI(\gamma^{a^\prime},\gamma^{b^\prime})}{|b^\prime-a^\prime|^\theta}\right)^{\!p\,}\right]\\
       &\ge c^{p(1-\theta)}\bE_{x}\!\left[\sup_{0\le a^\prime<b^\prime\le y}\left(\frac{d_\cI(\gamma^{a^\prime},\gamma^{b^\prime})}{|b^\prime-a^\prime|^\theta}\right)^{\!p\,}\right],
\end{align*}
so the expectation is finite for any initial mass, including unit initial mass $x=1$, and for any $y\in(0,\infty)$.       
\end{proof}

We conclude this section by returning to the $\cM^\circ$-valued processes that capture projected metric tree structure and projected mass measures, in preparation for a
$d_{\rm GHP}$-continuous evolution of continuum random trees as claimed in Theorem \ref{thm:intro:AD}.

\begin{corollary} Consider any pseudo-stationary $\cJ^\circ$-valued type-2 evolution $(\Gamma^y,y\ge 0)$. Then the associated $\cM^\circ$-valued evolution $(M_2(\Gamma^y),y\ge 0)$
  is almost surely $\theta$-H\"older continuous in $(\cM^\circ,d_{\rm HP})$ for all $\theta\in(0,\frac14)$. 
\end{corollary}
\begin{proof} This follows from Proposition \ref{type2holder}, 1-self-similar scaling by an independent initial mass and the Lipschitz property of $M_2\colon\cJ^\circ\rightarrow\cM^\circ$ that we noted above \eqref{eq:onebranch}.
\end{proof}

\chapter{Self-similar and unit-mass $k$-tree evolutions}
\label{ch:constr}

In this chapter, we generalize the type-2 evolutions of Chapter \ref{ch:type-2} to several variants of $k$-tree evolutions with fluctuating total mass processes, and we generalize the results about (pseudo-)stationarity and unit-mass 2-tree evolutions of Chapter \ref{dePoiss} to establish stationary unit-mass (resampling) $k$-tree evolutions. To do this, we pull together several threads, which we recall in an informal way here, leaving precise statements to the later sections in this chapter. 

\begin{figure}[b!]\centering
  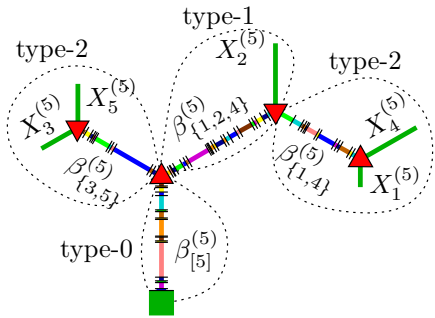
  \caption{$(\ft_5,(X_1^{(5)},\ldots,X_5^{(5)}),(\beta_{[5]}^{(5)},\beta_{\{1,2,4\}}^{(5)},\beta_{\{1,4\}}^{(5)},\beta_{\{3,5\}}^{(5)}))$ is a 5-tree with two type-2, one type-1 and one type-0 compounds.\label{fig:B_k_tree_proj_cons}}
\end{figure}

Informally, a $k$-tree for us is a tree of the sort shown in Figure \ref{fig:B_k_tree_proj_cons}, consisting of a binary combinatorial tree shape with $k$ labeled leaves, top masses for the leaf edges and interval partitions for the other (internal) edges. We further associate with each of the internal edges the number of leaf edges incident to their upper vertex (the vertex further away from the root) hence forming compounds of three types, where type $i=0,1,2$ has $i$ top masses.  

The fundamental idea is to use an independent type-$i$ evolution for each compound of type $i=0,1,2$. As noted in Chapters \ref{ch:prelim}--\ref{ch:type-2}, type-1 and type-2 evolutions degenerate in finite time. As a consequence, such $k$-tree evolutions also degenerate. This gives rise to a first kind of (self-similar) $k$-tree evolution that is killed when one of the constituent type-1 or type-2 evolutions reaches its degeneration time. 

It is natural to view degeneration as the loss of a label. In order to prepare for the consistency results of Chapter \ref{ch:consistency}, indeed of Theorem \ref{thm:intro:k_tree}\ref{main cons}, we may make a swap of two labels before reducing the tree shape, applying in the framework of $k$-trees with real-valued top masses and $\cI$-valued edge partitions the rules of \cite{Paper2} developed in a setting of combinatorial trees, as recalled in Section \ref{sec:discrete_to_cts}. In a \em non-resampling $k$-tree evolution\em, the evolution resumes as a succession of killed $j$-tree evolutions for $j=k-1,k-2,\ldots,2$ until becoming a final type-2 evolution, which eventually has one of the two remaining top masses continue as a $\besq(-1)$ evolution between the degeneration time and the lifetime of this final type-2 evolution. Here, this remaining $\besq(-1)$ process beyond the degeneration time of this type-2 evolution can be viewed as a 1-tree evolution.  

The next thread is \em resampling\em. In Definition \ref{def:resampling}, we resampled de-Poissonized type-2 evolutions to obtain unit-mass 2-tree evolutions. Here, we enhance the notion of resampling to handle more complex states at degeneration. Indeed, we develop this in the self-similar setting to obtain \em resampling $k$-tree evolutions\em. 

Obtaining unit-mass $k$-tree evolutions by \em de-Poissonization \em is then straightforward following the same steps and arguments as in Section \ref{sec:type2:dePoiss}, but 
\em (pseudo-) stationarity \em of $k$-tree evolutions requires refined arguments to handle the enhanced resampling. Specifically, we have obtained type-0/1/2 pseudo-stationarity results that were conditional on non-degeneration, and we will here establish $k$-tree pseudo-stationarity results that are unconditional.   

The structure of this chapter is as follows. In Section \ref{sec:type012}, we collect some combined results for type-$i$ evolutions, $i=0,1,2$, for ease of reference, and we record some further consequences. In Section \ref{sec:killed_def}, we introduce the spaces of $k$-trees, in which our $k$-tree evolutions take their values. In Section \ref{sec:resamp_def}, we introduce killed, non-resampling and resampling $k$-tree evolutions in the self-similar regime and identify the total mass processes of the latter two as $\besq(-1)$ processes (up to the accumulation time of resampling times, but we show in Chapter \ref{ch:consistency} that the total mass reaches zero continuously at this accumulation time). In Section \ref{sec:pseudo}, we establish pseudo-stationarity properties of $k$-tree evolutions. In Section \ref{sec:ktree:dePoiss} we de-Poissonize to obtain unit-mass processes.

\section{Summary of properties of type-0, type-1, and type-2 evolutions}\label{sec:type012}

Recall from Chapters \ref{ch:prelim}--\ref{ch:type-2} that type-0, type-1, and type-2 evolutions are Markov processes introduced with pathwise constructions. In Section \ref{sec:discrete_to_cts}, we argued via a connection to ordered Chinese restaurant processes that the type-2 evolution is a continuum analogue of a certain 2-tree projection of the discrete Aldous chain discussed in \cite[Appendix A]{Paper2}. 
By the same argument, the $k$-tree projection of the Aldous chain also discussed there can be decomposed into parts whose evolutions are analogous to type-0/1/2 evolutions. In Figure \ref{fig:B_k_tree_proj_cons}, the dashed lines separate parts of the $k$-tree that evolve as type-0/1/2 evolutions; see Definition \ref{def:killed_ktree}. 

The aforementioned pathwise construction brings a lot of symmetry to light, and it makes many calculations accessible. In this chapter, we will not delve into this construction, and in fact, only a few key properties of these processes are needed. We refer to Section \ref{sec:prelim:type01_def} and Definition  \ref{def:type2:v1} for the definitions of type-0, type-1 and type-2 evolutions. For ease of reference, and in order to better exhibit some patterns that emerge, we re-group and re-state some results from Chapters \ref{ch:prelim}--\ref{ch:type-2} and record some elementary consequences in this section.

Type-$i$ evolutions, for $i=0,1,2$, are valued in (subsets of) the product space $[0,\infty)^i\times \cI$. We refer to the real-valued first coordinates of type-1 and type-2 evolutions as \emph{top blocks} or \emph{top masses}. Each block in a type-$i$ evolution, including these top blocks, has mass that fluctuates as a squared Bessel diffusion \BESQ[-1]. 
Informally, when a top block of a type-1 or type-2 evolution hits mass zero, the leftmost blocks of the interval partition component of the evolution are successively (informally speaking, as the blocks are not well-ordered) pulled out of the interval partition to serve as new top blocks, until their \BESQ[-1] masses are absorbed at zero.

\begin{proposition}[cf.\ Proposition \ref{type1totalmass} and Theorem \ref{thm:type2:total_mass}]\label{prop:012:mass}
 For $i=0,1,2$, the total mass process for a type-$i$ evolution is a \BESQ[1-i]. Moreover, this total mass process is a strong Markov process in the filtration of the type-$i$ evolution.
\end{proposition}

The strong Markov property in the larger filtration noted above follows from Dynkin's criterion; see Theorem \ref{thm:Dynkin}.

Recall that 0 is instantaneous reflecting for \BESQ[1], and a type-0 evolution is similarly reflecting at $\emptyset$ and it will be useful to say that it has 
an infinite degeneration time. On the other hand, we consider 0 as absorbing for \BESQ[0] and \BESQ[-1]. A type-1 evolution is said to \emph{degenerate} when it is absorbed at $(0,\emptyset)$ at a time that we also refer to as its \em lifetime\em. While a type-2 evolution $((m_1^y,m_2^y,\beta^y),y\ge 0)$ is 
eventually absorbed in $(0,0,\emptyset)$, and we refer to this time $\zeta$ as the \emph{lifetime} of the type-2 evolution, a first degeneration happens at the earlier time $D<\zeta$ at which either $m_1^y+\|\beta^y\|$ or $m_2^y+\|\beta^y\|$ hits zero and is absorbed. This is the \emph{degeneration time} of the type-2 evolution. See Corollary \ref{type2:degenlife}. 

\begin{proposition}[cf.\ Propositions \ref{prop:type01:diffusion}, \ref{prop:type2:IPvalued}, Corollary \ref{corpairtype1}, Theorem \ref{thm:diffusion}]\label{prop:012:pred}$\,$
 \begin{enumerate}[label=(\roman*),ref=(\roman*)]
  \item Type-0/1/2 evolutions and $\cI$-valued type-1 and $\cI^\circ$-valued type-2 evolutions are self-similar Borel right Markov processes.\label{item:pred:Markov}
  \item Type-0 evolutions and $\cI$-valued type-1 evolutions are path-continuous. \label{item:pred:0}
  \item If $\big(\big(m_1^y,m_2^y,\beta^y\big),y\!\ge\!0\big)$ is a type-2 evolution and $I(y) :\equiv \max\{n\!\ge\! 0\colon Y_n\!\le\! y\}\linebreak\text{mod }2$ is $\{1,2\}$-valued, where $(Y_n,n\ge0)$ is as in Definition \ref{def:type2:v1}, then the \emph{$\cI^\circ$-valued type-2 evolution} $\big(\big(0,m_{3-I(y)}^y\big)\concat \big(0,m_{I(y)}^y\big)\concat\beta^y,y\ge0\big)$ is a diffusion. Also, 
  	each of $m_1^y$ and $m_2^y$ can only equal zero when $\beta^y$ has no leftmost block, and they can only both equal zero if $\beta^y = \emptyset$.\label{item:pred:2}
 \end{enumerate}
\end{proposition}

\begin{proposition}[Concatenation properties; Proposition \ref{prop:type1:construction}, Corollary \ref{type01plustype1}]\label{prop:012:concat}
 Consider a type-1 evolution $((m^y,\beta^y),y\ge0)$.
 \begin{enumerate}[label=(\roman*),ref=(\roman*)]
  \item \label{item:012concat:BESQ+0}
  	Let $\zeta$ denote the first time that $m^y$ hits zero. Then $(m^y\cf\{y\le\zeta\},y\ge0)$ is a \BESQ[-1] and $(\beta^y,y\in [0,\zeta])$ distributed as an independent type-0 evolution stopped at $\zeta$.
  \item \label{item:012concat:0+1}
  	If $(\widetilde\beta^y,y\ge0)$ is an independent type-0 evolution, then $(\widetilde\beta^y\concat(0,m^y)\concat\beta^y,y\ge0)$ is a type-0 evolution.
  \item \label{item:012concat:1+1}
  	Suppose instead that $((\widetilde m^y,\widetilde\beta^y),y\ge0)$ is an independent type-1 evolution and let $\widetilde D$ denote its degeneration time. Then the following process is a type-1 evolution:
  \begin{equation}\label{eq:012concat:1+1}
   \left\{\begin{array}{ll}
    (\widetilde m^y,\widetilde\beta^y\concat(0,m^y)\concat\beta^y)	& \text{for }y\in [0,\widetilde D),\\
    (m^y,\beta^y)	& \text{for }y\ge \widetilde D.
   \end{array}\right.
  \end{equation}
  \item \label{item:012concat:2+1}
  	Suppose instead that $((\widetilde m_1^y,\widetilde m_2^y,\widetilde\beta^y),y\ge0)$ is an independent type-2 evolution. Let $\widetilde D$ denote its degeneration time. Let $(\widehat x_1,\widehat x_2)$ equal $(\widetilde m_1^{\widetilde D},m^{\widetilde D})$ if $\widetilde{m}_2^{\widetilde D} = 0$ (i.e.\ if label 2 is the label that degenerates at time $\widetilde{D}$), or equal $(m^{\widetilde D},\widetilde m_2^{\widetilde D})$ otherwise (if label 1 degenerates). Let $((\widehat m_1^y,\widehat m_2^y,\widehat\beta^y),y\ge0)$ be a type-2 evolution with initial state $(\widehat x_1,\widehat x_2, \beta^{\widetilde D})$, conditionally independent of the other processes given its initial state. The following is a type-2 evolution:
  \begin{equation}\label{eq:012concat:2+1}
   \left\{\begin{array}{ll}
    (\widetilde m_1^y,\widetilde m_2^y,\widetilde\beta^y\concat(0,m^y)\concat\beta^y)	& \text{for }y\in [0,\widetilde D),\\
    (\widehat m_1^{y-\widetilde D},\widehat m_2^{y-\widetilde D},\widehat\beta^{y-\widetilde D})	& \text{for }y\ge \widetilde D.
   \end{array}\right.
  \end{equation}
 \end{enumerate}
 Moreover, the concatenated evolutions constructed in \ref{item:012concat:0+1}, \ref{item:012concat:1+1}, and \ref{item:012concat:2+1} each possess the strong Markov property in the larger filtrations generated by their constituent parts.
\end{proposition}

\begin{corollary}\label{cor:decomp_ker}
 Consider an initial state $(m_1,m_2,\beta)\in \cJ^\circ$ with $\beta\neq\emptyset$ and a distinguished block $(a,b)\in\beta$, and suppose $(\Gamma^y,\,y\ge0)$ is a type-2 evolution with this initial state. 
 Then there exist, possibly on an enlarged probability space, independent type-2 and type-1 evolutions $(\widetilde m_1^y,\widetilde m_2^y,\widetilde\beta^y)$ and $(m^y,\beta^y)$, $y\ge0$, such that their concatenation in the sense of \eqref{eq:012concat:2+1} equals $(\Gamma^y)$, up until the degeneration time $\widetilde D$ of $(\widetilde m_1^y,\widetilde m_2^y,\widetilde\beta^y)$, with $(a,b)\in\beta$ corresponding to the block $(0,m^0)$ in the initial concatenated process.
 
 There exists a regular conditional distribution $\kappa$ for the joint law of $(\widetilde m_1^y,\widetilde m_2^y,\widetilde\beta^y)$ and $(m^y,\beta^y)$, $y\ge0$, given $(\Gamma^y,\,y\ge0)$ and the block $(a,b)$.
 
 Corresponding claims hold relating to assertions \ref{item:012concat:0+1} and \ref{item:012concat:1+1} of Proposition \ref{prop:012:concat}.
\end{corollary}

\begin{proof}
 The claimed existence of regular conditional distributions follows from the properties that $(\cI,d_{\cI})$ is Lusin (Proposition \ref{prop:IPspace:Lusin}) and that these evolutions have c\`adl\`ag paths. The remaining assertions are immediate consequences of Proposition \ref{prop:012:concat}.
\end{proof}

Recall from Proposition \ref{prop:PDIP} that a Poisson--Dirichlet interval partition with parameters $\big(\frac12,\frac12\big)$, called $\PDIP\big(\frac12,\frac12\big)$, is an interval partition whose ranked block sizes have law $\PoiDir[\frac12,\frac12]$, with the blocks exchangeably ordered from left to right. Let $A\sim \BetaDist[\frac12,\frac12]$, $(A_1,A_2,A_3)\sim\distribfont{Dirichlet}\big(\frac12,\frac12,\frac12\big)$, and $\bar\beta\sim\PDIP\big(\frac12,\frac12\big)$ independent of each other. Recall that we refer to a probability distribution on $\cI$ as a \emph{pseudo-stationary law for the type-0 evolution} if it is the law of $M\bar\beta$, i.e.\ $\bar\beta$ scaled by $M$, for an independent random mass $M>0$. Likewise a law on $[0,\infty)\times\cI$,
respectively $[0,\infty)^2\times\cI$, is a \emph{pseudo-stationary law for the type-1, {\rm resp.} type-2, evolution} if it is the law of any independent multiple of $(A,(1-A)\bar\beta)$, resp.\ $(A_1,A_2,A_3\bar\beta)$. This language is in reference to the following proposition.

\begin{proposition}[Propositions \ref{prop:01:pseudo} and \ref{prop:type2:pseudo}--\ref{prop:type2:pseudo_g}]\label{prop:012:pseudo}
 For $i=0,1,2$, if a type-$i$ evolution has a pseudo-stationary initial distribution, then given that it does not degenerate prior to time $y$, its conditional law at time $y$ is also pseudo-stationary. In the special case that its initial mass has law $\GammaDist[\frac{1+i}{2},\lambda]$, then its mass at time $y$ has conditional law $\GammaDist[\frac{1+i}{2},\lambda/(2\lambda y+1)]$.
\end{proposition}

\section{State spaces of $k$-trees}\label{sec:killed_def}

In this section we formalize the notion of a $k$-tree introduced in Section \ref{sec:intro:B_k_tree}, which we write as a tree shape equipped with top masses and edge partitions: 
$$T_k = \left(\ft_k,(x_j^{(k)}\!,j\!\in\![k]),(\beta_E^{(k)}\!,E\!\in\!\ft_k)\right).$$

\begin{figure}
 \centering
 \includegraphics[scale=0.4]{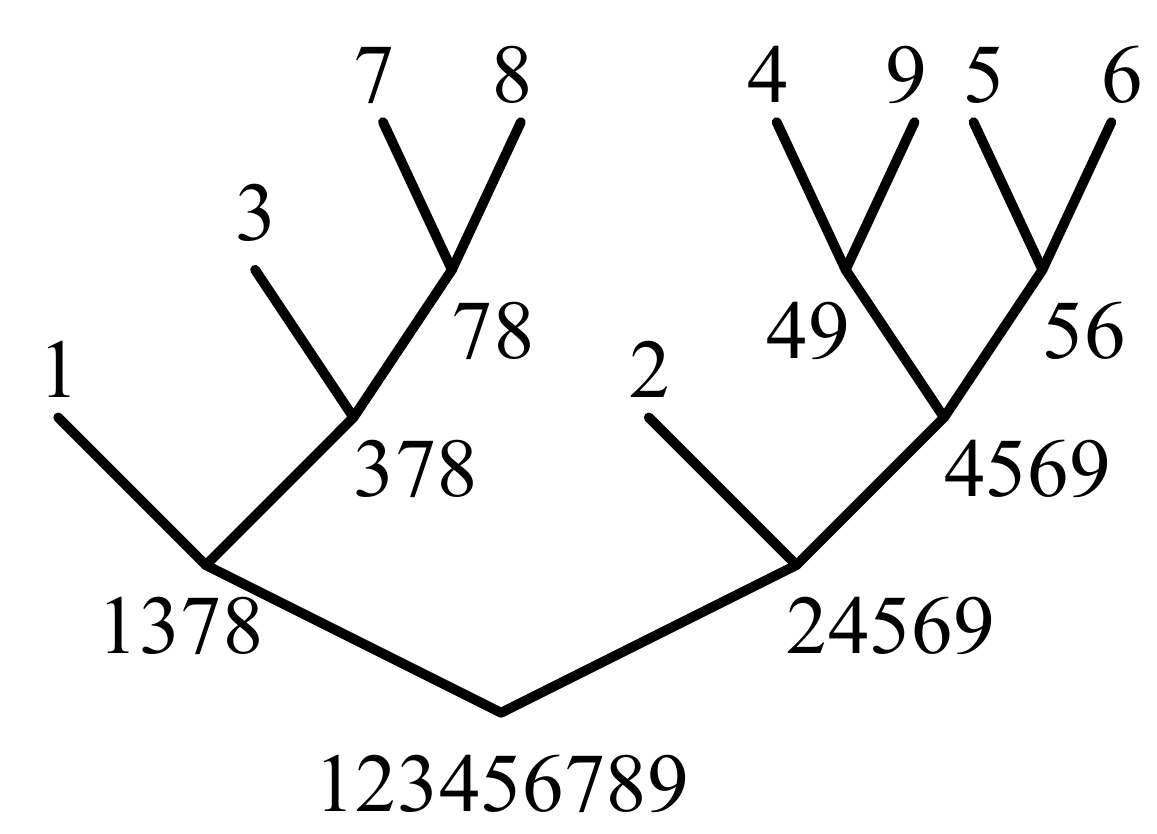}\quad\includegraphics[scale=0.4]{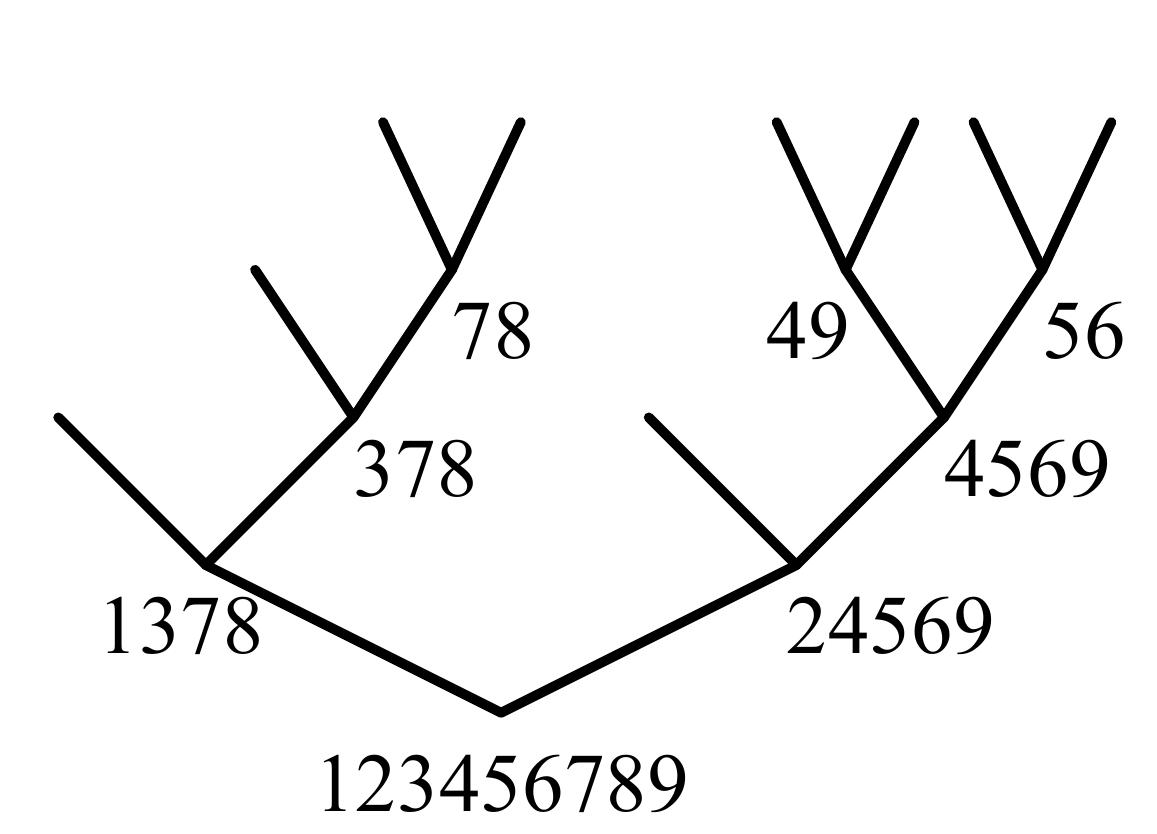}
 \caption{Left: A tree-representation of a binary hierarchy on $[9]$. Right: The same tree with leaf labels omitted. Tree shape notation: $\ft =\{\{1,2,3,4,5,6,7,8,9\},$ $\{1,3,7,8\},$ $\{2,4,5,6,9\},$ $\{3,7,8\},$ $\{4,5,6,9\},$ $\{7,8\},$ $\{4,9\},$ $\{5,6\}\}$.\label{fig:fragmentation}}
\end{figure}

Our notion of a tree shape is a variation of similar notions that capture a combinatorial tree structure within a richer model; see e.g.\ Pitman \cite[Chapter 7]{CSP}. Indeed, our tree shape $\ft_k$ is equivalent to a leaf-labeled combinatorial tree also known as a \em cladogram \em \cite{Aldous96,For-05}, \em fragmentation \em \cite{Bertoin06,HMPW,MPW}, \em hierarchy \em \cite{FlajSedg09} or \em total partition \em \cite{StanleyVol2} of $[k]$. 
We visualize tree shapes as rooted binary trees, but we formalize them as sets of subsets of a leaf set $A$, rather than as graphs $G=(V,E)$. Before developing this formally, see Figure \ref{fig:fragmentation} for an example of a binary tree and its tree shape. Below, we list all tree shapes with leaf set $[3]=\{1,2,3\}$,
$$\bT_{[3]}^{\rm shape}\!\!=\!\{\{[3],\{2,3\}\},\,\{[3],\{1,3\}\},\,\{[3],\{1,2\}\}\},\ \mbox{formalizing }\,
	\parbox{0.9cm}{\includegraphics[height=1cm]{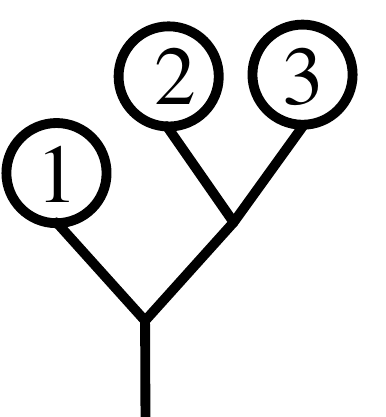}}\!,\,%
    \parbox{0.9cm}{\includegraphics[height=1cm]{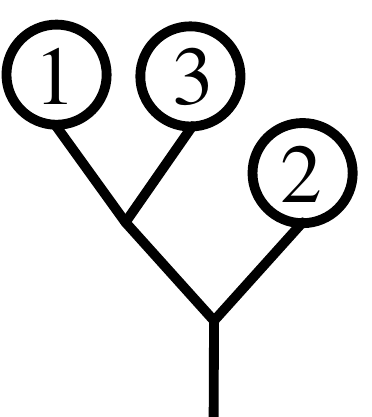}},%
    \parbox{0.9cm}{\includegraphics[height=1cm]{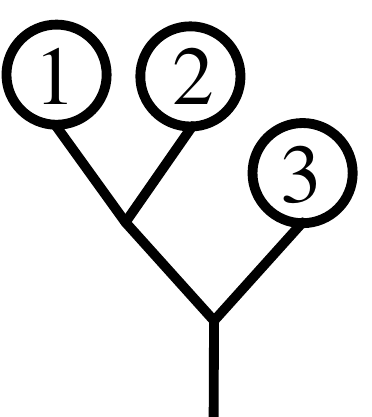}},$$
which we have illustrated as \em planted trees \em \cite[Chapter 7]{CSP} that have both a base vertex and a degree-1 root vertex connected to the base vertex. This allows us to associate the label of each non-root vertex also with the edge below the vertex. Indeed, we think of top masses and edge partitions as being associated with leaf edges and internal edges (including the root edge), respectively. 

Formally, a subset $\overline{\ft}$ of the power set of $A\subset\bN$ is a \em binary hierarchy on A \em if 
\begin{itemize}
  \item $A\in\overline{\ft}$ and $\emptyset \notin \overline{\ft}$, 
  \item each $B\!\in\!\overline{\ft}$ with $\#B\!\ge\! 2$ has a non-trivial partition into $C_1,C_2\!\in\!\overline{\ft}$, 
  \item for all $B,C\in\overline{\ft}$, we have $B\cap C=\emptyset$ or $B\subseteq C$ or $C\subseteq B$,
  \item $\{j\}\in\overline{\ft}$ for all $j\in A$.   
\end{itemize}
Note that in the setting of (ii), the pair $\{C_1,C_2\}$ is unique. 
We call $B$ the \em parent \em of $C_1$ and $C_2$, writing $\parent{C_1}=\parent{C_2}=B$. We call $C_1$ and $C_2$ \em siblings\em, writing
$C_1={\rm sib}(C_2)$ and $C_2={\rm sib}(C_1)$. The sibling of the parent of a set is its \emph{uncle}. 
We denote the set of binary hierarchies on $A$ by $\bT_A^{\rm hierarchy}$.

We can associate with any binary hierarchy $\overline{\ft}\in\bT_A^{\rm hierarchy}$ a graph-theoretic rooted binary tree 
$({\rm vert}(\overline{\ft}),{\rm edge}(\overline{\ft}))$ with vertex and edge sets
$${\rm vert}(\overline{\ft})=\overline{\ft}\cup\{\textsc{root}\}\quad\mbox{and}\quad {\rm edge}(\overline{\ft})=\{\{B,C\}\subset {\rm vert}(\overline{\ft})\colon B=\parent{C}\},$$
with the convention that $\parent{A}=\textsc{root}$. Here, \em binary \em means that apart from the leaves and the root of degree 1, all vertices have degree 3, i.e.\ are binary branch points. 
We denote the set of graph-theoretic rooted binary trees with leaves labeled by $A$ as $\bT_A^{\rm graph}$.

Now, consider the injective map that sends $\overline{\ft}\in\bT_A^{\rm hierarchy}$ to $\ft = \overline{\ft}\setminus\{\{j\}\colon j\in A\}$. We define the image of $\bT_A^{\rm hierarchy}$ under this map to be $\bT^{\rm shape}_A$, the set of \emph{tree shapes} with leaf labels in $A$. The elements of a tree shape $\ft$ correspond to branch points in ${\rm vert}(\overline{\ft})$. We also take these elements to represent the \emph{parent edges} of the branch points, i.e.\ we use label $C\in\ft$ to refer to the \em internal edge \em $\{\parent{C},C\}\in {\rm edge}(\overline{\ft})$.

We emphasize that $\bT_A^{\rm hierarchy}$, $\bT_A^{\rm graph}$ and $\bT_A^{\rm shape}$ are all in natural one-to-one correspondence. We will mainly use the notation $\ft\in\bT^{\rm shape}_A$, but we regard $\ft$ as representing the corresponding binary tree, including its leaves. For example, if $\ft = \{[3],\{1,2\}\}$, we say that the internal edge $\{1,2\}$ is the sibling of the leaf edge $\{3\}$ in (the graph-theoretic tree represented by) $\ft$, even though $\{3\}\notin\ft$. 
We refer to the members of $\ft$ as edges and use graph-theoretic descriptions of operations on $({\rm vert}(\overline{\ft}),{\rm edge}(\overline{\ft}))\in\bT^{\rm graph}_A$ that induce operations on tree shapes. 

Indeed, by removing $\{j\}$, $j\in A$, from $\overline{\ft}$, this notion of a tree shape $\ft\in\bT_A^{\rm shape}$ offers simple notation for both leaf edges (members of $A$) and internal edges (members of $\ft$). Our aim is to study evolutions of interval partitions $\beta_E\in\cI$ associated with each $E\in\ft$ together with those top masses $x_j\in[0,\infty)$ associated with label $j\in A$ for which $\parent{\{j\}}=E$. In our binary setting, this gives rise to \em type-$i$ edges \em with $i$ top masses, $i=0,1,2$. 
E.g., in Figure \ref{fig:B_k_tree_proj_cons}, edges $\{3,5\}$ and $\{1,4\}$ are type-2 edges, label 2 and its top mass $X_2^{(5)}$ are associated with the type-1 edge $\{1,2,4\}$, while edge $[5]$ is a type-0 edge. 

More precisely, an edge $E\in\ft$ with $\#E=2$ is called a \em type-2 edge\em. If $E=\{j,j^\prime\}$, then $x_j$ and $x_{j^\prime}$ are the two top masses associated with $E$. If $j$ is not associated with a type-2 edge, it is associated with the \em type-1 edge \em $E:=\parent{\{j\}}\in\ft$, which satisfies $\#E\ge 3$ since $\{j\}$ has a sibling ($E\setminus\{j\}$) in $\ft$. Edges that are neither type-1 nor type-2 edges are called \em type-0 edges\em. They have no top masses. Instead, such edges $E$ can be written as $E_1\cup E_2$ for two edges $E_1,E_2\in\ft$ with disjoint label sets: $E_1\cap E_2=\varnothing$. See Figure \ref{fig:B_k_tree_proj_cons} for labeled examples of type-0, type-1, and type-2 edges.

For a finite, non-empty set $A\subset\bN$, an \emph{$A$-tree} is a tree shape $\ft\in\bT_{A}^{\rm shape}$ equipped with non-negative weights on leaf edges and interval partitions marking the internal edges:
\begin{equation}\label{eq:k_tree_def}
 \bTInt_{A}=\bigcup_{\ft\in\bT^{\rm shape}_{A}}\{\ft\}\times[0,\infty)^A\times\cI^\ft.
\end{equation}
For $k\!\ge\!1$, we call elements of $\bTInt_{k} \!:=\! \bTInt_{[k]}$ \emph{$k$-trees}. For $T\!=\!(\ft,(x_j,j\!\in\! A),(\beta_E,E\!\in\!\ft))$ \linebreak $\in \bTInt_{A}$, we write $\|T\|=\sum_{j\in A}x_j+\sum_{E\in\ft}\|\beta_E\|$ for its total mass. Think of this representation in connection with Figure \ref{fig:B_k_tree_proj} and the description of the $k$-tree projection of a Brownian CRT in the introduction. The $x_j$ represent masses of subtrees corresponding to leaves of the tree represented by $\ft$, while the $\beta_E$ represent totally ordered collections of subtree masses. 
In this interpretation, the intervals in $\beta_E$ that are closer to $0$ represent subtrees that are farther from the root of the CRT.

We refer to each top mass $x_j$, $j\in A$, and each interval in each of the partitions $\beta_E$, $E\in\ft$, as a \emph{block} of $T$. Formally, we denote the set of blocks by
\begin{equation}
 \block(\ft,(x_j,j\in A),(\beta_E,E\in\ft)) := A\cup\{(E,a,b)\colon E\in\ft,\,(a,b)\in \beta_E\}.
\end{equation}
We will write $\|\ell\|$ for the \emph{mass} of $\ell\in \block(T)$; i.e.\ for the top masses $\|\ell\|:=x_\ell$, $\ell\in A$, for the other blocks $\|\ell\|=\|(E,a,b)\|:=b-a$. 
Then $\sum_{\ell\in \block(T)}\|\ell\|=\|T\|$.

For each label set $A$ and each $\ft\in\TShape_A$, we topologize the set of $A$-trees with shape $\ft$ by the product over the topologies in the components. This can be metrized by setting
\begin{equation}\label{eq:ktree:metric_1}
 d_{\bT}(T,T') = \sum_{j\in A} |x_j-x'_j| + \sum_{E\in\ft}\dI(\beta_E,\beta'_E)
\end{equation}
for $T,T'\in\bTInt_A$ with shapes $\ft=\ft'$. Within the set of trees with a given label set $A$ and shape $\ft$, there is a single $A$-tree $0_\ft:=(\ft,(0)_{j\in A},(\emptyset)_{E\in\ft})$ of 
zero total mass; we topologize the space of all $A$-trees, for all finite label sets $A$, by identifying all of these trees of zero mass, thereby gluing these spaces together. This is metrized by
\begin{equation}\label{eq:ktree:metric_2}
 d_{\bT}(T,T') = \sum_{j\in A} x_j + \sum_{j\in A'}x'_j + \sum_{E\in\ft}\dI(\beta_E,\emptyset) + \sum_{E\in\ft'}\dI(\beta'_E,\emptyset)
\end{equation}
for $T\in\bTInt_A$, $T'\in\bTInt_{A'}$ with differing tree shapes. We note that $\dI(\beta,\emptyset) = \max\{\|\beta\|,\sD_\beta(\infty)\}$ for any $\beta\in\cI$. We will also 
write $0$ for a zero-mass tree.


\begin{proposition}\label{prop:Lusin}
 $\big(\big(\bigcup_A \big(\bTInt_A\setminus\{0_\ft,\ft\in\bT_A^{\rm shape}\}\big) \big)\cup\{0\},d_{\bT}\big)$ is a Lusin space.
\end{proposition}

\begin{proof}
 From Proposition \ref{prop:IPspace:Lusin}, $(\cI,\dI)$ is Lusin. Thus, so are the product topologies on the set of trees in $\bTInt_A$ with a given shape $\ft\in\TShape_A$, for each non-empty finite $A\subset\mathbb{N}$. The countable union of these sets is equipped with a metric that corresponds to gluing metric spaces of trees by identifying the unique zero-mass tree for each tree shape. This entails the claimed Lusin property.
\end{proof}

We are interested in $k$-tree-valued Markov processes that avoid certain degenerate states. For example, states with multiple zero top masses will be inaccessible by our evolutions. We also exclude states having a zero top mass with an empty partition on its parent edge. Such states will arise as left limits but force jumps ``away from the boundary.'' Specifically, for finite $A\subset\bN$ with $\#A\ge 2$, we define
\begin{align}
  \tdTInt_A \!:=\! \left\{ T \!=\! (\ft,(x_j,j\!\in\! A),(\beta_E,E\!\in\!\ft))\!\in\!\bTInt_A\, \middle|\!\!
    \begin{array}{l}
      x_i\!+\!x_j\!>\!0\mbox{ for all }E\!=\!\{i,j\}\!\in\!\ft\mbox{ and}\\
      x_i\!+\!\big\|\beta_{\parent{\{i\}}}\big\| \!=\! 0\mbox{ for at most one }i\!\in\! A
    \end{array}\!\!\!\right\}&\!\notag\\ \label{eq:k_tree_spaces}
  \TInt_{A} \!:= \!\left\{ T \!=\! \left(\ft,(x_j,j\!\in\! A),(\beta_E,E\!\in\!\ft)\right)\!\in\!\tdTInt_A \,\middle|\, 
  		x_j+\big\|\beta_{\parent{\{j\}}}\big\| \!>\! 0\mbox{ for all }j\!\in\! A\right\}\!.&
\end{align}
Let $I\colon\tdTInt_A\rightarrow A\cup\{\infty\}$ record $I(T)=i$ if $x_i+\|\beta_{\parent{\{i\}}}\| = 0$ and set $I(T)=\infty$ if $T\in\TInt_A$. In the former case, we say that \emph{label $i$ is degenerate in $T$}.

Because we will only ever consider single-leaf trees in the case where the leaf has label 1, we take the convention that $\tdTInt_{1} = [0,\infty)$ and $\TInt_{1} = (0,\infty)$, with this real number representing the mass on the leaf 1 component, which is then the total mass of the tree. We also define $\TInt_\emptyset = \{0\}$. 
As noted above Proposition \ref{prop:Lusin}, we identify all trees of zero mass. We take the convention of writing $0$ to denote such a tree. 

\section{Self-similar non-resampling and resampling $k$-tree evolutions}\label{sec:non_resamp_def}\label{sec:resamp_def}

The key building blocks for both non-resampling and resampling $k$-tree evolutions are killed $A$-tree evolutions:

\begin{definition}[Killed $A$-tree evolution]\label{def:killed_ktree}
 Consider some finite $A\subset\bN$ with $\#A\ge 2$ and an $A$-tree $T = (\ft,(m^0_j,j\in A),(\beta^0_E,E\in\ft))\in\TInt_A$. 
 \begin{itemize}
  \item For each type-2 edge $E = \{i,j\}\in \ft$, let $((m_i^y,m_j^y,\beta_E^y),y\ge0)$ denote a type-2 evolution from initial state $(m_i^0,m_j^0,\beta_E^0)$, and let $D_E$ denote its degeneration time (when one top mass and the edge partition vanish).
  \item For each type-1 edge $E = \parent{\{i\}}\in \ft$, let $((m_i^y,\beta_E^y),y\ge0)$ denote a type-1 evolution from initial state $(m_i^0,\beta_E^0)$, and let $D_E$ denote its degeneration time (when the top mass and edge partition both vanish).
  \item For each type-0 edge $E\in\ft$, let $(\beta_E^y,y\ge0)$ denote a type-0 evolution from initial state $\beta^0_E$ and define $D_E=\infty$.
 \end{itemize}
 We take these evolutions to be jointly independent. Let $D = \min_{E\in\ft}D_E$. Define $\cT^y = (\ft,(m_i^y,i\in A),(\beta_E^y,E\in\ft))$ for $y\in [0,D)$ and $\cT^y = 0$ for $y\ge D$. This is the \emph{killed $A$-tree evolution from initial state $T$}. We call $D$ the \emph{degeneration time} of the evolution.
 
For $\#A=1$, and $T\in\bT_1^{\rm int}=(0,\infty)$, define $(\cT^y)$ to be a \BESQ[-1] starting from $T$, killed upon hitting zero.
\end{definition}

In light of this construction, in an $A$-tree $T = (\ft,(x_j,j\in A),(\beta_E,E\in\ft))\in\TInt_A$, we refer to each type-2 edge partition with its two top masses, $(x_i,x_j,\beta_{\{i,j\}})$, as a \emph{type-2 compound}. Likewise, for a type-1 edge $E = \parent{\{j\}}$, we call $(x_j,\beta_E)$ a \emph{type-1 compound}, and for each type-0 edge $F$, the partition $\beta_F$ is a \emph{type-0 compound}. In Figure \ref{fig:B_k_tree_proj_cons}, $\beta_{[5]}^{(5)}$ is a type-0 compound, $\big(X_2^{(5)},\beta_{\{1,2,4\}}^{(5)}\big)$ is a type-1 compound, and $\big(X_3^{(5)},X_5^{(5)},\beta_{\{3,5\}}^{(5)}\big)$ and $\big(X_1^{(5)},X_4^{(5)},\beta_{\{1,4\}}^{(5)}\big)$ are type-2 compounds.


Recall from Chapter \ref{ch:prelim} that a Markov process $(\cT^y,y\ge0)$ is said to be \emph{self-similar} if it has the same semigroup as $(c\cT^{y/c},y\ge0)$ for all $c>0$.

\begin{proposition}\label{prop:killed:Markov} Killed $A$-tree evolutions are self-similar Borel right Markov processes, but they are not Hunt pocesses.
\end{proposition}
\begin{proof}  Killed $A$-tree evolutions are Borel right Markov processes as they are effectively tuples of independent type-0/1/2 evolutions, which are themselves self-similar Borel right Markov processes as noted in Proposition \ref{prop:012:pred}\ref{item:pred:Markov}, killed at a stopping time. 
We use stopping times $$S_n := \inf\left\{y\ge0\colon \min_{j\in A}(m_j^y+\|\beta_{\parent{\{j\}}}^y\|) < 1/n\right\},\quad n\ge1,$$
 to show that these processes fail to be quasi-left-continuous, hence fail to be Hunt processes.  
 These times are eventually strictly increasing (as soon as they exceed 0) and they converge to the killing time. Thus, the killing time is an increasing limit of stopping times, so it is visible in the left-continuous filtration and is a time at which the killed $A$-tree evolution is discontinuous.
\end{proof}

In the theory of Borel right Markov processes, \emph{branch states} are states that are not visited by the right-continuous Markov process but may be attained as a left limit, triggering an instantaneous jump. We will now define non-resampling $k$-tree evolutions with branch states in $\tdTInt_A\setminus\TInt_A$. When a type-1 or type-2 compound in an $A$-tree degenerates in a non-resampling evolution, we project this compound down and the evolution proceeds with one fewer leaf label.

However, in \cite{Paper2} we found that in the discrete regime, in order to construct a family of projectively consistent Markov processes, it was necessary to have degenerate labels sometimes swap places with other nearby, higher labels before dropping the degenerate component and its label with it. The following two definitions 
lead to an analogous construction in the present setting. The role of this mechanism in preserving consistency will be evident in the proof of Proposition \ref{prop:Dynkin:killed}.

\begin{figure}
 \centering
 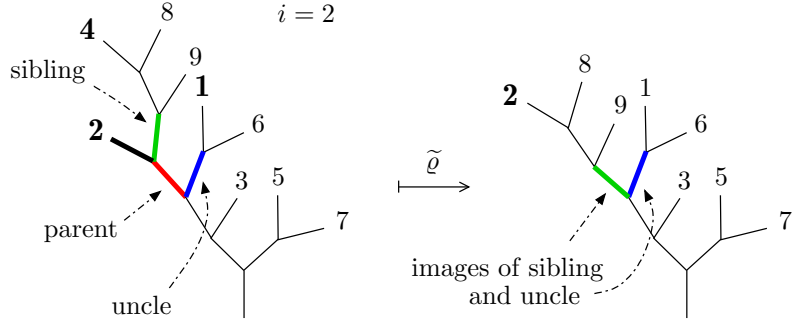
 \caption{Example of the swap-and-reduce map on a tree shape. Least labels in the two subtrees descended from sibling and uncle of leaf edge $\{i\}$ are shown in bold.\label{fig:label_swap}}
\end{figure}
\medskip

\textbf{Swap-and-reduce map for tree shapes. }
 Consider a tree shape $\ft\in\TShape_A$ on some label set with $\#A\ge2$ and label $i\in A$. Let 
 \begin{equation}\label{eq:swapred:iab}
  J(\ft,i) := \max\{i,a,b\} \quad \text{where} \quad a = \min \!\big({\rm sib}(\{i\})\big), \quad 
  b = \min \!\big({\rm sib}\big(\parent{\{i\}}\big)\big);
 \end{equation}
 i.e.\ $a$ and $b$ are the respective least elements in the label sets on the sibling and uncle of leaf edge $\{i\}$. In the special case that the parent $\parent{\{i\}} = A$, in which case $\{i\}$ has no uncle, we define $b=0$. This is illustrated in Figure \ref{fig:label_swap}, where
 $$\ft = \big\{ [9], \{5,7\}, \{1,2,3,4,6,8,9\}, \{1,2,4,6,8,9\}, \{1,6\}, \{2,4,8,9\}, \{4,8,9\}, \{4,8\}  \}.$$
 Leaf edge $\{2\}$ has sibling $\{4,8,9\}$ and uncle $\{1,6\}$, so $a=4$, $b=1$, $J(\ft,2) = \max\{2,4,1\} = 4$.
 
 We define a \emph{swap-and-reduce map on tree shapes}, 
  \begin{equation}\label{eq:varrhotilde}
      \widetilde\varrho\colon \TShape_A\times A \to \bigcup_{j\in A}\TShape_{A\setminus\{j\}}
  \end{equation}
 mapping $(\ft,i)$ to the tree shape $\ft'$ obtained from $\ft$ by first swapping labels $i$ and $j=J(\ft,i)$, then deleting the leaf subsequently labeled $j$ and contracting away its parent branch point. Formally, $\ft'$ is the image of $\ft\setminus\big\{\parent{\{i\}}\big\}$ under the map $\phi_{\ft,i}$ that modifies label sets $E\in\ft$ by first deleting label $i$ from the sets, and then replacing label $j$ by $i$. 
 In the example in Figure \ref{fig:label_swap}, with $i=2$ and $j=4$, $\phi_{\ft,i}$ is the map
 \newcommand{\mapsdown}{\rotatebox[origin=c]{-90}{\ensuremath{\mapsto}}}
 \newcommand{\pcom}{\hphantom{,}}
 \definecolor{mygreen}{rgb}{0,0.75,0}
\begin{align*}
\ft \!= &\big\{\ \ [9], \ \{5,7\},\{1,2,3,4,6,8,9\},\{1,2,4,6,8,9\}, &\!\!\!\!\!\textcolor{blue}{\{1,6\}},&\textcolor{red}{\{2,4,8,9\}},\textcolor{mygreen}{\{4,8,9\}},&\!\!\!\!\!\!\{4,8\}\!\big\}\\
&\quad\ \ \mapsdown\qquad\ \ \mapsdown\quad\qquad\qquad\mapsdown\qquad\qquad\qquad\mapsdown&\mapsdown\quad\ &\quad\ \qquad\qquad\mapsdown&\mapsdown\quad\ \\
 \ft' \!=&\big\{[9]\!\setminus\!\{4\},\{5,7\},\{1,2,3,6,8,9\}, \ \{1,2,6,8,9\}, &\!\!\!\!\!\textcolor{blue}{\{1,6\}},&\qquad\quad \textcolor{mygreen}{\{2,8,9\}},&\!\!\!\!\!\!\!\{2,8\}\!\big\}\!.
\end{align*}


Note that in the preceding definition, $\phi_{\ft,i}(E_1) = \phi_{\ft,i}(E_2)$ if and only if $E_1\setminus\{i\} = E_2\setminus\{i\}$. But the only distinct edges $E_1\neq E_2$ in $\ft$ with this relationship are the sibling and parent of leaf edge $\{i\}$. Thus, by excluding $\parent{\{i\}}$ from its domain, we render $\phi_{\ft,i}$ injective and ensure that the range of this map is an element of $\TShape_{A\setminus\{J(\ft,i)\}}$.

This swap-and-reduce map is the same as the down-move of the modified Aldous chain of Definition \ref{def:modified_AC}. This map on tree shapes induces a corresponding map for degenerate $A$-trees, where labels are swapped and the degenerate component is projected away, but everything else remains unchanged.
\medskip

\textbf{Swap-and-reduce map for $A$-trees. }
 Let $T = (\ft,(x_h,h\in A),(\beta_E,E\in\ft))$ $\in\tdTInt_A\setminus\TInt_A$. Recall that for such an $A$-tree, $I(T)$ denotes the unique index $i\in A$ for which $x_i+\big\|\beta_{\parent{\{i\}}}\big\| = 0$. We define $J\colon\tdTInt_A\setminus\TInt_A\rightarrow A$ by  
 $J(T) = J(\ft,I(T))$, as defined above. The \emph{swap-and-reduce map on $A$-trees} is the map
 \begin{equation}\label{eq:varrho}
   \varrho\colon \tdTInt_A\setminus \TInt_A \rightarrow \bigcup_{j\in A}\TInt_{A\setminus\{j\}}
\end{equation}
 that sends $T$ to 
 $(\widetilde\varrho(\ft,I(T)),(x'_h, h\in A\setminus\{J(T)\}),(\beta'_E, E\in \widetilde\varrho(\ft,I(T))))$ where 
 \begin{enumerate}\item[(i)] $x'_h = x_h$ for $h\neq I(T)$, and $x'_{I(T)} = x_{J(T)}$ if $I(T)\neq J(T)$,
  \item[(ii)] $\beta'_{E} = \beta_{\phi_{\ft,I(T)}^{-1}(E)}$ for each $E\in \widetilde\varrho(\ft,I(T))$, where $\phi_{\ft,I(T)}$ is the injective map defined in the
definition of the swap-and-reduce map for tree shapes.
  \end{enumerate} 

\begin{definition}[Non-resampling $k$-tree evolution]\label{def:nonresamp_1}
 Set $A_1 = [k]$ and fix some $\cT^0_{(1)} =T\in \TInt_{A_1}$. Inductively for $1\le n\le k-1$, let $(\cT^y_{(n)},y\in [0,\Delta_n))$ denote a killed $A_{n}$-tree evolution from initial state $\cT_{(n)}^0$, run until its degeneration time $\Delta_n$, conditionally independent of $(\cT_{(j)},j < n)$ given its initial state. If $n\le k-2$, we then set $A_{n+1} = A_n\setminus \{J(\cT_{(n)}^{\Delta_n-})\}$, let $\cT_{(n+1)}^0 = \varrho(\cT_{(n)}^{\Delta_n-})$ and repeat. For $n=k$, let $(\cT^y_{(k)},y\in[0,\Delta_k))$ denote a 
$\besq(-1)$ process from initial state $\|\cT_{(k-1)}^{\Delta_{k-1}-}\|$.
 
 For $1\le n\le k$ we define $D_n = \sum_{j=1}^n\Delta_j$ and set $D_0=0$. For $y\in [D_{n-1},D_{n})$ we define $\cT^y = \cT_{(n)}^{y-D_{n-1}}$. 
For $y\ge D_{k}$ we set $\cT^y = 0 \in \TInt_{\emptyset}$. Then $(\cT^y,y\ge0)$ is a \emph{non-resampling $k$-tree evolution from initial state $T$}. We say that at each time $D_n$, label $I(\cT^{D_n-})$ has \emph{caused degeneration} and label $J(\cT^{D_n-})$ is \emph{dropped in degeneration}.
\end{definition}


We now define a resampling $k$-tree evolution in which at degeneration times we first apply $\varrho$ and then jump into a random state according to a resampling kernel, which reinserts the label lost in degeneration, so that the evolution always retains all $k$ labels.
\medskip

\textbf{Label insertion operator $\oplus$. }
 \emph{For tree shapes.} Consider $\ft\in\TShape_A$. Given an edge $F\in\ft\cup\{\{h\}\colon h\in A\}$, we define $\ft\oplus (F,j)$ to be the tree shape with labels $A\cup\{j\}$ formed by replacing edge $F$ by a path of length 2, and inserting label $j$ as a child of the new branch point in the middle of the path. Formally, for each $E\in\ft$ we define (a) $\phi(E)=E\cup\{j\}$ if $F\subsetneq E$ and (b) $\phi(E)=E$ otherwise. Then $\ft\oplus (F,j)$ equals $\phi(\ft)\cup \{F\cup\{j\}\}$.
 
\begin{figure}
  \centering
  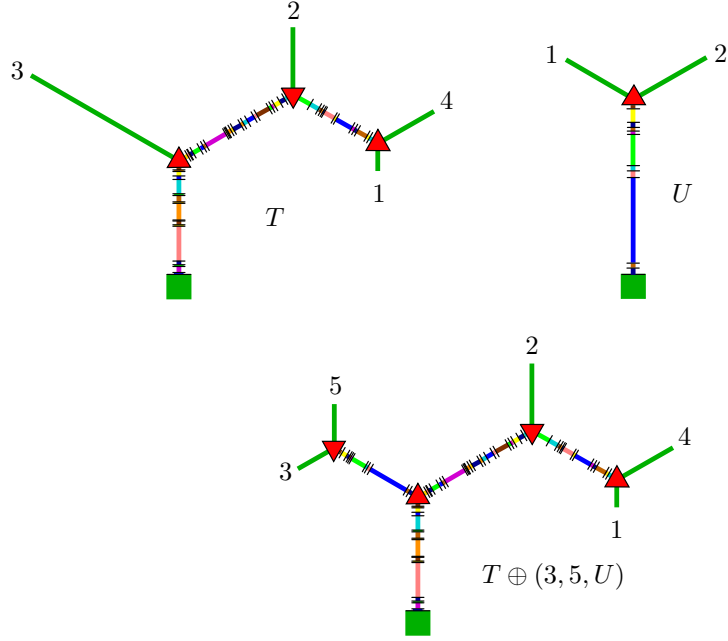
  \caption{The label insertion operator, inserting label $5$ into leaf block $3$ in a $4$-tree.\label{fig:label_insertion}}
 \end{figure}

 \emph{For $A$-trees.} Consider an $A$-tree $T = (\ft,(x_h,h\in A),(\beta_E,E\in\ft))$, a label $i\in A$, and a 2-tree 
 $U = (y_1,y_2,\gamma)\in \TInt_{2}$ with $\|U\| = 1$, where we have dropped the tree shape because all elements of $\TInt_{2}$ have the same shape. We define $T\oplus (i,j,U)$ to be the $(A\cup\{j\})$-tree formed by replacing the leaf block $i$ and its weight $x_i$ by the rescaled 2-tree in which label $i$ gets weight $x_iy_1$, a new label $j$ gets weight $x_iy_2$, and their new parent edge bears partition $x_i\gamma$. This operation is illustrated in Figure \ref{fig:label_insertion}. 
 Formally, 
  \begin{equation}\label{eq:insert1}T\oplus (i,j,U) = (\ft\oplus(\{i\},j),(x'_h,h\in A\cup\{j\}),(\beta'_E,E\in\ft\oplus(\{i\},j))),
  \end{equation}
where: (i) $(x'_i,x'_j,\beta'_{\{i,j\}}) = x_i U$, (ii) $x'_h=x_h$ for $h\notin \{i,j\}$, and (iii) $\beta'_E = \beta_{\phi^{-1}(E)}$ for $E\neq \{i,j\}$, where $\phi$ is as for tree shapes.
 
 Now consider a block $\ell = (F,a,b)\in \block(T)$. This block splits $\beta_F$ into $\beta_{F,0}\concat (0,b-a)\concat \beta_{F,1}$. We define $T\oplus (\ell,j,U)$ to be the $(A\cup\{j\})$-tree formed by inserting label $j$ into block $\ell$. In this definition, $U$ is redundant. Formally, 
 \begin{equation}\label{eq:insert2}T\oplus (\ell,j,U) = (\ft\oplus (F,j),(x'_h,h\in A\cup\{j\}),(\beta'_E,E\in\ft\oplus(F,j))),
 \end{equation} 
where: (i) $x'_h=x_h$ for $h\neq j$, (ii) $\beta'_E = \beta_{\phi^{-1}(E)}$ for $E\notin \{F,F\cup\{j\}\}$, and (iii) $(\beta'_F,x'_j,\beta'_{F\cup\{j\}}) = (\beta_{F,0},b-a,\beta_{F,1})$.
\medskip

\textbf{Resampling kernel for $A$-trees. }
 For finite non-empty $A\subset\bN$ and $j\in\bN\setminus A$, we define the \emph{resampling kernel} as the distribution of the tree obtained by inserting label $j$ into a block chosen at random according to the masses of blocks and, if the chosen block is a top mass $x_i$, then replacing the block by a rescaled Brownian reduced $2$-tree. More formally, we define a kernel $\Lambda_{j,A}$ from $\TInt_A$ to $\TInt_{A\cup\{j\}}$ by
 \begin{equation}
 \label{eq:resampling_kernel:def}
  \int_{T^\prime\in\bT_{A\cup\{j\}}^{\rm int}}\!\!\varphi(T^\prime)\Lambda_{j,A}(T,dT^\prime)=\sum_{\ell\in \block(T)}\frac{\|\ell\|}{\|T\|}\!\!\int_{U\in\TInt_{2}}\varphi(T\oplus(\ell,j,U))Q(dU),
 \end{equation}
 where $Q$ denotes the distribution of a Brownian reduced $2$-tree with leaf labels $\{1,2\}$, as defined in the introduction.

In \eqref{eq:B_ktree_resamp}, we will describe how these resampling kernels can be used to generate a Brownian reduced $k$-tree for $k\ge 3$.

\begin{definition}[Resampling $k$-tree evolution]\label{def:resamp_1}
 Fix some $\cT^0_{(1)} =T\in \TInt_{k}$. Inductively for $n\ge1$, let $(\cT^y_{(n)},y\in [0,\Delta_n))$ denote a killed $k$-tree evolution from initial state $\cT_{(n)}^0$, run until its degeneration time $\Delta_n$, conditionally independent of $(\cT_{(j)},j< n)$ given its initial state. We define $\cT_{(n+1)}^0$ to have conditional distribution $\Lambda_{J_n,[k]\setminus\{J_n\}}\big(\varrho(\cT_{(n)}^{\Delta_n-}),\cdot\,\big)$ given $(\cT_{(j)},j\le n)$, where $J_n = J(\cT_{(n)}^{\Delta_n-})$.
 
 We set $D_0=0$ and define $D_n = \sum_{j=1}^n\Delta_j$, $n\ge 1$. For $y\in [D_{n-1},D_{n})$ we define $\cT^y = \cT_{(n)}^{y-D_{n-1}}$. For $y\ge D_\infty:=\sup_{n\ge 0}D_n$ we set $\cT^y= 0 \in \TInt_{\emptyset}$. Then $(\cT^y,\,y\ge 0)$ is a \emph{resampling $k$-tree evolution with initial state $T$}.
\end{definition}

\begin{theorem}\label{thm:Markov}
 Non-resampling and resampling $k$-tree evolutions are self-\linebreak similar Borel right Markov processes, but they are not Hunt processes.
\end{theorem}
\begin{proof}
 As noted in Proposition \ref{prop:killed:Markov}, killed $A$-tree evolutions are self-similar Borel right Markov processes. Note that because these evolutions in the various type-0/1/2 compounds in the tree are independent and their degeneration times are continuous random variables, almost surely one of them degenerates before all of the others.  Since each type-$i$ compound has $i$ positive top masses and positive interval partition mass at almost all times before its degeneration time, $\cT^{D_1-}\in\tdTInt_{k}$ a.s. Therefore, the non-resampling and resampling $k$-tree evolutions are well-defined. Moreover, the type of construction undertaken in Definitions \ref{def:nonresamp_1} and \ref{def:resamp_1} of non-resampling and resampling $k$-tree evolutions is well-studied; it yields a Borel right Markov process by Th\'eor\`eme 1 and the Remarque on p.\ 474 of Meyer \cite{Mey75}.

Recall from the proof of Proposition \ref{prop:killed:Markov} the sequence $(S_n)$ of stopping times that increases to the killing time, which here is the first degeneration/resampling time. Indeed, for non-resampling and resampling $k$-tree evolutions, the jump to $0$ has been replaced by a swap-and-reduce and/or resampling jump, so the discontinuity along $(S_n)$ prevails, hence the Hunt property still fails for non-resampling and resampling $k$-tree evolutions.
\end{proof}

\begin{theorem}\label{thm:total_mass}
 Let $(\cT^y,y\ge0)$ be a non-resampling or resampling $k$-tree evolution with initial state with mass $\|\cT^0\| = m$. Then the total mass process $(\|\cT^y\|,y\ge0)$ has law $\BESQ_m(-1)$.
\end{theorem}

At this stage, we can only prove a partial form of Theorem \ref{thm:total_mass}, as follows.

\begin{proposition}\label{prop:total_mass_0}
 The total mass process of a (self-similar) non-resampling $A$-tree evolution is a \BESQ[-1]. The total mass process of a resampling $k$-tree evolution is a \BESQ[-1], killed at the random time $D_\infty := \sup_nD_n$.
\end{proposition}
\begin{proof}
 Let $(\cT^y,y\ge0)$ denote a non-resampling $k$-tree evolution. Up until its first degeneration, its total mass $\|\cT^y\|$ is the sum of the total masses of $k-1$ type-0/1/2 evolutions -- one compound for each internal edge $E\in\ft$ of the tree shape. In particular, the sum of the ``type numbers'' of these compounds is $k$: if we let $n_i$ denote the number of type-$i$ compounds, $i\in\{0,1,2\}$, then
 $$k-1 = n_0+n_1+n_2 \qquad\text{and} \qquad k = 0\times n_0 + 1\times n_1 + 2\times n_2 = n_1 + 2n_2.$$ 
 This gives $n_2 = n_0+1$. By Proposition \ref{prop:012:mass}, 
 the total mass process of a type-$i$ evolution is a \BESQ[1-i]. Then, by the (generalized) additivity of squared Bessel processes \cite[Proposition 1.1]{PW18}, the sum of these total masses, $(\|\cT^y\|,y\in [0,D_1])$, evolves as a squared Bessel process with parameter $1\times n_0 + 0\times n_1 - 1\times n_2 = -1$, stopped at a stopping time in a filtration to which the squared Bessel process is adapted. Moreover, the same argument and the strong Markov property show that the total mass continues to evolve as a \BESQ[-1] between the first and second degeneration times, and so on. Thus, the process evolves as a \BESQ[-1] until its absorption at $0$. The same argument proves the assertion as stated for the resampling $k$-tree evolution.
\end{proof}

In Section \ref{sec:non_acc} (and Appendix \ref{sec:non_acc_2}), we will complete the proof of Theorem \ref{thm:total_mass} by establishing the following result.

\begin{proposition}\label{prop:non_accumulation}
 For resampling $k$-tree evolutions with degeneration times $D_n$, $n\ge 1$, the limit $D_\infty = \lim_{n\rightarrow\infty}D_n$ equals $\inf\{y\ge0\colon \|\cT^{y-}\|=0\}$, and this is a.s.\ finite.
\end{proposition}

\section{Pseudo-stationarity of self-similar $k$-tree evolutions}\label{sec:pseudo}

Recall the law of a Brownian reduced $k$-tree stated in Proposition \ref{prop:B_ktree} in terms of a uniform random tree shape $\tau$ from $\bT_{[k]}^{\rm shape}$, a 
${\tt Dirichlet}(\frac12,\ldots,\frac12)$ mass split onto the $k$ top masses and the $k-1$ edges, and independent ${\tt PDIP}(\frac12,\frac12)$ proportions to further split the edge
masses into interval partitions. Denote by $Q_{z,A}(dU)$ the distribution on $\TInt_A$ of the $A$-tree obtained from the distribution $Q_{z,[k]}(dU)$ of a Brownian reduced $k$-tree scaled to have total mass $z$, with leaves then relabeled by the increasing bijection $[k]\rightarrow A$, for $k=\#A$. The resampling kernel $\Lambda_{j,A}$ of Definition \ref{def:resamp_1} satisfies
\begin{equation}\label{eq:B_ktree_resamp}
 \int_{T\in\TInt_{k}}\!\!Q_{z,[k]}(dT)f(T) =\! \int_{(T_i,i\in [2,k])}\!\!\Lambda_{2,[1]}(z,dT_2)\cdots\Lambda_{k,[k-1]}(T_{k-1},dT_k)f(T_k),\!\!
\end{equation}
where $z\in\TInt_{1}$ denotes the $1$-tree with weight $z$ on its sole component, leaf 1. This formula indicates that the Markov chain that begins with $z$ and at each step, adds a successive label via the resampling kernel, has as its path a consistent system of Brownian reduced $k$-trees, $k\ge1$, each scaled to have total mass $z$. Like Proposition \ref{prop:B_ktree}, this formula follows from the development in \cite[Section 3.3]{PW13}.

\begin{proposition}\label{prop:pseudo:resamp}
 Let $(\cT^y,y\ge 0)$ be a resampling $k$-tree evolution starting from an independently scaled Brownian reduced $k$-tree of arbitrary total mass $M$, and let $B\sim\besq_M(-1)$. Then at any fixed time $y\ge0$, $\cT^y$ has the distribution of an independently scaled Brownian reduced $k$-tree of mass $B(y)$.  
\end{proposition}

In light of this result, we refer to the laws of independently scaled Brownian reduced $k$-trees as the \emph{pseudo-stationary laws for resampling $k$-tree evolutions.} 
Before we prove Proposition \ref{prop:pseudo:resamp}, recall that type-0 evolutions do not degenerate (and are reflected when reaching zero total mass), while we say that type-1 evolutions degenerate when they reach the (absorbing) state of zero total mass and type-2 evolutions degenerate when they reach a single-top-mass state on an empty interval partition. In particular, total mass evolutions conditioned on no degeneration up to time $y$ are unaffected by the conditioning for type-0 evolutions as we are conditioning on an event of probability 1, while they are conditioned to be positive for type-1 evolutions and conditioned on an event that depends on the underlying type-2 evolution for type 2.

\begin{proposition}\label{prop:pseudo:pre_D}
 Let $(\cT^y,y\ge0)$ be a killed/non-resampling/resampling $k$-tree evolution starting from a Brownian reduced $k$-tree scaled by an independent random initial mass. Then for $y\ge0$, given $\{D_1>y\}$, the tree $\cT^y$ is again conditionally a Brownian reduced $k$-tree scaled by an independent random mass. In the special case that $\|\cT^0\|\sim\GammaDist[k-\frac12,\lambda]$, given $\{D_1>y\}$, $\|\cT^y\|$ has conditional law $\GammaDist[k-\frac12,\lambda/(2\lambda y+1)]$. 
\end{proposition}
\begin{proof}
 First, suppose $M\sim \GammaDist[k-\frac12,\lambda]$. Note that
 \begin{equation}
  \bP(D_1>y) = (2y\lambda+1)^{-k},
 \end{equation}
 since each type-1 compound contributes $(2y\lambda+1)^{-1}$ by \eqref{eq:pseudo:type_1_survival} and each type-2 compound contributes $(2y\lambda+1)^{-2}$, by Proposition \ref{prop:degeneration}, all independently, with $k$ top masses altogether. Conditioning on non-degeneration means conditioning each independent type-$i$ evolution, $i=1,2$, not to degenerate; thus, this conditioning does not break the independence of these evolutions. By Proposition \ref{prop:012:pseudo}, the conditional distribution of each edge partition and top mass at time $y$ is the same as the initial distribution, but with each mass and partition scaled up by a factor of $2\lambda y+1$, as claimed. 
 
 The result for deterministic initial total mass follows by Laplace inversion, and for general random mass by integration. We leave the details to the reader and refer to Propositions \ref{prop:type2:pseudo} and \ref{prop:pseudo_f} or indeed to \cite[Proposition 4.3 and the proof of Theorem 1.5]{Paper1-2} for similar arguments.
\end{proof}




\begin{proposition}\label{prop:pseudo:degen}
 Let $(\cT^y,y\ge 0)$ be a killed/non-resampling/resampling $k$-tree evolution starting from a Brownian reduced $k$-tree, scaled by any independent initial mass $M$, and let $y\ge 0$. Then the following hold.
 \begin{enumerate}[label=(\roman*), ref=(\roman*)]
  \item \label{item:pseudoD:J}
  	The label $J = J(\cT^{D-})$ dropped at the first degeneration time $D = D_1$ has law $\bP(J=2)=2/k(2k-3)$ and $\bP(J=j)=(4j-5)/k(2k-3)$, $j\in\{3,\ldots,k\}$.
  \item \label{item:pseudoD:swapred}
  	Conditionally given $J\!=\!j$, the normalized tree $\varrho\big(\cT^{D-}\big) \big/ \big\|\cT^{D-}\big\|$, which is simply $\cT^D / \big\|\cT^D\big\|$ in the non-resampling case, is a Brownian reduced $([k]\!\setminus\!\{j\})$-tree.
  \item \label{item:pseudoD:indep}
  	The pair $\big(J\big(\cT^{D-}\big),\varrho\big(\cT^{D-}\big) \big/ \big\|\cT^{D-}\big\|\big)$ is independent of $\big(M,D,\big\|\cT^{D-}\big\|\big)$.
  \item \label{item:pseudoD:gamma_mass}
  	In the special case that $M\sim\GammaDist[k-\frac12,\lambda]$, conditionally given $D=y$,
  	$$\|\cT^{y-}\|\sim{\tt Gamma}\Big(k-\textstyle\frac32,\lambda/(2y\lambda+1)\Big).$$
  \item \label{item:pseudoD:resamp}
  	In the resampling case, properties \ref{item:pseudoD:J}, \ref{item:pseudoD:swapred}, and \ref{item:pseudoD:indep} also hold at each subsequent degeneration time $D=D_n$, $n\ge1$. Moreover, $\cT^{D_n}/\big\|\cT^{D_n}\big\|$ is a Brownian reduced $k$-tree.
 \end{enumerate}
\end{proposition}
\begin{proof}
 First, we derive \ref{item:pseudoD:resamp} as a consequence of the other assertions. Equation \eqref{eq:B_ktree_resamp}, along with exchangeability of labels in the Brownian reduced $k$-tree, implies that taking a Brownian reduced $([k]\!\setminus\!\{j\})$-tree and inserting label $j$ via the resampling kernel results in a Brownian reduced $k$-tree. Thus, \ref{item:pseudoD:swapred} gives us $\cT^{D_1}/\big\|\cT^{D_1}\big\| \stackrel{d}{=} \cT^0/\big\|\cT^0\big\|$. Assertion \ref{item:pseudoD:resamp} for subsequent degeneration times then follows by induction and the strong Markov property of resampling $k$-tree evolutions at degeneration times.
 
 It remains to prove \ref{item:pseudoD:J}, \ref{item:pseudoD:swapred}, \ref{item:pseudoD:indep}, and \ref{item:pseudoD:gamma_mass} for killed evolutions. We begin with the special case $M\!\sim\!\GammaDist[k\!-\!\frac12,\lambda]$. In this case, by Proposition \ref{prop:B_ktree}, each type-$i$ compound has initial mass $\GammaDist[(i+1)/2,\lambda]$, with all initial masses being independent. For $y>0$, each type-1 compound avoids degeneration prior to time $y$ with probability $(2y\lambda+1)^{-1}$ by \eqref{eq:pseudo:type_1_survival}. For type-2 the corresponding probability is $(2y\lambda+1)^{-2}$ by Proposition \ref{prop:degeneration}. Moreover, when a type-2 compound degenerates, each of the two labels is equally likely to be the one to cause degeneration Proposition \ref{degdist}. Thus, the first label $I$ to cause degeneration is uniformly random in $[k]$ and is jointly independent with the initial tree shape $\tau_k$ and the time of degeneration $D$. But recall that this does not necessarily mean that label $I$ is dropped at degeneration; we must account for the swapping part of the swap--and-reduce map $\varrho$. 
 
 This places us in the setting of our study of the modified Aldous chain of Definition \ref{def:modified_AC}, where we begin with a uniform random tree shape $\tau_k$ and select a uniform random leaf $I$ for removal, with the same label-swapping dynamics as in the definition of $\widetilde{\varrho}$ in Section \ref{sec:non_resamp_def}. In particular, \cite[Corollary 4]{Paper2} gives the distribution of $J:=J(\tau_k,I)$ as $p_1:=\bP(J=1)=0$, $p_2:=\bP(J=2)=2/k(2k-3)$ and $p_j:=\bP(J=j)=(4j-5)/k(2k-3)$ for $j\in\{3,\ldots,k\}$; and \cite[Lemma 3]{Paper2} says that given $\{J=j\}$, the tree shape $\tau_{k-1}$ after swapping and reduction is conditionally uniformly distributed on $\TShape_{[k]\setminus\{j\}}$. 
 
 Since $D$ is independent of $(\tau_k,I)$ and since $\tau_{k-1} = \widetilde{\varrho}(\tau_k,I)$, if we additionally condition on $D$, then the above conclusion still holds: given $\{D=y,J=j\}$, the resulting tree shape $\tau_{k-1}$ is still conditionally uniform on $\TShape_{[k]\setminus\{j\}}$. Moreover, by the independence of the evolutions of the type-0/1/2 compounds in the tree (prior to conditioning), each type-1 compound that does not degenerate is conditionally distributed as a type-1 evolution in pseudo-stationarity, conditioned not to die up to time $y$, and correspondingly for type-2 and type-0 compounds. As noted in Proposition \ref{prop:012:pseudo}, the law at time $y$ is the same as the initial distribution, but with total mass scaled up by a factor of $2y\lambda+1$, meaning that each top mass $m_h^y$ in these compounds is conditionally independent with law $\GammaDist[\frac12,\lambda/(2y\lambda+1)]$, and each internal edge partition $\beta_E^y$ is a conditionally independent \PDIP[\frac12,\frac12] scaled by a $\GammaDist[\frac12,\lambda/(2y\lambda+1)]$ mass. Similarly, if the degeneration occurs in a type-2 compound in $\tau_k$, then the remaining top mass in that compound also has conditional law $\GammaDist[\frac12,\lambda/(2y\lambda+1)]$ by Proposition \ref{degdist}. Thus, $\varrho\big(\cT^{D-}\big)/\big\|\cT^{D-}\big\|$ is conditionally a Brownian reduced $([k]\!\setminus\!\{j\})$-tree, as claimed, and is conditionally independent of $\|\cT^D\| \sim \GammaDist[k-\frac32,\lambda/(2y\lambda+1)]$.
 
  This completes the proof of \ref{item:pseudoD:gamma_mass} as well as of \ref{item:pseudoD:J} and \ref{item:pseudoD:swapred} in the special case when the initial total mass is $M\sim\GammaDist[k-\frac12,\lambda]$. Moreover, since the above conditional law of the normalized tree does not depend on the particular value of $D$, we find that the pair in \ref{item:pseudoD:indep} is independent of $\big(D,\big\|\cT^D\big\|\big)$ in this case; it remains to show independence from $M$.
 
 Now consider a $k$-tree evolution $(\wbT^{y},y\ge 0)$ starting from a unit-mass Brownian reduced $k$-tree, with degeneration time $\wbD$. Let $M$ be independent of this evolution with law $\GammaDist[k-\frac12,\lambda]$. By the self-similarity noted in Proposition \ref{prop:killed:Markov}, $\cT^y = M\wbT^{y/M}$, $y\ge0$, is a $k$-tree evolution with initial mass $M$, as studied above. In particular, $  \big(D,\big\|\cT^{D-}\big\|\big) = \big(M\wbD,M\big\|\wbT^{\wbD-}\big\|\big)$ and
 \begin{equation*}
  \left(J\big(\cT^{D-}\big),\frac{\varrho\big(\cT^{D-}\big)}{\big\|\cT^{D-}\big\|}\right) = \left(J\big(\wbT^{\wbD-}\big),\frac{\varrho\big(\wbT^{\wbD-}\big)}{\big\|\wbT^{\wbD-}\big\|}\right).
 \end{equation*}
 %
 %
 We showed that
 \begin{equation*}
 \begin{split}
  &\int_0^\infty 
  \bE \left[ f\left( J\big(\wbT^{\wbD-}\big) , \frac{\varrho\big(\wbT^{\wbD-}\big)}{\big\|\wbT^{\wbD-}\big\|}\right)
  		g\left(x\left\|\wbT^{\wbD-}\right\|,x\wbD\right)\right] \frac{\lambda^{k-\frac12}}{\Gamma(k-\frac12)}x^{k-\frac32}e^{-\lambda x} dx\\
  	&\ \ = \bE [ f(J^*,\wbT^*) ] \int_0^\infty \bE\left[g\left(x\left\|\wbT^{\wbD-}\right\|,x\wbD\right)\right] \frac{\lambda^{k-\frac12}}{\Gamma(k-\frac12)}x^{k-\frac32}e^{-\lambda x} dx,
 \end{split}
 \end{equation*}
 where $J^*$ and $\wbT^*$ have the laws described in \ref{item:pseudoD:J} and \ref{item:pseudoD:swapred} for the dropped label and the normalized tree. If we cancel out the constant factors of $\lambda
^{k-\frac12}/\Gamma(k-\frac12)$ on each side, appeal to Laplace inversion, and then cancel out factors of $x^{k-\frac32}$, then we find
 \begin{equation*}
  \bE\!\left[ f\!\left(\! J\big(\wbT^{\wbD-}\big) , \frac{\varrho\big(\wbT^{\wbD-}\big)}{\big\|\wbT^{\wbD-}\big\|}\right)
  		g\left(x\left\|\wbT^{\wbD-}\right\|,x\wbD\right)\right]\! = \bE [ f(J^*,\wbT^*) ] \bE\!\left[g\left(x\left\|\wbT^{\wbD-}\right\|,x\wbD\right)\right]
 \end{equation*}
 for Lebesgue-a.e.\ $x>0$. By the self-similarity observed in Proposition \ref{prop:killed:Markov}, it follows that this holds for every $x>0$, thus proving \ref{item:pseudoD:J} and \ref{item:pseudoD:swapred} for fixed initial mass, or for any independent random initial mass by integration. This formula also demonstrates the independence of the dropped label and normalized tree from the degeneration time and mass at degeneration. Since the laws that we find for the dropped index and the normalized tree do not depend on the initial mass $x$, this also proves \ref{item:pseudoD:indep}.
\end{proof}

This means that for any scaled Brownian reduced $k$-tree, the $k$-tree evolution without resampling runs through independent multiples of $\widebar{\cT}_k^{*},\widebar{\cT}_{k-1}^{*},\ldots,\widebar{\cT}_2^{*},\widebar{\cT}_1^{*},0$, where each of the $n$-trees for $1\le n\le k$ has as its distribution the appropriate mixture of Brownian reduced trees with label sets of size $n$. We now combine the previous results to establish the laws of independently scaled Brownian reduced $k$-trees as pseudo-stationary laws for resampling $k$-tree evolutions.

\begin{proof}[Proof of Proposition \ref{prop:pseudo:resamp}]
 By Proposition \ref{prop:pseudo:degen} and \eqref{eq:B_ktree_resamp}, conditional on $D_1 = z>0$, the tree $\cT^{z}$ is distributed as a Brownian reduced $k$-tree scaled by an independent random mass. By the strong Markov property at degeneration times and induction, the same holds conditional on $D_n=z>0$, for any $n\ge1$.
 
 While we have not yet proved Proposition \ref{prop:non_accumulation}, that $D_\infty$ is the hitting time at zero for the total mass process, it is easier to prove it in this special setting. Indeed, the rescaled inter-degeneration times $(D_{n+1}-D_n)/\|\cT^{D_n}\|$ are independent and identically distributed. As shown in Proposition \ref{prop:total_mass_0}, the total mass $\|\cT^y\|$ evolves as a \BESQ[-1] up until time $D_\infty$, so this time must be a.s.\ finite so as to not exceed the time of absorption for the \BESQ. From this we conclude that the masses $\|\cT^{D_n}\|$ must tend to zero almost surely.
 
 Let $(\cT^x_n,x\ge0)$ denote a resampling $k$-tree evolution with $\cT^0_n=\cT^{D_n}$, conditionally independent of $(\cT^y,y\ge0)$ given $\cT^{D_n}$. Now, we condition on $D_n = z\le y < D_{n+1}$. Then by the strong Markov property at time $D_n$, the tree $\cT^y$ is conditionally distributed according to the conditional law of $\cT^{y-z}_n$, given that $(\cT^x_n)$ does not degenerate prior to this time. By Proposition \ref{prop:pseudo:pre_D}, this too is a Brownian reduced $k$-tree scaled by an independent random mass. Integrating out this conditioning preserves the property of $\cT^y$ being a Brownian reduced $k$-tree scaled by an independent mass.
\end{proof}

As in Lemma \ref{pretotal}, we can strengthen the pseudo-stationarity of Proposition \ref{prop:pseudo:resamp}, at fixed times $y$, to certain stopping times $Y$.

\begin{corollary}[Strong pseudo-stationarity]\label{cor:pseudo:resamp}
 Consider a resampling $k$-tree evolution $(\cT^y,\,y\ge 0)$, whose initial state is an independent multiple of a random state with unit-mass pseudo-stationary distribution $Q_{1,[k]}$. 
  Denote by $M(y) = \|\cT^y\|$, $y\ge 0$, the associated total mass process and by $(\cF^y_{\rm mass},\,y\ge 0)$ the right-continuous filtration it generates. Let $Y$ be a stopping time in this filtration. Then for all $\cF^Y_{\rm mass}$-measurable $\eta\colon\Omega\!\rightarrow\![0,\infty)$ and measurable $H\colon\bT_k^{\rm int}\!\rightarrow\![0,\infty)$,
  $$\bE\left[\eta H(\cT^Y)\right]=\bE\left[\eta Q_{M(Y),[k]}[H]\right].$$
\end{corollary}
\begin{proof}
 The proof of Lemma \ref{pretotal} is easily adapted, omitting all indicators such as $\mathbf{1}\{D>y\}$, and using Proposition \ref{prop:pseudo:resamp} instead of Propositions  \ref{prop:type2:pseudo} or \ref{prop:pseudo_f}. 
\end{proof}

\section{Unit-mass $k$-tree evolutions}\label{sec:ktree:dePoiss}

We proceed as in the case $k=2$ in \eqref{eq:dePoi_time_change} and Definition \ref{defdePoiss}, noting that 
$$\cJ^\circ=\{(m_1,m_2,\beta_{\{1,2\}})\colon (\ft_2,(m_j,j\in[2]),(\beta_E,E\in\ft_2))\in\tdTInt_{\{1,2\}}\}\cup\{(0,0,\emptyset)\},$$
where $\ft_2\in\mathbb{T}^{\rm shape}_2$ is the unique 2-tree shape, which is irrelevant for the total mass of $T=(\ft_2,(m_j,j\in[2]),(\beta_E,E\in\ft_2))\in\tdTInt_{\{1,2\}}$ as
$\|T\|=m_1+m_2+\|\beta_{\{1,2\}}\|=\|(m_1,m_2,\beta_{\{1,2\}})\|$. In general, given a \cadlag\ path $\fT = (T^y,y\!\ge\!0)$ in $\bigcup_A\tdTInt_{A}$, consider the de-Poissonization time-change function $\rho_\fT\colon[0,\infty)\rightarrow[0,\infty]$,
\begin{equation}\label{eq:dePoi:time_change}
 \rho_\fT(u) := \inf\left\{y\ge 0\colon\int_0^y\|T^x\|^{-1}dx>u\right\},\quad u\ge 0.
\end{equation}
If the total mass process $(\|T^y\|,y\ge0)$ evolves as a $\besq(-1)$, as in Theorem \ref{thm:total_mass}, then $\rho(\fT)$ is bijective onto $[0,\zeta)$ a.s., where $\zeta = \inf\{y\ge0\colon \|T^y\| = 0\}$ is a.s.\ finite; see e.g.\ \cite[p.\ 314-5]{GoinYor03}. 

Let $\TInt_{k,1} := \big\{T\in\TInt_{k}\colon \|T\| = 1\big\}$. We define unit-mass $k$-tree evolutions by de-Poissonizing self-similar $k$-tree evolutions.

\begin{definition}\label{def:dePoi}
 Let $\fT = (\cT^y,y\ge0)$ denote a self-similar resampling (respectively, non-resampling) $k$-tree evolution from initial state $T\in\TInt_{k,1}$. Then
 $$\widebar\cT^u := \big\|\cT^{\rho_u(\fT)}\big\|^{-1} \cT^{\rho_u(\fT)},\quad u\ge0$$
 is a \emph{unit-mass resampling} (resp.\ \emph{non-resampling}) \emph{$k$-tree evolution from initial state $T$}. 
 %
\end{definition}

\begin{theorem}\label{thm:dePoi}
Unit-mass resampling and non-resampling $k$-tree evolutions are Borel right Markov processes. The former are stationary with the laws of the Brownian reduced $k$-trees. The latter are eventually absorbed at the state $1\in\TInt_{1,1}$ of the degenerate tree consisting of only one top mass of unit weight.
\end{theorem}
\begin{proof}
 The proofs of the $k=2$ case in Theorems \ref{absorption}--\ref{thm:stationary} are easily adapted using the Markov property, total mass, and pseudo-stationarity results of the self-similar $k$-tree evolutions obtained in Theorems \ref{thm:Markov} and \ref{thm:total_mass} and Proposition \ref{prop:pseudo:resamp}.
\end{proof}

We also obtain the following result for (unit-mass resampling) $k$-tree evolutions.

\begin{corollary}\label{cor:WF}
Let $\widebar \fT = (\widebar \cT^u,u\!\ge\!0)= ((\widebar \ft^u_k,(\widebar X^u_j ,j\!\in\![k]),(\widebar\beta^u_E,E\!\in\!\widebar\ft^u_k)), u\!\geq\! 0)$ denote a $k$-tree evolution started from $T= (\ft_k,(X_j,j\!\in\![k]),(\beta_E,E\!\in\!\ft_k))$ $\in\TInt_{k,1}$, and let $\tau$ be the first time either a top mass or an interval partition has mass $0$.  Observe that $\tau \leq \widebar D_1$, where $\widebar D_1$ is the first time $\widebar \fT$ resamples, so for $u<\tau$, $\widebar \ft^u_k=\ft_k$.  Then $((\widebar X^{u/4}_j,j\!\in\![k]),(\|\widebar\beta^{u/4}_E\|,E\!\in\!\ft_k)), 0\leq u <\tau)$ is a Wright--Fisher diffusion, killed when one of the coordinates vanishes, with parameters $-\frac12$ respectively $\frac12$ for coordinates corresponding to top masses respectively masses of interval partitions.
\end{corollary}

\begin{proof}
 Let $\fT$ be a self-similar $k$-tree evolution started from $T$.  By Propositions \ref{prop:012:mass} and \ref{prop:012:concat}(i), up until the first time a top mass or the mass of an interval partition is zero, the top masses evolve as \BESQ[-1] processes and the masses of internal interval partitions evolve as \BESQ[1] processes, and all of these are independent.  The effect of de-Poissonization procedure in Definition \ref{def:dePoi} on these evolutions is identical to Pal's de-Poissonization procedure \cite{Pal11, Pal13} used to construct Wright--Fisher diffusions. See \eqref{eq:WF_constr}. The result follows.\qedhere
\end{proof}


\chapter{Projective consistency of $k$-tree evolutions}
\label{ch:consistency}

Recall that the ultimate goal of this memoir is to construct a path-continuous continuum-tree-valued Markov process. 
Our strategy, as indicated in the statements of Theorems \ref{thm:intro:k_tree} and \ref{thm:intro:AD}, is to obtain this as a projective limit of $k$-tree-valued processes. We think of these $k$-trees as projections of a Brownian CRT, as described in Section \ref{sec:intro:B_k_tree}. In this chapter, we prove projective consistency results for the (self-similar) $k$-tree evolutions and their unit-mass variants, as defined in Chapter \ref{ch:constr}, thereby proving Theorem \ref{thm:intro:k_tree}. In Chapter \ref{ch:properties}, we study the projective limits.

\begin{definition}[Projection maps for $A$-trees]\label{def:proj}
 For $j\in\bN$ and finite $A\subset\bN$ with $\#(A\setminus\{j\})\ge 1$, we define a projection map
 $$\pi_{-j}\colon \bTInt_A\to\bTInt_{A\setminus\{j\}}$$
 to remove label $j$ from an $A$-tree, as follows. Let $T = (\ft,(x_i,i\in A),(\beta_E,E\in\ft))\in \bTInt_A$. If $j\notin A$, then $\pi_{-j}(T) = T$. Otherwise, let $\phi$ denote the map $E\mapsto E\setminus\{j\}$ with domain $\ft\setminus\big\{\parent{\{j\}}\big\}$. As noted for a similar map in Section \ref{sec:non_resamp_def}, this map is injective. Then $\pi_{-j}(T) := (\ft',(x'_i,i\in A\setminus\{j\}),(\beta'_E,E\in\ft'))$, where
 \begin{enumerate}[label=(\roman*),ref=(\roman*)]
  \item $\displaystyle\ft' = \phi(\ft) = \big\{E\setminus\{j\}\colon E\in \ft\setminus\big\{\parent{\{j\}}\big\}\big\}$,
  \item if $E = \parent{\{j\}}$ is a type-1 edge in $\ft$, then $\beta'_{E\setminus\{j\}} = \beta_{E\setminus\{j\}}\concat (0,x_j)\concat\beta_E$,\label{item:proj:merge}
  \item if $\parent{\{j\}} = \{a,j\}$ is a type-2 edge in $\ft$, then $x_a' = x_a+x_j+\|\beta_{\{a,j\}}\|$, \label{item:proj:add}
  \item if $i\in A\setminus\{j\}$ is \emph{not} the sibling of $\{j\}$ in $\ft$, then $x'_i = x_i$, and 
  \item if $E\in \ft'$ is \emph{not} the sibling of $\{j\}$ in $\ft$, then $\beta'_E = \beta_{\phi^{-1}(E)}$.
 \end{enumerate}
 
 For $k\ge1$ and any finite $A\subseteq\bN$ with $A\cap[k]\neq\emptyset$, we define $\pi_k\colon \bTInt_A\to\bTInt_{A\cap [k]}$ to be the composition $\pi_{-(k+1)}\circ\pi_{-(k+2)}\circ\cdots\circ\pi_{-\max(A)}$ to project onto trees labeled by $A\cap[k]$. It is straightforward to check that this composition of projection maps commutes. 
%
\end{definition}

\begin{figure}
 \centering
 \scalebox{.95}{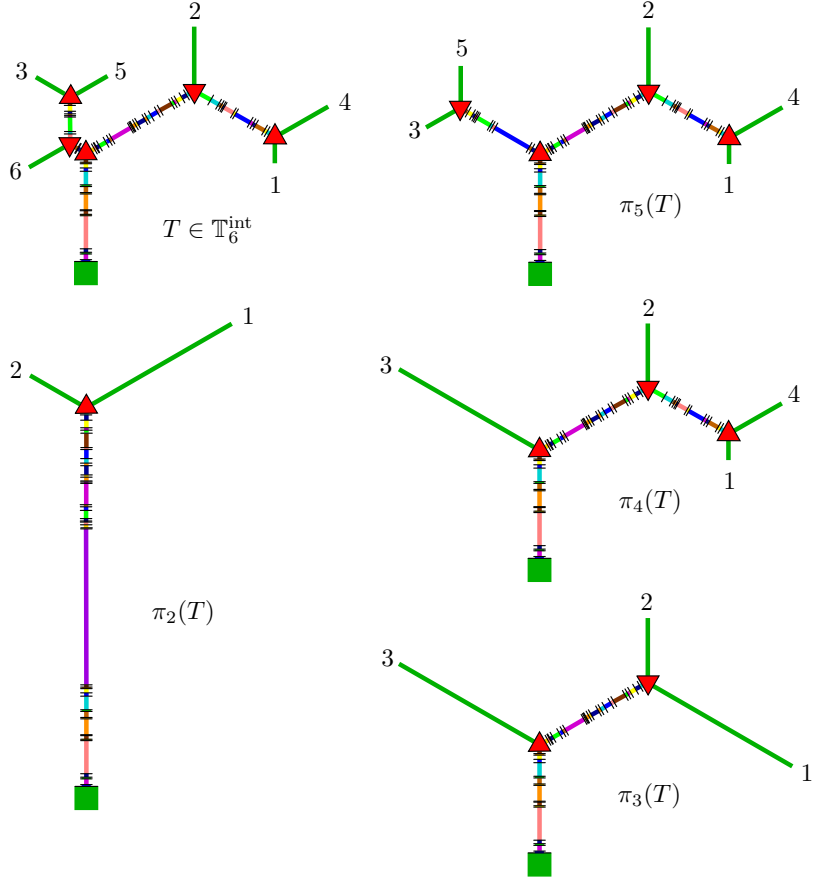}
 \caption{Projections of a $6$-tree.\label{fig:k_tree_proj}}
\end{figure}

These projections are illustrated in Figure \ref{fig:k_tree_proj}. 
Note how, in that example, in passing from $T$ to $\pi_5(T)$, the condition of item \ref{item:proj:merge} of the above definition applies, whereas in passing from $\pi_5(T)$ to $\pi_4(T)$, item \ref{item:proj:add} applies. By comparing these maps to the label insertion operator introduced prior to Definition \ref{def:resamp_1}, we see that the projection neatly undoes label insertion: if $j\notin A\subset\BN$, then
\begin{equation}
\label{eq:projection_insertion}
 \pi_{-j}(T\oplus (\ell,j,U)) = T \quad \text{for any }T\in\bTInt_A,\ \ell\in\block(T),\ U\in \TInt_{2,1}.
\end{equation}

It is easily seen from \cite[Lemma 2.4]{Paper1-1} that these are weak contraction maps:
\begin{equation*}
 d_{\bT}(\pi_{k}(T),\pi_{k}(T^*)) \le d_{\bT}(T,T^*) \quad \text{for }T\in\bTInt_A,\ T^*\in\bTInt_B,
\end{equation*}
where $d_{\bT}$ is defined in \eqref{eq:ktree:metric_1}-\eqref{eq:ktree:metric_2}.

\begin{theorem}\label{thm:consistency}
 \begin{enumerate}[label=(\roman*),ref=(\roman*)]
  \item \label{item:cnst:nonresamp}
    Let $2\le j<k$. For $(\cT^y,\,y\ge0)$ any (self-similar) non-resampling $k$-tree evolution, $(\pi_j(\cT^y),\,y\ge0)$ is a non-resampling $j$-tree evolution.
  \item \label{item:cnst:resamp}
    If $(\cT^y,\,y\ge0)$ is a resampling $k$-tree evolution and $\cT^0$ satisfies
   \begin{equation}\label{eq:cnst:init}
    \bE [\varphi(\cT^0)] = \int_{\TInt_{j+1}}\Lambda_{j+1,[j]}(T_j,dT_{j+1})\cdots\int_{\TInt_k}\Lambda_{k,[k-1]}(T_{k-1},dT_k)\varphi(T_k)
   \end{equation}
   for some $T_j\in \TInt_{j}$, then $(\pi_j(\cT^y),\,y\ge0)$ is a resampling $j$-tree evolution.
  \item \label{item:cnst:dePoi} These same results hold for unit-mass versions of these processes.
 \end{enumerate}
\end{theorem}

For $j\ge1$, we say that $(T_k,k\ge j)$ is a consistent family of $k$-trees if $\pi_{k-1}(T_k)=T_{k-1}$ for all $k > j$. We say that a family of $k$-tree evolutions $(\cT_k^y,y\ge 0)$, $k\ge j$, is consistent if $(\cT_k^y,k\ge j)$ is consistent for each $y\ge 0$. This next result follows from Theorem \ref{thm:consistency} by the Kolmogorov consistency theorem.

\begin{corollary}\label{cor:consistent_fam}
 \begin{enumerate}[label=(\roman*), ref=(\roman*)]
    \item For every consistent family $T_k\in\TInt_{k}$, $k\ge 1$, there are consistent families of non-resampling $k$-tree evolutions $(\cT_k^{y},y\ge 0)$, $k\ge 1$, with $\cT_k^0 = T_k$ for each $k$.\label{item:cnstfam:nonresamp}
    \item For any fixed $j\ge 1$ and $T\in\TInt_{j}$, there exists a consistent family of resampling $k$-tree evolutions $(\cT_k^y,y\ge0)$, $k\ge j$, $\cT_j^0 = T$. If $j=1$, each process in this family is pseudo-stationary, as in Proposition \ref{prop:pseudo:resamp}.\label{item:cnstfam:resamp}
    \item Assertions \ref{item:cnstfam:nonresamp} and \ref{item:cnst:resamp} hold for unit-mass versions of these processes; in particular, there is a consistent family of stationary unit-mass resampling $k$-tree evolutions, $(\widebar{\cT}_k^{u},u\ge 0)$, $k\ge 1$.\label{item:cnstfam:dePoi}
  \end{enumerate}
\end{corollary}

We note one more consistency result. For any $B\subseteq \bN$ with $A\cap B\neq\emptyset$ we define $\pi_B\colon \bTInt_A\to\bTInt_{A\cap B}$ analogously to $\pi_k$, to be the composition of projection maps $\pi_{-j}$ dropping each successive label $j\in A\setminus B$.

\begin{proposition}
 \label{prop:resamp_to_non}
 Suppose $(\cT^y,y\ge0)$ is a resampling $k$-tree evolution. Then there exists a process $((A_y,B_y,\sigma_y),y\ge0)$ that is constant between degeneration times, where for each $y$, $\sigma_y\colon A_y\to B_y$ is a bijection between two subsets of $[k]$, with the property that $\sigma_y\circ\pi_{A_y}(\cT^y)\in\TInt_{B^y}$, $y\ge0$, is a non-resampling $k$-tree evolution.
\end{proposition}

In \cite[Theorem 2]{Paper2} we proved the discrete analogue to Theorem \ref{thm:consistency}\ref{item:cnst:resamp} for the label-swapping variant of Aldous's Markov chain on cladograms described in Definition \ref{def:modified_AC}. Our approach was to find an intermediate process in between the ``upper'' process -- in the present setting, $(\cT^y_{k+1},y\ge0)$ -- and its projection, and give a two-step proof. In the first step, we used the Rogers--Pitman intertwining criterion to show that the intermediate process is Markovian. In the second step, we used Dynkin's criterion to identify the lower, projected process. We refer to Appendix \ref{sec:DynkinIntertwining} for a discussion of these two criteria. 
%


We take an approach inspired by this strategy. As in \cite{Paper2}, neither Dynkin's criterion nor the intertwining criterion holds between the resampling $(k\!+\!1)$- and $k$-tree evolutions. Due to the obstacles presented by degeneration times, we prefer coupling arguments over appeals to Dynkin's criterion between times at which a lower label ($\le k$) is dropped in degeneration and has to resample. At those degeneration times, we argue again by a combination of intertwining and Dynkin arguments connecting the $(k\!+\!1)$-tree, an intermediate object, and the $k$-tree.

This intermediate process and associated intertwining property are introduced in Section \ref{sec:const:intertwin}. In Section \ref{sec:proj_at_degen}, we establish some lemmas about projections of degenerate trees, including the aforementioned Dynkin arguments. Then, in Section \ref{sec:const:Dynkin1}, we put these pieces together with our coupling arguments to show that the consistency of Theorem \ref{thm:consistency}\ref{item:cnst:resamp} holds up until the accumulation time $D_\infty$ of degeneration times. In Section \ref{sec:non_acc}, we use this partial result to prove Proposition \ref{prop:non_accumulation}, which states that $D_\infty$ is a.s.\ the time at which the total mass process converges to 0. This completes the proofs of Theorems \ref{thm:total_mass} and \ref{thm:consistency}\ref{item:cnst:resamp}. The remaining results of Theorem \ref{thm:consistency} and Proposition \ref{prop:resamp_to_non} are then proved in Section \ref{sec:const:other}, thereby completing the proof of Theorem \ref{thm:intro:k_tree}. We conclude in Section \ref{sec:partial_resamp} with a definition and consistency results for partially resampling $k$-tree evolutions.

\section{Intermediate process intertwined below a resampling $(k\!+\!1)$-tree}\label{sec:const:intertwin}

\begin{figure}[t]
 \centering
 \scalebox{.965}{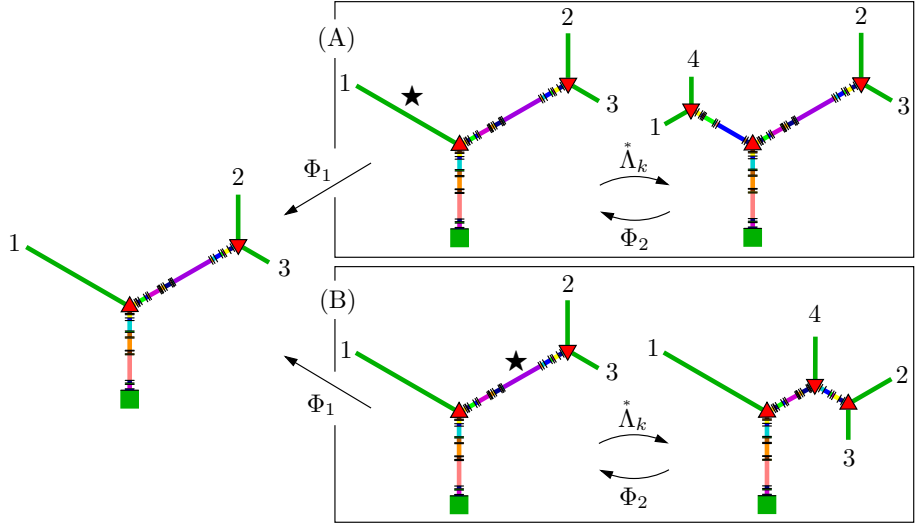}
 \caption{Two marked $k$-trees ($k=3$) based on the same $k$-tree, with marked blocks indicated with a ``$\bigstar$.'' In example (A), the marking is on a leaf block. Then $\Ast\Lambda_k$ splits the marked block into a Brownian reduced 2-tree. In (B), the marking is on an internal block. Then the kernel $\Ast\Lambda_k$ inserts label $k\!+\!1$ into the block. \label{fig:mark_k_tree}}
\end{figure}

Let $k\ge 1$. We define marked $k$-trees as $k$-trees with one block of the tree ``marked.'' In particular, we are interested in projecting from $(k\!+\!1)$-trees and marking the block of the resulting $k$-tree into which label $k\!+\!1$ must be inserted to recover the $(k\!+\!1)$-tree from the $k$-tree. See Figure \ref{fig:mark_k_tree}. 

\begin{definition}\label{def:mark}
 Let $k\ge 1$. We define the set of \emph{marked $k$-trees}
 \begin{equation}\label{eq:marked}
 \begin{split}
  \bTMarkk &:= \Big\{ (T,\ell) \;|\; T\in\bTInt_{k}\setminus\{0\},\ \ell \in \overline{\block}(T)\Big\}\cup\{0\}\text{, where}\\
  \overline{\block}(T) &:= \block(T)\cup\left\{(F,a,a)\colon F\in\ft,\,a\in \left(\big[0,\|\beta_F\|\big]\setminus\bigcup\nolimits_{V\in\beta_F}V\right)\right\}
 \end{split}
 \end{equation}
 for $T = (\ft,(x_j,j\in [k]),(\beta_E,E\in\ft))$. 
 We view marked $k$-trees as intermediate objects between $(k\!+\!1)$-trees and $k$-trees, via a pair of projection maps. First, $\phi_1\colon \bTMarkk\to \bTInt_{k}$ is the projection $\phi_1(T,\ell) = T$. The map $\phi_2\colon \bTInt_{k+1}\to\bTMarkk$ is illustrated in Figure \ref{fig:mark_k_tree} and defined as follows.
 \begin{enumerate}[label=(\roman*),ref=(\roman*)]
  \item If in $T\in\bTInt_{k+1}$ we have $\longparent{\{k\!+\!1\}} = \{j,k+1\}$ for some $j\in [k]$, then $\phi_2(T) = (\pi_{k}(T),j)$. This is the case in Figure \ref{fig:mark_k_tree}(A).
  \item Otherwise, if $E = \longparent{\{k\!+\!1\}}$ is \emph{not} a type-2 edge, then recall part \ref{item:proj:merge} of Definition \ref{def:proj} of $\pi_k(T)$, in which the interval partitions marking the parent edge $E$ and sibling edge $F := E\setminus\{k\!+\!1\}$ are combined with the top mass $x_{k+1}$ to form the partition $\beta'_F = \beta_F\concat (0,x_{k+1})\concat\beta_{E}$ marking $F$ in the projected tree. In this case we define $\phi_2(T) = \big(\pi_{k}(T),(F,\|\beta_F\|,\|\beta_F\|+x_{k+1})\big)$, where the marked block is the block in $\pi_{k}(T)$ corresponding to the top mass $x_{k+1}$ in $T$. This is the case in Figure \ref{fig:mark_k_tree}(B).
 \end{enumerate} 
 We also define a stochastic kernel from $\bTMarkk$ to $\bTInt_{k+1}$. Recall the label insertion operator, $\oplus$, of Section \ref{sec:resamp_def}. Let $\Ast\Lambda_k$ denote the kernel from $\bTMarkk$ to $\bTInt_{k+1}$ that associates with each $(T,\ell)\in \bTMarkk$ the law of $T\oplus (\ell,k\!+\!1,U)$, where $U\sim Q$ is a Brownian reduced $2$-tree of unit mass.
 
 We adopt the convention that $\phi_2(0) = \phi_1(0) = 0$ and $\Ast\Lambda_k(0,\cdot\,) = \delta_0(\,\cdot\,)$.
 
 We denote corresponding Markov kernels by $\Phi_1$ and $\Phi_2$:
 $$\Phi_1((T,\ell),\cdot\,) = \delta_T(\,\cdot\,),\qquad \Phi_2(T,\cdot\,) = \delta_{\phi_2(T)}(\,\cdot\,).$$
\end{definition}

The term in \eqref{eq:marked} in which we allow a marking $\ell$ of the form $(F,a,a)$ allows the description of a $(k\!+\!1)$-tree in which the top mass $x_{k+1}$ equals $0$ and sits in a type-1 compound. The special case $a=\|\beta_F\|$ corresponds to this type-1 compound being degenerate, with a null interval partition below the zero mass leaf component. We write $\TMarkk^{\rm int}:=\phi_2(\TInt_{k+1})$ and $\tdTMarkk:=\phi_2(\tdTInt_{k+1})$ for the spaces of marked trees with respectively no degenerate labels or at most one degenerate label, which could be label $k+1$, as discussed above.

We may think of the resampling kernel $\Lambda_{k+1,[k]}$ of Section \ref{sec:resamp_def} as representing a two-step transition in which a block is first selected at random and then, if a leaf block was chosen, it is split into a scaled Brownian reduced 2-tree. Then $\Ast\Lambda_k$ represents the second step: 
for $(T,\ell)\in\TMarkk^{\rm int}$ with $\|\ell\|>0$, 
\begin{equation}\label{eq:resamp_vs_mark}
 \Ast\Lambda_k\left(\left(T,\ell\right),\cdot\,\right) = \Lambda_{k+1,[k]}\left(T,\cdot\ \middle|\ k\!+\!1\text{ is inserted into }\ell\right),
\end{equation}

In Appendix \ref{sec:marked_tree_metric}, we introduce a natural metric $d_{\Ast\bT}$ on $\bTMarkk$ that possesses the following properties.

\begin{lemma}\label{lem:Lstar_cont}
 The projection maps $\phi_1\colon\big(\bigcup_{k\ge 1}\bTInt_{k+1},d_\bT\big)\rightarrow\big(\bigcup_{k\ge 1}\bTMarkk,d_{\Ast\bT}\big)$ and 
                $\phi_2\colon\big(\bigcup_{k\ge 1}\bTMarkk,d_{\Ast\bT}\big)\rightarrow\big(\bigcup_{k\ge 1}\bTInt_k,d_\bT\big)$ are continuous. The kernel $\Ast\Lambda_k$ is weakly continuous in its first coordinate.
\end{lemma}

When composing stochastic kernels, we adopt the standard convention that sequential transitions are ordered from left to right:
\begin{equation}\label{eq:composition_of_kernels}
 \int PQ(x,dz)f(z) = \int P(x,dy)\int Q(y,dz)f(z).
\end{equation}
We discuss Rogers and Pitman's \cite{RogersPitman} notion of intertwining in a continuous-time regime in Appendix \ref{sec:DynkinIntertwining}; for now, we note the property in the discrete-time regime.

\begin{definition}
\label{def:intertwining:MC}
 Suppose $P$ is a Markov kernel on a state space $(S,\cS)$, $\phi\colon S\to R$ is a surjective measurable map to $(R,\cR)$, and $\Lambda\colon R\times\cS\to [0,1]$ is a stochastic kernel. Let $\Phi\colon S\times\cR\to [0,1]$ denote the kernel $\Phi(s,\cdot\,) = \delta_{\phi(s)}(\,\cdot\,)$. Let $Q := \Lambda P\Phi$. We say that \emph{$Q$ is intertwined below $P$} (via $(\phi,\Lambda)$) if
 \begin{enumerate}[label=(I\arabic*),ref=(I\arabic*)]
  \item \label{item:intertwin:up:MC} $\Lambda\Phi$ is the identity kernel on $(R,\cR)$ and
  \item \label{item:intertwining:MC} $\Lambda P = Q\Lambda$.
 \end{enumerate}
\end{definition}

\begin{proposition}[Theorem 2 in \cite{RogersPitman}, discrete-time regime]\label{RogePitmdiscr}
 Let $P$, $\phi$, $\Lambda$, and $Q$ be as in Definition \ref{def:intertwining:MC}, with $Q$ intertwined below $P$ via $(\phi,\Lambda)$. If $(X_n,\,n\ge0)$ is a discrete-time Markov process on $S$ with transition kernel $P$ and
 \begin{enumerate}[start=3,label=(I\arabic*),ref=(I\arabic*)]
  \item[\rm (I3)] \label{item:intertwin:init:MC} $X_0 \sim \int_R \Lambda(y,\cdot\,)\mu(dy)$ for some probability measure $\mu$ on $R$,
 \end{enumerate}
 then $Y_n := \phi(X_n)$, $n\ge0$, is a discrete-time Markov process on $R$ with $Y_0 \sim \mu$. We say that this process is \emph{intertwined below $(X_n,\,n\ge0)$} via $\Lambda$.
\end{proposition}

For $(T,\ell)\in\TMarkk^{\rm int}$, let $\big(\cT_{k+1}^y,y\ge0\big)$ denote a resampling $(k\!+\!1)$-tree evolution with initial distribution $\Ast\Lambda_k((T,\ell),\cdot\,)$. Let $\Ast\cT_{k}^y = (\cT_{k}^y,\ell^y) := \phi_2(\cT^y_{k+1})$, $y\ge0$. 

\begin{remark}
 In fact, $\big(\Ast\cT_{k}^y,\,y\ge0\big)$ is a (continuous-time) Markov process intertwined below $\big(\cT_{k+1}^y,y\ge0\big)$, in the sense described in Appendix \ref{sec:DynkinIntertwining}. However, as we do not require this full result, we only prove a pair of partial results in this direction.
\end{remark}

Recall from Section \ref{sec:killed_def} our language around tree shapes and from Section \ref{sec:non_resamp_def} the definition of the swap-and-reduce map, $\varrho$. When a label $i$ in a resampling or non-resampling $(k\!+\!1)$-tree evolution degenerates, it swaps places with $\max\{i,a,b\}$, where $a$ is the least label descended from the sibling of $i$ and $b$ is the least label descended from its uncle. Since $k\!+\!1$ is the greatest label in the tree, there are three cases in which it will resample:
\begin{enumerate}[label=(D\arabic*), ref= (D\arabic*), topsep=4pt, itemsep=3pt]
 \item\label{case:degen:type2} $k\!+\!1$ belongs to a type-2 compound, and either it (in the case $k\!+\!1=i$) or its sibling (in the case $k\!+\!1=a$) causes degeneration;
 \item\label{case:degen:self} $k\!+\!1$ belongs to a type-1 compound and causes degeneration, so $k\!+\!1=i$; or
 \item\label{case:degen:nephew} $k\!+\!1=b$, as leaf $k\!+\!1$ belongs to a type-1 compound and its sibling in the tree shape is an internal edge that belongs to a type-1 or type-2 compound that degenerates.
\end{enumerate}

Let $D_1,D_2,\ldots$ denote the sequence of degeneration times of $\big(\cT^y_{k+1},\,y\ge0\big)$, with $D_1^{\le k},D_2^{\le k},\ldots$ the subsequence of degeneration times at which a label other than $k+1$ is dropped and resampled.

\begin{proposition}\label{prop:inter:skel_1}
 $\Big(\Ast\cT^{D_n}_{k}\cf\{D_n<D^{\le k}_1\},\,n\ge0\Big)$ is a discrete-time Markov process intertwined below $\Big(\cT^{D_n}_{k+1}\cf\{D_n<D^{\le k}_1\},\,n\ge0\Big)$ via $\big(\phi_1,\Ast\Lambda_k\big)$.
\end{proposition}

\begin{proof} First note that the upper process is a discrete-time Markov process by construction. 
 Conditions \ref{item:intertwin:up:MC} and (I3) are satisfied by definition of the relevant kernels and our setup. By Proposition \ref{RogePitmdiscr} and as 
 discussed in Appendix \ref{sec:DynkinIntertwining}, it now suffices to check that for all $(T,\ell)\in\TMarkk^{\rm int}$, the r.c.d.\ of $\cT^{D_1}_{k+1}\cf\{D_1<D^{\le k}_1\}$ is $\Ast\Lambda_k(\Ast\cT^{D_1}_k\cf\{D_1<D^{\le 1}_1\},\cdot)$, if $\cT^{0}_{k+1}$ has distribution $\Ast\Lambda_k((T,\ell),\cdot)$.
For brevity, let $\cT^{\circ,y}_{k+1} := \cT^{y}_{k+1}\cf\{y<D^{\le k}_1\}$ and $\Ast\cT^{\circ,y}_{k} := \Ast\cT^{y}_{k}\cf\{y<D^{\le k}_1\}$, $y\ge0$.
 
 Let $(T,\ell)\in\TMarkk^{\rm int}$ and let $\cT^{0}_{k+1}$ have distribution $\Ast\Lambda_k((T,\ell),\cdot)$. The claimed r.c.d.\ holds trivially where $\Ast\cT^{\circ,D_1}_{k} = 0$ and where $\Ast\cT^{\circ,D_n}_{k} \in \TMarkk^{\rm int}\setminus (\{0\}\cup(\TInt_k\times [k]))$, as in both cases, 
$\Ast\Lambda_k$ maps the marked $k$-tree to a point mass at the unique $(k\!+\!1)$-tree that projects down to it via $\phi_2$. The latter of those cases is that depicted in Figure \ref{fig:mark_k_tree}(B). This leaves the only non-trivial case, where $\Ast\cT_k^{\circ,D_1}\in\bT_k^{\rm int}\times[k]$, i.e.\ when $k\!+\!1$ is the label to resample at $D_1$ and the mark is put into a leaf block $i\in[k]$, forming a type-2 compound, i.e.\ $D_1 < D^{\le k}_1$ and $\ell^{D_1} = i\in [k]$. Hence, when $\Ast\cT^{\circ,D_1}_{k} \in \TInt_k\times [k]$, the claimed r.c.d.\ follows from \eqref{eq:resamp_vs_mark} and the definition of resampling in \eqref{eq:resampling_kernel:def}.
\end{proof}

\begin{lemma}\label{lem:inter:preswap}
 Given the event $\{D_1^{\le k}<D_\infty\}$ and the path $\big(\Ast\cT_k^y,\, y\in [0,D^{\le k}_1)\big)$, the tree prior to resampling at this time, $\cT_{k+1}^{D^{\le k}_1-}$, has regular conditional distribution $\Ast\Lambda_k\big(\Ast\cT_k^{D^{\le k}_1-},\cdot\,\big)$.
\end{lemma}

\begin{proof} First note that $\mathbb{P}\big(\Ast\cT_k^{D^{\le k}_1-} \in \tdTMarkk\setminus\TMarkk^{\rm int}\big)=1$ since $D_1^{\le k}$ is a degeneration time.
 As in the preceding proof, the claimed r.c.d.\ holds trivially where $\Ast\cT_k^{D^{\le k}_1-} \in \tdTMarkk\setminus \big(\TMarkk^{\rm int}\cup(\tdTInt_k\times [k])\big)$. It remains to check the case where $\Ast\cT_k^{D^{\le k}_1-} \!\in\! \big(\tdTInt_k\setminus\TInt_k\big)\times [k]$, in which the marked block at time $D_1^{\le k}-$ is a top mass.
 
 Fix $i\in [k]$, $n\ge0$, and consider the event $E_{i,n} := \big\{\ell^{D_n} = i, D^{\le k}_1=D_{n+1}\big\}$. On $E_{i,n}$, indeed on $\big\{\ell^{D_n} = i, D^{\le k}_1>D_n\big\}$, labels $i$ and $k+1$ are in a type-2 compound $\Gamma_{i,k+1}^y := \big(m^{y}_i,m^{y}_{k+1},\beta^{y}_{\{i,k+1\}}\big)$, $y\in [D_n,D_{n+1})$. We denote the associated total mass and normalized 2-tree by $M_{i,k+1}^y := \big\|\Gamma_{i,k+1}^y\big\|$ and $\bar\Gamma_{i,k+1}^y := \Gamma_{i,k+1}^y / M_{i,k+1}^y$, $y\in [D_n,D_{n+1})$. 
 Recall that $\Ast\Lambda_k$ acts on a $k$-tree with a marked leaf block by splitting that leaf block into a suitably scaled Brownian reduced 2-tree. 
 Thus, it will suffice to prove that, given the event $E_{i,n}$ for any choice of $i$ and $n$, the tree $\bar\Gamma_{i,k+1}^{D^{\le k}_1-}$ is a Brownian reduced 2-tree and is  conditionally independent of $\big(\Ast\cT_{k}^{y},\,y\in \big[0,D^{\le k}_1\big)\big)$.
 
 Let us first prove the claimed conditional distribution. By Proposition \ref{prop:inter:skel_1}, $\bar\Gamma^{D_n}_{i,k+1}$ is a Brownian reduced 2-tree independent of $M_{i,k+1}^y$, conditionally given $\big\{\ell^{D_n}=i,D_1^{\le k}>D_n\big\}$.  
Further conditioning on $\big\{D_1^{\le k}=D_{n+1}\big\}$ is conditioning the type-2 evolution $\big(\Gamma_{i,k+1}^{D_n+y}\big)$ to not degenerate prior to some other (independent) type-1 or type-2 compound in $\big(\cT_{k+1}^{D_n+y},\,y\ge 0\big)$. By the strong pseudo-stationarity observed in Lemma \ref{pretotal}\ref{item:pretotal:1}, $\bar\Gamma_{i,k+1}^{D^{\le k}_1-}$ is again a Brownian reduced 2-tree, conditionally given $E_{i,n}$.
 
 Now, we show the claimed conditional independence. 
 Lemma \ref{pretotal}\ref{item:pretotal:1} observes further that $\bar\Gamma_{i,k+1}^{D^{\le k}_1-}$ is 
 conditionally independent of its past total mass process $\big(M_{i,k+1}^{D_n+y},\,y\in \big[0,D^{\le k}_1-D_n\big)\big)$ given $E_{i,n}$. By the independence of the constituent type-0/1/2 compounds that make up the evolving tree in between degenerations, $\bar\Gamma_{i,k+1}^{D^{\le k}_1-}$ is furthermore conditionally independent of $\big(\Ast\cT_{k}^{D_n+y},\,y\in \big[0,D^{\le k}_1-D_n\big)\big)$ given $E_{i,n}$. Finally, if $n\ge1$ then by the iterative construction of resampling $(k\!+\!1)$-tree evolutions in Definition \ref{def:resamp_1}, $\bar\Gamma_{i,k+1}^{D^{\le k}_1-}$ is conditionally independent of $\big(\cT^y_{k+1},\,y\in [0,D_n)\big)$ given $E_{i,n}$ and $\Ast\cT_k^{D_n}$. Thus, $\bar\Gamma_{i,k+1}^{D^{\le k}_1-}$ is conditionally independent of $\big(\Ast\cT_{k}^{y},\,y\in \big[D_n,D^{\le k}_1\big)\big)$ given $E_{i,n}$, as desired.
\end{proof}

\section{Projections of degenerate trees}
\label{sec:proj_at_degen}

In this section we prove three lemmas regarding projection maps applied to degenerate trees, in further preparation to prove the projective consistency of resampling $k$-tree evolutions. 
Recall from Section \ref{sec:non_resamp_def} that for a tree with one degenerate label, $T\in \tdTInt_{k}\setminus\TInt_{k}$, $I(T)$ denotes the degenerate label and $J(T)$ denotes the label dropped when applying the swap-and-reduce map: $\varrho(T)\in \TInt_{[k]\setminus\{J(T)\}}$.
 


\begin{figure}[t]
 \centering
 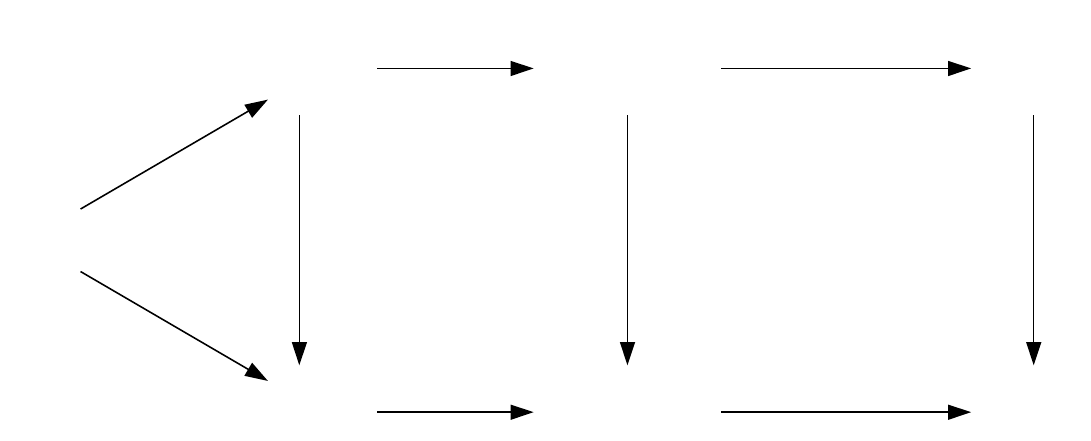
 \caption{Here, $J_n$ denotes $J\big(\cT^{D^{\le k}_n-}_{k}\big)$. Lemma \ref{lem:inter:preswap} asserts that the conditional law of $\cT^{D^{\le k}_n-}_{k+1}$ (upper left) given $\Ast{\cT}^{D^{\le k}_n-}_k$ (left) is as in the upper left arrow in this diagram. By construction, the conditional law of $\cT^{D^{\le k}_n}_k$ (in the lower right) given $\cT^{D^{\le k}_n-}_{k+1}$ (upper left) is as in the upper path. Lemma \ref{lem:projection_at_degen} asserts that the left box of the diagram commutes; Lemma \ref{lem:Dynkin:preswap} claims that the rest of the diagram commutes.\label{fig:preswap}}
\end{figure}

\begin{lemma}
\label{lem:projection_at_degen}
 For a degenerate $(k\!+\!1)$-tree $T\in \tdTInt_{k+1}\setminus\TInt_{k+1}$, if $J(T)\in [k]$, then $\pi_k(T)$ is also degenerate: $\pi_k(T)\in \tdTInt_{k}\setminus\TInt_{k}$ with
 \begin{equation}
 \label{eq:projection_at_degen}
  I(T) = I(\pi_k(T)),\ \ J(T) = J(\pi_k(T)),\ \ \text{and}\ \ \pi_k(\varrho(T))= \varrho(\pi_k(T)).
 \end{equation}
\end{lemma}

\begin{proof}
 Let $T$ be as in the statement of the lemma. We denote its coordinates by $T = (\ft,(x_i,\,i\in [k]),(\beta_E,\,E\in\ft))$. Let $i := I(T)$, so $x_i = 0 = \big\|\beta_{\parent{\{i\}}}\big\|$. 
 By Definition \ref{def:proj} of $\pi_k$, this projection merges the masses and partitions on leaf $k+1$, its parent, and its sibling. Specifically, if the parent edge is of type 1, the mass of $k+1$ and the partitions of the parent and sibling edges form a new edge partition, and if the parent edge is of type 2, the masses of $k+1$, the entire parent partition and the sibling top mass merge to form a new top mass. Neither leaf $i$ nor its parent edge can be the sibling or parent of $k+1$, as if they were, then $J(T)$ would equal $k+1$, in violation of our hypothesis. Thus, weight $x_i$ and partition $\beta_{\parent{\{i\}}}$ are unchanged in this projection. We further remark that each mass or partition in $T$ either is preserved or is involved in the merger described above, so no additional compound can become degenerate in $\pi_k(T)$. Hence, this projected tree has a unique degenerate label $I(\pi_k(T))= i = I(T)$, and $\pi_k(T)\in \tdTInt_{k}\setminus\TInt_{k}$. 
 
 Recall that $J(T) = \max\{i,\min(A),\min(B)\}$, where $A$ and $B$ denote the label sets on the sibling and uncle of leaf $i$, respectively. By our hypothesis, neither $A$ nor $B$ can equal $\{k+1\}$. The projection can only possibly change these label sets by removing label $k+1$, which does not affect $\min(A)$ or $\min(B)$. This proves that $J(\pi_k(T)) = J(T)=:j$.
 
 Finally, the swap-and-reduce map $\varrho$ swaps the places of labels $i$ and $j$, then applies the projection map $\pi_{-j}$ to remove the weightless leaf $j$ and its equally weightless parent edge. As noted in Definition \ref{def:proj}, the projection maps $\pi_{-j}$ and $\pi_{-(k+1)}$ commute. Moreover, $\pi_{-(k+1)}$ clearly commutes with the operation of swapping two lower labels. This proves that $\pi_k(\varrho(T)) = \varrho(\pi_k(T))$.
\end{proof}


\begin{lemma}
\label{lem:Dynkin:preswap}
 Fix $j\in [k]$ and $(T,\ell)\in\tdTMarkk\setminus\TMarkk^{\textnormal{int}}$ with $J(T) = j$. Let $T^\prime\sim \Ast\Lambda_k((T,\ell),\cdot\,)$. Then the following two $k$-trees have the same distribution:
  \begin{itemize}\item $\pi_k(T^{\prime\prime})$, where $T^{\prime\prime}$ has conditional law $\Lambda_{j,[k+1]\setminus\{j\}}(\varrho(T^\prime),\cdot\,)$ given $T^\prime$, and
   \item $T^{\prime\prime\prime}\sim\Lambda_{j,[k]\setminus\{j\}}(\varrho(T),\cdot\,)$.
  \end{itemize}
\end{lemma}

When applied to $(T,\ell)=\Ast\cT_k^{D_n^{\le k}-}$, this has the flavor of a Dynkin's criterion companion to the intertwining-like assertion of Lemma \ref{lem:inter:preswap}, claiming that when a label resamples, the added information in the marked tree, namely the marked block, does not inform the transition of the projected $k$-tree evolution $\cT^y_k = \pi_k\big(\cT^y_{k+1}\big) = \phi_1\big(\Ast\cT^y_k\big)$. This assertion and its relationship with Lemma \ref{lem:inter:preswap} are described via a commutative diagram of stochastic kernels in Figure \ref{fig:preswap}.

\begin{proof}
 First we observe that, as noted in Definition \ref{def:mark} of $\Ast\Lambda_k$, $\pi_k(T^\prime) = T$. Thus, by Lemma \ref{lem:projection_at_degen}, $\pi_k(\varrho(T^\prime)) = \varrho(T)$.
 
 This lemma is trivial in the case $\ell\notin [k]$ that the marked block is internal: in that case, the kernel $\Ast\Lambda_k$ acts trivially, mapping $(T,\ell)$ to a Dirac point mass at the unique deterministic $(k\!+\!1)$-tree $T_1$ that satisfies $\phi_2(T_1) = (T,\ell)$, as in Figure \ref{fig:mark_k_tree}(B), and there is a natural tree-structure- and block-mass-preserving bijection between the blocks of $T$ and those of $T_1 = T^\prime$, allowing us to couple $\Lambda_{j,[k+1]\setminus\{j\}}(\varrho(T_1),\cdot\,)$ with $\Lambda_{j,[k]\setminus\{j\}}(\varrho(T),\cdot\,)$ so that the $\pi_k$-projection of the former equals the latter.
 
 Henceforth, we assume that $\ell = i\in [k]$, so labels $i$ and $k\!+\!1$ are in a type-2 compound in $T^\prime$. Let $H$ denote the event that the kernel $\Lambda_{j,[k+1]\setminus\{j\}}$ inserts label $j$ somewhere into this type-2 compound. From the standpoint of the $k$-tree $T$, this is the event that label $j$ is inserted into the marked leaf block labeled $i$. Again, the assertion is trivial on the event $H^c$, as then there is a tree-structure- and block-mass-preserving bijection between the \emph{remaining} blocks of the trees, i.e.\ the unmarked blocks in $T$ and the blocks outside of the type-2 compound containing $i$ and $k\!+\!1$ in $T^\prime$. Thus it remains only to prove the assertion conditional on the event $H$.
 
 This event $H$ has probability $x_i / \|T\|$, where $x_i$ is the top mass labeled $i$ in $T$. Moreover, $H$ is independent of $T^\prime$: it does not depend on the normalized ``internal structure'' $U$ of the type-2 compound containing $i$ and $k+1$, which is the only random part of $T^\prime$. 
 By definition of $\Ast\Lambda_k$, this type-2 compound in $T^\prime$ is distributed as a Brownian reduced 2-tree scaled to have mass $x_i$. By \eqref{eq:B_ktree_resamp} and the exchangeability of labels in Brownian reduced $m$-trees noted in Proposition \ref{prop:B_ktree}, this means that after inserting label $j$, blocks $i$, $j$, and $k\!+\!1$ in $T^{\prime\prime}$, along with the partitions marking their parent edges, comprise a Brownian reduced 3-tree of total mass $x_i$. Thus, the $\pi_{-(k+1)}$-projection of this 3-tree is another Brownian reduced 2-tree of mass $x_i$, with leaf labels $i$ and $j$; i.e.\ $\pi_k(T^{\prime\prime})$ is distributed as it would be under $\Lambda_{j,[k]\setminus\{j\}}(T,\cdot)$ conditioned on label $j$ resampling into block $i$.
\end{proof}

\section{Consistent resampling $k$-tree evolutions}
\label{sec:const:Dynkin1}

This section is devoted to proving the following proposition, that the consistency of Theorem \ref{thm:consistency}\ref{item:cnst:resamp} holds up until time $D_\infty$.

\begin{proposition}\label{prop:consistency_0}
 Fix $T\in\TInt_{k}$. If $(\cT^y_{k+1},y\ge0)$ is a resampling $(k\!+\!1)$-tree evolution with $\cT^0_{k+1}\sim \Lambda_{k+1,[k]}(T,\cdot\,)$, then $\cT^y_k := \pi_k\big(\cT^y_{k+1}\big)$, $y\ge0$, evolves as a resampling $k$-tree evolution prior to the limit $D_\infty$ of degeneration times in $\big(\cT_{k+1}^y\big)$.
\end{proposition}

Since it is not clear at this stage whether $D_\infty$ is a stopping time in the natural filtration of $(\cT^y_k)$, we will also construct a bigger filtration in which degeneration times
of $(\cT^y_{k+1})$, and $D_\infty$, are stopping times, while $(\cT^y_k)$ is still a resampling $k$-tree evolution with respect to this bigger filtration. We begin with a series of intermediate results starting from killed $(k\!+\!1)$-tree evolutions and successively extending past degeneration times.

\begin{lemma}\label{lem:Dynkin:killed:1}
 For any initial distribution $\mu$ on $\TInt_{k+1}$, it is possible to define a pair of coupled processes such that:
\begin{itemize}
  \item $\big(\cT^y_{k+1},\,y\ge0\big)$ is a killed $(k\!+\!1)$-tree evolution with $\cT^0_{k+1}\sim\mu$;
  \item $\big(\cT^y_{k},\,y\ge0\big)$ is a killed $k$-tree evolution;
  \item $\cT^y_k = \pi_k\big(\cT^y_{k+1}\big)$ for all $y$ less than the degeneration time of $\big(\cT^y_{k+1}\big)$; and
  \item $\big(\cT^y_{k},\,y\ge0\big)$ is strongly Markovian in the filtration generated by both processes.
\end{itemize}
\end{lemma}

\begin{proof}
 We will prove this for a generic fixed initial state $T'\in \TInt_{k+1}$; the extension to general initial distributions follows by mixing. Let $T := \pi_k(T')$. We denote the coordinates of $T$ by $(\ft,(x_i,i\in[k]),(\beta_E,E\in\ft))$. Let $(\Omega_{(0)},\cF_{(0)},\bP_{(0)})$ denote a probability space on which we have defined an independent type-$d$ evolution corresponding to each type-$d$ compound in $T$, for $d=0,1,2$, with the initial state of each evolution equal to the corresponding compound in $T$. We denote the top mass and interval partition evolutions corresponding to each leaf $i$ and each internal edge $E$ by $\big(m_i^y,y\ge0\big)$ and $\big(\beta_E^y,y\ge0\big)$, respectively. Let $D^{\le k}$ denote the minimum of the degeneration times of these type-$d$ evolutions. As in Definition \ref{def:killed_ktree}, $\cT^y_k := \big(\ft,\big(m_j^y,j\in[k]\big), \big(\beta_E^y,E\in\ft\big)\big)$, $y\in [0,D^{\le k})$, with $\cT^y_k := 0$ for $y\ge D^{\le k}$, is a killed $k$-tree evolution. We will extend this to include a construction of a $(k\!+\!1)$-tree evolution in two cases.
 \medskip
 
 Case 1: Leaf $k+1$ is in a type-2 compound in $T'$, $U' \!=\! (x'_i,x'_{k+1},\beta'_{\{i,k+1\}}) \!\in\! \TInt_2$, with some sibling leaf $i\in [k]$. We extend our probability space to $(\Omega_{(1)},\cF_{(1)},\bP_{(1)})$ to include a process $\big(U_{(1)}^y,y\ge0\big)$ so that, given the sub-$\sigma$-algebra of $\cF_{(1)}$ corresponding to $\cF_{(0)}$, it is conditionally distributed as a type-2 evolution with initial state $U'$, conditioned to have total mass evolution $\big\|U_{(1)}^y\big\| = m_i^y$ for $0\le y \le \inf\big\{z\ge 0\colon m_i^z = 0\big\}$. Such a conditional distribution exists as, by Proposition \ref{prop:IPspace:Lusin}, $(\cI,d_{\cI})$ is Lusin, and type-2 evolutions are c\`adl\`ag. 
 If leaf $i$ belongs to a type-1 compound in $T$ (not in $T'$, where it shares a type-2 compound with $k\!+\!1$), then Proposition \ref{prop:012:concat}\ref{item:012concat:BESQ+0} indicates that this total mass process $\big(m_i^y,\,y\ge0\big)$ evolves as a $\besq(-1)$ on this time interval. 
 The same holds if $i$ belongs to a type-2 compound in $T$, by Definition \ref{def:type2:v1} of type-2 evolutions and the symmetry asserted in Lemma \ref{lem:type2_symm}. 
 Thus, by Proposition \ref{prop:012:mass}, after integrating out this conditioning, $\big(U_{(1)}^y,y\ge0\big)$ is a type-2 evolution. We define $\Delta_1$ to be the degeneration time of $\big(U_{(1)}^y,y\ge0\big)$ and set
 \begin{equation}\label{eq:Dynkin:D1}
  D_1 := D^{\le k}\wedge \Delta_1.
 \end{equation}
 Recall the label insertion operator $\oplus$ defined in Section \ref{sec:resamp_def}. We define 
 \begin{equation}\label{eq:Dynkin:insert_ext}
  \cT^y_{k+1} := \cT^y_k\oplus \left(i,\, k\!+\!1,\, U_{(1)}^y/\left\|U_{(1)}^y\right\|\right), \quad y\in [0,D_1),
 \end{equation}
 with $\cT^y_{k+1} := 0$ for $y \ge D_1$. Then this is a killed $(k\!+\!1)$-tree evolution with initial state $T'$, in which the type-2 compound containing label $k\!+\!1$ equals $\big(U_{(1)}^y,y\ge0\big)$. 
 \medskip
 
 Case 2: Leaf $k+1$ belongs to a type-1 compound in $T'$. Then, following Definition \ref{def:proj}\ref{item:proj:merge} of $\pi_k$ in this case, leaf $k+1$ corresponds to a block $(a,b)\in\beta_E$ along some internal edge partition in $T$, if we exclude for a moment the subcase where the leaf mass of $k+1$ vanishes. 
We consider the case that this internal edge $E\in\ft$ belongs to a type-2 compound in $T$, $E = \{i,j\}$ for some $i,j\in [k]$; the other cases can be handled similarly. 

 For the purpose of the following, we denote the type-2 evolution on this compound by $\Gamma_{E}^y := \big(m_i^y,m_j^y,\beta_E^y \big)$, $y\ge0$. 
 Corollary \ref{cor:decomp_ker} notes the existence of a kernel $\kappa$ that takes the path of a type-2 compound and a block in the interval partition component at time zero, and yields the conditional joint law for a type-2 and type-1 evolution, conditioned to concatenate to equal the specified type-2 path, up until a degeneration time. More specifically, let 
 \begin{equation}
 \label{eq:DK:type2_split}
  \left( \Gamma_{(1)}^y, \left(m_{(1)}^y,\beta_{(1)}^y\right),\, y\in \big[0,D_{(1)}\big)\right) \sim
  	\kappa\Big( \big( \big( \Gamma_E^y,y\ge0\big),\,(a,b)\big),\cdot\,\Big).
 \end{equation}
 We recall three properties from Corollary \ref{cor:decomp_ker}:
 \begin{itemize}
  \item up until the degeneration time of $\Gamma^{(1)}$, these two processes concatenate to equal $\Gamma_E^y$, in the sense of \eqref{eq:012concat:2+1};
  \item after mixing over the law of $\big(\Gamma_E^y,y\ge0\big)$, these constituent evolutions are independent; and
  \item the top mass $m_{(1)}^0$ corresponds to the marked block $(a,b)\in \beta_E^0$.
 \end{itemize}
 We extend our probability space to $(\Omega_{(1)},\cF_{(1)},\bP_{(1)})$ to include a pair with this conditional law, as in \eqref{eq:DK:type2_split}, given the sub-$\sigma$-algebra of $\cF_{(1)}$ corresponding to $\cF_{(0)}$.
 
 We define $a_{(1)}^y$ to equal the mass of the interval partition component of $\Gamma_{(1)}^y$ and we set $b_{(1)}^y := a_{(1)}^y + m_{(1)}^y$ for $y\in [0,D_{(1)})$. We take $\Delta_1$ to denote the lesser of the two degeneration times of $\Gamma_{(1)}^y$ and $\left(m_{(1)}^y,\beta_{(1)}^y\right)$, $y\ge0$, and use that to define $D_1$ as in \eqref{eq:Dynkin:D1}. Then we set
 \begin{equation}\label{eq:Dynkin:insert_int}
  \cT^y_{k+1} := \cT^y_k\oplus \left(\left(E,a_{(1)}^y,b_{(1)}^y\right),\, k\!+\!1,\, U \right),\quad y\in [0,D_1),
 \end{equation}
 with $\cT^y_{k+1} := 0$ for $y \ge D_1$, where $U$ is an arbitrary 2-tree, say $(1/2,1/2,\emptyset)$, which, we recall from Section \ref{sec:resamp_def}, is redundant in the label insertion operator when inserting into an internal block. It is easily checked that this all works, too, when the mass of $k+1$ vanishes in $T^\prime$ or at exceptional times, i.e.\ when $m_{(1)}^y=0$.
 \medskip
 
 In each case, the constructed process $\big(\cT^y_{k+1},y\ge0\big)$ is a killed $(k\!+\!1)$-tree evolution with initial state $T'$, satisfying $\pi_k\big(\cT^y_{k+1}\big) = \cT^y_{k}$ for $y\in [0, D_1)$. The claim that $\big(\cT^y_{k},y\ge0\big)$ is a strong Markov process in the filtration generated by itself and $\big(\cT^y_{k+1},y\ge0\big)$ follows from the Proposition \ref{prop:012:mass}, which asserts that the total mass process is strongly Markovian in the filtration generated by a type-2 evolution, and Proposition \ref{prop:012:concat}, which asserts that a type-0/1/2 evolution is strongly Markovian in the filtration generated by two constituent evolutions into which it can be decomposed, as in \eqref{eq:DK:type2_split}.
\end{proof}

\begin{lemma}\label{lem:Dynkin:killed:2}
 We continue in the setting of Lemma \ref{lem:Dynkin:killed:1}, with $D_1$ denoting the degeneration time of $\big(\cT^y_{k+1}\big)$. On the event $\big\{J\big(\cT_{k+1}^{D_1-}\big) = k\!+\!1\big\}$, that label $k + 1$ would be dropped in the first degeneration, we additionally find
 \begin{equation}\label{eq:swapred:k+1}
  \varrho\big(\cT_{k+1}^{D_1-}\big) = \cT_k^{D_1}\ \ \text{a.s.}.
 \end{equation}
\end{lemma}
 

\begin{proof}
 Recall the three cases in which we can have $J\big(\cT_{k+1}^{D_1-}\big) = k\!+\!1$, listed as \ref{case:degen:type2}, \ref{case:degen:self}, and \ref{case:degen:nephew} above Proposition \ref{prop:inter:skel_1}.

 Case \ref{case:degen:type2}: $k+1$ and another label $i$ are in a type-2 compound that degenerates at time $D_1$. In this case, \eqref{eq:swapred:k+1} is clear: on both sides of the formula, this type-2 compound in $\cT_{k+1}^{D_1-}$ is reduced to a single leaf mass with leaf label $i$.
 
 Case \ref{case:degen:self}: $k+1$ is in a type-1 compound that degenerates at time $D_1$. Again, \eqref{eq:swapred:k+1} clear: this type-1 compound has zero mass, and it is contracted away on both sides of the formula.
 
 For the last case, we take up the notation in the proof of Lemma \ref{lem:Dynkin:killed:1}.
 
 Case \ref{case:degen:nephew}: $k+1$ is in a type-1 compound and (one of) its nephew label(s), meaning one of the labels in the process $\big(\Gamma_{(1)}^y\big)$ described in \eqref{eq:DK:type2_split}, causes degeneration at time $D_1$. We will address the case where the sibling edge $E = \longparent{\{k\!+\!1\}}\setminus\{k\!+\!1\}$ is in a type-2 compound, $E = \{i,j\}$, and say label $i$ causes degeneration; the type-1 case is similar. Then, in $\cT_{k+1}^{D_1-}$, block $i$ and edge $\{i,j\}$ both have mass zero; but in $\varrho\big(\cT_{k+1}^{D_1-}\big)$, label $i$ displaces label $k\!+\!1$, and the edge that was formerly $\longparent{\{k\!+\!1\}} = \{i,j,k\!+\!1\}$ gets relabeled as $\{i,j\}$, so that the newly labeled block $i$ has mass $m_{(1)}^{D_1}$ while edge $\{i,j\}$ bears the partition $\beta_{(1)}^{D_1}$. This is consistent with the second line of the formula in Proposition \ref{prop:012:concat}\ref{item:012concat:2+1}, which describes construction of a type-2 evolution by concatenating a type-2 and a type-1, 
 so we can conclude that $\big(m_{(1)}^{D_1},\beta_{(1)}^{D_1}\big) = \big(m_i^{D_1},\beta_E^{D_1}\big)$. Thus, again, \eqref{eq:swapred:k+1} holds.
\end{proof}


\begin{proposition}
\label{prop:Dynkin:killed}
 For any initial distribution $\mu$ on $\TInt_{k+1}$, it is possible to define a pair of coupled processes such that:
\begin{itemize}
  \item $\big(\cT^y_{k+1},\,y\ge0\big)$ is a resampling $(k\!+\!1)$-tree evolution with $\cT^0_{k+1}\sim\mu$;
  \item $\big(\cT^y_{k},\,y\ge0\big)$ is a killed $k$-tree evolution;
  \item $\cT^y_k = \pi_k\big(\cT^y_{k+1}\big)$ for all $y\in [0,D^{\le k}\wedge D_\infty)$, where $D^{\le k}$ is the degeneration time of $\big(\cT^y_{k},\,y\ge0\big)$, while $D_\infty$ is the accumulation point of the degeneration times of $\big(\cT^y_{k+1},\,y\ge0\big)$; and
  \item $\big(\cT^y_{k},\,y\ge0\big)$ is strongly Markovian in the filtration generated by both processes.
\end{itemize}
\end{proposition}

\begin{proof}
 We begin with the coupled killed evolutions of the previous two lemmas, defined on a probability space $(\Omega_{(1)},\cF_{(1)},\bP_{(1)})$. We will extend this construction recursively, one degeneration at a time, to obtain a resampling $(k\!+\!1)$-tree evolution with the claimed properties.
 
 Suppose that for some $n\ge 1$ we have defined $\big(\cT^y_{k+1},y\in [0, D_n)\big)$ on some extension $(\Omega_{(n)},\cF_{(n)},\bP_{(n)})$ of $(\Omega_{(1)},\cF_{(1)},\bP_{(1)})$ so that this is distributed as a resampling $(k\!+\!1)$-tree evolution stopped at its $n^{\text{th}}$ degeneration time. Suppose also that $\cT^y_k := \pi_k\big(\cT^y_{k+1}\big)$ for $y\in [0, D_n\wedge D^{\le k})$ and that, on the event $A_n := \{D_n < D^{\le k}\}$ in which $k\!+\!1$ is the label dropped in each of the first $n$ degenerations, we get 
$\varrho\big(\cT^{D_n-}_{k+1}\big) = \cT^{D_n}_k$ a.s., as in \eqref{eq:swapred:k+1}.
 
 We further extend our probability space to $(\Omega_{(n+1)},\cF_{(n+1)},\bP_{(n+1)})$ to include additional random objects with the following conditional distributions given the sub-$\sigma$-algebra of $\cF_{(n+1)}$ that corresponds to $\cF_{(n)}$.
 \begin{itemize}
  \item On $A_n^c$, we require a process $\big(\cT^{(n),y}_{k+1},\,y\ge0\big)$ conditionally distributed as a killed $(k\!+\!1)$-tree evolution with initial law $\Lambda_{k+1,[k]}\big(\varrho\big(\cT^{D_n-}_{k+1}\big),\cdot\,\big)$. On $A_n$ we instead define this to be the constant process at $0$.
  \item On $A_n$, we require a random block $L_{n}$ conditionally distributed as a size-biased pick from $\block\big(\cT_k^{D_n}\big)$. On $A_n^c$ we set $L_n = 0$.
  \item On $A_{n,1} := \{L_n\in [k]\} \subseteq A_n$, we require a process $\big(U_{(n+1)}^y,y\ge 0\big)$ that is conditionally distributed 
as a pseudo-stationary type-2 evolution, as in Proposition \ref{prop:012:pseudo}, conditioned to have total mass process $\big\|U_{(n+1)}^y\| = m_{L_n}^{D_n + y}$, $y\ge0$. On $A_{n,1}^c$ we define this to be the constant process at $0$.
  \item On $A_{n,2} := A_n\setminus A_{n,1}$, with $L_n = (E,a,b)\in \block\big(\cT^{D_n}_k\big)\setminus [k]$, we require a pair of processes with conditional law 
   \begin{equation}
   \label{eq:DK:split_v2}
    \big( \Gamma_{(n+1)}^y, \big(m_{(n+1)}^y,\beta_{(n+1)}^y\big), y\!\in\! [0,\wt D)\big) \sim
    	\kappa\left( \left( \left( \Gamma_E^{D_n + y},y\!\ge\!0\right)\!,(a,b)\right)\!,\cdot\,\right)\!,
   \end{equation}
   where $\big(\Gamma_E^y,y\ge0\big)$ denotes the evolution on the type-0/1/2 compound in $\big(\cT^y_k\big)$ containing edge $E$, and $\kappa$ denotes the kernel described in Corollary \ref{cor:decomp_ker}, allowing us to decompose this type-0/1/2 evolution into a type-0/1/2 evolution $\big(\Gamma_{(n+1)}^y\big)$ concatenated with a type-1 evolution $\big(m_{(n+1)}^y,\beta_{(n+1)}^y\big)$ up until the degeneration time $\wt D$ of $\big(\Gamma_{(n+1)}^y\big)$. On $A_{n,2}^c$ we define these to be constant $0$ processes.
 \end{itemize}
 On each of these events, we define $\Delta_{n+1}$ to be a different degeneration time. On $A_n^c$, it equals the degeneration time of $\big(\cT^{(n),y}_{k+1}\big)$; on $A_{n,1}$, it is the minimum of $D^{\le k}-D_n$ and the degeneration time of $\big(U_{(n+1)}^y\big)$; and on $A_{n,2}$, it is the minimum of $D^{\le k}-D_n$, the degeneration time $\wt D$ of $\big( \Gamma_{(n+1)}^y\big)$, and that of $\big(m_{(n+1)}^y,\beta_{(n+1)}^y\big)$.
 
 We define $D_{n+1} := D_n + \Delta_{n+1}$; thus, on the event $A_n$ that no lower label has degenerated prior to $D_n$, we get $D_{n+1}\le D^{\le k}$.
 
 On $A_n^c$ we define $\cT^y_{k+1} := \cT^{(n),y-D_n}_{k+1}$ for $y\in [D_n,D_{n+1})$. On $A_{n,1}$, we define $\big(\cT^y_{k+1},\, y\in [D_n,D_{n+1})\big)$ as in \eqref{eq:Dynkin:insert_ext}, in Case 1 in the proof of Lemma \ref{lem:Dynkin:killed:1}, inserting the type-2 evolution $U_{(n+1)}^y$ in place of the individual evolving mass on leaf $L_n$. On $A_{n,2}$, we define this process as in \eqref{eq:Dynkin:insert_int}, in Case 2 in that proof.
 
 By construction and by definition of the resampling kernel in \eqref{eq:resampling_kernel:def}, $\cT^{D_n}_{k+1}$ has conditional law $\Lambda_{k+1,[k]}\big(\varrho\big(\cT^{D_n-}_{k+1}\big),\cdot\,\big)$ given $\big(\cT^{y}_{k+1},\,y\in [0,D_n)\big)$. Hence, at this degeneration time, $\big(\cT^{y}_{k+1}\big)$ behaves like a resampling $(k\!+\!1)$-tree evolution. Moreover, by the argument in the proof of Lemma \ref{lem:Dynkin:killed:1}, $\big(\cT^{D_n+y}_{k+1},y\in [0,\Delta_{n+1})\big)$ is distributed as a killed $(k\!+\!1)$-tree evolution. Putting these pieces together, $\big(\cT^y_{k+1},y\in [0, D_{n+1})\big)$ is a resampling $(k\!+\!1)$-tree evolution stopped at its $(n+1)^{\text{st}}$ degeneration time, and $\cT^y_k = \pi_k\big(\cT^y_{k+1}\big)$ for $y\in [0,D_{n+1}\wedge D^{\le k})$. By the same arguments as in the proof of Lemma \ref{lem:Dynkin:killed:2}, equation \eqref{eq:swapred:k+1} holds at time $D_{n+1}$ on the event $A_{n+1}$ that label $k\!+\!1$ resamples an $(n\!+\!1)^{\text{st}}$ time before the first time that a lower label would resample. 
 
 By the Ionescu Tulcea theorem \cite[Theorem 6.17]{Kallenberg}, there is a probability space $(\Omega_\infty,\cF_{\infty},\bP_{\infty})$ on which we can define a resampling $(k\!+\!1)$-tree evolution $\big(\cT^y_{k+1},\,y\ge0\big)$, with $\pi_k\big(\cT^y_{k+1}\big) = \cT^y_k$ for $y\in [0, D_\infty)$ and $\cT^0_{k+1}\sim \mu$.
 As in the proof of Lemma \ref{lem:Dynkin:killed:1}, the claim that $\big(\cT^y_{k},y\in [0, D^{\le k})\big)$ is strongly Markovian in the filtration generated by itself and $\big(\cT^y_{k+1},y\ge 0\big)$ follows from the assertions concerning filtrations at the ends of Propositions \ref{prop:012:mass} and \ref{prop:012:concat}.
\end{proof}
 

\begin{proof}[Proof of Proposition \ref{prop:consistency_0}] 
 As in the statement of the proposition, let $(\cT^y_{k+1},y\ge0)$ denote a resampling $(k\!+\!1)$-tree evolution with $\cT^0_{k+1} \sim\Lambda_{k+1,[k]}(T,\cdot\,)$, and let $\cT^y_k:=\pi_k\big(\cT^y_{k+1}\big)$, $y\ge0$. 
 Let $\big(D^{\le k}_n,\,n\ge1\big)$ denote the sub-sequence of degeneration times for $\big(\cT^y_{k+1}\big)$ at which a label in $[k]$ drops and resamples. 
 By Proposition \ref{prop:Dynkin:killed}, $\cT^{y}_k$ evolves as a resampling $k$-tree evolution for $y\in \big[0,D_\infty\wedge D^{\le k}_1\big)$. By Lemma \ref{lem:projection_at_degen}, given $\big\{D^{\le k}_1 < D_\infty\big\}$, 
 $\cT^{D^{\le k}_1-}_k$ is degenerate with $J\big(\cT^{D^{\le k}_1-}_k\big) = J\big(\cT^{D^{\le k}_1-}_{k+1}\big) =: J_1$.
 
 For the purpose of the following, let $R$ and $\Pi_k$ denote the trivial stochastic kernels associated with $\varrho$ and $\pi_k$, i.e.\ $R(T,\cdot\,) = \delta_{\varrho(T)}(\,\cdot\,)$ and $\Pi_k(T,\cdot\,) = \delta_{\pi_k(T)}(\,\cdot\,)$. Recall from Definition \ref{def:mark} that $\phi_1$ and $\Phi_1$ are the map and associated stochastic kernel projecting a marked $k$-tree $(T,\ell)$ to an (unmarked) $k$-tree $T$. As in \eqref{eq:composition_of_kernels}, we follow the standard convention of reading compositions of stochastic kernels from left to right.
 
 Now recall Lemma \ref{lem:inter:preswap} and the definition of resampling: given the marked $k$-tree process $(\Ast\cT^y_{k},\,y\in[0,D^{\le k}_1))$ and the event $\big\{D^{\le k}_1 < D_\infty\big\}$, the projected tree after resampling, $\cT^{D^{\le k}_1}_k$, has conditional law
 \begin{equation*}
  \Ast\Lambda_k R \Lambda_{J_1,[k+1]\setminus\{J_1\}} \Pi_k \Big( \Ast\cT^{D^{\le k}_1-}_{k},\cdot\, \Big)
  = \Lambda_{J_1,[k]\setminus\{J_1\}}\Big( \varrho\Big(\cT^{D^{\le k}_1-}_k\Big),\cdot\, \Big),
 \end{equation*}
 where the second expression follows from Lemma \ref{lem:Dynkin:preswap}. See the commutative diagram in Figure \ref{fig:preswap}. Thus, at this degeneration time, the projected tree behaves the same as in Definition \ref{def:resamp_1} of resampling $k$-tree evolutions. We conclude by induction and the strong Markov property of $\big(\cT^y_{k+1},\,y\ge0\big)$ applied at the degeneration times $D^{\le k}_n$, $n\ge1$. 
\end{proof}

\section{Accumulation of degeneration times as mass hits zero}\label{sec:non_acc}

In this section we finally prove Proposition \ref{prop:non_accumulation}, which allows us to complete the proofs of Theorem \ref{thm:total_mass} and Theorem \ref{thm:consistency}\ref{item:cnst:resamp}. 
Recall, Proposition \ref{prop:non_accumulation} states that $D_\infty := \sup_nD_n$ equals $\inf\{y\ge0\colon \|\cT^{y-}\|=0\}$ for a resampling $k$-tree evolution $(\cT^y,y\ge 0)$ with degeneration times $(D_n,n\ge 1)$. 

\begin{lemma}\label{lem:degen_diff}
 Fix $k\ge 3$ and $\epsilon>0$. Let $T\in\TInt_{k-1}$ with $\|T\|>\epsilon$ and let $(\cT^y,y\!\ge\!0)$ be a resampling $k$-tree evolution with $\cT^0\!\sim\!\Lambda_{k,[k-1]}(T,\cdot\,)$. Let $(D^*_n,n\!\ge\!1)$\linebreak denote the subsequence of degeneration times at which label $k$ is dropped and resamples. Assume that with probability one we get $D_\infty > D^*_2$. Then there is some $\delta = \delta(k,\epsilon)>0$ that does not depend on $T$ such that $\bP(D^*_2>\delta)>\delta$.
\end{lemma}

We prove this lemma in Appendix \ref{sec:non_acc_2}.

\begin{proof}[Proof of Proposition \ref{prop:non_accumulation} using Lemma \ref{lem:degen_diff}]
 By Proposition \ref{prop:total_mass_0}, the total mass $\|\cT^y\|$ of a resampling $k$-tree evolution evolves as a \besq$(-1)$ stopped at a random stopping time $D_\infty$. 
 Since a \besq$(-1)$ a.s.\ hits zero in finite time, if there were no infinite accumulation of degenerations prior to the total mass hitting 0, then there would be multiple degenerations simultaneously at that time. 
 But this is impossible, by the independence of the type-$i$ evolutions in the compounds of the $k$-tree in Definition \ref{def:killed_ktree} and the continuity of the distributions of their degeneration times, cf.\ Lemma \ref{lem:BESQ:length}. Hence, $D_\infty < \infty$ a.s. In fact, Lemma \ref{lem:BESQ:length} entails $\bE[D_\infty] < \infty$.%
 
 We will prove the proposition by showing that for every $\epsilon\in (0,\|\cT^0\|)$ we get $H_\epsilon := \inf\{y\ge0\colon 0<\|\cT^y\|\le\epsilon\} < D_\infty$ a.s., or equivalently, $H_\epsilon < \infty$, as the resampling $k$-tree evolution is defined to jump to 0 
at time $D_\infty$. This would imply that these times $H_\epsilon$ have a limit in $[0,D_\infty]$ at which time $\|\cT^y\|$ converges to zero, by the continuity of the total mass process noted in Proposition \ref{prop:total_mass_0}. Thus, $D_\infty$ would be sandwiched between the supremum of the sequence of $H_{1/n}$ stopping times and the time when total mass hits 0, so all three would have to be equal, thereby completing the proof.
 
 Fix $\epsilon > 0$. We will proceed by induction on the number $k$ of leaves in our $k$-tree evolutions, beginning with $k=2$. Consider a (self-similar) resampling 2-tree evolution starting from any unit-mass 2-tree. After its first degeneration, when label 2 resamples, it takes a new state according to a pseudo-stationary distribution, as described in Proposition \ref{prop:pseudo:resamp}. There is some $\delta>0$ such that a pseudo-stationary (self-similar) 2-tree evolution with initial mass $\epsilon$ will \emph{not} degenerate prior to time $\delta$ with probability at least $\delta$. By the self-similarity noted in Theorem \ref{thm:Markov}, the same holds for any larger initial mass with the same $\delta$. This proves that
 \begin{equation}\label{eq:degen_diff}
  \bP\big(D^*_{n+2}-D^*_{n} > \delta\ \big|\ H_{\epsilon}>D^*_{n}\big) > \delta \quad \text{for all }n\ge 1,
 \end{equation}
 where $(D^*_n,\,n\ge1)$ is the sequence of times at which the highest label, in this case label 2, is dropped in degeneration. In this base case, this is a meaningless distinction, as label 2 is always dropped in degeneration, and the ``$D_{n+2}$'' in the formula could be replaced by ``$D_{n+1}$,'' but we are preparing to appeal to this display again in the inductive step. By linearity of expectation, \eqref{eq:degen_diff} implies
 \begin{equation}
 \label{eq:degen:BCL}
  \infty > \bE[D_\infty] > \sum_{n \ge 1}\delta^2 \bP\{H_{\epsilon} > D^*_{2n}\}.
 \end{equation}
 It follows by the Borel--Cantelli Lemma that $H_\epsilon < D_\infty$ a.s., as desired.
 
 Now, suppose for induction that the proposition holds for $k$-tree evolutions and consider a resampling $(k\!+\!1)$-tree evolution $(\cT^y,y\ge0)$. By Proposition \ref{prop:consistency_0}, $(\pi_k(\cT^y),y\ge0)$ is a resampling $k$-tree evolution up to the accumulation time $D_\infty$ of degenerations of the $(k\!+\!1)$-tree evolution. The degeneration times of $(\pi_k(\cT^y))$ are the times at which a label less than or equal to $k$ is dropped and resamples in $(\cT^y)$. By the inductive hypothesis, these degeneration times do not have an accumulation point prior to the extinction time of the \besq$(-1)$ total mass. Thus, $D_\infty$ must equal the accumulation point of degeneration times $(D^*_n,\,n\ge1)$ at which label $k\!+\!1$ resamples. After $D^*_1$, label $k\!+\!1$ resamples, so that by the strong Markov property, $(\cT^{D^*_1+y},\,y\ge0)$ satisfies the hypotheses of Lemma \ref{lem:degen_diff}. That lemma now implies \eqref{eq:degen_diff}, which in turn implies \eqref{eq:degen:BCL}, which again proves the proposition in this case, by the Borel--Cantelli Lemma. By induction, this completes the proof.
\end{proof}

We can now complete the proofs of the two other theorems mentioned at the beginning of this section.

\begin{proof}[Proof of Theorem \ref{thm:total_mass}]
 By Proposition \ref{prop:total_mass_0}, the total mass process of a resampling $k$-tree evolution is $\besq(-1)$ up until time $D_\infty$. By Proposition \ref{prop:non_accumulation}, at this time the total mass approaches 0 continuously, as desired.
\end{proof}

\begin{proof}[Proof of Theorem \ref{thm:consistency}$(\mathrm{ii})$]
 By Proposition \ref{prop:consistency_0} and induction on the difference $k-j$, the desired projective consistency between resampling $k$- and $j$-tree evolutions, for $k>j\ge2$, holds up until the accumulation time $D_\infty$ of degenerations in the $k$-tree evolution. By Theorem \ref{thm:total_mass}, $D_\infty$ is the time at which both processes are continuously absorbed at total mass 0.
\end{proof}

\section[Proofs of remaining consistency results, including Theorem 1.5]{Proofs of remaining consistency results, including Theorem \ref{thm:intro:k_tree}}\label{sec:const:other}

\begin{proof}[Proof of Theorem \ref{thm:consistency}$(\mathrm{i})$]
 It suffices to prove projective consistency between non-resampling $(k\!+\!1)$- and $k$-tree evolutions, as we can extend to general projections from $k$- to $j$-tree evolutions for $k>j\ge2$ by induction on the difference $k-j$.
 
 By Lemma \ref{lem:Dynkin:killed:1}, we can define a coupled pair $\big(\cT^y_{k+1},\,y\ge0\big)$ and $\big(\cT^y_{k},\,y\ge0\big)$ of a killed $(k\!+\!1)$- and $k$-tree evolution, respectively, that satisfy $\cT^y_{k} = \pi_k\big(\cT^y_{k+1}\big)$ for $y\in [0,D_1)$, where $D_1$ is the degeneration time of $\big(\cT^y_{k+1}\big)$. In fact, because these are killed evolutions, and therefore do not invoke the swap-and-reduce map, this lemma is not sensitive to the particular choice of labels; we can do the same for any two finite label sets $A$ and $B$ of cardinality at least 2 that differ by the addition of one label.
 
 As in the proof of Proposition \ref{prop:Dynkin:killed}, we extend our construction inductively past degeneration times. If label $k\!+\!1$ degenerates in $\big(\cT^y_{k+1}\big)$ at time $D_n$, then by Lemma \ref{lem:Dynkin:killed:2}, $\cT^{D_n}_{k+1} := \varrho\big(\cT^{D_n-}_{k+1}\big) = \cT^{D_n}_{k}$. Thereafter, the two processes can be defined to be equal. Otherwise, if a lower label degenerates in $\big(\cT^y_{k+1}\big)$ at time $D_n$, then by Lemma \ref{lem:projection_at_degen}, the same label degenerates in $\big(\cT^y_{k}\big)$ at that time, and 
 $$\pi_k\big(\cT^{D_n}_{k+1}\big) = \pi_k\big(\varrho\big(\cT^{D_n-}_{k+1}\big)\big) = \varrho\big(\cT^{D_n-}_{k}\big) = \cT^{D_n}_k.$$
 Then, by Lemma \ref{lem:Dynkin:killed:1}, we can extend our coupled construction until the next degeneration time.
 
 The end result of this construction is a non-resampling $(k\!+\!1)$-tree evolution coupled to a non-resampling $k$-tree evolution so that the latter equals the $\pi_k$-projection of the former at all times, as desired.
\end{proof}

\begin{proof}[Proof of Theorem \ref{thm:consistency}$(\mathrm{iii})$]
 Fix $1\le j<k$. Suppose $\fT_k := (\cT_k^y,y\ge0)$ is a resampling $k$-tree evolution with initial distribution as in \eqref{eq:cnst:init}, so that $\fT_j = (\cT_j^y,y\ge0) := (\pi_j(\cT_k^y),y\ge0)$ is a resampling $j$-tree evolution. Then because these evolutions have the same total mass process, they require the same time change for de-Poissonization: $(\rho_{\fT_k}(u),u\!\ge\!0) = (\rho_{\fT_j}(u),u\!\ge\!0)$. Thus, the associated unit-mass processes of Definition \ref{def:dePoi} are also projectively consistent. The same argument holds in the non-resampling case.
\end{proof}

\begin{proof}[Proof of Proposition \ref{prop:resamp_to_non}]
%
 Suppose $\big(\cT_{k,+}^y,y\ge0\big)$ is a resampling $k$-tree evolution. Let $(D_n,n\ge1)$ denote its sequence of degeneration times, and set $D_0 := 0$. 
 We see, for example from Theorem \ref{thm:consistency}$(\mathrm{ii})$, that in between degenerations of $\cT_{k,+}$, any projection of this process with permuted labels will evolve as a killed $k$-tree evolution. Thus, to prove that we can get a non-resampling evolution from a projection and permutation, as described in the proposition, it suffices to describe a suitable sequence of projections and permutations that change at the degenerations of $\cT_{k,+}$.
 
 Recall that we consider each edge in a tree shape to be labeled by the set of all labels of leaves in the subtree above that edge.  Recall our terminology around tree shapes and the definition of $\varrho$ in Sections \ref{sec:killed_def} and \ref{sec:non_resamp_def}, respectively: when a label $i_n := I\big(\cT_{k,+}^{D_n-}\big)$ causes degeneration, it swaps places with label $j_n := J\big(\cT_{k,+}^{D_n-}\big) = \max\{i_n,a_n,b_n\}$, where $a_n$ and $b_n$ are respectively the least labels on the sibling and uncle of leaf edge $\{i_n\}$ 
in the tree shape of $\cT_{k,+}^{D_n-}$, with the convention that $b_n = 0$ in the special case that the parent of $\{i_n\}$ is the root edge, $\parent{\{i_n\}} = [k]$. In the resampling evolution, label $j_n$ is resampled.
 
 We extend this notation. Let $E_{n}^{(a)}$ and $E_n^{(b)}$ denote the sets of labels on the sibling and uncle of edge $\{i_n\}$, with the convention $E_n^{(b)} = \emptyset$ when $\parent{\{i_n\}} = [k]$. Then $a_n = \min\!\big(E_n^{(a)}\big)$ and $b_n = \min\!\big(E_n^{(b)}\cup\{0\}\big)$. Let $\tau_n$ denote the transposition permutation that swaps $i_n$ with $j_n$.
 
 Set $A_0 := B_0 := [k]$ and let $\sigma_0$ denote the identity map on $[k]$. Now suppose for a recursive construction that we have defined $(A_{n-1},B_{n-1},\sigma_{n-1})$. For $1\le j < n$, $D_{j-1} \le y < D_j$, let $\cT^y_{k,-} := \sigma_j\circ\pi_{A_j}\big(\cT^y_{k,+}\big)$. 
 We consider six cases.
 \smallskip 
 
 \begin{enumerate}[label = Case \arabic*:, ref = \arabic*, topsep=5pt, itemsep=5pt, itemindent=1.82cm, leftmargin=0pt]
  \item  $i_n\notin A_{n-1}$ and $j_n\notin A_{n-1}$. In this case, the degeneration, swap-and-reduce map, and resampling in $\cT_k^{D_n}$ are invisible under $\sigma_{n-1}\circ \pi_{A_{n-1}}$, since the projection erases both labels involved. We set $(A_n,B_n,\sigma_n) := (A_{n-1},B_{n-1},\sigma_{n-1})$.
  
  \item $i_n\notin A_{n-1}$ and $j_n\in A_{n-1}$. In this case, the label $i_n$ that has caused degeneration is invisible under $\pi_{A_{n-1}}$, so there is no degeneration in the projected process, but $i_n$ displaces a label that \emph{is} visible. To maintain continuity in the projected process at this time, $i_n$ takes the place of $j_n$ in such a way that $\sigma_n(i_n) = \sigma_{n-1}(j_n)$. In particular, $A_n:= (A_{n-1}\setminus\{j_n\})\cup \{i_n\}$, $B_n := B_{n-1}$, and $\sigma_n := \sigma_{n-1}\circ\tau_n|_{A_n}$.\medskip
  
  In each of the remaining cases, $i_n\in A_{n-1}$. Let $\wi_n := \sigma_{n-1}(i_n)$.
  
  \item \label{case:r2n:degen} $i_n\in A_{n-1}$ and both $E_n^{(a)}$ and $E_n^{(b)}$ intersect $A_{n-1}$ non-trivially. In this case, the degeneration caused by $i_n$ in $\cT_{k,+}$ corresponds to a degeneration caused by $\wi_n$ in $\cT_{k,-}$.
  
  Let $\tilde a_n := \min\!\big(\sigma_{n-1}\big(E_n^{(a)}\cap A_{n-1}\big)\!\big)$ and $\tilde b_n := \min\!\big(\sigma_{n-1}\big(E_n^{(b)}\cap A_{n-1}\big)\!\big)$. Let $\wj_n := \max\{\wi_n,\tilde a_n,\tilde b_n\}$ and let $\tilde\tau_n$ denote the transposition permutation that swaps $\wi_n$ with $\wj_n$. If $j_n\in A_{n-1}$, then we set $A_n := A_{n-1}\setminus\{j_n\}$; otherwise, we set $A_n := A_{n-1}\setminus\{i_n\}$. In either case, we define $B_n := B_{n-1}\setminus\{\wj_n\}$ and $\sigma_n := \tilde\tau_n\circ\sigma_{n-1}\circ\tau_n|_{A_n}$.
 
  \item \label{case:r2n:degen2} $i_n\in A_{n-1}$ and $E_n^{(a)}$ intersects $A_{n-1}$ non-trivially, while $\parent{\{i_n\}} = [k]$ (recall that when $\parent{\{i_n\}} = [k]$, there is no uncle to edge $\{i_n\}$, so by convention, $E_n^{(b)} = \emptyset$ and $b_n=0$).
  
  As in Case \ref{case:r2n:degen}, the degeneration caused by $i_n$ in $\cT_{k,+}$ corresponds to a degeneration caused by $\wi_n$ in $\cT_{k,-}$, but in this case, both degenerations occur at the root of the respective trees. We define $(A_n,B_n,\sigma_n)$ as in the previous case, but with $\tilde b_n := 0$.
  
  \item \label{case:r2n:shrink} $i_n\in A_{n-1}$ and $E_n^{(a)}$ is disjoint from $A_{n-1}$. Then leaf block $i_n$ and the subtree that contains label set $E_n^{(a)}$ in $\cT_{k,+}^y$ project down to a single leaf block, $\sigma_{n-1}(i_n)$, in $\cT_{k,-}^y$ as $y$ approaches $D_n$. By leaving open the possibility that $E_n^{(b)}$ may be disjoint from $A_{n-1}$ as well, we include in this case the possibility that the subtree of $\cT_{k,+}^y$ with label set $E_n^{(b)}$ projects to this same leaf block as well. Regardless, this degeneration is ``invisible'' in $\cT_{k,-}$. In order to keep label $\wi_n$ in place in the projected process, if label $i_n$ resamples or swaps with a label in $E_n^{(b)}$, then we choose a label in $E_n^{(a)}$ to map to $\wi_n$ under $\sigma_n$.
  
  \begin{enumerate}[label = Case \theenumi.\arabic*: , ref=\theenumi.\arabic*, topsep=0pt, itemsep=0pt, itemindent=2cm, leftmargin=0pt]
   \item $j_n = a_n$. Then we define $(A_n,B_n,\sigma_n) := (A_{n-1},B_{n-1},\sigma_{n-1})$.\label{case:r2n:shrink:null}
   
   \item $j_n = i_n$ or $j_n = b_n$. Then let $\hat\tau_n$ denote the transposition that swaps $i_n$ with $a_n$. If $j_n\in A_{n-1}$, as is always the case when $j_n = i_n$, then we set $A_n := (A_{n-1}\setminus\{j_n\})\cup\{a_n\}$. Otherwise, if $j_n\notin A_{n-1}$, then we set $A_n := (A_{n-1}\setminus\{i_n\})\cup\{a_n\}$. In either case, we define $B_n := B_{n-1}$ and $\sigma_n := \sigma_{n-1}\circ\hat\tau_n\circ\tau_n|_{A_n}$.\label{case:r2n:shrink:swap}
  \end{enumerate}
  
  \item $i_n\in A_{n-1}$ while $E_n^{(a)}$ intersects $A_{n-1}$ non-trivially, $\parent{\{i_n\}} \neq [k]$, and $E_n^{(b)}$ is disjoint from $A_{n-1}$. This degeneration time in $\cT_{k,+}^y$ corresponds to a time at which the leaf block $\wi_n$ in $\cT_{k,-}^y$ has mass approaching zero (more precisely, it is a.s.\ an accumulation point of prior times at which this mass equals zero) while the interval partition on its parent edge has a leftmost block. The subtree of $\cT^{D_n-}_{k,+}$ that contains the leaf labels $E_n^{(b)}$ maps to a single \emph{internal} block, the aforementioned leftmost block, in $\cT_{k,-}^{D_n-}$. Therefore, we define $\sigma_n$ in such a way that some label in $E_n^{(b)}$ gets mapped to $\wi_n$, so that the label $\wi_n$ ``moves into'' the leftmost block in the projected process, as in a type-1 or type-2 evolution; see Proposition \ref{prop:012:pred}. In fact, we can accomplish this with the same definitions of $(A_n,B_n,\sigma_n)$ as in Cases \ref{case:r2n:shrink:null} and \ref{case:r2n:shrink:swap}, but with roles of $a_n$ and $b_n$ reversed.
%
%
 \end{enumerate}
 
 It follows from the consistency result of Theorem \ref{thm:consistency}\ref{item:cnst:nonresamp} that for each $n$, the projected process evolves as a stopped non-resampling $k$-tree evolution (or $B_n$-tree evolution) during the interval $[D_n,D_{n+1})$. By our construction, we have $\cT_{k,-}^{D_n} = \varrho\big(\cT_{k,-}^{D_n-}\big)$ in Case \ref{case:r2n:degen}, as in Definition \ref{def:nonresamp_1} of non-resampling evolutions. In the other cases, it follows from the arguments in the proof of Lemma \ref{lem:Dynkin:killed:2} that each type-0/1/2 compound in $\cT_{k,-}^{D_n}$ attains the value required by the type-0/1/2 evolution in that compound, given its left limit in $\cT_{k,-}^{D_n-}$. Thus, $\big(\cT_{k,-}^y,y\ge0\big)$ is a non-resampling $k$-tree evolution. 
\end{proof}

Finally, we have the ingredients needed to prove one of the main theorems of this memoir.

\begin{proof}[Proof of Theorem \ref{thm:intro:k_tree}]
 Resampling unit-mass $k$-tree evolutions satisfy all of the properties claimed in Theorem \ref{thm:intro:k_tree}.
 
 (i) These are $k$-tree-valued processes, by definition.
 
 (ii) They possess the required projective consistency, by Theorem \ref{thm:consistency}\ref{item:cnstfam:dePoi}.
 
 (iii) By Theorem \ref{thm:dePoi}, they are stationary with the laws of Brownian reduced $k$-trees, $k\ge 2$.
 
 (iv) By Corollary \ref{cor:WF}, they project to mixed-parameter Wright--Fisher diffusions, as claimed.
\end{proof}

\section{Consistent partially resampling $k$-tree evolutions}
\label{sec:partial_resamp}

Theorem \ref{thm:intro:k_tree}, which has now been proved, allows the construction of a rooted, weighted $\BR$-tree-valued process, as described above that theorem statement. In order to prove the main properties of that process, as listed in Theorem \ref{thm:intro:AD}, we require one more family of $k$-tree-valued processes. Recall that, for a degenerate tree $T\in\tdTInt_k \setminus \TInt_k$, $J(T)$ denotes the label that would be dropped by the swap-and-reduce map, $\varrho(T)\in \TInt_{[k]\setminus\{J(T)\}}$. Here, we define $m$-tree evolutions in which low labels (up to some threshold $k\le m$) do not resample when dropped in degeneration, but high labels do.

\begin{definition}[Partially resampling $m$-tree evolution]\label{def:partial_resamp}
 Fix $m\ge k \ge 1$ and some $\cT^0_{(1)} =T\in \TInt_{m}$. Let $A_1 := [m]$. Inductively for $n\ge1$, let $(\cT^y_{(n)},y\in [0,\Delta_n))$ denote a killed $A_n$-tree evolution from initial state $\cT_{(n)}^0$, run until its degeneration time $\Delta_n$, conditionally independent of $(\cT_{(j)},j< n)$ given its initial state. Let $J_n := J\big(\cT^{\Delta_n-}_{(n)}\big)$. We define $\cT_{(n+1)}^0$ to equal $\varrho\big(\cT^{\Delta_n-}_{(n)}\big)$ on the event $\{J_n\in [k]\}$, and on the event $\{J_n > k\}$ to have conditional distribution $\Lambda_{J_n,A_n\setminus\{J_n\}}\big(\varrho(\cT_{(n)}^{\Delta_n-}),\cdot\,\big)$ given $(\cT_{(j)},j\le n)$. In other words, labels in $[k]$ do not resample, but higher labels do.
 
 We set $D_0=0$ and define $D_n = \sum_{j=1}^n\Delta_j$, $n\ge 1$. For $y\in [D_{n-1},D_{n})$ we define $\cT^y = \cT_{(n)}^{y-D_{n-1}}$. For $y\ge D_\infty:=\sup_{n\ge 0}D_n$ we set $\cT^y= 0 \in \TInt_{\emptyset}$. Then $(\cT^y,\,y\ge 0)$ is an \emph{$(m+,k-)$-partially resampling (self-similar) $m$-tree evolution with initial state $T$}.
\end{definition}


\begin{proposition}
\label{prop:partial_resamp:stat}
 Fix $k\ge 1$, and let $\big(\cT_{m+,k-}^0,\, m\ge k\big)$ be a consistent family of Brownian reduced $m$-trees. Then, starting from this initial family, there are consistent $(m+,k-)$-partially resampling $m$-tree evolutions $(\mathcal{T}_{m+,k-}^y,\,y\ge 0)$, $m\ge k$. 
 
 Fix $y\ge 0$ and denote by $A^y$ the label set of $\mathcal{T}_{k+,k-}^y$. Conditionally given $A^y$, the tree $\mathcal{T}_{m+,k-}^y$ is a scaled Brownian reduced $(m-k+\#A^y)$-tree with label set $([m]\setminus[k])\cup A^y$. 
\end{proposition}
\begin{proof} Let $R_0=0$ and $A_0=[k]$. Suppose that, for some $j\ge 0$, we have constructed the consistent family up to a resampling time $R_j$, with remaining label set 
  $A_j\subseteq[k]$ of size $k-j$. 
 By Corollary \ref{cor:consistent_fam}(ii) and straightforward relabeling, there exists a consistent family of resampling $(m-k+\#A_j)$-tree evolutions starting from $\mathcal{T}_{m+,k-}^{R_j}$, $m\ge k$, which we consider up to but excluding the first time $\Delta_j$ that a label in $A_j$ first resamples. We define $\mathcal{T}_{m+,k-}^{R_j+z}$ to be this process for $0<z<\Delta_j$, and we define $R_{j+1} := R_j+\Delta_j$.
 
 In the case $j=k-1$, we recall that label 1 is never dropped (see e.g.\ Proposition \ref{prop:pseudo:degen}). Therefore, we have $A_j=\{1\}$ and $R_{j+1}=\infty$, and the construction of $((\mathcal{T}_{m+,k-}^y,\,y\ge 0),\,m\ge k)$ is complete.
 
 If $j<k-1$, then $R_{j+1}<\infty$. In this case, we define $A_{j+1} := A_j\setminus \{J_j\}$, where $J_j := J\big(\mathcal{T}_{k+,k-}^{R_{j+1}-}\big)$. By Lemma \ref{lem:projection_at_degen}, the family $\mathcal{T}_{m+,k-}^{R_{j+1}} := \varrho\big(\mathcal{T}_{m+,k-}^{R_{j+1}-}\big)$, $m\ge k$, is again projectively consistent; and by Proposition \ref{prop:pseudo:degen}, these are Brownian reduced $(m-k+\#A_{j+1})$-trees labeled by $([m]\setminus[k])\cup A_{j+1}$.  Hence, the induction proceeds and completes the construction.
  
Now let $A^y=A_{j-1}$ for $R_{j-1}\le y<R_j$, $j\in[k]$. The claimed time-$y$ marginal distributions follow from the aforementioned distributions at the stopping times $R_0,\ldots,R_{k-1}$, Proposition \ref{prop:pseudo:pre_D}, and straightforward relabeling.  
\end{proof}

\begin{proposition}\label{prop:partialfromresampling}
 Let $(\cT^y_m,y\ge 0)$, $m\ge 1$, be a consistent system of resampling $m$-tree evolutions starting from Brownian reduced $m$-trees. For each pair $m\ge k\ge 1$, there exists a process $((A_y^{m,k},B_y^k,\sigma_y^{m,k}),y\ge 0)$ that is constant between degeneration times of $(\cT^y_m,\,y\ge 0)$, such that $\sigma_y^{m,k}$ is a bijection between $A_y^{m,k}\subset[m]$ and $B_y^{k}\cup([m]\setminus[k])$ with $B_y^k\subseteq[k]$, and such that $\cT_{m+,k-}^y:=\sigma_y^{m,k}\circ\pi_{A_y^{m,k}}(\cT^y_m)$, $y\ge 0$, is an $(m+,k-)$-partially resampling $m$-tree evolution. Furthermore, for any $k\ge 1$, these processes can be chosen to be projectively consistent in $m$, $m\ge k$.
\end{proposition}
\begin{proof}
 The proof of Proposition \ref{prop:resamp_to_non} can be adapted, as follows. The construction of the bijections is the same, except that in
  Cases \ref{case:r2n:degen} and \ref{case:r2n:degen2}, we reduce the sizes of $A_n$ and $B_n$ only when $\wj_n\le k$ and set $A_n:=A_{n-1}$ and $B_n:=B_{n-1}$ otherwise, extending
  $\sigma_n:=\tilde{\tau}_n\circ\sigma_{n-1}\circ\tau_n$ to hold on the larger $A_n=A_{n-1}$. This achieves that the resampling of the higher label
  for the partially resampling evolution follows the corresponding resampling in the fully resampling evolution.
\end{proof}

\chapter{The Aldous diffusion as a projective limit of $k$-tree evolutions}
\label{ch:properties}

In this chapter,  
we construct a continuum-tree-valued Markov process and identify it as the process conjectured by Aldous. 
Specifically, we consider the consistent system of stationary unit-mass $k$-tree evolutions starting from Brownian reduced $k$-trees, $k\ge 1$, of Theorem \ref{thm:intro:k_tree} as established in Chapter \ref{ch:consistency}. We reverse the construction of Brownian reduced $k$-trees from a Brownian CRT and study the map $S$ 
that associates with suitable consistent families of $k$-trees an associated (${\rm GHP}$-isometry class of a) rooted, weighted $\mathbb{R}$-tree, and we prove 
\begin{enumerate}\item[1.] The map $S$ projects the consistent family of stationary unit-mass $k$-tree evolutions, $k\ge 1$, to a stationary continuum-tree-valued Markov process, 
  which possesses a continuous modification. The stationary distribution is the distribution of the Brownian CRT. Cf. Theorem \ref{thm:intro:AD}. 
\end{enumerate}
As explained in the introduction, this process solves a conjecture that David Aldous formulated in the late 1990s. Specifically, Aldous \cite{Aldous00} studied a Markov chain on unrooted binary trees with $n$ labelled leaves, where each transition consists of removing and reinserting a leaf uniformly at random. Our process relates to the rooted variant of this Markov chain, where the root is an additional degree-1 vertex. These Markov chains have uniform stationary distributions. Aldous \cite{AldousCRT1} showed that suitable 
representations of uniform binary $n$-tree shapes have as their $n\rightarrow\infty$ scaling limit the Brownian CRT. Aldous \cite{ADProb2, AldousDiffusionProblem} observed that there are induced three-mass (or 
$(2k-1)$-mass) Markov chains that record subtree sizes around one (or $k-1$) branch points, and that these Markov chains, scaled and suitably sped up to make $n^2$ steps per unit time, appear to converge to
Wright--Fisher-like diffusions. He conjectured that ``these diffusions are recording certain aspects of an underlying diffusion on continuum trees.'' We show the following, hence identifying the process in 1.\ as this conjectured process, which we call the Aldous diffusion.
\begin{enumerate}\item[2.] For each $k\ge 2$ and a sample of $k$ leaves of the initial tree of the (stationary) continuum-tree-valued process of Claim 1, consider the
  reduced subtree spanned by these $k$ leaves and the root, as time evolves. Remove from the continuum trees the $k-1$ branch points of this reduced subtree and record the $2k-1$ masses
  of the connected components. Then the process of Claim 1 induces an evolution of the $2k-1$ component masses, stopped when one mass vanishes. 
  This stopped process is a Wright--Fisher diffusion with parameter $\frac12$ for each of the $k-1$ components between two branch points (or a branch point and the root) and 
  parameter $-\frac12$ for each of the other $k$ components. 
\end{enumerate}  
In Chapter \ref{chap8}, we use this construction to study the resulting Aldous diffusion. 

The structure of this chapter is as follows. In Section \ref{sec:Rtrees}, we recall from the literature the Gromov--Hausdorff--Prokhorov space of weighted $\mathbb{R}$-trees and study $\mathbb{R}$-tree projections of $k$-trees and their $k\rightarrow\infty$ limits associated with consistent families of $k$-trees hence formalising the map $S$ of Claim 1.
In Section \ref{sec:defmarkov} we use $S$ to define the process that will be our Aldous diffusion, establish the Markov property and discuss why the strong Markov property fails, in general. In Section 
\ref{sec:dItoGHP} we derive general bounds in terms of the interval partition metric $d_\cI$ bounding the Gromov--Hausdorff and Gromov--Hausdorff--Prokhorov distances 
between the (weighted) $\mathbb{R}$-trees associated with $k$-trees. In Section \ref{sec:subtreedecomp} we enhance subtree decompositions of the Brownian CRT from \cite{CW,DuW,HPW,PW13}. In Section \ref{sec:continuity}, we establish a path-continuous modification of our continuum-tree-valued process. 
In Section \ref{sec:properties} we pull the threads together and hence establish Claims 1 and 2 and thereby identify our process as the Aldous diffusion.
In Section \ref{sec:gen} we revisit our arguments and establish general Markovianity and continuity theorems under assumptions that require the construction of a suitable consistent system of k-tree evolutions, and we give some context as to where 
                                                                          these general results may apply.

\section{Introduction to weighted $\mathbb{R}$-trees and $k$-tree projections}\label{sec:Rtrees}

The aim of this section is to give a formal definition of the map $S$ that associates with suitable consistent families of $k$-trees a weighted $\mathbb{R}$-tree. We formalize the discussion of $\BR$-trees from the introduction, following \cite{ADH13,Gromov1999,M09}. To this end, first recall the \emph{Hausdorff distance} $d_M^{\rm H}$ on the set $\mathcal{K}_M$ of compact subsets and the \emph{Prokhorov distance} $d_M^{\rm P}$ on the set $\mathcal{M}_M$ of 
finite Borel measures in a complete and separable metric space $(M,d_M)$. Specifically, for $x\in M$ and $C\in\mathcal{K}_M$, denote by $d_M(x,C)=\min\{d_M(x,y)\colon y\in C\}$ the closest distance from $x$ to $C$ and by $C^\varepsilon=\{x\in M\colon d_M(x,C)\le\varepsilon\}$ the $\varepsilon$-thickening of $C$. Then for all $C,C^\prime\in\mathcal{K}_M$ and
$\lambda,\lambda^\prime\in\mathcal{M}_M$, 
\begin{align*}
d_M^{\rm H}(C,C^\prime)
 &:=\inf\big\{\varepsilon\!>\!0\colon C\subseteq (C^\prime)^\varepsilon\mbox{ and }C^\prime\subseteq C^\varepsilon\big\}\\
d_M^{\rm P}(\lambda,\lambda^\prime)
 &:=\inf\big\{\varepsilon\!>\!0\colon \lambda(C)\le\lambda^\prime(C^\varepsilon)\!+\!\varepsilon\mbox{ and }\lambda^\prime(C)\le\lambda(C^\varepsilon)\!+\!\varepsilon\mbox{ for all }C\in\mathcal{K}_M\big\}.
\end{align*}

\begin{definition}\label{def:Rtree}
 An \emph{$\BR$-tree} (\emph{real tree}) is a complete, separable metric space $(T,d)$ with the property that: (i) for each $x,y\in T$, there is a unique non-self-intersecting path in $T$ from $x$ to $y$ $[\![x,y]\!]_{T}$, and (ii) each such path $[\![x,y]\!]_{T}$ is isometric to a real interval $[0,d(x,y)]$. We will only consider \em compact \em $\mathbb{R}$-trees.
 
 A \emph{rooted, weighted $\BR$-tree} is a quadruple $(T,d,\rho,\mu)$, where $(T,d)$ is an $\BR$-tree, $\rho\in T$ is a distinguished vertex called the \emph{root}, and $\mu$ is a finite measure on the $\sigma$-algebra of Borel sets of $(T,d)$. 

  The (rooted) \emph{Gromov--Hausdorff--Prokhorov distance} $d_{\rm GHP}(\mathrm{T},\mathrm{T}^\prime)$ between two rooted, weighted $\mathbb{R}$-trees 
  $\mathrm{T}=(T,d,\rho,\mu)$ and $\mathrm{T}^\prime=(T^\prime,d^\prime,\rho^\prime,\mu^\prime)$ is defined as 
  \[
  d_{\rm GHP}(\mathrm{T},\mathrm{T}^\prime) := \inf_{\phi,\phi^\prime}\max\Big\{d_M^{\rm H}(\phi(T),\phi^\prime(T^\prime)),\ d_M(\phi(\rho),\phi^\prime(\rho^\prime)),\ d_M^{\rm P}(\phi_*\mu,\phi^\prime_*\mu^\prime)\Big\},
  \]
  where the infimum is taken over all metric spaces $(M,d_M)$ and all injective isometries $\phi\colon(T,d)\rightarrow(M,d_M)$ and 
  $\phi^\prime\colon(T^\prime,d^\prime)\rightarrow(M,d_M)$. 
  We say that $\mathrm{T}$ and $\mathrm{T}^\prime$ are \em ${\rm GHP}$-isometric \em if there is a bijective isometry $\iota\colon(T,d)\rightarrow(T^\prime,d^\prime)$ such that $\iota(\rho)=\rho^\prime$ and $\iota_*\mu=\mu^\prime$.
  We denote by $\mathbb{T}^{\rm real}$ the set of ${\rm GHP}$-isometry classes of rooted, weighted, compact $\mathbb{R}$-trees.
\end{definition} 

\begin{proposition}[Theorem 2.5 of \cite{ADH13}] The distance $d_{\rm GHP}(\mathrm{T},\mathrm{T}^\prime)$ only depends on the ${\rm GHP}$-isometry classes of $\mathrm{T}$ and $\mathrm{T}^\prime$, and 
  induces a metric on $\mathbb{T}^{\rm real}$, also denoted by $d_{\rm GHP}$. Furthermore, $(\mathbb{T}^{\rm real},d_{\rm GHP})$ is separable and complete.
\end{proposition}

\begin{definition}
 A \emph{random rooted, weighted $\BR$-tree} is a $(\mathbb{T}^{\rm real},d_{\rm GHP})$-valued random variable.
\end{definition}

Now consider a consistent family $(R_k,\,k\ge 1)\in\prod_{k\ge 1}\widebar{\bT}_k^{\rm int}$, i.e.\ $k$-trees $R_k=\big(\ft_k,(x_j^{(k)},j\in[k]),(\beta_E^{(k)},E\in\ft_k)\big)\in\widebar{\bT}_k^{\rm int}$,
$k\ge 1$, such that $\pi_{-k}(R_k)=R_{k-1}$ for all $k\ge 2$, in the sense of Definition \ref{def:proj}. Since $\beta_E^{(k)}\in\cI$ for all $E\in\ft_k$ and $k\ge 2$, 
each edge has a diversity $\sD\big(\beta_E^{(k)}\big)$ that we can use as a branch length. Although we will eventually work with ${\rm GHP}$-isometry classes in $(\mathbb{T}^{\rm real},d_{\rm GHP})$, we first construct representations of $R_k$, $k\ge 2$, as 
actual rooted, weighted $\mathbb{R}$-trees. We will work in $[0,\infty)^{\ft_k}$, $k\ge 2$, equipped with the $\ell_1$-distance. We denote by $e_E$ the unit vector in direction 
$E\in\ft_k$ and, for $w\in [0,\infty)^{\ft_k}$ and $c>0$, we write $w+[0,c]e_E:=\{w+xe_E,x\in[0,c]\}$.

Recall \eqref{eq:metric2tree}, where we associated with a 2-tree $(a,b,\gamma)\in\cJ^\circ$ an interval $[0,\sD(\gamma)]$
equipped with a weight measure, $M_2(a,b,\gamma)\in\mathcal{M}^\circ$ that adds two atoms of masses $a$ and $b$ at the ``top'' at $\sD(\gamma)$ and an atom of mass 
${\rm Leb}(U)$ at distance $\sD_\gamma(U)$ from the top, for each $U\in\gamma$. Note that $\mathcal{M}^\circ$ can be seen as a set of weighted (one-branch) 
$\mathbb{R}$-trees rooted at 0. We can represent $R_2$ by $M_2\big(x_1^{(2)},x_2^{(2)},\beta_{\{1,2\}}^{(2)}\big)$. 

\begin{definition}\label{def:SkandSkcirc} Adapting \eqref{eq:metric2tree}, we associate with $\beta\in\cI$ the weighted interval
\begin{equation}\label{eq:metric0tree}
M_0(\beta):=\big(M^\circ_0(\beta),\mu_0(\beta)\big):=\bigg(\big[0,\sD(\beta)\big]\ ,\ \ \sum_{U\in\beta}\Leb(U)\delta\big(W_\beta(U)\big)\bigg),
\end{equation}  
where $W_\beta(U)=\sD(\beta)-\sD_\beta(U)$, $U\in\beta$. Similarly, we associate with a $k$-tree $R_k=(\ft_k,(x_j^{(k)},j\in[k]),(\beta_E^{(k)},E\in\ft_k))\in\widebar{\bT}_k^{\rm int}$, 
the compact set
\begin{equation}\label{eq:skcirc}
S^\circ_k(R_k):=\bigcup_{E\in\ft_k}\left(W_k(E)+\left[0,\sD\big(\beta_E^{(k)}\big)\right]e_E\right)\subset[0,\infty)^{\ft_k}
\end{equation}
equipped with the $\ell_1$-distance $d_{\ell_1}$, where $W_k([k])=0$ and,  for $E\in\ft_k\setminus{\{[k]\}}$, $W_k(E)=W_k\big(\parent{E}\big)+\sD\big(\beta_{\parent{E}}^{(k)}\big)e_{\parent{E}}$.
We further equip $S^\circ_k(R_k)$ with a measure and let
\begin{equation}\label{eq:sk}
S_k(R_k):=\big(S^\circ_k(R_k),d_{\ell_1},0,\mu_k\big),\quad\mbox{with }\mu_k:=\sum_{\ell\in\block(R_k)}\|\ell\|\delta\big(W_k(\ell)\big),
\end{equation}
where we define the location of block $\ell$ in $S^\circ_k(R_k)$ as $W_k(\ell)=W_k\big(\parent{\{j\}}\big)+\sD\big(\beta_{\parent{\{j\}}}^{(k)}\big)e_{\parent{\{j\}}}$ for top blocks labelled by $\ell=j\in[k]$
and $W_k(\ell)=W_k(E)+W_{\beta^{(k)}_E}(U)e_E$ for blocks $\ell=(E,a,b)$ for $(a,b)=U\in\beta_E^{(k)}$, $E\in\ft_k$.

  Let $\tau\colon\bigcup_{k\ge 1}\widebar{\bT}_k^{\rm int}\rightarrow\bT^{\rm real}$ be the function that assigns to $R_k\in\widebar{\bT}_k^{\rm int}$ the ${\rm GHP}$-isometry class of 
  $S_k(R_k)$. 
\end{definition}

\begin{definition}\label{def:S} Let $\widebar{\bT}_\infty^{\rm int}$ be the subset of $\prod_{k\ge 1}\widebar{\bT}_k^{\rm int}$ of all consistent families $\mathrm{R}=(R_k,k\ge 1)$. We define a function 
  $S\colon\widebar{\bT}_{\infty}^{\rm int}\rightarrow\mathbb{T}^{\rm real}$, 
  \[
  S(\mathrm{R})=\left\{\begin{array}{ll}\lim_{k\rightarrow\infty}\tau(R_k),&\mbox{if this limit exists in }(\bT^{\rm real},d_{\rm GHP})\\[0.1cm] 
					         \Upsilon&\mbox{otherwise,}\end{array}\right.
  \] 
  where $\Upsilon\in\bT^{\rm real}$ is the ${\rm GHP}$-isometry class of
  the one-point tree equipped with the zero measure, $(\{0\},0,0,0)$.
\end{definition}
 
\begin{proposition}\label{thm:S} The map $S\colon\widebar{\bT}_{\infty}^{\rm int}\rightarrow\mathbb{T}^{\rm real}$ is Borel measurable.
\end{proposition}
\begin{proof} By \cite[Theorem 2.5(a)--(b)]{Paper1-0}, the map $M_0\colon\cI\rightarrow\mathcal{M}$ is continuous, where the space $\mathcal{M}$ of 
  \eqref{eq:onebranch} is equipped with the Hausdorff--Prokhorov metric. The function that projects $(C,\nu)\in M_0(\cI)$ onto the ${\rm GHP}$-isometry class of 
  $(C,|\cdot|,0,\nu)$ is clearly Lipschitz continuous. An induction shows that 
  $\tau\colon\bigcup_{k\ge 1}\widebar{\bT}_k^{\rm int}\rightarrow\bT^{\rm real}$ is also continuous on each part $\{\ft\}\times[0,\infty)^k\times\cI^\ft$ of the partition of $\widebar{\bT}_k^{\rm int}$ according to tree
  shape. Then $S$ is Borel measurable as a limit of Borel measurable functions.
\end{proof}

We will apply $S$ to the consistent families of unit-mass and self-similar $k$-tree evolutions. As defined for $\cI$-valued Markov processes above Proposition \ref{prop:type01:diffusion} and noted for $k$-tree evolutions in Theorem \ref{thm:Markov}, self-similarity of $(\cT_{k,+}^y,y\ge 0)$ starting from $\cT_{k,+}^0=\big(\ft_k,(x_j^{(k)},j\in[k]),(\beta_E^{(k)},E\in\ft)\big)\in\mathbb{T}_k^{\rm int}$ means that $(c\cT_{k,+}^{y/c},y\ge 0)$ is also a $k$-tree evolutions, starting from the scaled initial tree
\begin{equation}\label{eq:ch7:scaling}
c\cT^0_{k,+}:=\big(\ft_k,(cx_j^{(k)},j\in[k]),(c\beta_E^{(k)},E\in\ft_k)\big),
\end{equation}
in which all block masses, i.e.\ both top masses $x_j^{(k)}$ and blocks $U\in\beta_E^{(k)}$ are scaled by $c>0$. This scaling of block masses and time is 
naturally consistent when applied to consistent families of $k$-tree evolutions. As noted in \cite[Lemma 3.3]{Paper1-0}, the effect of this scaling of masses by $c$ is easily seen from Definition \ref{def:diversity} to induce a scaling of diversities by $\sqrt{c}$. In particular, if $\mathrm{T}=(T,d,\rho,\mu)$ is a representative of 
$S(\cT_{k,+},k\ge 1)$, then $c\mathrm{T}:=(T,\sqrt{c}d,\rho,c\mu)$ is a representative of $S(c\cT_{k,+}^0,k\ge 1)$. This ties in with the natural notion of scaling of Brownian excursions (encoding Brownian CRTs as in Section \ref{sec:BCRT}) that scales Brownian motion space (distances in the CRT) by $\sqrt{c}$ when scaling Brownian motion time (masses in the CRT) by $c$. 
We further observe that for any two rooted, weighted $\mathbb{R}$-trees $
\mathrm{T}$ and $\mathrm{T}^\prime$, we have
\begin{equation}\label{eq:ch7:scaling2}
  d_{\rm GHP}(c\mathrm{T},c\mathrm{T}^\prime)\le\max\{c,\sqrt{c}\}d_{\rm GHP}(\mathrm{T},\mathrm{T}^\prime).
\end{equation}

The map $S$ is not one-to-one, for instance, because $S(\mathrm{R})$ is invariant under permutations of labels, in the sense that consistently permuting labels $1,\ldots,m$
in $R_k$, $k\ge m$, for some $m\ge 2$, does not change $\tau(R_k)$, $k\ge m$. However, if we suitably enrich $S(\mathrm{R})$ by a sequence of marked points, we will be able to establish a partial inverse of $S$, in Theorem \ref{thm:realtoint} below.   

We follow \cite{M09,RW2} and extend Definition \ref{def:Rtree} to consider rooted, weighted compact $\mathbb{R}$-trees equipped with a sequence of marked points and to define 
\begin{align*}&d_{\rm GHP}^{\infty}\Big(\big(T,d,\rho,\mu,(\sigma_j,j\ge 1)\big),\big(T^\prime,d^\prime,\rho^\prime,\mu^\prime,(\sigma_j^\prime,j\ge 1)\big)\Big)\\
   &=\sum_{k\ge 1}2^{-k}d_{\rm GHP}^{[k]}\Big(\big(T,d,\rho,\mu,(\sigma_1,\ldots,\sigma_k)\big),\big(T^\prime,d^\prime,\rho^\prime,\mu^\prime,(\sigma_1^\prime,\ldots,\sigma_k^\prime)\big)\Big),
\end{align*}
with
\begin{align*}&d_{\rm GHP}^{[k]}\Big(\big(T,d,\sigma_0,\mu,(\sigma_1,\ldots,\sigma_k)\big),\big(T^\prime,d^\prime,\sigma_0^\prime,\mu^\prime,(\sigma_1^\prime,\ldots,\sigma_k^\prime)\big)\Big)\\
  &=\inf_{\phi,\phi^\prime}\left\{\max\left\{d_M^{\rm H}(\phi(T),\phi^\prime(T^\prime)),\ d_M^{\rm P}(\phi_*\mu,\phi^\prime_*\mu^\prime),\ \max_{0\le i\le k}d_M(\phi(\sigma_i),\phi^\prime(\sigma_i^\prime))\right\}\right\},
\end{align*}
where the infimum is over all metric spaces $(M,d_M)$ and all injective isometries $\phi\colon(T,d)\rightarrow(M,d_M)$ and $\phi^\prime\colon (T^\prime,d^\prime)\rightarrow (M,d_M)$. A ${\rm GHP}$-isometry $\iota\colon(T,d)\rightarrow(T^\prime,d^\prime)$ is a ${\rm GHP}^\infty$-isometry if furthermore $\iota(\sigma_j)=\sigma_j^\prime$ for
all $j\ge 1$. Then $d_{\rm GHP}^\infty$ can be viewed as a metric on the set $\bT_\infty^{\rm real}$ of ${\rm GHP}^\infty$-isometry classes. 

Recall from Section \ref{sec:BCRT} the definition of a Brownian CRT. Specifically, the line-breaking construction and the construction from a Brownian excursion yield random 
$\mathbb{R}$-trees whose respective projections to their ${\rm GHP}$-isometry classes in $\mathbb{T}^{\rm real}$ have the same distribution \cite{AldousCRT3,LeGall05}. Taking a sample from the weight measure of the $\mathbb{R}$-tree in either representative is straightforward and gives rise to a random rooted weighted $\mathbb{R}$-tree equipped with a sequence of
random marked points, which are almost surely leaves. The following result formalizes the idea that the distribution on $\bT_\infty^{\rm real}$ of its ${\rm GHP}^\infty$-isometry class does not depend on the 
choice of representative.

\begin{proposition}\label{prop:sampleleaves} There is a natural stochastic kernel $\mathbf{m}_\infty$ from $\mathbb{T}^{\rm real}$ to $\mathbb{T}^{\rm real}_\infty$ such that 
  $\mathbf{m}_\infty(\mathrm{T},\cdot)$ can be considered as the distribution on $\mathbb{T}^{\rm real}_\infty$ of the space $\mathrm{T}=(T,d,\rho,\mu)$ equipped with a sequence 
  of independent identically $\mu/\|\mu\|$-distributed marked points, if $\mu\neq 0$ and where $\|\mu\|=\mu(T)$ is the (finite) total mass of $\mu$. To be definite, we use the convention that we sample from $\delta_\rho$ if $\mu=0$.
\end{proposition}
\begin{proof} The proof of Lemma 13 from Miermont \cite{M09} for any finite number of marked points applies mutatis mutandis. 
\end{proof}

\begin{theorem}\label{thm:realtoint} There is a natural measurable map $R\colon\mathbb{T}_\infty^{\rm real}\rightarrow\widebar{\mathbb{T}}_\infty^{\rm int}$ such that 
  \[
  S(R(\mathrm{T},\boldsymbol{\sigma}))=\mathrm{T}
  \]
  for $\mathbf{m}_\infty(\mathrm{T},d\boldsymbol{\sigma})\mathbb{P}(\cT\!\in\! d\mathrm{T})$-a.e.\ $(\mathrm{T},\boldsymbol{\sigma})\!=\!(\mathrm{T},(\sigma_j,j\!\ge\! 1))$, where $\cT$ is a Brownian CRT.
\end{theorem}

We make the map $R$ and the kernel $\mathbf{m}_\infty$ explicit and prove this theorem in Appendix \ref{sec:contproj}. In the following, we will use the image $\mathcal{L}$ of $\mathbf{m}_\infty(\mathrm{T},\cdot)$ under $R$ to sample a random system of consistent $k$-trees associated with $\mathrm{T}\in\mathbb{T}^{\rm real}$. Indeed, $\mathcal{L}$ is a kernel from $\mathbb{T}^{\rm real}$ to $\widebar{\bT}_\infty^{\rm int}$. In particular, we
can carry out the construction of a consistent family of Brownian reduced $k$-trees from a $\mathbb{T}^{\rm real}$-valued Brownian CRT $\cT$ by sampling from
$\mathcal{L}(\cT,\cdot)$, on a suitably enlarged probability space.



\section{The Markov property of the projective continuum tree limit}\label{sec:defmarkov}

Now that we have formally introduced all ingredients, let us make precise the less formal definition of the Aldous diffusion given in the introduction.

\begin{definition}\label{df:adformal}
  Let $(\widebar{\cT}_{k,+}^s,s\ge 0)$ be a consistent family of stationary unit-mass resampling $k$-tree evolutions, $k\ge 1$, as in Corollary 
  \ref{cor:consistent_fam}(iii). Then we define the \em Aldous diffusion \em as a ${\rm GHP}$-path-continuous modification of the process
  $\widebar{\cT}(s)=S\big(\widebar{\cT}_{k,+}^s,k\ge 1\big)$, $s\ge 0$, where 
  $S\colon\widebar{\mathbb{T}}_\infty^{\rm int}\rightarrow\mathbb{T}^{\rm real}$ is as defined in Definition \ref{def:S}. 

  Similarly, given consistent pseudo-stationary self-similar resampling $k$-tree evolutions $(\cT_{k,+}^y,y\ge 0)$, $k\ge 1$, as in Corollary
  \ref{cor:consistent_fam}(ii), we define the \em self-similar Aldous diffusion \em as a ${\rm GHP}$-path-continuous modification of
  ${\cT}(y)=S\big(\cT_{k,+}^y,k\ge 1\big)$, $y\ge 0$. 
\end{definition}

Indeed, we can view the consistent family of $k$-tree evolutions as a single evolution in the subset $\widebar{\mathbb{T}}_\infty^{\rm int}$ of the product space $\prod_{k\ge 1}\widetilde{\bT}_{k}^{\rm int}$. Since the stationary distribution is the distribution of a consistent family of Brownian reduced $k$-trees, $k\ge 1$, and Theorem \ref{thm:realtoint} confirms that mapping consistent families of Brownian reduced $k$-trees under $S$ returns Brownian CRTs, the 
process $(\widebar{\cT}(s),\,s\ge 0)$ is well-defined. In this section, we establish the Markov property of this process. To justify calling it a diffusion, we show in Corollary \ref{cor:ad:contversion} that it has a ${\rm GHP}$-path-continuous modification hence establishing the existence of the Aldous diffusion. We will, however, argue that for this Markov process, the strong Markov property fails, in general. 

To establish the Markov property, let us think about the transition mechanism. In the construction, the $k$-trees form, at all times, a consistent family of Brownian reduced $k$-trees. Given a 
Brownian continuum random tree as initial state, we obtain an evolution, as follows. First, we sample a random initial system of consistent reduced $k$-trees, then we use the consistent evolution of those, and finally we consider the limiting continuum tree induced by the consistent system of $k$-trees, at time $s$, as the state at time $s\ge 0$. A priori, such a construction may not yield the Markov property, since a time-$s$ transition followed by a time-$r$ transition 
(using a newly sampled consistent system of reduced $k$-trees from the same time-$s$ CRT), may not give the same time-$(s+r)$ distribution as a time-$(s+r)$ transition (without sampling new $k$-trees at time $s$). The following result expresses the idea that 
the target state as a continuum tree does not depend on the choice of sampled $k$-trees.

\begin{proposition}\label{couplelabels} Let $(\cT,d,\rho,\mu)$ be a Brownian CRT and $\widebar{\cT}_{\!\bullet}^0=(\widebar{\cT}_{\!k,\bullet}^0,\,k\ge 1)$ and $\widebar{\cT}_{\!\circ}^0=(\widebar{\cT}_{\!k,\circ}^0,\,k\ge 1)$ two
  families of reduced $k$-trees of the same Brownian CRT $\cT$,  independently sampled according to $\mu$, in the sense of Proposition \ref{prop:sampleleaves} and Theorem \ref{thm:realtoint} Then there are two coupled families $(\widebar{\cT}^s_{\!\bullet},\,s\ge 0)$ and $(\widebar{\cT}^s_{\!\circ},\,s\ge 0)$ of 
  consistent unit-mass resampling $k$-tree evolutions, $k\ge 1$, for which 
  $$\big(S(\widebar{\cT}_{\!\bullet}^s),s\ge 0\big)=\big(S(\widebar{\cT}_{\!\circ}^s),s\ge 0\big)\quad\mbox{a.s.}$$   
\end{proposition}

Before we prove this, we establish a more elementary lemma.

\begin{lemma}\label{lmperm} Consider a self-similar resampling $k$-tree evolution $(\cT_{k,+}^y,\,y\ge 0)$ with resampling times $D_m$, $m\ge 1$, and any permutation $p_0$ of $[k]$. Then there is a sequence 
  $(p_m,\,m\ge 1)$ of random permutations, on the same probability space, such that $\widehat{\cT}_{k,+}^y:=p_m\cT_{k,+}^y$, $D_m\le y<D_{m+1}$, $m\ge 0$, defines a resampling $k$-tree evolution 
  $\big(\widehat{\cT}_{k,+}^y,\,y\ge 0\big)$ with the same resampling times $D_m$, $m\ge 1$.

  The same holds for a non-resampling $k$-tree evolution $(\cT_{k,-}^y,\,y\ge 0)$ with degeneration times $D_1,\ldots,D_k$.
\end{lemma}
\begin{proof} It follows from elementary symmetry properties of killed $k$-tree evolutions of Definition \ref{def:killed_ktree} that $(p_0\cT_{k,\pm}^y,\,0\le y<D_1)$ is also a killed  
  $k$-tree evolution in the sense of that definition. At each resampling (or degeneration) time $D_m$, the swapping part of the swap-reduction function $\varrho$ defined in Section \ref{sec:non_resamp_def} may yield different transpositions, but  
  appropriately composing $p_{m-1}$ with these transpositions, if any, yields a new permutation $p_m$ with the desired properties. Specifically, symmetry properties of the resampling kernel defined in Section \ref{sec:resamp_def} are such
  that we can achieve that $\widehat{\cT}_{k,+}^{D_m}=p_m\cT_{k,+}^{D_m}$ has performed the resampling step as required. By induction, $\big(\widehat{\cT}_{k,+}^y,\,y\ge 0\big)$ is a  resampling  $k$-tree evolution with
  resampling times $D_m$, $m\ge 1$. We conclude similarly in the non-resampling case. 
\end{proof}
  
\begin{proof}[Proof of Proposition \ref{couplelabels}] Suppose that $\cT_\bullet^0:=\widebar{\cT}_{\!\bullet}^{0}$ and $\cT_\circ^0:=\widebar{\cT}_{\!\circ}^{0}$ are associated with 
  leaf samples  $(\Sigma_{k,\bullet},\,k\ge 1)$ and $(\Sigma_{k,\circ},\,k\ge 1)$. Define the merged sample by alternating $\Sigma_{2k-1}=\Sigma_{k,\bullet}$ and 
  $\Sigma_{2k}=\Sigma_{k,\circ}$, $k\ge 1$, and consider a family of consistent pseudo-stationary resampling $k$-tree evolutions $(\cT_k^y,y\ge 0)$, $k\ge 1$, starting from 
  the associated consistent system of Brownian reduced $k$-trees.
  
  Now fix $k\ge 2$. Consider the permutation $p_0^\circ$ of $[2k]$ given by $p_0^\circ(2i)=i$, $p_0^\circ(2i-1)=k+i$, $i\in[k]$, and the process $\big(\widehat{\cT}_{2k}^y,\,y\ge 0\big)$ constructed as in 
  Lemma \ref{lmperm}, with permutations $(p_m^\circ,\,m\ge 0)$ and resampling times $(D_m,m\ge 1)$. We define the projection 
  $\cT_{k,\circ}^y=\pi_{k}(\widehat{\cT}_{2k}^y)$, $y\ge 0$, to obtain a $k$-tree evolution. By the invariance of 
  Brownian reduced $2k$-trees under permutations, $\big(\widehat{\cT}_{2k}^y,\,y\ge 0\big)$ is a pseudo-stationary  resampling $2k$-tree evolution. By Theorem \ref{thm:consistency}, $\big(\cT_{k,\circ}^y,\,y\ge 0\big)$ is a 
  pseudo-stationary resampling $k$-tree evolution. Starting from the permutation $p_0^\bullet(2i-1)=i$, $p_0^\bullet(2i)=k+i$, $i\in[k]$, we similarly define $(\cT_{k,\bullet}^y,\,y\ge 0)$ with 
  permutations $(p_m^\bullet,\,m\ge 0)$.
  
  By further projection to $[m]\subseteq[k]$, we obtain consistent evolutions by Theorem \ref{thm:consistency}, and by Kolmogorov's consistency theorem and de-Poissonization, 
  we obtain three coupled families $\big(\widebar{\cT}^{s}_{\!\bullet},\,s\ge 0\big)$, $\big(\widebar{\cT}^{s}_{\!\circ},\,s\ge 0\big)$, and $\big(\big(\widebar{\cT}_{\!k}^s,\,s\ge 0\big),\,k\ge 1\big)$. Denote by $\widebar{\cT}(s)=S\big(\widebar{\cT}_{\!k}^s,\,k\ge 1\big)$ the Brownian CRT associated with the
  third family. As  $\widebar{\cT}^s_{\!k,\circ}$ and $\widebar{\cT}^s_{\!k,\bullet}$ are projections of $\widebar{\cT}^s_{\!2k}$, their weighted $\mathbb{R}$-tree representations are projections of $S_{2k}(\widebar{\cT}^s_{\!2k})$ for all $k\ge 1$, and can all be viewed (up to ${\rm GHP}$-isometry) as projected subtrees of (any representative of) 
$\widebar{\cT}(s)$. Hence, the projective systems of $\bR$-trees a.s.\ increase to subsets $S\big(\widebar{\cT}^s_{\!\circ}\big)$ and $S\big(\widebar{\cT}^s_{\!\bullet}\big)$ of 
  $\widebar{\cT}(s)$, equipped with projected mass measures. But the mass measures of these three Brownian CRTs are diffuse and charge all fringe subtrees \cite{AldousCRT3}, hence they must be equal.     
\end{proof}


\newcommand{\BCRT}{\ensuremath{\mathtt{BCRT}}}
\newcommand{\TReal}{\bT^{\textnormal{real}}}
\newcommand{\rR}{\mathrm{R}}

Before proceeding, we find it useful to reframe the preceding result in the language of stochastic kernels.
\begin{itemize}
 \item We denote the law on $\bT^{\rm real}$ of the Brownian CRT by \BCRT.
 \item Let $\widetilde S$ denote the stochastic kernel from $\widebar{\mathbb{T}}^{\rm int}_\infty$ to $\TReal$ associated with the map $S$, i.e.\ $\widetilde S(\mathrm{R},\cdot) = \delta_{S(\mathrm{R})}(\,\cdot\,)$. Thus, $\cL\widetilde S$ is the identity kernel on a \BCRT -a.s.\ subset of $\TReal$.
 \item Let $Q_\infty$ denote the law on $\widebar{\bT}^{\rm int}_\infty$ of a projectively consistent system of Brownian reduced $k$-trees, i.e.\ 
  $Q_\infty(\,\cdot\,) = \int \cL(\mathrm{T},\cdot\,)\BCRT(d\mathrm{T})$.
 \item 
 For $u\ge0$, let $\kappa_u$ denote the time-$u$ transition kernel for the Markov processes $\big(\widebar{\cT}^{s}_{\!\bullet},\,s\ge0\big)$ and $\big(\widebar{\cT}^{s}_{\!\circ},\,s\ge0\big)$ of Proposition \ref{couplelabels}. Following Proposition \ref{prop:Lusin}, $\widebar{\mathbb{T}}^{\rm int}_\infty$ is a subspace of a countable product of Borel spaces and hence is Borel; thus, following \cite[Theorem 6.3]{Kallenberg}, such transition kernels exist.

 Equivalently, these are the transition kernels for the stationary, consistent family of unit-mass resampling $k$-tree evolutions described in Corollary \ref{cor:consistent_fam}\ref{item:cnstfam:dePoi}.
\end{itemize}

Recall from \eqref{eq:composition_of_kernels} the convention of left-to-right composition of Markov kernels $\kappa_u\kappa_v$, as distinct from the right-to-left notation for composition of functions $g\circ f$:\vspace{-0.1cm}
\begin{equation}
 \int_{\widebar{\mathbb{T}}^{\rm int}_\infty} \cL\kappa_1(\mathrm{T},d\mathrm{R})f(\mathrm{R}) = \int_{\widebar{\mathbb{T}}^{\rm int}_\infty}\cL(\mathrm{T},d\mathrm{R}')\int_{\widebar{\mathbb{T}}^{\rm int}_\infty}\kappa_1(\mathrm{R}',d\mathrm{R})f(\mathrm{R}).\vspace{-0.1cm}
\end{equation}
Our goal is to study the process $\big(S\big(\overline{\cT}^{\,s}_{\!\bullet}\big),\,s\ge 0\big)$ of Proposition \ref{couplelabels}, which is our proposed Aldous diffusion. To that end, we reformulate Proposition \ref{couplelabels} as follows.

\begin{corollary}\label{cor:couple_labels_ker}
 For each $u\ge0$ there exists a $Q_\infty$-a.s.\ domain $B_u\subset \widebar{\mathbb{T}}^{\rm int}_\infty$ on which $\widetilde S\cL\kappa_u\widetilde S=\kappa_u \widetilde S$, i.e.\vspace{-0.1cm}
 \begin{equation}
   \int f(\mathrm{T})\widetilde S\cL\kappa_u \widetilde S(\mathrm{R},d\mathrm{T})=\int f(\mathrm{T})\kappa_u \widetilde S(\mathrm{R},d\mathrm{T})\vspace{-0.1cm}
 \end{equation}
 for all bounded, measurable functions $f\colon \TReal\to \bR$ and all $\mathrm{R}\in B_u$.\pagebreak
\end{corollary}

\begin{proof}
 We copy the notation of Proposition \ref{couplelabels}. By that proposition, for any $u\ge0$ and any bounded, measurable $g\colon \widebar{\mathbb{T}}^{\rm int}_\infty\times\TReal\to \BR$,
 \begin{equation*}
  \bE\big[ g(\widebar{\cT}_{\!\bullet}^{\,0},S\big( \widebar{\cT}_{\!\circ}^{\,u}\big)\big) \big]
  	= \bE\big[ g(\widebar{\cT}_{\!\bullet}^{\,0},S\big( \widebar{\cT}_{\!\bullet}^{\,u}\big)\big) \big].
 \end{equation*}
 Note that $\widebar{\cT}_{\!\bullet}^{\,0}$ and $\widebar{\cT}_{\!\circ}^{\,0}$ each have law $Q_\infty$, and $\widebar{\cT}_{\!\circ}^{\,0}$ has conditional law $\widetilde S\cL(\widebar{\cT}_{\!\bullet}^{\,0},\cdot\,)$ given $\widebar{\cT}_{\!\bullet}^{\,0}$. Hence, the above formula is equivalent to
 \begin{equation*}
  \iint g(\mathrm{R},\mathrm{T})\widetilde S\cL\kappa_u\widetilde S(\mathrm{R},d\mathrm{T})Q_\infty(d\mathrm{R})
  	= \iint g(\mathrm{R},\mathrm{T})\kappa_u\widetilde S(\mathrm{R},d\mathrm{T}) Q_\infty(d\mathrm{R}).
 \end{equation*}
 The corollary follows by the a.s.\ uniqueness of regular conditional distributions \cite[Theorem 6.3]{Kallenberg}.
\end{proof}

\begin{theorem}\label{thm:ad:markov}
The Aldous diffusion has the simple Markov property.
\end{theorem}

\begin{proof}
 Fix $0=s_0<s_1<\cdots <s_{k+1}$ with $\Delta_j := s_j - s_{j-1}$, $j\in [k+1]$, and let $f_0,\dots, f_k,f\colon\TReal\to \bR$ be bounded, measurable functions. 
 Kallenberg \cite[Corollary 8.3]{Kallenberg} observes the following \emph{a.s.\ semi-group property} for general Markov processes on Borel spaces: for every $u,v\ge 0$,
 \begin{equation}\label{eq:a.s._semigroup}
  \kappa_u\kappa_v(\mathrm{R},\cdot\,) = \kappa_{u+v}(\mathrm{R},\cdot\,)\quad \text{for }Q_\infty\text{-a.e }\mathrm{R}\in\widebar{\mathbb{T}}^{\rm int}_\infty.
 \end{equation}
 Hence, there exists a $Q_\infty$-a.s.\ set $C\subset\widebar{\mathbb{T}}^{\rm int}_\infty$ such that
 \begin{equation*}
  \kappa_{\Delta_1}\kappa_{\Delta_2}\ldots\kappa_{\Delta_j}(\mathrm{R},\cdot\,) = \kappa_{s_j}(\mathrm{R},\cdot\,)\quad \text{for all }j\in [k+1],\ \mathrm{R}\in C.
 \end{equation*}
 For the Aldous diffusion $\big(S\big(\overline{\cT}^{\,s}_{\!\bullet}\big),\,s\geq 0\big)$, as defined above, this yields  
 \begin{equation*}
 \begin{split}
  \bE\left[ f\big(S\big(\overline{\cT}^{\,s_{k+1}}_{\!\bullet}\big)\big)\prod_{j=0}^{k} f_j\big(S\big(\overline{\cT}^{\,s_j}_{\!\bullet}\big)\big) \right]  
  	&= \int_{C} Q_\infty(d\rR_0)\int_{\widebar{\mathbb{T}}^{\rm int}_\infty}\kappa_{\Delta_1}(\rR_0,d\rR_1)(f_1\circ S)(\rR_1)\\
  		&\qquad \cdots\int_{\widebar{\mathbb{T}}^{\rm int}_\infty}\kappa_{\Delta_k}(\rR_{k-1},d\rR_k)(f_k\circ S)(\rR_k)\\
		&\qquad \int_{\widebar{\mathbb{T}}^{\rm int}_\infty}\kappa_{\Delta_{k+1}}(\rR_k,d\rR_{k+1})(f\circ S)(\rR_{k+1}).
 \end{split}
 \end{equation*}
 Let $B_{\Delta_{k+1}} \subset \widebar{\mathbb{T}}^{\rm int}_\infty$ be as in Corollary \ref{cor:couple_labels_ker}. The stationarity of our Markov process implies that $B_{\Delta_{k+1}}$ has full measure under the $\kappa_{\Delta_1}\ldots\kappa_{\Delta_k}$-image of $Q_\infty$. Hence, in the context of that formula, the two innermost integrals can be rewritten as
 \begin{equation*}
 \begin{split}
  &\int_{B_{\Delta_{k+1}}}\kappa_{\Delta_k}(\rR_{k-1},d\rR_k)(f_k\circ S)(\rR_k)
  	\int_{\widebar{\mathbb{T}}^{\rm int}_\infty}\kappa_{\Delta_{k+1}}(\rR_k,d\rR_{k+1})(f\circ S)(\rR_{k+1})\\
  	&\quad = \int_{B_{\Delta_{k+1}}}\kappa_{\Delta_k}(\rR_{k-1},d\rR_k)(f_k\circ S)(\rR_k)
  	\int_{\TReal}\widetilde S\cL\kappa_{\Delta_{k+1}}\widetilde S(\rR_k,d\mathrm{T}_{k+1})f(\mathrm{T}_{k+1})\\
  	&\quad = \int_{\TReal}\kappa_{\Delta_k}\widetilde S(\rR_{k-1},d\mathrm{T}_k)f_k(\mathrm{T}_k)
  	\int_{\TReal}\cL\kappa_{\Delta_{k+1}}\widetilde S(\mathrm{T}_k,d\mathrm{T}_{k+1})f(\mathrm{T}_{k+1}),
 \end{split}
 \end{equation*}
 with the last line following by a superficial rearrangement of kernels and the observation, again, that $Q_\infty(B_{\Delta_{k+1}}) = 1$.
Plugging this back in, we get
 \begin{equation*}
 \begin{split}
  &\bE\left[ f\big(S\big(\overline{\cT}^{\,s_{k+1}}_{\!\bullet}\big)\big)\prod_{j=0}^{k} f_j\big( S\big(\overline{\cT}^{\,s_j}_{\!\bullet}\big) \big)\right]
  	= \bE\left[ \widetilde\kappa_{\Delta_{k+1}}\big(S\big(\overline{\cT}^{\,s_k}_{\!\bullet}\big),f\big)\prod_{j=0}^{k} f_j\big( S\big(\overline{\cT}^{\,s_j}_{\!\bullet}\big)\big) \right]\\
  	&\qquad\qquad\qquad = \bE\left[ \bE\Big[ f\big(S\big( \overline\cT^{\,s_{k+1}}_{\!\bullet} \big)\big)\, \Big|\, S\big( \overline{\cT}^{\,s_k}_{\!\bullet} \big)\Big]\prod_{j=0}^{k} f_j\big( S\big(\overline{\cT}^{\,s_j}_{\!\bullet}\big)\big) \right],
 \end{split}
 \end{equation*}
 where $\widetilde\kappa_{\Delta_{k+1}} := \cL\kappa_{\Delta_{k+1}}\widetilde S$. 
 Monotone class arguments allow to further extend this form of the simple Markov property, see e.g.\ \cite[Lemma 8.1]{Kallenberg}.
\end{proof}

\begin{remark} The simple Markov property of the self-similar Aldous diffusion follows by similar arguments.
\end{remark}


We now claim that the Aldous diffusion is not a strong Markov process. Informally, we see this by considering the first time that a ternary branch point with four 
large component masses is formed. This can be set up as a stopping time. Before this time, there was an edge separating the root component and three other components into two pairs of components. This edge, before de-Poissonization performing a type-0 evolutions with total mass process $\besq(1)$, has just shrunk to zero mass, as $\besq(1)$ does, but this is not a degeneration since edges only degenerate when one of their top masses also vanishes. After this time, the same arrangement into pairs of labels persists. However, the state in $\bT^{\rm real}$ at this time does not contain the information about the pairing. The following remark formalizes this.

\begin{remark}\label{prop:notstrongmarkov}
  The Aldous diffusion in the state space $(\bT^{\rm real},d_{\rm GHP})$ does not have the strong Markov property. Indeed, we will argue that the strong Markov property fails at a stopping time at which the Aldous diffusion $\big(\widebar{\cT}(s),\,s\ge 0\big)$ possesses a degree-4 branch point. 
  Specifically, recall that for each $s\ge 0$, the Brownian CRT $\widebar{\cT}(s)$ has only degree-3 branch points almost surely \cite{AldousCRT3}. In particular, a Brownian reduced 3-tree consists of 
  three top masses and two edge partitions that all have positive mass almost surely. The shape is necessarily a type-1 edge between the root and a branch point and a type-2 edge 
  above. When the mass of the edge partition of the type-2 edge vanishes, the corresponding branch points in the CRT coincide hence forming a degree-4 branch point.

  More precisely, it will be convenient to explore this in the context of the construction $\widebar{\cT}(s):=S\big(\widebar{\cT}^s_k,\,k\ge 1\big)$, $s\ge 0$, of the Aldous diffusion from consistent
  resampling unit-mass $k$-tree evolutions $\big(\widebar{\cT}^s_k,\,s\ge 0\big)$, $k\ge 1$. It is well-known (e.g.\ as a consequence of \cite[Theorem 2]{Aldous94} and the sampling properties of Dirichlet distributions) that the vector of five masses obtained from a Brownian reduced 3-tree has a 
  ${\tt Dirichlet}(\frac12,\frac12,\frac12,\frac12,\frac12)$ distribution. We consider the event $\widebar{A}_1^{\rm int}$ that the top masses of the type-2 edge of $\widebar{\cT}_3^0$ each exceed $\frac13$ and that the remaining top mass and the type-1 edge mass each exceed 
  $\frac19$. (We will later replace these thresholds by $\frac{12}{37},\frac{10}{37},\frac{8}{37},\frac{6}{37}$ for technical reasons, but the principle is the same.) This 
  event has positive probability. Similarly, consider the event $\widebar{A}_2^{\rm int}$ that two top masses including the one of the type-1 edge exceed 
  the larger thresholds, and that the remaining top mass and the type-1 edge mass exceed the smaller thresholds, respectively. For each permutation of the thresholds, this event has the same 
  positive probability.  

  By Corollary \ref{cor:WF}, the evolution of the five masses under the 3-tree evolution stopped when the first component vanishes is a Wright--Fisher diffusion with three parameters
  $-\frac12$ and two $\frac12$. In particular, Pal's \cite{Pal13} construction from squared Bessel processes, here ${\tt BESQ}(-1)$ and ${\tt BESQ}(1)$, easily yields that there is positive
  probability that the type-2 edge mass vanishes before any of the other four masses violates its constraint to lie above their respective thresholds. We denote by $\sigma^{\rm int}_i$ the first time when either the zero mass is attained or one of the four mass constraints fails, for the process starting in $\widebar{A}_i^{\rm int}$, respectively, $i=1,2$.

  If at $\sigma^{\rm int}_1$ or $\sigma^{\rm int}_2$, the type-2 edge mass vanishes, this is not a degeneration time for the resampling 3-tree evolution, and $\widebar{\cT}^{\sigma_1^{\rm int}}_3$ and $\widebar{\cT}^{\sigma_2^{\rm int}}_3$ retain the respective tree shape and the further evolution 
  preserves the position of the top masses exceeding the respective thresholds. We also oberve that this information is not retained in $S\big(\widebar{\cT}^{\sigma_1^{\rm int}}_k,\,k\ge 1\big)$ and 
  $S\big(\widebar{\cT}^{\sigma_2^{\rm int}}_k,\,k\ge 1\big)$,
  which will both have a degree-4 branch point with three subtrees corresponding to the top masses. The random times $\sigma_1^{\rm int}$ and $\sigma_2^{\rm int}$ are not stopping times in the
  filtration of the Aldous diffusion, so a formalisation of this argument will require us to define related times that are. This will also require some sample path regularity, so we postpone the further discussion to 
  Chapter \ref{chap8}.
\end{remark}

\section{The GHP-distance between $k$-trees is bounded by $d_\cI$-distances}\label{sec:dItoGHP}

Let $R_k\in\widebar{\bT}_{k}^{\rm int}=\widebar{\bT}_{[k]}^{\rm int}$ and $1\le i<j\le k$. The projection map $\pi_{\{i,j\}}$ defined in Definition \ref{def:proj} projects $\widebar{\bT}_{[k]}^{\rm int}$ to $\widebar{\bT}_{\{i,j\}}^{\rm int}$, which is isometric to the space $\cJ^\circ$ of \eqref{type2spaces}. As a slight variation, we denote by $\pi_{i,j}R_k\in\cI^\circ$ the interval partition representation $(0,x_i)\concat(0,x_j)\concat\beta$ associated with
$(x_i,x_j,\beta)=\pi_{\{i,j\}}R_k$ as in Remark \ref{prop:type2:IPvalued}. Note that $\pi_{j,i}R_k$ differs from $\pi_{i,j}R_k$ in the order of their two left-most blocks (if $x_i>0$ and $x_j>0$ are distinct). Let 
$S_k(R_k)$ and $S_k^\circ(R_k)$ denote the weighted $\bR$-tree and the $\mathbb{R}$-tree without the weight measure associated with $R_k$ in Definition \ref{def:SkandSkcirc}. Now let $R_k,R_k^\prime\in\widebar{\bT}_k^{\rm int}$. In this section we bound ${\rm GHP}$-distances between trees $S_k(R_k)$ and $S_k(R_k^\prime)$ by $d_\cI$-distances between interval partitions $\pi_{i,j}R_k$ and $\pi_{i,j}R_k^\prime$, $1\le i,j\le k$. 

Before turning to weighted $\mathbb{R}$-trees, we bound distances in the sense of the Gromov--Hausdorff distance without
Prokhorov component
\begin{equation}\label{eq:GH}
d_{\rm GH}(\mathrm{T},\mathrm{T}^\prime)=\inf_{\phi,\phi^\prime}\max\Big\{d_M^{\rm H}(\phi(T),\phi^\prime(T^\prime)),d_M(\phi(\rho),\phi^\prime(\rho^\prime))\Big\},
\end{equation}
where the infimum is taken over all metric spaces $(M,d_M)$ and all injective isometries $\phi\colon(T,d)\rightarrow(M,d_M)$ and 
$\phi^\prime\colon(T^\prime,d^\prime)\rightarrow(M,d_M)$. 
We say that $\mathrm{T}$ and $\mathrm{T}^\prime$ are \em ${\rm GH}$-isometric \em if there is a bijective isometry $\iota\colon(T,d)\rightarrow(T^\prime,d^\prime)$ such that $\iota(\rho)=\rho^\prime$.
Then $d_{\rm GH}$ can be viewed as a metric on the set $\mathbb{T}^{\rm real}_\circ$ of ${\rm GH}$-isometry classes of rooted compact $\mathbb{R}$-trees,
and there is a useful equivalent definition \cite{EPW,Gromov1999,M09}
\begin{equation}\label{eq:GH2}
d_{\rm GH}(\mathrm{T},\mathrm{T}^\prime)
=\inf\Big\{{\rm dis}_{\rm GH}(K)\colon K\,\mbox{ ${\rm GH}$-correspondence between $\mathrm{T}$ and $\mathrm{T}^\prime$}\Big\},
\end{equation}
where a \em ${\rm GH}$-correspondence \em between $\mathrm{T}=(T,d,\rho)$ and $\mathrm{T}^\prime=(T^\prime,d^\prime,\rho^\prime)$ is a subset 
$K\subseteq T\times T^\prime$ with $(\rho,\rho^\prime)\in K$, whose coordinate projections are surjective onto $T$ and $T^\prime$ respectively, and where the \em 
${\rm GH}$-distortion \em of $K$ is given by 
\[
{\rm dis}_{\rm GH}(K)=\frac{1}{2}\sup\Big\{|d(x,y)-d^\prime(x^\prime,y^\prime)|\colon (x,x^\prime),(y,y^\prime)\in K\Big\}.
\]
\begin{proposition}\label{propGH} Let $k\ge 2$. Consider $R_k,R_k^\prime\in\widebar{\bT}_{k}^{\rm int}$ with the same shape $\ft_k\in\bT_{k}^{\rm shape}$. 
  Consider the $\bR$-trees  $S_k^\circ(R_k)$ and $S_k^\circ(R_k^\prime)$ and for each $1\le i<j\le k$ the four interval partitions 
  $\pi_{i,j} R_k,\pi_{j,i}R_k, \pi_{i,j}R_k^\prime,\pi_{j,i}R_k^\prime\in\cI$.
    Then 
    $$d_{\rm GH}(S_k^\circ(R_k),S_k^\circ(R_k^\prime))\le 2\max_{1\le i<j\le k}\min\Big\{d_\cI(\pi_{i,j}R_k,\pi_{i,j}R_k^\prime),d_\cI(\pi_{j,i}R_k,\pi_{j,i}R_k^\prime)\Big\} .$$
\end{proposition}
\begin{proof} 
  We use notation $R_k\!=\!(\ft_k,(x_i,i\!\in\![k]),(\beta_E,E\!\in\!\ft_k))$ and simplify notation from Definition \ref{def:SkandSkcirc} for the locations
  $\Sigma_i=W_k(i)\in S_k^\circ(R_k)$ corresponding to the top mass labelled $i\in[k]$ and $b_E=W_k(E)+\sD(\beta_E)\in S_k^\circ(R_k)$ corresponding to 
  $E\in\ft_k$, i.e. the vertex at the top end (away from the root $\rho:=W_k([k])$) of the branch in $S_k^\circ(R_k)$ built from $\beta_E$. 
  This vertex typically has degree $3-g$ for a type-$g$ edge, $g=0,1,2$. Indeed, if $E=\{i,j\}\in\ft_k$ is a type-2 edge, then $b_E=\Sigma_i=\Sigma_j$; if $E=F\cup\{i\}\in\ft_k$, with $F\in\ft_k$ and $i\not\in F$, is a type-1 edge, then $b_E=\Sigma_i$; also, if $\beta_E=\varnothing$, then
  $b_E=b_{\parent{E}}$. We use similar notation $\Sigma_i^\prime,b_E^\prime,\rho^\prime\in S_k^\circ( R_k^\prime)$. In the following, we denote the metrics
  of $S_k^\circ(R_k)$ and $S_k^\circ(R_k^\prime)$ by $d$ and $d^\prime$, respectively. Now consider the 
  ${\rm GH}$-correspondence $K$ between $S_k^\circ(R_k)$ and $S_k^\circ( R_k^\prime)$ which consists of
  \begin{itemize}
  \item pairs of \em special \em vertices $(\rho,\rho^\prime)$ and $(b_E,b^\prime_E)$ for all $E\in\ft_k$,\vspace{0.1cm}
  \item and all pairs of points $(\lambda b_E+(1-\lambda)b_{\parent{E}},\lambda b_E^\prime+(1-\lambda)b_{\parent{E}}^\prime)$, $0<\lambda<1$, 
    on the branches between special vertices $b_E,b_{\parent{E}}$ and $b_E^\prime,b_{\parent{E}}^\prime$, associated with $E\in\ft_k$ 
    and its parent $\parent{E}$, with the convention that $b_{\parent{[k]}}=\rho$ and $b_{\parent{[k]}}^\prime=\rho^\prime$.
  \end{itemize}
  For any parent edge, i.e.\ any edge of the form $\parent{E}\in\ft_k$ for some $E\in\ft_k$, there are $i\in E$ and $j\in\parent{E}\setminus E$ for which $d(\rho,b_{\parent{E}})=\sD(\pi_{i,j}R_k)=\sD(\pi_{j,i}R_k)$. Any non-parent edge $E\in\ft_k$ is of the form $E=\{i,j\}\in\ft_k$, and similarly $d(\rho,b_E)=\sD(\pi_{i,j}R_k)=\sD(\pi_{j,i}R_k)$. Hence, if 
  $$\max_{1\le i<j\le k}\min\Big\{d_\cI(\pi_{i,j}R_k,\pi_{i,j}R_k^\prime),d_\cI(\pi_{j,i}R_k,\pi_{j,i}R_k^\prime)\Big\}<\varepsilon,$$
  then we have 
  \begin{equation}\label{eq:spinedist}
    \max_{E\in\ft_k}\left|d(\rho,b_E)-d^\prime(\rho^\prime,b_E^\prime)\right|<\varepsilon.
  \end{equation}
  This also constrains other distances. Specifically, for $A,B\in\ft_k$, consider the ``most recent common ancestor'' $C=\bigcap_{D\in\ft_k\colon A,B\subseteq D}D$ of $A$ and 
  $B$ in $\ft_k$. As the shortest path $b_A$ to $b_B$ in $S_k^\circ(R_k)$ passes through $b_C$, we have $d(b_A,b_B)=d(b_A,b_C)+d(b_C,b_B)$ in $(S_k^\circ(R_k),d)$, and likewise in 
  $(S_k^\circ( R_k^\prime),d^\prime)$. The triangular inequality (in $\mathbb{R}$) yields
  $$|d(b_A,b_B)-d^\prime(b_A^\prime,b_B^\prime)|\le|d(b_A,b_C)-d^\prime(b_A^\prime,b_C^\prime)|+|d(b_C,b_B)-d^\prime(b_C^\prime,b_B^\prime)|.$$
  Similarly, as $C$ is an ancestor of $A$ in $\ft_k$, we have $d(b_A,b_C)=d(\rho,b_A)-d(\rho,b_C)$ in $(S_k^\circ(R_k),d)$ and likewise in $(S_k^\circ( R_k^\prime),d^\prime)$, and this
  yields 
  $$|d(b_A,b_C)-d^\prime(b_A^\prime,b_C^\prime)|\le|d(\rho,b_A)-d^\prime(\rho^\prime,b_A^\prime)|+|d(\rho,b_C)-d^\prime(\rho^\prime,b_C^\prime)|,$$
  and likewise for $|d(b_B,b_C)-d^\prime(b_B^\prime,b_C^\prime)|$. Combining these inequalities with \eqref{eq:spinedist} yields
  \begin{equation}\label{eq:bpdist}\max_{A,B\in\ft_k}\left|d(b_A,b_B)-d^\prime(b_A^\prime,b_B^\prime)\right|<4\varepsilon.
  \end{equation}
  Apart from the special vertices, the ${\rm GH}$-correspondence $K$ includes pairs of points on branches, say $(v,v^\prime)$ of the form 
  $v=\lambda b_{A}+(1-\lambda)b_{\parent{A}}$ and $v^\prime=\lambda b_{A}^\prime+(1-\lambda)b_{\parent{A}}^\prime$. For this pair 
  $(v,v^\prime)$ and another pair $(w,w^\prime)$ obtained by replacing $A$ by $B$ and $\lambda$ by $\mu$, we can write $|d(v,w)-d(v^\prime,w^\prime)|$ as
  \begin{align*}
    \Big|&\lambda\mu\,\Big(d(b_A,b_B)-d^\prime(b_A^\prime,b_B^\prime)\Big)
           +\lambda(1-\mu)\Big(d(b_A,b_{\parent{B}})-d^\prime(b_A^\prime,b_{\parent{B}}^\prime)\Big)\\
           &+(1-\lambda)\mu\,\Big(d(b_{\parent{A}},b_B)-d^\prime(b_{\parent{A}}^\prime,b_B^\prime)\Big)
           +(1-\lambda)(1-\mu)\Big(d(b_{\parent{A}},b_{\parent{B}})-d^\prime(b_{\parent{A}}^\prime,b_{\parent{B}}^\prime)\Big)
    \Big|
  \end{align*}
  and, by (the triangular inequality in $\mathbb{R}$ and) \eqref{eq:bpdist}, this is also bounded above by $4\varepsilon$. We conclude that the ${\rm GH}$-distortion of $K$ is at most 
  $4\varepsilon$. By \eqref{eq:GH2},  the GH-distance is just half the infimum of ${\rm GH}$-distortions among all ${\rm GH}$-correspondences. The ${\rm GH}$-correspondence
  we have found therefore establishes this proposition.
\end{proof}

Our next aim is to derive a similar result for the Gromov--Hausdorff--Prohorov metric. Recall that we write $S_k(R_k)=(S_k^\circ(R_k),d_{\ell_1},0,\mu_k)$ for the rooted $\bR$-tree $S_k^\circ(R_k)$ further equipped with the mass measure $\mu_k$ that has an atom in each point of
$S_k^\circ( R_k)$ corresponding to a block in an edge partition, whose size is the block mass (or the sum of masses if there are two or more blocks associated with the same point of $S_k^\circ( R_k)$). The following example demonstrates that a ${\rm GHP}$-version of Proposition \ref{propGH} will
 need constants that increase at least linearly with $k$.

\begin{example} We begin by designing some edge partitions. Let $k\in 2\mathbb{N}$ and $\varepsilon>2\delta>0$. Consider an interval partition $\beta_0$ with diversity 
$\ell=\sD(\beta_0)=3k\varepsilon$ and such that any stretch $[a,b]\subset[0,\|\beta_0\|]$ of diversity $\sD_{\beta_0}(b)-\sD_{\beta_0}(a)\le k\varepsilon$ has total length $b-a<\delta$. In particular, all intervals of $\beta_0$ are of length less than $\delta<\varepsilon/2$. Let $\beta$ be $\beta_0$ with a interval of length $\varepsilon$ inserted at diversity $k\varepsilon$ from one end and let $\beta^\prime$ be $\beta_0$ with an interval of length $\varepsilon$ inserted at diversity $k\varepsilon$ from the other end. Finally, let $\gamma$ be an interval partition with two intervals of length $\ell$, top masses, followed by $\beta$, and construct 
$\gamma^\prime$ from $\beta^\prime$ in the same way. Then $d_\cI(\gamma,\gamma^\prime)=\varepsilon$, as 
\begin{itemize}
  \item this $d_\cI$-distortion is approached by $d_\cI$-correspondences that do not match the intervals of length $\varepsilon$, in the limit towards including all other intervals, 
    which are naturally matched as both $\gamma$ and $\gamma^\prime$ are built from $\beta_0$;
  \item the diversities of the two intervals of length $\varepsilon$ differ by $k\varepsilon>\varepsilon$ so the $d_\cI$-distortion of any $d_\cI$-correspondence that matches these two intervals 
    with each other will have $d_\cI$-distortion exceeding $\varepsilon$; matching an interval of length $\varepsilon$ with a top mass has a similar effect; 
  \item other intervals are of length at most $\delta<\varepsilon/2$, so if both intervals of length $\varepsilon$ are matched with those, the combined mass difference exceeds
    $2(\varepsilon-\delta)>\varepsilon$, and if only one of them is matched in this way, there is further unmatched mass of $\varepsilon$ in addition to the discrepancy 
    $\varepsilon-\delta$. 
\end{itemize}
Now let $\ft_k$ be a line of $k/2-1$ type-0 edges with $k/2$ type-2 edges to form a comb. Let $R_k$ have $\ell$ for all top masses, $\beta_0$ on all type-0 edges and $\beta$ on all type-2 edges, while $ R_k^\prime$ has $\ell$ for all top masses, $\beta_0$ on all type-0 edges and $\beta^\prime$ on all type-2 edges. Then $S_k(R_k)=(S_k^\circ(R_k),d,\rho,\mu)$ and $S_k(R_k^\prime)=(S_k^\circ(R_k),d^\prime,\rho^\prime,\mu^\prime)$ satisfy
$$d_{\rm GHP}(S_k( R_k),S_k (R_k^\prime))>k(\varepsilon-\delta)/2\qquad\mbox{and}\qquad
  \max_{1\le i<j\le k}d_{\cI}(\pi_{i,j} R_k,\pi_{i,j} R_k^\prime)=\varepsilon,$$
since
\begin{itemize}
  \item each $\pi_{i,j} R_k$ only contains a single interval of length $\varepsilon$, which leads to a distance $d_\cI(\pi_{i,j}R_k,\pi_{i,j}R_k^\prime)=\varepsilon$ as above for
    $d_\cI(\gamma,\gamma^\prime)$; 
  \item the $k/2$ locations $A_j$ and $A_j^\prime$, $1\le j\le k/2$, of atoms of size $\varepsilon$ in $S_k(R_k)$ and $S_k(R_k^\prime)$ are all, respectively, at distances 
    $(3i+1)k\varepsilon$ and $(3i+2)k\varepsilon$ from the respective roots, $1\le i\le k/2-1$; in particular, 
    $|d(A_j,\rho)-d^\prime(A_{j^\prime}^\prime,\rho^\prime)|\ge k\varepsilon$; 
  \item for any injective isometries $\phi\colon S_k^\circ(R_k)\rightarrow M$ and $\phi^\prime\colon S_k^\circ(R_k^\prime)\rightarrow M$ with 
$d_M(\phi(\rho_k),\phi^\prime(\rho_k^\prime))\le k(\varepsilon-\delta)/2$, the triangle inequality of $d_M$ yields
    \begin{align*}
    &d_M(\phi^\prime(A_{j^\prime}^\prime),\phi(A_j))\\
    &\ge|d_M(\phi^\prime(A_{j^\prime}^\prime),\phi^\prime(\rho^\prime))-d_M(\phi(A_j),\phi(\rho))|-d_M(\phi^\prime(\rho^\prime),\phi(\rho))\\
    &>k(\varepsilon-\delta)/2.
    \end{align*}
    Hence, none of the atoms $\phi^\prime(A_{j^\prime}^\prime)$, $1\le j^\prime\le k/2$, on $\phi^\prime(S_k^\circ(R_k^\prime))$ are in the set $C^{k(\varepsilon-\delta)/2}$ when $C=\{\phi(A_j),1\le j\le k/2\}$ is the set of atoms of size $\varepsilon$ on $\phi(S_k^\circ(R_k))$. The same is
    true for the branch points of $\phi^\prime(R_k^\prime)$. Hence, $C^{k(\varepsilon-\delta)/2}$ consists of, at most, $k/2$ stretches of length $k(\varepsilon-\delta)<k\varepsilon$, with mass in $\phi^\prime(S_k(R_k^\prime))$ of at most $\delta$. Therefore,
   \[
   \phi_*\mu(C)-\phi^\prime_*\mu^\prime(C^{k(\varepsilon-\delta)/2})>k\varepsilon/2-k\delta/2=k(\varepsilon-\delta)/2.
   \]
   Hence, $d_M^{\rm P}(\phi_*\mu,\phi^\prime_*\mu^\prime)>k(\varepsilon-\delta)/2$.
\end{itemize}
As $\delta\downarrow 0$, we find that the least upper bound on $d_{\rm GHP}(S_k( R_k),S_k(R_k^\prime))$ of the form
$$c_k\max_{1\le i<j\le k}d_{\cI}(\pi_{i,j} R_k,\pi_{i,j} R_k^\prime),$$
must have $c_k\ge k/2$, at least in any generality that includes these examples.
\end{example}

In this example, we have a block/atom of size $\varepsilon$ near the top or near the bottom of each of $k/2$ edges too far apart and too large to allow them to be matched for $d_\cI$ and to allow them to be embedded close together in any $(M,d_M)$ for $d_{\rm GHP}$. The key difference between the metrics is that only one such block affects each $d_\cI$, while $k/2$ affect $d_{\rm GHP}$. The tree shape plays very little role in this example. What matters is the number of type-2 edges, and a typical uniform tree shape has $k/4$ type-2 edges, see e.g.\ \cite[Proposition 59]{For-05}. The main
general observation is that $d_{\rm GHP}$ takes into account the spines under $\pi_{i,j}$, $1\le i<j\le k$, in a more additive
way than the maximum over $d_\cI$. We could replace the max by a sum in the claimed upper bound, but it is a sum over $\binom{k}{2}=O(k^2)$ terms, and the example only shows that we need a factor at least linear in $k$. A factor linear in $k$ is in fact sufficient. We do not chase the best constant. The following result will suffice for our purposes. 

\begin{proposition}\label{propGHP} Let $k\ge 2$. Consider two trees $R_k,R_k^\prime\in\widebar{\bT}_{k}^{\rm int}$ with the same shape $\ft_k$. Consider the weighted $\bR$-trees
  $S_k(R_k)$ and $S_k(R_k^\prime)$ and for each $1\le i<j\le k$ the four interval partitions $\pi_{i,j}R_k,\pi_{j,i}R_k,\pi_{i,j}R_k^\prime,\pi_{j,i}R_k^\prime\in\cI$.
  Then 
  $$d_{\rm GHP}(S_k(R_k),S_k(R_k^\prime))\le 3k\max_{1\le i<j\le k}\min\Big\{d_\cI(\pi_{i,j} R_k,\pi_{i,j} R_k^\prime),d_\cI(\pi_{j,i}R_k,\pi_{j,i}R_k^\prime)\Big\}.$$
\end{proposition}

The proof, while similar to the proof of Proposition \ref{propGH}, is somewhat lengthy and can be found in Appendix \ref{app:pfpropGHP}. The added complexity derives from the fact that the GH-correspondence of the proof of Proposition \ref{propGH} needs to be adjusted to facilitate the construction of a coupling of 
weight measures. In principle, this can be done using the $d_\cI$-correspondences between interval partitions to couple (large) atoms. To do this
consistently for all (large) atoms, we revert to an induction on tree shape in which we encounter a variety of cases.  
To apply Proposition \ref{propGHP} to non-resampling $k$-tree evolutions, we need to address shape changes. In the following, we will use the convention 
$\pi_{0,i}R_k=\pi_{i,0}R_k=\pi_{0,0}R_k=\{(0,\|R_k\|)\}$ for any $R_k\in\widebar{\mathbb{T}}_k^{\rm int}$ and $1\le i\le k$.

\begin{lemma}\label{lm:allpairs} Let $(\mathcal{T}_{k,-}^y,y\ge 0)$ be a non-resampling $k$-tree evolution starting from any initial $k$-tree 
  $T\in\mathbb{T}^{\rm int}_k$ with tree shape $\mathbf{t}$. We denote by $D_1,\ldots,D_k$ the degeneration times and by $A_m\subseteq[k]$ the label sets on 
  $[D_m,D_{m+1})$, $0\le m\le k$ where $D_0=0$ and $D_{k+1}=\infty$. Then there are a $\widebar{\mathbb{T}}_k^{\rm int}$-valued process  
  $(\widetilde{\mathcal{T}}^y_{k,-},y\ge 0)$ and maps $g_m\colon[k]^2\rightarrow (A_m\cup\{0\})^2$,  
  $0\le m\le k$, including the identity $g_0$, such that 
  \begin{itemize}
    \item $\tau(\widetilde{\mathcal{T}}^y_{k,-})=\tau(\mathcal{T}^y_{k,-})$, and $\widetilde{\mathcal{T}}^y_{k,-}$ has tree shape $\mathbf{t}$ for all $y\ge 0$,
    \item $\pi_{i,j}\widetilde{\mathcal{T}}^y_{k,-}\!=\pi_{g_m(i,j)}\mathcal{T}_{k,-}^y$ for all $1\le i\neq j\le k$, $D_m\le y<D_{m+1}$, $0\le m\le k$, 
    \item $(\pi_{i,j}\widetilde{\mathcal{T}}^y_{k,-},y\ge 0)$ is a $\cI^\circ$-valued type-2 evolution for all $1\le i\neq j\le k$.
  \end{itemize}
\end{lemma}
\begin{proof} To specify $(\widetilde{\mathcal{T}}^y_{k,-},y\ge 0)$, let us revisit the construction of $(\mathcal{T}^y_{k,-},y\ge 0)$, where the tree shape only changes at degeneration times, 
  through a combination of first (possibly) swapping two labels and then (always) reducing tree shape by dropping a label, see \eqref{eq:varrho}.  On the one hand, preserving tree shape calls for simpler 
  label dynamics by not swapping labels (which we did to obtain consistent $k$-tree evolutions). On the other hand, preserving tree shape means not to reduce tree shape (which we 
  did by removing edges with a zero top mass on an empty edge partition, as this naturally reassigns any other top mass of the removed parent edge to the 
  grandparent edge; this is relevant for the continuation of the $k$-tree evolution as it can increase the type of the grandparent edge from 1 to 2 or from 0 to 1). 

  In the light of this construction and as we are given $(\mathcal{T}^y_{k,-},y\ge 0)$, it is \vspace{-0.1cm} natural to construct $\widetilde{\mathcal{T}}^y_{k,-}$ from $\mathcal{T}^y_{k,-}$
  for each fixed $y$, by first inserting the empty edge partitions and zero top masses removed in the evolution up to time $y$, and then permuting labels to undo any label swapping and attain tree shape $\mathbf{t}$. Specifically,
  for $0\le y<D_1$, we have label set $A_0=[k]$ and $\mathcal{T}^y_{k,-}$ has tree shape $\mathbf{t}$. We \vspace{-0.1cm} set $\widetilde{p}_0={\rm id}\colon[k]\rightarrow[k]$ and 
  $\widetilde{\mathcal{T}}_{k,-}^y=\mathcal{T}_{k,-}^y$, $0\le y<D_1$. Inductively, assume we \vspace{-0.1cm} have constructed  
  $(\widetilde{\mathcal{T}}^y_{k,-},0\le y<D_m)$ and permutations $\widetilde{p}_0,\ldots,\widetilde{p}_{m-1}$ of $[k]$ such that\vspace{-0.1cm}
  \[
  \mathcal{T}_{k,-}^y=\pi_{A_n}\widetilde{p}_n\widetilde{\mathcal{T}}^y_{k,-},\quad D_n\le y<D_{n+1},\quad 0\le n\le m-1.\vspace{-0.1cm}
  \]
  By definition of non-resampling $k$-tree evolutions, $\mathcal{T}_{k,-}^{D_m}$ is obtained \vspace{-0.05cm} from $\mathcal{T}_{k,-}^{D_m-}$ by swap-reduction as in \eqref{eq:varrho}. 
  Specifically, this means first \vspace{-0.1cm} swapping $i_m:=I(\mathcal{T}_{k,-}^{D_m-})$ and $j_m:=J(\mathcal{T}_{k,-}^{D_m-})$ via the transposition $\sigma_{i_m,j_m}$ and then 
  removing $j_m$ so \vspace{-0.1cm} that $\mathcal{T}_{k,-}^{D_m}=\pi_{A_{m-1}\setminus\{j_m\}}\sigma_{i_m,j_m}\mathcal{T}_{k,-}^{D_m-}=\sigma_{i_m,j_m}\pi_{A_{m-1}\setminus\{i_m\}}\mathcal{T}_{k,-}^{D_m-}$ and \vspace{-0.1cm}
  \[
  \mathcal{T}_{k,-}^{D_m-}=\sigma_{i_m,j_m}\left(\mathcal{T}_{k,-}^{D_m}\oplus(\ell_m,j_m)\right)=\left(\sigma_{i_m,j_m}\mathcal{T}_{k,-}^{D_m}\right)\oplus(\sigma_{i_m,j_m}\ell_m,i_m),\vspace{-0.1cm}
  \]
  where we use a simplified form of the insertion operator notation of \eqref{eq:insert1}--\eqref{eq:insert2} so that $T\oplus(\ell,j):=T\oplus(\ell,j,U)$ with 
  $U=(1,0,\emptyset)\in\mathbb{T}_2^{\rm int}$ a degenerate 2-tree, and where $\ell_m=a_m$ if $i_m$ is part of a type-2 edge $\{i_m,a_m\}$ in the tree shape of $\mathcal{T}_{k,-}^{D_m-}$, or $\ell_m=(F,0,0)$ if the parent $F\cup\{i_m\}$ of $i_m$ is a 
  type-1 edge in the tree shape of $\mathcal{T}_{k,-}^{D_m-}$. We then define $\widetilde{p}_m=\sigma_{i_m,j_m}\circ\widetilde{p}_{m-1}$, 
  $\widetilde{p}_{n:m}=\widetilde{p}_m\circ\widetilde{p}_n^{-1}$, $1\le n\le m$, and\vspace{-0.1cm}
  \[
  \widetilde{\mathcal{T}}^y_{k,-}:=\widetilde{p}_m^{-1}\left(\cdots\left(\mathcal{T}_{k,-}^y\oplus(\widetilde{p}_{m:m}\ell_m,j_m)\right)\cdots\oplus(\widetilde{p}_{1:m}\ell_1,j_1)\right),\quad D_{m}\le y<D_{m+1}.\vspace{-0.1cm}
  \]
  Then the induction proceeds, since we have $A_m=A_{m-1}\setminus\{j_m\}$, so that\vspace{-0.1cm}
  \[
  \mathcal{T}^y_{k,-}=\pi_{A_m}\widetilde{p}_m\widetilde{\mathcal{T}}_{k,-}^y,\qquad D_m\le y<D_{m+1}.\vspace{-0.1cm}
  \]
  This means that $(\widetilde{\mathcal{T}}^y_{k,-},y\ge 0)$ satisfies the first bullet point.

  For the other two bullet points, we construct maps beginning with $g_0={\rm id}$. Now specifically, consider the pair $(i,j)=(1,2)$. Then  
  $\pi_{1,2}\widetilde{\cT}^y_{k,-}=\pi_{1,2}\cT^y_{k,-}$ for all $D_m\le y<D_{m+1}$ on the event $\{2\in A_m\}$, and we set $g_m(1,2)=(1,2)$. We have 
  $\pi_{1,2}\widetilde{\cT}^y_{k,-}=\{(0,\|\cT^y_{k,-}\|)\}=:\pi_{1,2}\cT^y_{k,-}$ for all $D_m\le y<D_{m+1}$, after label 2 has degenerated, i.e.\ on the event $\{2\not\in A_m\}$, and we set 
  $g_m(1,2)=(1,0)$. 
  For other label pairs $(i,j)\neq(1,2)$, we consider the setting of Lemma \ref{lmperm} for an initial permutation that satisfies $p_0(i)=1$ and $p_0(j)=2$. Changing 
  permutations as in the conclusion of the lemma, the $\cI^\circ$-projection $\pi_{1,2}$ yields a type-2 evolution starting from 
  $\pi_{i,j}\cT^0_{k,-}\in\cI^\circ$.
  We claim that this $\pi_{1,2}$-projection is equal to $(\pi_{i,j}\widetilde{\cT}_{k,-}^y,\,y\ge 0)$. Indeed, the two $\cI^\circ$-valued evolutions clearly coincide between
  shape changes of the underlying $k$-tree evolutions, and they are both continuous at shape changes, so they are equal. We accordingly define
  $g_m(i,j)=(p_m^{-1}(1),p_m^{-1}(2))$ on $\{2\in p_m A_m\}$ and $g_m(i,j)=(p_m^{-1}(1),0)$ on $\{2\not\in p_mA_m\}$.
\end{proof}

\section{Subtree decompositions of the Brownian CRT}\label{sec:subtreedecomp}
In this section we study the decomposition of a Brownian CRT along a reduced $k$-tree. To begin with, we sample a sequence of leaves from $\mu$. This can be done in $\mathbb{T}^{\rm real}$ using the kernel $\mathbf{m}_\infty$ of Proposition \ref{prop:sampleleaves}, but it will be easier to use a Brownian CRT 
$(\mathcal{T},d_{\ell_1},\mathbf{0})$ embedded in $\ell_1$, which in the setting of Section \ref{sec:BCRT} is naturally equipped with points 
$\Sigma_j=\mathbf{x}^{(j)}+D_j\mathbf{e}_j$, $j\ge 1$, and with the weak limit $\mu$ of the empirical measures $\mu_k=k^{-1}\sum_{j\in[k]}\delta(\Sigma_j)$, as $k\rightarrow\infty$. Recall also Aldous's observation that subtrees obtained when sampling from a Brownian excursion $2B^{\rm ex}$ have the same joint distribution, when suitably represented, as the trees in the line-breaking construction. In the current context, this entails that the  
${\rm GHP}^\infty$-isometry class of $(\mathcal{T},d_{\ell_1},\mathbf{0},\mu,(\Sigma_j,j\ge 1))$ is $\mathbf{m}_\infty(\mathrm{T},d(\sigma_j,j\ge 1))\mathbb{P}(\mathcal{T}\in d\mathrm{T})$.
 
We consider the subtrees $\cR_k^+=\bigcup_{j\in[k]}[\![\mathbf{0},\Sigma_j]\!]\subset\mathcal{T}$ spanned by the root $\mathbf{0}$ and leaves $\Sigma_1,\ldots,\Sigma_k$ and the 
projection of $\mu$ onto $\cR_k^+$, i.e.\ the image $\mu_k^+$ of $\mu$ under the natural projection onto the first $k$ coordinates in $\ell_1$. 
Subtrees like $\cR_k^+$ sampled from a CRT have played an important role ever since Aldous \cite{AldousCRT1,AldousCRT3} and others \cite{EPW,EW} initiated the study of CRTs. An important feature of $\cR_k^+$ is that it has
a discrete branching structure captured by a discrete tree shape $\ft_k$ in the space $\bT_{k}^{\rm shape}$ of binary combinatorial trees with $k$ leaves (and leaf edges) labeled by $[k]$, and with $k-1$ internal edges that are the elements of $\ft_k$ in the formalism of Section \ref{sec:killed_def}. Then $\cR_k^+$ can be decomposed into branches $(\cB_E,d_{\ell_1},\mu_k^+|_E)$, $E\in [k]\cup\ft_k$.
More recently, spinal and subtree decomposition 
theorems \cite{CW,DuW,HPW,PW13} have decomposed $(\cT,d_{\ell_1},\mathbf{0},\mu)$ along $\cR_k^+$. Vice versa, the joint distribution of the CRT $(\cT,d_{\ell_1},\mathbf{0},\mu)$ and $(\cR_k^+,d_{\ell_1},\mathbf{0},\mu_k^+)$, can be described, as follows. 

\begin{theorem}[Subtree decomposition I]\label{subdecI}\begin{enumerate}\item[\rm(a)] Let $(\cT,d_{\ell_1},\mathbf{0},\mu)$ be a BCRT embedded in $\ell_1$ and $\Sigma_j$, $j\ge 1$, a sequence
    sampled from $\mu$. Then $(\cR_k^+,d_{\ell_1},\mathbf{0},\mu_k^+)$ is determined by the following independent random variables:
    \begin{itemize}
      \item a tree shape $\ft_k$ that is uniformly distributed on $\bT_{k}^{\rm shape}$;
      \item a vector $(M_1,\ldots,M_{2k-1})\sim{\tt Dirichlet}(\frac12,\ldots,\frac12)$, representing the masses 
        $(\mu_k^+(\cB_E),E\in [k]\cup\ft_k)$, when listed in an order such as depth first search;
      \item interval partitions $\beta_1,\ldots,\beta_{2k-1}\!\sim\!\PDIP(\frac12,\frac12)$ representing the atoms of 
        $\mu_k^+|_{\cB_E}$, $E\!\in\! [k]\cup\ft_k$, in spinal order and normalized by $M_1,\ldots,M_{2k-1}$.\pagebreak
    \end{itemize}
  \item[\rm(b)] Conditionally given $(\cR_k^+,d_{\ell_1},\mathbf{0},\mu_k^+)$, the distribution of $(\cT,d_{\ell_1},\mathbf{0},\mu)$ is that of the tree
    obtained from $(\cR_k^+,d_{\ell_1},\mathbf{0},\mu_k^+)$ by grafting at each atom $x\in\cR_k^+$ of mass $m_x=\mu_k^+(x)>0$ a tree 
    $(\cT_x,\sqrt{m_x}d_x,\rho_x,m_x\mu_x)$ at its root $\rho_x$, where
    \begin{itemize}\item $(\cT_x,d_x,\rho_x,\mu_x)$ is an independent BCRT for each atom $x$ of $\mu_k^+$.
    \end{itemize}
  \end{enumerate}
\end{theorem}
\begin{proof} The case $k=1$ is essentially in \cite[Proposition 4(ii)]{HPW}, which shows that the collection of normalized
  spinal subtrees is a family of independent copies $\cT_i$, $i\ge 1$, of $\cT$, independent of the spinal mass partition. 
  The construction there of the Brownian CRT $\cT$ from a homogeneous partition-valued fragmentation process \cite{Bertoin06} and the 
  stopping line argument to decompose along the block containing 1 actually show that the independence of normalized spinal subtrees can be strengthened
  to include a joint independence from the spinal interval partition. See also \cite{RW}.
 
  Similarly, the stopping line argument applied to the stopping line at times when higher labels leave the blocks containing
  $1,\ldots,k$, shows that again the collection of normalized spinal subtrees is a family of independent copies $\cT_i^{(k)}$,
  $i\ge 1$, of $\cT$, independent of the $k$-tree $(\cR_k^+,d_{\ell_1},\mathbf{0},\mu_k^+)$, as required for (b). The $k$-tree can be represented
  as in \cite[Proposition 26]{PW13} in terms of the independent random variables listed in (a).  
\end{proof}


In the terminology of \cite{PitmWink09,PW13}, the weighted interval $M_0(\beta)$ associated with $\beta\sim{\tt PDIP}(\frac12,\frac12)$, is called a 
$(\frac12,\frac12)$-string of beads. Let us denote by $\phi\colon\cI\times (\bT^{\rm real})^\bN\rightarrow\bT^{\rm real}$ the map that assigns to 
$\beta\in\cI$ and $(\cT_i,i\ge 1)\in(\bT^{\rm real})^\bN$ the ($\mathrm{GHP}$-isometry class of the) tree that consists of a spine $[0,\sD(\beta)]$, to which the 
tree $\cT_i$ with mass rescaled by the size $P_i=b_i-a_i$ of the $i$th-largest block $U_i=(a_i,b_i)$ of $\beta$ and with distances rescaled by $\sqrt{P_i}$ is grafted at 
$\sD_\beta(U_i)$, provided that the resulting tree is a well-defined compact $\bR$-tree. To be definite, we assign the trivial tree $\Upsilon$ otherwise, and also, we break mass ties by spinal order. It was shown in \cite{RW} that this map is measurable and a.s.\ non-trivial when applied to the above random arguments. The case $k=1$ of Theorem \ref{subdecI} says that $\phi(\beta,(\cT_i,i\ge 1))$ is a Brownian CRT.

%
The $k=1$ case of Theorem \ref{subdecI} also yields that the branches $\cB_j$ of $\cR_k^+$ that have a leaf $\Sigma_j$ as an end point, with their subtrees grafted are rescaled Brownian CRTs $\cT_j$, $j\in[k]$. Constructions such as \eqref{eq:sk}, when applied to Brownian reduced $k$-trees, yield subtrees that we can represent as $(\cR_k,d_{\ell_1},\mathbf{0},\mu_k)$, $k\ge 1$, where $\cR_k=\cR_k^+\setminus\bigcup_{j\in[k]}\cB_j$ and $\mu_k$ is obtained by projecting $\mu$ or $\mu_k^+$ onto $\cR_k$. Then the decomposition of $(\cT,d_{\ell_1},\mathbf{0},\mu)$ along 
$(\cR_k,d_{\ell_1},\mathbf{0},\mu_k)$ is naturally obtained from the decomposition along $\cR_k^+$. To make a precise statement, we return to the notion of  a Brownian reduced $k$-tree of the form $\big(\ft_k,(X_j^{(k)},j\in[k]),(\beta_E^{(k)},E\in\ft_k)\big)\in\mathbb{T}_k^{\rm int}$ and deduce the following variant of Theorem \ref{subdecI}.

Recall from Definition \ref{def:SkandSkcirc} notation $S_k^\circ(R_k)$ and $W_k(\ell)$ for an $\mathbb{R}$-tree representation and atom locations 
derived from $R_k\in\widebar{\mathbb{T}}_k^{\rm int}$. 

\begin{corollary}[Subtree decomposition II]\label{subdecII}\begin{enumerate}\item[\rm(a)] Consider $(\cT,d_{\ell_1},\mathbf{0},\mu)$, a Brownian CRT embedded in $\ell_1$, and a sequence  $\Sigma_j$, $j\ge 1$, sampled from $\mu$. 
    Then the Brownian reduced $k$-tree $R_k=(\ft_k,(X_j^{(k)},j\in[k]),(\beta_E^{(k)},E\in T_k))$ is determined by the following independent random variables:
    \begin{itemize}\item a tree shape $\ft_k$ that is uniformly distributed on $\bT^{\rm shape}_{[k]}$;
         \item a vector $(M_1,\ldots,M_{2k-1})\sim{\tt Dirichlet}(\frac12,\ldots,\frac12)$ listing $k$ top masses $\big(X_j^{(k)},j\in[k]\big)$ followed by $k-1$ edge masses $\big(\|\beta_E^{(k)}\|,E\in\ft_k\big)$,
           the latter in an order such as depth first search; 
      \item interval partitions $\beta_{k+1},\ldots,\beta_{2k-1}\sim{\tt PDIP}(\frac12,\frac12)$
        representing the edge partitions $\big(\beta_E^{(k)},E\in\ft_k\big)$ normalized by $M_{k+1},\ldots,M_{2k-1}$.
    \end{itemize}
  \item[\rm(b)] Conditionally given $R_k$, the distribution of $(\cT,d_{\ell_1},\mathbf{0},\mu)$ is that of the tree
    obtained from $S_k^\circ(R_k)$ by grafting at $W_k(\ell)$ for each block $\ell\in\block(R_k)$ a tree $\big(\cT_\ell,\sqrt{\|\ell\|}d_\ell,\rho_\ell,\|\ell\|\mu_\ell\big)$ at its root $\rho_\ell$, where
    \begin{itemize}\item $(\cT_\ell,d_\ell,\rho_\ell,\mu_\ell)$ is an independent Brownian CRT for each $\ell\in\block(R_k)$.
    \end{itemize}
  \end{enumerate}
\end{corollary}
\begin{proof} (a) is a direct consequence of Theorem \ref{subdecI}. For (b), we apply the further argument there to stopping lines at heights
  when higher labels leave the blocks containing two or more labels $1,\ldots,k$.
\end{proof}

\section{${\rm GHP}$-path-continuity of the Aldous diffusion}\label{sec:continuity}

In this section we establish the existence of a continuous modification of the process $\widebar{\cT}(s)=S\big(\widebar{\cT}_{\!k,+}^s,k\ge 1\big)$, $s\ge 0$, of 
Definition \ref{df:adformal} and hence complete the proof of the existence of the Aldous diffusion as a ${\rm GHP}$-path-continuous Markov process. First, we study the $\bT^{\rm real}$-valued processes 
$\cT(y):=S(\cT^y_{k,+},k\ge 1)$, $y\ge 0$, constructed from a consistent system of \em self-similar \em $\widebar{\bT}^{\rm int}_k$-valued resampling $k$-tree evolutions $(\cT^y_{k,+},y\ge 0)$, $k\ge 1$, starting from a consistent family of
unit-mass Brownian reduced $k$-trees, i.e.\ in the case $j=1$ and $T=1\in\mathbb{T}_1^{\rm int}$ of Corollary \ref{cor:consistent_fam}(ii). Recall from Theorem \ref{thm:total_mass} that the total mass $B(y):=\|\cT^y_{k,+}\|$, $y\ge 0$, is fluctuating according to ${\tt BESQ}_1(-1)$. Since projections preserve total mass, 
this is indeed the same ${\tt BESQ}_1(-1)$ process for all $k\ge 1$ and for the limiting process $(\cT(y),\,y\ge 0)$. Furthermore, the marginal distributions of $\cT(y)$ are 
Brownian reduced $k$-trees with masses scaled by $B(y)$ and distances by $\sqrt{B(y)}$, by consistency and by the pseudo-stationarity of Proposition \ref{prop:pseudo:resamp}.   

Specifically, we will apply the Kolmogorov--Chentsov theorem to what we naturally call the \em self-similar Aldous diffusion \em $(\cT(y),y\ge 0)$, starting from a unit-mass Brownian CRT. We start by estimating the distance between $\cT(0)$ and the representation $\tau(\cT^0_{k,+})$ in $\bT^{\rm real}$ of $\cT^0_{k,+}\in\bT^{\rm int}_k$.

\begin{lemma}\label{prop1} For all $p>2$, there is a constant $K_{p,1}>0$ such that for all $k\ge 2$
  $$\bE\left[\left(d_{\rm GHP}(\cT(0),\tau(\cT^0_{k,+})\right)^{p \,}\right]\le K_{p ,1}k^{1-p /2}.$$
\end{lemma}
\begin{proof} In a unit-mass Brownian reduced $k$-tree $\cT^0_{k,+}$, mass is split into $2k-1$ \linebreak parts according to ${\tt Dirichlet}\big(\frac12,\ldots,\frac12\big)$, as we noted
  in Proposition \ref{prop:B_ktree} and Theorem \ref{subdecI}. If we aggregate the $k-1$  parts associated with edge partitions into a single mass and keep the $k$ 
  top masses separate, the distribution of these $k+1$ masses $(M_0,M_1,\ldots,M_k)$ is ${\tt Dirichlet}\big(\frac12(k-1),\frac12,\frac12,\ldots,\frac12\big)$, by 
  aggregation properties of Dirichlet vectors. The combined mass of the edge partitions is further split into block masses according to an independent ${\tt PD}\big(\frac12,\frac12(k-1)\big)$ distribution. 
We denote its parts in a size-biased random order by $P_j$, $j\ge 1$. From the subtree decomposition of Corollary \ref{subdecII}, each block mass is associated with a Brownian CRT. Denote these by $\cS_1,\ldots,\cS_k$ and $\cS^{(j)}$, $j\ge 1$. Denote the $p $-moment of the height ${\rm ht}(\cT)$ of a Brownian CRT $\cT$ by $h_p $ and note that we can bound the $(p /2)$-moment of a ${\tt Beta}\big(\frac12,k-1\big)$-variable above by $K_p  k^{-p /2}$, the $(p /2-1)$-moment of a ${\tt Beta}\big(\frac12,k/2\big)$ variable by $K_p  k^{1-p /2}$, for some $K_p >0$.
\begin{align*}&\bE\left[\left(d_{\rm GHP}(\cT(0),\tau(\cT_{k,+}^0))\right)^{p \,}\right]\\
  &\le\bE\left[\max\left\{\max_{i\in[k]}M_i^{p /2}({\rm ht}(\cS_i))^p ,\sup_{j\ge 1}M_0^{p /2} P_j^{p /2}({\rm ht}(\cS^{(j)}))^p \right\}\right]\\
  &\le\sum_{i=1}^k\bE\left[M_i^{p /2}\right]h_p +\sum_{j\ge 1}\bE\left[M_0^{p /2}\right]\bE\left[P_j^{p /2}\right]h_p \\
  &\le K_p  k^{1-p /2}h_p +\bE\left[P_1^{p /2-1}\right]h_p \le K_{p ,1}k^{1-p /2},
\end{align*}
as required. (See e.g.\ \cite{CSP} for details on the properties of Dirichlet and Poisson--Dirichlet distributions that we have used, and specifically \cite[Equation (2.23)]{CSP} for the penultimate step.)
\end{proof}

For the passage from time $0$ to time $y$ we will use a $k$-tree evolution $(\cT^{y}_{k,-},y\ge 0)$ without resampling (as resampling causes discontinuities in the pre-limiting processes) and note that for all $k\ge 2$, we can bound $d_{\rm GHP}(\mathcal{T}(0),\mathcal{T}(y))$ above by
\begin{equation}\label{eq:3GHP}
d_{\rm GHP}\big(\mathcal{T}(0),\tau(\cT^0_{k,-})\big)+d_{\rm GHP}\big(\tau(\cT^0_{k,-}),\tau(\cT^y_{k,-})\big)
                                                                     +d_{\rm GHP}\big(\tau(\cT^y_{k,-}),\mathcal{T}(y)\big).
\end{equation}
We will ultimately choose larger $k$ for smaller $y>0$, specifically, $k>y^{-\delta}$ for some $\delta\in(0,1)$. The first term in \eqref{eq:3GHP} can be bounded by Lemma \ref{prop1}, and we rephrase the lemma with our intended choice of $k$ in the following corollary. Indeed Propositions \ref{prop:partial_resamp:stat}--\ref{prop:partialfromresampling} will allow us to apply this corollary also for the third term in \eqref{eq:3GHP}, by separately taking into account the total mass at time $y$.

\begin{corollary}\label{cor:time0} Fix $\delta\in(0,1)$. Then for all $p>2$, $y\in(0,1]$ and $k>y^{-\delta}$
$$\bE\left[\left(d_{\rm GHP}(\cT(0),\tau(\cT_{k,-}^0)\right)^{p \,}\right]\le K_{p ,1}y^{p \delta/2-\delta}.$$
\end{corollary}
In order to obtain a good bound from the lemma applied to $\cT(y)$, we will need to make sure that the random number $J$ of remaining top masses in a non-resampling evolution is large with high probability, in the following sense.

\begin{lemma}\label{lm2} Fix $\varepsilon,\delta\in(0,1)$. Denote by $J_k^y$ the 
  random number of top masses of $\cT^{y}_{k,-}$, the non-resampling evolution at time $y$. Then for any $p>2$ there is $K_{p ,2}>0$ such that $\bP(J_k^y\le k(1\!-\!\varepsilon))\le K_{p ,2}y^{2(p \delta/2-\delta)}$ for all $y>0$ and $k=\lceil y^{-\delta}/(1-\varepsilon)\rceil$.
\end{lemma}
\begin{proof} Consider a ${\tt BESQ}(-1)$ process $(M(y),\,y\ge 0)$ starting from $M(0)=M\sim{\tt Gamma}(\frac{1}{2},1)$. By Lemma \ref{lem:BESQ:length}, the events of survival beyond times $x>0$ can be 
  given in terms of an independent ${\tt Gamma}(\frac32,1)$-variable $G$ as $\{M/2G>x\}$. Then 
  $$\bP(M/2G>2ky)\rightarrow 1\qquad\mbox{as $y\downarrow 0$, since $ky\rightarrow 0$.}$$
  In a non-resampling $k$-tree evolution starting from a Brownian reduced $k$-tree with $M_k\sim{\rm Gamma}(k-\frac12,1)$ initial mass, the $k$ top mass evolutions 
  (each stopped when hitting 0) are independent and distributed as $(M(y),\,y\ge 0)$. Denote by $J_{M_k}^{(k)}$ the number of top mass evolutions surviving to time $2ky$. 
  Then for all $y$ sufficiently small, $\bE\big[J_{M_k}^{(k)}\big]\ge k(1-\varepsilon/2)$ and hence by Hoeffding's 
  inequality \cite[Theorem 1]{Hoeffding1963}
  $$\bP\left(J_{M_k}^{(k)}\le k(1-\varepsilon)\right)\le \bP\left(J_{M_k}^{(k)}\le\bE\left[J_{M_k}^{(k)}\right]-k\varepsilon/2\right)\le\exp(-k\varepsilon^2/2).$$
  By the Weak Law of Large Numbers, $\bP\big(M_k/k<2\big)\ge\frac12$ for all $k$ sufficiently large, i.e.\ for all $y$ sufficiently small. By self-similarity, survival probabilities of
  ${\tt BESQ}(-1)$ are increasing functions of the initial mass, and $\mathbb{P}\big(J_m^{(k)}\le k(1-\varepsilon)\big)$ is a decreasing function of the initial mass $m$ of the 
  non-resampling $k$-tree evolution. In particular, for $k$ such that $\mathbb{P}\big(M_k/k<2\big)\ge\frac12$,
  \begin{align*}\frac{1}{2}\bP\left(J_{2k}^{(k)}\le k(1-\varepsilon)\right)&\le \int_{0}^{2k}\bP\left(J_m^{(k)}\le k(1-\varepsilon)\right)\mathbb{P}(M_k\in dm)\\
                                                                                                    &\le\bP\left(J_{M_k}^{(k)}\le k(1-\varepsilon)\right)\le\exp(-k\varepsilon^2/2).
  \end{align*}
  Again by self-similarity, scaling from initial mass $2k$ to unit initial mass yields
  $$\bP(J_k^y\le k(1-\varepsilon))=\bP\left(J_{2k}^{(k)}\le k(1-\varepsilon)\right)\le 2\exp(-k\varepsilon^2/2),$$
  and since $p >2$ is fixed and $k=\lceil y^{-\delta}/(1-\varepsilon)\rceil$, we can bound this probability by $K_{p,2}y^{2(p \delta/2-\delta)}$ for all $y>0$, by choosing $K_{p ,2}$ 
  sufficiently large. 
\end{proof}

With this estimate on the number of top masses in the non-resampling evolution established, we now derive the following corollary of Lemma \ref{prop1}.
\begin{corollary}\label{cor:timey} Fix $\varepsilon,\delta\in(0,1)$. For all $p >2$, there is a constant $K_{p ,3}>0$ such 
  that for all $y\in(0,1]$ and $k=\lceil y^{-\delta}/(1-\varepsilon)\rceil$, 
  \[
  \bE\left[\left(d_{\rm GHP}\left(\cT(y),\tau(\cT^{y}_{k,-})\right)\right)^{p \,}\right]\le K_{p ,3}y^{p \delta/2-\delta}.
  \]
\end{corollary}
\begin{proof}  Recall that $\cT(y)=S(\cT^y_{m,+},\,m\ge 1)$ a.s. By Proposition \ref{prop:partialfromresampling}, we can associate with $((\cT^y_{m,+},\,y\ge 0),\,m\ge 1)$, on the 
  same probability space, a consistent family $(\cT^y_{m+,k-},\,y\ge 0)$, $m\ge k\ge 1$, of partially resampling $m$-tree evolutions, in which labels in $[k]$ do not resample, while
  higher labels do. This includes a non-resampling $k$-tree evolution $\cT_{k,-}^y:=\cT_{k+,k-}^y$, $y\ge 0$. More precisely, this family can be chosen in such a way that for fixed $y$ and $k$, the trees $\cT^y_{m+,k-}$, $m\ge k$, are
  projectively consistent in the sense that $\pi_{-m}\cT_{m+,k-}^y=\cT_{(m-1)+,k-}^y$, $m\ge k+1$. By Proposition \ref{prop:partial_resamp:stat}, when conditioning on the label set $A_k^y$ of $\cT_{k,-}^y:=\cT_{k+,k-}^y$, on the event $\{\#A_k^y\ge 1\}$, these trees form a consistent family of independently scaled 
  Brownian reduced $(m-k+\#A_k^y)$-trees with label set $([m]\setminus[k])\cup A_k^y$. 

  Since furthermore, for the consistent family constructed in Proposition \ref{prop:partialfromresampling}, the tree $\cT_{m+,k-}^y$ is obtained as a projection from $\cT_{m,+}^y$ 
  up to relabelling, for each $m\ge k$, we can couple representatives of 
  $S(\cT^y_{m+,k-},m\ge k)$ and $S(\cT^y_{m,+},m\ge 1)=\cT(y)$ a.s., so that the former is a subset of the latter, equipped with the projected mass measure. But since 
  $S(\cT^y_{m+,k-},m\ge k)$ is itself a scaled BCRT, the 
  inclusion is an equality a.s. In the notation of Lemma \ref{lm2}, we have $J_k^y=\#A_k^y$. We split the 
  expectation according to the number $J_k^y$ of surviving top masses
  \begin{align*}\bE\left[\left(d_{\rm GHP}\left(\tau(\cT^{y}_{k,-}),\cT(y)\right)\!\right)^{\!p \,}\right]
				&=\bE\left[\left(d_{\rm GHP}\left(\tau(\cT^{y}_{k,-}),\cT(y)\right)\!\right)^{\!p }\cf\big\{J_k^y\le k(1-\varepsilon)\big\}\right]\\
					&\quad+\bE\left[\left(d_{\rm GHP}\left(\tau(\cT^{y}_{k,-}),\cT(y)\right)\!\right)^{p \,}\cf\big\{J_k^y>k(1-\varepsilon)\big\}\right]\!.
  \end{align*}
  Now denote by $B(y)$ the total mass of $\cT(y)$. By Theorem \ref{thm:total_mass} and since $(\cT^0_{j,+},j\ge 1)$ has unit mass, this sequence together with $\cT(0)=S(\cT^0_{j,+},j\ge 1)$
  is independent of $B(y)$. We remark that $B(y)$ and $J_k^y$ are not independent and that $(\cT^0_{j,+},j\ge 1)$ is not independent of $(B(y),J_k^y)$. We
  extend our probability space to support an independent copy $(\widehat{\cT}_{j},j\ge 1)$ of $(\cT^0_{j,+},j\ge 1)$, and we let
  $\widehat{\cT}=S(\widehat{\cT}_{j},j\ge 1)$. The following argument does not depend on the conditional distribution of $B(y)$ given $J_k^y$.

  On $\{J_k^y\le k(1-\varepsilon)\}$, we bound $d_{\rm GHP}(\tau(\cT_{k,-}^y),\cT(y))$ by ${\rm ht}(\cT(y))$.    
  By \cite[equation (18)]{Paper0}, we have $\bE\big[(B(y))^{p }\big]\le(1+2y(1+2(p -1)))^{p }\le (4p -1)^{p }$ for all $y\in(0,1]$, 
  and hence, using the Cauchy--Schwarz inequality, the first term can be bounded by Lemma \ref{lm2} as	
  \begin{align*}&\left(\bE\left[\left({\rm ht}\left(\cT(y)\right)\right)^{2p }\right)\right]^{1/2}\Big(\bP\left\{J_k^y\le k(1-\varepsilon)\right\}\Big)^{\!1/2}\\
	&\le(4p -1)^{p /2}\left(\bE\left[\left({\rm ht}(\cT(0))\right)^{2p \,}\right]\right)^{1/2}K_{p ,2}y^{p \delta/2-\delta},
  \end{align*}	
  On $\{J_k^y=j\}$ for $j>k(1-\varepsilon)$, we note that we can first replace $(\tau(\cT_{k,-}^y),\cT(y))$ by $(\tau(\widetilde{\cT}_{j}),\widehat{\cT})$,
  both scaled by (the independent) $B(y)$; then we increase the GHP-distance by reducing $j$ to $\lceil k(1-\varepsilon)\rceil$; then we drop the indicator; then
  we can replace $(\tau(\widetilde{\cT}_{j}),\widehat{\cT})$ by $(\tau(\cT^0_{\lceil k(1-\varepsilon)\rceil,+}),\cT(0))$ maintaining the scaling by 
  $B(y)$, which is an independent factor. Using Lemma \ref{prop1}, the second term is bounded above by
  \begin{align*}
  &\left(\bE\left[\left(B(y)\right)^p \right]\right)^{1/2}\left(\bE\left[\left(d_{\rm GHP}\left(\tau(\cT^0_{\lceil k(1-\varepsilon)\rceil,+}),\cT(0)\right)\right)^{2p }\right]\right)^{\!1/2}\\
    &\le(4p -1)^{p /2}\sqrt{K_{2p ,1}}(1-\varepsilon)^{1/2-p/2}\,y^{p \delta/2-\delta/2}.\qedhere
  \end{align*} 
\end{proof}

\medskip

Finally, we compare $\tau(\cT^0_{k,-})$ and $\tau(\cT^y_{k,-})$.

\begin{lemma}\label{lm:0toy} Fix $\varepsilon,\delta\in(0,1)$ and $\theta\in(0,\frac14)$. For each $p>2$, there is a constant $K_{p,4}$ such that for all $y\in(0,1]$ and 
  $k=\lceil y^{-\delta}/(1-\varepsilon)\rceil$, we have 
  \begin{equation}\label{eq:0toy}
     \bE\left[\left(d_{\rm GHP}\left(\tau(\cT_{k,-}^0),\tau(\cT_{k,-}^{y})\right)\right)^{p \,}\right]\le K_{p,4}y^{(\theta-\delta)p-2\delta}.
  \end{equation}
\end{lemma}
\begin{proof}
 In order to apply the bounds of Proposition \ref{propGHP}, we invoke Lemma \ref{lm:allpairs}, which allows us to replace $(\cT_{k,-}^y,y\ge 0)$ by a
  $\widebar{\mathbb{T}}^{\rm int}_k$-valued evolution $(\widetilde{\cT}_{k,-}^y,y\ge 0)$ whose tree shape remains constant and that has projections 
  $(\pi_{i,j}\widetilde{\cT}_{k,-}^y,y\ge 0)$ that are $\cI^\circ$-valued type-2 evolutions for all $1\le i\neq j\le k$, all starting from identically distributed initial states. 
  Furthermore, their initial distribution, which is obtained by concatenating the two top masses of a Brownian reduced 2-tree at the left end of the interval partition, was denoted by 
  $\widetilde{\mu}$ in Proposition \ref{type2holder}. 
  By Proposition \ref{propGHP}, this entails for all $y\in(0,1]$, $k\ge 2$,
  \begin{align*}\bE\left[\left(d_{\rm GHP}\left(\tau(\cT_{k,-}^0),\tau(\cT_{k,-}^{y})\right)\right)^{p \,}\right]
  &\le 3^p  k^p \bE\left[\max_{1\le i<j\le k}\left(d_{\cI}\left(\pi_{i,j}\cT_{k,-}^0,\pi_{i,j}\widetilde{\cT}^{y}_{k,-}\right)\right)^{p \,}\right]\\
  &\le 3^p  k^p \sum_{1\le i<j\le k}\bE\left[\left(d_{\cI}\left(\pi_{i,j}\cT_{k,-}^0,\pi_{i,j}\widetilde{\cT}^{y}_{k,-}\right)\right)^{p \,}\right]\\
  &\le 3^p  k^{p +2}\bE\left[\left(d_{\cI}\left(\pi_{1,2}\cT_{2,-}^0,\pi_{1,2}\cT^{y}_{2,-}\right)\right)^{p \,}\right],
  \end{align*}
We now take  $k=\lceil y^{-\delta}/(1-\varepsilon)\rceil$ and apply Proposition \ref{type2holder} to the $\cI^\circ$-valued type-2 evolution
$(\pi_{1,2}\cT_{2,-}^y,\,y\ge 0)$ to obtain \eqref{eq:0toy}.
Specifically, we established the existence of random H\"older constants $L=L_{\theta,y}$, with moments of all orders, for any $\theta\in(0,\frac14)$ and $y>0$, such that
  $$d_\cI(\gamma^a,\gamma^b)\le L|b-a|^\theta\qquad\mbox{for all }0\le a<b\le y,$$ 
  for an $\cI^\circ$-valued type-2 evolution  $(\gamma^y,y\ge 0)$ starting from $(0,A)\concat(0,B)\concat C\widebar{\beta}$ where 
  $(A,B,C)\sim{\tt Dirichlet}(\frac12,\frac12,\frac12)$ and $\widebar{\beta}\sim{\tt PDIP}(\frac12,\frac12)$ are independent.
  This entails that there is a constant $C_{\theta,p}$ such 
  that for all $y\in[0,1]$
  \[
  \bE\left[\left(d_{\cI}\left(\gamma^0,\gamma^y\right)\right)^{p\,}\right]\le C_{\theta,p}y^{\theta p}.\qedhere
  \]
\end{proof}

Now, we choose $\delta$ so as to get the best overall bound out of the three estimates of Corollaries \ref{cor:time0} and \ref{cor:timey} and Lemma \ref{lm:0toy}. Since the first two want $\delta$ large and the last one wants $\delta$ small, this easily gives $\delta=\frac{2}{3}\theta\in(0,\frac16)$ as the best choice for large 
$p $, and hence the proof of Theorem \ref{holder} will give GHP-H\"older continuity of index up to (but excluding) $\frac{1}{12}$. We do not claim the optimality of this index. We also establish GH-H\"older continuity up to (but excluding) $\frac{1}{4}$, and this will be optimal like for local times of $\Stable(\frac32)$-processes \cite{Boylan} where these bounds originate, see Section \ref{sec:Holder}.

\begin{theorem}\label{holder} The self-similar Aldous diffusion admits a continuous modification. This modification is a.s.\ 
  ${\rm GHP}$-H\"older continuous of index $\alpha$ for all $\alpha\in(0,\frac{1}{12})$ and  
  ${\rm GH}$-H\"older continuous of index $\alpha$ for all $\alpha\in(0,\frac14)$.
\end{theorem}
\begin{proof} Let $\theta\in(0,1/4)$, $\delta=\frac23\theta$ and $p>2$. By \eqref{eq:3GHP}, Corollaries \ref{cor:time0} and \ref{cor:timey} and Lemma \ref{lm:0toy}, there is $K_{p ,5}>0$ such that for all $y\in(0,1]$  
  \begin{align*}&\bE\left[\left(d_{\rm GHP}(\cT(0),\cT(y))\right)^p\right]\\
		&\le 3^p \bE\left[\left(d_{\rm GHP}\left(\cT(0),\tau(\cT^0_{k,-})\right)\right)^{p\,}\right]+
		      3^p \bE\left[\left(d_{\rm GHP}\left(\tau(\cT^0_{k,-}),\tau(\cT^{y}_{k,-})\right)\right)^{p\,}\right]\\ 
		 &\qquad+     3^p \bE\left[\left(d_{\rm GHP}\left(\tau(\cT^{y}_{k,-}),\cT(y)\right)\right)^{p\,}\right]\\
		&\le K_{p ,5}y^{\theta p /3-4\theta/3}.
  \end{align*}
  Denote by $\bP_x$ the law of the self-similar Aldous diffusion $(\cT(y),\,y\ge 0)$ starting from a Brownian CRT with initial 
  mass $x$, defined in the natural way as $S(\cT^y_{k,+},k\ge 1)$, $y\ge 0$, from a consistent system of pseudo-stationary resampling $k$-tree evolutions with initial mass $x$, as in Corollary \ref{cor:consistent_fam}(ii).  
  Denoting by $Z(y)=\|\cT(y)\|$, $y\ge 0$, the total mass evolution ({\tt BESQ}($-1$)), we have for all $0\le a<b\le 1$
\begin{align*}
  &\bE_1\left[(d_{\rm GHP}(\cT(a),\cT(b)))^p  \right]\\
    &\le \bE_1\left[(d_{\rm GHP}(\cT(a),\cT(b)))^p  \cf{\{Z(a)\ge b-a\}}\right]\\
    &\quad+2^{p  -1}\bE_1\left[({\rm ht}(\cT(a)))^p  \cf{\{Z(a)\le b-a\}}\right]+2^{p  -1}\bE_1\left[\|\cT(a)\|^p  \cf{\{Z(a)\le b-a\}}\right]\\
    &\quad+2^{p  -1}\bE_1\left[({\rm ht}(\cT(b)))^p  \cf{\{Z(a)\le b-a\}}\right]+2^{p  -1}\bE_1\left[\|\cT(b)\|^p  \cf{\{Z(a)\le b-a\}}\right],
\end{align*}
where we have 
\begin{itemize}
  \item split the term for $Z(a)\le b-a$ by ($(s+t)^p  \le 2^{p  -1}(s^p  +t^p  )$ and) the triangular inequality for $d_{\rm GHP}$ to compare 
  $\cT(a)$ and $\cT(b)$ with the degenerate one-point tree $\Upsilon$ with zero mass,
  \item and split the resulting two terms using $d_{\rm GHP}(\mathrm{T},\Upsilon)\le\max\{{\rm ht}(T),\mu(T)\}$. 
\end{itemize}
Now we consider each of these terms separately. We now apply the Markov property, pseudo-stationarity and self-similarity of $((\cT_{k,+}^y,k\ge 1),y\ge 0)$, as well as the effect of applying $S$ as noted around \eqref{eq:ch7:scaling} and \eqref{eq:ch7:scaling2}. Then the first term is bounded by
\begin{align*}&\int_{b-a}^\infty\max\{x^p  ,x^{p  /2}\}\bE_1\Big[\big(d_{\rm GHP}(\cT(0),\cT((b-a)/x))\big)^p  \Big]\bP_1(Z(a)\in dx)\\
&\le K_{p ,5}|b-a|^{\theta p /3-4\theta/3}\bE_1\left[\max\left\{(Z(a))^{p  /2-\theta p /3+4\theta/3},(Z(a))^{p  -\theta p /3+4\theta/3}\right\}\right],
\end{align*}
where the latter expectation is bounded uniformly in $a\in[0,1]$ for $p  $ sufficiently large. Next, clearly,
$$\bE_1\big[\|\cT(a)\|^p  \cf{\{Z(a)\le b-a\}}\big]\le|b-a|^p  ,$$
and by pseudo-stationarity and since distances scale by $\sqrt{x}$ when masses are scaled by $x$,
\begin{align*}\bE_1\Big[({\rm ht}(\cT(a)))^p  \cf{\{Z(a)\le b-a\}}\Big]
  &=\bE_1\left[(Z(a))^{p  /2}\cf{\{Z(a)\le b-a\}}\right]\bE_1\big[({\rm ht}(\cT(0)))^p  \big]\\
  &\le|b-a|^{p  /2}\bE_1\big[({\rm ht}(\cT(0)))^p  \big].
\end{align*}
To study the analogous time-$b$ quantities, it will be useful to first calculate moments of the squared Bessel process, using the transition 
density identified by G\"oing-Jaeschke and Yor \cite{GoinYor03} and well-known series representations of Bessel functions: 
\begin{align*}
  \bE_1\big[(Z(y))^q\big]
	&=\int_0^\infty\frac{1}{2y}x^{q-3/4}\exp\left(-\frac{x+1}{2y}\right)
		\sum_{m=0}^\infty\frac{1}{m!\Gamma(m+5/2)}\left(\frac{\sqrt{x}}{2y}\right)^{2m+3/2}dx\\
	&=\sum_{m=0}^\infty\frac{1}{m!\Gamma(m+5/2)}\left(\frac{1}{2y}\right)^{2m+5/2}\exp\left(-\frac{1}{2y}\right)
	    \int_0^\infty x^{q+m}e^{-x/2y}dx\\
	&=\sum_{m=0}^\infty\frac{\Gamma(q+m+1)}{m!\Gamma(m+5/2)}\left(\frac{1}{2y}\right)^{m+3/2-q}\exp\left(-\frac{1}{2y}\right)\\
	&\le \sum_{m=0}^\infty\frac{\Gamma(q+m+1)}{m!\Gamma(m+5/2)}\left(\frac{1}{2y}\right)^{m+3/2-q},
\end{align*}
where the series is finite for all $q>0$ and $y>0$, by the ratio test. Applying the Markov property of ${\tt BESQ}_1(-1)$ at time $a$ and then self-similarity, we calculate
\begin{align*}
  &\bE_1\big[\|\cT(b)\|^p  \cf{\{Z(a)\le b-a\}}\big]\\
    &=\bE_1\big[\cf{\{Z(a)\le b-a\}}\bE_{Z(a)}[(Z(b-a))^p  ]\big]\\
    &=\int_0^{b-a}x^p  \bE_1\left[(Z((b-a)/x))^p  \right]\bP_1(Z(a)\in dx)\\
    &\le\int_0^{b-a}x^p  \sum_{m=0}^\infty\frac{\Gamma(p  +m+1)}{m!\Gamma(m+5/2)}\left(\frac{x}{2(b-a)}\right)^{m+3/2-p  }\bP_1(Z(a)\in dx)\\
    &\le|b-a|^p  \sum_{m=0}^\infty\frac{\Gamma(p  +m+1)}{m!\Gamma(m+5/2)}\left(\frac{1}{2}\right)^{m+3/2-p  },
\end{align*}    
and similarly,
\begin{align*}
  &\bE_1\left[({\rm ht}(\cT(b)))^p  \cf{\{Z(a)\le b-a\}}\right]\\
    &=\bE_1\left[\cf{\{Z(a)\le b-a\}}\bE_{Z(a)}[{\rm ht}(\cT(b-a))^p  ]\right]\\
    &=\int_0^{b-a}x^{p  /2}\bE_1\left[{\rm ht}(\cT((b-a)/x))^p  \right]\bP_1(Z(a)\in dx)\\
    &\le\!\!\int_0^{b-a}\!\!\!x^{p  /2}\!\sum_{m=0}^\infty\!\frac{\Gamma(p  /2+m+1)}{m!\Gamma(m+5/2)}\!\left(\!\frac{x}{2(b-a)}\!\right)^{\!m+3/2-p  /2}\!\bP_1(Z(a)\!\in\! dx)\bE_1\!\left[({\rm ht}(\cT(0)))^p  \right]\\
    &\le|b-a|^{p  /2}\bE_1\left[({\rm ht}(\cT(0)))^p  \right]\sum_{m=0}^\infty\frac{\Gamma(p  /2+m+1)}{m!\Gamma(m+5/2)}\left(\frac{1}{2}\right)^{m+3/2-p  /2},
\end{align*}    
as required.

  Therefore, the Kolmogorov--Chentsov criterion \cite[Theorem I.(2.1)]{RevuzYor} applies. Specifically, if we write $\varepsilon(p )=\theta p /3-4\theta/3-1$, we get H\"older continuity for all indices in $(0,\varepsilon(p )/p )$. As $p \rightarrow\infty$, this includes all $\alpha\in(0,\frac{1}{12})$, as required, since $\theta\in(0,\frac14)$ was arbitrary.
  
  For GH, we can improve the bound on the distance between $\tau(\cT^0_{k,-})$ and $\tau(\cT^{y}_{k,-})$ using Proposition \ref{propGH} to obtain
  $$\bE\left[\left(d_{\rm GH}(\tau(\cT^0_{k,-}),\tau(\cT^{y}_{k,-}))\right)^p \right]\le 2^p  k^2C_{\theta,p }y^{\theta p }\le K_{p ,6}y^{\theta p -2\delta}.$$
  Since the coefficient of $p $ does not depend on $\delta$, we can choose $\delta=\frac12$ and obtain
  $$\bE\left[(d_{\rm GH}(\cT(0),\cT(y)))^p \right]\le K_{p ,7}y^{\theta p -1}.$$
  The same argument as for GHP, now with $\varepsilon(p )=\theta p -2$ yields H\"older continuity of all indices $\alpha\in(0,\frac14)$, since
  $\theta\in(0,\frac14)$ was arbitrary.
\end{proof}

\begin{corollary}\label{cor:ad:contversion}
  The Aldous diffusion has a GHP-path-continuous modification.
\end{corollary}		
\begin{proof} This follows from Theorem \ref{holder} because the de-Poissonization time-change is differentiable and because the total mass process
  of the self-similar process is (almost) $\frac12$-H\"older. 
\end{proof}

\section[Proof of Theorem 1.6 and resolution of Conjecture 1.1]{Proof of Theorem \ref{thm:intro:AD} and resolving the first part of Conjecture \ref{conj:Aldous}}\label{sec:properties}

In this section we pull the threads together and complete our proof of Claims 1 and 2 made at the beginning of this chapter. Specifically, we prove Theorem \ref{thm:intro:AD}, which claims that our construction yields a continuum-tree-valued process that (i) is stationary with the law ${\tt BCRT}$ (ii) is ${\rm GHP}$-path-continuous and (iii) has the simple Markov property. We further note that the process that we call the Aldous diffusion reduces to Wright--Fisher diffusions when decomposing around finitely many branch points, as Aldous \cite{ADProb2, AldousDiffusionProblem} stipulated, hence resolving the first part of Conjecture \ref{conj:Aldous}. 

\begin{proof}[Proof of Theorem \ref{thm:intro:AD}] Consider the process  $\widebar{\cT}(s)=S\big(\widebar{\cT}_{\!k,+}^s,k\!\ge\! 1\big)$, $s\!\ge\! 0$, constructed 
  by mapping the stationary $\widebar{\mathbb{T}}_\infty^{\rm int}$-valued consistent family of unit-mass resampling $k$-tree evolutions of
  Corollary \ref{cor:consistent_fam}(iii) by the map $S\colon\widebar{\mathbb{T}}_\infty^{\rm int}\rightarrow\mathbb{T}^{\rm real}$ of
  Definition \ref{def:S}, which was shown to be measurable in Proposition \ref{thm:S}. 

  (i) By Corollary \ref{cor:consistent_fam}(iii), the stationary distribution of the $\widebar{\mathbb{T}}_\infty^{\rm int}$-valued process is a consistent system
  of Brownian reduced $k$-trees. By Theorem \ref{thm:realtoint}, the image under $S$ is a Brownian CRT. Hence $(\widebar{\cT}(s),\,s\ge 0)$ is stationary with
  the law ${\tt BCRT}$.

  (ii) A GHP-path-continuous modification exists by Corollary \ref{cor:ad:contversion}.

  (iii) The simple Markov property was shown in Theorem \ref{thm:ad:markov}.
\end{proof}

This also establishes Claim 1 as stated at the beginning of this chapter. Specifically, this establishes the Aldous diffusion as a path-continuous Markov process in the Gromov--Hausdorff--Prokhorov space of weighted $\mathbb{R}$-trees that is stationary with the law of the BCRT. As a consequence of the construction and of Corollary \ref{cor:WF}, we can indeed conclude that this process resolves the first part of Conjecture \ref{conj:Aldous}.

\begin{corollary}\label{cor:conjproved} 
  Consider the stationary $\widebar{\mathbb{T}}_\infty^{\rm int}$-valued consistent family of unit-mass resampling $k$-tree evolutions 
  \[
  \widebar{\cT}_{\!k,+}^{\,s}=\Big(\widebar{\ft}_k^s,\big(\widebar{X}_j^{(k)}(s),\,j\in[k]\big),\big(\widebar{\beta}_E^{(k)}(s),\,E\in\widebar{\ft}_k^s\big)\Big),\quad s\ge 0,\, k\ge 1,
  \]
  of Corollary \ref{cor:consistent_fam}(iii) and weighted $\mathbb{R}$-trees $S_k\big(\widebar{\cT}_{\!k,+}^{\,s}\big)$, $k\ge 1$, and 
  $S\big(\widebar{\cT}_{\!k,+}^{\,s},\,k\ge 1\big)$ associated as in \eqref{eq:sk} and Definition \ref{def:S}, $s\ge 0$. Then the Wright--Fisher diffusions
  \[
  \Big(\big(\widebar X^{(k)}_j(u/4),\,j\in[k]\big),\big(\|\widebar\beta^{(k)}_E(u/4)\|,\,E\in\widebar{\ft}_k^{u/4})\big)\Big),\quad 0\leq u/4 <\tau_k,
  \]
  of Corollary \ref{cor:WF}, where $\tau_k:=\inf\big\{s\ge 0\colon\min_{j\in[k]}\widebar{X}^{(k)}_j(s)=0\big\}$, are embedded in the Aldous diffusion 
  $\widebar{\cT}(s)=S\big(\widebar{\cT}_{k,+}^{\,s},k\!\ge\! 1\big)$, $s\!\ge\! 0$, of Definition \ref{df:adformal} 
  as masses of connected components separated in any representative of
  $\widebar{\cT}(s)$ by the branch points corresponding to $W_k^{s}(E)+\sD\big(\beta^{(k)}_E\big)e_E\in S_k^\circ\big(\widebar{\cT}_{\!k,+}^{\,s}\big)$, 
  $E\in\widebar{\ft}_k^s$, defined for each $s\ge 0$ as specified before \eqref{eq:sk}.
\end{corollary}

This also establishes Claim 2 since $\widebar{\cT}_{\!k,+}^{\,0}$ is a Brownian reduced $k$-tree associated with the BCRT $\widebar{\cT}(0)$ and hence jointly distributed as the reduced $k$-tree obtained by sampling $k$ leaves from the mass measure of $\widebar{\cT}(0)$.

\section{General Markovianity and Continuity theorems}\label{sec:gen}

In this chapter, we have so far constructed the (unit-mass and self-similar) Aldous diffusion(s) and derived their Markovianity and GHP-path-continuity from  consistent systems of $k$-tree evolutions and their properties (Definition \ref{df:adformal} and Theorems \ref{thm:ad:markov} and \ref{holder} and their consequences). The philosophy has been that in the same way as continuum (random) trees are entirely characterized by properties of their (consistent system of) reduced $k$-trees, continuum-tree-valued evolutions can be obtained from consistent systems of evolutions of reduced $k$-trees. Furthermore, we have demonstrated that as the continuum structure (diffuse mass measure carried by a dense set of leaves) of a continuum tree is a feature not present in the reduced $k$-trees, the GHP-path-continuity does not have to be present in evolutions of reduced $k$-trees. 

The aim of this section is to generalize the setting away from BCRTs, the Aldous chain and the specific consistent system of $k$-tree evolutions towards more 
general CRTs and tree-valued Markov chains and to explore what properties of associated $k$-tree evolutions we require to establish more general continuum-tree-valued diffusions. We do this by revisiting the developments of this chapter and formulate a set of assumptions on the CRT and on the evolutions of reduced $k$-trees, under which we can prove Markovianity and continuity theorems for associated continuum-tree-valued processes.

\subsection*{Interval partitions with diversity and reduced $k$-trees of CRTs}

Recall 
from Section \ref{sec:IP} the notion of an interval partition with diversity and note that the set $\cI$ of such interval partitions, the
total diversity function $\sD$, the block diversity functions $\sD_\beta$, $\beta\in\cI$, the notions of total mass $\|\beta\|$ of $\beta\in\cI$ and mass
$\|U\|$ of $U\in\beta$, $\beta\in\cI$, and the notion of concatenation $\concat$ satisfy the following.  
\begin{enumerate}
  \item[A.] $(\cI,d_\cI)$ is a Lusin space and $\sD,\|\cdot\|\colon\cI\rightarrow[0,\infty)$ are continuous functions. Each $\beta\in\cI$ is a
    countable set equipped with functions $\sD_\beta\colon\beta\rightarrow[0,\sD(\beta)]$ and $\|\cdot\|\colon\beta\rightarrow[0,\|\beta\|]$ so that 
    $\sum_{U\in\beta}\|U\|=\|\beta\|$ for all $\beta\in\cI$. There is an associative operation of concatenation $\concat\colon\cI^2\rightarrow\cI$ such 
    that $\sD(\beta\concat\gamma)=\sD(\beta)+\sD(\gamma)$ and $\|\beta\concat\gamma\|=\|\beta\|+\|\gamma\|$.
\end{enumerate} 
For any space $(\cI,d_\cI)$ with these properties, we can proceed as in Section \ref{sec:killed_def}, let $\widebar{\mathbb{T}}^{\rm int}_{k}=\bigcup_{\ft\in\mathbb{T}^{\rm shape}_k}\{\ft\}\times[0,\infty)^k\times\cI^\ft$ be the space of $k$-trees, and note that projections $\pi_{-k}$ defined there, still give rise to a notion of consistent family $(R_k,k\ge 1)\in\prod_{k\ge 1}\widebar{\mathbb{T}}_{k}^{\rm int}$. In particular, this includes the spaces
$(\cI_\alpha,d_\alpha)$ of \cite{Paper1-0}, in which the ($\frac12$-)diversity of Definition \ref{def:diversity} is replaced by the notion of $\alpha$-diversity 
obtained by replacing $\sqrt{\pi}\sqrt{h}$ by $\Gamma(1-\alpha)h^\alpha$ in \eqref{eq:diversity}.
\begin{remark}\label{rem:mult} Indeed, 
we can include multifurcating trees, at some cost. Let us abuse notation and redefine $\mathbb{T}^{\rm shape}_k$ as the set of subsets $\ft$ of the power set of $[k]$ such 
that $[k]\in\ft$, $\#B\ge 2$ for all $B\in\ft$, and for all $B,C\in\ft$, we have $B\cap C=\emptyset$ or $B\subseteq C$ or $C\subseteq B$. Then we write $C=\parent{B}$ if $B\subset C$ and there is no $D\in\ft$ with $B\subset D\subset C$, for any $C\in\ft$ and 
$B\in\ft\cup\{\{j\},j\in[k]\}$. We say $C$ is multifurcating if there are three or more $B\in\ft\cup\{\{j\},j\in[k]\}$ with $\parent{B}=C$. Of course, we still set 
$\widebar{\mathbb{T}}^{\rm int}_{k}=\bigcup_{\ft\in\mathbb{T}^{\rm shape}_k}\{\ft\}\times[0,\infty)^k\times\cI^\ft$ now including multifurcating tree shapes. 
In the multifurcating setting, the definition of $\pi_{-k}$ in Definition \ref{def:proj} will have to be adjusted by including $\parent{\{j\}}$ in the domain of $\phi$, if this is a multifurcating branch point in $\ft$ so that this branch point is not removed. The form of this projection function will depend on the structure of $\cI$, which may not be a set of interval partitions.
\end{remark}
Definitions \ref{def:SkandSkcirc} and \ref{def:S} of $M_0$, $S_k$ and $S$ apply verbatim, but we do require some regularity of the map $M_0$ that is not implied by Assumption A, but is known to hold, for instance, when $(\cI,d_\cI)=(\cI_\alpha,d_\alpha)$, by \cite[Theorem 2.5(a)--(b)]{Paper1-0}.
\begin{enumerate}
\item[B.] $M_0\colon\cI\rightarrow\mathcal{M}$ defined as in Definition \ref{def:SkandSkcirc} is continuous, where the space $\mathcal{M}$ of \eqref{eq:onebranch} is equipped with the 
  Hausdorff--Prokhorov metric. 
\end{enumerate}
This is the key ingredient that makes the proof of Proposition \ref{thm:S} apply and establish the following generalization for the function $S$, now defined using any
$(\cI,d_\cI,\sD,\|\cdot\|,\concat)$ that satisfies Assumptions A--B.
\begin{proposition}\label{thm:S2} Under Assumptions A--B, the map $S\colon\widebar{\bT}_{\infty}^{\rm int}\rightarrow\mathbb{T}^{\rm real}$ is Borel measurable.
\end{proposition}
This proposition allows us to associate a continuum tree with a family of reduced $k$-trees, via the map $S$. In Theorem \ref{thm:realtoint}, we established a
map $R$ that allows us to associate a family of reduced $k$-trees with a continuum tree (with an infinite sample from its mass measure) in a generality sufficient for the BCRT. In fact, this map is also sufficient
for larger classes of CRTs.
We make 
the following assumption. 
\begin{enumerate}\item[C.] $(\cT,d,\rho,\mu)$ is a self-similar CRT in the sense of \cite{HM04}, binary with diffuse weight measure whose support is the set of leaves 
  and dense in $\cT$. Furthermore, its reduced 2-tree $(\ft,X_1,X_2,\beta):=(\ft_2,X_1^{(2)},X_2^{(2)},\beta_{[2]}^{(2)})$, constructed in the same way as 
  explained above Proposition \ref{prop:B_ktree}, is such that almost surely, $\cR_2$ is isometric to $[0,\sD(\beta)]$, every block $U\in\beta$ 
  corresponds to a connected component of $\cT\setminus\cR_2$ at distance $\sD(\beta)-\sD_{\beta}(U)$ from the root $\rho$, while blocks $X_1$ and $X_2$ 
  correspond to connected components at distance $\sD(\beta)$ from $\rho$. 
\end{enumerate}   

\begin{theorem}\label{thm:realtoint2} The measurable map 
  $R\colon\mathbb{T}_\infty^{\rm real}\rightarrow\widebar{\mathbb{T}}_\infty^{\rm int}$ of Theorem \ref{thm:realtoint} is such that, whenever Assumptions A--C are satisfied, we have
  \[
  S(R(\mathrm{T},\boldsymbol{\sigma}))=\mathrm{T}
  \]
  for $\mathbf{m}_\infty(\mathrm{T},d\boldsymbol{\sigma})\mathbb{P}(\cT\!\in\! d\mathrm{T})$-a.e.\ $(\mathrm{T},\boldsymbol{\sigma})\!=\!(\mathrm{T},(\sigma_j,j\!\ge\! 1))$.
\end{theorem}

\subsection*{The Markovianity theorem}

Recall that our aim is  to define continuum-tree-valued stochastic processes using the map $S$. The missing ingredient for this is a consistent system of 
$k$-tree evolutions. We therefore make the following assumption, modeled closely on the developments of
Chapters \ref{ch:constr}--\ref{ch:consistency}, requiring a generalization of the consistent (pseudo-)stationary families and other unit-mass and/or self-similar resampling/non-resampling/killed $k$-tree evolutions $(\widebar{\cT}_{k,\pm/\dagger}^s,s\!\ge\! 0)$ and/or $(\cT_{k,\pm/\dagger}^y,y\!\ge\! 0)$, $k\ge 1$, constructed there. In the following, we use unified notation $(\cT_{k,\pm/\dagger}^s,s\!\ge\! 0)$ for a process that may be unit-mass or 1-self-similar, depending on what the assumptions provide, but we no longer require to have both frameworks nor the independent evolution of type-$i$ compounds, $i=0,1,2$.
\begin{enumerate}
\item[D.] $((\cT_{k,+}^s,k\ge 1),s\ge 0)$  
  is a unit-mass (or 1-self-similar) $\widebar{\mathbb{T}}_\infty^{\rm int}$-valued Borel right Markov process that is (pseudo-)stationary with 
  marginal distributions obtained as consistent reduced $k$-trees of $(\cT,d,\rho,\mu)$ (scaled by independent random constants in the self-similar case). Processes $(\cT_{k,\pm/\dagger}^s,s\ge 0)$ for fixed $k\ge 1$ are Markovian and evolve continuously and according to the same dynamics until a top mass vanishes together with its parent edge partition. At 
  such times, killing or swap-reduction and/or resampling take place as specified in Section \ref{sec:resamp_def}, using the distribution $Q$ of a reduced 2-tree 
  associated with $(\cT,d,\rho,\mu)$. Killed $k$-tree evolutions are invariant under the permutation of labels. 
\end{enumerate}
This assumption is such that, under Assumptions A--D, Proposition \ref{couplelabels}, Lemma \ref{lmperm} and Corollary \ref{cor:couple_labels_ker} remain valid and their proofs are easily adapted, both as a pure unit-mass and as a pure 1-self-similar argument. We now deduce the first main theorem of this section, retracing the proof of Theorem \ref{thm:ad:markov}.

\begin{theorem}[Markovianity theorem]\label{thm:markov2} Under Assumptions A--D, the (unit-mass or 1-self-similar) process $\cT(s)=S\big(\cT_{k,+}^s,k\ge 1\big)$, $s\ge 0$, has the simple Markov property.
\end{theorem}

\subsection*{Representing non-resampling $k$-tree evolutions by $2$-tree evolutions}

Lemma \ref{lm:allpairs} naturally generalizes, as follows.

\begin{lemma}\label{lm:allpairs2} Let $(\mathcal{T}_{k,-}^s,s\ge 0)$ be a (unit-mass or 1-self-similar) non-resampling $k$-tree evolution under Assumptions A--D, starting from any initial $k$-tree 
  $T\in\mathbb{T}^{\rm int}_k$ with tree shape $\mathbf{t}$. We denote by $D_1,\ldots,D_k$ the degeneration times and by $A_m\subseteq[k]$ the label sets on 
  $[D_m,D_{m+1})$, $0\le m\le k$, where $D_0=0$ and $D_{k+1}=\infty$. Then there are a $\widebar{\mathbb{T}}_k^{\rm int}$-valued process  
  $(\widetilde{\mathcal{T}}^s_{k,-},s\ge 0)$ and maps $g_m\colon[k]^2\rightarrow (A_m\cup\{0\})^2$,  
  $0\le m\le k$, including the identity function $g_0$, such that 
  \begin{itemize}
    \item $\tau(\widetilde{\mathcal{T}}^s_{k,-})=\tau(\mathcal{T}^s_{k,-})$, and $\widetilde{\mathcal{T}}^s_{k,-}$ has tree shape $\mathbf{t}$ for all $s\ge 0$,
    \item $\pi_{i,j}\widetilde{\mathcal{T}}^s_{k,-}\!=\pi_{g_m(i,j)}\mathcal{T}_{k,-}^s$ for all $1\le i\neq j\le k$, $D_m\le s<D_{m+1}$, $0\le m\le k$, 
    \item $(\pi_{i,j}\widetilde{\mathcal{T}}^s_{k,-},s\ge 0)$ is an $\cI^\circ$-valued type-2 evolution for all $1\le i\neq j\le k$, i.e.\ a non-resampling self-similar $2$-tree evolution represented in $\cI^\circ$ as explained at the beginning of Section \ref{sec:dItoGHP}.
  \end{itemize}
\end{lemma}

\subsection*{Subtree decompositions of self-similar CRTs}

Part (b) of the subtree decomposition theorem, Corollary \ref{subdecII}, of the BCRT holds for all self-similar CRTs:

\begin{theorem}[Subtree decomposition -- general]\label{subdecgen} Consider any $\alpha$-self-similar CRT $(\cT,d_{\ell_1},\mathbf{0},\mu)$ embedded in $\ell_1$, a 
  sequence $\Sigma_j$, $j\ge 1$, sampled from $\mu$, and the associated reduced $k$-tree $R_k=(\ft_k,(X_j^{(k)},j\in[k]),(\beta_E^{(k)},E\in T_k))$. Then
  conditionally given $R_k$, the distribution of $(\cT,d_{\ell_1},\mathbf{0},\mu)$ is that of the tree
  obtained from $S_k^\circ(R_k)$ by grafting at $W_k(\ell)$ for each block $\ell\in\block(R_k)$ a tree 
  $\big(\cT_\ell,\|\ell\|^\alpha d_\ell,\rho_\ell,\|\ell\|\mu_\ell\big)$ at its root $\rho_\ell$, where $(\cT_\ell,d_\ell,\rho_\ell,\mu_\ell)$, $\ell\in\block(R_k)$, are i.i.d.\ copies of $(\cT,d_{\ell_1},\mathbf{0},\mu)$. 
\end{theorem}
\begin{proof} \cite[Proposition 4(ii)]{HPW} and the stopping-line arguments in the proofs of
  Theorem \ref{subdecI} and Corollary \ref{subdecII} apply to all self-similar CRTs.
\end{proof}

Instead of an exact distribution of $R_k$, $k\ge 1$, as in part (a) of Corollary \ref{subdecII} in the general case, we just assume the following.
\begin{enumerate}
  \item[E.] In the setting of Theorem \ref{subdecgen}, $R_k$ has block sizes $(M_{k,i},i\ge 1)$, $k\ge 2$, such that for all $q>1$ there is 
    $K_{q}>0$ such that $\sum_{i\ge 1}\mathbb{E}[M_{k,i}^q]\le K_{q}k^{1-q}$ for all $k\ge 2$.
\end{enumerate}
In fact, the check that the BCRT satisfies Assumption E is part of Lemma \ref{prop1} and works just the same with $\frac12$ suitably replaced by $\alpha\in(0,1)$. 
More generally, this is closely related to a property of Bertoin's conservative 1-self-similar fragmentation chains \cite[Corollary 3]{Bertoin03}, which have 
similar asymptotics for sums of powers of blocks. Under Assumption C, we expect that Assumption E is always satisfied.

\subsection*{The continuity theorem}

Finally, let us retrace our steps towards establishing the existence of a path-continuous modification of the Aldous diffusion, from Section \ref{sec:continuity}, in the present generality. 
Lemma \ref{prop1} generalizes, as follows.

\begin{lemma}\label{prop1_2} Under Assumptions A--E, for all $p>2$, there is a constant $K_{p,1}>0$ such that for all $k\ge 2$
  $$\bE\left[\left(d_{\rm GHP}(\cT(0),\tau(\cT^0_{k,+})\right)^{p \,}\right]\le K_{p ,1}k^{1-p /2}.$$
\end{lemma}
\begin{proof} From the subtree decomposition of Theorem \ref{subdecgen}, each block mass is associated with a CRT. Denote these by 
  $\cS^{(i)}$, $i\ge 1$, and the $p $-moment of the height ${\rm ht}(\cT)$ of the CRT $\cT$ of Assumption C by $h_p $. Now note that, by Assumption E,
\begin{align*}
\bE\Big[\big(d_{\rm GHP}(\cT(0),\tau(\cT_{k,+}^0))\big)^{p}\Big]
 &\le\bE\Big[\sup_{i\ge 1}M_{k,i}^{p /2}({\rm ht}(\cS^{(i)}))^p\Big]\\
 &\le\sum_{i\ge 1}\bE\Big[M_{k,i}^{p /2}\Big]h_p 
\le K_{p ,1}k^{1-p /2}.\qedhere
\end{align*}
\end{proof}
For the passage from time $0$ to time $s$ we will again use a $k$-tree evolution $(\cT^{s}_{k,-},s\ge 0)$ without resampling and bound 
$d_{\rm GHP}(\mathcal{T}(0),\mathcal{T}(s))$ above by
\begin{equation}\label{eq:3GHP2}
d_{\rm GHP}\big(\mathcal{T}(0),\tau(\cT^0_{k,-})\big)+d_{\rm GHP}\big(\tau(\cT^0_{k,-}),\tau(\cT^s_{k,-})\big)
                                                                     +d_{\rm GHP}\big(\tau(\cT^s_{k,-}),\mathcal{T}(s)\big).
\end{equation}
Choosing $k>s^{-\delta}$, the first term in \eqref{eq:3GHP2} can be bounded by Lemma \ref{prop1_2}, as in Corollary \ref{cor:time0}, which we restate here under Assumptions A--E.

\begin{corollary}\label{cor:time0_2} Fix $\delta\in(0,1)$. Then for all $p>2$, $s\in(0,1]$ and $k>s^{-\delta}$
$$\bE\left[\left(d_{\rm GHP}(\cT(0),\tau(\cT_{k,-}^0)\right)^{p \,}\right]\le K_{p ,1}s^{p \delta/2-\delta}.$$
\end{corollary}

In  
order to obtain a good bound from the lemma applied to $\cT(s)$, we will again need to make sure that the random number $J$ of remaining top masses in a non-resampling evolution is large with high probability, which we formulate here as an assumption, to play the role that Lemma \ref{lm2} plays in Section \ref{sec:continuity}.

\begin{enumerate}
  \item[F.] Fix $\varepsilon,\delta\in(0,1)$. Denote by $J_k^s$ the 
  random number of top masses of $\cT^{s}_{k,-}$, the non-resampling evolution at time $s$. Then for any $p>2$ there is $K_{p ,2}>0$ such that $\bP(J_k^s\le k(1\!-\!\varepsilon))\le K_{p ,2}s^{2(p \delta/2-\delta)}$ for all $s>0$ and $k=\lceil s^{-\delta}/(1-\varepsilon)\rceil$.
\end{enumerate}

In the setting of Section \ref{sec:continuity}, we use Propositions \ref{prop:partial_resamp:stat}--\ref{prop:partialfromresampling} to deduce from the consistent system of (pseudo-)stationary resampling $m$-tree
evolutions the existence of the following further processes on the same probability space. As the consistent system of Assumption D is less explicit, we also 
formulate the existence of the larger system of partially resampling evolutions as an assumption.
\begin{enumerate}
  \item[G.] Given a consistent system of (pseudo-)stationary resampling $m$-tree evolutions  $(\cT^s_m,s\ge 0)$, $m\ge 1$, as in Assumption D, there exists for each pair 
    $m\ge k\ge 1$, a process $((A_s^{m,k},B_s^k,\sigma_s^{m,k}),y\ge 0)$ that is constant between resampling times of $(\cT^s_m,\,s\ge 0)$, such that 
    $\sigma_s^{m,k}$ is a bijection between $A_s^{m,k}\subset[m]$ and $B_s^{k}\cup([m]\setminus[k])$ with $B_s^k\subseteq[k]$, and such that 
    $\cT_{m+,k-}^s:=\sigma_s^{m,k}\circ\pi_{A_s^{m,k}}(\cT^s_m)$, $s\ge 0$, is an $(m+,k-)$-partially resampling $m$-tree evolution in the sense of
    Definition \ref{def:partial_resamp}, i.e.\ in which labels in $[k]$ do not resample, while higher labels do. Furthermore, for any $k\ge 1$, these processes can be chosen to be projectively consistent in $m$, $m\ge k$.
    Finally,  for each $s\ge 0$, denote by $A^s_k$ the label set of $\mathcal{T}_{k+,k-}^s$, then conditionally given $A^s_k$, the tree 
    $\mathcal{T}_{m+,k-}^s$ is distributed like the (scaled) reduced $(m-k+\#A^s_k)$-tree with label set $([m]\setminus[k])\cup A^s_k$ constructed from
    the CRT $(\mathcal{T},d,\rho,\mu)$ of Assumption C. 
\end{enumerate}
In the self-similar case, we also need some control of the moments of the total mass evolution, which is trivial in the unit-mass case. In the self-similar case, we 
obtain this from \cite[equation (18)]{Paper0}, since all continuous positive 1-self-similar Markov processes are linear time-changes of ${\tt BESQ}(c)$ for some $c\in\mathbb{R}$:
\begin{lemma}\label{lm:gentot} Under assumptions A--D, the total mass evolution $(B(s),s\ge 0)$ is such that for all $p\ge 0$, we have  
  $\sup_{s\in[0,1]}\bE\big[(B(s))^{p }\big]<\infty$.
\end{lemma}

From the estimate on the number of top masses in the non-resampling evolution of Assumption F,  the existence of partially resampling evolutions as in  Assumption G and the bounds on the total mass evolution Lemma \ref{lm:gentot}, we now derive the following corollary of Lemma \ref{prop1_2}.

\begin{corollary}\label{cor:timey2} Under Assumptions A--G, fix $\varepsilon,\delta\in(0,1)$. For all $p >2$, there is a constant $K_{p ,3}>0$ such 
  that for all $s\in(0,1]$ and $k=\lceil s^{-\delta}/(1-\varepsilon)\rceil$, 
  $$\bE\left[\left(d_{\rm GHP}\left(\cT(s),\tau(\cT^{s}_{k,-})\right)\right)^{p \,}\right]\le K_{p ,3}s^{p \delta/2-\delta}.$$
\end{corollary}
Finally, 
we compare $\tau(\cT^0_{k,-})$ and $\tau(\cT^s_{k,-})$. To this end, we need a bound as in Proposition \ref{type2holder}.
\begin{enumerate}
  \item[H.] For $\cI^\circ$-valued non-resampling (unit-mass or 1-self-similar) 2-tree evolution $(\widetilde{\gamma}^s,s\ge 0)$ starting from an initial state, which is obtained by 
    concatenating the two top masses of a unit-mass reduced 2-tree at the left end of the interval partition, and for all $\theta\in(0,\frac{\alpha}{2})$, there is a random 
    constant $L=L_{\theta,s}$ with moments of all orders such that $d_{\cI}(\widetilde{\gamma}^a,\widetilde{\gamma}^b)\le L|b-a|^\theta$ for all 
    $0\le a<b\le s$. 
\end{enumerate}

\begin{lemma}\label{lm:0toy2} Under
Assumptions A--H, fix $\varepsilon,\delta\in(0,1)$ and $\theta\in(0,\frac{\alpha}{2})$. For each $p>2$, there is a constant $K_{p,4}$ such that for all $s\in(0,1]$ and 
  $k=\lceil s^{-\delta}/(1-\varepsilon)\rceil$, we have 
  \begin{equation*}
     \bE\left[\left(d_{\rm GHP}\left(\tau(\cT_{k,-}^0),\tau(\cT_{k,-}^{s})\right)\right)^{p \,}\right]\le K_{p,4}s^{(\theta-\delta)p-2\delta}.
  \end{equation*}
\end{lemma}
Again, we choose the $\delta$ that gives the best overall bounds from Corollaries \ref{cor:time0_2} and \ref{cor:timey2} and Lemma \ref{lm:0toy2}, and we conclude, as follows.
\begin{theorem}[Continuity theorem]\label{holder2} Under Assumptions A--H, the (unit-mass or 1-self-similar) continuum-tree-valued diffusion admits a continuous modification. This modification is a.s.\ 
  ${\rm GHP}$-H\"older continuous of index $\varrho$ for all $\varrho\in(0,\frac{\alpha}{6})$ and  
  ${\rm GH}$-H\"older continuous of index $\varrho$ for all $\varrho\in(0,\frac{\alpha}{2})$.
\end{theorem}

\subsection*{Examples} 
The  
reader may sense that the self-similar and
unit-mass Aldous diffusions could be the only processes for which all steps of the construction work. They certainly are the most significant such processes, which is why we have presented our results as results about these specific processes in the first instance. However, there are several more general settings, in which either substantial 
partial results are available or there are other reasons to believe such processes now become accessible and will substantially benefit from the generality of some of the methods presented here.
\begin{itemize}\item Arguably, closest to the Aldous chain is a chain built from Ford's alpha model \cite{For-05}, where up-steps are not insertions into a uniform random edge, 
  but according to weights $\alpha\in(0,1)$ for internal edges and $1-\alpha$ for external edges. This tree growth model yields binary Markov branching trees  
  that are weakly sampling consistent \cite{For-05} and have $\alpha$-self-similar CRTs as their scaling limits \cite{HMPW}. The meaning and relevance of weak sampling 
  consistency are that a down-step triggered by selecting a uniform leaf (and appropriate relabelling \cite{Soerensen}) yields a tree distributed like the tree with one
  fewer leaf. As a consequence, we can define down-up Markov chains with Ford trees as their stationary distributions, and the question of 
  a continuum-tree-valued diffusion arises in the same way as for the Aldous chain. 

  Consistent combinatorial $k$-tree chains have been constructed \cite{Soerensen}, Poissonization yields the same decoupling of the evolutions of type-0,
  type-1 and type-2 compounds. The theory of type-0 and type-1 evolutions has been developed \cite{Paper1-2} in a generality that covers cases relevant to the 
  chains associated with Ford's model in state spaces $(\cI_\alpha,d_\alpha)$ of \cite{Paper1-0} alluded to earlier. 

  The remaining challenge for Ford's model is that labels are not exchangeable. The above proofs of the Markovianity and continuity theorems both depend on exchangeable labels 
  obtained by sampling from the self-similar CRT. An alternative approach to the diffusion is to start from the strongly sampling consistent Markov branching 
  model associated 
  with the Ford CRT, as discussed more generally in the following bullet point -- in the special case of the Ford CRT, the above further structure may still be relevant. 
\item Consider any binary self-similar CRT, as in Assumption C. Sampling a sequence of independent leaves from its mass measure always yields tree shapes that form a strongly
    sampling consistent Markov branching model in the sense of \cite{HMPW}. This yields up-steps according to some regenerative tree growth
    rule \cite{PRW}. Strong sampling consistency means that the Aldous chain with its up-step replaced by the regenerative tree growth rule yields
    a stationary Markov chain. In general, there is no reason why any form of Poissonization should lead to the decoupling of type-$i$ compounds $i=0,1,2$, 
    but the regenerative structure appears to be useful in a unit-mass setting.
   
    The challenges here are captured by Assumptions A--H. Under what further conditions is there a notion of diversity that captures edge lengths from subtree 
    masses and a metric on interval partitions so that Assumptions A--B hold? What do $k$-tree evolutions look like? Can we set them up as autonomous Markov
    processes that satisfy the consistency requirements of Assumptions D and G? And can we then also establish Assumptions E, F and H? Note that in the absence
    of Poissonization, this is a programme that may be carried out entirely in the unit-mass framework.
\item Another setting, where the parallels reach further, is Marchal's growth process \cite{Marchal08} for stable CRTs. These are a one-parameter family of 
    multifurcating CRTs, so they do not satisfy Assumption C in the above wording, but if we relax the assumption of binary self-similar CRTs to general self-similar CRTs in the sense of \cite{HM04}, 
    they are included, indeed they form arguably the most important one-parameter family of self-similar CRTs \cite{DuquLeGall02,HM04,HMPW}. 
    As indicated in Remark \ref{rem:mult}, this can be handled by a suitable space $(\cI,d_\cI)$ generalizing the notion of an interval partition to record sizes of 
    subtrees that have the same distance from the root in blocks that are not necessarily totally ordered. Associated tree-valued up-down chains and induced
    consistent $k$-tree chains have been studied in \cite{Soerensen} in the combinatorial setting.

    The notion of $\alpha$-diversity in spaces $(\cI_\alpha,d_\alpha)$ again captures edge lengths from (coarse) spinal partitions and hence    
    provide again useful building blocks to capture the full coarse-fine spinal partition identified in \cite[Corollary 10]{HPW}. Some progress towards constructing
    relevant evolutions has been made in the interval partition evolutions with two-sided immigration and in the nested interval partition evolutions of 
    \cite{ShiWinkel-2}. Setting up a consistent system of $k$-tree evolutions is the subject of ongoing research, which aims to connect to the theory
    of this section to obtain (self-similar and unit-mass variants of) a continuum-tree valued \em stable Aldous diffusion\em.  
\end{itemize}


\chapter{Further properties of the Aldous diffusion}\label{chap8}

In this chapter, we study the Aldous diffusion constructed in Chapter \ref{ch:properties}. Specifically, we address the following points.
\begin{enumerate}\item[1.] We show that, although the Brownian CRT (and hence the Aldous diffusion at any fixed time) is almost surely binary, there is a dense Lebesgue-null set of times when the Aldous diffusion is a CRT with a non-binary branch point. Indeed, we argue, but stop short of a rigorous proof that at such times, there is precisely one degree-4 branch point and no branch points of any degree 5 or higher. The latter properties distinguish our process from Zambotti's excursion-valued process, which we discussed in Section \ref{sec:related_work}.  
  \item[2.] We continue our discussion of the failure of the strong Markov property of the Aldous diffusion by showing that there are stopping times in the natural filtration of the Aldous diffusion at which the 
   Aldous diffusion is non-binary, and that at such stopping times, the strong Markov property fails. 
  \item[3.] We embed the stationary Aldous chain into the Aldous diffusion in such a way that the steps of the Aldous chain occur at times that are spaced by independent exponentially distributed random 
    variables. We show that these embedded continuous-time Aldous chains, suitably rescaled, converge to the Aldous diffusion hence completing the resolution of Conjecture \ref{conj:Aldous}. We further deduce the reversibility of the Aldous diffusion. 
\end{enumerate}
These three points are covered, respectively, in Sections \ref{sec:nonbin}, \ref{sec:smp} and \ref{sec:emb}. We conclude this chapter by stating open problems related to the Aldous diffusion in Section \ref{sec:openprob}.

\section{Non-binary branch points}\label{sec:nonbin}

Recall that all branch points of the Brownian CRT are binary almost surely. More precisely, it is almost surely the case that for all points in a Brownian CRT, their removal disconnects the Brownian CRT into three connected components (around countably many branch points) or  two connected components (around other non-leaf vertices) or the tree remains connected (when removing a leaf). In our 
discussion of the failure of the strong Markov property of the Aldous diffusion in Remark \ref{prop:notstrongmarkov} we argued that there are nevertheless random times at which the Aldous diffusion has a ternary branch point (whose removal disconnects into four connected components). Before we return to this discussion in the next section, we consider non-binary branch points (whose removal disconnects into at least four connected components) more systematically. 

\begin{theorem} \label{thm:exceptionaldegrees}
  On an event of probability 1, the Aldous diffusion $\widebar{\cT}(s)$, $s\ge 0$, has a dense null set of times $s\ge 0$ at which $\widebar{\cT}(s)$ has a non-binary branch point.
\end{theorem}

For the proof of Theorem \ref{thm:exceptionaldegrees}, it will be important to have $k$-trees embedded in the Aldous diffusion at all times. While the Aldous diffusion, as defined in Definition \ref{df:adformal}, has been constructed from $k$-tree evolutions in a way that naturally entails an embedding at fixed times almost surely, the  Kolmorogov--Chentsov argument of Theorem \ref{holder} passes to a continuous modification that, a priori, only preserves the embedding property at a countable dense set of times (for instance all dyadic times in the proof of \cite{RevuzYor}). The following lemma allows us to strengthen this,
with the help of path properties of $k$-tree evolutions. Recall the definition of the Gromov--Hausdorff space $(\mathbb{T}_\circ^{\rm real},d_{\rm GH})$ of \eqref{eq:GH}.

\begin{proposition}\label{prop:embcont}  Let $(\widebar{\cT}_{k,+}^s,s\ge 0)$ be a consistent family of stationary unit-mass resampling $k$-tree evolutions, $k\ge 1$, as in Corollary 
  \ref{cor:consistent_fam}(iii) and $\big(\widebar{\cT}(s),s\ge 0\big)$ the ${\rm GHP}$-path-continuous modification of 
  $\big(S\big(\widebar{\cT}_{k,+}^s,k\ge 1\big), s\ge 0\big)$, where $S\colon\widebar{\mathbb{T}}_\infty^{\rm int}\rightarrow\mathbb{T}^{\rm real}$ is as defined in 
  Definition \ref{def:S}. Then it is almost surely the case that for all $s\ge 0$, the $k$-trees $\tau(\widebar{\cT}_{k,+}^s)$, $k\ge 1$, can be embedded isometrically into (any representative of) 
  $\widebar{\cT}(s)$ as a nested family of subsets. 
\end{proposition}

We prove this proposition in Appendix \ref{appx:markedGH}.

\begin{corollary}\label{cor:allemb} It is almost surely the case that
  \[
   S\big(\widebar{\cT}_{k,+}^s,k\ge 1\big)=\lim_{k\rightarrow\infty}\tau(\widebar{\cT}_{k,+}^s)\subseteq\widebar{\cT}(s),\qquad\mbox{for all }\ s\ge 0,
  \]
  in the sense that any representative of the left-hand side can be embedded isometrically into any representative of the right-hand side as rooted $\mathbb{R}$-trees.
\end{corollary}

While we believe that this subset property can be strengthened to an equality, we have not been able to prove this. We also emphasize that the embedding is as
a rooted $\mathbb{R}$-tree, not as a weighted rooted $\mathbb{R}$-tree.

\begin{proof}[Proof of Theorem \ref{thm:exceptionaldegrees}] Since the Aldous diffusion is stationary with the distribution of the binary BCRT, the set of times where any non-binary branch points exist is a Lebesgue null set. 
  Since de-Poissonization does not affect the relevant 
  statements about non-binary branch points, it suffices to establish the corresponding results for the self-similar Aldous diffusion (restricting the claim that the set of exceptional times is 
  dense to the lifetime of the self-similar process). 

  To show that the set of times with a non-binary branch point is dense almost surely, it suffices to show that there is a sequence of such times accumulating at time 0, since this then extends 
  from time 0 to all rational times, by pseudo-stationarity of the self-similar processes of Corollary \ref{cor:consistent_fam}(ii), during the lifetime of the self-similar Aldous diffusion $(\mathcal{T}(y),y\ge 0)$. Equivalently, we can consider $\eta_0=\inf\{y\ge 0\colon\mathcal{T}(y)\mbox{ non-binary}\}$ and show that $\eta_0=0$ almost surely.

  Consider a killed $3$-tree evolution $\big(\cT_{3,\dagger}^y,\,y\ge 0\big)$ starting according to the pseudo-stationary distribution with ${\tt Gamma}(\frac52,\lambda)$ initial
  mass. As in the proof of Proposition \ref{prop:pseudo:degen}, this process is composed of two 
  independent evolutions of types 1 and 2, respectively. By Constructions \ref{deftype2} and \ref{type1:construction} and Proposition 
  \ref{type1totalmass}, the three top masses and total masses of the two interval partitions form independent ${\tt BESQ}(-1)$ and ${\tt BESQ}(1)$ processes starting from 
  ${\tt Gamma}(\frac12,\lambda)$ initial mass, until the first of these vanishes. With positive probability $p_3>0$, say, we observe the event $A_3$ that this is the ${\tt BESQ}(1)$ 
  process associated with the type-2 
  evolution. At this time $\eta$, all three top masses are projected to $W_3^\eta(1)=W_3^\eta(2)=W_3^\eta(3)$ in $S_3(\cT_{3,\dagger}^\eta)$. This should correspond to a   
  ternary 
  branch point of the self-similar Aldous diffusion. 

  Rather than confirming this rigorously in the $3$-tree setup, we consider a similar event $A_6$ associated with a killed $6$-tree evolution $\big(\cT_{6,\dagger}^y,\,y\ge 0\big)$ starting 
  according to the pseudo-stationary distribution with ${\tt Gamma}(\frac{11}{2},\lambda)$ initial mass. Consider tree shapes that have three type-2 edges with label 
  sets $\{1,4\}$, $\{2,5\}$ and $\{3,6\}$ and two type-0 edges. The edge of interest is the type-0 edge that is between two branch points of the tree shape, not the one adjacent to the root. Consider
  the event $A_6$ that we observe such a tree shape and that of the five interval partitions and six top masses, the first to reach zero mass, or zero diversity in the case of an edge, is the edge 
  of interest. Then $A_6$ has positive probability $p_6>0$. At this random time $\eta^\prime$, the $\mathbb{R}$-tree $S_6(\cT_{6,\dagger}^{\eta^\prime})$ is a star 
  tree that has a ternary branch point with four branches respectively leading to the root and three leaves $W_6^{\eta^\prime}(1)=W_6^{\eta^\prime}(4)$,  $W_6^{\eta^\prime}(2)=W_6^{\eta^\prime}(5)$ and  $W_6^{\eta^\prime}(3)=W_6^{\eta^\prime}(6)$, and on the event $A_6$, these branches have positive lengths given by the diversities of the other four
  interval partitions.  
     
  By scaling, it is easy to see that this probability $p_6=\mathbb{P}(A_6)$ is unchanged if we change the initial total mass to 
  any other positive initial total mass, deterministic or random. Now fix $\varepsilon>0$. Starting from unit initial mass, the probability 
  $p_6^\prime$ that $A_6$ holds with $\eta^\prime\le\varepsilon$ is positive, and by scaling, the probability of this event is bounded below by 
  $p_6^\prime$ when the initial mass is bounded above by 1.

  Intuitively, events like $A_6$ happen at all scales. More formally, we will consider a consistent sequence $(\mathcal{T}^y_{k,+},y\ge 0)$, $k\ge 1$, of 
  resampling self-similar $k$-tree evolutions starting from a consistent family of unit-mass Brownian reduced $k$-trees, as in Corollary 
  \ref{cor:consistent_fam}(ii) for $T=1\in\mathbb{T}_1^{\rm int}$. For each $n\ge 1$ we will choose a
  random $K_n$ so that the initial tree can be projected to $n$ disjoint scaled Brownian reduced $6$-trees. Specifically, consider the initial $n$-tree and let
  $K_n$ be the first $k\ge n$ such that the tree shape of $\mathcal{T}_{k,+}^0$ contains edges with label sets $\widetilde{E}^{(n)}_1,\ldots\widetilde{E}^{(n)}_n$ such that
  $i\in \widetilde{E}^{(n)}_i$ and $\#\widetilde{E}^{(n)}_i\ge 6$ for all $i\in[n]$. For each $i\in[n]$ let $E^{(n)}_i\subseteq\widetilde{E}^{(n)}_i$ be the subset of six lowest labels. Then the
  subtrees of $\cT^0_{k,+}$ above $\widetilde{E}_i^{(n)}$, $i\in[n]$, are disjoint and, projected to $E^{(n)}_i$, $i\in[n]$, give rise to $n$ 6-trees
  $\cS_i^{(n)}$, $i\in[n]$.  

  Now we consider 6-tree evolutions starting from $\cS_i^{(n)}$, $i\in[n]$, $n\ge 1$, induced by $(\cT_{k,+}^y,\,y\ge 0)$, $k\ge 1$. Specifically, for each 
  $i\in[n]$, consider a permutation $p_i^{(n)}$ of $[K_n]$ that maps
  $E_i^{(n)}$ to $\{1,2,3,4,5,6\}$. By Lemma \ref{lmperm}, Proposition \ref{prop:resamp_to_non} and Corollary \ref{cor:consistent_fam}, there are an 
  induced resampling $K_n$-tree evolution starting from the initial tree with labels permuted by $p_i^{(n)}$, an associated non-resampling $K_n$-tree evolution, a non-resampling 
  $6$-tree evolution obtained by projection via $\pi_6$, and finally a non-resampling $E_i^{(n)}$-tree evolution obtained by relabeling via $p_n^{-1}|_{[6]}$, in which the subtree with root
  edge labeled $E_i^{(n)}$ performs a 6-tree evolution starting from $\cS_i^{(n)}$ and run until its first degeneration time. 
  The purpose of passing to non-resampling evolutions is to avoid resampling into $\cS_i^{(n)}$. The purpose of the relabeling is to prevent lower labels from 
  swapping with labels in $E_i^{(n)}$ in a swap-reduction step of the non-resampling $K_n$-tree evolution. 

  Now fix $n$ and still consider the non-resampling $K_{n}$-tree evolution (conditionally given $K_{n}$). Then the tree shape of the projection
  $\cN_{n}^0=\pi_{E_1^{(n)}\cup\cdots\cup E_j^{(n)}}\cT_{K_{n},+}^0$, includes, for each $i\in[n]$, a 6-tree $\cS_i^{(n)}$, and these 6-trees are disjoint subtrees of $\cN_{n}^0$. It follows from Proposition \ref{prop:B_ktree} and aggregation properties of Dirichlet vectors that conditionally given their masses 
  $\|\cS_i^{(n)}\|$, $i\in[n]$, the scaled Brownian reduced 6-trees $\cS_i^{(n)}$, $i\in[n]$, are independent. Similarly projecting the non-resampling $K_n$-tree evolution gives rise to 
  a non-resampling $E_1^{(n)}\cup\cdots\cup E_n^{(n)}$-tree evolution $(\cN_{n}^y,\,y\ge 0)$ in which the induced evolutions 
  $\big((\cS_i^{(n)})^y,\,y\ge 0\big)$ are non-resampling 6-tree evolutions starting from $(\cS_i^{(n)})^0=\cS_i^{(n)}$, $i\in[n]$, and they are conditionally independent given their initial masses 
  $\|\cS_i^{(n)}\|$, $i\in[n]$.  

  Now let $\eta_n=\inf\{y\ge 0\colon S_{K_n}(\mathcal{T}^y_{K_n,+})\mbox{ non-binary}\}$. Denote by $B_i^{(n)}$ the event that 
  $\big((\cS_i^{(n)})^y,\,y\ge 0\big)$ rescaled to unit mass satisfies $A_6$ and $\eta^\prime\le\varepsilon$. Then we find that
  $\mathbb{P}(\eta_n\le\varepsilon)\ge\mathbb{P}\big(\bigcup_{i\in[n]}B_i^{(n)}\big)=1-(1-p_6^\prime)^n\rightarrow 1$ as $n\rightarrow\infty$.  

  We now translate these statements into statements about the self-similar Aldous diffusion $\cT(y)=S\big(\cT_{k,+}^{\,y},k\!\ge\! 1\big)$, $y\ge 0$. By Corollary \ref{cor:allemb},  $(S_k(\cT_{k,+}^y),k\ge 1)$ is, almost surely,  
  a GHP-convergent tree growth process for all $y\ge 0$.  By Proposition \ref{prop:embcont}, any pre-limiting tree can be embedded in the limiting tree. In particular, any non-binary branch point in a pre-limiting tree is also a non-binary branch point in the limit. Hence, $\mathbb{P}(\eta_0\le\varepsilon)\ge\mathbb{P}(\eta_n\le\varepsilon)\rightarrow 1$ as $n\rightarrow\infty$. But this holds for all $\varepsilon>0$. Therefore, $\eta_0=0$ almost surely. This completes the proof.\end{proof}

\begin{remark}\label{rmk:higherdegree} 
  It is natural to ask about branch points of degree 5 or higher or several simultaneous ternary branch points. We believe that neither is possible, and indeed the reasoning in the above
  proof can be extended. 

  Specifically, the argument of a ${\tt BESQ}(1)$ vanishing in a $3$-tree is easily adapted to the situation of $4$-trees, as follows. In order for there to exist 
  $\eta\ge 0$ with $W_4^\eta(j)$ equal for all $j\in[4]$, we need two independent ${\tt BESQ}(1)$ to vanish simultaneously. But they are distributed like squared Brownian 
  motions, and as the origin is polar for planar Brownian motion, this event has probability zero. The same argument applies to any two adjacent edges in $k$-tree evolutions for 
  any $k\ge 4$ and to the evolutions started from any rational time. Similarly, for $6$-trees, the probability vanishes that there is $\eta\ge 0$ with $W_6^\eta(j)$, $j\in[6]$, forming 
  two groups of three equal points.

  Vice versa, any branch point of degree $d$ in the limiting tree $\cT(y)$ has $d-1$ subtrees of positive height and it is an easy consequence of the GHP metric that it must be 
  present in $S_k(\cT_{k,+}^y)$ for large enough $k$, as long as $S(\cT^y_{k,+},k\ge 1)=\cT(y)$. Similarly, if its degree in the limiting tree is finite, this degree is then attained for 
  large enough $k$, and if its degree is infinite, this is approached as $k\rightarrow\infty$. However, we have not been able to prove that Corollary \ref{cor:allemb}, or the 
  corresponding statement in the self-similar setting, can be strengthened to an equality. Furthermore, we would have to consider that there may be other ways of forming 
  higher-degree branch points. Specifically, while each block $\ell$ of $\cT_{k,+}^{y}$ corresponds to a connected component of $\cT(y)\setminus S_k^\circ(\cT_{k,+}^y)$ 
  grafted at $W_k^y(\ell)$, non-binary branch points are formed if the locations $W_k^y(\ell)$ are not distinct. We have discussed when this happens to top masses due to vanishing 
  interval partition mass. It remains to discuss the possibility of vanishing interval partition diversity. 

  As recalled at the end of Section \ref{sec:prelim:SSS}, we showed in \cite{Paper0,Paper1-1,Paper1-2} that diversity equals local time almost surely simultaneously at all levels for the 
  construction of Lemma \ref{lem:type1:pseudo_constr} of a pseudo-stationary type-0 evolution. However, while it is well-known that local times are strictly increasing on the level 
  set almost surely at any fixed level, at exceptional levels this is known to fail at isolated points of level sets (decrease times \cite{Bertoin91,Marsalle98}). We are not aware of any literature that ensures that   
  this cannot happen for non-isolated points. While we believe that this is true, we have been unable to prove this. However, the above arguments (or indeed the proof of 
  \cite[Theorem 1.10]{PaperFV}) would entail that the interval partition evolution does not visit states with finitely many blocks. But we would also like to 
  rule out the possibility that there are exceptional times with infinitely many blocks with zero diversity, which would correspond to branch points of infinite degree.  
\end{remark}

\begin{remark}\label{prop:Zambotti} We would like to remark that Zambotti's process \eqref{eq:Zambotti} does not yield the Aldous diffusion when mapped to a 
  continuum-tree-valued process via the association explained around \eqref{eq:tree_from_path}. While we stop short of proving this, an intriguing way to approach a proof is to 
  look at branch point degrees. 

  Specifically, by \cite[Theorems 7.2--7.3]{Zambotti17}, there are exceptional times when Zambotti's process has at least three zeroes in $(0,1)$, but not when 
  it has five or more. In the continuum tree encoded by such a function, having $j$ zeroes means the root has degree $j+1$, with one 
  subtree encoded by each of the $j+1$ excursions separated by the $j$ zeroes. It therefore suffices to show that the degree of the root is never 4 or 5 under the 
  Aldous diffusion or equivalently under the self-similar Aldous diffusion.

  In the setting of the proof of Theorem \ref{thm:exceptionaldegrees}, recall the argument used to find a ternary branch points in a killed 3-tree evolution. Applying this to the
  ${\tt BESQ}(1)$ process associated with the type-1 evolution yields a root with degree 2 at exceptional times with positive probability. Similarly, considering the root edge and
  an adjacent edge, we find that with probability one, the associated masses do not vanish simultaneously hence preventing the creation of a root with degree 3 or higher in this 
  way. A complete proof would again have to rule out other ways in which the root could have higher degree, as discussed in the preceding remark.   

  As an alternative approach to a proof, we could also use the Markov chain of \cite{EtheLabb15} to see that also the limit can resolve ternary branch points in a way other than how 
  they were formed. Indeed, this sheds some light on the possible degrees of the root since the Aldous diffusion dynamics appear to create blockage through the 
  memory in the binary tree shape that fixes the order of subtrees on any spine. 
  In Zambotti's process, the only reflection is at level 0 and several local minima can descend to zero simultaneously with hardly any interaction. Indeed, in 
  \cite{EtheLabb15}, flipping a few maxima into minima can make a subtree move past a branch point.
\end{remark}

\section{3-sided states and the failure of the Strong Markov property}\label{sec:smp}

In Remark \ref{prop:notstrongmarkov}, we argued that the strong Markov property of the Aldous diffusion fails at times where there is a ternary branch point. We have now developed the tools to establish this rigorously. Specifically, we will study here examples of states with ternary branch points and explain how they can be viewed as states with three reflecting sides that ensure that the Aldous diffusion can
only exit such states on the side from which they were approached. In the following, we consider the Aldous diffusion in the ${\rm GHP}$-closed subset 
$\mathbb{T}^{\rm real}_{\rm unit}\subset\mathbb{T}^{\rm real}$ of ${\rm GHP}$-isometry classes of unit-mass rooted, weighted $\mathbb{R}$-trees, where it takes its values.

\begin{proposition}\label{prop:notstrongmarkov2}
   The Aldous diffusion in the state space $(\bT^{\rm real}_{\rm unit},d_{\rm GHP})$ does not have the strong Markov property. 
\end{proposition}

In the proof, we will consider stopping times, at which the Aldous diffusion hits, with the same positive probability, but from different ``sides'', the following set of (${\rm GHP}$-isometry classes $\mathrm{T}$ of) trees $(T,d,\mu,\rho)$ with a ternary branch point:
\begin{align*}
A=\Big\{\mathrm{T}\in\mathbb{T}^{\rm real}_{\rm unit}\colon\exists_{v\in T}&\exists_{\overset{C_0,C_1,C_2,C_3\subseteq T}{\text{closed connected}}}\forall_{i\neq j}C_i\cap C_j=\{v\},
  \rho\in C_0,\\
   &\ \ \forall_{j}\ \mu(C_j)\ge\textstyle\frac{12-2j}{37},\  \sup\{d(v,x)\colon x\in C_j\}\ge 1\Big\}.
\end{align*} 
In other words, for trees in $A$ there is a point $v\in T$ at which the tree can be split into at least four components of height at least 1 and masses at least $\frac{12}{37}$,
$\frac{10}{37}$, $\frac{8}{37}$, $\frac{6}{37}$, respectively. These \em mass thresholds\em, as well as the \em height threshold \em 1, are of technical help to avoid ambiguities and degeneracies, 
but are of no intrinsic significance beyond their role in finding times in a dense set of times where a ternary branch point exists. In particular, these mass thresholds add to $\frac{36}{37}$ close to 1, and the point $v$ is unique (unless $\mu$ has atoms, which will be addressed appropriately). The set $A$ can be approached in three different ways, distinguishing trees that may have a 
branch separating four components into pairs of mass thresholds, either $(\frac{12}{37},\frac{10}{37})$ and $(\frac{8}{37},\frac{6}{37})$, or $(\frac{12}{37},\frac{8}{37})$ and $(\frac{10}{37},\frac{6}{37})$, or $(\frac{12}{37},\frac{6}{37})$ and $(\frac{10}{37},\frac{8}{37})$:
\begin{align*}
A_1=\Big\{\mathrm{T}\in\mathbb{T}^{\rm real}_{\rm unit}\colon\exists_{v,w\in T}&\exists_{\overset{C_0,C_1,C_2,C_3\subseteq T\setminus]\!]v,w[\![}{\text{closed connected}}}C_0\cap C_1=\{v\},C_2\cap C_3=\{w\}, \rho\in C_0,\\
 &\ \ \forall_{j}\ \mu(C_j)\ge\textstyle\frac{12-2j}{37},
   \ \sup\{d(v,x)\wedge d(w,x)\colon x\in C_j\}\ge 1\Big\},\\
A_2=\Big\{\mathrm{T}\in\mathbb{T}^{\rm real}_{\rm unit}\colon\exists_{v,w\in T}&\exists_{\overset{C_0,C_1,C_2,C_3\subseteq T\setminus]\!]v,w[\![}{\text{closed connected}}}C_0\cap C_2=\{v\},C_1\cap C_3=\{w\},\rho\in C_0,\\
 &\ \ \forall_{j}\ \mu(C_j)\ge\textstyle\frac{12-2j}{37}, 
   \ \sup\{d(v,x)\wedge d(w,x)\colon x\in C_j\}\ge 1\Big\},\\
A_3=\Big\{\mathrm{T}\in\mathbb{T}^{\rm real}_{\rm unit}\colon\exists_{v,w\in T}&\exists_{\overset{C_0,C_1,C_2,C_3\subseteq T\setminus]\!]v,w[\![}{\text{closed connected}}}C_0\cap C_3=\{v\},C_1\cap C_2=\{w\},\rho\in C_0,\\
 &\ \ \forall_{j}\ \mu(C_j)\ge\textstyle\frac{12-2j}{37}, 
   \ \sup\{d(v,x)\wedge d(w,x)\colon x\in C_j\}\ge 1\Big\}.
\end{align*}
We will be interested in the statement of the strong Markov property at the first hitting time of $A$ on the events that the Aldous diffusion is starting and staying in $A_1$, $A_2$ and $A_3$ respectively until hitting $A$. To this end, we first note some properties of these sets of $\mathbb{R}$-trees.

\begin{lemma}\label{lm:smp1}
  \begin{enumerate}\item[(i)] For any $\mathbb{R}$-tree $(T,d,\rho,\mu)$ with diffuse $\mu$ and isometry class in $A$ (respectively $A_i$, $i=1,2,3$), the vertex $v\in T$ (respectively $v,w\in T$) and the components $C_0,C_1,C_2,C_3\subseteq T$ that satisfy all constraints are unique. 
  \item[(ii)] The subsets $A,A_1,A_2,A_3\subset\mathbb{T}^{\rm real}_{\rm unit}$ are closed. 
  \item[(iii)] Consider subsets $A_i^\circ$ of $A_i$, $i=1,2,3$, obtained by replacing all weak inequalities by strict inequalities, removing $v$ and $w$ from $C_0,C_1,C_2,C_3$ for the purpose of satisfying the connectedness constraint and for exceeding the mass thresholds, and also requiring $v\neq w$. Then $A_i^\circ$ is open in $(\mathbb{T}_{\rm unit}^{\rm real},d_{\rm GHP})$, and its closure is a subset of $A_i$, $i=1,2,3$.
  \item[(iv)] Let $\gamma\in(0,\frac{1}{4\times 37})$ and consider subsets $\widetilde{A}_i(\gamma)$ of $A_i$, $i=1,2,3$, obtained by increasing all 
     mass and height thresholds of $A_i$ by $\gamma$ and replacing $C_j$ by the $\gamma$-thickening $C_j^\gamma$ for the purposes of satisfying the intersection 
     constraints. Then $\widetilde{A}_i(\gamma)$ is closed in $\mathbb{T}_{\rm unit}^{\rm real}$ and a subset of $A_i^\circ\cup A$.
  \item[(v)] We have $A_i^\circ\subseteq\bigcup_{\gamma\in(0,1/(4\times 37))}\widetilde{A}_i(\gamma)\subseteq A_i^\circ\cup A$. 
\end{enumerate}
\end{lemma}
 
We prove this in Appendix \ref{appx:smp}.

\begin{lemma}\label{lm:smp2} Consider the Aldous diffusion $(\widebar{\mathcal{T}}(s),s\ge 0)$ of Definition \ref{df:adformal}. Then $\mathcal{T}:=\widebar{\mathcal{T}}(0)$ is a Brownian CRT and $\mathbb{P}\{\mathcal{T}\!\in\! A_1^\circ\}=\mathbb{P}\{\mathcal{T}\!\in\! A_2^\circ\}=\mathbb{P}\{\mathcal{T}\!\in\! A_3^\circ\}>0$. 
  Furthermore, for $i=1,2,3$, consider the exit time $\sigma_{A_i^\circ}=\inf\{s\ge 0\colon\widebar{\mathcal{T}}(s)\not\in A_i^\circ\}$ from $A_i^\circ$. Then $\mathbb{P}(\sigma_{A_i^\circ}>0\,|\,\mathcal{T}\in A_i^\circ)=1$. 
\end{lemma}
\begin{proof} For the first claim, recall from Proposition \ref{prop:B_ktree} the distribution of a Brownian reduced 6-tree $R_6$. In particular, with positive probability for $i=1,2,3$,
  \begin{itemize}\item the uniform tree shape has three type-2 edges,
     \item the $\distribfont{Dirichlet}\big(\frac12,\ldots,\frac12\big)$ mass split assigns mass $M_0>\frac{12}{37}$ to the (type-0) root edge, $M_1+M_2+M_3>\frac{12-2i}{37}$ 
          to the adjacent type-2 edge including its top masses, and masses $M_4+M_5+M_6$, $M_7+M_8+M_9$ exceeding the other two thresholds 
          to the other two type-2 edges including their top masses,
     \item the $\PDIP\big(\frac12,\frac12\big)$ partitions have diversities exceeding $1/\sqrt{M_{3j}}$ for the four edges with masses $M_{3j}$ that are subject to a mass threshold, $j=0,1,2,3$.
  \end{itemize}
  Then the associated weighted $\mathbb{R}$-tree $S_6(R_6)$ has two branch points, which we call $v$ and $w$, and the decomposition of $S_6(R_6)$ around $v$ and $w$ 
  yields branches of weights exceeding the respective mass thresholds and lengths exceeding $\sqrt{M_{3j}}/\sqrt{M_{3j}}=1$, $j=0,1,2,3$. We may assume that $S_6(R_6)$ can be embedded in $\mathcal{T}$, with the mass measure of $\mathcal{T}$ projecting to the mass measure of $S_6(R_6)$, and this easily entails that 
  $\mathbb{P}\{\mathcal{T}\in A_i^\circ\}\ge \mathbb{P}\{S_6(R_6)\in A_i^\circ\}>0$.

  While the above is not the only way that $\{\mathcal{T}\in A_i^\circ\}$ can occur, note that $\mathbb{P}\{S_6(R_6)\in A_i^\circ\}$ does not depend on $i=1,2,3$, by symmetries
  of Brownian reduced 6-trees. For Brownian reduced $k$-trees for larger $k$, we can similarly consider all tree shapes and selections of $v$ and $w$ in which $\{S_k(R_k)\in A_i^\circ\}$ can occur. Then the transformation that,
  on $\{S_k(R_k)\in A_1^\circ\cup A_2^\circ\}$, disconnects, swaps and regrafts the unique (by Lemma \ref{lm:smp1}(i)) pair of components of $S_k^\circ(R_k)\setminus\{v,w\}$ with masses between $\frac{8}{37}$ and $\frac{11}{37}$, leaves the distribution of the Brownian 
  reduced $k$-tree invariant, and so do similar swaps on $\{S_k(R_k)\in A_1^\circ\cup A_3^\circ\}$. Hence, these probabilities do not depend on $i=1,2,3$, either. 

  Conditionally given the event $\{\mathcal{T}\in A_i\}$, the mass and height thresholds are strictly exceeded almost surely, by Corollary \ref{subdecII}. Since the mass measure of $\mathcal{T}$ is diffuse almost surely, there are no 
  atoms at $v$ or $w$, so the mass thresholds are still strictly exceeded when $v$ and $w$ are removed from $C_0,C_1,C_2,C_3$. Hence $\mathcal{T}\in A_i^\circ$ almost surely, i.e.
  $\mathbb{P}\{\mathcal{T}\in A_i\}=\mathbb{P}\{\mathcal{T}\in A_i^\circ\}$. 
  Since the boundary of $A_i^\circ$ is a subset of $A_i\setminus A_i^\circ$, by Lemma \ref{lm:smp1}(iii), and $\mathbb{P}\{\mathcal{T}\in A_i\setminus A_i^\circ\}=0$, and 
  $d_{\rm GHP}(S_k(R_k),\mathcal{T})\rightarrow 0$, the Portmanteau theorem yields $\mathbb{P}\{\mathcal{T}\in A_i^\circ\}=\lim_{k\rightarrow\infty}\mathbb{P}\{S_k(R_k)\in A_i^\circ\}$. Hence, $\mathbb{P}\{\mathcal{T}\in A_1^\circ\}=\mathbb{P}\{\mathcal{T}\in A_2^\circ\}=\mathbb{P}\{\mathcal{T}\in A_3^\circ\}$.

   Since the
  Aldous diffusion is $d_{\rm GHP}$-path-continuous and $A_i^\circ$ is open by Lemma \ref{lm:smp1}(iii), the exit time from the open set $A_i^\circ$ is almost surely positive, conditionally given $\{\mathcal{T}\in A_i^\circ\}$, for each $i=1,2,3$. 
\end{proof}

\begin{proof}[Proof of Proposition \ref{prop:notstrongmarkov2}]
  Consider the Aldous diffusion $(\widebar{\mathcal{T}}(s),s\ge 0)$ of Definition \ref{df:adformal} and Corollary \ref{cor:ad:contversion}, constructed from a consistent family
  $(\widebar{\mathcal{T}}_{k,+}^s,s\!\ge\! 0)$, $k\ge 1$, of stationary unit-mass resampling $k$-tree evolutions, which in turn have been obtained by de-Poissonizing a consistent
  family of self-similar resampling $k$-tree evolutions $(\mathcal{T}_{k,+}^y,y\ge 0)$, $k\ge 1$. Recall the exit times
  \[
   \sigma_{A_i^\circ}=\inf\{s\ge 0\colon \widebar{\mathcal{T}}(s)\not\in A_i^\circ\},\qquad i=1,2,
  \]
  introduced in Lemma \ref{lm:smp2}. For any $k\ge 6$, on the event $\{S_k(\widebar{\mathcal{T}}_{k,+}^0)\in A_1^\circ\}$, we can identify unique 
  $v,w\in S_k^\circ(\widebar{\mathcal{T}}_{k,+}^0)$ such that components of $S_k^\circ(\widebar{\mathcal{T}}_{k,+}^0)\setminus\{v,w\}$ satisfy the mass and height
  constraints for membership in $A_1^\circ$, by Lemma \ref{lm:smp1}(i). Conditionally given $\mathcal{T}_{k,+}^0$, each component gives rise to a number of independent type-$d$ evolutions in 
  $(\mathcal{T}_{k,+}^y,y\ge 0)$, which eventually trigger resampling events. Disconnecting, swapping and regrafting the pair of components of 
  $S_k^\circ(\widebar{\mathcal{T}}_{k,+}^0)\setminus\{v,w\}$ with masses between $\frac{8}{37}$ and $\frac{11}{37}$, as in the proof of Lemma \ref{lm:smp2}, corresponds
  to a change of tree shape of $\mathcal{T}_{k,+}^0$. Inductively, we can use the same type-$d$ evolutions and couple the resampling events to obtain a self-similar $k$-tree
  evolution starting from the initial tree with components swapped, hence in $A_2^\circ$. Furthermore, this construction is naturally consistent in $k$ and, via de-Poissonization,
  mapping under $S$ and passing to a continuous modification, we can associate, on the event $\{\widebar{\mathcal{T}}(0)\in A_1^\circ\}$, an Aldous diffusion 
  $(\widebar{\mathcal{T}}^*(s),s\ge 0)$ starting within $A_2^\circ$. Then the two processes are coupled so that
  \[
  \sigma_{A_2^\circ}^*:=\inf\{s\ge 0\colon \widebar{\mathcal{T}}^*(s)\not\in A_2^\circ\}=\sigma_{A_1^\circ}
  \]
  and for $s\in[0,\sigma_{A_1^\circ}]$, the trees $\widebar{\mathcal{T}}(s)$ and $\widebar{\mathcal{T}}^*(s)$ differ only in that the unique components around the appropriate branch 
  points with masses between $\frac{8}{37}$ and $\frac{11}{37}$ are swapped. In particular, we have
  \[
  \widebar{\mathcal{T}}(\sigma_{A_1^\circ})=\widebar{\mathcal{T}}^*(\sigma_{A_2^\circ})\qquad\mbox{on the event }\big\{\widebar{\mathcal{T}}(\sigma_{A_1^\circ})\in A\big\}.
  \]
  Now we claim that $\mathbb{P}\{\widebar{\mathcal{T}}(\sigma_{A_1^\circ})\in A\}>0$. 
  To show this, we first observe (as a consequence of
  \eqref{eq:bpdist}) the continuity of the map that associates with $\mathrm{T}\in A_1$ the distance $d(v,w)$ between $v$ and $w$, which are unique by 
  Lemma \ref{lm:smp1}(i). Next consider the unit-mass 6-tree evolution 
  $(\widebar{\mathcal{T}}_{6,+}^s,s\ge 0)$ used in the construction of $(\widebar{\mathcal{T}}(s),s\ge 0)$. Recall from the discussion around Proposition \ref{prop:embcont} and Corollary \ref{cor:allemb} that the relationship between $k$-tree evolutions and the Aldous diffusion is rather subtle, due to the passage to a continuous 
  modification, but also that modifications preserve the original process at a dense set of times and at those times, $S_6(\widebar{\mathcal{T}}_{6,+}^s)$ is
  not just embedded into $\widebar{\mathcal{T}}(s)$, but also carries the projected weight measure of $\widebar{\mathcal{T}}(s)$. In the following we fix 
  $\gamma\in(0,\frac{1}{4\times 37})$ and recall the closed subset $\widetilde{A}_1(\gamma)$ of $A_1^\circ\cup A$ from Lemma \ref{lm:smp1}(iv), which 
  we here write as 
  \begin{align*}
  \widetilde{A}_1=\Big\{\mathrm{T}\in\mathbb{T}^{\rm real}_{\rm unit}\colon&\exists_{v,w\in T}\exists_{\overset{C_0,C_1,C_2,C_3\subseteq T\setminus]\!]v,w[\![}{\text{closed connected}}}\,C_0^\gamma\cap C_1^\gamma\!=\!\{v\},C_2^\gamma\cap C_3^\gamma\!=\!\{w\},\rho\in C_0,\\
      &\ \forall_{j}\ \mu(C_j)\ge\textstyle\frac{12-2j}{37}+\gamma,
   \ \sup\{d(v,x)\wedge d(w,x)\colon x\in C_j\}\ge 1+\gamma\Big\}.
  \end{align*} 
  Then having $\widebar{\mathcal{T}}(s)\in\widetilde{A}_1$ on a set of times that is dense in 
  $[0,\sigma_{A_1^\circ}]$ already entails that $\widebar{\mathcal{T}}(s)\in\widetilde{A}_1$ for all $s\in[0,\sigma_{A_1^\circ}]$. Hence,  
  we have
  \[
  \Big\{S_6(\widebar{\mathcal{T}}_{6,+}^s)\in\widetilde{A}_1\mbox{ for all }s\in[0,\sigma_{A_1^\circ}]\Big\}\subseteq\big\{\widebar{\mathcal{T}}(\sigma_{A_1^\circ})\in A\big\}
  \]
  since the membership of the closed set $\widetilde{A}_1\subseteq A_1^\circ\cup A$ transfers from $S_6(\widebar{\mathcal{T}}_{6,+}^s)$ to $\widebar{\mathcal{T}}(s)$. Indeed, on the left-hand event, we have $\sigma_{A_1^\circ}=\inf\{s\ge 0\colon S_6(\widebar{\mathcal{T}}_{6,+}^s)\in A\}$. But this event 
  corresponds, in the self-similar setting before de-Poissonization, to the ${\tt BESQ}(1)$ process of mass in the component between $v$ and $w$ vanishing 
  while independent type-$d$ evolutions of other parts of the process maintain the mass and height constraints (after de-Poissonization). This has positive
  probability. 

  Now assume that the Aldous diffusion satisfies the strong Markov property at $\sigma_{A_1^\circ}$ and at $\sigma_{A_2^\circ}$. Then 
  \begin{align*}
  \mathbb{P}\big(\widebar{\mathcal{T}}(\sigma_{A_1^\circ})\in\cdot\,|\,\widebar{\mathcal{T}}(0)\!\in\! A_1^\circ,\widebar{\mathcal{T}}(\sigma_{A_1^\circ})\!\in\! A\big)
   &=\mathbb{P}\big(\widebar{\mathcal{T}}^*(\sigma_{A_2^\circ})\in\cdot\,|\,\widebar{\mathcal{T}}^*(0)\!\in\! A_2^\circ,\widebar{\mathcal{T}}^*(\sigma_{A_2^\circ})\!\in\! A\big)\\
   &=\mathbb{P}\big(\widebar{\mathcal{T}}(\sigma_{A_2^\circ})\in\cdot\,|\,\widebar{\mathcal{T}}(0)\!\in\! A_2^\circ,\widebar{\mathcal{T}}(\sigma_{A_2^\circ})\!\in\! A\big)
  \end{align*}
  implies that the processes $(\widebar{\mathcal{T}}(\sigma_{A_i^\circ}+u),u\ge 0)$ under 
  $\mathbb{P}(\,\cdot\,|\,\widebar{\mathcal{T}}(0)\in A_i^\circ,\widebar{\mathcal{T}}(\sigma_{A_i^\circ})\in A)$, $i=1,2$,
  have the same distribution. But this is false since the ${\tt BESQ}(1)$ process reflects at 0, while the mass and height constraints are almost surely satisfied strictly at time 
  $\sigma_{A_i^\circ}$ and, by path-continuity, will continue to hold for a positive amount of time after $\sigma_{A_i^\circ}$ so that for $u$ sufficiently small, we will have
  \[
  \mathbb{P}(\widebar{\mathcal{T}}(\sigma_{A_1^\circ}+u)\in A_1^\circ\,|\,\widebar{\mathcal{T}}(\sigma_{A_1^\circ})\in A)=
   \mathbb{P}(\widebar{\mathcal{T}}(\sigma_{A_2^\circ}+u)\in A_2^\circ\,|\,\widebar{\mathcal{T}}(\sigma_{A_2^\circ})\in A)>\frac{1}{2},
   \]
  and this contradicts the equality in distribution since $A_1^\circ$ and $A_2^\circ$ are disjoint.
\end{proof}
 
While we focussed on $A_1^\circ$ and $A_2^\circ$ in the proof of Proposition \ref{prop:notstrongmarkov2}, it should be clear that either can be swapped with $A_3^\circ$, and the same effect is observed. Also the reflection at $A$ observed there when $(S_6(\widebar{T}_{6,+}^s),s\ge 0)$ stays within 
$\widetilde{A}_i=\widetilde{A}_i(\gamma)$ holds for all $i=1,2,3$ and as long as any
$(S_k(\widebar{T}_{k,+}^s),s\ge 0)$, $k\ge 6$, stays within $\bigcup_{\gamma\in(0,1/(4\times 37))}\widetilde{A}_i(\gamma)=A_i^\circ\cup A$, by Lemma \ref{lm:smp1}(v). 
We conclude, as follows.

\begin{corollary} The Aldous diffusion starting in $A_i^\circ$ and stopped when first hitting $A_i\setminus(A_i^\circ\cup A)$ has $A$ as a reflecting boundary for all $i=1,2,3$. 
\end{corollary}

We can therefore think of states in $A$ as having three sides from which they can be reached. Indeed, the behaviour of the Aldous diffusion at $A$ bears some similarities with 
the behaviour of reflecting planar Brownian motion in a disk with a slit removed, $\mathbb{D}=\{z\in\mathbb{C}\colon|z|<1\}\setminus[0,1)$, for which each of the boundary states in $(0,1]$ has two sides,
one in the upper right quarter-disk, one in the lower right quarter-disk. In order to move between these quarter-disks, the reflecting planar Brownian motion has to enter the left half-disk. For the Aldous diffusion
to move between any two of the three $A_i^\circ$, it has to breach the height or mass thresholds. 

Reflecting planar Brownian motion in $\mathbb{D}$ fails to be strongly Markovian in the Euclidean closure, but becomes strongly Markovian if the topology is changed to a compactification that
effectively contains two disjoint copies of $(0,1]$. More generally, the study of Brownian motion in bounded domains with reflection on the boundary is a classical problem that was studied by Fukushima
\cite{Fuku67}. See also \cite[Section 3]{BurdChen98}. One way to equip the compactification of $\mathbb{D}$ with a metric is to use the natural extension of the intrinsic metric on $\mathbb{D}$, which
assigns any two points in $\mathbb{D}$ as their distance the infimum of the Euclidean length of paths between the two points. 

We believe it is worth exploring the generalization of this idea, where the roles of $\mathbb{C}$ and $\mathbb{D}$ are taken by $\mathbb{T}^{\rm real}$ and a suitable subset of binary trees. Such a 
generalization is delicate since suitable sets of binary trees will not be $d_{\rm GHP}$-open and $A$ will only be a small part of the boundary. Indeed, the Aldous diffusion will exhibit similar behaviour at 
every tree with a ternary branch point, and it visits trees with ternary branch points on a dense set of times, by Theorem \ref{thm:exceptionaldegrees}. We will
return to this idea in Section \ref{sec:openprob}.

\section{The modified Aldous chain embedded in the Aldous diffusion}\label{sec:emb}

In this section we show that tree shapes in a unit-mass resampling $k$-tree evolution in stationarity are Markovian and we provide an embedding of the stationary modified Aldous chain of Definition \ref{def:modified_AC}. We use this to show that the Aldous diffusion is the scaling limit of these continuous-time Markov chains and that the Aldous diffusion is reversible.  

\begin{theorem}\label{thm:embedding}
  Consider the consistent system $\big(\widebar{\cT}_{\!k,+}^s,\,s\ge 0\big)$, $k\ge 1$, of stationary unit-mass $k$-tree evolutions. Then the
  associated tree shape evolutions $\big(\widebar{\ft}^s_{k,+},\,s\ge 0\big)$, $k\ge 1$, form a consistent family of stationary continuous-time Markov chains, which have the
  same distribution as the modified Aldous chain embedded into continuous time by independent ${\tt Exponential}(\widetilde{c}_k)$ times between steps, where $\widetilde{c}_k=k(2k-3)$, $k\ge 2$.
  
  More precisely, let $\Lambda(\ft_k,\,\cdot\,)$ denote the kernel from $\bT_k^{\rm shape}$ to $\bT_k^{\rm int}$ that associates with $\ft_k\in\bT_k^{\rm shape}$ the distribution of a unit-mass Brownian reduced $k$-tree conditioned to have shape $\ft_k$. Then the continuous-time modified Aldous chain on $k$-tree shapes is intertwined below the stationary unit-mass $k$-tree evolution via $\Lambda(\ft_k,\,\cdot\,)$, in the sense of Theorem \ref{thm:RogersPitman}.
\end{theorem}

Before we prove this theorem, we need some auxiliary results. We begin in the self-similar regime and study degeneration times, mass evolutions and pseudo-stationarity conditioned on tree shapes. 

\begin{lemma}\label{lm:massindepshape} Consider a killed $k$-tree evolution $(\cT^y,\,y\ge 0)$ starting from a Brownian reduced $k$-tree 
  with any initial mass distribution $\mu$ and conditioned to have tree shape $\ft\in\mathbb{T}^{\rm shape}_k$. Denote its degeneration time by $D$. Then the distribution of the
  total mass process $(\|\cT^y\|,\,0\le y<D)$ does not depend on $\ft$.
\end{lemma}
\begin{proof} First consider $\mu={\tt Gamma}(k-\frac12,\lambda)$, cf.\ the proof of Proposition \ref{prop:pseudo:degen} for related arguments. Recall that Definition 
  \ref{def:killed_ktree} builds such an evolution from independent type-$d$ evolutions for each type-$d$ edge of $\ft$, $d=0,1,2$. More specifically, we may use Construction 
  \ref{interweaving} for all type-2 edges and therefore use as building blocks $k$ independent identically distributed pseudo-stationary type-1 evolutions to obtain the type-$d$ 
  evolutions of our construction for all edges of types $d=1,2$. Furthermore, by Remark \ref{rem:type0fromtype2}, the remainder of this construction for each type-2 evolution also 
  gives rise to an evolution that is a type-0 evolution starting from a ${\tt Gamma}(\frac12,\lambda)$-multiple of ${\tt PDIP}(\frac12,\frac12)$, up to the lifetime of the type-2 
  evolution, in such a way that given this lifetime, the type-0 evolution is independent of the type-2 evolution. These type-0 evolutions are just what is needed in our construction for 
  type-0 edges. The lifetime of the killed $k$-tree evolution is the minimum of the $k$ independent type-1 evolutions. Since every $k$-tree shape has one more type-2 edge than 
  type-0 edges, the joint distribution of this lifetime with the total mass process of the killed $k$-tree evolution is the same for all $k$-tree shapes.

  To deduce the claim for general $\mu$, we proceed as in the proof of Proposition \ref{prop:pseudo:degen}. Specifically, here is a sketch of the argument. For any two tree 
  shapes, we here express the distributional identity of total mass processes when $\mu={\tt Gamma}(k-\frac12,\lambda)$ in terms of expectations of functions of the total mass 
  process. We invert the Laplace transform in $\lambda$ to obtain the result for fixed initial mass. We then integrate the fixed mass result against $\mu$ for the general case.
\end{proof}

Let us investigate pseudo-stationarity results of resampling $k$-tree evolutions conditioned on their initial tree shape. Specifically, we have the
following variants of Propositions \ref{prop:pseudo:pre_D} and \ref{prop:pseudo:degen} when conditioning on an initial tree shape $\ft_0\in\mathbb{T}^{\rm shape}_k$, all in
the following setting. We recall from the beginning of Section \ref{sec:pseudo} notation $Q_{1,[k]}$ for the distribution of a unit-mass Brownian reduced $k$-tree and write $\textsc{shape}\colon\mathbb{T}^{\rm int}_k\rightarrow\mathbb{T}_k^{\rm shape}$ for the map that assigns with $T\in\mathbb{T}^{\rm int}_k$ 
its shape in $\mathbb{T}^{\rm shape}_k$.
\begin{enumerate}
  \item[$\mathbf{(S)}$] Let $\ft_0\in\mathbb{T}^{\rm shape}_k$ and $(\cT^y,\,y\ge 0)$ a resampling $k$-tree evolution, whose initial state is 
  an independent multiple $M$ of a random state with unit-mass distribution $Q_{1,[k]}(\,\cdot\,|\,\textsc{shape}=\ft_0)$. 
\end{enumerate}

We will denote the distribution of $(\cT^y,\,y\ge 0)$ in setting $\mathbf{(S)}$ by $\mathbb{P}_{Q_{\mu,\ft_0}}$, if $M\sim\mu$, and by $\mathbb{P}_{Q_{m,\ft_0}}$ if
$M=m$. We also denote by $\mathbb{P}_T$ the distribution of a resampling $k$-tree evolution starting from $T\in\mathbb{T}_k^{\rm int}$.

\begin{proposition}\label{prop:condpseudo:pre_D} Consider setting $\mathbf{(S)}$. Given $\{D_1>y\}$, the tree $\mathcal{T}^y$ is conditionally an independently scaled Brownian reduced $k$-tree 
  conditioned to have tree shape $\ft_0$.
\end{proposition} 
\begin{proof} 
 Since the tree shape in $(\cT^y,\,y\ge 0)$ does not change before $D_1$, the proof of 
  Proposition \ref{prop:pseudo:pre_D} applies verbatim. 
\end{proof} 

\begin{proposition}\label{prop:condpseudo:degen} Consider setting $\mathbf{(S)}$. Then the following hold for all $n\ge 1$. 
  \begin{enumerate}
    \item[(i)] Given successive labels $I_r=i_r$ causing degeneration, $r\in[n-1]$, and shapes $\ft_{r-1}$ at $D_r-$, $r\in[n]$, we have 
       $I_n=I(\cT^{D_n-})\sim{\tt Unif}([k])$ and the tree shape of $\cT^{D_n}$ is uniformly distributed on the $2k-3$ shapes in $\bT^{\rm shape}_k$ of the form 
       $\ft_{n-1}\oplus(F,J(\ft_{n-1},I_n))$, $F\in\widetilde{\varrho}(\ft_{n-1},I_n)\cup\{\{h\}\colon h\in[k]\setminus\{J(\ft_{n-1},I_n)\}\}$.
    \item[(ii)] Under the conditioning of (i) and further conditioning on $I_n=i_n$ and on the tree shape at $D_n$ being $\ft_n$, the normalized tree $\cT^{D_n}/\|\cT^{D_n}\|$ is a Brownian reduced 
       $k$-tree conditioned to have tree shape $\ft_n$.
    \item[(iii)] Under the conditioning of (ii), the normalized tree $\cT^{D_n}/\|\cT^{D_n}\|$ is independent of $(M,\|\cT^{D_1}\|,\ldots,\|\cT^{D_n}\|,D_1,\ldots,D_n)$
  \end{enumerate}
\end{proposition}
\begin{proof} We refine the proof of Proposition \ref{prop:pseudo:degen}. We denote the random tree shape of $\cT^y$ by $\ft^y$ so that the conditioning on a fixed 
  sequence $(\ft_0,\ldots,\ft_{n-1})$ of tree shapes entails that $\ft^y=\ft_r$ for $D_r\le y<D_{r+1}$, $0\le r\le n-1$. To establish the joint distributions claimed in (i)--(iii), we consider measurable test 
  functions $f_r,g_r\colon[0,\infty)\rightarrow[0,\infty)$, $r\ge 0$, and $H\colon\mathbb{T}_k^{\rm int}\rightarrow[0,\infty)$. Then it suffices to show that
  \begin{align*}
    &\bE_{Q_{\mu,\ft_0}}\Bigg[f_0(\|\cT^0\|)\prod_{r=1}^n\Big(g_r(D_r)f_r(\|\cT^{D_r}\|)\mathbf{1}\{I_r=i_r,\ft^{D_r}=\ft_r\}\Big)H\left(\frac{\cT^{D_n}}{\|\cT^{D_n}\|}\right)\Bigg]\\
    &=\left(\frac{1}{k(2k-3)}\right)^{\!n}\!Q_{1,[k]}[H\,|\,\textsc{shape}=\ft_n]\,\bE_{Q_{\mu,\ft_0}}\Bigg[f_0(\|\cT^0\|)\prod_{r=1}^n\Big(g_r(D_r)f_r(\|\cT^{D_r}\|)\Big)\!\Bigg].
  \end{align*}
  We will prove this by induction on $n$ together with the further claim that the final expectation in this display does not depend on $\ft_0$.  

  First, we reduce the claim to $n=1$ by the strong Markov and self-similarity properties of 
  resampling $k$-tree evolutions at degeneration times. Indeed, if (i)--(iii) hold for $n=1$, then conditionally given the tree shape $\cT^{D_1}$ is $\ft_1$, the post-$D_1$ process 
  satisfies $\mathbf{(S)}$ with $\ft_0$ replaced by $\ft_1$. Furthermore, by (iii), this post-$D_1$ process, after self-similar scaling to start from unit mass, is conditionally independent 
  of $(M,\|\cT^{D_1}\|,D_1)$ given $I_1=i_1$ and given the tree shapes of $\cT^0$ and $\cT^{D_1}$. Inductively, if (i)--(iii) hold with $n$ replaced by $n^\prime=n-1\ge 1$ and 
  $(\ft_{r},i_r,f_r,g_r)$ replaced by $(\ft_{r}^\prime,i_r^\prime,f_r^\prime,g_r^\prime)=(\ft_{r+1},i_{r+1},f_{r+1},g_{r+1})$, $r\in[n]$, then
  \begin{align*}
    &\bE_{Q_{\mu,\ft_0}}\Bigg[f_0(\|\cT^0\|)\prod_{r=1}^n\Big(g_r(D_r)f_r(\|\cT^{D_r}\|)\mathbf{1}\{I_r=i_r,\ft^{D_r}=\ft_r\}\Big)H\left(\frac{\cT^{D_n}}{\|\cT^{D_n}\|}\right)\Bigg]\\
    &=\bE_{Q_{\mu,\ft_0}}\Bigg[f_0(\|\cT^0\|)g_1(D_1)f_1(\|\cT^{D_1}\|)\mathbf{1}\{I_1=i_1,\ft^{D_1}=\ft_1\}\\[-0.05cm]
    &\qquad\quad\bE_{\cT^{D_1}}\Bigg[\prod_{r=1}^{n^\prime}\Big(g_{r}^\prime(D_r)f_{r}^\prime(\|\cT^{D_r}\|)\mathbf{1}\{I_r=i_r^\prime,\ft^{D_r}=\ft_r^\prime\}\Big)H\left(\frac{\cT^{D_{n^\prime}}}{\|\cT^{D_{n\prime}}\|}\right)\Bigg]\Bigg]\\
    &\underset{(*)}{=}\bE_{Q_{\mu,\ft_0}}\Bigg[f_0(\|\cT^0\|)g_1(D_1)f_1(\|\cT^{D_1}\|)\mathbf{1}\{I_1=i_1,\ft^{D_1}=\ft_1\}\\[-0.05cm]
    &\qquad\quad\bE_{Q_{\|\cT^{D_1}\|,\ft_1}}\Bigg[\prod_{r=1}^{n^\prime}\Big(g_{r}^\prime(D_r)f_{r}^\prime(\|\cT^{D_r}\|)\mathbf{1}\{I_r=i_r^\prime,\ft^{D_r}=\ft_r^\prime\}\Big)H\left(\frac{\cT^{D_{n^\prime}}}{\|\cT^{D_{n^\prime}}\|}\right)\Bigg]\Bigg]
  \end{align*}
  allows us to apply first the first part of the induction hypothesis, and then the $n=1$ result (with $H$), to deduce that this further equals
  \begin{align*}
   &\bE_{Q_{\mu,\ft_0}}\Bigg[f_0(\|\cT^0\|)g_1(D_1)f_1(\|\cT^{D_1}\|)\mathbf{1}\{I_1=i_1,\ft^{D_1}=\ft_1\}\\[-0.05cm]
   &\qquad\quad\left(\frac{1}{k(2k-3)}\right)^{\!n^\prime}\!Q_{1,[k]}[H\,|\,\textsc{shape}=\ft_{n^\prime}^\prime]\,\bE_{Q_{\|\cT^{D_1}\|,\ft_1}}\Bigg[\prod_{r=1}^{\!n^\prime}\!\Big(g_r^\prime(D_r)f_r^\prime(\|\cT^{D_r}\|)\Big)\!\Bigg]\Bigg]\\
   &=\left(\frac{1}{k(2k-3)}\right)^{\!n}Q_{1,[k]}[H\,|\,\textsc{shape}=\ft_n]\,\\[-0.05cm]
   &\qquad\quad\bE_{Q_{\mu,\ft_0}}\Bigg[f_0(\|\cT^0\|)g_1(D_1)f_1(\|\cT^{D_1}\|)\bE_{Q_{\|\cT^{D_1}\|,\ft_1}}\Bigg[\prod_{r=1}^{\!n^\prime}\!\Big(g_r^\prime(D_r)f_r^\prime(\|\cT^{D_r}\|)\Big)\!\Bigg]\Bigg].
  \end{align*}

  By the second part of the induction hypothesis, the inner conditional expectation of the product does not depend on $\ft_1$, and by the corresponding statement for $n=1$, the 
  outer expectation does not depend on $\ft_0$. In particular, we can replace $\ft_1$ by $\ft^{D_1}$ and apply the Markov property. 
  More precisely, we argue as follows. Firstly, we insert $\sum_{i\in[k],\ft\in\mathbb{T}^{\rm shape}_k}\mathbf{1}\{I_1=i,\ft^{D_1}=\ft\}$, which equals 1. Secondly, we replace 
  $\ft_1$ by $\ft$ in the conditional expectation of the product. Thirdly, we apply the step $(*)$ of the previous display in reverse for each $i$ and $\ft$ (taking the roles of $i_1$ and $\ft_1$). Fourthly, with no more 
  dependence on $(i,\ft)$, we can remove the sums again. Finally, an application of the strong Markov property completes the induction step.

  It remains to prove the case $n=1$. As in the proof of Proposition \ref{prop:pseudo:degen}, it suffices to consider a killed $k$-tree evolution and a resampling step carried out on
  the left limit at the killing time. By the invariance of the resampling kernel under permutations of labels, we obtain a variant of  \eqref{eq:B_ktree_resamp} that includes shapes. 
  Specifically, we find that resampling $j$
  into a Brownian reduced $(k-1)$-tree with shape $\underline{\ft}\in\mathbb{T}^{\rm shape}_{[k]\setminus\{j\}}$ inserts label $j$ into an edge uniformly chosen from the 
  $2k-3$ edges of $\underline{\ft}$ and then yields a Brownian reduced $k$-trees conditioned on this shape:
  \begin{align*}
 &\int_{T_k\in\TInt_{k}}\int_{T_{k-1}\in\TInt_{k-1}}\!\!Q_{z,[k]\setminus\{j\}}(dT_{k-1}\,|\,\textsc{shape}=\underline{\ft})\Lambda_{j,[k]\setminus\{j\}}(T_{k-1},dT_k)f(T_k)\\
 &=\frac{1}{2k-3}\sum_{F\in\underline{\ft}\cup\{\{h\}\colon h\in[k]\setminus\{j\}\}}\int_{\TInt_k}Q_{z,[k]}(dT\,|\,\textsc{shape}=\underline{\ft}\oplus(F,j))f(T).
  \end{align*} 
  This identifies the claimed conditional distribution of $\ft^{D_1}$ given $I_1$ and means it suffices to establish (i)--(iii) with $\varrho(\cT^{D_1-})/\|\cT^{D_1-}\|$ and 
  $\widetilde{\varrho}(\ft_0,i_1)$ instead of $\cT^{D_1}/\|\cT^{D_1}\|$ and $\ft_1$. The remaining claim, including the second claim that the distribution of $(\|\cT^0\|,D_1,\|\cT^{D_1-}\|)$ does not depend on $\ft_0$, only depends on the killed $k$-tree evolution. 

  With fixed initial tree shape
  $\ft_0=\ft$, the argument of Proposition \ref{prop:pseudo:degen} still yields $I=I_1\sim{\tt Unif}([k])$ independent of $D=D_1$. The distribution of $J=J(\ft,I)$ identified there 
  changes, but we leave this implicit here. More importantly, on the event $\{I=i\}$, we have an induced tree shape 
  $\underline{\ft}_i=\widetilde{\varrho}(\ft,i)\in\mathbb{T}^{\rm shape}_{[k]\setminus\{J(\ft,i)\}}$ after swap-reduction. The further arguments of 
  Proposition \ref{prop:pseudo:degen} now yield that, conditionally given this tree shape, 
  $\varrho(\cT^{D-})/\|\cT^{D-}\|\sim Q_{1,[k]\setminus\{J(\ft,i)\}}(\,\cdot\,|\,\textsc{shape}=\underline{\ft}_i)$ that is conditionally independent of
  $(\|\cT^0\|,D,\|\cT^{D-}\|)$. The second claim follows from Lemma \ref{lm:massindepshape}.
\end{proof}

In particular, we read off the distribution of the tree shapes at resampling times.

\begin{corollary}\label{cor:ac} Consider setting $\mathbf{(S)}$. Then the tree shapes $(\ft^{D_n},\,n\ge 0)$ at $D_0=0$ and at resampling times $D_n$, $n\ge 1$, evolve according to the
  modified Aldous chain of Definition \ref{def:modified_AC}. 
\end{corollary}

\begin{corollary}[Conditional strong pseudo-stationarity]\label{cor:condpseudo}
 Consider setting $\mathbf{(S)}$. Denote by $M(y)$ and $\ft(y)$ the total mass and shape of 
  $\cT^y$ and by $\cF^y_{\rm mass+shape}$, $y\ge 0$, the filtration they generate.  Let $Y$ be a stopping time in this filtration. Then for all $\ft^\prime\in\mathbb{T}^{\rm shape}_{[k]}$ and all 
  $\cF^Y_{\rm mass+shape}$-measurable $\eta\colon\Omega\rightarrow[0,\infty)$ and measurable $H\colon\mathbb{T}^{\rm int}_k\rightarrow[0,\infty)$,
  \begin{align}\label{eq:condpseudo}
    &\bE\left[\eta\mathbf{1}\{M(Y)>0,\ft(Y)=\ft^\prime\}H(\cT^Y)\right]\\
    &=\bE\left[\eta\mathbf{1}\{M(Y)>0,\ft(Y)=\ft^\prime\} Q_{M(Y),[k]}[H\,|\,\textsc{shape}=\ft^\prime]\right].\nonumber
  \end{align} 
\end{corollary}
\begin{proof} As in the proof of Proposition \ref{prop:pseudo:resamp}, we deduce from Propositions \ref{prop:condpseudo:pre_D} and \ref{prop:condpseudo:degen} 
  corresponding statements for any fixed time by conditioning on the resampling times. Specifically, we find that for all $y\ge 0$ and $\ft^\prime\in\mathbb{T}_k^{\rm shape}$, 
  conditionally given $\{M(y)>0,\ft(y)=\ft^\prime\}$, the tree $\cT^y$ is an independently scaled Brownian reduced $k$-tree conditioned to have tree shape $\ft^\prime$.

  Based on these conditional distributions, we adapt the proof of Lemma \ref{pretotal}, as follows. First suppose that $Y=y$ is non-random. An induction yields that for all $0=y_0<y_1<\cdots<y_n=y$ 
  \begin{align*}  
    &\mathbb{E}_{Q_{m,\ft}}\Bigg[\prod_{r=0}^nf_r(\ft(y_r),M(y_r))\mathbf{1}\{M(y)>0,\ft(y)=\ft^\prime\}H(\cT^y)\Bigg]\\
    &=\mathbb{E}_{Q_{m,\ft}}\Bigg[\prod_{r=0}^nf_r(\ft(y_r),M(y_r))\mathbf{1}\{M(y)>0,\ft(y)=\ft^\prime\}Q_{M(y),[k]}[H\,|\,\textsc{shape}=\ft^\prime]\Bigg]
  \end{align*}
  and a monotone class theorem establishes \eqref{eq:condpseudo} when $Y=y$ is non-random. The generalization to stopping times by discretization and right-continuity is 
  again standard.
\end{proof}

\begin{proof}[Proof of Theorem \ref{thm:embedding}]  In the following, we denote by $\widebar{\mathbb{P}}_{1,\ft}$ (respectively $\mathbb{P}_{m,\ft}$) the distribution of a unit-mass (respectively self-similar) resampling $k$-tree evolution starting from a Brownian reduced $k$-tree (of mass $m$ in the self-similar case) conditioned to have tree shape $\ft\in\mathbb{T}^{\rm shape}_k$ for some $k\ge 1$. Let $\big(\cT^y,\,y\ge 0\big)$ $\sim\mathbb{P}_{\mu,\ft}=\int_{(0,\infty)}\mathbb{P}_{m,\ft}\mu(dm)$ with induced tree shape evolution $(\ft^y,\,y\ge 0)$. Recall that the de-Poissonization stopping times $\rho_{\boldsymbol{\mathcal{T}}}(s)$, $s\ge 0$, satisfy $\|\cT^{\rho_{\boldsymbol{\mathcal{T}}}(s)}\|>0$ a.s., as noted below \eqref{eq:dePoi_time_change}. By Corollary \ref{cor:condpseudo} applied to $Y=\rho_{\boldsymbol{\mathcal{T}}}(s)$, we find for  
  $\cF^Y_{\rm mass+shape}$-measurable $\eta\colon\Omega\rightarrow[0,\infty)$, measurable $h\colon\mathbb{T}^{\rm int}_k\rightarrow[0,\infty)$ and
  $\ft^\prime\in\mathbb{T}^{\rm shape}_k$
  \[
  \mathbb{E}_{m,\ft}\Big[\eta\cf\{\ft^{\rho_{\boldsymbol{\mathcal{T}}}(s)}=\ft^\prime\}h\big(\cT^{\rho_{\boldsymbol{\mathcal{T}}}(s)}/\|\cT^{\rho_{\boldsymbol{\mathcal{T}}}(s)}\|\big)\Big]=\mathbb{E}_{m,\ft}\Big[\eta\cf\{\ft^{\rho_{\boldsymbol{\mathcal{T}}}(s)}=\ft^\prime\}\Big]\mathbb{E}_{1,\ft^\prime}\Big[h\big(\cT^0\big)\Big].
  \]
  This readily entails 
  \[
  \widebar{\mathbb{E}}_{1,\ft}\Big[\widebar{\eta}\cf\{\widebar{\ft}^s=\ft^\prime\}h\big(\widebar{\cT}^s\big)\Big]=\widebar{\mathbb{E}}_{1,\ft}\Big[\widebar{\eta}\cf\{\widebar{\ft}^s=\ft^\prime\}\Big]\mathbb{E}_{1,\ft^\prime}\Big[h\big(\cT^0\big)\Big],
  \]
  where $\widebar\eta\colon\Omega\rightarrow[0,\infty)$ is measurable in $\widebar\cF^s_{\rm mass+shape} = \cF^{\rho_{\boldsymbol{\mathcal{T}}}(s)}_{\rm mass+shape}$. 
  In the notation of intertwining of Theorem \ref{thm:RogersPitman}, this means that $\Lambda P_s=\Lambda P_s\Phi\Lambda$, where $P_s$ is the transition kernel of the 
  unit-mass resampling $k$-tree evolution, $\Phi$ is the kernel associated with projection to tree shape and $\Lambda$ is the kernel stated in the statement of the theorem.
  Since also $\Lambda\Phi$ is the identity kernel on $\mathbb{T}^{\rm shape}_k$ and $\widebar{\cT}^0$ has distribution $\Lambda(\ft,\,\cdot\,)$, Theorem
  \ref{thm:RogersPitman} applies and yields that $(\widebar{\ft}^s,\,s\ge 0)$ is Markovian.

  To make the transition kernel of $(\widebar{\ft}^s,\,s\ge 0)$ more explicit, note that by construction, the tree shape stays constant between the times 
  $\widebar{D}_r$ such that $\rho_{\boldsymbol{\mathcal{T}}}(\widebar{D}_r)=D_r$, $r\ge 0$, where we write $D_0:=0$. By Lemma \ref{lm:massindepshape}, the distribution of 
  $\widebar{D}_1$ under $\widebar{\bP}_{1,\ft}$ does not depend on the initial tree shape $\widebar{\ft}^0=\ft$, since $\widebar{D}_1$ only depends on 
  $(\|\cT^y\|,\,0\le y<D_1)$. 

  In a continuous-time Markov chain whose holding times $\widebar{D}_r-\widebar{D}_{r-1}$, $r\ge 1$, are identically distributed, they are furthermore independent exponentially 
  distributed and independent of the jump chain $(\ft^{D_r},\,r\ge 0)$. This also follows inductively from the independence noted in Proposition \ref{prop:condpseudo:degen}(iii) 
  in conjunction with the further independence from tree shapes. To summarize,
  \[
  \widebar{\mathbb{E}}_{1,\ft}\!\Bigg[\!\prod_{r=1}^{n}\!g_r\big(\widebar{D}_{r}-\widebar{D}_{r-1}\big)\cf\big\{I(\widebar{\cT}^{\widebar{D}_r-})\!=\!i_r,\widebar{\ft}^{\widebar{D}_r}\!=\!\ft_r\big\}\!\Bigg]
  \!=\!\Big(\frac{1}{k(2k-3)}\Big)^{\!n}\prod_{r=1}^n\!\widebar{\mathbb{E}}_{1,\ft}\Big[g_r\big(\widebar{D}_1\big)\Big].
  \]
  While it does not seem straightforward to determine the distribution of $\widebar{D}_1$ directly, we obtained from the Markov property of $(\widebar{\ft}^s,\,s\ge 0)$ that it is 
  exponential.  We denote the rate parameter by $\widetilde{c}_k$. We further identify the transition probabilities of $(\widebar{\ft}^{\widebar{D}_r},\,r\ge 0)$ as the ones of the 
  modified Aldous chain, cf.\ Corollary \ref{cor:ac}.

  To relate the rate parameters, we note that in a stationary unit-mass resampling $(k+1)$-tree evolution with rates $\widetilde{c}_{k+1}$, the resampling label $J$ in the sense of 
  Proposition \ref{prop:pseudo:degen}(i) is $k+1$ with probability $(4k-1)/(k+1)(2k-1)$. By Poisson thinning, this means that 
  \[
  \widetilde{c}_k=\widetilde{c}_{k+1}\Big(1-\frac{4k-1}{(k+1)(2k-1)}\Big)=\widetilde{c}_{k+1}\frac{k(2k-3)}{(k+1)(2k-1)}.
  \]
  Hence, there is $c\in(0,\infty)$ such that $\widetilde{c}_k=ck(2k-3)$. Since $\widetilde{c}_2=2c$, the following proposition entails that $c=1$, and this completes the proof.
\end{proof}

\begin{proposition}\label{lm:degrate} The first resampling time of a stationary unit-mass 2-tree evolution as in Theorem \ref{thm:stationary} is exponentially distributed with rate parameter 2. 
\end{proposition}

We prove this lemma in Appendix \ref{appx:degrate}.

We further note that with probability $1/(2k-3)$, the state of the modified Aldous chain does not change. Hence, the actual jump rate of the continuous-time Markov chain is further thinned to $c_k=k(2k-4)=2k(k-2)$, which is reminiscent of, but not the same as the rates appearing in Kingman's coalescent. But still, in a non-resampling evolution, these rates are such that it should be possible to show 
that the Aldous diffusion comes down from infinity in the sense that, starting in stationarity and with suitable labelling conventions, the labels perform a variant of Kingman's coalescent.

The embedding of Theorem \ref{thm:embedding} allows us to deduce scaling limits. The following corollary restates Theorem \ref{thm:scalinglim} in the context of Theorem \ref{thm:embedding} and completes the resolution of Conjecture \ref{conj:Aldous}.

\begin{corollary}\label{cor:embconv8}
 This continuous-time (modified) Aldous chain, represented as a process of $\mathbb{R}$-trees with edge lengths $1/\sqrt{k}$ and uniform weight measure on the leaves, 
  converges to the Aldous diffusion as $k\to\infty$, in the sense of finite-dimensional distributions on $(\mathbb{T}^{\rm real},d_{\rm GHP})$. Indeed, in the coupling provided by Theorem \ref{thm:embedding}, the convergence holds almost surely.
\end{corollary}
\begin{proof}[Proof of Corollary \ref{cor:embconv8} and Theorem \ref{thm:scalinglim}] In the setting of Theorem \ref{thm:embedding}, consider the Aldous diffusion $\widebar{\cT}(s)=S\big(\widebar{\cT}^s_{k,+},k\ge 1\big)$, $s\ge 0$. Then for 
  each $s\ge 0$, the joint distribution of $(\widebar{\ft}^s_k,\,k\ge 1)$ with the limiting CRT $\widebar{\cT}(s)$ is as in Curien and Haas \cite[Theorem 5(iii)]{CuriHaas13}. Hence, the
  convergence as claimed holds almost surely for each $s\ge 0$. This entails the almost sure convergence at any finite number of times, which also entails the convergence in the 
  sense of finite-dimensional distributions, as claimed. 
\end{proof}

\begin{remark} The holding times are i.i.d.\ exponential variables that are coupled for different $k$ so that they form a superposition of Poisson processes of rates 
$\widetilde{c}_j-\widetilde{c}_{j-1}$, $j\ge 2$. The points of a Poisson process in any interval are distributed like a Poisson number of i.i.d.\ 
uniform random variables. Now consider the time change that maps equi-distant times at multiples of $1/k^2$ to the times of the points of the superposition up to level $k$. It follows from the strong law of large numbers for Poisson variables and the Glivenko--Cantelli theorem for the uniform random variables that these time changes converge uniformly in any time interval to a linear function. If we had uniform rather than finite-dimensional convergence in Corollary \ref{cor:embconv8}, we could deduce the convergence of the discrete Aldous chain.
\end{remark}

Finally, we turn to the question of reversibility of the Aldous diffusion. 
Specifically, recall the (unmodified) Aldous chain of Section \ref{sec:adconj}. We note as Schweinsberg \cite{Schweinsberg02} did in the unrooted case, that the symmetry of the transition rules entails the reversibility with the uniform stationary distribution. The modified label dynamics of Definition \ref{def:modified_AC}, on the other hand, sacrifice reversibility for $k\ge 4$ by allowing certain forward moves that cannot be reversed in one step. However, as noted before,
the two chains share the same dynamics of the underlying unlabeled trees (and representations in $\mathbb{T}^{\rm real}$), and indeed, this will allow us to apply the reversibility of the (unmodified) Aldous chain. None of these observations is affected by
passing to continuous time, with steps separated by independent exponential times. We can now prove Corollary \ref{cor:rev:intro}, which we restate here.

\begin{corollary}\label{cor:reversible} The Aldous diffusion is reversible with the distribution of the Brownian CRT.
\end{corollary}
\begin{proof}[Proof of Corollaries \ref{cor:reversible} and \ref{cor:rev:intro}] For each $k\ge 1$, denote by $(\ft_k(s),s\ge 0)$ a stationary continuous-time (unmodified) Aldous chain with steps separated by independent exponential times with the rates 
  $\widetilde{c}_k=k(2k-3)$ of Theorem \ref{thm:embedding}. This continuous-time chain is reversible and therefore $(\ft_k(0),\ft_k(s))$ and $(\ft_k(s),\ft_k(0))$ have the same distribution for each 
  $s\ge 0$. 

  Adapting ideas from Definition \ref{def:SkandSkcirc} to the simpler setting without edge partitions, we denote by 
  $\tau_{\sqrt{k}}\colon\mathbb{T}_k^{\rm shape}\rightarrow\mathbb{T}^{\rm real}$ the map that assigns with a tree shape the ${\rm GHP}$-isometry class of an $\mathbb{R}$-tree representative of 
  $\mathbf{t}$ with edge lengths $1/\sqrt{k}$ and uniform weight measure on the leaves. 

  Now fix $s\ge 0$. In the setting of Theorem \ref{thm:embedding}, Corollary \ref{cor:embconv8} yields that 
  $d_{\rm GHP}(\tau_{\sqrt{k}}(\widebar{\mathbf{t}}_{k,+}^0),\widebar{\mathcal{T}}(0))\rightarrow 0$ and 
  $d_{\rm GHP}(\tau_{\sqrt{k}}(\widebar{\mathbf{t}}_{k,+}^s),\widebar{\mathcal{T}}(s))\rightarrow 0$ almost surely. Then for any bounded continuous function 
  $f\colon\mathbb{T}^{\rm real}\times\mathbb{T}^{\rm real}\rightarrow[0,\infty)$, we have
  \begin{align*}
  \mathbb{E}\left[f\left(\widebar{\mathcal{T}}(0),\widebar{\mathcal{T}}(s)\right)\right]
  &=\lim_{k\rightarrow\infty}\mathbb{E}\left[f\left(\tau_{\sqrt{k}}(\widebar{\mathbf{t}}_{k,+}^0),\tau_{\sqrt{k}}(\widebar{\mathbf{t}}_{k,+}^s)\right)\right]\\
  &=\lim_{k\rightarrow\infty}\mathbb{E}\left[f\left(\tau_{\sqrt{k}}(\ft_k(0)),\tau_{\sqrt{k}}(\ft_k(s))\right)\right]\\
  &=\lim_{k\rightarrow\infty}\mathbb{E}\left[f\left(\tau_{\sqrt{k}}(\ft_k(s)),\tau_{\sqrt{k}}(\ft_k(0))\right)\right]\\
  &=\lim_{k\rightarrow\infty}\mathbb{E}\left[f\left(\tau_{\sqrt{k}}(\widebar{\mathbf{t}}_{k,+}^s),\tau_{\sqrt{k}}(\widebar{\mathbf{t}}_{k,+}^0)\right)\right]
  =\mathbb{E}\left[f\left(\widebar{\mathcal{T}}(s),\widebar{\mathcal{T}}(0)\right)\right]
  \end{align*}
  and this completes the proof.
\end{proof}

\section{Open problems}\label{sec:openprob}

The definition of the Aldous diffusion, Definition \ref{df:adformal}, is in a stationary setting. The Markov property established in Theorem \ref{thm:ad:markov} 
gives rise to transition kernels $\widetilde{\kappa}_s$, $s\ge 0$, that are defined ${\tt BCRT}$-almost everywhere, satisfy a semi-group property ${\tt BCRT}$-almost everywhere, and can be paraphrased, as follows. 
\begin{itemize}
\item Sample a sequence of leaves from the mass measure of the initial unit-mass (binary) continuum tree, 
\item use the root and the first $k$ leaves to build an initial reduced $k$-tree in $\widebar{\mathbb{T}}^{\rm int}_{k}$, for all $k\ge 1$, forming an 
  initial consistent family in $\widebar{\mathbb{T}}^{\rm int}_\infty$,
\item run a $\widebar{\mathbb{T}}^{\rm int}_\infty$-valued evolution of unit-mass $k$-tree evolutions for time $s\ge 0$,
\item use the function $S$ of Definition \ref{def:S} to project the consistent family at time $s$ into $\mathbb{T}^{\rm real}$.
\end{itemize}
This is sufficient to establish the stationary process as a simple Markov process and to derive some other properties, as we have demonstrated, but raises 
further questions, whose answers may open up a more direct analytic study of the Aldous diffusion without relying on the delicate consistency in stationarity of resampling $k$-tree evolutions in Corollary \ref{cor:consistent_fam}, which gives rise to the $\widebar{\mathbb{T}}^{\rm int}_\infty$-valued evolution.\pagebreak

\begin{problem}\label{prob1} Identify a state space of (binary) continuum trees from which the Aldous diffusion with transition kernels $\widetilde{\kappa}_s$, $s\ge 0$, can 
  start. Provide an explicit description of $\widetilde{\kappa}_s$, $s\ge 0$, as a family of kernels on this state space that satisfies the semi-group property everywhere. Identify the generator of the Aldous diffusion.
\end{problem}

Indeed, a restriction to binary continuum trees is necessary here, because for any continuum tree with ternary or higher-degree branch points and any sequence of leaves sampled from its mass measure, there will not be a unique way 
to associate a consistent family in $\widebar{\mathbb{T}}^{\rm int}_\infty$. More precisely, a reduced $k$-tree constructed from labels in all four 
components around a ternary branch point must split the four labels into two pairs separated by an empty edge partition. Each of the three ways of pairing up four labels will typically lead to a different continuum-tree-valued evolution in much the same way as the Aldous diffusion resolves ternary branch points instantaneously into two binary branch points after the stopping times explored in Section \ref{sec:smp}. 

On the other hand, a state space of binary trees is insufficient for a continuous modification, by Theorem \ref{thm:exceptionaldegrees}. As indicated in Remark \ref{rmk:higherdegree}, we
believe it is sufficient to allow one ternary and no higher-degree branch points.

\begin{problem} Show that the Aldous diffusion almost surely has no times at which there is any branch point of degree 5 or higher or more than one ternary
  branch point.
\end{problem}

We discussed in Remark \ref{prop:Zambotti} that this would also be one way, but not the only way, to approach the following problem.

\begin{problem} Prove rigorously that Zambotti's process \eqref{eq:Zambotti} on a space of excursions does not yield the Aldous diffusion when mapped to a 
  continuum-tree-valued process via the association explained around \eqref{eq:tree_from_path}.  
\end{problem}

This clearly leaves open Zambotti's problem \cite[Section 5.6.4]{Zambotti17} of providing a description of his process as an evolution of trees, which would 
naturally take place in $(\mathbb{T}^{\rm real},d_{\rm GHP})$, or in a subspace. The following is the complementary problem.

\begin{problem} Construct and study an \em excursion-valued Aldous diffusion \em that projects to the Aldous diffusion via the association explained around \eqref{eq:tree_from_path}. 
\end{problem}

From the perspective of our construction, this involves the construction of consistent planar structures. Intuitively, this can be done using independent Bernoulli variables 
to make tree shapes planar. To construct excursions, it seems useful to also associate Bernoulli variables with every block in an interval partition, to indicate whether the corresponding subtree is to the left or to the right of the branch represented by the interval partition. This also helps set up consistent evolutions where naturally the Bernoulli variables are associated with atoms of the Poisson random measures of Section \ref{sec:prelim:SSS}. Following this route would involve revisiting many developments of this memoir with added structure. Could  
other techniques make use of the less abstract state space of excursions, using some of the insights from this memoir that translate easily between frameworks?

An excursion-valued Aldous diffusion is a richer object as it encodes planar structure. Planar structure plays an important role in some applications of continuum trees, notably 
to random planar maps and Liouville quantum gravity \cite{MatingTrees,LeGallMier12,MarcMokk}. The planar order of the BCRT was already considered by Aldous \cite{AldousCRT3}, as was the general coding of ordered graph-theoretic trees as walks. See also Le Gall \cite{LeGall91}. There is also a more recent literature on the representation of order structure of $\mathbb{R}$-trees. See e.g. Evans et al.\ \cite{EvanGrueWako2017,EvanWako2020}.

We do not believe that an excursion-valued Aldous diffusion would be strongly Markovian, but the additional order structure appears to reduce the 3-sided
nature of states explored in Section \ref{sec:smp} to 2-sided states. Specifically, a continuum tree with a ternary branch point corresponds to an excursion
with three adjacent sub-excursions above the same level. A resolution into two binary branch points corresponds to one of the touch points moving above the 
other, which can happen in only two ways, the third being ruled out by the planar order -- the subtrees corresponding to the left and right excursion cannot form
a third pair. 

As far as the strong Markov property is concerned, other approaches seem more promising. One possibility is to discard the metric structure and work on a state space of rooted
algebraic trees building on the work of L\"ohr et al.\ \cite{LohrMytnWint20,LohrWint21} in the unrooted case.
\begin{problem} Define a topology on a space of rooted algebraic trees. Show that the Aldous diffusion projected to rooted algebraic trees is a strong Markov process, a \em rooted algebraic Aldous diffusion\em.
\end{problem}
This problem may alternatively be addressed in Forman's space of interval-partition trees \cite{IPTrees}. Either way, the metric structure of the continuum tree  would be removed from the state space. A potential alternative may be changing the metric on a state space of binary continuum trees to a 
metric that measures the distance between two binary continuum trees as the infimum of ${\rm GHP}$-lengths of paths in this space of binary continuum trees. 
This space is not complete, but approaching the same continuum tree with a ternary branch point from the three sides indicated in Section \ref{sec:smp} appears to correspond to three distinct
points in a completion, which we denote by $(\widebar{\mathbb{T}}^{\rm real}_{\rm binary},d_{\rm GHP}^{\rm path})$. Following \cite{EPW,EW}, it is easy to see that this completion is a separable metric space. 
\begin{problem} Show that the Aldous diffusion is a path-continuous strong Markov process in $(\widebar{\mathbb{T}}^{\rm real}_{\rm binary},d_{\rm GHP}^{\rm path})$.
\end{problem}
Intuitively, the self-similar Aldous diffusion of Section \ref{sec:continuity} inherits the independence of the evolution of subtrees in subtree decompositions that is expressed in the underlying $k$-tree evolutions (modulo labels, which only play an auxiliary role when mapping into $(\mathbb{T}^{\rm real},d_{\rm GHP})$).   
\begin{problem} Describe the evolution of the subtree decomposition of Corollary \ref{subdecII} under the self-similar Aldous diffusion, until one of the subtrees corresponding to the top masses vanishes.
\end{problem}
It would be particularly interesting to identify a $\sigma$-finite measure that describes the evolution of new subtrees created during this evolution. In the following, we refer to this $\sigma$-finite measure as the \em excursion measure \em of the self-similar Aldous diffusion. 
\begin{problem} Study the self-similar Aldous diffusion under its excursion measure. Study de-Poissonization under the excursion measure of the self-similar Aldous diffusion.
\end{problem} 
Following any fixed time, with an entrance law of the excursion measure of the self-similar Aldous diffusion, the forward evolution is a self-similar Aldous diffusion and its de-Poissonization must give a unit-mass
Aldous diffusion, and the backward evolution can be approached using the reversibility of Corollary \ref{cor:reversible}. Returning to the context of limit theorems of Section \ref{sec:emb}, we also pose the following problem.
\begin{problem}\label{prob:funcconv} Strengthen the convergence of the embedded continuous-time Aldous chain in Corollary \ref{cor:embconv8} to functional convergence. Show the convergence of the discrete-time Aldous chain to the Aldous diffusion. Identify other Markov chains that converge to the Aldous diffusion.
\end{problem}
In this memoir, we have approached problems about continuum trees via embedded $k$-trees, spinal decompositions and Poissonization. It is 
instructive to do the same to approach Problem \ref{prob:funcconv}. Recall the Poissonized ${\tt oCRP}(\frac12,0)$ from Section \ref{sec:discrete_to_cts}. 
In the Poissonized (modified) Aldous chain, every leaf (and adjacent branch point) is deleted at rate 1, into every edge a new branch point (and adjacent leaf) is inserted at rate $\frac12$. Given the subtree spanned by two leaves and the root of the initial tree, the associated 2-tree has two top masses evolving as birth-and-death chains and an evolving vector of leaf counts in spinal subtrees until one of the top masses vanishes. We can view this as a first top mass, and a second top mass followed by the vector of spinal leaf counts. The so-extended vector evolves as an ${\tt oCRP}(\frac12,0)$ independently of the first top mass, jointly stopped when the first top mass vanishes. This is the discrete analogue of the type-2 evolution of Definition \ref{def:type2:v1} run as a type-1 evolution and an independent ${\tt BESQ}(-1)$-top mass until one of the two top masses vanishes.

It was shown in \cite{RogeWink20,ShiWinkel-2} that this ${\tt oCRP}(\frac12,0)$ has a type-1 evolution as its scaling limit. This scaling limit holds as functional 
convergence \cite[Theorem 3.12]{ShiWinkel-2} if represented in a space of interval partitions equipped with the Hausdorff distance, or equivalently the distance obtained in Definition 
\ref{def:IP:metric} if ${\rm dis}(\beta,\gamma,(U_j,V_j)_{j\in[n]})$ is replaced by ${\rm dis}_H(\beta,\gamma,(U_j,V_j)_{j\in[n]})$ defined to be the 
maximum of just (iii) and (iv) in Definition \ref{def:IP:metric} hence ignoring diversities. This is proved by first showing the convergence of the discrete 
scaffolding-and-spindles construction of Section \ref{sec:discrete_to_cts} to the continuous scaffolding-and-spindles construction of Section \ref{sec:prelim:SSS}. Distances in the continuum trees are diversities of interval partitions and local times of the scaffolding L\'evy processes. Establishing the functional convergence of 
distances in the Poissonized Aldous chain is therefore closely related to the functional convergence of L\'evy process local times in the following sense.

\begin{problem} Consider a sequence $X^n$ of spectrally positive compound Poisson processes compensated to have zero mean by adding a negative drift.
  Suppose that $X^n\rightarrow X$ weakly, where $X$ is an unbounded variation L\'evy process with bi-continuous occupation density local time process $L$.
  Show that the occupation density local times $L^n$ of $X^n$ converge weakly to $L$.
\end{problem}

For a Brownian motion limit, this problem was addressed by Khoshnevisan \cite{Khoshnevisan1993} and Lambert et al.\ \cite{LambSimaZwar13}.
A general finite-dimensional convergence result (even without assuming bi-continous limiting local times) was proved by Lambert and Simatos 
\cite[Theorem 2.4]{LambSima15}. They also illustrate for a specific heavy-tailed jump distribution in the relevant domain of attraction of a stable process that tightness holds. The argument is very technical, but any obstacles to proving the corresponding result in the setting of \cite[Theorem 1.5]{RogeWink20} appear
to be technical in nature rather than any suspected lack of tightness. In any case, this is only a first step or practice step towards establishing tightness in Problem 
\ref{prob:funcconv}. 

Finally, 
recall from Section 7.7 the discussion of generalizations of the Aldous diffusions to other (pseudo-)stationary continuum-tree-valued evolutions. Let us here formulate the three examples as open problems.

\begin{problem} Construct a unit-mass and self-similar continuum-tree-valued evolutions for Ford's CRT \cite{For-05,HMPW} that relate to consistent systems of $k$-tree evolutions with non-exchangeable labels, using as building blocks the type-0 and type-1 evolutions of \cite{Paper1-2}.
\end{problem}

\begin{problem} Construct unit-mass continuum-tree-valued evolutions associated with down-up Markov chains derived from strongly sampling consistent
  Markov branching models in the domain of attraction of a binary self-similar CRT of \cite{HaasMier12,HMPW}.
\end{problem}

\begin{problem} Construct unit-mass and self-similar continuum-tree-valued evolutions for the stable CRTs of \cite{DuquLeGall02,HM04} based on a down-up chain
  whose up-steps are Marchal's growth procedure \cite{Marchal08} and building on nested interval partition evolutions of \cite{ShiWinkel-2}.
\end{problem}

\appendix
\chapter*{Appendix}

\stepcounter{chapter}

This appendix is a collection of material mostly of a technical nature. Section \ref{sec:DynkinIntertwining} reviews Dynkin's criterion and intertwining, which are used throughout to show where functions of Markov processes are Markovian. The remainder contains technical proofs of results stated in the main chapters,  Sections \ref{sec:marked_tree_metric}--\ref{sec:non_acc_2} are relevant in Chapter \ref{ch:consistency}, Sections \ref{sec:contproj}--\ref{app:pfpropGHP} in Chapter \ref{ch:properties} and Sections \ref{appx:markedGH}--\ref{appx:degrate} in Chapter \ref{chap8}.

\section{Dynkin's criterion and intertwining}
\label{sec:DynkinIntertwining}

Throughout this section: $(X(t),\,t\ge0)$ is a continuous-time Markov process on a state space $(S,\cS)$, $(P_t,t\ge0)$ is the family of transition kernels for $X$, $\phi\colon S\to T$ is a surjective measurable map to $(T,\cT)$, and $Y(t) = \phi(X(t))$, $t\ge0$. We will discuss two different sufficient criteria for $(Y(t),\,t\ge0)$ to also be Markovian.

\begin{theorem}[Theorem 10.13 of \cite{Dynkin}]\label{thm:Dynkin}
 Let $\phi^{-1}$ denote the pre-image under $\phi$. If $\phi$ satisfies \emph{Dynkin's criterion} that for all $A\in \cT$ and $(x,y)\in S^2\text{ with }\phi(x)=\phi(y)$,
 \begin{equation}\label{eq:Dynkin}
  P_t(x, \phi^{-1}(A)) = P_t(y,\phi^{-1}(A)),
 \end{equation}
 then $(Y(t),\,t\ge0)$ is a Markov process in the filtration generated by $(X(t),\,t\ge0)$.
\end{theorem}

See \cite[Lemma I.14.1]{RogersWilliams} for another version of this result. This is also sometimes credited as the Kemeny--Snell criterion, after \cite[Theorem 6.3.2]{KemenySnell}.

The second criterion that we discuss is stated in terms of compositions of stochastic kernels. We adopt the standard convention that sequential transitions are ordered from left to right, unlike the notation for compositions of functions:
$$\int PQ(x,dz)f(z) = \int P(x,dy)\int Q(y,dz)f(z).$$

\begin{definition}\label{def:intertwining}
 Consider a stochastic kernel $\Lambda\colon T\times\cS\to [0,1]$ and let $Q_t := \Lambda P_t\Phi$, $t\ge0$, where $\Phi$ denotes the kernel associated with the map $\phi$, $\Phi(x,\cdot\,) = \delta_{\phi(x)}(\,\cdot\,)$. We say $(Q_t,\,t\ge0)$ is \emph{intertwined below} $(P_t,\,t\ge0)$ via $\Lambda$ if
 \begin{enumerate}[label=(\roman*),ref=(\roman*)]
  \item \label{item:intertwin:up} $\Lambda\Phi$ equals the identity kernel on $(T,\cT)$ and
  \item \label{item:intertwining} $\Lambda P_t = Q_t\Lambda$, $t\ge0$.
 \end{enumerate}
\end{definition}

\begin{theorem}[Theorem 2 of \cite{RogersPitman}]\label{thm:RogersPitman}
 If $(Q_t,\,t\ge0)$ is intertwined below $(P_t,\,t\geq 0)$ via $\Lambda$ and additionally,
 \begin{enumerate}[start=3,label=(\roman*),ref=(\roman*)]
  \item \label{item:intertwin:init} $X(0)$ has regular conditional distribution (r.c.d.) $\Lambda(Y(0),\cdot\,)$ given $Y(0)$,
 \end{enumerate}
 then $(Y(t),\,t\ge0)$ is a Markov process. We then say that $(Y(t),\,t\ge0)$ is \emph{intertwined below} $(X(t),\,t\ge0)$ via $\Lambda$.
\end{theorem}

If conditions \ref{item:intertwin:up} and \ref{item:intertwin:init} are satisfied, then \ref{item:intertwining} is equivalent \cite[Remark (ii)]{RogersPitman} to
\begin{enumerate}[start=2, label=(\roman*'), ref=(\roman*')]
 \item \label{item:intertwining_v2} For all $t\ge 0$ and $y\in T$, if $X(0)$ has distribution $\Lambda(y,\cdot)$, then the r.c.d.\ of $X(t)$ given $Y(t)$ is $\Lambda(Y(t),\cdot\,)$.
\end{enumerate}

\section[A metric on marked $k$-trees and the proof of Lemma 6.6]{A metric on marked $k$-trees and the proof of Lemma \ref{lem:Lstar_cont}}\label{sec:marked_tree_metric}

In this section we discuss a metric on the space $\bTMarkk$ of marked $k$-trees that was introduced in Section \ref{sec:const:intertwin}. Recall from Definition \ref{def:mark} that a marked $k$-tree is an ordered pair $\Ast T = (T,\ell)$, where $T$ is a $k$-tree in the sense of Section \ref{sec:killed_def}  (i.e.\ with leaf masses and internal edge partitions) and $\ell$ is a distinguished block in $T$, either a leaf $i\in [k]$ or one of the blocks along one of the internal edge partitions. Such marked trees arise as the $\phi_2$-projections of $(k\!+\!1)$-trees; recall from Definition \ref{def:mark} and Figure \ref{fig:mark_k_tree} that this map contracts away leaf $k\!+\!1$ and leaves a marking in the block of the resulting $k$-tree where that leaf would need to be inserted in order to recover our initial $(k\!+\!1)$-tree.

We metrize $\bTMarkk$ by
\begin{equation}\label{eq:markk:dist_def}
 d_{\Ast\bT}\left(\Ast T_{k,1},\Ast T_{k,2}\right) := \inf\left\{ d_{\bT}(T_{k+1,1},T_{k+1,2})\colon \phi_2(T_{k+1,1}) = \Ast T_{k,1},\,\phi_2(T_{k+1,2}) = \Ast T_{k,2}\right\}.
\end{equation}

Note that if $T_{k,1}$ and $T_{k,2}$ have the same tree shape as each other and both are marked in corresponding leaf blocks $i\in [k]$, then
\begin{equation}\label{eq:markk:dist_1}
 d_{\Ast\bT}\left((T_{k,1},i),(T_{k,2},i)\right) = d_{\bT}(T_{k,1},T_{k,2}).
\end{equation}
Indeed, $d_{\bT}(T_{k,1},T_{k,2})$ is a general lower bound for distances between marked $k$-trees. In this special case, the bound can be seen to be sharp by splitting block $i$ in each of the marked $k$-trees into a very small block $k\!+\!1$, a small edge partition with little diversity, and a massive block $i$, in order to form $(k\!+\!1)$-trees that project down as desired. In the limit as block $k\!+\!1$ and the edge partition on its parent approach mass and diversity zero, the $d_{\bT}$-distance between the resulting $(k\!+\!1)$-trees converges to $d_{\bT}(T_{k,1},T_{k,2})$.

On the other hand, if two marked $k$-trees have equal tree shape but the marked blocks lie in different leaf components or internal edge partitions, then
\begin{equation}\label{eq:markk:dist_2}
 d_{\Ast\bT}\left((T_{k,1},\ell_1),(T_{k,2},\ell_2)\right) = d_{\bT}(T_{k,1},0) + d_{\bT}(0,T_{k,2}).
\end{equation}
If $T_{k,1}$ and $T_{k,2}$ have different tree shapes, then both \eqref{eq:markk:dist_1} and \eqref{eq:markk:dist_2} hold, as the right hand sides are then equal, by \eqref{eq:ktree:metric_2}.

This leaves only the case where the two marked $k$-trees have the same shape and the marked blocks each lie in corresponding internal edge partitions in the two trees. Then each marked $k$-tree is as in example (A) in Figure \ref{fig:mark_k_tree}: for $i=1,2$, there is a unique $(k\!+\!1)$-tree $T_{k+1,i}$ for which $\phi_2(T_{k+1,i}) = (T_{k,i},\ell_i)$. Then
\begin{equation}\label{eq:markk:dist_3}
 d_{\Ast\bT}\big((T_{k,1},(E,a_1,b_1)),\, (T_{k,2},(E,a_2,b_2))\big) = d_{\bT}(T_{k+1,1},T_{k+1,2}).
\end{equation}

\begin{proof}[Proof of Lemma \ref{lem:Lstar_cont}]
 First, we note that for $T_1,T_2\in\bTInt_{k+1}$ and $\Ast T_1,\Ast T_2\in\bTMarkk$,
 \begin{equation}
  d_{\Ast\bT}(\phi_1(T_1),\phi_1(T_2)) \le d_{\bT}(T_1,T_2) \quad \text{and} \quad d_{\bT}\Big(\phi_2\big(\Ast T_1\big),\phi_2\big(\Ast T_2\big)\Big) \le d_{\Ast T}\Big(\Ast T_1,\Ast T_2\Big).
 \end{equation}
 The first of these inequalities follows immediately from the definition of $d_{\Ast\bT}$. The second follows from \eqref{eq:markk:dist_1}, \eqref{eq:markk:dist_2}, and \eqref{eq:markk:dist_3}, with the added note that the projection map $\pi_k$ satisfies
 $$d_{\bT}(T_{k+1,1},T_{k+1,2}) \le d_{\bT}(\pi_k(T_{k+1,1}),\pi_k(T_{k+1,2})).$$
 This proves the continuity of $\phi_1$ and $\phi_2$.
 
 We now prove that $\Ast\Lambda_k$ is weakly continuous in its first coordinate. We separately check continuity at zero, at $k$-trees with a marked leaf, and at $k$-trees with the mark in a block of an interval partition. In each case, we consider a sequence $((T_n,\ell_n),n\ge1)$ of marked $k$-trees converging to a limit $\Ast T_\infty$ of that type.
 
 Case 1: $\Ast T_\infty = 0$. Then the total mass $\|T_n\|$ and the diversities of all interval partition components of the $T_n$ must go to zero. Let $U = (m_1,m_2,\beta)\sim Q$ denote a Brownian reduced 2-tree of unit mass. For $n\ge1$, let $\cT_n := T_n\oplus (\ell_n,U)$. Then $\cT_n$ has law $\Ast\Lambda_k((T_n,\ell_n),\cdot\,)$. We recall that, as noted in \cite[Equation (3.5)]{Paper1-0}, scaling an interval partition by $c$, causes its diversity to scale by $\sqrt{c}$. Thus,
 $$d_{\bT}(\cT_n,0) \le d_{\bT}(T_n,0) + \sqrt{\|\ell_n\|}\sD_{\beta}(\infty) \le d_{\bT}(T_n,0) + \sqrt{\|T_n\|}\sD_{\beta}(\infty),$$
 which goes to zero as $n$ tends to infinity. We conclude that $\Ast\Lambda_k$ is weakly continuous at 0.
 
 Case 2: $\Ast T_\infty = (T_\infty,i)$ for some $i\in [k]$. Then by \eqref{eq:markk:dist_2}, for all sufficiently large $n$, $\ell_n = i$ and $T_n$ has the same tree shape as $T_\infty$; call this tree shape $\ft$. Let $U$ and $(\cT_n,n\ge1)$ be as in Case 1. Let $x_{n,i}$ denote the mass of block $i$ in $T_n$. As noted in \cite[Equation (3.7)]{Paper1-0},
 $$d_{\bT}(T_n,T_m) \le d_{\bT}(\cT_n,\cT_m) \le d_{\bT}(T_n,T_m) + \big|\sqrt{x_{n,i}} - \sqrt{x_{m,i}}\big|\sD_{\beta}(\infty).$$
 Since the sequence $(x_{n,i},n\ge1)$ is Cauchy, it follows from the above bounds that $(\cT_n,n\ge1)$ is a.s.\ Cauchy as well. Thus, we conclude that $\Ast\Lambda_k(T_n,\cdot\,)$ converges weakly.
 
 Case 3: $\ell_n = (E,a_n,b_n)$ for some $E\in \ft$, for all sufficiently large $n$. Then for each such large $n$, there is some $(k\!+\!1)$-tree $T_{k+1,n}$ such that $\Ast\Lambda_k((T_n,\ell_n),\cdot\,) = \delta_{T_{k+1,n}}(\,\cdot\,)$. By \eqref{eq:markk:dist_3}, if the marked $k$-trees $(T_n,\ell_n)$ converge then so do the $(k+1)$-trees $T_{k+1,n}$.
 
 This proves that $\Ast\Lambda_k$ is weakly continuous in its first coordinate everywhere on $\bTMarkk$.
\end{proof}

\section[Proof of Lemma 6.18]{Proof of Lemma \ref{lem:degen_diff}}\label{sec:non_acc_2}

In this section, we prove Lemma \ref{lem:degen_diff}, which we restate here for convenience.

\begin{lemma}(Lemma \ref{lem:degen_diff})
 Fix $k\ge 3$ and $\epsilon>0$. Let $T\in\TInt_{k-1}$ with $\|T\|>\epsilon$ and let $(\cT^y,y\ge0)$ be a resampling $k$-tree evolution with $\cT^0\sim\Lambda_{k,[k-1]}(T,\cdot\,)$. Let $(D_n,n\ge1)$ denote the sequence of all degeneration times of this evolution and $(D^*_n,n\ge1)$ the subsequence of degeneration times at which label $k$ drops and resamples. Assume that with probability one we get $D_\infty:=\sup_n D_n > D^*_2$. Then there is some $\delta = \delta(k,\epsilon)>0$ that does not depend on $T$ such that $\bP(D^*_2>\delta)>\delta$.
\end{lemma}

We prove this lemma in three cases.
\begin{enumerate}[label = Case \arabic*:, ref = \arabic*, itemindent=1.82cm, leftmargin=0pt]
 \item\label{case:degdif:leaf} $T$ contains a leaf block of mass $x_i > \|T\|/2k$.
 \item\label{case:degdif:IP_big} $T$ contains an edge partition $\beta$ of mass at least $\|T\|/2k$, and $\beta$ contains a block of mass at least $\|\beta\|/2k^2$.
 \item\label{case:degdif:IP_small} $T$ contains an edge partition $\beta$ of mass at least $\|T\|/2k$, and each block in $\beta$ has mass less than $\|\beta\|/2k^2$.
\end{enumerate}

\begin{proof}[Proof of Lemma \ref{lem:degen_diff}, Case \ref{case:degdif:leaf}]
 With probability at least $1/2k$, the kernel $\Lambda_{k,[k-1]}(T,\cdot\,)$ inserts label $k$ into the large leaf block $i$, splitting it into a Brownian reduced 2-tree $(x_i,x_k,\beta_{\{i,k\}})$.  This type-2 compound $(\cU^y,y\in [0,D_1))$ will then evolve in pseudo-stationarity, as in Proposition \ref{prop:012:pseudo}, until the first degeneration time $D_1$ of $(\cT^y,y\ge0)$.
 
 Let $A_1$ denote the event that $\cU^{D_1-}$ is not degenerate, i.e.\ some other compound degenerates at time $D_1$. On $A_1$, some outside label may swap places with $i$ and cause $i$ to resample. However, as noted in the discussion of cases \ref{case:degen:type2}, \ref{case:degen:self}, and \ref{case:degen:nephew} in Section \ref{sec:const:intertwin}, no label will swap places with label $k$ at time $D_1$ on the event $A_1$. Let $\mathcal{R}_1$ denote the subtree of $\varrho(\cT^{D_1-})$ corresponding to $\cU^{D_1-}$. This equals $\cU^{D_1-}$ if no label swaps with $i$. By Proposition \ref{prop:012:pseudo} and exchangeability of labels in Brownian reduced 2-trees, on the event $A_1$ the tree $\cR_1$ is a Brownian reduced 2-tree.
 
 Let $B_1$ denote the event that the label dropped at $D_1$ resamples into a block of $\cR_1$.  Let $\cU^{D_1}$ denote the resulting subtree after resampling. On the event $A_1\cap B_1^c$, $\cU^{D_1} = \cR_1$ is again a Brownian reduced 2-tree, by definition of the resampling $k$-tree evolution. On the event $A_1\cap B_1$, the tree $\cU^{D_1}$ is a Brownian reduced 3-tree, by \eqref{eq:B_ktree_resamp} and the exchangeability of labels.
 
 We extend this construction inductively. Suppose that on the event $\bigcap_{m=1}^n A_m$, the tree $\cU^{D_n}$ is a Brownian reduced $M$-tree, for some (random) $M$. We define $(\cU^y,y\in [D_n,D_{n+1}))$ to be the $M$-tree evolution in this subtree during this time interval. Let $A_{n+1}$ denote the event that $\cU^{D_{n+1}-}$ is non-degenerate, $\cR_{n+1}$ the corresponding subtree in $\varrho(\cT^{D_{n+1}-})$, $B_{n+1}$ the event that the dropped label resamples into a block in $\cR_{n+1}$, and $\cU^{D_{n+1}}$ the corresponding subtree in $\cT^{D_{n+1}}$. Then on $\bigcap_{m=1}^{n+1} A_m$ the tree $\cR_{n+1}$ is again a Brownian reduced $M$-tree, by the same arguments as above, with Proposition \ref{prop:pseudo:pre_D} in place of Proposition \ref{prop:012:pseudo}. On $B_{n+1}^c\cap \bigcap_{m=1}^{n+1} A_m$, the tree $\cU^{D_{n+1}} = \cR_{n+1}$ is a Brownian reduced $M$-tree, and on $B_{n+1}\cap \bigcap_{m=1}^{n+1} A_m$ the tree $\cU^{D_{n+1}}$ is a Brownian reduced $(M\!+\!1)$-tree.
 
 In this manner, we define $(\cU^y,y\in [0,D_N))$ where $D_N$ is the first time that $\cU^{y-}$ attains a degenerate state as a left limit. Let $A^y$ denote the label set of $\cU^y$ for $y\in [0,D_N)$; by the preceding argument and Proposition \ref{prop:pseudo:pre_D}, $\cU^y$ is conditionally a Brownian reduced $(\# A^y)$-tree given $\{y < D_N\}$. Let $(\sigma^y,y\in [0,D_N))$ denote the evolving permutation that composes all label swaps due to the swap-and-reduce map, $\sigma^y = \tau_n\circ\tau_{n-1}\circ\cdots\circ\tau_1$ for $y\in [D_n,D_{n+1})$, where $\tau_m$ is the label swap permutation that occurs at time $D_m$. 
 
 We can simplify this account by considering a pseudo-stationary killed $k$-tree evolution $(\cV^y,y\in [0,D''))$ coupled so that $\pi_{A^y}\circ\sigma^y(\cV^y) = \cU^y$, $y\in [0,D'')$, where $D''$ is the degeneration time of $(\cV^y)$. Such a coupling is possible due to the consistency result of Proposition \ref{prop:consistency_0} and the exchangeability of labels evident in Definition \ref{def:killed_ktree} of killed $k$-tree evolutions. Note that, in particular, $D''$ precedes the first time at which a label in $(\cU^y)$ degenerates. Moreover, following the discussion of cases \ref{case:degen:type2}, \ref{case:degen:self}, and \ref{case:degen:nephew} in Section \ref{sec:const:intertwin}, label $k$ cannot be dropped in degeneration until a label within $(\cU^y)$ degenerates.
 
 There is some $\delta>0$ sufficiently small so that a pseudo-stationary $k$-tree evolution with initial mass $\epsilon/2k$ will avoid degenerating prior to time $\delta$ with probability at least $2k\delta$. By the self-similarity noted in Theorem \ref{thm:Markov}, this same $\delta$ bound holds for pseudo-stationary $k$-tree evolutions with greater initial mass. Applying this bound to $(\cV^y,y\in [0,D''))$ proves the lemma in this case.
\end{proof}

\begin{proof}[Proof of Lemma \ref{lem:degen_diff}, Case \ref{case:degdif:IP_big}]
 In this case, with probability at least $1/4k^3$, label $k$ is inserted into a ``large'' block in $\beta$ of mass at least $\|\beta\|/2k^2$. If another label resamples into this same block prior to time $D^*_1$, then we are in the regime of Case \ref{case:degdif:leaf}, and the same argument applies, albeit with smaller initial mass proportion. However, if no other label resamples into this block then, although it is unlikely for this block to vanish quickly, it is possible for label $k$ to be dropped in degeneration if a label that is a nephew of $k$ causes degeneration (case \ref{case:degen:nephew} in Section \ref{sec:const:intertwin}). In this latter case, however, that label swaps into the block in which label $k$ was sitting. Then, label $k$ resamples and may jump back into this large block with probability bounded away from zero. This, again, puts us in the regime of Case \ref{case:degdif:leaf}. In this case, $D^*_1$ may be small with high probability, but not $D^*_2$.
 
 More formally, a version of the argument for Case \ref{case:degdif:leaf} yields $\delta>0$ for which, with probability at least $\delta$: (i) the kernel $\Lambda_{k,[k-1]}(T,\cdot\,)$ inserts label $k$ into a block in $\beta$ with mass at least $\|\beta\|/2k^2$; (ii) this block, or a subtree created within this block survives to time $\delta$ with its mass staying above $\|\beta\|/3k^2$; (iii) the total mass stays below $2\|T\|$ and either (iv) $D^*_1 > \delta$; or (v) $D^*_{1} \le \delta$ but at time $D^*_{1}$, label $k$ resamples back into this same block, which only holds a single other label at that time; and then (vi) $D^*_{2} - D^*_{1} > \delta$.
\end{proof}

To prove Case \ref{case:degdif:IP_small}, we require two lemmas, one of which recalls additional properties of type-0/1/2 evolutions.

\begin{lemma}\label{lem:012:clade_ish}
 Fix $(x_1,x_2,\beta)\in [0,\infty)^2\times\cI$ with $x_1+x_2>0$. There exist a type-0 evolution $(\beta_0^y,y\ge0)$, a type-1 evolution $((m^y,\beta_1^y),y\ge0)$, and a type-2 evolution $((m_1^y,m_2^y,\beta_2^y),y\ge0)$ with respective initial states $\beta$, $(x_1,\beta)$, and $(x_1,x_2,\beta)$, coupled in such a way that for every $y$, there exists an injective, left-to-right order-preserving and mass-preserving map sending the blocks of $\beta_2^y$ to blocks of $\beta_1^y$, and a map with these same properties sending the blocks of $\beta_1^y$ to blocks of $\beta_0^y$.
\end{lemma}

These assertions are immediate from the pathwise constructions of type-0/1/2 evolutions in Constructions \ref{type0:construction}, \ref{type1:construction} and \ref{deftype2}. 

\begin{lemma}\label{lem:type2:smallblocks}
 Fix $c\in (0,1/2)$ and $x>0$. Consider $u_1,u_2\ge0$ with $u_1+u_2>0$ and $\beta\in\cI$ with $\|\beta\|>x$ and none of the blocks of $\beta$ having mass greater than $c\|\beta\|$. For every $\epsilon>0$ there exists some $\delta = \delta(x,c) > 0$ that does not depend on $(u_1,u_2,\beta)$ such that with probability at least $1-\epsilon$, a type-2 evolution with initial state $(u_1,u_2,\beta)$ avoids degenerating prior to time $\delta$.
\end{lemma}

\begin{proof}
 Fix a block $(a,b)\in\beta$ with $a\in [c\|\beta\|,2c\|\beta\|]$ and let
 \begin{equation*}
  \beta_0 := \{(a',b')\in\beta\colon a'<a\},\quad \beta_1 := \{(a'-b,b'-b)\colon (a',b')\in\beta, a'\ge b\}
 \end{equation*}
 so that $\beta = \beta_0\concat (0,b-a)\concat\beta_1$. We follow Proposition \ref{prop:012:concat}\ref{item:012concat:2+1}, in which a type-2 evolution is formed by concatenating a type-2 with a type-1. In particular, let $\wh\Gamma^y := \big( \wh m_1^y, \wh m_2^y,\wh \beta^y\big)$ and $\widetilde\Gamma^y := (\widetilde m^y,\widetilde\beta^y)$, $y\ge0$, denote a type-2 and a type-1 evolution with respective initial states $(u_1,u_2,\beta_0)$ and $(b-a,\beta_1)$. Let $\wh D$ denote the degeneration time of $(\wh\Gamma^y,y\ge0)$ and let $\wh Z$ denote the time at which $\|\wh\Gamma^y\|$ hits zero. Let $I$ equal 1 if $\wh m_1^{\wh D}>0$ or 2 if $\wh m_2^{\wh D}>0$, and set $(X_I,X_{3-I}) := \big(\wh m_I^{\wh D},\widetilde m^{\wh D}\big)$. Finally, let $\big(\widebar m_1^y,\widebar m_2^y,\widebar\beta^y\big)$, $y\ge 0$ denote a type-2 evolution with initial state $(X_1,X_2,\widetilde\beta^{\widehat D})$, conditionally independent of $((\wh\Gamma^y,\widetilde\Gamma^y)),y\in [0,\widehat D]))$ given this initial state, but coupled to have $\widebar m_I^y = \wh m_I^{\widehat D + y}$ for $y\in [0,\widehat Z-\widehat D]$. By Proposition \ref{prop:012:concat}\ref{item:012concat:2+1}, the following is a type-2 evolution:
 \begin{equation}\label{eq:SB:concat_1}
  \left\{\begin{array}{ll}
   (\wh m_1^y,\;\wh m_2^y,\;\wh\beta^y\concat(0,\widetilde m^y)\concat\widetilde\beta^y)	& \text{for }y\in [0,\wh D),\\
   (\widebar m_1^{y-\wh D},\;\widebar m_2^{y-\wh D},\;\widebar\beta^{y-\wh D})	& \text{for }y\ge \wh D.
  \end{array}\right.
 \end{equation}
 Moreover, by the Markov property of type-1 evolutions, Definition \ref{def:type2:v1} of type-2 evolutions, and the symmetry noted in Lemma \ref{lem:type2_symm}, the following is a stopped type-1 evolution:
 \begin{equation}\label{eq:SB:concat_2}
  \left\{\begin{array}{ll}
   (\widetilde m^y,\widetilde\beta^y)	& \text{for }y\in [0,\wh D),\\
   (\widebar m_{3-I}^{y-\wh D},\widebar\beta^{y-\wh D})	& \text{for }y\in [\wh D,\wh Z].
  \end{array}\right.
 \end{equation}
 
 Let $\delta>0$ be sufficiently small so that, with probability at least $\sqrt{1-\epsilon}$, a $\besq_{cx}(-1)$ avoids hitting zero prior to time $\delta$, and likewise for a $\besq_{(1-2c)x}(0)$. Then with probability at least $1 - \epsilon = (\sqrt{1-\epsilon})^2$, both $(\|\wh\Gamma^y\|,y\ge0)$ and the $\besq(0)$ total mass of the type-1 evolution of \eqref{eq:SB:concat_2} avoid hitting zero prior to time $\delta$. On this event, the type-2 evolution of \eqref{eq:SB:concat_1} does not degenerate prior to time $\delta$.
\end{proof}

\begin{proof}[Proof of Lemma \ref{lem:degen_diff}, Case \ref{case:degdif:IP_small}]
  Informally, we proved that degenerations of $k$ may take a long time in Cases \ref{case:degdif:leaf} and \ref{case:degdif:IP_big} by controlling 
  the degeneration times of pseudo-stationary structures inserted into large blocks in repeated resampling events. In Case 
  \ref{case:degdif:IP_small}, there are no large blocks, and indeed large blocks may never form. Instead, there must be a large interval partition,
  which we can cut rather evenly into $2k-1$ sub-partitions. We will control the degeneration of evolving sub-partitions and the probability
  that insertions are into distinct non-adjacent internal sub-partitions.
  
  Specifically, let us follow the notation introduced at the start of this appendix. We can decompose $\beta=\beta_1\star\cdots\star\beta_{2k-1}$
  into sub-partitions $\beta_i$ with $\|\beta_1\|\ge \|\beta\|/k$ and $\|\beta\|/4k\le\|\beta_i\|\le\|\beta\|/2k$ for $2\le i\le 2k-1$, since no 
  block exceeds mass $\|\beta\|/2k^2\le\|\beta\|/4k$. 
  
  With probability at least $1/8k^2$, the kernel $\Lambda_{k,[k-1]}(T,\cdot\,)$ inserts label $k$ into $\beta_2$, splitting 
  $\beta_2=\beta_2^-\star(0,x_k)\star\beta_2^+$. Then label $k$ is in the type-1 compound 
  $\cU^0=(x_k,\beta_2^+\star\beta_3\star\ldots\star\beta_{2k-1})$, while $\beta_1\star\beta_2^-$ is the interval partition of the 
  compound associated with the sibling edge of $k$. Consider the concatenation $\cV^y:=\cV^y_3\star\cV^y_4\star\cdots\star\cV^y_{2k-1}$
  of type-1 evolutions $(\cV^y_i,y\ge 0)$, $3\le i\le 2k-1$, starting respectively from $\beta_2^+\star\beta_3,\beta_4,\ldots,\beta_{2k-1}$, for
  times $y$ up to the first time $D_\cV$ that one of them reaches half or double its initial mass. We denote by $D_\cW$ the degeneration time of
  the sibling edge of $k$ as part of this resampling $k$-tree evolution. If this sibling edge of $k$ is a type-2 edge, 
  denote by $u_1$ and $u_2$ its top masses and consider a type-2 evolution $(\cW^y,y\ge 0)$ starting from $(u_1,u_2,\beta_1\star\beta_2^-)$, with
  degeneration time $D_\cW$. 
  Otherwise, this edge has three or more labels, so one or both children of this edge have more than one label. For each of these children, we 
  choose as $u_1$ or $u_2$, respectively, the top mass of its smallest label. Then the degeneration time of a type-2 evolution 
  $(\widetilde{\cW}^y,y\ge 0)$ starting from $(u_1,u_2,\beta_1\star\beta_2^-)$ is stochastically dominated by the time $D_\cW$ at which the sibling edge of $k$ degenerates. We denote by $D_\cT$ the first time that $\|\cT^y\|$ reaches half or double its initial mass.
  
  Since $\|\beta_1\star\beta_2^-\|\ge\|\beta_1\|\ge\|\beta\|/k\ge\|T\|/2k^2\ge\epsilon/2k^2$, Lemma \ref{lem:type2:smallblocks} yields 
  $\delta_\cW=\delta(\epsilon/2k^2,1/2k)>0$ such that $\bP(D_\cW>\delta_1)\ge 1-1/(32k^2)^k$. 
  Let $\delta_\cT>0$ be such that a $\besq_{\epsilon}(-1)$ stays in $(\epsilon/2,2\epsilon)$ up to time $\delta_\cT$ with probability at least
  $1-1/(32k^2)^k$. Then $\bP(D_\cT>\delta_\cT)\ge 1-1/(32k^2)^k$. Finally, let $\delta_\cV>0$ be such that the probability that 
  $\besq_{\epsilon/8k^2}(0)$ does not exit $(\epsilon/16k^2,\epsilon/4k^2)$ before time $\delta_\cV$ exceeds $(1/2)^{1/2k}$. Since the $2k-3$ 
  independent type-1 evolutions are starting from greater initial mass, we obtain from Proposition \ref{prop:012:mass} and the self-similarity assertion of Theorem \ref{thm:Markov} that $\bP(D_\cV>\delta_\cV)>1/2$. 
  
  We proceed in a way similar to Case 1 and inductively construct a subtree evolution $(\cU^y,y\in[0,D_\cU))$ coupled to $(\cV^y,y\in[0,D_\cV))$ on 
  events $A_{n+1}$, $n\ge 0$, on which $D_{n+1}<\min\{D_\cV,D_\cW,D_\cT\}$ and any resampling of a label at time $D_{n+1}$ when $\cU^{D_{n+1}-}$ 
  has $j-1$ labels occurs into a block of $\cV^{D_{n+1}}_{2j}$, $j=2,\ldots,k-1$. Given that $D_{n+1}<\min\{D_\cV,D_\cW,D_\cT\}$, such a block is 
  chosen by the resampling kernel with (conditional) probability exceeding $\|\beta_{2j}\|/(4\|T\|)\ge 1/32k^2$. 
  
  Thus, with probability at least $\delta:=\min\{\delta_\cV,\delta_\cW,\delta_\cT,1/(32k^2)^{k-2}\}$, we have $D_1^*>\delta$ and so $D_2^*>\delta$.
\end{proof}

\section{Representation of $\mathbb{R}$-trees by consistent families of $k$-trees}\label{sec:contproj}

Recall from Section \ref{sec:Rtrees} 
\begin{itemize}
  \item the subset $\widebar{\mathbb{T}}^{\rm int}_\infty$ of projectively consistent families $\mathrm{R}=(R_k,k\ge 1)$ of $k$-trees 
    $R_k=\big(\ft_k,(x_j^{(k)},j\in[k]),(\beta_E^{(k)},E\in\ft_k)\big)\in\widebar{\mathbb{T}}_k^{\rm int}$, $k\ge 1$, equipped with the subset topology of the product topology formed by the metric topologies induced by 
    \eqref{eq:ktree:metric_1}--\eqref{eq:ktree:metric_2},
  \item and the metric space $\big(\mathbb{T}^{\rm real}_\infty,d_{\rm GHP}^\infty\big)$ of ${\rm GHP}_\infty$-isometry classes of rooted, weighted compact $\mathbb{R}$-trees with a 
    sequence of marked points. 
\end{itemize}
In this section, we make explicit the kernel $\mathbf{m}_\infty$ from $\mathbb{T}^{\rm real}$ to $\mathbb{T}^{\rm real}_\infty$ that samples a sequence of marked points from the normalized weight measure, and the map $R\colon\mathbb{T}_\infty^{\rm real}\rightarrow\widebar{\mathbb{T}}^{\rm int}_\infty$ that uses the marked points of
the $\mathbb{R}$-tree (representative) to build a corresponding projectively consistent family in $\widebar{\mathbb{T}}^{\rm int}_\infty$, and we prove Theorem \ref{thm:realtoint}. 

Let $\big(\bT^{\rm real}_k,d_{\rm GHP}^{[k]}\big)$ be Miermont's \cite{M09} space of ${\rm GHP}^{[k]}$-isometry classes of rooted, weighted compact $\mathbb{R}$-trees with $k$ marked points.  In this section, we adapt Miermont's notation and write $\mathrm{T}=\big[T,d,\mu,(\sigma_i)_{i=0}^k\big]=\big[T,d,\mu,(\sigma_0,\boldsymbol{\sigma})\big]\in\bT^{\rm real}_k$ for the ${\rm GHP}^{[k]}$-isometry class of $(T,d,\mu,(\sigma_i)_{i=0}^k)$, a rooted, weighted compact $\mathbb{R}$-tree $(T,d,\sigma_0,\mu)$ with $k$ marked points $\sigma_1,\ldots,\sigma_k\in T$. 
We extend this notation to similarly write $\mathrm{T}=\big[T,d,\mu,(\sigma_i)_{i=0}^\infty\big]=\big[T,d,\mu,(\sigma_0,\boldsymbol{\sigma})\big]\in\mathbb{T}^{\rm real}_\infty$.

\begin{definition} For any $1\le k\le \infty$ and any rooted, weighted compact $\mathbb{R}$-tree $(T,d,\sigma_0,\mu)$ with $\mu\neq 0$, let
  \[ 
  \mathbf{m}_k\big((T,d,\sigma_0,\mu),\mathbf{A}\big) 
  = \int_{T^k} \big(\mu/\|\mu\|\big)^{\otimes k}(d\boldsymbol{\sigma}) \mathbf{1}_{\mathbf{A}}\Big( \big[T,d,\mu,(\sigma_0,\boldsymbol{\sigma})\big]\Big).
  \]
\end{definition}

As noted in Section \ref{sec:Rtrees}, Miermont's arguments for finite $k$ extend to $k=\infty$. Specifically, $\mathbf{m}_k((T,d,\sigma_0,\mu),\cdot)$ only depends on the
${\rm GHP}$-isometry class $\big[T,d,\sigma_0,\mu\big]$ and induces a kernel from $\mathbb{T}^{\rm real}$ to $\mathbb{T}^{\rm real}_k$, for all $1\le k\le \infty$. See
\cite[Lemma 13]{M09}.

Denote by $L(T,d)$ the set of leaves of the $\mathbb{R}$-tree $(T,d)$, i.e.\ the set of members of $T$ whose removal does not disconnect $T$. For $1\le k<\infty$, let
\[ 
\bT^{{\rm real},L}_k = \Big\{\big[T,d,\mu, (\sigma_i)_{i=0}^k\big] \in \bT_k^{\rm real}\colon L(T,d) \subseteq \{\sigma_0,\sigma_1,\dots, \sigma_k\} \Big\}
\]
be the set of ${\rm GHP}^{[k]}$-isometry classes of trees that are spanned by the root and their other $k$ marked points. 

For any $k$-marked rooted, weighted compact $\mathbb{R}$-tree $\big(T,d,\mu,(\sigma_i)_{i=0}^k\big)$ and $0\le i\neq j\le k$, we consider the isometries
$\rho_{i,j}^T\colon[0,d(\sigma_i,\sigma_j)]\rightarrow T$ with $\rho_{i,j}^T(0)=\sigma_i$ and $\rho_{i,j}^T(d(\sigma_i,\sigma_j))=\sigma_j$ and the associated paths 
$[\![\sigma_i,\sigma_j]\!]_T=\rho_{i,j}^T([0,d(\sigma_i,\sigma_j)])$ in $T$. We write $]\!]\sigma_i,\sigma_j[\![_T=[\![\sigma_i,\sigma_j]\!]_T\setminus\{\sigma_i,\sigma_j\}$ and define
$[\![\sigma_i,\sigma_j[\![_T$ and $]\!]\sigma_i,\sigma_j]\!]_T$ similarly. Define the subset
\begin{equation}\label{eq:reduce}
T_k^+=\bigcup_{1\le j\le k}[\![\sigma_0,\sigma_j]\!]_T
\end{equation}
of $T$ spanned by $\sigma_0,\ldots,\sigma_k$ and the projection map $\pi_k^+\colon T\rightarrow T_k^+$ given by $\pi_k^+(\sigma)=\arg\min_{x\in T_k^+}d(x,\sigma)$. We note that this is 
well-defined because $(T,d)$ is an $\mathbb{R}$-tree and $T_k^+\subseteq T$ is a closed connected subset. We consider $T_k^+$ as a $k$-marked rooted $\mathbb{R}$-tree by equipping it with the metric and marked points inherited from $T$, and
we further equip $T_k^+$ with the measure $\mu_k^+=(\pi_k^+)_*\mu$, the image of $\mu$ by $\pi_k^+$. Since the ${\rm GHP}^{[k]}$-isometry class $\mathrm{T}_k^+$ of 
$\big(T_k^+,d,\mu_k^+,(\sigma_i)_{i=0}^k\big)$ only depends on the ${\rm GHP}^{[k]}$-isometry class $\mathrm{T}$ of $\big(T,d,\mu,(\sigma_i)_{i=0}^k\big)$, this induces a map
${\textsc{reduce}}_k^+\colon\mathbb{T}^{\rm real}_k\rightarrow\mathbb{T}^{{\rm real},L}_k$ given by $\textsc{reduce}_k^+(\mathrm{T})=\mathrm{T}_k^+$  . 

\begin{lemma}\label{lm:kprojection}
The map ${\textsc{reduce}}_k^+\colon\mathbb{T}^{\rm real}_k\rightarrow\mathbb{T}^{{\rm real},L}_k$ is Lipschitz continuous with Lipschitz constant less than or equal to 17. 
\end{lemma}
\begin{proof} For the purposes of this proof, we omit the superscript $+$ from reduced trees and projection maps. Consider a metric space $(M,\delta)$ and two embedded $k$-marked rooted, weighted compact $\mathbb{R}$-trees $\big(T,\delta,\mu,(\sigma_i)_{i=0}^k\big)$
  and $\big(T^\prime,\delta,\mu^\prime,(\sigma_i^\prime)_{i=0}^k\big)$ with 
  \[
  \delta^{\mathrm{H}}(T,T^\prime)<\varepsilon,\qquad \delta^{\rm P}(\mu,\mu^\prime)<\varepsilon,\qquad \delta(\sigma_i,\sigma_i^\prime)<\varepsilon,\quad 0\le i\le k,
  \]
  and their reduced subtrees $T_k$ and $T_k^\prime$ defined as in \eqref{eq:reduce}, with associated projection maps $\pi_k\colon T\rightarrow T_k$ and
  $\pi_k^\prime\colon T^\prime\rightarrow T_k^\prime$ and measures $\mu_k$ and $\mu_k^\prime$. It suffices to show that 
  $\delta^{\mathrm{H}}(T_k,T_k^\prime)<17\varepsilon$ and $\delta^{\mathrm{P}}(\mu_k,\mu_k^\prime)<17\varepsilon$.

  \medskip

  Let $\sigma\in T_k$. Then $\sigma\in[\![\sigma_0,\sigma_i]\!]_T$ for some $i\in[k]$. Since $\delta(\sigma_0,\sigma_0^\prime)<\varepsilon$ and 
  $\delta(\sigma_i,\sigma_i^\prime)<\varepsilon$, the triangular inequality yields $|\delta(\sigma_0,\sigma_i)-\delta(\sigma_0^\prime,\sigma_i^\prime)|<2\varepsilon$. Therefore,
  we can find $\sigma^\prime\in[\![\sigma_0^\prime,\sigma_i^\prime]\!]_{T^\prime}$ such that
  \[
  \big|\delta(\sigma_0^\prime,\sigma^\prime)-\delta(\sigma_0,\sigma)\big|<2\varepsilon\quad\mbox{and}\quad
  \big|\delta(\sigma_i^\prime,\sigma^\prime)-\delta(\sigma_i,\sigma)\big|<2\varepsilon.
  \]
  A priori, $\sigma^\prime$ may be far from $\sigma$ in $M$. However, since $\delta^{\rm H}(T,T^\prime)<\varepsilon$, we can also find $\widetilde{\sigma}^\prime\in T^\prime$ with
  $\delta(\sigma,\widetilde{\sigma}^\prime)<\varepsilon$, and since $T^\prime$ is an $\mathbb{R}$-tree, we have  
  $\sigma^\prime\in[\![\sigma_i^\prime,\widetilde{\sigma}^\prime]\!]_{T^\prime}$ or $\sigma^\prime\in[\![\sigma_0^\prime,\widetilde{\sigma}^\prime]\!]_{T^\prime}$. In the first case,
  \begin{align*}
  \delta(\sigma,\sigma^\prime)&\le\delta(\sigma,\widetilde{\sigma}^\prime)+\delta(\widetilde{\sigma}^\prime,\sigma^\prime)
                                          =\delta(\sigma,\widetilde{\sigma}^\prime)+\delta(\widetilde{\sigma}^\prime,\sigma_i^\prime)-\delta(\sigma^\prime,\sigma_i^\prime)\\
                                       &\le\varepsilon+\delta(\widetilde{\sigma}^\prime,\sigma)+\delta(\sigma,\sigma_i)+\delta(\sigma_i,\sigma_i^\prime)-\delta(\sigma,\sigma_i)+2\varepsilon
	\ <\  5\varepsilon.
  \end{align*}
  In the second case, the same argument with $i$ replaced by $0$ yields the same conclusion. Reversing the roles of $T_k$ and $T_k^\prime$, we conclude that 
  $\delta^{\rm H}(T_k,T_k^\prime)<5\varepsilon$. 

\medskip

  Turning to the measures, consider any closed $C\subseteq M$ and recall notation $C^\varepsilon$ for its $\varepsilon$-thickening. Then $\pi_k^{-1}(C\cap T_k)\subseteq M$ is also  closed and so
  \[
  \mu_k(C)=\mu(\pi_k^{-1}(C\cap T_k))\le\mu^\prime\big((\pi_k^{-1}(C\cap T_k))^\varepsilon\big)+\varepsilon.
  \]
  Our aim is to show that this can be further bounded by $\mu_k^\prime(C^{17\varepsilon})+17\varepsilon$. To this end, let 
  $\sigma^\prime\in\big(\pi_k^{-1}(C\cap T_k)\big)^\varepsilon\cap T^\prime$. Then there is $\sigma\in\pi_k^{-1}(C\cap T_k)$ such that 
  $\delta(\sigma,\sigma^\prime)\le\varepsilon$. Now consider the projections $\pi_k(\sigma)$ into $T_k$ and $\pi_k^\prime(\sigma^\prime)$ into $T_k^\prime$. 
  Since $\delta^{\rm H}(T_k,T_k^\prime)<5\varepsilon$, there is $p\in T_k$ such that $\delta(\pi_k^\prime(\sigma^\prime),p)<5\varepsilon$. Then
  \[
  \delta(\sigma,\pi_k(\sigma))\le\delta(\sigma,p)
                                            \le \delta(\sigma,\sigma^\prime)+\delta(\sigma^\prime,\pi_k^\prime(\sigma^\prime))+\delta(\pi_k^\prime(\sigma^\prime),p)
                                            \le \delta(\sigma^\prime,\pi_k^\prime(\sigma^\prime))+6\varepsilon.
  \]
  Similarly, $\delta(\sigma^\prime,\pi_k^\prime(\sigma^\prime))\le\delta(\sigma,\pi_k(\sigma))+6\varepsilon$. Since $\pi_k(\sigma)\in[\![\sigma,p]\!]_T$, we find
  \begin{align*}
  \delta(\pi_k^\prime(\sigma^\prime),\pi_k(\sigma))
               &\le\delta(\pi_k^\prime(\sigma^\prime),p)+\delta(p,\pi_k(\sigma))
               =\delta(\pi_k^\prime(\sigma^\prime),p)+\delta(\sigma,p)-\delta(\sigma,\pi_k(\sigma))\\
               &\le 5\varepsilon+\delta(\sigma^\prime,\pi_k^\prime(\sigma^\prime))+6\varepsilon-\delta(\sigma^\prime,\pi_k^\prime(\sigma^\prime))+6\varepsilon
               \ =\ 17\varepsilon.
  \end{align*}
  We conclude that $\big(\pi_k^{-1}(C\cap T_k)\big)^\varepsilon\cap T^\prime\subseteq\big(\pi_k^\prime\big)^{-1}\big(C^{17\varepsilon}\cap T_k^\prime\big)$, and this entails that
  \[
  \mu^\prime\big((\pi_k^{-1}(C\cap T_k))^\varepsilon\big)+\varepsilon\le\mu^\prime\big(\big(\pi_k^\prime\big)^{-1}\big(C^{17\varepsilon}\cap T_k^\prime\big)\big)+17\varepsilon,
  \]
  as required. 
\end{proof}

For $\mathrm{T}=\big[T,d,\mu,(\sigma_i)_{i=0}^k\big]\in\mathbb{T}_k^{{\rm real},L}$, we construct the tree shape in $\mathbb{T}_{k}^{\rm shape}$, denoted by
$\textsc{shape}_k(\mathrm{T})$, via its graph-theoretic representation $(V,E)$, by the following procedure. The vertex set is 
$V=\{\sigma_i,1\le i\le k\}\cup\{b_{i,j},1\le i<j\le k\}\subseteq T$, where $b_{i,j}$ is the unique point in $T$ that satisfies $[\![\sigma_0,b_{i,j}]\!]_T=[\![\sigma_0,\sigma_i]\!]_T\cap[\![\sigma_0,\sigma_j]\!]_T$. The edge set $E$ includes an edge between $u\in V$ and $v\in V$ if and only if $]\!]u,v[\![_T\cap V=\emptyset$. With each $v\in V$, 
associate $\textsc{labels}(v)=\{i\in[k]\colon\sigma_i=v\}\cup\bigcup_{1\le i<j\le k\colon b_{i,j}=v}\{i,j\}$. Since $\big\{\textsc{labels}(v),v\in V\big\}$ does not depend on 
the choice of representative of the isometry class $\mathrm{T}$, we can define the tree shape of $\mathrm{T}$ as
\[
\textsc{shape}_k(\mathrm{T})
  =\left\{
      \begin{array}{ll}
         \big\{\textsc{labels}(v),v\in V\big\}\setminus\big\{\{j\}\colon j\in[k]\big\}&\mbox{if this is in }\mathbb{T}^{\rm shape}_{k},\\[0.2cm]
         \Delta&\mbox{otherwise,}
      \end{array}
    \right.
\] 
where $\Delta\not\in\mathbb{T}_{k}^{\rm shape}$ is a cemetery state. If $B=\textsc{labels}(v)$ and $v=b_{i,j}$, we also write $v=b_B$. 
In fact, it is not hard to see that while $\overline{\mathbf{t}}:=\big\{\textsc{labels}(v),v\in V\big\}$ is not necessarily a binary hierarchy in the sense of Section \ref{sec:killed_def}, it is always a hierarchy in the generalised sense where each 
$B\in\overline{\mathbf{t}}$ with $\#B\ge 2$ has a unique partition into a minimal number of parts $C_1,\ldots,C_k\in\overline{\mathbf{t}}$, for some $k\ge 2$, without further requiring that $k=2$.
For the purposes of this section, it suffices to consider $\textsc{shape}_k(\mathrm{T})$ when this hierarchy is binary and hence gives rise to a (binary!) tree shape in $\mathbb{T}^{\rm shape}_{k}$. 

\begin{proposition}
  The function $\textsc{shape}_k\colon\mathbb{T}^{{\rm real},L}_k\rightarrow\mathbb{T}_{k}^{\rm shape}\cup\{\Delta\}$ is measurable.
\end{proposition}
\begin{proof} Since $\mathbb{T}_{k}^{\rm shape}\cup\{\Delta\}$ is a finite set, it suffices to show that for all $\ft\in\bT_{k}^{\rm shape}$, the set 
  $\textsc{shape}_k^{-1}(\ft)$ is measurable. We will argue that this set is open in $(\mathbb{T}^{{\rm real},L},d_{\rm GHP}^{[k]})$. To this end, consider any 
  $\mathrm{T}=\big[T,d,\mu,(\sigma_i)_{i=0}^k\big]\in\mathbb{T}^{{\rm real},L}_k$ with tree shape $\textsc{shape}_k(\mathrm{T})=\ft$. The construction of its tree shape via 
  the graph-theoretic representation $(V,E)$ is such that any two points in the finite vertex set $V\subseteq T$ are at a strictly positive $d$-distance in $T$. From this and \eqref{eq:bpdist}, 
  we can find a  $d_{\rm GHP}^{[k]}$-ball in $\mathbb{T}_{k}^{{\rm real},L}$ centered at $\mathrm{T}$ in which the tree shape remains constant, as 
  required.
\end{proof}

\begin{remark} As we are applying \eqref{eq:bpdist} from the proof of Proposition \ref{propGH} in what we will use to prove Theorem \ref{thm:realtoint}, we point out that Proposition \ref{propGH}   
  does not depend on Theorem \ref{thm:realtoint} nor indeed on any material from Section \ref{sec:defmarkov}.
\end{remark}

For $\mathrm{T}=\big[T,d,\mu,(\sigma_0,\sigma_1,\sigma_2)\big]\in\mathbb{T}_2^{{\rm real},L}$, recall $[\![\sigma_0,b_{1,2}]\!]_T=[\![\sigma_0,\sigma_1]\!]_T\cap[\![\sigma_0,\sigma_2]\!]_T$. The (first) top mass $\mu\big(]\!]b_{1,2},\sigma_1]\!]_T\big)$ does not depend on the choice of representative of the isometry class $\mathrm{T}$ and will be denoted by $\textsc{top}(\mathrm{T})$.

More generally, let $\mathrm{T}=\big[T,d,\mu,(\sigma_i)_{i=0}^k\big]\in\mathbb{T}_k^{{\rm real},L}$ and $j\in[k]$. Then we define the $j^{\rm th}$ top mass
\[
\textsc{top}_j^{(k)}(\mathrm{T})=\mu\left(\bigcap_{i\in[k]\setminus\{j\}}]\!]b_{i,j},\sigma_j]\!]_T\right).
\]
\begin{proposition}
  The functions $\textsc{top}_j^{(k)}\colon\mathbb{T}_k^{{\rm real},L}\rightarrow[0,\infty)$ are measurable for all $1\le j\le k$.
\end{proposition}
\begin{proof} First consider the function $\textsc{top}\!=\!\textsc{top}_1^{(2)}$. Consider a sequence $\mathrm{T}_n\!\rightarrow\!\mathrm{T}=[T,d,\mu,(\sigma_0,\sigma_1,\sigma_2)]$ in the space 
  $(\mathbb{T}^{{\rm real},L}_2,d_{\rm GHP}^{[2]})$. Then the distances between the branch point and each of the marked points converge (cf. \eqref{eq:bpdist}) 
  and since $]\!]b_{1,2},\sigma_1]\!]_T$ 
  is open in $T$, we have $\liminf_{n\rightarrow\infty}\textsc{top}(\mathrm{T}_n)\ge\textsc{top}(\mathrm{T})$, i.e.\ $\textsc{top}$ is lower semi-continuous and thus 
  measurable. The argument is easily adapted to handle $\textsc{top}_j^{(k)}$, $1\le j\le k$.
\end{proof}

For any $\mathrm{T}=\big[T,d,\mu,(\sigma_0,\sigma_1,\sigma_2)\big]\in\mathbb{T}_2^{{\rm real},L}$, we consider the ``restriction'' $\big([\![\sigma_0,b_{1,2}]\!]_T,d,\mu|_{[\![\sigma_0,b_{1,2}]\!]_T},(\sigma_0,b_{1,2})\big)$. This space is ${\rm GHP}^{[1]}$-isometric to the interval $[0,d(\sigma_0,b_{1,2})]$ equipped with the root $0$, the marked point $d(\sigma_0,b_{1,2})$, and the image weight measure. We denote the weighted interval by $\textsc{edge}(\mathrm{T})\in\mathcal{M}$, using the space introduced in 
\eqref{eq:onebranch}. 

More generally, let $\mathrm{T}=\big[T,d,\mu,(\sigma_i)_{i=0}^k\big]\in\mathbb{T}_k^{{\rm real},L}$ and $B\in\textsc{shape}_k(\mathrm{T})\setminus\{[k]\}$, we consider
the restriction $\big([\![b_{\parent{B}},b_B]\!]_T,d,\mu|_{]\!]b_{\parent{B}},b_B]\!]_T},(b_{\parent{B}},b_B)\big)$. This space is ${\rm GHP}^{[1]}$-isometric to the interval
$[0,d(b_{\parent{B}},b_B)]$ equipped with the root $0$, the marked point $d(b_{\parent{B}},b_B)$, and the image weight measure.  
We similarly handle the case $B=[k]$ using the convention that $b_{\parent{[k]}}=\sigma_0$, and restricting $\mu$ to the left-closed geodesic $[\![\sigma_0,b_{[k]}]\!]_T$ to include a potential atom at the root. In all cases, we denote the weighted interval by $\textsc{edge}_B^{(k)}(\mathrm{T})\in\mathcal{M}$.

\begin{lemma}\label{lm:edge}
  For each $\ft\in\mathbb{T}_{k}^{\rm shape}$, the functions $\textsc{edge}_B^{(k)}\colon\textsc{shape}^{-1}_k(\ft)\rightarrow\mathcal{M}$, $B\in\ft$, are measurable.
\end{lemma}
\begin{proof} First consider the function $\textsc{edge}=\textsc{edge}_{[2]}^{(2)}$. Let $\mathrm{T}_n\rightarrow\mathrm{T}$ as in the proof of the previous lemma. 
  Here, the closed interval that forms the first component in $\textsc{edge}(\mathrm{T}_n)$ converges in the Hausdorff sense, as $n\rightarrow\infty$, while the weight 
  measure component in $\textsc{edge}(\mathrm{T}_n)$ when evaluated on closed subsets of the closed limiting interval exhibits upper semi-continuity, as $n\rightarrow\infty$. This 
  entails (vague convergence when restricted to the interior of the limiting closed interval and hence) measurability of function $\textsc{edge}$. Again this argument can be adapted
  to handle $\textsc{edge}_B^{(k)}$, $B\in\ft$, $\ft\in\mathbb{T}_{k}^{\rm shape}$.
\end{proof}

The following results are straightforward consequences of the definitions.

\begin{lemma} The map $\mu\mapsto\big(\mu([0,t]),t\ge 0\big)$ from the space of finite measures on $[0,\infty)$ to the Skorokhod space $\mathbb{D}^\uparrow_b([0,\infty),[0,\infty))$ of bounded increasing functions is
  measurable. The map $\mu\mapsto\mu([0,\infty))$ is continuous, and the map $(f,\sup f)\mapsto\textsc{range}(f)=\overline{f([0,\infty))}$ from  $\{(g,M)\in\mathbb{D}_b^\uparrow([0,\infty),[0,\infty))\colon\sup g=M\}$ to the space of compact subsets of $[0,\infty)$ equipped with the Hausdorff metric is also continuous.
\end{lemma}

\begin{corollary} The map $\textsc{ip}(\mu)=\textsc{range}\big(\mu([0,t]),t\ge 0\big)$ is measurable.
\end{corollary}

We will slightly abuse notation and consider $\textsc{ip}\colon\cM\rightarrow\overline{\cI}_H$ as a map on $\cM$ that only depends on the measure component in $\cM$ and that takes values in the space $\overline{\cI}_H$ of all interval partitions in a sense generalizing Definition \ref{def:intro:IP}, where we associate with any compact subset $C\subseteq[0,\infty)$ the set of bounded connected components of $[0,\infty)\setminus C$ and equip $\overline{\cI}_H$ with a metric induced by the Hausdorff metric. See \cite{Paper1-0} for a fuller discussion of this induced topology on the most relevant subspace $\cI_H\subset\overline{\cI}_H$ of interval partitions $\beta$ whose partition points $C=[0,\|\beta\|]\setminus\bigcup_{U\in\beta}U$ have zero Lebesgue measure, and also \cite[Section 4]{IPPAT} for the relevance of generalized interval partitions for interval partition evolutions. In the setting of Lemma \ref{lm:edge}, we now introduce maps $\textsc{ip}_B^{(k)}:=\textsc{ip}\circ\textsc{edge}_B^{(k)}$.   

\begin{proposition} For all $\ft\in\mathbb{T}_{k}^{\rm shape}$, the functions $\textsc{ip}_B^{(k)}\colon\textsc{shape}^{-1}_k(\ft)\rightarrow\overline{\cI}_H$, $B\in\ft$, are measurable.
\end{proposition}

Now recall that the ultimate aim of this section is the proof of Theorem \ref{thm:realtoint}, which we can now phrase more precisely, as follows. The claim is we can construct a function 
$R\colon\mathbb{T}_\infty^{\rm real}\rightarrow\widebar{\mathbb{T}}_\infty^{\rm int}$ such that  
$S(R([T,d,\mu,(\sigma_i)_{i=0}^\infty]))=[T,d,\sigma_0,\mu]$
for $\mathbf{m}_\infty([T,d,\sigma_0,\mu],d[T,d,\mu,(\sigma_i)_{i=0}^\infty]){\tt BCRT}(d[T,d,\sigma_0,\mu])$-a.e. $[T,d,\mu,(\sigma_i)_{i=0}^\infty]\in\mathbb{T}_\infty^{\rm real}$, where ${\tt BCRT}$ denotes the distribution on $\mathbb{T}^{\rm real}$ of a Brownian CRT. We will slightly abuse notation and abbreviate this as saying 
\[ 
S(R(\mathrm{T},(\sigma_i)_{i=1}^\infty))=\mathrm{T}\quad\mbox{for }{\tt BCRT}_\infty\mbox{-a.e. }(\mathrm{T},(\sigma_i)_{i=1}^\infty)\in\mathbb{T}_\infty^{\rm real}.
\]

In the following we will further abuse notation and consider $\textsc{reduce}_k^+\colon\mathbb{T}^{\rm real}_\infty\rightarrow\mathbb{T}_k^{{\rm real},L}$, naturally defined
by projecting away the redundant marks beyond the first $k$ via the natural 1-Lipschitz map from $\mathbb{T}^{\rm real}_\infty$ to $\mathbb{T}^{\rm real}_k$.

\begin{definition}\label{def:mapR} Let $\mathrm{T}=[T,d,\mu,(\sigma_i)_{i=0}^\infty]\in\mathbb{T}_\infty^{\rm real}$. Then we define 
$R(\mathrm{T})=(R_k(\textsc{reduce}_k^+(\mathrm{T})),k\ge 1)$, if this is in $\widebar{\mathbb{T}}^{\rm int}_\infty$, where
  \[
  R_k(\mathrm{T}_k^+)=\Big(\textsc{shape}_k(\mathrm{T}_k^+),\big(\textsc{top}_j^{(k)}(\mathrm{T}_k^+),j\in[k]\big),\big(\textsc{ip}_B^{(k)}(\mathrm{T}_k^+),B\in\textsc{shape}_k(\mathrm{T}_k^+)\big)\Big)
  \]
  if $\textsc{shape}_k(\mathrm{T}_k^+)\in\mathbb{T}_{k}^{\rm shape}$. If $\textsc{shape}_k(\mathrm{T}_k^+)=\Delta$ or if $\textsc{ip}_B^{(k)}(\mathrm{T}_k^+)\not\in\cI$ for any $B\in\textsc{shape}_k(\mathrm{T}_k^+)$, we set 
  $R_k(\mathrm{T}_k^+)=0\in\widebar{\mathbb{T}}_k^{\rm int}$. If $(R_k(\textsc{reduce}_k^+(\mathrm{T})),k\ge 1)$ is not consistent, we define $R(\mathrm{T})=0\in\widebar{\mathbb{T}}_\infty^{\rm int}$.
\end{definition}

\begin{proof}[Proof of Theorem \ref{thm:realtoint}] First note that a Brownian CRT is almost surely binary with a diffuse weight measure supported by the leaves \cite{AldousCRT1}. For any 
  rooted, weighted compact $\mathbb{R}$-tree $(T,d,\rho,\mu)$ with these properties, points sampled from $\mu$ will be distinct almost surely, and any finite number of distinct 
  leaves $\sigma_1,\ldots,\sigma_k$ gives rise to a binary tree shape. Hence, $\textsc{shape}_k(\textsc{reduce}_k^+([T,d,\mu,(\sigma_i)_{i=0}^k]))$ is well-defined as an element 
  of $\mathbb{T}_{k}^{\rm shape}$, for each $k$, and given a sequence $(\sigma_i)_{i=1}^\infty$ of distinct leaves, also $R([T,d,\mu,(\sigma_i)_{i=0}^\infty])$ is well-defined 
  as a member of $\prod_{k\ge 1}\widebar{\mathbb{T}}_k^{\rm int}$, since Brownian reduced $k$-trees almost surely have interval partitions that are members of $\cI$. And also
  the consistency is a consequence of the construction, as noted previously when constructing Brownian reduced $k$-trees in Section \ref{sec:intro:B_k_tree}.

  We now turn to the application of $S$ to $R([T,d,\mu,(\sigma_i)_{i=0}^\infty])$. In Definition \ref{def:S}, we defined $S(R_k,k\ge 1)=\lim_{k\rightarrow\infty}\tau(R_k)$, where
  $\tau(R_k)$ is the ${\rm GHP}$-isometry class of the weighted $\mathbb{R}$-tree $S_k(R_k)$ constructed in \eqref{eq:skcirc}--\eqref{eq:sk} branch by branch using diversities
  of interval partitions as locations and block sizes as sizes of atoms. 

  On the other hand, 
  $\lim_{k\rightarrow\infty}\textsc{reduce}_k^+(\mathrm{T},(\sigma_i)_{i=1}^\infty)$ exists in $(\mathbb{T}^{\rm real},d_{\rm GHP})$ for any $(\mathrm{T},(\sigma_i)_{i=1}^\infty)\in\mathbb{T}^{\rm real}_\infty$, since representatives of reduced trees are naturally embedded in 
  any representative of $\mathrm{T}$ and form an increasing sequence of closed sets in a compact metric space that converges to the closure of their union in the Hausdorff sense. As the weight 
  measure of the Brownian CRT has dense support \cite[Theorem 3]{AldousCRT1}, this Hausdorff limit is actually $\mathrm{T}$ for ${\tt BCRT}_\infty$-a.e.\ $(\mathrm{T},(\sigma_i)_{i=1}^\infty)\in\mathbb{T}_\infty^{\rm real}$. As the weight measures of the reduced trees are just projections of the weight measure of the limiting tree, the Hausdorff convergence further entails the Hausdorff--Prokhorov convergence of embedded weighted $\mathbb{R}$-trees.

  As a slight variation of the above argument, we can consider the function $\textsc{reduce}_k\colon\mathbb{T}^{\rm real}_\infty\rightarrow \mathbb{T}^{{\rm real},L}_k$ that
  associates with $\big[T,d,\mu,(\sigma_i)_{i=0}^\infty\big]$ the ${\rm GHP}^{[k]}$-isometry class $\mathrm{T}_k$ of 
  \[
  T_k=\bigcup_{1\le i<j\le k}[\![\sigma_0,b_{i,j}]\!]_T
  \]
  of $T$ with root $\sigma_0$ inherited from $\big[T,d,\mu,(\sigma_i)_{i=0}^\infty\big]$ and now equipped with the marked points $\pi_k(\sigma_i)$, $i\in[k]$, and the weight measure $(\pi_k)_*\mu$, both
  projected by the natural projection $\pi_k\colon T\rightarrow T_k$. Then similarly $\lim_{k\rightarrow\infty}\textsc{reduce}_k(\mathrm{T},(\sigma_i)_{i=1}^\infty)$ exists in $(\mathbb{T}^{\rm real},d_{\rm GHP})$ and equals $\mathrm{T}$
  for ${\tt BCRT}_\infty$-a.e.\ $(\mathrm{T},(\sigma_i)_{i=1}^\infty)\in\mathbb{T}_\infty^{\rm real}$.

  To show $S(R(\mathrm{T},(\sigma_i)_{i=1}^\infty))=\mathrm{T}$, it thus suffices to show
  that $\textsc{reduce}_k(\mathrm{T},(\sigma_i)_{i=1}^\infty)$ is the ${\rm GHP}^{[k]}$-isometry class of 
  $S_k(R_k(\textsc{reduce}_k^+(\mathrm{T},(\sigma_i)_{i=1}^\infty)))$, where 
  $R_k\colon\mathbb{T}_k^{{\rm real},L}\rightarrow\widebar{\mathbb{T}}_k^{\rm int}$ is as in Definition \ref{def:mapR}.

  Specifically, for $k=2$ and any binary $(\mathrm{T},(\sigma_i)_{i=1}^\infty)\in\mathbb{T}_\infty^{\rm real}$ , write 
  $\mathrm{T}_2^+=\textsc{reduce}_2^+(\mathrm{T},(\sigma_i)_{i=1}^\infty)$ and 
  suppose that $\textsc{edge}(\mathrm{T}_2^+)\in\cI$. Then the weighted one-branch $\mathbb{R}$-tree $S_2(R_2(\mathrm{T}_2^+))$ of length 
  $\sD(\textsc{edge}(T_2^+))$ has an atom at the end whose size is a sum that includes $\textsc{top}_1^{(2)}(\mathrm{T}_2^+)$ and 
  $\textsc{top}_2^{(2)}(\mathrm{T}_2^+)$, which is precisely the mass projected to the end of 
  $\mathrm{T}_2=\textsc{reduce}_2(\mathrm{T},(\sigma_i)_{i=1}^\infty)$. Apart from top masses, $S_2(R_2(\mathrm{T}_2^+))$ has further atoms built from 
  sizes and locations (not necessarily distinct and possibly including the ends, in general). By construction, these are the block sizes and associated diversities of 
  $\textsc{edge}(\mathrm{T}_2^+)$, hence replicating the order of the atoms in (any representative of) $\mathrm{T}_2$. The subtle point is that the atom 
  locations on $S_2(R_2(\mathrm{T}_2^+))$ and $\mathrm{T}_2$ coincide for ${\tt BCRT}_\infty$-a.e.\ 
  $(\mathrm{T},(\sigma_i)_{i=1}^\infty)\in\mathbb{T}_\infty^{\rm real}$. This is a property established in \cite{PitmWink09}, as discussed in the context of 
  Proposition \ref{prop:PDIP}.  

  For $k\ge 3$, using analogous notation, we note that $S_k(R_k(\mathrm{T}_k^+))$ and $\mathrm{T}_k$ have the same tree shape, by construction. To 
  complete their identification, we apply the argument for the case $k=2$ to projections onto subtrees spanned by any two labels $1\le i<j\le k$. An induction 
  beginning with the edge to the branch point adjacent to the root (in the graph-theoretic tree shape) and proceeding to adjacent branch points in subtrees can be 
  used to complete the proof. We leave the details to the reader.
\end{proof}

\section[Proof of Proposition 7.18]{Proof of Proposition \ref{propGHP}}\label{app:pfpropGHP}

Recall that Proposition \ref{propGHP} claims that the GHP-distance between the weighted $\bR$-trees $S_k(R_k)$ and $S_k(R_k^\prime)$ associated with two trees $R_k,R_k^\prime\in\widebar{\bT}_{[k]}^{\rm int}$ with the same shape $\ft_k$ can be bounded above, as follows,
  $$d_{\rm GHP}(S_k(R_k),S_k(R_k^\prime))\le 3k\max_{1\le i<j\le k}\min\Big\{d_\cI(\pi_{i,j} R_k,\pi_{i,j} R_k^\prime),d_\cI(\pi_{j,i}R_k,\pi_{j,i}R_k^\prime)\Big\},$$
where for each $1\le i<j\le k$ the interval partitions $\pi_{i,j}R_k,\pi_{j,i}R_k,\pi_{i,j}R_k^\prime,\pi_{j,i}R_k^\prime\in\cI$ were defined at the beginning of Section \ref{sec:dItoGHP} to capture interval partition representations of projected 2-trees that have the top masses as left-most intervals in the order indicated by the indices $i$ and $j$.

Before we begin the proof, let us recall from \cite[Proposition 6]{M09} Miermont's representation of $d_{\rm GHP}$ that extends the similar representation \eqref{eq:GH2} of $d_{\rm GH}$ as infimum of ${\rm GH}$-distortions of ${\rm GH}$-correspondences.
For two unit-mass weighted $\mathbb{R}$-trees $\mathrm{T}=(T,d,\rho,\mu)$ and 
$\mathrm{T}^\prime=(T^\prime,d^\prime,\rho^\prime,\mu^\prime)$, we consider pairs $(K,\nu)$, where $K\subseteq T\times T^\prime$ is a ${\rm GH}$-correspondence and $\nu$ a coupling of
$\mu$ and $\mu^\prime$, i.e.\ a probability measure on $T\times T^\prime$ whose marginal distributions are $\mu$ and $\mu^\prime$. Then 
\begin{equation}\label{eq:ghpunit2}
d_{\rm GHP}(\mathrm{T},\mathrm{T}^\prime)=\inf\left\{{\rm dis}_{\rm GHP}(K,\nu)\colon\begin{array}{c}\text{$K\subseteq T\times T^\prime$ ${\rm GH}$-correspondence}\\
                                                                                                                                                                              \text{$\nu$ coupling of $\mu$ and $\mu^\prime$}
                                                                                                                                                 \end{array}\right\},
\end{equation}
where  ${\rm dis}_{\rm GHP}(K,\nu):=\min\{{\rm dis}_{\rm GH}(K),1-\nu(K)\}$ is the ${\rm GHP}$-distortion of $(K,\nu)$. We will handle weighted $\mathbb{R}$-trees that are not necessarily of unit mass and to allow masses of $\mu$ and $\mu^\prime$ to differ, will involve a partial coupling that leaves mass outside the correspondence unallocated rather than assigned to $(T\times T^\prime)\setminus K$. 
Instead of establishing a general representation now, we will indicate at the end of the proof how the partial coupling can be used to bound $d_{\rm GHP}(S_k(R_k),S_k(R_k^\prime))$. 

In this context, let us also explain some of the main ideas. Firstly, in the case $k=2$, a bound 
$d_\cI(\pi_{1,2} R_k,\pi_{1,2} R_k^\prime)<\varepsilon$ means there is a correspondence $(U_j,V_j)_{j\in[n]}$ from 
$\beta:=\pi_{1,2}R_k$ to $\gamma:=\pi_{1,2}R_k^\prime$ of distortion less than $\varepsilon$, in the sense of Definition \ref{def:IP:metric}. In 
$S_2(R_2)$ and $S_2(R_2^\prime)$, which we can represent as measures on intervals of the form $[0,\sD_\beta(\infty)]$ and $[0,\sD_\gamma(\infty)]$, each pair 
$(U_j,V_j)$ gives rise to a pair ${\rm Leb}(U_j)\delta(\sD_\beta(U_j))$ and ${\rm Leb}(V_j)\delta(\sD_\gamma(V_j))$ of atoms. The natural associated partial 
coupling is $\nu=\sum_{j\in[n]}\min\{{\rm Leb}(U_j),{\rm Leb}(V_j)\}\delta(\sD_\beta(U_j),\sD_\gamma(V_j))$. By the definition of the distortion of 
$(U_j,V_j)_{j\in[n]}$, this leaves mass at most $\varepsilon$ unmatched, and all the coupled mass is in the subset 
$K\subseteq[0,\sD_\beta(\infty)]\times[0,\sD_\gamma(\infty)]$ of points within $\varepsilon$ of the diagonal, which can help build a ${\rm GH}$-correspondence of ${\rm dis}_{\rm GH}(K)=\varepsilon$. 

Secondly, in the cases $k\ge 3$, we will have to build such a partial coupling consistently from the various (overlapping!) interval partitions $\pi_{i,j}R_k$ and 
$\pi_{i,j}R_k^\prime$, $i,j\in[k]$, $i\neq j$. The challenge is that some of the blocks of $\pi_{i,j} R_k$ and $\pi_{i,j} R_k^\prime$ correspond to subtree 
masses rather than atoms of $S_k(R_k)$ and $S_k(R_k^\prime)$, and the $d_\cI$-correspondence of blocks in the definition of 
$d_\cI(\pi_{i,j} R_k,\pi_{i,j} R_k^\prime)$ does not take into account such ``internal structure'' of the block, so a block of $\pi_{i,j}R_k$ that gives rise to a 
single atom of $S_k(R_k)$ may correspond to a block of $\pi_{i,j}R_k^\prime$ whose mass in $S_k( R_k^\prime)$ is spread over a subtree. 

Here is a lemma that makes some elementary observations about correspondences of small distortion.
\begin{lemma}\label{lm:match} Let $\varepsilon>0$ and $\beta,\gamma\in\cI$ with $d_{\cI}(\beta,\gamma)<\varepsilon$. Let $(U_i,V_i)_{i\in[n]}$ and 
  $(U_j^\prime,V_j^\prime)_{j\in[m]}$ be two correspondences from $\beta$ to $\gamma$ with distortion at most $\varepsilon$. Then for every $U\in\beta$ with 
  ${\rm Leb}(U)>\varepsilon$, there are $i\in[n]$ and $j\in[m]$ such that $U=U_i=U_j^\prime$. If furthermore ${\rm Leb}(U)>2\varepsilon$, then $V_i=V_j^\prime=:V$ with 
  ${\rm Leb}(V)>\varepsilon$, and $\beta_{<U}:=\{W\in\beta\colon W<U\}$ and $\gamma_{<V}:=\{W\in\gamma\colon W<V\}$ satisfy 
  $d_{\cI}(\beta_{<U},\gamma_{<V})<\varepsilon$.  
\end{lemma}

%

\begin{proof}[Proof of Proposition \ref{propGHP}] As in the proof of Proposition \ref{propGH}, we proceed by setting up a ${\rm GH}$-correspondence, but now including pairs of big atoms as well as endpoints of branches among the pairs of special vertices, to construct a partial coupling of the weight measures, as well as a ${\rm GH}$-correspondence between the $\bR$-trees. 
  Now suppose that 
  \begin{equation}\label{eq:appxdIbound}
    \max_{1\le i<j\le k}\min\Big\{d_\cI(\pi_{i,j} R_k,\pi_{i,j} R_k^\prime),d_\cI(\pi_{j,i} R_k,\pi_{j,i}R_k^\prime)\Big\}<\varepsilon.
  \end{equation}
  For the purposes of this proof, we will call a $d_\cI$-correspondence 
  of blocks of $\pi_{i,j} R_k$ and $\pi_{i,j} R_k^\prime$ in the sense of the definition of $d_\cI$ an $(i,j)$-matching.  For each $1\le i<j\le k$, consider an $(i,j)$-matching, or a $(j,i)$-matching, for 
  which the $d_\cI$-distortion is less than $\varepsilon$. This exists by the definition of $d_\cI$. Recall that the right end of $\pi_{i,j}R_k$ corresponds to the root of $S_k(R_k)$,
  for each $1\le i,j\le k$, $i\neq j$. 
  
  In this setting, we prove by strong induction on $k$ that the $(i,j)$-matchings, $1\le i<j\le k$, induce a finite collection of corresponding pairs of atoms of $S_k(R_k)$ and 
  $S_k(R_k^\prime)$, which leaves at most mass $(3k-1)\varepsilon$ unmatched, where we call unmatched mass of $S_k( R_k)$ the full mass of any atom of $S_k(R_k)$ that 
  is not included in a pairing and also the residual mass of any atom of $S_k(R_k)$ that has been paired with an atom of $S_k(R_k^\prime)$ of smaller mass. More precisely, we
  show a corresponding claim for $R_A,R_A^\prime\in\widebar{\mathbb{T}}_A^{\rm int}$, where $S_A(R_A)$ and $S_A^\circ(R_A)$ etc.\ are defined as in Definition 
  \ref{def:SkandSkcirc}, which applies verbatim with $k$ and $[k]$ replaced by $A$.

\medskip

  {\it Induction hypothesis: for all $A\subseteq\bN$, $2\le\#A\le k-1$ and trees $R_A,R_A^\prime\in\widebar{\bT}_A^{\rm int}$ with the same tree shape $\ft\in\bT^{\rm shape}_A$, we 
  have the following. If there is an $(i,j)$-matching or a $(j,i)$-matching, $M_{i,j}$, for each $i,j\in A$, $i<j$, whose $d_\cI$-distortion is less than $\varepsilon$, 
  then there is a finite collection of pairs of atoms corresponding either to matched blocks $(X,X^\prime)$, or to $(X,Y^\prime)$ where $(X,X^\prime)$, $(Y,X^\prime)$ and $(Y,Y^\prime)$ are all matched blocks (for different $i,j$), that leaves at most mass $(3\#A-1)\varepsilon$ unmatched.}

\medskip

  For $A=\{b,c\}$ and $R_A,R_A^\prime\in\widebar{\bT}_A^{\rm int}$, there is only one pair $(i,j)=(b,c)$. By possibly swapping roles of $b$ and $c$, we may suppose 
  we have a $(b,c)$-matching, i.e.\ a $d_\cI$-correspondence of 
  blocks $U_1,\ldots,U_n$ of $\pi_{b,c} R_A$ and $U_1^\prime,\ldots,U_n^\prime$ of $\pi_{b,c} R_A^\prime$, with $d_\cI$-distortion less than $\varepsilon$. 
  We use this $d_\cI$-correspondence to pair associated atoms of 
  $S_A(R_A)$ and $S_A(R_A^\prime)$. By definition of $d_\cI$, we have 
  $$\sum_{r=1}^n\left|{\rm Leb}(U_r)-{\rm Leb}(U_r^\prime)\right|+||\pi_{b,c} R_A||-\sum_{r=1}^n{\rm Leb}(U_r)<\varepsilon<(3\#A-1)\varepsilon$$
  and
  $$\sum_{r=1}^n\left|{\rm Leb}(U_r)-{\rm Leb}(U_r^\prime)\right|+||\pi_{b,c} R_A^\prime||-\sum_{r=1}^n{\rm Leb}(U_r^\prime)
  <\varepsilon<(3\#A-1)\varepsilon.$$
  These quantities are precisely unmatched masses of $S_A(R_A)$ and $S_A(R_A^\prime)$ associated with this pairing.
  
  For $k=\#A\ge 3$ and $R_A,R_A^\prime\in\widebar{\bT}_A^{\rm int}$ with the same shape $\ft_A$, consider the first branch point 
  of $\ft_A$, which splits subtree labels $A=A(1)=B(1)\cup C(1)$, say. The choice of this branch point is preliminary and we will adjust this in some cases to a 
  different branch point that splits subtree labels $A(m)=B(m)\cup C(m)$ for some $m\ge 2$ in a way that we explain later. For the purposes of this proof, we distinguish five types of branch point according 
  to the sizes $\#B(m)$ and $\#C(m)$ of subtree label sets and corresponding subtree masses. For the latter, it will be convenient to consider $b(m)=\min B(m)$ and $c(m)=\min C(m)$ so that 
  $d_{\cI}(\pi_{b(m),c(m)}R_A,\pi_{b(m),c(m)}R_A^\prime)<\varepsilon$ and the subtree masses of $B(m)$, respectively $C(m)$, are the first block sizes $f_1(m)$ and $f_1^\prime(m)$, 
  respectively second block sizes $f_2(m)$ and $f_2^\prime(m)$, of $\pi_{b(m),c(m)}R_A$ and $\pi_{b(m),c(m)}R_A^\prime$. Here we use the convention that $f_2(m)=0$ if $\pi_{b(m),c(m)}R_A$ does
  not have a second block, similarly for $f_2^\prime(m)$. 
  \begin{itemize}
   \item Type I: $\#B(m)=1$ or $\#C(m)=1$, but not both. 
   \item Type II: $\#B(m)\ge 2$, $\#C(m)\ge 2$, and (i) or (ii) or both hold, where\\ (i) $f_1(m),f_1^\prime(m)\in(0,3\varepsilon]$ with at least one in $(0,2\varepsilon]$,\\
      (ii) $f_2(m),f_2^\prime(m)\in[0,3\varepsilon]$ with at least one in $[0,2\varepsilon]$.
   \item Type III: $\#B(m)=\#C(m)=1$. 
   \item Type IV: $\#B(m)\ge 2$, $\#C(m)\ge 2$, $f_1(m),f_1^\prime(m),f_2(m),f_2^\prime(m)\in(2\varepsilon,\infty)$.
    \item Type V: $\#B(m)\ge 2$, $\#C(m)\ge 2$ and one of (i)--(iv) hold, where\\ 
           (i) either $f_1(m)\in(0,2\varepsilon]$ and $f_2^\prime(m)\in[0,2\varepsilon]$ and $f_1^\prime(m),f_2(m)\in(3\varepsilon,\infty)$, 
                or $f_1^\prime(m)\in(0,2\varepsilon]$ and $f_2(m)\in[0,2\varepsilon]$ and $f_1(m),f_2^\prime(m)\in(3\varepsilon,\infty)$,\\
	(ii) one of $f_1(m),f_1^\prime(m)$ is in $(0,2\varepsilon]$, the other in $(3\varepsilon,\infty)$, while
      $f_2(m),f_2^\prime(m)$ are in $(2\varepsilon,\infty)$ with at least one in $(3\varepsilon,\infty)$,\\  
           (iii) one of $f_2(m),f_2^\prime(m)$ is in $[0,2\varepsilon]$, the other in $(3\varepsilon,\infty)$, 
      while $f_1(m),f_1^\prime(m)$ are in $(2\varepsilon,\infty)$ with at least one in $(3\varepsilon,\infty)$,\\
           (iv) either $f_1(m)\in(0,2\varepsilon]$ and $f_2(m)\in[0,2\varepsilon]$ and $f_1^\prime(m),f_2^\prime(m)\in(3\varepsilon,\infty)$, 
                or $f_1^\prime(m)\in(0,2\varepsilon]$ and $f_2^\prime(m)\in[0,2\varepsilon]$ and $f_1(m),f_2(m)\in(3\varepsilon,\infty)$.
  \end{itemize}
  For $A(m)$ of types I or II, we will say that the subtree
  labeled $B(m)$ is \em larger \em than the subtree labeled $C(m)$ if the branch point $A(m)$ is of type I with $\#C(m)=1$ or of type II(ii). 
  If we are now given $A(m)=B(m)\cup C(m)$ with $\#A(m)\ge 3$, we let $A(m+1)=B(m)$ if the subtree labeled $B(m)$ is larger than the subtree labeled 
  $C(m)$, and we let $A(m+1)=C(m)$ otherwise. Let $m_0$ be the first $m\ge 1$ for which the branch point $A(m)=B(m)\cup C(m)$ has type III, IV or V. We 
  will now identify pairs of atoms by considering a $(b(m_0),c(m_0))$-matching of distortion less than $\varepsilon$.\medskip

  First suppose that $A(m_0)=B(m_0)\cup C(m_0)$ has type III. In this case, the top masses of $\pi_{b(m_0),c(m_0)}R_A$ and $\pi_{b(m_0),c(m_0)}R_A^\prime$ have atom locations in $S_A(R_A)$ and 
  $S_A(R_A^\prime)$, but there may be other blocks in the matching that correspond to subtrees in $S_A(R_A)$ and/or $S_A(R_A^\prime)$. These are necessarily associated with the respective smaller subtree labeled $B(m)$ or $C(m)$ for $1\le m<m_0$.
  These $m_0-1$ subtrees arise from branch points of types I and II. We will also refer to these blocks as blocks of types I and II. Let us denote by $m_1$, $m_2$ the numbers of blocks of types
  I, II, respectively. Then $m_1+m_2=m_0-1$, while $m_1+2m_2\le\#A-2$ since the smaller subtree has one label for type I and at least two labels for type II. 
  \begin{itemize}
    \item Type-I blocks have atom locations in $S_A(R_A)$ and $S_A(R_A^\prime)$. Indeed, they are the top masses of type-1 edges in $R_A$ and $R_A^\prime$.
    \item Type-II blocks typically do not have atom locations in $S_A(R_A)$ and $S_A(R_A^\prime)$. Indeed, they are total masses of a type-2 edge (if they correspond to subtrees with precisely two labels) 
      or of several edges (if there are three or more labels). But they have sizes at most $3\varepsilon$. We remove from the $(b(m_0),c(m_0))$-matching all type-II blocks and all blocks matched with type-II 
      blocks. This increases the distortion by at most $6\varepsilon m_4$: in each pair, one type-II block has mass at most $2\varepsilon$, their match (if any) has mass at most 
      $3\varepsilon$, while the other has mass at most $3\varepsilon$ and their match (if any) at most $4\varepsilon$. Note that of these, the ones with bounds $2\varepsilon$ and 
      $4\varepsilon$ are in $\pi_{b(m_0),c(m_0)}R_A$, and the two with bound $3\varepsilon$ are in $\pi_{b(m_0),c(m_0)}R_A^\prime$, or vice versa, in each case summing to
      at most $6\varepsilon$ per type-II branch point.
   \end{itemize}
   The resulting matching has distortion at most 
   $\varepsilon+6\varepsilon m_2\le \varepsilon+3(\#A-2)\varepsilon\le (3\#A-1)\varepsilon$. Also, all blocks that remain in the matching have atom locations. \medskip
  
   Next suppose that $A(m_0)=B(m_0)\cup C(m_0)$ has type IV. Then $\pi_{b(m_0),c(m_0)}R_A$ and $\pi_{b(m_0),c(m_0)}R_A^\prime$ each have two top masses exceeding 
   $2\varepsilon$. By Lemma \ref{lm:match}, the $(i,j)$-matchings $M_{i,j}$, $i,j\in B(m_0)$, have matched blocks corresponding to the top masses of 
   $\pi_{b(m_0),c(m_0)}R_A$ and $\pi_{b(m_0),c(m_0)}R_A^\prime$ labeled $c(m_0)$ and induce $(i,j)$-matchings $M_{i,j}^{B(m_0)}$ of $\pi_{i,j}R_{B(m_0)}$ and
   $\pi_{i,j}R_{B(m_0)}^\prime$ of distortion at most $\varepsilon$. By the induction hypothesis, there is an associated finite collection of pairs of blocks that 
   leave at most mass $(3\#B(m_0)-1)\varepsilon$ unmatched. This holds similarly for $C(m)$ leaving mass at
       most $(3\#C(m_0)-1)\varepsilon$ unmatched. Finally, the non-top-mass part of the $(b(m_0),c(m_0))$-matching can be handled like the
   type-III case, here leaving at most $\varepsilon+3(\#A-\#A(m_0))\varepsilon$ unmatched. This adds to $(3\#A-1)\varepsilon$, as required.\smallskip

   Now suppose that $A(m_0)=B(m_0)\cup C(m_0)$ has type V. We argue in each case that for each label set it suffices to remove any matched top masses in $[0,2\varepsilon]$ and their match, and to locate 
   an atom in the corresponding subtree for the one in $(3\varepsilon,\infty)$, or to apply the induction hypothesis if for one label set, both top masses are in $(2\varepsilon,\infty)$. We discuss the four cases 
   (i)--(iv) separately.
   \begin{enumerate}
      \item[(i)] Suppose $f_1(m_0)\in(0,2\varepsilon]$, $f_2^\prime(m_0)\in[0,2\varepsilon]$, $f_1^\prime(m_0),f_2(m_0)\in(3\varepsilon,\infty)$. The other subcase then follows by symmetry. Then the blocks $U^\prime$ labeled 
         $B(m_0)$ and $U$ labeled $C(m_0)$ of top masses $f_1^\prime(m_0)$ and $f_2(m_0)$ are the first pair in the $(b(m_0),c(m_0))$-matching. The top mass of size $f_1(m_0)$ must be unmatched due to 
         the order constraints of $d_\cI$-correspondences. The top block $V^\prime$ of size $f_2^\prime(m_0)\le 2\varepsilon$ may be part of the $(b(m_0),c(m_0))$ matching, and if so, is matched to a block $V$ of size at 
         most $3\varepsilon$. 

         If $f_2(m_0)\le 7\varepsilon$, then $f_1^\prime(m_0)\le 8\varepsilon$. Removing the one or two matched pairs involving top masses from the $(b(m_0),c(m_0))$-matching increases its distortion by
         at most $\max\{7\varepsilon+3\varepsilon,8\varepsilon+2\varepsilon\}=10\varepsilon$. Proceeding as for type III further increases the distortion by at most 
        $6\varepsilon m_2\le 3(\#A-\#A(m_0))\varepsilon\le 3(\#A-4)\varepsilon$ summing to a total distortion of at most $(3\#A-1)\varepsilon$.

         If $f_2(m_0)>7\varepsilon$, consider any $i,j\in B(m_0)$. Since the block $U$ is also a block in $\pi_{i,j}R_A$, it is matched in the $(i,j)$-matching to a block $W^\prime$ in $\pi_{i,j}R_A^\prime$ with 
         ${\rm Leb}(W^\prime)>6\varepsilon$, which corresponds to part of the subtree labeled $B(m_0)$ of mass ${\rm Leb}(U^\prime)$. Since the sizes of both $U^\prime$ and $W^\prime$ differ from $f_2(m_0)={\rm Leb}(U)$ by at most
         $\varepsilon$, we have $|{\rm Leb}(W^\prime)-{\rm Leb}(U^\prime)|\le 2\varepsilon$, and $i,j\in B(m_0)$ can be chosen so that $W^\prime$ has an atom location in $S_A(R_A^\prime)$. Similarly,
         $U^\prime$, shifted by $f_2^\prime(m_0)$, is a block of $\pi_{i,j}R_A^\prime$ for any $i,j\in C(m_0)$, and we can choose $i,j$ and a block with atom location in $S_A(R_A)$ so that 
         $(W,U^\prime+f_2^\prime(m_0))$ is a pair in this $(i,j)$-matching. Then the masses of $W$ and $W^\prime$ differ by at most
         $3\varepsilon$. Hence removing from the $(b(m_0),c(m_0))$-matching pairs involving top masses and adding the pair $(W,W^\prime)$ increases the unmatched mass by at most 
         $\max\{2\varepsilon+3\varepsilon,2\varepsilon+2\varepsilon\}=5\varepsilon$. Proceeding as for type III yields total unmatched mass bounded above by $(3\#A-1)\varepsilon$. 
      \item[(ii)] Suppose $f_1(m_0)\in(0,2\varepsilon]$, $f_1^\prime(m_0)\in(3\varepsilon,\infty)$, $f_2(m_0),f_2^\prime(m_0)\in(2\varepsilon,\infty)$. The other subcase follows again by symmetry. With notation as in (i), $(U,U^\prime)$ is 
          again a pair in the $(b(m_0),c(m_0))$-matching and the top mass of size $f_1(m_0)$ is unmatched. Here, the second top mass $V^\prime$ of $\pi_{b(m_0),c(m_0)}R_A^\prime$ of size $f_2^\prime(m_0)>2\varepsilon$ will be matched to a block
          $V$ in $\pi_{b(m_0),c(m_0)}R_A$, of size at least $\varepsilon$. Now consider $i,j\in C(m_0)$. Then $U^\prime+{\rm Leb}(V^\prime)$, is a block in the $(i,j)$-matching that must be 
          matched to $V$: on the one hand $f_2(m_0)<f_1^\prime(m_0)+\varepsilon$ since $(U,U^\prime)$ is a pair in a matching; on the other hand,
          if $U^\prime+{\rm Leb}(V^\prime)$ were matched to the left of $V$, it would be matched to a part of $f_2(m_0)$, so $f_2(m_0)>f_1^\prime(m_0)+f_2^\prime(m_0)-\varepsilon$; but then $f_2^\prime(m_0)<2\varepsilon$, which is false. In particular, we find that ${\rm Leb}(V)>2\varepsilon$. Also, for $i,j\in B(m_0)$, $(V,V^\prime)$ is in the $(i,j)$-matching, by Lemma \ref{lm:match}, and $U$ is a block in the 
         $(i,j)$-matching that must be matched to a block $W^\prime$ to the left of $V^\prime$. In particular, ${\rm Leb}(U^\prime)\ge{\rm Leb}(W^\prime)\ge{\rm Leb}(U^\prime)-2\varepsilon$.

         If ${\rm Leb}(U^\prime)\le 4\varepsilon$, then the matched pairs identified above imply that ${\rm Leb}(U)\le 5\varepsilon$, ${\rm Leb}(V)\le 5\varepsilon$ and 
         ${\rm Leb}(V^\prime)\le 6\varepsilon$. If we remove both pairs $(U,U^\prime)$ and $(V,V^\prime)$ from the $(b(m_0),c(m_0))$-matching, we increase the distortion by at most 
         $\max\{5\varepsilon+5\varepsilon,6\varepsilon+4\varepsilon\}=10\varepsilon$. We can argue as in case (i) that the resulting total distortion is at most $(3\#A-1)\varepsilon$.
         
        If ${\rm Leb}(U^\prime)>4\varepsilon$, then ${\rm Leb}(V)>3\varepsilon$, so $V$ has an atom location in $S_A(R_A)$. Also, the constraint 
        ${\rm Leb}(W^\prime)\ge{\rm Leb}(U^\prime)-2\varepsilon$ makes $W^\prime$ unique, for any $i,j\in B(m_0)$, and we may choose $i,j\in B(m_0)$ so that $W^\prime$ has an atom location in $S_A(R_A^\prime)$.
        Replacing $U^\prime$ by $W^\prime$ increases the unmatched mass by at most $2\varepsilon$. Since the $(i,j)$-matchings $M_{i,j}$, $i,j\in C(m_0)$, have matched blocks 
        $V$, and $U^\prime+{\rm Leb}(V^\prime)$, Lemma \ref{lm:match} ensures that they induce $(i,j)$-matchings 
        $M_{i,j}^{C(m_0)}$ of $\pi_{i,j}R_{C(m_0)}$ and $\pi_{i,j}R_{C(m_0)}^\prime$ of distortion at most $\varepsilon$. By the induction hypothesis, there is an associated finite collection of pairs of
        blocks that leave at most mass $(3\#C(m_0)-1)\varepsilon$ unmatched. Proceeding as for type III yields total unmatched mass of at most
        $\varepsilon+(3\#C(m_0)-1)\varepsilon+2\varepsilon+6\varepsilon m_2\le (3\#A-1)\varepsilon$.  
     \item[(iii)] Suppose $f_2(m_0)\in[0,2\varepsilon]$, $f_2^\prime(m_0)\in(3\varepsilon,\infty)$, $f_1(m_0),f_1^\prime(m_0)\in(2\varepsilon,\infty)$. We adapt the argument of (ii), as follows. Again,
        the first two pairs $(U,U^\prime)$ and $(V,V^\prime)$ of the $(b(m_0),c(m_0))$-matching consist of the three large top masses and the non-top mass $V$. In a suitable $(i,j)$-matching for
        $i,j\in C(m_0)$ we find $V$ matched with $U^\prime+{\rm Leb}(V^\prime)$, while $U+f_2(m_0)$ is matched with a block $W^{\prime\prime}$ to the left of $U^\prime+{\rm Leb}(V^\prime)$,
        where $W^{\prime\prime}$ has an atom location in $S_A(R_A^\prime)$.

        If ${\rm Leb}(V^\prime)\le 4\varepsilon$, then ${\rm Leb}(V)\le 5\varepsilon$, and also ${\rm Leb}(W^{\prime\prime})\le 4\varepsilon$ so that ${\rm Leb}(U)\le 5\varepsilon$ and ${\rm Leb}(U^\prime)<6\varepsilon$. If ${\rm Leb}(V^\prime)>4\varepsilon$, then ${\rm Leb}(U^\prime)+{\rm Leb}(V^\prime)\le{\rm Leb}(U)+{\rm Leb}(V)+\varepsilon\le{\rm Leb}(W^{\prime\prime})+{\rm Leb}(U^\prime)+2\varepsilon$, and ${\rm Leb}(V^\prime)-{\rm Leb}(W^{\prime\prime})\le 2\varepsilon$. In particular, ${\rm Leb}(W^{\prime\prime})>2\varepsilon$. The remainder of the argument is easily
 adapted, here applying the induction hypothesis to $R_{B(m_0)}$, $R_{B(m_0)}^\prime$.
    \item[(iv)] Suppose $f_1(m)\in(0,2\varepsilon]$ and $f_2(m)\in[0,2\varepsilon]$ and $f_1^\prime(m),f_2^\prime(m)\in(3\varepsilon,\infty)$. We combine the arguments of (ii) and (iii). Indeed, we 
      consider the first two pairs $(U,U^\prime)$ and $(V,V^\prime)$ of the $(b(m_0),c(m_0))$-matching, with neither $U$ nor $V$ top mass. When ${\rm Leb}(U^\prime)\le 4\varepsilon$ or
      ${\rm Leb}(V^\prime)\le 4\varepsilon$, we can drop both pairs and increase the distortion by at most $10\varepsilon$. If both ${\rm Leb}(U^\prime)>4\varepsilon$ and 
      ${\rm Leb}(V^\prime)>4\varepsilon$, we can replace $U^\prime$ and $V^\prime$ by $W^\prime$ and $W^{\prime\prime}$ and increase the unmatched mass by at most $2\varepsilon+2\varepsilon=4\varepsilon$. In either
      case, the type-III argument completes this last case hence completing the induction step.
   \end{enumerate}
  For $A=[k]$, we obtain a partial coupling of atom mass, which gives a finite collection $K_0\subset S_k^\circ(R_k)\times S_k^\circ(R_k^\prime)$ of pairs of atoms. To build a  
  ${\rm GH}$-correspondence
  between $S_k^\circ( R_k)$ and $S_k^\circ( R_k^\prime)$ we define a set $K_1$ of special pairs that contains
  \begin{itemize}
  \item pairs of vertices $(b_E,b_E^\prime)$ for all $E\in\ft_k$
    \item and the pairs of coupled atoms $(v,v^\prime)\in K_0$.
  \end{itemize}
  We now define the ${\rm GH}$-correspondence $K$ generated by $K_1$ as containing
  \begin{itemize}
    \item pairs of points $\big(\lambda v+(1-\lambda)\overset{\leftarrow}{v},\lambda v^\prime+(1-\lambda)\big(\overset{\leftarrow}{v}\big)^\prime\big)$, 
    $0\le\lambda\le 1$, on the branch between adjacent corresponding special points $(v,v^\prime),\big(\overset{\leftarrow}{v},\big(\overset{\leftarrow}{v}\big)^\prime\big)\in K_1$. 
  \end{itemize} 
  We now proceed as in the proof of Proposition \ref{propGH}, where we note that \eqref{eq:appxdIbound} ensures that distances from the root for pairs of 
  matched atoms are still bounded by $\varepsilon$, except where we included $(W,W^\prime)$, $(U,W^\prime)$, $(V,W^\prime)$ and/or $(V,W^{\prime\prime})$ for type-V blocks. In those cases, which
  only arise for $k\ge 4$, the bound is $3\varepsilon$. We can therefore extend \eqref{eq:spinedist} to a maximum over our extended list of special pairs, subject to a factor 3 when $k\ge 4$. 
  Then the ${\rm GH}$-distortion of $K$ is at most $4\varepsilon\le 3k\varepsilon$ when $k=2,3$ and at most $12\varepsilon\le 3k\varepsilon$ when $k\ge 4$, by the same argument as in the proof of 
  Proposition \ref{propGH}. Furthermore, the ${\rm GH}$-correspondence gives rise to a partial coupling of the
  weight measures $\mu$ of $S_k(R_k)$ and $\mu^\prime$ of $S_k(R_k^\prime)$, which is given by
  $$\nu=\sum_{(U,U^\prime)}\min\{{\rm Leb}(U),{\rm Leb}(U^\prime)\}\delta(v_U,v_{U^\prime}^\prime),$$
  where the sum is taken over all matched blocks $(U,U^\prime)$ chosen for the construction of $K_0$, and where $v_U\in S_k^\circ(R_k)$ and 
  $v_{U^\prime}^\prime\in S_k^\circ(R_k^\prime)$ are the atom locations corresponding to these blocks. Specifically, the fact that the coupling of atom mass leaves mass at most
  $(3k-1)\varepsilon$ uncovered now implies that
  $$\mu(S_k^\circ( R_k))-\nu(K_0)\le(3k-1)\varepsilon\qquad\mbox{and}\qquad\mu^\prime(S_k^\circ( R_k^\prime))-\nu(K_0)\le(3k-1)\varepsilon.$$
  By the argument of \cite[Proposition 6]{M09}, adapted from the case of probability measures to the case of finite measures, this permits the construction of a   
  metric space $(M,d_M)$ and injective isometries $\phi\colon S_k^\circ(R_k)\rightarrow M$ and $\phi^\prime\colon S_k^\circ(R_k^\prime)\rightarrow M$ that
  show that $d_{\rm GHP}(S_k(R_k),S_k(R_k^\prime))\le 3k\varepsilon$, as required. 
\end{proof}

\section[Embedding of $\mathbb{R}$-trees and the proof of Proposition 8.2]{Embedding of $\mathbb{R}$-trees and the proof of Proposition \ref{prop:embcont}}\label{appx:markedGH}

\begin{lemma}\label{lm:embconv} Consider a sequence $\mathrm{T}_n=(T_n,d_n,\rho_n)$, $n\ge 0$, of compact rooted $\mathbb{R}$-trees and compact connected subsets $S_n\subseteq T_n$ with $\rho_n\in S_n$, $n\ge 1$. Suppose that $d_{\rm GH}(\mathrm{T}_n,\mathrm{T}_0)\rightarrow 0$ as $n\rightarrow\infty$ and that $\mathrm{S}_n=(S_n,d_n,\rho_n)$, $n\ge 1$, converges in $(\mathbb{T}_\circ^{\rm real},d_{\rm GH})$. Then there is a compact connected subset $S_0\subseteq T_0$ with $\rho_0\in S_0$ such that $d_{\rm GH}(\mathrm{S}_n,\mathrm{S}_0)\rightarrow 0$, where $\mathrm{S}_0=(S_0,d_0,\rho_0)$.
\end{lemma}
\begin{proof}
Let $(E_n,\delta_n)$ be a sequence of metric spaces and $\phi_n, \psi_n$ sequences of isometric embeddings of $T_n$ and $T_0$ into $E_n$ such that $\delta_n^{\rm H}(\phi_n(T_n),\psi_n(T_0)) \to 0$ and $\delta_n(\phi_n(\rho_n),\psi_n(\rho_0))\to 0$.  Let $a_n =\max\{\delta_n^{\rm H}(\phi_n(T_n),\psi_n(T_0)),\delta_n(\phi_n(\rho_n),\psi_n(\rho_0))\}+1/n$ and define $\widebar{S}_n = \{t\in T_0\colon \delta_n(\psi_n(t),\phi_n(S_n))\le a_n\}$ and $\widebar{\mathrm{S}}_n=(\widebar{S}_n,d_0,\rho_0)$.  By construction, $\rho_0\in\widebar{S}_n$ and $d_{\rm GH}(\mathrm{S}_n,\widebar{\mathrm{S}}_n) \leq a_n\to 0$.  

Since $\mathrm{T}_0$ is compact, the Hausdorff topology on compact subsets of $T_0$ is compact as well.  Therefore there is a subsequence $\widebar{S}_{n_k}$ converging to some compact 
$S_0\subseteq T_0$ in the Hausdorff metric.  Hence $d_{\rm GH}(\mathrm{S}_{n_k},\mathrm{S}_0) \to 0$.  However, since $\mathrm{S}_n$ converges in the GH topology, $d_{\rm GH}(\mathrm{S}_n,\mathrm{S}_0)\to 0$ as well.

Finally, since the set of isometry classes of $\mathbb{R}$-trees is $d_{\rm GH}$-closed, $\mathrm{S}_0$ is an $\mathbb{R}$-tree and hence connected.
\end{proof}

\begin{corollary}\label{cor:embconv} Suppose that in the setting of the Lemma \ref{lm:embconv}, we have compact connected $S_n^{(k)}\subseteq T_n$ with $\rho_n\in S_n^{(k)}$ for all $k\ge 1$, and that they are nested,
  i.e. $S_n^{(k)}\subseteq S_n^{(k+1)}$ for all $k\ge 1$. If $\mathrm{S}_n^{(k)}=(S_n^{(k)},d_n,\rho_n)$ converges in $(\mathbb{T}_\circ^{\rm real},d_{\rm GH})$ for all $k\ge 1$, then there is a
  nested family of compact connected $S_0^{(k)}\subseteq T_0$ with $\rho_0\in S_0^{(k)}$, $k\ge 1$, such that $d_{\rm GH}(\mathrm{S}_n^{(k)},\mathrm{S}_0^{(k)})\rightarrow 0$ for all $k\ge 1$. 
\end{corollary}
\begin{proof} We argue as in the proof of the lemma, for $k=1$ to find a subsequence along which $\widebar{S}_n^{(1)}$, in obvious notation, converges to some $S_0^{(1)}\subseteq T_0$. For each 
  $k\ge 2$, we inductively pass to a further subsequence to identify subsequential limits $S_0^{(k)}\subseteq T_0$ of $\widebar{S}_n^{(k)}$ while maintaining previous convergences. By definition 
  $S_n^{(k)}\subseteq S_n^{(k+1)}$ implies $\widebar{S}_n^{(k)}\subseteq\widebar{S}_n^{(k+1)}$ and hence $S_0^{(k)}\subseteq S_0^{(k+1)}$, as claimed. The claimed Gromov--Hausdorff 
  convergences follow as in the proof of Lemma \ref{lm:embconv}.
\end{proof}

We remark that while the Hausdorff limit of $S_0^{(k)}$ as $k\rightarrow\infty$ is well-defined as a subset of $T_0$, it may be a strict subset of $T_0$ even if all $T_n$ are the Hausdorff limits of 
$S_n^{(k)}$ as $k\rightarrow\infty$. This is because there may be subtrees of $T_0$ that correspond to subtrees of $T_n$ that are included in $S_n^{(k_n+1)}$, but do not intersect $S_n^{(k_n)}$ for a 
sequence $k_n\rightarrow\infty$. 

\begin{proof}[Proof of Proposition \ref{prop:embcont}] Denote by $D=\{k2^{-m},k\ge 0,m\ge 0\}$ the set of dyadic rationals. Recall the 
  Definition \ref{def:S} of the map $S$. First consider the setting of the self-similar Aldous diffusion in Section \ref{sec:continuity} and let $z\in D$. The map $S$ takes a Gromov--Hausdorff--Prokhorov limit of rooted, weighted $\mathbb{R}$-trees 
  $S_k(\cT_{k,+}^z)$, which in turn are trees with a discrete tree shape and edge lengths, further equipped with a measure. By consistency of the 
  $k$-tree evolutions, these trees $S_k(\cT_{k,+}^z)$, $k\ge 1$, can be constructed as a projectively consistent sequence of weighted 
  $\mathbb{R}$-trees, for instance embedded in $\ell_1(\mathbb{N})$ as proposed by Aldous \cite{AldousCRT1}, and this construction also yields a representative of 
  the limit $S\big(\cT_{k,+}^z,k\ge 1\big)$, $z\in D$, in which $\tau(\cT_{k,+}^z)$, $k\ge 1$, naturally have embedded representatives that form a nested family.
  By the proof of the Kolmogorov--Chentsov theorem \cite[Proof of Theorem I.(2.1)]{RevuzYor}, there is an almost sure event on which we may define 
  \[
  \cT(y)=\lim_{z\rightarrow y, z\in D}S\big(\cT_{k,+}^z,k\ge 1\big),\qquad y\ge 0,
  \]
  to obtain a path-continuous process. These limits are in the Gromov--Hausdorff--Prokhorov sense. We work on this almost sure event. Let $k\ge 1$ and
  $y\in[0,\infty)$. Then there is a sequence $z_n\in D$, $n\ge 1$, with $z_n\downarrow y$. Recall that the resampling $k$-tree evolution 
  $(\cT_{k,+}^y,y\ge 0)$ is right-continuous and continuous between the resampling times, which form a strictly increasing sequence that accumulates at the
  random time when the total mass vanishes. In particular, we may assume without loss of generality that the sequence $(z_n,n\ge 1)$ is between 
  two resampling times (but allowing $y$ to be a resampling time). For the purposes of the remainder of this proof, we use notation $\tau^\circ=\pi\circ\tau$, where 
  $\pi\colon\bT^{\rm real}\rightarrow\bT_\circ^{\rm real}$ is the natural projection that associates with a (rooted, weighted isometry class of a) rooted, weighted $\mathbb{R}$-tree the (rooted isometry class
  of the) rooted $\mathbb{R}$-tree. Since diversities of edge partitions evolve continuously between resampling times, we deduce that $\tau^\circ(\cT_{k,+}^{z_n})$
  tends to $\tau^\circ(\cT_{k,+}^y)$ in the Gromov--Hausdorff sense.

  We can now apply Lemma \ref{lm:embconv} taking as $\mathrm{S}_n$ the representative of $\tau^\circ(\cT_{k,+}^{z_n})$ embedded into a
  representative $\mathrm{T}_n$ of $\cT(z_n)$, and $u_0=s$. Then the lemma entails that $\tau^\circ(\cT_{k,+}^y)$ can be isometrically 
  embedded into (any representative of) $\cT(y)$. By Corollary \ref{cor:embconv}, this can be achieved for all $k\ge 1$ simultaneously, as a nested family, as required. 

  Since the Aldous diffusion is obtained by scaling and time-changing the self-similar Aldous diffusion, the conclusion also holds for the Aldous diffusion.    
\end{proof}

\section[Proof of Lemma 8.7]{Proof of Lemma \ref{lm:smp1}}\label{appx:smp}

Recall that we defined $A_1$ to be the set of (${\rm GHP}$-isometry classes of) unit-mass rooted, weighted $\mathbb{R}$-trees $(T,d,\rho,\mu)$ that have two branch points $v,w\in T$ at which $T$ can
be decomposed into connected subsets $C_0,C_1,C_2,C_3$ of $T\setminus]\!]v,w[\![$ in such a way that $C_0\cap C_1=\{v\}$, $C_2\cap C_3=\{w\}$, where $C_0$ contains the root $\rho$, and such
that the heights and masses of $C_j$, $j=0,1,2.3$, are all greater than or equal to some threshold values. Specifically, the connectedness and intersection properties imply that $C_j$ can be viewed as an
$\mathbb{R}$-tree rooted at $u_j=v$ for $j=0,1$ and at $u_j=w$ for $j=2,3$, and the height constraints are $\sup\{d(u_j,x)\colon x\in C_j\}\ge 1$, $j=0,1,2,3$, while the mass constraints are
\[
\textstyle\mu(C_0)\ge\frac{12}{37},\ \mu(C_1)\ge\frac{10}{37},\ \mu(C_2)\ge\frac{8}{37},\ \mu(C_3)\ge\frac{6}{37}.
\] 
In this setting, Lemma \ref{lm:smp1}(i) states the following in the case of $A_1$.

\begin{lemma}\label{lm1:smp1} If in the above setting $\mu$ is diffuse, then $v,w,C_0,C_1,C_2,C_3$ are unique. 
\end{lemma}
\begin{proof} Let us fix a choice of $v,w,C_0,C_1,C_2,C_3$. Since $\mu$ is diffuse, there is no atom in $v$ or $w$. The mass decomposition around $v$ includes three 
  components that are subject to mass thresholds $\frac{10}{37}$, $\frac{12}{37}$ and $\frac{14}{37}$, and as there is no atom in $v$ or $w$, none of these is exceeded 
  by more than $\frac{1}{37}$, so the components cannot swap roles and none of them is big enough to contain two sufficiently heavy components in a decomposition around another branch point that also 
  exceeds these mass thresholds. Therefore, $v,C_0,C_1$ are unique. The further 
  branch point $w$ splits the component of mass between $\frac{14}{37}$ and $\frac{15}{37}$ further into two components of sizes at least $\frac{6}{37}$ and $\frac{8}{37}$, and this is possible only in this 
  component, not in any other component around either $v$ or the present $w$. Hence, $w,C_2,C_3$ are also unique.    
\end{proof}

We now consider three of the claims to Lemma \ref{lm:smp1}(ii)--(iii), and start by studying
  \begin{align*}
A_1^\circ=\Big\{\mathrm{T}\in\mathbb{T}^{\rm real}_{\rm unit}\colon\exists_{\overset{v,w\in T}{v\neq w}}&\exists_{\overset{U_0,U_1,U_2,U_3\subseteq T\setminus[\![v,w]\!]}{\text{open connected}}}\widebar{U}_0\cap\widebar{U}_1=\{v\},\widebar{U}_2\cap\widebar{U}_3=\{w\}, \rho\in U_0,\\
 &\ \ \forall_{j}\ \mu(U_j)>\textstyle\frac{12-2j}{37},
   \ \sup\{d(v,x)\wedge d(w,x)\colon x\in U_j\}>1\Big\}.
  \end{align*}
\begin{lemma}\label{lm2:smp1} The set $A_1^\circ$ is open in $(\mathbb{T}^{\rm real}_{\rm unit},d_{\rm GHP})$
\end{lemma}
\begin{proof}
  Let $\mathrm{T}\in A_1^\circ$ be the weighted isometry class of $(T,d,\rho,\mu)$ and suppose that $v,w,U_0,U_1,U_2,U_3$ are such that the constraints for membership in $A_1^\circ$ hold. We will write
  $u_j=v$ for $j=0,1$ and $u_j=w$ for $j=2,3$. By possibly making the $U_j$, $j=0,1,2,3$, larger, we may assume that they are entire connected components of $T\setminus[\![v,w]\!]$, since the
  constraints in combination with the $\mathbb{R}$-tree property of $(T,d)$ already guarantee that they are subsets of distinct connected components. 

  We denote by $B_\delta(v)=\{x\in T\colon d(x,v)<\delta\}$ the open ball of radius $\delta$. By the regularity of measures, there is $\delta>0$ such that $d(v,w)>10\delta$ and the compact sets $V_j:=U_j\setminus B_{2\delta}(u_j)$, $j=0,1,2,3$, still exceed the 
  respective mass thresholds by more than $2\delta$ and stricter height thresholds of $1+2\delta$. By compactness, only finitely many connected components of $U_j\setminus\overline{B_\delta(u_j)}$
  intersect $V_j$. As $U_j$ is connected, the most recent common ancestor $v_j$ of $V_j$ is in $U_j$. We denote by $W_j=\{x\in U_j\colon v_j\in[\![u_j,x]\!]\}$ the subtree above 
  $v_j$ and also set $\varepsilon:=\min\{d(u_j,v_j),j=0,1,2,3\}$. Then $\varepsilon\in(0,2\delta]$ and so $W_j$ is a closed connected subset of $U_j$ that contains $V_j$ and satisfies the stricter mass and height constraints.   
  We note that any pair of points in different subsets $W_j$ has distance at least $2\varepsilon$ in $(T,d)$. 

  Now suppose that the weighted isometry class $\mathrm{T}^\prime\in\mathbb{T}^{\rm real}_{\rm unit}$ of $(T^\prime,d^\prime,\rho^\prime,\mu^\prime)$ satisfies 
  $d_{\rm GHP}(\mathrm{T},\mathrm{T}^\prime)<\varepsilon/12$. Recall Miermont's representation \eqref{eq:ghpunit2} of $d_{\rm GHP}$. Consider a correspondence with ${\rm GH}$-distortion below $\varepsilon/12$ and a coupling   of $\mu$ and $\mu^\prime$ that assigns
  all but mass strictly below $\varepsilon/12$ to pairs of points in the correspondence. We define the sets $W_j^\prime$ of points in 
  $T^\prime$ that are in correspondence with points in $W_j$, $j=0,1,2,3$.  Then the lower bound on the distances of $W_j$ implies that they are disjoint and at distance at least $5\varepsilon/3$. We claim that while $W_j^\prime$ may not be 
  connected, the $\varepsilon/6$-thickening $(W_j^\prime)^{\varepsilon/6}$ will be connected (and at distance at least $4\varepsilon/3$ from each other, and they will satisfy the respective mass and
  height constraints). This is because any two points $x_1^\prime,x_2^\prime\in W_j^\prime$ are in correspondence with two points $x_1,x_2\in W_j$ that are connected by a path 
  $[\![x_1,x_2]\!]\subseteq W_j$. Each point on the path is corresponding to a point in $W_j^\prime$, so travelling along $[\![x_1,x_2]\!]$ at step sizes of $\gamma$ corresponds to a
  sequence of points in $W_j^\prime$ at step size strictly below $\gamma+\varepsilon/6$, for any $\gamma>0$, and for $\gamma$ sufficiently small, this will be below $\varepsilon/3$, so the unique path in 
  the $\mathbb{R}$-tree $T^\prime$ between any two adjacent steps is in $(W_j^\prime)^{\varepsilon/6}$.

  Let us show that there are no $y_j^\prime\in(W_j^\prime)^{\varepsilon/6}$, $j=0,2,3$, such that $y_2^\prime\in[\![y_0^\prime,y_3^\prime]\!]$. Assume for contradiction that this were the case. For corresponding points
  $y_j\in T$, we have $y_j\in(W_j)^{\varepsilon/3}\subseteq U_j\setminus B_{\varepsilon/2}(u_j)$, and this does give a contradiction:
  \begin{align*}
  \frac{\varepsilon}{2}\le d(y_2,w)&=\frac{1}{2}\big(d(y_0,y_2)+d(y_2,y_3)-d(y_0,y_3)\big)\\
                                     &<\frac{1}{2}\Big(d^\prime(y_0^\prime,y_2^\prime)+d^\prime(y_2^\prime,y_3^\prime)-d^\prime(y_0^\prime,y_3^\prime)+3\frac{\varepsilon}{6}\Big)=\frac{\varepsilon}{4}.
  \end{align*}
  Similarly, we cannot have $y_i^\prime\in[\![y_j^\prime,y_k^\prime]\!]$ for any distinct $i,j,k\in\{0,1,2,3\}$. Hence, there 
  are unique $v^\prime,w^\prime\in T^\prime$ such that 
  \begin{align*}
  [\![y_0^\prime,v^\prime]\!]&=[\![y_0^\prime,y_1^\prime]\!]\cap[\![y_0^\prime,y_2^\prime]\!]\quad\mbox{for all }y_j^\prime\in(W_j^\prime)^{\varepsilon/6},j=0,1,2,\\
  [\![y_0^\prime,w^\prime]\!]&=[\![y_0^\prime,y_2^\prime]\!]\cap[\![y_0^\prime,y_3^\prime]\!]\quad\mbox{for all }y_j^\prime\in(W_j^\prime)^{\varepsilon/6},j=0,2,3.
  \end{align*}
  Note that $v_j\in W_j\cap\overline{B_{2\delta}(u_j)}$, $j=0,1,2,3$. In particular, corresponding $v_j^\prime\in W_j^\prime$ satisfy $d^\prime(v_0^\prime,v_1^\prime)\le 4\delta+\varepsilon/6<5\delta$ and similarly $d^\prime(v_2^\prime,v_3^\prime)<5\delta$. On the other hand,
$d^\prime(v_1^\prime,v_j^\prime)\ge d(v_1,v_j)-\varepsilon/6\ge d(v,w)+11\varepsilon/6>10\delta$ for $j=2,3$, and similarly
$d^\prime(v_0^\prime,v_j^\prime)>10\delta$. An elementary argument considering the possible shapes of the tree spanned by $v_j^\prime$, $j=0,1,2,3$,
entails that $v^\prime\in]\!]v_j^\prime,w^\prime[\![$ for $j=0,1$ and $w^\prime\in]\!]v^\prime,v_j^\prime[\![$ for $j=2,3$.   We conclude that the
  connected components $U_j^\prime$ of $T^\prime\setminus[\![v^\prime,w^\prime]\!]$ containing $W_j^\prime$, $j=0,1,2,3$, together with $v^\prime$ and $w^\prime$ satisfy all the 
  constraints to imply that $\mathrm{T}^\prime\in A_1^\circ$, as required.
\end{proof}
\begin{lemma}\label{lm3:smp1} The set $A_1$ is closed.
\end{lemma}
\begin{proof} Consider a sequence $(\mathrm{T}_n,n\ge 1)$ in $A_1$ that converges in $(\mathbb{T}^{\rm real},d_{\rm GHP})$. By \cite[Lemmas 5.8 and A.1]{GrevPfafWint09}, there is a compact
  metric space $(M,d_M)$ and embeddings $T_n\subseteq M$, $n\ge 1$, such that representatives $(T_n,d_M,\rho_n,\mu_n)$ of $\mathrm{T}_n$ converge in the sense that 
  \[
  d_M(\rho_n,\rho)\rightarrow 0,\qquad d_M^{\rm H}(T_n,T)\rightarrow 0,\qquad\mbox{and }d_M^{\rm P}(\mu_n,\mu)\rightarrow 0.
  \]
  Since $\mathrm{T}_n\in A_1$, we may take $v_n,w_n\in T_n$ and $C_j^{(n)}\subset T_n$, $j=0,1,2,3$, satisfying all constraints for membership in $A_1$, for all $n\ge 1$. Since the space of compact 
  subsets of $M$ equipped with the Hausdorff distance $d_M^{\rm H}$ is compact, we may assume, by taking successive subsequences that we suppress notationally, that also
  \[
  d_M(v_n,v)\rightarrow 0,\qquad d_M(w_n,w)\rightarrow 0,\qquad d_M^{\rm H}(C_j^{(n)},C_j)\rightarrow 0,\quad j=0,1,2,3,
  \]
  for some $v,w\in T$, $C_j\subset T$, $j=0,1,2,3$. Then $C_j$, $j=0,1,2,3$, are connected as Hausdorff limits of $\mathbb{R}$-trees. Also clearly $v\in C_0\cap C_1$ and $w\in C_2\cap C_3$. To show 
  that these intersections contain no further points, assume for contradiction that $x\in C_0\cap C_1$ with $x\neq v$. Then there are sequences $x_n\in C_0^{(n)}$ and $x_n^\prime\in C_1^{(n)}$ such 
  that $d_M(x_n,x)\rightarrow 0$ and $d_M(x_n^\prime,x)\rightarrow 0$. But also, $\delta=d(v,x)>0$. We write $u_j^{(n)}=v_n$, $j=0,1$, and $u_j^{(n)}=w_n$, $j=2,3$, $n\ge 1$. Now consider 
  $\widetilde{C}_j^{(n)}:=C_j^{(n)}\setminus B_{\delta/2}(u_j^{(n)})$, $j=0,1,2,3$, $n\ge 1$. Then passing to further subsequences, we may assume that 
  \[ 
  d_M^{\rm H}(\widetilde{C}_j^{(n)},\widetilde{C}_{j})\rightarrow 0,\qquad j=0,1,2,3,
  \]
  for some $\widetilde{C}_j\subseteq C_j$, $j=0,1,2,3$. Furthermore, any two points in distinct $\widetilde{C}_j^{(n)}$, $j=0,1,2,3$, are at least 
  $\delta$ apart and this property is maintained in the limit. But we will have $x_n\in\widetilde{C}_0^{(n)}$
  for $n$ sufficiently large, hence $x\in\widetilde{C}_{0}$. Similarly $x\in \widetilde{C}_{1}$, but as $\widetilde{C}_{0}\cap\widetilde{C}_{1}=\emptyset$, this is a contradiction. So $C_0\cap C_1=\{v\}$. 
  Similarly $C_2\cap C_3=\{w\}$.
  That $C_j$ satisfies the height constraint is a consequence of Hausdorff convergence as points in $C_j^{(n)}$ at distance 1 from $u_j^{(n)}$ have limit points in $C_j$ at distance 1 from $u_j$, where 
  $u_j=v$ for $j=0,1$ and $u_j=w$ for $j=2,3$. That $C_j$ satisfies the mass constraint is a consequence of Prokhorov convergence: first for all $\varepsilon$ and $n$ sufficiently large,
  \[
  \mu(C_j^\varepsilon)\ge\mu_n((C_j)^\varepsilon)-\varepsilon\ge\mu_n(C_j^{(n)})-\varepsilon\ge\frac{12-2j}{37}-\varepsilon,
  \]
  then the regularity of measures means the $\lim_{\varepsilon\downto 0}$ of the left-hand side is $\mu(C_j)$ while the right-hand side tends to the required threshold.
\end{proof}
\begin{corollary}\label{lm4:smp1} The closure of $A_1^\circ$ is a subset of $A_1$.
\end{corollary}
\begin{proof} It suffices to note that $A_1^\circ\subset A_1$, and that $A_1$ is closed, by Lemma \ref{lm3:smp1}.
\end{proof}
\begin{proof}[Proof of Lemma \ref{lm:smp1}] 
(i) This was proved for $A_1$ in Lemma \ref{lm1:smp1}, while the argument is easily adapted for $A_2$, $A_3$ and $A$.

(ii)--(iii) We showed in Lemma \ref{lm3:smp1} that $A_1$ is closed, in Lemma \ref{lm2:smp1} that $A_1^\circ$ is open and in Corollary \ref{lm4:smp1} that the closure of $A_1^\circ$ is a subset of $A_1$. Again, 
these arguments are easily adapted for $A_2$ and $A_3$, and the argument of Lemma \ref{lm3:smp1} also to $A$.

(iv) Since $C_j$ and $C_j^\gamma$ are both closed subsets of $T$ for any $\mathrm{T}=[T,d,\rho,\mu]\in\widetilde{A}_i(\gamma)$, the argument of
Lemma \ref{lm3:smp1} again applies to show that $\widetilde{A}_i(\gamma)$ is closed. Given $v$, $w$ and $C_j$, $j=0,1,2,3$, that satisfy the constraints of
$\widetilde{A}_i(\gamma)$, we can consider as $U_j$ the connected component of $T\setminus[\![v,w]\!]$ containing $C_j$, $j=0,1,2,3$. Then all constraints for membership in $A_i^\circ$ hold, with the possible exception of the requirement that $v\neq w$, but if $v=w$, then $\mathrm{T}\in A$. We conclude that $\widetilde{A}_i(\gamma)\subseteq A_i^\circ\cup A$. 

(v) Now let $\mathrm{T}=[T,d,\rho,\mu]\in A_i^\circ$. In the notation of the proof of Lemma \ref{lm2:smp1}, we consider $v_j\in W_j\subseteq U_j$, $j=0,1,2,3$,
and $\varepsilon=\min\{d(u_j,v_j),j=0,1,2,3\}$. This minimum is attained by at least one $j=0,1,2,3$. If we let $\gamma=\varepsilon$ and define 
$\widetilde{v}_j\in[\![u_j,v_j]\!]$ with $d(u_j,\widetilde{v}_j)=\varepsilon$, then $C_j=W_j\cup[\![v_j,\widetilde{v}_j]\!]$ inherits the closure and connectedness properties from $W_j$. This together with the intersection constraints of $U_j$ further entails that $C_j$, $j=0,1,2,3$, satisfy
the intersection constraints for membership in $\widetilde{A}_i(\gamma)$. Since $W_j$ also exceeds the mass and height thresholds of $A_i^\circ$ by more
than $2\delta\ge\varepsilon=\gamma$, in the notation of the proof of Lemma \ref{lm2:smp1}, $C_j$ also satisfies the mass and height constraints for
membership of $\mathrm{T}$ in $\widetilde{A}_i(\gamma)$.      
\end{proof}

\section[Proof of Proposition 8.16]{Proof of Proposition \ref{lm:degrate}}\label{appx:degrate}

By the intertwining argument at the beginning of the proof of Theorem \ref{thm:embedding}, the resampling time $\widebar{D}$ of a unit-mass 2-tree evolution is exponentially distributed. It therefore suffices to identify the rate parameter, which can be obtained as $\lim_{\varepsilon\downarrow 0}\frac{1}{\varepsilon}\mathbb{P}(\overline{D}\le\varepsilon)$. Since de-Poissonization involves a time change that has continuously changing speed starting from 1, Proposition \ref{lm:degrate} follows immediately from the following.

\begin{proposition}\label{prop:2-tree_degen_rate} Consider the degeneration time $D$ of a type-2 evolution starting from the unit-mass pseudo-stationary distribution. Then
  \[
  \lim_{\varepsilon\downarrow 0}\frac{1}{\varepsilon}\mathbb{P}(D\le\varepsilon)=2.
  \]
\end{proposition}

This appendix is devoted to the proof of this proposition as a culmination of intermediate results. Let  $((m_1^y,m_2^y,\beta^y),\, y\ge0)$ be a 
type-2 evolution starting from the unit-mass pseudo-stationary distribution and $D$ its degeneration time.
For the purpose of the following, for $i=1,2$, let $T_i$ denote the first time that the mass on label $i$ approaches 0:
\begin{equation*}
 T_i := \inf\{y\ge0 \colon m_i^{y-} = 0\} = \sup_{h>0} \inf\{y\ge0 \colon m_i^{y} < h\}.
\end{equation*}
For $\epsilon\in \big(0,\frac14\big)$, we cover the event $\{D\le\epsilon\}$ with a union of 6 events.
\begin{align}
 A_{\epsilon,1} &:= \left\{m_1^0 > 1-\sqrt{\epsilon}\text{\ and\ }m_2^\epsilon+\|\beta^\epsilon\| = 0\text{, while\ }T_1\ge \epsilon\right\},\\
 B_{\epsilon,1} &:= \left\{m_1^0 \le 1-\sqrt{\epsilon}\text{\ and\ }m_2^\epsilon+\|\beta^\epsilon\| = 0\text{, while\ }T_1\ge \epsilon\right\},\\
 C_\epsilon &:= \{\max\{T_1,T_2,D\}\le \epsilon\text{\ and\ }m_1^\epsilon + m_2^\epsilon \ge 1-\epsilon^{1/3}\},\\
 D_\epsilon &:= \{m_1^\epsilon + m_2^\epsilon + \|\beta^\epsilon\| < 1-\epsilon^{1/3}\}.
\end{align}
We define $A_{\epsilon,2}$ and $B_{\epsilon,2}$ in the same manner as $A_{\epsilon,1}$ and $B_{\epsilon,1}$, respectively, but with the roles of labels 1 and 2 switched. In particular,  $A_{\epsilon,1}$ is the event that the process degenerates prior to time $\epsilon$ and that label 1 starts with a large initial mass and avoids converging to mass 0 prior to time $\epsilon$. By virtue of our choice that $m_1^0 > 1-\sqrt{\epsilon} > \frac12$, this is disjoint from $A_{\epsilon,2}$. Thus, by symmetry,
\begin{equation}\label{eq:cover_degen_event}
 2\bP(A_{\epsilon,1}) \le \bP\{D\le \epsilon\} \le 2\bP(A_{\epsilon,1}) + 2\bP(B_{\epsilon,1}) + \bP(C_\epsilon) + \bP(D_\epsilon).
\end{equation}

First, we prove the following.

\begin{lemma}\label{lem:PS_degen_rate:A_lim}
 $\displaystyle
 \lim_{\epsilon\to 0+} \frac{\bP(A_{\epsilon,1})}{\epsilon} 
 	= \lim_{\epsilon\to 0+}\frac{1}{\epsilon}\bP\left\{ \begin{array}{c}
 		m_1^0 > 1-\sqrt{\epsilon},\ T_1 \ge \epsilon,\\[3pt]
 		\textnormal{and }m_2^\epsilon + \|\beta^\epsilon\| = 0
 		\end{array}\right\} 
 	= 1$.
\end{lemma}
Afterwards, we will show that each of the remaining probabilities goes to 0 faster than $\epsilon$.

\begin{proof}
 Fix $\epsilon \in (0,\frac14)$.
 The top mass $m_1^y$ evolves as a $\besq(-1)$ up until it approaches 0. If it starts from mass $1-x$ then by Lemma \ref{lem:BESQ:length}, $T_1\sim\distribfont{InvGamma}\big(\frac32,\frac{1-x}{2}\big)$. Conditional on $\{T_1>\epsilon\}$, the process $(m_2^y,\beta^y)$, $y\in [0,\epsilon]$, is a type-1 evolution, per Definition \ref{def:type2:v1}. Thus, by Proposition \ref{type1totalmass}, its total mass evolves as a $\besq(0)$, so its degeneration time $D$ from initial mass $x$ is $\distribfont{InvExp}\big(\frac{x}{2}\big)$ \cite[equation (13)]{GoinYor03}. Combining this with the $\BetaDist\big(1,\frac12\big)$ distribution of $m_2^0+\|\beta^0\|$ gives
 \begin{equation}\label{eq:PS_degen_rate:A_bound}
  \bP(A_{\epsilon,1}) = \int_{x=0}^{\sqrt{\epsilon}}\frac{1}{2\sqrt{1-x}} e^{-x/2\epsilon} \left(1 - \int_{y=(1-x)/2\epsilon}^\infty \frac{1}{\Gamma\left(\frac32\right)} \sqrt{y}e^{-y} dy\right)dx,
 \end{equation}
 with the leftmost term under the outer integral being the density of $\BetaDist\big(1,\frac12\big)$, the middle term being the cumulative distribution function of $\distribfont{InvExp}\big(\frac{x}{2}\big)$, and the inner integral calculating one minus the cumulative distribution function of $\distribfont{InvGamma}\big(\frac32,\frac{1-x}{2}\big)$. To get a lower bound, we bound the first term below by $\frac12$ and reduce the lower bound of the inner integral down to $(1-\sqrt{\epsilon})/2\epsilon$.
 \begin{equation*}
 \begin{split}
  \bP(A_{\epsilon,1}) &\ge \int_{x=0}^{\sqrt{\epsilon}}\frac{1}{2} e^{-x/2\epsilon} dx\left(1 - \int_{y=(1-\sqrt{\epsilon})/2\epsilon}^\infty \frac{2}{\sqrt{\pi}} \sqrt{y}e^{-y} dy\right)\\
  	&\ge [-\epsilon e^{-x/2\epsilon}]_0^{\sqrt{\epsilon}} \left(1 - \frac{2}{\sqrt{\pi}}\int_{y=(1-\sqrt{\epsilon})/2\epsilon}^\infty ye^{-y}dy\right).
 \end{split}
 \end{equation*}
 The bound $y \ge \sqrt{y}$ in the inner integral is justified as $(1-\sqrt{\epsilon})/2\epsilon \ge (1-\frac12)/\frac12 = 1$. Thus,
 \begin{equation*}
 \begin{split}
  \bP(A_{\epsilon,1}) &\ge \epsilon \left(1 - e^{-1/2\sqrt{\epsilon}}\right) \left(1 - \frac{2}{\sqrt{\pi}}[-ye^{-y} - e^{-y}]_{(1-\sqrt{\epsilon})/2\epsilon}^\infty\right)\\
  	&= \epsilon \left(1 - e^{-1/2\sqrt{\epsilon}}\right) \left( 1 - \frac{2}{\sqrt{\pi}} e^{-(1-\sqrt{\epsilon})/2\epsilon}\left(\frac{1-\sqrt{\epsilon}}{2\epsilon} + 1\right) \right).
 \end{split}
 \end{equation*}
 Both terms after the initial $\epsilon$ converge to 1 as epsilon tends to 0.
 
 We now derive the upper bound. Continuing from \eqref{eq:PS_degen_rate:A_bound} and bounding the inner integral below by 0,
 \begin{equation*}
  \bP(A_{\epsilon,1}) \le \int_{x=0}^{\sqrt{\epsilon}}\frac{1}{2\sqrt{1-x}} e^{-x/2\epsilon}dx
  	\le \int_{x=0}^{\infty} \frac{1}{2\sqrt{1-\sqrt{\epsilon}}} e^{-x/2\epsilon}dx
  	= \frac{\epsilon}{\sqrt{1-\sqrt{\epsilon}}}.
 \end{equation*}
 Combining this with our lower bound proves the limit.
\end{proof}

\begin{lemma}\label{lem:PS_degen_rate:B_lim}$\displaystyle
  \lim_{\epsilon\to 0+} \frac{\bP(B_{\epsilon,1})}{\epsilon}
  	= \lim_{\epsilon\to 0+}\frac{1}{\epsilon}\bP\left\{\begin{array}{c}
  		m_1^0 \le 1-\sqrt{\epsilon},\ T_1\ge \epsilon,\\[3pt]
  		\textnormal{and\ }m_2^\epsilon+\|\beta^\epsilon\| = 0
  	\end{array}\right\}
 	= 0.$
\end{lemma}

\begin{proof}
 Fix $\epsilon\in \big(0,\frac14\big)$. We get a formula for $\bP(B_{\epsilon,1})$ via the same argument as that giving rise to \eqref{eq:PS_degen_rate:A_bound}, just changing the bounds on the outer integral:
 \begin{equation}\label{eq:PS_degen_rate:B_bound}
  \bP(B_{\epsilon,1}) = \int_{x=\sqrt{\epsilon}}^1\frac{1}{2\sqrt{1-x}} e^{-x/2\epsilon} \left(1 - \int_{y=(1-x)/2\epsilon}^\infty\frac{1}{\Gamma\left(\frac32\right)} \sqrt{y}e^{-y} dy\right)dx.
 \end{equation}
 Thus,
 \begin{equation*}
 \begin{split}
  \bP(B_{\epsilon,1}) &\le \int_{x=\sqrt{\epsilon}}^1 \frac{1}{2\sqrt{1-x}} e^{-x/2\epsilon}dx\\
  	&\le \int_{x=\sqrt{\epsilon}}^{3/4} \frac{1}{2\sqrt{1/4}} e^{-x/2\epsilon}dx + \int_{x=3/4}^1 \frac{1}{2\sqrt{1-x}} e^{-3/8\epsilon}dx\\
  	&= [-2\epsilon e^{-x/2\epsilon}]_{\sqrt{\epsilon}}^{3/4} + e^{-3/8\epsilon}[-\sqrt{1-x}]_{3/4}^1\\
  	&\le 2\epsilon e^{-1/2\sqrt{\epsilon}} + \frac12 e^{-3/8\epsilon}.
 \end{split}
 \end{equation*}
 Dividing by $\epsilon$ and taking the limit proves the result.
\end{proof}

\begin{lemma}\label{lem:PS_degen_rate:C_lim}
 $\displaystyle
  \lim_{\epsilon\to 0+}\frac{\bP(C_\epsilon)}{\epsilon} = 
 \lim_{\epsilon\to 0+}\frac{1}{\epsilon}\bP\left\{\begin{array}{c}\max\{T_1,T_2,D\}\le \epsilon,\\[3pt]\textnormal{and\ }m_1^\epsilon + m_2^\epsilon \ge 1-\epsilon^{1/3}\end{array}\right\}=0.
 $
\end{lemma}


\begin{proof}
 For $j=1,2$, let $E_j := \{T_j \le \epsilon\}$. Lemma \ref{lem:BESQ:length} entails that
 \begin{equation}
  \label{eq:PS_degen_rate:clock_change}
  \bP(E_j \mid m_j^0) = \bP(G \ge m_j^0 / 2\epsilon\mid m_j^0),
 \end{equation}
 where $G \sim \GammaDist\big(\frac32,1\big)$. Plugging in the $\BetaDist\big(\frac12,1\big)$ distribution of $m_1^0$,
 \begin{equation*}
  \bP(E_1) = \int_0^1 \frac{1}{2\sqrt{x}}\bP\left\{G \ge \frac{x}{2\epsilon}\right\}dx
  	\le \int_0^\epsilon \frac{1}{2\sqrt{x}}dx + \int_\epsilon^1 \frac{1}{2\sqrt{x}} \frac{3\epsilon}{x}dx,
 \end{equation*}
 by Markov's inequality. Evaluating these integrals gives
 \begin{equation}\label{eq:PS_degen_rate:dbl_clock_change}
  \bP(E_1) = \bP(E_2) \le 4\sqrt{\epsilon} - 3\epsilon = O(\sqrt{\epsilon})
 \end{equation}
 in big-$O$ notation.
 
 Note that $E_1$ and $E_2$ are conditionally independent given the initial mass split $(m_1^0,m_2^0,\|\beta^0\|)$. By \eqref{eq:PS_degen_rate:clock_change}, the conditional probabilities of these two events are monotone decreasing in $m_1^0$ and $m_2^0$, respectively. Under the $\distribfont{Dirichlet}\big(\frac12,\frac12,\frac12\big)$ distribution of the initial mass split, there is a strong negative stochastic relationship between these masses: the conditional law of $m_1^0$ given $m_2^0 = a$ stochastically dominates that of $m_1^0$ given $m_2^0 = b$ for any $0\le a < b\le 1$. Thus,
 \begin{equation}
  \bP(E_1\cap E_2) \le \bP(E_1)\bP(E_2) = O(\epsilon).
 \end{equation}
 
 Construction \ref{deftype2} of type-2 evolutions shows that we may view the blocks of the interval partition component $\beta^y$ of our type-2 evolution as a subset of the blocks of a type-0 evolution $(\widehat \beta^z,\,z\ge0)$ with the same initial state, $\widehat\beta^0 = \beta^0$. Moreover, after the first times $T_1$ and $T_2$ when each of the initial top masses $m_1^y$ and $m_2^y$ converges to 0, these top masses also correspond to blocks in this type-0 evolution.
 
 Now, define $E^* := E_1\cap E_2 \cap\{D\le \epsilon\}$. On $E^*$, at most one of the top masses is non-zero at time $\epsilon$, so to bound the probability that the total mass of the process exceeds $1-\epsilon^{1/3}$ at time $\epsilon$ on this event, we need only bound the conditional probability of
 $$E_3 := \big\{\widehat \beta^\epsilon\text{ has a block of mass }1-\epsilon^{1/3}\big\}$$
 given $E_1$ and $E_2$. Note that, by \eqref{eq:PS_degen_rate:clock_change}, under this conditioning, $m_1^0$ and $m_2^0$ are biased to be small, and thus $\|\beta^0\|$ is biased to be large, thus making $E_3$ more probable. Therefore, we will bound this probability in the extreme event that $\|\beta^0\| = 1$.
 
 The total mass process $\big(\big\|\widehat\beta^y\big\|,\,y\ge0\big)$ is a $\besq(1)$, per Proposition \ref{type1totalmass}. Thus, if we take $(B(t),\,t\ge0)$ to denote standard one-dimensional Brownian motion,
 \begin{align*}
  \bP\big(\big\|\widehat\beta^\epsilon\big\| > 1+\epsilon^{1/3}\ \big|\ \big\|\widehat\beta^0\big\| = 1 \big) &\le 2\bP\big(B(\epsilon) > \sqrt{1+\epsilon^{1/3}} - 1\big),\\
  	&\le 2\bP\big(B(\epsilon) > \epsilon^{1/3}/4\big) \le 2 e^{-\epsilon^{-1/3}/32},
 \end{align*}
 bounding the square-root function below by a first-order approximation and applying the Chernoff bound to the Gaussian $B(\epsilon)$. Next, by the pseudo-stationarity of $\big(\widehat\beta^y,\,y\ge0\big)$ as described in Proposition \ref{prop:01:pseudo},
 \begin{align*}
  \bP(E_3 \mid E_1\cap E_2) &\le \bP\left(\frac{\widehat\beta^\epsilon}{\|\widehat\beta^\epsilon\|}\text{ has a block }>\frac{1-\epsilon^{1/3}}{1+\epsilon^{1/3}}\right) + e^{-\epsilon^{-1/3}/32}\\
  	&\le \bP\left(\parbox{5cm}{\centering A size-biased random block from a $\PoiDir(\frac12,\frac12)$ is $> 1-2\epsilon^{1/3}$}\right)\frac{1}{1-2\epsilon^{1/3}} + o(\epsilon).\notag
 \end{align*}
 If we only consider $\epsilon \in (0,2^{-6})$ then $1/(1-2\epsilon^{1/3}) < 2$. A size-biased random block of a $\distribfont{PD}\left(\frac12,\frac12\right)$ has $\BetaDist\big(\frac12,1\big)$ distribution, which has cumulative distribution function $F(x) = \sqrt{x}$ for $x\in [0,1]$. Thus,
 \begin{equation*}
  \bP(E_3 \mid E_1\cap E_2) \le 2(1 - \sqrt{1-2\epsilon^{1/3}}) + o(\epsilon) = O(\epsilon^{1/3}).
 \end{equation*}
 Plugging in \eqref{eq:PS_degen_rate:dbl_clock_change} gives
 $$\bP(C_\epsilon) \le \bP(E_1\cap E_2\cap E_3) \le O(\epsilon) O(\epsilon^{1/3}) = o(\epsilon),$$
 as desired.
\end{proof}

\begin{lemma}\label{lem:PS_degen_rate:D_lim}
 $\displaystyle \lim_{\epsilon\to 0+}\frac{\bP(D_\epsilon)}{\epsilon} = 
 \lim_{\epsilon\to 0+}\frac{1}{\epsilon}{\bP\big(m_1^\epsilon + m_2^\epsilon + \|\beta^\epsilon\| < 1-\epsilon^{1/3}\big)} = 0.$
\end{lemma}

\begin{proof}
 By Theorem \ref{thm:type2:total_mass}, the total mass process $M(y) := m_1^y + m_2^y + \|\beta^y\|$, $y\ge0$, is a $\besq_1(-1)$. There are many ways to bound change in a squared Bessel process; we apply \cite[Lemma 33]{Paper0} with $s=0$, $t=\epsilon < 1$ to get $\mathbb{E}[|M(\epsilon)-M(0)|^p] = O(\epsilon^{p/2})$ for every $p\ge2$. By Markov's inequality,
 $$\bP\big(M(\epsilon) < 1-\epsilon^{1/3}\big) \le \mathbb{E}[|M(\epsilon)-M(0)|^p] / \epsilon^{p/3} = O(\epsilon^{p/6}).$$
 Taking $p=7$, for example, completes the proof.
\end{proof}

We can now prove the main result of this appendix.

\begin{proof}[Proof of Proposition \ref{prop:2-tree_degen_rate}]
 We conclude from Lemmas \ref{lem:PS_degen_rate:A_lim}, \ref{lem:PS_degen_rate:B_lim}, \ref{lem:PS_degen_rate:C_lim}, and \ref{lem:PS_degen_rate:D_lim} together with equation \eqref{eq:cover_degen_event} that
 $$\hfill\lim_{\epsilon\to0+}\epsilon^{-1}\bP\{D \le \epsilon\} = 2.\vspace{-20pt}$$
\end{proof}

\backmatter
\bibliographystyle{abbrv}
\bibliography{AldousDiffusion5_NF_23-05-03,addbib_MW_7May2023}
\printindex

\end{document}
